\documentclass[10pt]{book}
\usepackage{latexsym}
\usepackage{amsfonts}
\usepackage{enumerate}
\usepackage{multicol}
\usepackage{graphicx}
\usepackage{amssymb}
\usepackage{amsmath}
\usepackage{epic}
\usepackage{makeidx}

\usepackage[bookmarks=false]{hyperref}

\newtheorem{thm}{Theorem}[chapter]
\newtheorem{lem}[thm]{Lemma}
\newtheorem{cor}[thm]{Corollary}
\newtheorem{prop}[thm]{Proposition}
\newtheorem{conj}[thm]{Conjecture}

\newtheorem{prob}[thm]{Problem}

\makeindex

\begin{document}

\frontmatter

\thispagestyle{empty}

\begin{center}

\vspace*{.8in}

{\bf {\huge ADDITIVE COMBINATORICS} 

\vspace*{.2in}

{\large A MENU OF RESEARCH PROBLEMS}

\vspace*{1.3in}
{\bf {\large AUGUST, 2017}

\vspace*{2.5in}

B\'{E}LA BAJNOK}

Department of Mathematics

Gettysburg College

Gettysburg, PA 17325-1486 USA

bbajnok@gettysburg.edu}

\end{center}

\newpage

\thispagestyle{empty}

\copyright 2000-17 B\'{e}la Bajnok

\tableofcontents

\chapter{Preface}

{\large {\bf Additive Combinatorics}}

\vspace*{.1in}

This book deals with {\em additive combinatorics}, a vibrant area of current mathematical research.  Additive combinatorics---an offspring of combinatorial number theory and additive number theory---can be described as the study of combinatorial properties of {\em sumsets} (collections of sums with terms from given subsets) in additive structures.  For example, given a subset $A$ in an abelian group $G$, one can consider the sumset $A+A$ consisting of all two-term sums of the elements of $A$, and then ask how small this sumset may be; furthermore, given that $A+A$ is small, one can examine what can be said about the structure of $A$.  

Additive combinatorics is a rather new field within mathematics that is just now coming to its own; although some of its results have been known for a very long time, many of its most fundamental questions have only been settled recently or are still unsolved.  For example, the question about the minimum size of $A+A$ mentioned above was first determined in groups of prime order by Cauchy in\index{Cauchy, A--L.} 1813 (then re-discovered by Davenport\index{Davenport, H.} over a hundred years later), but it has only been determined in the general setting in the twenty-first century.  For this and many other reasons, additive combinatorics provides an excellent area for research by students of any background: it has intriguing and promising questions for everyone.

\vspace*{.25in}

\noindent {\large {\bf Student Research}}

\vspace*{.1in}

Student research in mathematics has been increasing dramatically at all levels over the past several decades.  Producing research is an expectation at all doctoral programs in mathematics, and it appears that the number of publications before students graduate has increased substantially in recent years.  The National Science Foundation (NSF) is sponsoring the popular Research Experiences for Undergraduates (REU) program, and undergraduate research is indeed blossoming at many colleges and universities.  There are even organized mathematics research programs for students in high school, such as PRIMES, held at the Massachusetts Institute of Technology and online nationally, and PROMYS, held at Boston University and at  the University of Oxford. 

The benefits and costs of student research in mathematics have been described elsewhere---see, for example, the report \cite{MAA:2006a}  of the Committee on the Undergraduate Program (CUPM) of the Mathematical Association of America (MAA).  While the need for an ample inventory of questions for student research is clear, it is often noted how challenging it is to produce an appropriate supply.  Indeed, there are quite a few demands on such questions; in the view of this author (see \cite{Baj:2006a}), they are most appropriate when they are:
\begin{itemize}

\item based on substantial topics -- students should be engaged in the study of non-trivial and not-too-esoteric mathematics;

\item challenging at a variety of different levels -- students with different backgrounds and interests should be able to engage in the projects;

\item approachable with a variety of different methods -- students interested in theoretical, computational, abstract, or concrete work should be able to choose their own approaches;

\item incrementally attainable, where at least partial results are within reach yet complete solutions are not easy -- no student would want to spend long hours of hard work and not feel productive; and

\item new and unsolved -- the results attained by the students would have to be, at least in theory, publishable.

\end{itemize}
I believe that the problems in the book abide by these objectives.

 \vspace*{.25in}

\noindent {\large {\bf About This Book}}

\vspace*{.1in}

In this book we interpret additive combinatorics, somewhat narrowly, as the study of sumsets in finite abelian groups.  (A second volume, focusing on the set of integers and other infinite groups, is in the planning phase.)   This book offers an extensive menu of research projects to any student interested in pursuing investigations in this area.  It contains five parts that we briefly describe, as follows.

{\bf Ingredients.}  We make no assumptions on the backgrounds of students wishing to engage in research projects based on this collection.  In order to equip beginners with the necessary background, we provide brief introductions to the relevant branches of number theory, combinatorics, and group theory.  The sections in each chapter contain exercises aimed to solidify the understanding of the material discussed.

{\bf Appetizers.}
The short articles in this part are meant to invite everyone to the main entrees of the menu; the appetizers are carefully chosen so that they provide bite-size representative samples of the research projects, as well as make connections to other parts of mathematics that students might have encountered.  (One article describes how this author's research in spherical geometry lead him to additive combinatorics.)

{\bf Sides.}  Here we present and study some auxiliary functions that appear in several different chapters of the text.

{\bf Entrees.}  Each chapter in this main part of the book discusses a particular family of open questions in additive combinatorics.  At the present time, these questions are about 
  \\ \hspace*{.5in} -- maximum sumset size, 
  \\ \hspace*{.5in} -- spanning sets, 
  \\ \hspace*{.5in} -- Sidon sets, 
  \\ \hspace*{.5in} -- minimum sumset size,  
  \\ \hspace*{.5in} -- critical numbers, 
  \\ \hspace*{.5in} -- zero-sum-free sets, and 
  \\ \hspace*{.5in} -- sum-free sets.   \\  (We plan to include additional chapters in the near future.)  Each chapter is divided into the same four sections, depending on whether the sums may contain repeated terms or not and whether terms can be both added and subtracted or can only be added; each section is then further divided into the same three subsections depending on the number of terms in the sums: whether it is fixed, limited, or is arbitrary.  Within most subsections, we investigate both {\em direct questions} (e.g. what is the maximum size of a set with the required sumset property) and {\em inverse questions} (how can we characterize all subsets that achieve this maximum size). 

Many of the questions we discuss have been investigated extensively already, in which case we include all relevant results that are available to the best of our knowledge.  We find it essential that our expositions are complete, and the author is committed to continually survey the literature to assure that the material remains current.  (For a variety of reasons, this is not always an easy task: as is often the case with relatively young fields, results appear in a diverse set of outlets, are often presented using different notations and terminology, and, occasionally, are even incorrect.)   

While we aim to present each topic as thoroughly as possible, this book is not a historical survey.  In particular, in most cases, we only include the best results known on any given question; previously achieved special cases and weaker or partial results can then usually be found in the citations we provide.  Furthermore, any reader with an interest in a more contextual treatment of our topics or a desire to see them from different perspectives is encouraged to turn to one of the available books on the subject, such as \cite{Fre:1973a} by Freiman,\index{Freiman, G. A.}  \cite{GerRuz:2009a} by Geroldinger\index{Geroldinger, A.} and Ruzsa,\index{Ruzsa, I.} \cite{Gry:2013a} by Grynkiewicz,\index{Grynkiewicz, D. J.} \cite{Guy:2004a} by Guy,\index{Guy, R.} \cite{HalRot:1983a} by Halberstam\index{Halberstam, H.} and Roth,\index{Roth, K. F.} \cite{Nat:1996a} by Nathanson,\index{Nathanson, M. B.} or \cite{TaoVu:2006a} by Tao\index{Tao, T.} and Vu.\index{Vu, V. H.}

{\bf Pudding.}  The proofs are in the pudding!---well, not quite.  The book contains several hundred stated results; most of these results have been published, in which case we give citations where proofs (as well as further material) can be read.  In other cases, results have not appeared elsewhere; proofs of these statements---unless rather brief and particularly instructive, in which case they are presented where the results are stated---are separated into this part of the book.

 \vspace*{.25in}

\noindent {\large {\bf Acknowledgments}}

\vspace*{.1in}

The author is most grateful for the immense support he has received from Gettysburg College.  The frequent full-year sabbaticals and the many very generous travel grants are much appreciated.  Special thanks go to my colleagues in the Mathematics Department who allowed me to teach a research course based on this book in virtually every semester during the past fifteen years.  

I am grateful to the many friends and colleagues who have provided useful comments, feedback, and encouragement, including Noga Alon, L\'aszl\'o Babai, Mikl\'os B\'ona, Jill Dietz, R\'obert Freud, Ronald Graham, Ben Green, David Grynkiewicz, Mark Kayll, Vsevolod Lev, Ryan Matzke, Steven J. Miller, Wolfgang Schmid, Neil Sloane, Zhi--Wei Sun, James Tanton, Terence Tao, and Paul Zeitz.  I am particularly thankful to Samuel Edwards for his help with the bibliography and to Ivaylo Ilinkin for the computer program \cite{Ili:2017a}  that he created for this book.

But my biggest appreciation goes to the two hundred or so students at Gettysburg College who have carried out research projects based on this text thus far.  Over thirty of these students attained results that are mentioned and cited in this book, at least twenty have given conference presentations on their work, and several have even published their results in refereed journals---I congratulate them full-heartedly!  I am convinced, however, that every student benefited from the opportunities for perfecting a variety of skills by engaging in a research experience based on this book.

\chapter{Notations}

Following is the list of the most commonly used notations throughout this book.

\vspace{.4in}

\noindent {\Large {\bf Chapter 1: Number theory}}

\vspace{.2in}

{\bf Number sets:}

Integers: $\mathbb{Z}=\{0, \pm 1, \pm 2, \pm 3, \dots\}$

Nonnegative integers: $\mathbb{N}_0=\{0, 1, 2, 3, \dots\}$

Positive integers: $\mathbb{N}=\{1, 2, 3, \dots\}$

Interval: $[a,b]=\{ c \in \mathbb{Z} \mid a \leq c \leq b\}$  for given $a, b \in \mathbb{Z}$ with $a \leq b$ 

\vspace{.1in}

{\bf Divisors:}

$D(n)=\{ d \in \mathbb{N} \mid d|n \}$ for a given $n \in \mathbb{N}$

$d(n)=|D(n)|$ for a given $n \in \mathbb{N}$

$\gcd (a,b) = \max (D(a) \cap D(b))$ for given $a, b \in \mathbb{N}$

$\gcd (A) = \max \{ d \in \mathbb{N} \mid \forall a \in A, d|n \}$ for a given finite $A \subset \mathbb{N}$

\vspace{.2in}

\noindent {\Large{\bf Chapter 2: Combinatorics}}

\vspace{.2in}

{\bf Counting functions:}

$n^{\overline{m}}=n(n+1)\cdots (n+m-1)$ and $n^{\overline{0}}=1$ for $n \in \mathbb{N}_0$ and $m \in \mathbb{N}$ 

$n^{\underline{m}}=n(n-1)\cdots (n-m+1)$ and $n^{\underline{0}}=1$ for $n \in \mathbb{N}_0$ and $m \in \mathbb{N}$

${n \choose m}= \frac{n^{\underline{m}}}{m!} = \frac{n(n-1)\cdots(n-m+1)}{m!}$ for $n, m \in \mathbb{N}_0$

$\left[ {n \atop m} \right] = \frac{n^{\overline{m}}}{m!} = \frac{n(n+1)\cdots(n+m-1)}{m!}$ for $n, m \in \mathbb{N}_0$

\vspace{.1in}

{\bf Some other useful counting functions:}

$a(j,k)=\sum_{i \geq 0} {j \choose i} {k \choose i} 2^i=\left\{
\begin{array}{cl}
1 & \mbox{if $j=0$ or $k=0$} \\
a(j-1,k)+a(j-1,k-1)+a(j,k-1) & \mbox{if $j \geq 1$ and $k \geq 1$}
\end{array}
\right.$

$c(j,k)=\sum_{i \geq 0} {j-1 \choose i-1} {k \choose i} 2^i=\left\{
\begin{array}{cl}
1 & \mbox{if $j=0$ and $k \geq 0$} \\
0 & \mbox{if $j \geq 1$ and $k =0$} \\
c(j-1,k)+c(j-1,k-1)+c(j,k-1) & \mbox{if $j \geq 1$ and $k \geq 1$}
\end{array}
\right.$

\vspace{.1in}

{\bf Layers of the integer lattice:}

$\mathbb{N}_0^m (h) = \{ (\lambda_1, \dots, \lambda_m) \in \mathbb{N}_0^m \mid \Sigma_{i=1}^m \lambda_i = h\}$

$\mathbb{Z}^m (h) = \{ (\lambda_1, \dots, \lambda_m) \in \mathbb{Z}^m \mid \Sigma_{i=1}^m |\lambda_i| = h\}$

$\hat{\mathbb{N}}_0^m (h) = \{ (\lambda_1, \dots, \lambda_m) \in \{0,1\}^m \mid \Sigma_{i=1}^m \lambda_i = h\}$

$\hat{\mathbb{Z}}^m (h) = \{ (\lambda_1, \dots, \lambda_m) \in \{-1,0,1\}^m \mid \Sigma_{i=1}^m |\lambda_i |= h\}$

\vspace{.2in}

\noindent {\Large{\bf Chapter 3: Group theory}}

\vspace{.2in}

{\bf Finite abelian groups:}

$\mathbb{Z}_n$: cyclic group of order $n$ using additive notation

$G$: (arbitrary) abelian group of order $n$ using additive notation

$r$: rank of $G$ (number of factors in the invariant decomposition of $G$)

$\kappa$: exponent of $G$ (order of largest factor in the invariant decomposition of $G$)

\vspace{.1in}

{\bf Elements and subsets:}

$\mathrm{ord}_G(a)=\min \{d \in \mathbb{N} \mid da=0\}$

$\mathrm{Ord}(G,d)=\{ a \in G \mid \mathrm{ord}_G(a)=d\}$

$\langle a \rangle =\{\lambda a \mid \lambda \in [1,\mathrm{ord}_G(a)]\}$ 

$A=\{a_1,\dots,a_m\}$: (arbitrary) $m$-subset of $G$

$\langle A \rangle =\{\Sigma_{i=1}^m \lambda_i a_i \mid \lambda_i \in [1,\mathrm{ord}_G(a_i)]\}$ 

\vspace{.1in}

{\bf Unrestricted sumsets:}

$hA=\{\Sigma_{i=1}^m \lambda_i a_i \mid \lambda_i \in \mathbb{N}_0, \Sigma \lambda_i=h\}$ 

$[0,s]A=\cup_{h=0}^s hA=\{\Sigma_{i=1}^m \lambda_i a_i \mid \lambda_i \in \mathbb{N}_0, \Sigma \lambda_i \leq s\}$ 

$\langle A \rangle =\cup_{h=0}^{\infty} hA=\{\Sigma_{i=1}^m \lambda_i a_i \mid \lambda_i \in \mathbb{N}_0\}$

\vspace{.1in}

{\bf Unrestricted signed sumsets:}

$h_{\pm}A=\{\Sigma_{i=1}^m \lambda_i a_i \mid \lambda_i \in \mathbb{Z}, \Sigma |\lambda_i|=h\}$

$[0,s]_{\pm}A=\cup_{h=0}^s h_{\pm}A =\{\Sigma_{i=1}^m \lambda_i a_i \mid \lambda_i \in \mathbb{Z}, \Sigma |\lambda_i| \leq s\}$

$\langle A \rangle =\cup_{h=0}^{\infty} h_{\pm}A=\{\Sigma_{i=1}^m \lambda_i a_i \mid \lambda_i \in \mathbb{Z}\}$

\vspace{.1in}

{\bf Restricted sumsets:}

$h \hat{\;} A=\{\Sigma_{i=1}^m \lambda_i a_i \mid \lambda_i \in \{0,1\}, \Sigma \lambda_i=h\}$

$[0,s] \hat{\;} A =\cup_{h=0}^s h \hat{\;} A=\{\Sigma_{i=1}^m \lambda_i a_i \mid \lambda_i \in \{0,1\}, \Sigma \lambda_i \leq s\}$

$\Sigma A =\cup_{h=0}^{\infty} h \hat{\;} A=\{\Sigma_{i=1}^m \lambda_i a_i \mid \lambda_i \in \{0,1\}\}$

$\Sigma^* A =\cup_{h=1}^{\infty} h \hat{\;} A=\{\Sigma_{i=1}^m \lambda_i a_i \mid \lambda_i \in \{0,1\}, \Sigma \lambda_i \geq 1\}$

\vspace{.1in}

{\bf Restricted signed sumsets:}

$h \hat{_{\pm}} A=\{\Sigma_{i=1}^m \lambda_i a_i \mid \lambda_i \in \{-1,0,1\}, \Sigma |\lambda_i|=h\}$

$[0,s] \hat{_{\pm}} A =\cup_{h=0}^s h \hat{_{\pm}} A =\{\Sigma_{i=1}^m \lambda_i a_i \mid \lambda_i \in \{-1,0,1\}, \Sigma |\lambda_i| \leq s\}$

$\Sigma_{\pm} A =\cup_{h=0}^{\infty} h \hat{_{\pm}} A=\{\Sigma_{i=1}^m \lambda_i a_i \mid \lambda_i \in \{-1,0,1\}\}$

$\Sigma_{\pm}^* A =\cup_{h=1}^{\infty} h \hat{_{\pm}} A=\{\Sigma_{i=1}^m \lambda_i a_i \mid \lambda_i \in \{-1,0,1\}, \Sigma |\lambda_i| \geq 1\}$

\vspace{.2in}

\noindent {\Large{\bf Sides}}

\vspace{.2in}

$v_g(n,h)=\max \left\{ \left( \left \lfloor \frac{d-1-\mathrm{gcd} (d, g)}{h} \right \rfloor +1  \right) \cdot \frac{n}{d}  \mid d \in D(n) \right\}$  for given $n, h, g \in \mathbb{N}$

$v_{\pm}(n,h)=\max \left\{ \left( 2 \cdot \left \lfloor \frac{d-2}{2h} \right \rfloor +1  \right) \cdot \frac{n}{d}  \mid d \in D(n) \right\}$  for given $n, h \in \mathbb{N}$

$u(n,m,h)=\min \left\{ \left( h \cdot \left \lceil \frac{m}{d} \right \rceil -h +1 \right) \cdot d  \mid d \in D(n) \right\}$ for given $n, m, h \in \mathbb{N}$

\vspace{.2in}

\noindent {\Large{\bf Chapter A: Maximum sumset size}}

\vspace{.2in}

{\bf A.1:  Unrestricted sumsets:}

$\nu(G,m,h)= \max \{ |hA| \mid A \subseteq G, |A|=m\}$

$\nu(G,m,[0,s])= \max \{ |[0,s]A| \mid A \subseteq G, |A|=m\}$

$\nu(G,m,\mathbb{N}_0)= \max \{ |\langle A \rangle | \mid A \subseteq G, |A|=m\}$

\vspace{.1in}

{\bf A.2: Unrestricted signed sumsets:}

$\nu_{\pm} (G,m,h)= \max \{ |h_{\pm} A| \mid A \subseteq G, |A|=m\}$

$\nu_{\pm} (G,m,[0,s])= \max \{ |[0,s]_{\pm} A| \mid A \subseteq G, |A|=m\}$

$\nu_{\pm} (G,m,\mathbb{N}_0)= \max \{ |\langle A \rangle | \mid A \subseteq G, |A|=m\}$

\vspace{.1in}

{\bf A.3: Restricted sumsets:}

$\nu \hat{\;} (G,m,h)= \max \{ |h\hat{\;} A| \mid A \subseteq G, |A|=m\}$

$\nu\hat{\;} (G,m,[0,s])= \max \{ |[0,s]\hat{\;} A| \mid A \subseteq G, |A|=m\}$

$\nu \hat{\;}  (G,m,\mathbb{N}_0)= \max \{ |\Sigma A | \mid A \subseteq G, |A|=m\}$

\vspace{.1in}

{\bf A.4: Restricted signed sumsets:}

$\nu \hat{_{\pm}} (G,m,h)= \max \{ |h \hat{_{\pm}} A| \mid A \subseteq G, |A|=m\}$

$\nu \hat{_{\pm}} (G,m,[0,s])= \max \{ |[0,s] \hat{_{\pm}} A| \mid A \subseteq G, |A|=m\}$

$\nu \hat{_{\pm}} (G,m,\mathbb{N}_0)= \max \{ |\Sigma _{\pm} A | \mid A \subseteq G, |A|=m\}$

\vspace{.2in}

\noindent {\Large{\bf Chapter B: Spanning sets}}

\vspace{.2in}

{\bf B.1:  Unrestricted sumsets:}

$\phi(G,h)= \min \{ |A| \mid A \subseteq G, hA=G\}$

$\phi(G,[0,s])= \min \{ |A| \mid A \subseteq G, [0,s]A=G\}$

$\phi(G,\mathbb{N}_0)= \min \{ |A | \mid A \subseteq G, \langle A \rangle=G\}$

\vspace{.1in}

{\bf B.2: Unrestricted signed sumsets:}

$\phi_{\pm} (G,h)= \min \{ | A| \mid A \subseteq G, h_{\pm} A=G\}$

$\phi_{\pm} (G,[0,s])= \min \{ | A| \mid A \subseteq G, [0,s]_{\pm} A=G\}$

$\phi_{\pm} (G,\mathbb{N}_0)= \min \{ |A | \mid A \subseteq G, \langle A \rangle=G\}$

\vspace{.1in}

{\bf B.3: Restricted sumsets:}

$\phi \hat{\;} (G,h)= \min \{ | A| \mid A \subseteq G, h\hat{\;} A=G\}$

$\phi\hat{\;} (G,[0,s])= \min \{ | A| \mid A \subseteq G, [0,s]\hat{\;} A = G\}$

$\phi \hat{\;}  (G,\mathbb{N}_0)= \min \{ | A | \mid A \subseteq G, \Sigma A = G\}$

\vspace{.1in}

{\bf B.4: Restricted signed sumsets:}

$\phi \hat{_{\pm}} (G,h)= \min \{ | A| \mid A \subseteq G, h \hat{_{\pm}} A = G\}$

$\phi \hat{_{\pm}} (G,[0,s])= \min \{ | A| \mid A \subseteq G, [0,s] \hat{_{\pm}} A = G\}$

$\phi \hat{_{\pm}} (G,\mathbb{N}_0)= \min \{ | A | \mid A \subseteq G, \Sigma _{\pm} A = G\}$

\vspace{.2in}

\noindent {\Large{\bf Chapter C: Sidon sets}}

\vspace{.2in}

{\bf C.1:  Unrestricted sumsets:}

$\sigma(G,h)= \max \left \{ |A| \mid A \subseteq G, |hA|=|\mathbb{N}_0^m (h)|={m+h-1 \choose h}\right \}$

$\sigma(G,[0,s])= \max \left \{ |A| \mid A \subseteq G, |[0,s]A|=\Sigma_{h=0}^s |\mathbb{N}_0^m (h)|={m+s \choose s} \right \}$

\vspace{.1in}

{\bf C.2: Unrestricted signed sumsets:}

$\sigma_{\pm} (G,h)= \max \left \{ | A| \mid A \subseteq G, |h_{\pm} A|=|\mathbb{Z}^m (h)|=c(h,m)\right \}$

$\sigma_{\pm} (G,[0,s])= \max \left \{ | A| \mid A \subseteq G, |[0,s]_{\pm} A|=\Sigma_{h=0}^s |\mathbb{Z}^m (h)|=a(m,s) \right \}$

\vspace{.1in}

{\bf C.3: Restricted sumsets:}

$\sigma \hat{\;} (G,h)= \max \left \{ | A| \mid A \subseteq G, |h\hat{\;} A|=|\hat{\mathbb{N}}_0^m (h)|={m \choose h}\right \}$

$\sigma\hat{\;} (G,[0,s])= \max \left \{ | A| \mid A \subseteq G, |[0,s]\hat{\;} A| = \Sigma_{h=0}^s |\hat{\mathbb{N}}_0^m (h)|= \Sigma_{h=0}^s  {m \choose h} \right \}$

$\sigma \hat{\;}  (G,\mathbb{N}_0)= \max \left \{ | A | \mid A \subseteq G, |\Sigma A |= \Sigma_{h=0}^m |\hat{\mathbb{N}}_0^m (h)|= 2^m \right \}$

\vspace{.1in}

{\bf C.4: Restricted signed sumsets:}

$\sigma \hat{_{\pm}} (G,h)= \max \left \{ | A| \mid A \subseteq G, |h \hat{_{\pm}} A| = |\hat{\mathbb{Z}}^m (h)|={m \choose h} 2^h\right \}$

$\sigma \hat{_{\pm}} (G,[0,s])= \max \left \{ | A| \mid A \subseteq G, |[0,s] \hat{_{\pm}} A| = \Sigma_{h=0}^s |\hat{\mathbb{Z}}^m (h)|=  \Sigma_{h=0}^s {m \choose h} 2^h  \right \}$

$\sigma \hat{_{\pm}} (G,\mathbb{N}_0)= \max \left \{ | A | \mid A \subseteq G, |\Sigma _{\pm} A |=\Sigma_{h=0}^m |\hat{\mathbb{Z}}^m (h)|= 3^m \right \}$

\vspace{.2in}

\noindent {\Large{\bf Chapter D: Minimum sumset size}}

\vspace{.2in}

{\bf D.1:  Unrestricted sumsets:}

$\rho(G,m,h)= \min \{ |hA| \mid A \subseteq G, |A|=m\}$

$\rho(G,m,[0,s])= \min \{ |[0,s]A| \mid A \subseteq G, |A|=m\}$

$\rho(G,m,\mathbb{N}_0)= \min \{ |\langle A \rangle | \mid A \subseteq G, |A|=m\}$

\vspace{.1in}

{\bf D.2: Unrestricted signed sumsets:}

$\rho_{\pm} (G,m,h)= \min \{ |h_{\pm} A| \mid A \subseteq G, |A|=m\}$

$\rho_{\pm} (G,m,[0,s])= \min \{ |[0,s]_{\pm} A| \mid A \subseteq G, |A|=m\}$

$\rho_{\pm} (G,m,\mathbb{N}_0)= \min \{ |\langle A \rangle | \mid A \subseteq G, |A|=m\}$

\vspace{.1in}

{\bf D.3: Restricted sumsets:}

$\rho \hat{\;} (G,m,h)= \min \{ |h\hat{\;} A| \mid A \subseteq G, |A|=m\}$

$\rho\hat{\;} (G,m,[0,s])= \min \{ |[0,s]\hat{\;} A| \mid A \subseteq G, |A|=m\}$

$\rho \hat{\;}  (G,m,\mathbb{N}_0)= \min \{ |\Sigma A | \mid A \subseteq G, |A|=m\}$

\vspace{.1in}

{\bf D.4: Restricted signed sumsets:}

$\rho \hat{_{\pm}} (G,m,h)= \min \{ |h \hat{_{\pm}} A| \mid A \subseteq G, |A|=m\}$

$\rho \hat{_{\pm}} (G,m,[0,s])= \min \{ |[0,s] \hat{_{\pm}} A| \mid A \subseteq G, |A|=m\}$

$\rho \hat{_{\pm}} (G,m,\mathbb{N}_0)= \min \{ |\Sigma _{\pm} A | \mid A \subseteq G, |A|=m\}$

\vspace{.2in}

\noindent {\Large{\bf Chapter E: The critical number}}

\vspace{.2in}

{\bf E.1:  Unrestricted sumsets:}

$\chi(G,h)= \min \{ m \mid A \subseteq G, |A|=m \Rightarrow hA=G  \}$

$\chi(G,[0,s])= \min \{ m\mid A \subseteq G, |A|=m  \Rightarrow [0,s]A=G \}$

$\chi(G,\mathbb{N}_0)= \min \{ m \mid A \subseteq G, |A|=m  \Rightarrow \langle A \rangle =G \}$

\vspace{.1in}

{\bf E.2: Unrestricted signed sumsets:}

$\chi_{\pm} (G,h)= \min \{ m \mid A \subseteq G, |A|=m  \Rightarrow h_{\pm}A=G \}$

$\chi_{\pm} (G,[0,s])= \min \{ m \mid A \subseteq G, |A|=m  \Rightarrow [0,s]_{\pm}A=G \}$

$\chi_{\pm} (G,\mathbb{N}_0)= \min \{ m \mid A \subseteq G, |A|=m  \Rightarrow \langle A \rangle=G \}$

\vspace{.1in}

{\bf E.3: Restricted sumsets:}

$\chi \hat{\;} (G,h)= \min \{ m \mid A \subseteq G, |A|=m  \Rightarrow h\hat{\;}A=G \}$

$\chi\hat{\;} (G,[0,s])= \min \{ m \mid A \subseteq G, |A|=m  \Rightarrow [0,s] \hat{\;}A=G \}$

$\chi \hat{\;}  (G,\mathbb{N}_0)= \min \{ m \mid A \subseteq G, |A|=m  \Rightarrow \Sigma A=G \}$

\vspace{.1in}

{\bf E.4: Restricted signed sumsets:}

$\chi \hat{_{\pm}} (G,h)= \min \{ m \mid A \subseteq G, |A|=m  \Rightarrow h\hat{_{\pm}} A=G \}$

$\chi \hat{_{\pm}} (G,[0,s])= \min \{ m \mid A \subseteq G, |A|=m  \Rightarrow [0,s] \hat{_{\pm}}  A=G \}$

$\chi \hat{_{\pm}} (G,\mathbb{N}_0)= \min \{ m \mid A \subseteq G, |A|=m  \Rightarrow \Sigma \hat{_{\pm}} A=G \}$

\vspace{.2in}

\noindent {\Large{\bf Chapter F: Zero-sum-free sets}}

\vspace{.2in}

{\bf F.1:  Unrestricted sumsets:}

$\tau(G,h)= \max \{ |A| \mid A \subseteq G, 0 \not \in hA\}$

$\tau(G,[1,t])= \max \{ |A| \mid A \subseteq G, 0 \not \in  [1,t]A\}$

\vspace{.1in}

{\bf F.2: Unrestricted signed sumsets:}

$\tau_{\pm} (G,h)= \max \{ | A| \mid A \subseteq G, 0 \not \in  h_{\pm} A\}$

$\tau_{\pm} (G,[1,t])= \max \{ | A| \mid A \subseteq G, 0 \not \in  [1,t]_{\pm} A\}$

\vspace{.1in}

{\bf F.3: Restricted sumsets:}

$\tau \hat{\;} (G,h)= \max \{ | A| \mid A \subseteq G, 0 \not \in  h\hat{\;} A\}$

$\tau\hat{\;} (G,[1,t])= \max \{ | A| \mid A \subseteq G, 0 \not \in [1,t]\hat{\;} A \}$

$\tau \hat{\;}  (G,\mathbb{N})= \max \{ | A | \mid A \subseteq G, 0 \not \in [1,m]\hat{\;} A \}$

\vspace{.1in}

{\bf F.4: Restricted signed sumsets:}

$\tau \hat{_{\pm}} (G,h)= \max \{ | A| \mid A \subseteq G, 0 \not \in h \hat{_{\pm}} A \}$

$\tau \hat{_{\pm}} (G,[1,t])= \max \{ | A| \mid A \subseteq G, 0 \not \in  [1,t] \hat{_{\pm}} A \}$

$\tau \hat{_{\pm}} (G,\mathbb{N})= \max \{ | A | \mid A \subseteq G, 0 \not \in  [1,m] \hat{_{\pm}} A \}$

\vspace{.2in}

\noindent {\Large{\bf Chapter G: Sum-free sets}}

\vspace{.2in}

{\bf G.1:  Unrestricted sumsets:}

$\mu(G,\{k,l\})= \max \{ |A| \mid A \subseteq G, kA \cap lA = \emptyset\}$

$\mu(G,[0,s])= \max \{ |A| \mid A \subseteq G, 0 \leq l < k \leq s \Rightarrow kA \cap lA = \emptyset\}$

\vspace{.1in}

{\bf G.2: Unrestricted signed sumsets:}

$\mu_{\pm} (G,\{k,l\})= \max \{ | A| \mid A \subseteq G, k_{\pm}A \cap l_{\pm}A = \emptyset\}$

$\mu_{\pm} (G,[0,s])= \max \{ | A| \mid A \subseteq G, 0 \leq l < k \leq s \Rightarrow k_{\pm}A \cap l_{\pm}A = \emptyset\}$

\vspace{.1in}

{\bf G.3: Restricted sumsets:}

$\mu \hat{\;} (G,\{k,l\})= \max \{ | A| \mid A \subseteq G, k\hat{\;}A \cap l\hat{\;}A = \emptyset\}$

$\mu\hat{\;} (G,[0,s])= \max \{ | A| \mid A \subseteq G, 0 \leq l < k \leq s \Rightarrow k\hat{\;}A \cap l\hat{\;}A = \emptyset \}$

$\mu \hat{\;}  (G,\mathbb{N}_0)= \max \{ | A | \mid A \subseteq G, 0 \leq l < k  \Rightarrow k \hat{\;} A \cap l \hat{\;} A = \emptyset \}$

\vspace{.1in}

{\bf G.4: Restricted signed sumsets:}

$\mu \hat{_{\pm}} (G,\{k,l\})= \max \{ | A| \mid A \subseteq G, kA \cap lA = \emptyset \}$

$\mu \hat{_{\pm}} (G,[0,s])= \max \{ | A| \mid A \subseteq G, 0 \leq l < k \leq s \Rightarrow k\hat{_{\pm}}A \cap l\hat{_{\pm}}A = \emptyset \}$

$\mu \hat{_{\pm}} (G,\mathbb{N}_0)= \max \{ | A | \mid A \subseteq G, 0 \leq l < k  \Rightarrow k\hat{_{\pm}}A \cap l\hat{_{\pm}}A = \emptyset \}$

\vspace{.2in}

\mainmatter

\part{Ingredients} \label{Chapter0}

This book deals with {\em additive combinatorics}, a vibrant area of current mathematical research.  Additive combinatorics---which grew out of combinatorial number theory and additive number theory---here is interpreted, somewhat narrowly, as the study of combinatorial properties of {\em sumsets} in abelian groups.  In Chapters \ref{IntroNTh}, \ref{IntroComb}, and \ref{IntroGrTh} we provide brief introductions to the relevant branches of {\em number theory}, {\em combinatorics}, and {\em group theory}, respectively.  The sections in each chapter contain exercises aimed to solidify the understanding of the material discussed.

\chapter{Number theory} \label{IntroNTh}

Number theory, at least in its most traditional form, is the branch of mathematics that studies the set of integers:
$$\mathbb{Z}=\{\dots,-2,-1,0,1,2,3,\dots\},$$ or one of its subsets, such as the set of nonnegative integers:
$$\mathbb{N}_0=\{0,1,2,3,4,5,\dots\}$$ or the set of positive integers:
$$\mathbb{N}=\{1,2,3,4,5,6,\dots\}.$$

The field of number theory occupies a distinguished spot within mathematics: the German mathematician Carl Friedrich Gauss (1777--1855) even dubbed it the ``Queen of Mathematics.''  Some of the reasons for this distinction is that it (or she?) manages to simultaneously possess the following seemingly contradictory attributes:

\begin{itemize}

\item although a substantial amount of the terminology and the methods in number theory, at least at the introductory level, are quite simple and need very little background, number theory reaches into some of the deepest and most complex areas of mathematics;

\item even though many of its questions and problems are easy to present, number theory has a cornucopia of impossibly difficult unsolved questions, some hundreds of years old; 

\item though number theory may be the oldest branch of mathematics and has always had a large number of devotees, it remains one of the most active and perplexing fields within mathematics with new developments and results being published every day.  

\end{itemize}

Below we review some of the foundations of number theory that we will rely on later.

\section{Divisibility of integers} \label{divis}

The most fundamental concept in number theory is probably {\em divisibility}: Given two integers $a$ and $b$, we say that $a$ is a \emph{divisor} of $b$ (or $b$ is \emph{divisible} by $a$) whenever there is an integer $c$ for which $a \cdot c=b$.  If $a$ is a divisor of $b$, we write $a | b$.

For example, $3 | 6$ and $6 | 6$, but $6\not | 3$. (We must be clear with our terminology: 3 is not divisible by 6, but of course 3 can be divided by 6; in the set $\mathbb{Q}$ of rational numbers---fractions of integers---the concept of divisibility is trivial in that every rational number is divisible by every nonzero rational number.) Also, $5 | 0$, since $5 \cdot 0 =0$; in fact, $0 | 0$, since (for example) $0 \cdot 7 =0$.  But $0\not | 5$, since there is no integer $c$ for which $0 \cdot c=5$ because for every real number $c$, $0 \cdot c=0$.
As it is often the case with mathematical definitions, one needs to be careful: saying that $a$ is a divisor of $b$ is not quite equivalent to saying that the fraction $b/a$ is an integer: $0$ is a divisor of $0$, but $0/0$ is not an integer (it's not even a number)!

In theory, for each integer $n$, we can easily find the set of its positive divisors, denoted by $D(n)$.  For example, we have
$$D(18)=D(-18)=\{1,2,3,6,9,18\}$$ and
$$D(19)=D(-19)=\{1,19\}.$$
Finding the divisors of large positive integers can be very difficult.  Cryptography, the study of encoding and decoding information, takes advantage of the dichotomy that multiplying two large integers (as sometimes used in encoding) is easy, but factoring the product without knowing any of the factors (i.e. decoding) could be, if the factors are chosen carefully, very hard!

Given an integer $n$, we denote the number of its positive divisors by $d(n)$; that is, we set $$d(n)=|D(n)|.$$  For instance, as the examples above show, we have $d(18)=d(-18)=6$ and $d(19)=d(-19)=2$.  Clearly, $d(n)$ is always positive as $1|n$ for every integer $n$.  The function $d(n)$ allows us to separate the set of integers into four classes:
\begin{itemize}
\item {\em units} have exactly one positive divisor;
\item {\em prime numbers} have exactly two positive divisors;
\item {\em composite numbers} have three or more positive divisors;
 and
\item {\em zero} has infinitely many positive divisors.

\end{itemize}
 
There are two units among the integers: 1 and $-1$ divide all integers, but have only one positive divisor.  The set of prime numbers $$P=\{\pm 2, \pm 3, \pm 5, \pm 7, \pm 11, \pm 13, \pm 17, \pm 19, \dots\},$$ as we explain shortly, forms the basic building block of the set of integers.  Primes have been studied for thousands of years.  They have many intriguing attributes, some of which are not fully understood to this day.  For those interested in more on primes may start their investigations with \cite{OEIS}, the On-Line Encyclopedia of Integer Sequences; the positive primes appear as the sequence A000040.

Given two positive integers $a$ and $b$, we define their {\em greatest common divisor}, denoted by $\mathrm{gcd}(a,b)$, to be the greatest integer that is a divisor of both $a$ and $b$; in other words, $$\mathrm{gcd}(a,b) = \max (D(a) \cap D(b)).$$  A pair of integers is said to be relatively prime if their greatest common divisor is 1.  The dual concept is the {\em least common multiple} of two positive integers: it is the smallest positive integer that is a multiple of both of them.  The least common multiple of integers $a$ and $b$ is denoted by $\mathrm{lcm}(a,b)$.  It is not hard to see that all pairs of positive integers have a unique greatest common divisor as well as a least common multiple.

\vspace{.17in}

{\bf Exercises}

\begin{enumerate}

\item 

Characterize all integers $n$ for which 

\begin{enumerate}

\item
$d(n)=3$;

\item
$d(n)=4$;

\item
$d(n)=5$.

\end{enumerate}

\item 

\begin{enumerate}

\item

The concepts of greatest common divisor and least common multiple of two positive integers can be extended to three or more positive integers.  Find 
$\mathrm{gcd}(24,32,60)$ and $\mathrm{lcm}(24,32,60)$.

\item Find positive integers $a$, $b$, $c$, and $d$ that are relatively prime, but no two of them are relatively prime (that is, $\mathrm{gcd}(a,b,c,d)=1$, but $\mathrm{gcd}(a,b)>1$, $\mathrm{gcd}(a,c)>1$, etc.).

\end{enumerate}

\end{enumerate}

\section{Congruences}

As a generalization of divisibility, we may consider situations where an integer leaves a remainder when divided by another integer.  More precisely, given a positive integer $m$, we say that an integer $a$ leaves a {\em remainder} $r$ when divided by $m$, if there is an integer $k$ for which $a=m \cdot k + r$ and $0 \leq r \leq m-1$; in this case we sometimes also say that $a$ is {\em congruent} to $r$ {\em mod} $m$ and write $a \equiv r$ mod $m$.  For example, both 13 and 1863 are congruent to 3 mod 10, and so is $-57$ as it can be written as $10 \cdot (-6)+3$.   (Recall that any remainder mod $10$ must be between $0$ and $9$, inclusive.)  It is not hard to see that for any positive integer $m$ and for any integer $a$, one can determine a unique remainder $r$ that $a$ is congruent to mod $m$.  

A bit more generally, given a positive integer $m$, we say that two integers  are {\em congruent mod} $m$ if they leave the same remainder when divided by $m$.  For example, 13 and 1863 are congruent mod $10$, since they both leave a remainder of 3 when divided by 10; we denote this by writing $13 \equiv 1863$ mod $10$.  In fact, all positive integers with a last decimal digit of 3 and all negative integers with a last digit of 7 are congruent to each other, and they form the {\em congruence class} of $$[3]_{10}=\{10k+3  \mid k \in \mathbb{Z} \}=\{\dots, -27, -17, -7, 3, 13, 23, 33, \dots\}.$$

Congruence classes allow us to partition the set of integers in a natural way; for example, we have
$$\mathbb{Z}=[0]_{10} \cup [1]_{10} \cup [2]_{10} \cup [3]_{10} \cup [4]_{10} \cup [5]_{10} \cup [6]_{10} \cup [7]_{10} \cup [8]_{10} \cup [9]_{10}.$$
Since, for a given $m \in \mathbb{N}$, the remainder $r$ may be any integer value between 0 and $m-1$, we have exactly $m$ congruence classes mod $m$; these classes are disjoint and their union contains every integer.  For example, the congruence classes $[0]_2$ and $[1]_2$ contain the even and the odd integers, respectively; the fact that any integer must be either even or odd but not both can be expressed by saying that $[0]_2$ and $[1]_2$ form a {\em partition} of $\mathbb{Z}$.

Congruence classes play a prominent role in additive combinatorics, and we will return to their study shortly.

\vspace{.17in}

{\bf Exercises}

\begin{enumerate}

\item Find all integers between $-100$ and 100 that are in the congruence class $[7]_{33}$.

\item \begin{enumerate}

\item Find integers $a$, $b$, $c$, and $d$ for which the congruence classes $[a]_{2}$, $[b]_{4}$, $[c]_{8}$, and $[d]_{8}$ partition $\mathbb{Z}$.  (To be a partition, every integer must belong to exactly one congruence class.)

\item Partition $\mathbb{Z}$ into exactly five congruence classes where four of the five moduli are distinct.  (It is a well-known result that the moduli cannot be all distinct.)

\end{enumerate}

\end{enumerate}

\section{The Fundamental Theorem of Number Theory} \label{FTNT}

In general, there may be many ways to factor an integer into a product of other integers.  In fact, if we allow 1 or $-1$ to appear as factors, then every integer has infinitely many different factorizations.  For example, factorizations of $18$ include $1 \cdot 18$, $1 \cdot (-1) \cdot (-18)$, $2 \cdot 3 \cdot 3$, and $3 \cdot 6$.  Here and in the next section we briefly discuss two of these factorizations (namely, generalizations of the last two factorizations of $18$ that we listed). 

First, there is what is referred to as the {\em prime factorization} of an integer.  According to the {\em Fundamental Theorem of Number Theory}, every integer $n$ with $n \geq 2$ is either a (positive) prime or can be expressed as a product of positive primes; furthermore, this factorization into primes is essentially unique (that is, there is only one factorization if we ignore the order of the prime factors or the possibility of using their negatives).  So, for example, the prime factorization of $18$ is $2 \cdot 3 \cdot 3$.  In general, the prime factorization of an integer $n$ with $n \geq 2$ can be written as $$n=\underbrace{p_1 \cdot p_1 \cdots p_1}_{\alpha_1} \cdot \underbrace{p_2 \cdot p_2 \cdots p_2}_{\alpha_2} \cdot \cdots \cdot \underbrace{p_k \cdot p_k \cdots p_k}_{\alpha_k},$$
where $p_1, p_2, \dots, p_k$ are the distinct prime factors of $n$ and $\alpha_1, \alpha_2, \dots, \alpha_k$ are positive integers.  This prime factorization is often turned into the {\em prime-power factorization} $$n=p_1^{\alpha_1} \cdot p_2^{\alpha_2} \cdot \cdots \cdot p_k^{\alpha_k};$$ for example, the prime-power factorization of $18$ is $2 \cdot 3^2$.
(If $n$ has only a single prime factor $p$ that appears $\alpha$ times in its factorization, we simply write $n=p^{\alpha}$.)  Note that $n=1$, of course, has no prime or prime-power factorization.

The prime and prime-power factorizations of integers are very useful and appear often in discussions.  For example, they allow us to quickly see if an integer $b$ is divisible by another integer $a$: this is the case if, and only if, each of the prime factors of $a$ appears in $b$ as well and at least as many times as it does in $a$.  For example, $$18=2 \cdot 3^2$$ is a factor of $$252=2^2 \cdot 3^2 \cdot 7,$$ but not of $$840000=2^7 \cdot 3 \cdot 5^5 \cdot 7.$$  
It also helps us compute the greatest common divisor and the least common multiple of two integers explicitly; for instance, given the prime-power factorizations of 252 and 840000 above, we immediately see that $$\mathrm{gcd}(252, 840000)=2^2 \cdot 3 \cdot 7$$ and $$\mathrm{lcm}(252, 840000)=2^7 \cdot 3^2 \cdot 5^5 \cdot 7.$$  (We mention, in passing, that $$\mathrm{gcd}(252, 840000) \cdot \mathrm{lcm}(252, 840000) =2^9 \cdot 3^3 \cdot 5^5 \cdot 7^2=252 \cdot 840000,$$ exemplifying the general fact that the product of the gcd and the lcm of two integers equals the product of the two integers.)

\vspace{.17in}

{\bf Exercises}

\begin{enumerate}

\item 
\begin{enumerate}

\item Find $d(1225)$.  

\item Suppose that $p$ and $k$ are positive integers and $p$ is a prime.  Find, in terms of $p$ and $k$, $d(p^k)$.

\item Find a formula for $d(n)$ for an arbitrary $n \in \mathbb{N}$ in terms of its prime factorization. 

\end{enumerate}

\item  Let us define, for a given $n \in \mathbb{N}$, $m \in \mathbb{N}$, and $i=0,1,\dots,m-1$, the set $$D_{m,i}(n)=\{ d \in D(n)  \mid d \equiv i \; (m) \}$$ and, if $D_{m,i}(n) \neq \emptyset$, let $$f_{m,i}(n)=\mathrm{min} D_{m,i}(n).$$

\begin{enumerate}

\item Find $f_{3,1}(1225)$ and $f_{3,2}(1225)$.

\item Explain why $f_{3,2}(n)$ is a prime number for every $n \in \mathbb{N}$ for which $D_{3,2}(n) \neq \emptyset$.

\item Is $f_{3,1}(n)$ also a prime number for every $n \in \mathbb{N}$ for which $D_{3,1}(n) \neq \emptyset$?

\item For each value of $i=0,1,2,3$, decide if $f_{4,i}(n)$ must be a prime or not whenever $D_{4,i}(n) \neq \emptyset$.

\item For each value of $i=0,1,2,3,4$, decide if $f_{5,i}(n)$ must be a prime or not whenever $D_{5,i}(n) \neq \emptyset$.

\end{enumerate}

\end{enumerate}

\section{Multiplicative number theory} \label{multnth}

The branch of number theory commonly referred to as {\em multiplicative number theory} deals with the various ways that integers can be factored into products of other integers.  The Fundamental Theorem of Number Theory, discussed above, plays a key role.  The prime factorization and the prime-power factorization of an integer are only two of the many different factorizations; here we discuss some others that we will use later.

First, a common generalization of the prime factorization and the prime-power factorization: the so-called {\em primary factorization}.  Indeed, the prime factorization $$n=\underbrace{p_1 \cdot p_1 \cdots p_1}_{\alpha_1} \cdot \underbrace{p_2 \cdot p_2 \cdots p_2}_{\alpha_2} \cdot \cdots \cdot \underbrace{p_k \cdot p_k \cdots p_k}_{\alpha_k}$$ and prime-power factorization $$n=p_1^{\alpha_1} \cdot p_2^{\alpha_2} \cdot \cdots \cdot p_k^{\alpha_k}$$ of an integer $n \geq 2$ can be considered the two extremes of the (potentially) many different  primary factorizations of the form
$$n=n_1 \cdot n_2 \cdot \cdots \cdot n_k$$ where each factor $n_i$ is a product of (one or more) prime  powers with base $p_i$ (here $i=1,2,\dots,k$).  For example, the number $n=18$ has only two primary factorizations, the prime factorization $(2) \cdot (3 \cdot 3)$ and the prime-power factorization $(2) \cdot (3^2)$, but a number such as $n=840000=2^7 \cdot 3 \cdot 5^5 \cdot 7$ has many, for example, $$840000=(2 \cdot 2^3 \cdot 2^3) \cdot (3) \cdot (5 \cdot 5^4) \cdot (7)$$ and  $$840000=(2^2 \cdot 2^2 \cdot 2^3) \cdot (3) \cdot ( 5 \cdot 5 \cdot 5^3) \cdot (7).$$  (Our parentheses indicate factors $n_1$, $n_2$, etc.)   

One can enumerate the number of primary factorizations of a given integer $n$, as follows.  Let us first consider the case when $n$ is a prime power itself; for example, let us examine $n=32=2^5$.  It is easy to see that $2^5$ has seven primary factorizations: $$2 \cdot 2 \cdot 2 \cdot 2 \cdot 2=2 \cdot 2 \cdot 2 \cdot 2^2=2 \cdot 2 \cdot 2^3=2 \cdot 2^2 \cdot 2^2=2 \cdot 2^4=2^2 \cdot 2^3=2^5.$$
More generally, if $n=p^{\alpha}$ for some prime $p$ and positive integer $\alpha$, then the number of primary factorizations of $n$ agrees with the number of ways that $\alpha$ can be written as the sum of positive integers (where the order of the terms is irrelevant).  Denoting this quantity by $p(\alpha)$, we find the following values. \label{p(n)values}
$$\begin{array}{|c||c|c|c|c|c|c|c|c|c||} \hline
\alpha & 1 & 2 & 3 & 4 & 5 & 6 & 7 & 8 & 9  \\ \hline
p(\alpha) & 1 & 2 & 3 & 5 & 7 & 11 & 15 & 22 & 30 \\ \hline
\end{array}$$
For example, $p(5)=7$ as the ways to write 5 as the sum of positive integers are $$1+1+1+1+1=1+1+1+2=1+1+3=1+2+2=1+4=2+3=5.$$  Consequently,  any number of the form $p^5$ with prime base $p$ has seven primary factorizations. 

The function $p(\alpha)$ is called the {\em partition function}; it plays an important role in mathematics in various ways, but, unfortunately, there is no closed formula for it.  For more values and information, see sequence A000041 in the On-Line Encyclopedia of Integer Sequences \cite{OEIS}. \label{p(n)}

In general, we can see that if the prime-power factorization of $n$ is $$n=p_1^{\alpha_1} \cdot p_2^{\alpha_2} \cdot \cdots \cdot p_k^{\alpha_k},$$ then the number of primary factorizations of $n$ equals $$p(\alpha_1) \cdot p(\alpha_2) \cdot \cdots \cdot p(\alpha_k).$$  Therefore, $$n=840000=2^7 \cdot 3 \cdot 5^5 \cdot 7$$ has $$p(7) \cdot p(1) \cdot p(5) \cdot p(1)=15 \cdot 1 \cdot 7 \cdot 1=105$$ different primary factorizations.

Related to primary factorizations, we have the so-called invariant factorizations.  We say that $$n=n_1 \cdot n_2 \cdot \cdots \cdot n_r$$ is an {\em invariant factorization} of the integer $n \geq 2$, if either $r=1$, or $r \geq 2$ with $n_1 \geq 2$, and  $n_i|n_{i+1}$ holds for each $i=1,2,\dots,r-1$.  For example, $18$ has two invariant factorizations, $3 \cdot 6$ and $18$ itself; $840000$, however, has many more, for instance,  
$$840000=(2) \cdot (2^3 \cdot 5) \cdot (2^3 \cdot 3 \cdot 5^4 \cdot  7)$$ and  $$840000=(2^2 \cdot 5) \cdot (2^2 \cdot 5) \cdot (2^3 \cdot 3 \cdot 5^3 \cdot 7).$$  (As before, our parentheses indicate factors $n_1$, $n_2$, etc.)

There is a nice one-to-one correspondence between primary factorizations and invariant factorizations.  Producing a primary factorization from an invariant factorization is easy: one can simply order and group the prime powers involved according to their prime bases.  For example, the primary factorization that we get from 
$$(2) \cdot (2^3 \cdot 5) \cdot (2^3 \cdot 3 \cdot 5^4 \cdot  7)$$ is $$(2 \cdot 2^3 \cdot 2^3) \cdot (3) \cdot (5 \cdot 5^4) \cdot (7).$$

To get an invariant factorization from a primary factorization, start by setting  the largest factor of the invariant factorization equal to the product of the largest factors of each of the factors in the primary factorization, follow that by the second largest factors, and so on.  For example, for the primary factorization $$(2^2 \cdot 2^2 \cdot 2^3) \cdot (3) \cdot ( 5 \cdot 5 \cdot 5^3) \cdot (7),$$ we see that the largest invariant factor equals $2^3 \cdot 3 \cdot 5^3 \cdot 7$, the next one is $2^2 \cdot 5$, and then $2^2 \cdot 5$ again, yielding the invariant factorization $$(2^2 \cdot 5) \cdot (2^2 \cdot 5) \cdot (2^3 \cdot 3 \cdot 5^3 \cdot 7).$$

Primary factorizations and invariant factorizations will enable us to classify all finite abelian groups---see Chapter \ref{IntroGrTh}.

\vspace{.17in}

{\bf Exercises}

\begin{enumerate}

\item Find $p(10)$.

\item  \begin{enumerate}

\item  How many primary factorizations and how many invariant factorizations does 72 have?

\item Find all primary factorizations and invariant factorizations of 72. 

\end{enumerate}

\end{enumerate}

\section{Additive number theory} \label{addnth}

In contrast to multiplicative number theory that deals with ways in which positive integers factor into products of other positive integers, {\em additive number theory} is concerned with ways that positive integers can be expressed as sums of certain other positive integers.   

In a typical setting, one is given a set $A$ of positive integers, and asks whether every positive integer can be written as a sum of terms all in $A$.  In further variations of this question, one may restrict the total number of terms in the sum or the number of times that a particular element of $A$ may occur in a sum.  

Regarding this latter restriction: the two most frequent variations are where we allow any element of $A$ to appear an arbitrary number of times and when each element may occur only at most once.  For a given $n \in \mathbb{N}$, $A \subseteq \mathbb{N}$ and $H\subseteq \mathbb{N}$, we introduce the notation $p(n,A,H)$ to denote the number of ways that $n$ can be written as a sum of elements of $A$, where the total number of terms in the sum must be an element of $H$, but there is no restriction on the number of times that elements may appear; similarly, $p\hat{\;}(n,A,H)$ will denote the number of those sums where the total number of terms in the sum must be an element of $H$, but where each element of $A$ may appear at most once.  

For instance, one can verify that there are three ways to write 11 as a sum of either three or four terms, where each term is 2, 3, 6, or 7:
$$7+2+2, \; \; 6+3+2, \; \; 3+3+3+2,$$
thus $$p(11, \{2,3,6,7\}, \{3,4\})=3,$$ but  $$p\hat{\;}(11, \{2,3,6,7\}, \{3,4\})=1,$$ as only the sum $6+3+2$ has distinct terms.

A more familiar example is the case when $A=\mathbb{N}$ and $H=\mathbb{N}$; that is, the terms may be arbitrary positive integers, and we have no restrictions on the number of terms in the sum or the number of times that a particular term may appear.  In this case, we get the partition function $p(n)$ introduced on page \pageref{p(n)}, so $$p(n,\mathbb{N}, \mathbb{N})=p(n).$$  For instance, $$p(5,\mathbb{N}, \mathbb{N})=p(5)=7,$$ as 5 can be written as 
\begin{eqnarray*}
5 & = & 5 \\
& = & 4+1 \\
& = & 3+2 \\
& = & 3+1+1 \\
& = & 2+2+1 \\
& = & 2+1+1+1 \\
& = & 1+1+1+1+1 \\
\end{eqnarray*}
If we allow only sums where the terms are distinct, we get $p\hat{\;}(n,\mathbb{N}, \mathbb{N})$; in the case of $n=5$, we see that $p\hat{\;}(5,\mathbb{N}, \mathbb{N})=3$ as only the first three sums above contain distinct terms.  As we noted on page \pageref{p(n)}, no closed formula exists for $p(n)$, though its values for small $n$ and various estimates for higher $n$ are known.  There is also no closed formula for $p\hat{\;}(n,\mathbb{N}, \mathbb{N})$, although we should mention the remarkable fact that the number of ways to partition a positive integer into distinct parts equals the number of its partitions into odd parts; that is, $$p\hat{\;}(n,\mathbb{N}, \mathbb{N})=p(n,\mathbb{O}, \mathbb{N}).$$ (We denote the set of odd positive integers here by $\mathbb{O}$.)  Indeed, we see that $5$ has three partitions comprised of only odd terms. 

Keeping $A=\mathbb{N}$ and putting no restrictions on the number of times that elements of $A$ may appear in a sum, but limiting the total number of possible terms to at most $s$ (for some $s \in \mathbb{N}$), we get the function $p(n,\mathbb{N}, [1,s])$.  (As usual, we denote the set $H=\{1,2,\dots,s\}$ by $[1,s]$.)  For example, we see that $p(5,\mathbb{N}, [1,3])=5$ since five of the sums above contain no more than three terms.  We have no closed formula for $p(n,\mathbb{N}, [1,s])$, though we have another remarkable identity:
$$p(n,\mathbb{N}, [1,s])=p(n, [1,s],\mathbb{N});$$   that is, for all positive integers $n$ and $s$, the number of partitions of $n$ into at most $s$ parts equals the number of partitions of $n$ into parts that do not exceed $s$.  For example, just like we had five ways to partition 5 into at most three parts, we have five partitions where each part is at most 3.  These are just some of the many amazing identities involving the partition function.  

In closing this section, we mention some other interesting examples.  The famous {\em Goldbach Conjecture} asserts that every positive even number is the sum of at most two positive primes; denoting the set of positive primes by $\mathbb{P}$ and the set of positive even integers by $\mathbb{E}$, this conjecture can be presented to say that 
$$p(n,P,[1,2]) \geq 1$$ holds for all $n \in \mathbb{E}$.  While several partial results have been achieved, the Goldbach Conjecture remains as one of the oldest and most famous unsolved problems in mathematics.   Very recently, Helfgott\index{Helfgott, H. A.} solved the related {\em Weak Goldbach Conjecture}: that $$p(n,P,[1,3]) \geq 1$$ for all $n \in \mathbb{O}$ and $n>1$ (see \cite{Hel:2015a}; this work has not been published in a refereed journal yet).

For a fixed positive integer $k$, let $S_k$ denote the set $\{1^k, 2^k, 3^k, \dots \}$ of all $k$-th powers of positive integers.  {\em Waring's Problem} asks for the smallest positive integer $s$, usually denoted by $g(k)$, for which $$p(n,S_k,[1,s]) \geq 1$$ for all $n \in \mathbb{N}$; that is, the smallest $s$ for which it is true that every positive integer can be written as at most $s$  $k$-th powers.  (Equivalently, we may include $0^k$ as an element of $S_k$ in which case $g(k)$ is the minimum value of $h$ for which $$p(n,S_k,\{h\}) \geq 1$$ for all $n \in \mathbb{N}$.)  For example, $g(2)=4$: as the {\em Four Squares Theorem} or {\em Lagrange's Theorem} asserts, every positive integer can be written as the sum of at most four positive squares (or, equivalently, as the sum of exactly 4 nonnegative squares).  To see that $g(2)$ cannot be less than 4, observe that $n=7$ (and infinitely many other $n$) indeed requires four squares.  We also know that $g(3)=9$, $g(4)=19$, and many other values of $g(k)$, but the question of finding all values is still open today.

Our final famous example is the so-called {\em Money Changing Problem}.  Suppose that we are given a finite set $A$ of relatively prime positive integers; we want to find the largest positive integer $n$, denoted by $f(A)$, for which $$p(n,A,\mathbb{N})=0;$$ that is, the largest positive integer $n$ that cannot be written as a sum of elements of $A$.  (The reason for the name of the problem should be obvious.)  Note that if the elements of $A$ are not relatively prime, then $f(A)$ does not exist; this is the case, for example, if all elements of $A$ are even, as only even numbers can be partitioned into even terms.  On the other hand, one can show that if the elements of $A$ are relatively prime, then $f(A)$, called the {\em Frobenius number of $A$}, exists.  For example, with $A=\{5,8\}$, we find that $f(A)=27$ as 27 cannot be written as a sum of 5's and 8's, but every number greater than 27 can be:
\begin{eqnarray*}
28 & = & 4 \cdot 5 + 1 \cdot 8, \\
29 & = & 1 \cdot 5 + 3 \cdot 8, \\
30 & = & 6 \cdot 5 + 0 \cdot 8, \\
31 & = & 3 \cdot 5 + 2 \cdot 8, \\
32 & = & 0 \cdot 5 + 4 \cdot 8, \\
33 & = & 5 \cdot 5 + 1 \cdot 8, \\
34 & = & 2 \cdot 5 + 3 \cdot 8,
\end{eqnarray*}  and so on.  In general, if $A=\{a,b\}$ (and $\mathrm{gcd}(a,b)=1$) then $f(A)=ab-a-b$; similar formulas for the case when $|A|>2$ are not known.

\vspace{.17in}

{\bf Exercises}

\begin{enumerate}

\item 
\begin{enumerate}

\item Verify that, as listed on page \pageref{p(n)values}, $p(7,\mathbb{N}, \mathbb{N})=p(7)=15.$

\item Verify that $p\hat{\;}(7,\mathbb{N}, \mathbb{N})=p(7,\mathbb{O}, \mathbb{N}).$

\item Verify that $p(7,\mathbb{N}, [1,3])=p(7, [1,3],\mathbb{N}).$

\end{enumerate}

\item

\begin{enumerate}
  \item Prove that the Weak Goldbach Conjecture follows from the Goldbach Conjecture.
  \item Prove that the Weak Goldbach Conjecture implies that $$p(n,P,[1,4]) \geq 1$$ holds for all $n \in \mathbb{N}$ and $n>1$.
\end{enumerate}

\item \begin{enumerate}

\item Prove that $g(3) \geq 9$ by finding a positive integer $n$ for which $p(n,S_3,[1,8])=0$.

\item Prove that $g(4) \geq 19$ by finding a positive integer $n$ for which $p(n,S_4,[1,18])=0$.

\end{enumerate}

\item Suppose that a fast-food chain sells chicken nuggets in packages of 6, 9, and 20.  Find the largest positive integer $n$ for which we are not able to buy exactly $n$ pieces.  

\end{enumerate}

\chapter{Combinatorics} \label{IntroComb}

Explaining what combinatorics is about may be simple since it deals with objects and techniques that are quite familiar to most people.   Yet, combinatorics is not easy to define precisely.  At the fundamental level, combinatorics deals with the questions related to counting the number of elements in a given set; a bit more precisely, combinatorics deals with {\em discrete structures}: sets---with, perhaps, some specific characteristics---whose elements can be listed and enumerated, as opposed to sets whose elements vary continuously and cannot be put in a list. 

  Moreover, beyond its object of study, combinatorics can be characterized by its methods.  Typically, by the {\em combinatorial method} we mean a relatively basic---but, perhaps, surprisingly deep and far reaching---argument using some relatively elementary tools, rather than the application of sophisticated and elaborately developed machinery.  That is what makes combinatorics so highly applicable and why it serves as a very elegant and accessible branch of study in the mathematics curriculum.

In this section we introduce some of the concepts and methods of combinatorics that we will need later.

\section{Basic enumeration principles} \label{0.2.1.1}

Enumeration---or, simply, counting---is probably one of our earliest intellectual pursuits, and it is a ubiquitous task in everyday life.  The principles of enumeration are also what several branches of mathematics are based on, especially probability theory and statistics.  In this section we briefly discuss elementary enumeration techniques.

A typical enumeration problem asks us to determine the {\em size} of a set: the size of a set $A$, denoted by $|A|$, is the number of elements in $A$.  Clearly, each set has either finite or infinite size.  Here we focus on finite sets only.

Most enumeration questions can be reduced to one of two fundamental principles: the Addition Rule and the Multiplication Rule.  According to the {\em Addition Rule}, if $A$ and $B$ are disjoint finite sets, then we have $$|A \cup B|=|A| + |B|.$$  More generally, if $A_1$, $A_2, \dots, A_n$ are pairwise disjoint finite sets ($n \in \mathbb{N}$), then we have $$\left| A_1 \cup \cdots \cup A_n \right|=|A_1|+ \cdots + |A_n|.$$  The {\em Multiplication Rule} says that for arbitrary finite sets $A$ and $B$, we have $$|A \times B|=|A| \cdot |B|,$$  and more generally, for arbitrary finite sets $A_1$, $A_2, \dots, A_n$ ($n \in \mathbb{N}$), we have $$\left| A_1 \times \cdots \times A_n \right|=|A_1|\cdot \cdots \cdot |A_n|.$$

Observe that the Addition Rule---unlike the Multiplication Rule---requires that the sets be pairwise disjoint.   A more general formula treats the case when our sets are not (or not known to be) pairwise disjoint: for finite sets $A$ and $B$ we can verify that $$|A \cup B|  = |A| + |B|- |A \cap B|;$$ indeed, to count the elements in the union of $A$ and $B$, adding the sizes of $A$ and $B$ together  would double-count the elements that are in both $A$ and $B$, so we need to subtract the number of elements in the intersection of $A$ and $B$.  
Similarly, for finite sets $A$, $B$, and $C$, we have 
$$|A \cup B \cup C|  = |A| + |B| + |C|- |A \cap B| - |A \cap C| - |B \cap C| + |A \cap B \cap C|.$$
The situation gets more complicated as the number of sets increases; while a precise statement (called the {\em Inclusion--Exclusion Rule}) is readily available, we will not state it here.  Instead, we just point out that, in general, for arbitrary sets $A_1$, $A_2, \dots, A_n$ we have  $$\left| A_1 \cup \cdots \cup A_n \right| \leq |A_1|+ \cdots + |A_n|.$$

An often-used consequence of this inequality is the {\em Pigeonhole Principle}, which says that, if we have $$\left| A_1 \cup \cdots \cup A_n \right| > kn$$ for some nonnegative integer $k$, then there must be at least one index $i \in \{1,\dots,n\}$ for which $|A_i| \geq k+1$.  To paraphrase: if more than $kn$ pigeons happen to sit in $n$ holes, then at least one hole must have at least $k+1$ pigeons in it.   Consequently, for example, we see that in a set of 101 positive integers one can always find 11 (or more) that share their last digits; similarly, among a group of 3000 people there is always a group of at least nine that share the same birthday.

While the counting principles we reviewed here may seem rather elementary, they have far-reaching consequences.  We present some in the exercises below.

\vspace{.17in}

{\bf Exercises}

\begin{enumerate}

\item Exhibit the Inclusion--Exclusion formula for four sets; that is, find an expression for the size of the union of four sets in terms of the sizes of their various intersections.

\item Prove that however we place seven points inside an 8-by-9 rectangle, we can always find

\begin{enumerate}

\item a pair whose distance is at most 5, and
\item three that form a triangle of area at most 12.
\end{enumerate}

\item

Let $S$ be a set of 100 distinct positive integers.  Which of the following statements are true?

\begin{enumerate}

\item If each element of $S$ is at most 198, then $S$ must contain two elements that are relatively prime.
\item If each element of $S$ is at most 199, then $S$ must contain two elements that are relatively prime.
\item If each element of $S$ is at most 198, then $S$ must contain two elements so that one is divisible by the other.
\item If each element of $S$ is at most 199, then $S$ must contain two elements so that one is divisible by the other.

\end{enumerate}

\end{enumerate}

\section{Counting lists, sequences, sets, and multisets} \label{0.2.1.2}

Before we discuss the four main counting questions in mathematics, we review some familiar terminology and notations, and introduce some new ones.  Recall that, for a given set $A$ and positive integer $m$, an element $(a_1,a_2,\dots,a_m)$ of $A^m$ is called a sequence of length $m$.  The order of the terms in the sequence matters; for example, the sequence $(2,3,4,5)$ of integers is different from $(3,2,4,5)$.  On the other hand, a subset of $A$ of size $m$  is simply a collection of $m$ of its elements, where two subsets are considered equal without regard of the order in which the terms are listed; for example, $\{2,3,4,5\}$ and $\{3,2,4,5\}$ are equal subsets of the set of integers.  Recall also that a set remains unchanged if we choose to list some of its elements more than once; for example, the sets $\{2,3,3,5\}$, $\{2,3,5,5\}$, and $\{2,3,5\}$ are all equal, while the sequences $(2,3,3,5)$, $(2,3,5,5)$, and $(2,3,5)$ are all different.  Thus, we can consider sets as two-fold relaxations of sequences: we don't care about the order in which the elements are listed, nor do we care how many times the elements are listed.

It will be useful for us to introduce two other objects.  First, we say that a sequence $(a_1,a_2,\dots,a_m)$ of elements of a set $A$ is a {\em list}, if the $m$ terms are pairwise distinct.  Thus, in a list, the order of the elements still matters, but each element is only allowed to appear once.  For example, the sequence $(2,3,4,5)$ is a list, but $(2,3,3,5)$ is not.  Conversely, in a so-called {\em multiset} $[a_1,a_2,\dots,a_m]$ of size $m$, the order of the elements $a_1, a_2, \dots, a_m$ of $A$ does not matter (as it is the case with sets), but elements may appear repeatedly (as they may in sequences).  For example, the multisets $[2,3,3,5]$, $[2,3,5,5]$,  and $[2,3,5]$ are all different, but $[2,3,3,5]$ is still the same as $[2,5,3,3]$.   

Given a set $A$ and a positive integer $m$, we are interested in counting the number of $m$-sequences (sequences of length $m$), $m$-lists  (lists of length $m$), $m$-multisubsets  (multisubsets of size $m$), and $m$-subsets   (subsets of size $m$) of $A$.  The schematic summary of these four terms is given in the following table.
\label{combintableterms}
\begin{center}
\begin{tabular}{l||c|c|}
& order matters & order does not matter \\ \hline \hline
& & \\ 
elements distinct & $m$-lists & $m$-sets \\ 
& & \\ \hline
& & \\ 
elements may repeat & $m$-sequences & $m$-multisets  \\  
& & \\ \hline
\end{tabular}
\end{center}

Obviously, if $|A|<m$, then $A$ has neither $m$-lists nor $m$-subsets.  If $|A|=m$, then the (only) $m$-subset of $A$ is $A$ itself, while, as we will soon see, if $|A|=m$, then $A$ has $m!$ $m$-lists.  For other situations, we introduce the following notations.  

Suppose that $n$ is a nonnegative integer and $m$ is a positive integer.  We define the {\em rising factorial $m$-th power} and the {\em falling factorial $m$-th power} of $n$ to be  
$$n^{\overline{m}}=n(n+1)\cdots (n+m-1)$$ and $$n^{\underline{m}}=n(n-1)\cdots (n-m+1),$$
respectively.  For example, we have $$10^{\overline{3}}=10 \cdot 11 \cdot 12=1320$$ and $$10^{\underline{3}}=10 \cdot 9 \cdot 8=720.$$  
Analogously to $n^0=1$ and $0!=1$, we extend these notations with $$n^{\underline{0}}=1 \; \mbox{ and } \; n^{\overline{0}}=1$$ for arbitrary nonnegative integers $n$.

Furthermore, we introduce the notations ${n \choose m}$ (pronounced ``$n$ choose $m$'') and $\left[ {n \atop m} \right]$ (pronounced ``$n$ multichoose $m$''): For nonnegative integers $m$ and $n$, $${n \choose m}= \frac{n^{\underline{m}}}{m!} = \frac{n(n-1)\cdots(n-m+1)}{m!}$$ and $$\left[ {n \atop m} \right] = \frac{n^{\overline{m}}}{m!} = \frac{n(n+1)\cdots(n+m-1)}{m!}.$$  

It is well known that these quantities always denote integers.  The values of ${n \choose m}$, also known as {\em binomial coefficients}, are exhibited in {\em Pascal's Triangle}; here we tabulate some of these values in a table format.  (Note that, when $m>n$, the formula above yields ${n \choose m}=0$; keeping the traditional shape of Pascal's Triangle, we omitted these entries from the table below.)

\label{Pascaltable}
\begin{tabular}{|c||c|c|c|c|c|c|c|c|} \hline
${n \choose m}$ & m=0 & m=1 & m=2 & m=3 & m=4 & m=5 & m=6 & m=7 \\ \hline \hline
n=0 & 1&  &  &  &  &  &  &   \\ \hline
n=1 & 1& 1&  &  &  &  &  &   \\ \hline
n=2 & 1& 2& 1&  &  &  &  &   \\ \hline
n=3 & 1& 3& 3& 1&  &  &  &   \\ \hline
n=4 & 1& 4& 6& 4& 1&  &  &   \\ \hline
n=5 & 1& 5& 10 & 10 & 5& 1&  &   \\ \hline
n=6 & 1& 6& 15& 20& 15& 6& 1&   \\ \hline
n=7 & 1& 7& 21& 35& 35& 21& 7& 1 \\ \hline
\end{tabular}

Observe that, since $$\frac{n(n-1)\cdots(n-m+1)}{m!}=\frac{n(n-1)\cdots(m+1)}{(n-m)!}$$ (which we can check by cross-multiplying), we have the identity $${n \choose m}={n \choose n-m},$$ expressing the fact that the rows in Pascal's Triangle are ``palindromic.''  The explanation for the term ``binomial coefficient'' will be clear once we discuss the Binomial Theorem  below.

The first few values of $\left[ {n \atop m} \right] $ are as follows.

\begin{tabular}{|c||c|c|c|c|c|c|c|c|} \hline
$\left[ {n \atop m} \right] $ & m=0 & m=1 & m=2 & m=3 & m=4 & m=5 & m=6 & m=7 \\ \hline \hline
n=1 & 1& 1& 1& 1& 1& 1& 1& 1   \\ \hline
n=2 & 1& 2& 3& 4& 5& 6& 7& 8 \\ \hline
n=3 & 1& 3& 6& 10& 15& 21& 28& 36 \\ \hline
n=4 & 1& 4& 10& 20& 35& 56& 84& 120  \\ \hline
n=5 & 1& 5& 15& 35& 70& 126& 210& 330  \\ \hline
n=6 & 1& 6& 21& 56& 126& 252& 462& 792  \\ \hline
n=7 & 1& 7& 28& 84& 210& 462& 924& 1716 \\ \hline
\end{tabular}
 
As we can see, the two tables contain the same data---values are just shifted: the entries in column $m$ in the first table are moved up by $m-1$ rows in the second table.  Indeed, since for integers $n$ and $m$ we clearly have $$n^{\overline{m}}=n(n+1) \cdots (n+m-1)=(n+m-1)(n+m-2) \cdots n=(n+m-1)^{\underline{m}},$$ we see that values of $\left[ {n \atop m} \right]$ can be expressed via the more-often used binomial coefficients as $$\left[ {n \atop m} \right]={n +m-1 \choose m}.$$

We are now ready to ``size up'' our four main configurations.  The {\em Enumeration Theorem} says that, if  $A$ is a set of size $n$ and $m$ is a positive integer,  then
\begin{itemize}
  \item the number of $m$-sequences of $A$ is $n^m$;
  \item the number of $m$-lists of $A$ is $n^{\underline{m}}$;
  \item the number of $m$-multisubsets of $A$ is $\left[ {n \atop m} \right]$; and  
  \item the number of $m$-subsets of $A$ is ${n \choose m}$.
\end{itemize}

Note that, when $n<m$, then $n^{\underline{m}}=0$ and ${n \choose m}=0$, in accordance with the fact that $A$ has no $m$-lists and no $m$-subsets in this case.  If $n=m$, then $n^{\underline{m}}=m!$ and ${n \choose m}=1$; indeed, in this case $A$ has $m!$ $m$-lists while its only $m$-subset is itself.

The enumeration techniques discussed above are often employed to determine the number of choices one has for selecting or arranging a given number of elements from a given set or collection of sets.  For example, the Addition Rule and the Multiplication Rule can be interpreted to say that, given boxes labeled $A_1$, $A_2, \dots, A_n$, if box $A_i$ contains $m_i$ distinct objects ($i=1,2,\dots,n$), then there are $$m_1+m_2+ \cdots + m_n$$ ways to choose an object from one of the boxes, and there are $$m_1 \cdot m_2 \cdot  \cdots \cdot  m_n$$ ways to choose an object from each of the boxes.  In a similar manner, the four basic enumeration functions of the Enumeration Theorem are sometimes called ``choice functions;'' the following table summarizes our results for the number of ways to choose $m$ elements from a given set of $n$ elements.

\label{combintable}
\begin{center}
\begin{tabular}{l||c|c|}
& order matters & order does not matter \\ \hline \hline
& & \\ 
elements distinct & $n^{\underline{m}}$ & ${n \choose m}$ \\ 
& & \\ \hline
& & \\ 
elements may repeat & $n^m$ & $\left[ {n \atop m} \right]$  \\  
& & \\ \hline
\end{tabular}
\end{center}

An important example for enumeration problems, one that we will refer to often, is to count the number of positive integer solutions to an equation of the form 
$$x_1+x_2+\cdots+x_m = h;$$ that is, to find, for a given $h \in \mathbb{N}$, the number of $m$-sequences of $\mathbb{N}$ with the property that the entries in the sequence add up to $h$.  (Note that we are counting sequences: the order of the terms does matter.)  We can visualize this question by imagining a segment of length $h$ inches with markings at all integer inches (that is, at 1, 2, and so on, all the way to $h-1$); our task is then to find the number of ways this segment can be broken into $m$ pieces at $m-1$ distinct markings: the lengths of the $m$ parts created will correspond, in order, to $x_1, x_2, \dots, x_m$.  By the Enumeration Theorem, the number of ways that this can be done, and therefore the 
number of positive integer solutions to our equation,  is \label{seqaddtoh} ${h-1 \choose m-1}$.  As a variation, one can easily prove (see one of the exercises below) that the number of nonnegative integer solutions to the same equation equals \label{seqaddtohnonneg} $\left[ {h+1 \atop m-1} \right]={m+h-1 \choose h}$.

\vspace{.17in}

{\bf Exercises}

\begin{enumerate}

\item Find the number of 

\begin{enumerate}

\item 4-sequences

\item 4-lists

\item 4-multisubsets

\item 4-subsets

\end{enumerate}

of a set of size 7.  Exhibit one example for each question.

\item Above we have shown that the number of positive integer solutions to an equation 
$$x_1+x_2+\cdots+x_m = h$$ equals ${h-1 \choose m-1}$.  
  Here we use three different approaches to find the number of nonnegative solutions to the equation.  

\begin{enumerate}
\item Modify the argument used for the number of positive integer solutions to prove that the number of nonnegative integer solutions equals $\left[ {h+1 \atop m-1} \right]$.

\item Explain why the number of nonnegative integer solutions to $$x_1+x_2+\cdots+x_m = h$$ is the same as the number of positive integer solutions to $$x_1+x_2+\cdots+x_m = h+m,$$ and use this fact to get that the result is ${h+m-1 \choose m-1}$.

\item Explain why the number of nonnegative integer solutions to $$x_1+x_2+\cdots+x_m = h$$ is the same as the number of ways one can place $h$ identical objects into $m$ distinct boxes, and use this fact to get that the result is $\left[ {m \atop h} \right]$.

\item Verify algebraically that the results of the three previous parts are the same.

\end{enumerate}

\end{enumerate}

\section{Binomial coefficients and Pascal's Triangle}  \label{0.2.2}

In this section we discuss some of the many famous and interesting properties of the so-called {\em binomial coefficients}, that is, the quantities ${n \choose m}$.  First, let us explain the reason for the name.  

A closer look at the rows of the table of entries for ${n \choose m}$, exhibited earlier, reveals, in order, the coefficients of the various terms in the expansion of the power $(a+b)^n$ (here $a$ and $b$ are arbitrary real numbers and $n$ is a nonnegative integer).  For example, the entries in row 4 of the table are 1, 4, 6, 4, and 1; indeed, the power $(a+b)^4$ expands as
$$(a+b)^4=a^4+4a^3b+6a^2b^2+4ab^3+b^4.$$  

We can easily explain this coincidence, as follows.  When using the distributive law to expand the expression $$(a+b)^n=(a+b) \cdots (a+b),$$ we arrive at a sum of products of $n$ factors, where each factor is either $a$ or $b$.  Using the commutative property of multiplication, each term can be arranged so that the $a$s (if any) all come before the $b$s (if any).  We can then collect ``like'' terms; that is, terms of the form $a^{n-m} b^{m}$ for the same $m=0,1,\dots,n$.  The number of such terms clearly equals the number of those $n$-sequences of the set $\{a,b\}$ that contain exactly $n-m$ $a$s and $m$ $b$s, which, by the Enumeration Theorem, is exactly ${n \choose m}$.  

This result is known as (Newton's) Binomial Theorem, and can be stated in general as the identity 
$$(a+b)^n=\sum_{m=0}^n {n \choose m} a^{n-m} b^{m}.$$
As the name implies, ${n \choose m}$ is indeed a ``binomial coefficient.''

The first few entries for ${n \choose m}$ (and, therefore, for $\left[ {n \atop m} \right]$ as well) are tabulated in {\em Pascal's Triangle} below.
\begin{center}
\begin{tabular}{ccccccccccccccc}
&&&&& & & 1 & & & & &&&\\
&&&&& & 1 & & 1 & & & &&&\\
&&&&& 1& & 2& & 1& & &&&\\
&&&&1& & 3& & 3& & 1& &&&\\
&&&1&& 4& & 6& & 4& & 1&&&\\
&&1&&5& & 10& & 10& & 5& &1&&\\
&1&&6&& 15& & 20& & 15& & 6&&1& \\
1&&7&&21& & 35& & 35& & 21& &7&&1\\
\end{tabular}
\end{center}
\label{pascalt} 
We can read off values of ${n \choose m}$, as follows.  If we label the rows, the ``left'' diagonals, and the ``right'' diagonals 0, 1, 2, etc. (we start with 0), then ${n \choose m}$ appears as the entry where row $n$ and right diagonal $m$ intersect.  For example, we see that ${6 \choose 3}=20$.

The binomial coefficients possess many interesting properties.  We have already mentioned the fact that the rows are `palindromic':
$${n \choose m} = {n \choose n-m}.$$  Another important property is known as Pascal's Identity: $${n \choose m}={n-1 \choose m}+{n-1 \choose m-1}.$$  This identity provides us, actually, with the easiest way to enumerate binomial coefficients: each entry in Pascal's Triangle is simply the sum of the two entries above it.  We can, thus, quite quickly find the next row:
$$1   \; \; \;  \; \; \;  \; \; \;  \; \; \; 8   \; \; \;  \; \; \;  \; \; \;  \; \; \; 28    \; \; \;  \; \; \;  \; \; \;  \; \; \; 56   \; \; \;  \; \;  \; \; \;  \; \; \; \; 70    \;  \; \; \; \; \;  \; \; \;  \; \; \; 56   \;  \; \; \; \; \;  \; \; \;  \; \; \; 28   \; \; \; \; \; \;  \; \; \;  \; \; \; 8    \; \; \; \; \; \;  \; \; \;  \; \; \; 1$$ 

Another interesting property of Pascal's Triangle is that the sum of the entries in each row add up to a power of 2: $${\displaystyle {n \choose 0}+{n \choose 1}+\cdots+{n \choose n}=2^n}.$$  Note that this identity follows directly from the Binomial Theorem (take $a=b=1$).  Similarly, (by taking $a=1$ and $b=2$) we have  
$${\displaystyle {n \choose 0}\cdot 2^0+{n \choose 1}\cdot 2^1+\cdots+{n \choose n}\cdot 2^n=3^n}.$$

Of the numerous other interesting properties of Pascal's Triangle, we list only two more: \label{binomialidentities}
$${\displaystyle {n-1 \choose m}+{n-2 \choose m-1}+\cdots+{n-m-1 \choose 0}={n \choose m}},$$ expressing the fact that the entries in each NW-SE diagonal, above a certain row, add to an entry in the next row.  Adding up numbers on NE-SW diagonals yields
$${\displaystyle {n-1 \choose m-1}+{n-2 \choose m-1}+\cdots+{m-1 \choose m-1}={n \choose m}}.$$  (See the Exercises below for proofs.)

We will use each of these identities later.

\vspace{.17in}

{\bf Exercises}

\begin{enumerate}

\item What identity of binomial coefficients arises from using the Binomial Theorem for evaluating $(1-1)^n$?  Verify your identity for $n=6$ and $n=7$.

\item Suppose that $n$ and $m$ are positive integers, and suppose that $m \leq n$.   
 
\begin{enumerate}

\item We set $A=\{1,2,\dots,n\}$; furthermore, we let $A_0$ be the set of $m$-subsets of $A$ that do not contain 1, let $A_1$ be the set of $m$-subsets of $A$ that contain 1 but do not contain 2, let $A_2$ be the set of $m$-subsets of $A$ that contain 1 and 2 but do not contain 3, and so on.  Prove the identity $${\displaystyle {n-1 \choose m}+{n-2 \choose m-1}+\cdots+{n-m-1 \choose 0}={n \choose m}}$$ by considering $A_0 \cup A_1 \cup A_2 \cup \cdots \cup A_m$.

\item Prove the identity $${\displaystyle {n-1 \choose m-1}+{n-2 \choose m-1}+\cdots+{m-1 \choose m-1}={n \choose m}}$$ using similar techniques as in part (a).

\end{enumerate}

\end{enumerate}

\section{Some recurrence relations} \label{0.2.3}

Let us return to the binomial coefficients discussed in the previous section.  Rather than looking at Pascal's Triangle, let's arrange their values in a more convenient table format: 

\begin{tabular}{|c||c|c|c|c|c|c|c|c|} \hline
p(j,k) & k=0 & k=1 & k=2 & k=3 & k=4 & k=5 & k=6 & k=7 \\ \hline \hline
j=0 & 1& 1& 1& 1& 1& 1& 1& 1   \\ \hline
j=1 & 1& 2& 3& 4& 5& 6& 7& 8 \\ \hline
j=2 & 1& 3& 6& 10& 15& 21& 28& 36 \\ \hline
j=3 & 1& 4& 10& 20& 35& 56& 84& 120  \\ \hline
j=4 & 1& 5& 15& 35& 70& 126& 210& 330  \\ \hline
j=5 & 1& 6& 21& 56& 126& 252& 462& 792  \\ \hline
j=6 & 1& 7& 28& 84& 210& 462& 924& 1716 \\ \hline
j=7 & 1& 8& 36& 120& 330& 792& 1716& 3432 \\ \hline
\end{tabular}

Here $p(j,k)$ denotes the entry in row $j$ and column $k$; we have $$p(j,k)=\left[ j+1 \atop k \right]= {j+k \choose k}.$$

As we noted before, Pascal's Identity, together with the values in the top row and the left-most column in the table, determine all entries {\em recursively}: we simply need to add the (previously determined) values directly above and directly to the left of the desired entry.  In fact, we can {\em define} the function $p(j,k)$ recursively by the {\em recurrence relation} $$p(j,k)=p(j-1,k)+p(j,k-1)$$ and the {\em initial conditions} that $p(j,0)=1$ for all $j \in \mathbb{N}_0$ and $p(0,k)=1$ for all $k \in \mathbb{N}_0$. 

Let us consider a variation where the function $a(j,k)$ is defined by the initial conditions $a(j,0)=1$ for all $j \in \mathbb{N}_0$ and $a(0,k)=1$ for all $k \in \mathbb{N}_0$ and by the recursive relation $$a(j,k)=a(j-1,k)+a(j-1,k-1)+a(j,k-1).$$  The first few values of the function are as follows.    

\begin{tabular}{c||c|c|c|c|c|c|c|} \label{a(j,k)table}
$a(j,k)$ & $k=0$  & $k=1$  & $k=2$  & $k=3$  & $k=4$  & $k=5$  & $k=6$ \\ \hline \hline
$j=0$ &1& 1& 1& 1& 1& 1& 1 \\ \hline
$j=1$ &1& 3& 5& 7& 9& 11& 13 \\ \hline
$j=2$ &1& 5& 13& 25& 41& 61& 85 \\ \hline
$j=3$ &1& 7& 25& 63& 129& 231& 377 \\ \hline
$j=4$ &1& 9& 41& 129& 321& 681& 1289 \\ \hline
$j=5$ &1& 11& 61& 231& 681& 1683& 3653 \\ \hline
$j=6$ &1& 13& 85& 377& 1289& 3653& 8989 \\ \hline
\end{tabular}

The numbers in this table are called {\em Delannoy numbers},\index{Delannoy, H.} named after the French amateur mathematician who introduced them in the nineteenth century in \cite{Del:1895a}.
According to its recurrence relation, the Delannoy number $a(j,k)$ is the sum of not only the entries directly above and to the left, but the entry in the ``above-left'' position as well.  Delannoy numbers---like any two-dimensional array---can be turned into a sequence by listing entries by its anti-diagonals; this sequence is given as A008288 in \cite{OEIS}.  (The sequence of entries in various columns can be found in \cite{OEIS} as well.)  In the next section we shall see an interesting interpretation of Delannoy numbers.

Changing the initial conditions, we next define the function $c(j,k)$ by the initial conditions $c(j,0)=0$ for all $j \in \mathbb{N}$ and $c(0,k)=1$ for all $k \in \mathbb{N}_0$ and by the (same) recursive relation $$c(j,k)=c(j-1,k)+c(j-1,k-1)+c(j,k-1).$$  The first few values of this function are as follows.    

\begin{tabular}{c||c|c|c|c|c|c|c|}
$c(j,k)$ & $k=0$  & $k=1$  & $k=2$  & $k=3$  & $k=4$  & $k=5$  & $k=6$ \\ \hline \hline
$j=0$ &1& 1& 1& 1& 1& 1& 1 \\ \hline
$j=1$ &0& 2& 4& 6& 8& 10& 12 \\ \hline
$j=2$ &0& 2& 8& 18& 32& 50& 72 \\ \hline
$j=3$ &0& 2& 12& 38& 88& 170& 292 \\ \hline
$j=4$ &0& 2& 16& 66& 192& 450& 912 \\ \hline
$j=5$ &0& 2& 20& 102& 360& 1002& 2364 \\ \hline
$j=6$ &0& 2& 24& 146& 608& 1970& 5336 \\ \hline
\end{tabular}    

The numbers given by this table can be found in sequence form at A266213 in \cite{OEIS}.

The functions $a(j,k)$ and $c(j,k)$ are strongly related; it is not difficult to reduce each one to the other, as we now show.  

Consider first the function $c(j,k)$.  A quick glance at the tables above suggests that for $j \geq 1$ and $k \geq 1$, the entry $c(j,k)$ is the sum of entries $a(j,k-1)$ and $a(j-1,k-1)$:
$$c(j,k)=a(j,k-1)+a(j-1,k-1).$$  We will prove this by induction.  We see that the equation holds for $j=1$ and for $k=1$; we then use the defining recursions for both $c$ and $a$, as well as our inductive hypothesis, for $j \geq 1$ and $k \geq 1$ to write
\begin{eqnarray*}
c(j,k) & = & c(j-1,k)+c(j-1,k-1)+c(j,k-1) \\
& = & a(j-1,k-1)+a(j-2,k-1) + \\
&  &  + a(j-1,k-2)+a(j-2,k-2) + \\
&  &  + a(j,k-2)+a(j-1,k-2) \\
& = & a(j,k-1)+a(j-1,k-1),
\end{eqnarray*}  
as claimed.

Using the recursion for $a(j,k)$ once more, we may rewrite this identity as
$$c(j,k)=a(j,k)-a(j-1,k),$$ from which we get
$$a(j,k)=c(j,k)+a(j-1,k).$$

We can then use this identity to express $a$ in terms of $c$:
\begin{eqnarray*}
a(j,k) & = & c(j,k)+a(j-1,k) \\
& = & c(j,k)+c(j-1,k)+a(j-2,k) \\
& = & c(j,k)+c(j-1,k)+c(j-2,k)+a(j-3,k) \\
& = & \dots \\
& = & c(j,k)+c(j-1,k)+\cdots+ c(1,k)+a(0,k) \\
& = & c(j,k)+c(j-1,k)+\cdots+ c(1,k)+c(0,k). 
\end{eqnarray*}
In summary, we have:

\begin{prop} \label{functionsac}
For the functions $a(j,k)$ and $c(j,k)$, defined recursively above for all nonnegative integers $j$ and $k$, we have
$$c(j,k)=a(j,k-1)+a(j-1,k-1)=a(j,k)-a(j-1,k)$$ for all $j, k \in \mathbb{N}$  
and  
$$a(j,k)=c(j,k)+c(j-1,k)+\cdots+ c(1,k)+c(0,k)$$ for all $j, k \in \mathbb{N}_0$.
\end{prop}

While recursive expressions are quite helpful, direct formulae, if they exist, would be even more useful, particularly when the variables are large.  For example, the binomial coefficients can be easily computed via the function $p(j,k)$ defined above, but it is good to know that $$p(j,k)={j+k \choose k}=\frac{(j+k)!}{j! \cdot k!}.$$   Although formulae for $a(j,k)$ and $c(j,k)$ are not so direct, we have the following expressions.

\begin{prop} \label{functionsacdirect}  For all nonnegative integers $j$ and $k$ we have
$$a(j,k)=\sum_{i \geq 0} {j \choose i} {k \choose i} 2^i$$ and
$$c(j,k)=\sum_{i \geq 0} {j-1 \choose i-1} {k \choose i} 2^i.$$
\end{prop}

For a proof of Proposition \ref{functionsacdirect}, see page \pageref{proofoffunctionsacdirect}.  Note that, while the summations seem to include infinitely many terms, all but finitely many are zero.  Note also that including $i=0$ in the last sum is only relevant if $j=0$.  (Here we use the convention that ${j-1 \choose -1}$ equals $1$ for $j=0$ and 0 if $j>0$.)  

It may seem a bit strange that, while $p(j,k)$ and $a(j,k)$ are similarly defined---with the only difference being that $p$ relies on a double recursion while $a$ uses a triple recursion---they yield very different formulae.   We gain some insight by the following consideration.  Suppose that we are given $j+k$ distinct (for example, numbered) balls, and that $j$ of them are green and $k$ are yellow.  Recall that $p(j,k)={j+k \choose k}$, a quantity expressing the number of ways one can select $k$ balls from the collection of these $j+k$ balls.  Now if $i$ of the $k$ balls selected are green and the other $k-i$ are yellow, then the number of such choices, by the Multiplication Rule, equals ${j \choose i} \cdot {k \choose i}$; summing over all possible values of $i$, we get $$p(j,k)={j+k \choose j}= \sum_{i \geq 0} {j \choose i} {k \choose i},$$ a form more closely reminiscent to the one for $a(j,k)$ in Proposition \ref{functionsacdirect}.

We will make frequent use of the quantities $a(j,k)$ and $c(j,k)$. Without going into detail here, we mention, for example, that
\begin{itemize}
\item for an $s$-spanning set of size $m$ in a group of order $n$ (see Chapter \ref{ChapterSpanning}), we have $n \leq a(m,s)$,
\item for a $B_h$ set over $\mathbb{Z}$ of size $m$ in a group of order $n$  (see Chapter \ref{ChapterSidon}), we have $n \geq c(h,m)$, and
\item for a $t$-independent set of size $m$ in a group of order $n$  (see Chapter \ref{ChapterZerosumfree}), we have 
\begin{itemize} \item $n \geq c(m,\tfrac{t+1}{2})$ if $t$ is odd and $t>1$, and 
\item $n \geq a(m,\tfrac{t}{2})$ if $t$ is even.
\end{itemize}
\end{itemize}

We will discuss each of these bounds in the relevant chapters of the book.

\vspace{.17in}

{\bf Exercises}

\begin{enumerate}

\item Find $a(7,7)$ and $c(7,7)$ using 

\begin{enumerate}

\item their recursive definitions;

\item Proposition \ref{functionsacdirect}.

\end{enumerate}

\item Suppose that $m \in \mathbb{N}$.  Express $a(m,3)$ and $c(m,3)$ as polynomial functions in $m$.

\end{enumerate}

\section{The integer lattice and its layers}  \label{0.2.4}

One of the most often discussed combinatorial objects---and one that we will frequently rely on---is the $m$-dimensional {\em integer lattice} $$\mathbb{Z}^m=\underbrace{\mathbb{Z} \times \mathbb{Z} \times \cdots \times \mathbb{Z}}_m,$$ consisting of all points (or vectors) $(\lambda_1,\lambda_2,\dots,\lambda_m)$ with integer coordinates. (Here $m \in \mathbb{N}$ is called the {\em dimension} of the lattice.)  Of course, the 1-dimensional integer lattice is simply the set of integers $\mathbb{Z}$, while the points in $\mathbb{Z}^2$ are arranged in an infinite (2-dimensional) grid in the plane, and $\mathbb{Z}^3$ can be visualized as an infinite grid in 3-space.  (For $m \geq 4$, geometric visualization of $\mathbb{Z}^m$ is not convenient.)

A {\em layer} of the integer lattice is defined as the collection of points with a given fixed {\em norm}; that is, for a given nonnegative integer $h$, the $h$th layer of $\mathbb{Z}^m$ is defined as
$$\mathbb{Z}^m (h) = \{ (\lambda_1,\lambda_2,\dots,\lambda_m) \in \mathbb{Z}^m  \mid |\lambda_1|+|\lambda_2|+\cdots+|\lambda_m|=h\}.$$  Obviously, $\mathbb{Z}^m (0)$ consists of a single point, the origin.  We can also easily see that $\mathbb{Z}^m (1)$ consists of the points of the $m$-dimensional lattice with all but one coordinate equal to $0$ and the remaining coordinate equal to $1$ or $-1$; there are exactly $2m$ such points.

Describing $\mathbb{Z}^m (h)$ explicitly gets more complicated as $h$ increases, however.  For instance, we find that
$$\mathbb{Z}^2 (3)=\{(0,\pm 3), (\pm 1, \pm 2), (\pm 2, \pm 1), (\pm 3, 0)\},$$ and
$$\mathbb{Z}^3 (2)=\{(0,0,\pm 2), (0,\pm 2, 0), (\pm 2, 0,0), (0,\pm 1, \pm 1), (\pm 1, 0,\pm 1), (\pm 1, \pm 1,0)\}.$$  The twelve points of $\mathbb{Z}^2 (3)$ lie on the boundary of a square in the plane (occupying the four vertices and  two points on each edge), and the eighteen points of $\mathbb{Z}^3 (2)$ are on the surface of an octahedron (the six vertices and the midpoints of the twelve edges).

We can derive a formula for the size of $\mathbb{Z}^m (h)$, as follows.  For $i=0,1,2,\dots,m$, let $I_i$ be the set of those elements of $\mathbb{Z}^m (h)$ where exactly $i$ of the $m$ coordinates are nonzero.  How many elements are in $I_i$?  We can choose which $i$ of the $m$ coordinates are nonzero in ${m \choose i}$ ways.  Next, we choose the absolute values of these nonzero coordinates: since the sum of these $i$ positive integers equals $h$, we have ${h-1 \choose i-1}$ choices (see page \pageref{seqaddtoh}).  Finally, each of these $i$ coordinates can be positive or negative, and therefore $$|I_i|= {m \choose i} {h-1 \choose i-1} 2^i.$$ Summing now for $i$ yields $$|\mathbb{Z}^m (h)|= \sum_{i=0}^h  {m \choose i} {h-1 \choose i-1} 2^i;$$ since for $i >h$ the terms vanish, we may write this as 
$$|\mathbb{Z}^m (h)|= \sum_{i \geq 0}  {m \choose i} {h-1 \choose i-1} 2^i.$$
Here we recognize the expression for the size of $\mathbb{Z}^m (h)$ as the quantity $c(h,m)$, discussed in detail in Section \ref{0.2.3}.

In our investigations later, we will also consider certain restrictions of $\mathbb{Z}^m (h)$.  In some cases, we will look at the part of $\mathbb{Z}^m (h)$ that is in the ``first quadrant;'' that is, the subset $\mathbb{N}_0^m (h)$ of $\mathbb{Z}^m (h)$ that contains only those points that contain no negative coordinates:
$$\mathbb{N}_0^m (h) = \{ (\lambda_1,\lambda_2,\dots,\lambda_m) \in \mathbb{N}_0^m  \mid \lambda_1+\lambda_2+\cdots+\lambda_m=h\}.$$  So, for example, 
$$\mathbb{N}_0^2 (3)=\{(0,3), (1, 2), (2, 1), (3, 0)\}$$ and
$$\mathbb{N}_0^3 (2)=\{(0,0,2), (0,2, 0), (2, 0,0), (0,1, 1), ( 1, 0, 1), ( 1, 1,0)\}.$$  

We can enumerate $\mathbb{N}_0^m (h)$ by observing that it is nothing but the set of $m$-sequences of $\mathbb{N}_0$ with the property that the entries in the sequence add up to $h$; as we have seen on page \pageref{seqaddtohnonneg}, this set has size
$$|\mathbb{N}_0^m (h)| = {m+h-1 \choose h}.$$

Two other, frequently appearing, cases occur when we restrict $\mathbb{Z}^m (h)$ or $\mathbb{N}_0^m (h)$ to those points where the absolute value of the coordinates are not more than 1.  (These points lie within a cube of side length 2 centered at the origin.)
We denote these sets by  
$\hat{\mathbb{Z}}^m (h)$ and $\hat{\mathbb{N}}_0^m (h)$, respectively; namely, we have
$$\hat{\mathbb{Z}}^m (h)=\{ (\lambda_1,\lambda_2,\dots,\lambda_m) \in \{-1,0,1\}^m  \mid |\lambda_1|+|\lambda_2|+\cdots+|\lambda_m|=h\}$$ and 
$$\hat{\mathbb{N}}_0^m (h)=\{ (\lambda_1,\lambda_2,\dots,\lambda_m) \in \{0,1\}^m  \mid \lambda_1+\lambda_2+\cdots+\lambda_m=h\}.$$
So, for example,
$$\hat{\mathbb{Z}}^3 (2)=\{(0,\pm 1, \pm 1), (\pm 1, 0,\pm 1), (\pm 1, \pm 1,0)\}$$
and
$$\hat{\mathbb{N}}_0^3 (2)=\{(0,1, 1), ( 1, 0, 1), ( 1, 1,0)\},$$ but we have $\hat{\mathbb{Z}}^2 (3)=\emptyset$ and $\hat{\mathbb{N}}_0^2 (3)=\emptyset$.
It is easy to see that the sizes of these sets are given by
$$|\hat{\mathbb{Z}}^m (h)|={m \choose h} 2^h$$ and 
$$|\hat{\mathbb{N}}_0^m (h)|={m \choose h}.$$

Often, rather than considering a single layer $\mathbb{Z}^m (h)$ of the integer lattice, we will study the union of several of them.  
Since the layers are pairwise disjoint, for a given range $H \subseteq \mathbb{N}_0$ of norms we have $$\left|\bigcup_{h \in H} \mathbb{Z}^m (h)\right|= \sum_{h \in H} |\mathbb{Z}^m (h)|;$$ we can similarly just add the sizes of $\mathbb{N}_0^m (h)$, $\hat{\mathbb{Z}}^m (h)$, and $\hat{\mathbb{N}}_0^m (h)$ for all $h \in H$.  Most often, we will consider $H$ consisting of 
\begin{itemize}

\item a single norm $h$ (with $h \in \mathbb{N}_0$),
\item a range $[0,s]=\{0,1,2,\dots,s\}$ (with $s \in \mathbb{N}_0$), or
\item allow all possible norms (i.e. have $H=\mathbb{N}_0$).

\end{itemize}
  The following table summarizes what we can say, using some of the identities we have seen earlier, about the size of $$\Lambda^m (H) = \{ (\lambda_1,\lambda_2,\dots,\lambda_m) \in \Lambda^m  \mid |\lambda_1|+|\lambda_2|+\cdots+|\lambda_m| \in H \}$$ for these choices of $H \subseteq \mathbb{N}_0$ and our four exemplary sets $\Lambda \subseteq \mathbb{Z}$.

$$\begin{array}{l||c|c|c|}  \label{sumsizeformula}
 |\Lambda^m (H)| & H=\{h\}  & H=[0,s] & H=\mathbb{N}_0 \\ \hline \hline

& && \\

\Lambda=\mathbb{N}_0  & {m +h -1 \choose h} & {m + s \choose s} & \infty  \\ 

& && \\ \hline 

& && \\
 
\Lambda=\mathbb{Z} & c(h,m)=\sum_{i \geq 0}  {m \choose i} {h-1 \choose i-1} 2^i & a(m,s)=\sum_{i \geq 0}  {m \choose i} {s \choose i} 2^i &  \infty \\ 

& && \\ \hline 

& && \\

\Lambda=\{0,1\} & {m \choose h} & \sum_{h \in H} {m \choose h}  & 2^m \\ 

& && \\ \hline 

& && \\

\Lambda=\{-1,0,1\}&  {m \choose h} 2^h & \sum_{h \in H} {m \choose h} 2^h & 3^m \\ 

& && \\ \hline

\end{array}$$

We need to compute the size of an additional set that we use later.  Namely, we want to find the number of lattice points that are strictly on one side of one of the coordinate planes of the $m$-dimensional space; that is, the size of the set
$$\mathbb{Z}^m (h)_{k+} = \{ (\lambda_1,\lambda_2,\dots,\lambda_m) \in \mathbb{Z}^m  \mid |\lambda_1|+|\lambda_2|+\cdots+|\lambda_m|=h, \lambda_k >0\}.$$
For example, the lattice points of the layer $\mathbb{Z}^2 (3)$ that are to the right of the $y$-axis are $$\mathbb{Z}^2 (3)_{1+}=\{(1, \pm 2), (2, \pm 1), (3, 0)\}.$$
Note that $|\mathbb{Z}^m (h)_{k+}|$ is the same for any $k=1,2,\dots,m$; here we calculate $|\mathbb{Z}^m (h)_{1+}|$.

As before, we let, for each $j=1,\dots,m$, $I_j$ denote the set of those elements of $\mathbb{Z}^m (h)_{1+}$ where exactly $j$ of the $m$ coordinates are nonzero.  (Note that $I_0=\emptyset$.)  How many elements are in $I_j$?  Here we can choose which $j$ of the $m$ coordinates are nonzero in ${m-1 \choose j-1}$ ways (since we must have $\lambda_1>0$).  Next, we choose the absolute values of these nonzero coordinates: since the sum of these $j$ positive integers equals $h$, we have ${h-1 \choose j-1}$ choices.  Finally, $j-1$ of these coordinates can be positive or negative, and therefore $$|I_j|= {m-1 \choose j-1} {h-1 \choose j-1} 2^{j-1}.$$ Summing now for $j$ yields $$|\mathbb{Z}^m (h)_{1+}|= \sum_{j \geq 1}  {m-1 \choose j-1} {h-1 \choose j-1} 2^{j-1}.$$  We can replace $j-1$ by $i$; this yields
$$|\mathbb{Z}^m (h)_{1+}|= \sum_{i \geq 0}  {m-1 \choose i} {h-1 \choose i} 2^{i},$$ and therefore
$$|\mathbb{Z}^m (h)_{k+}|= \sum_{i \geq 0}  {m-1 \choose i} {h-1 \choose i} 2^{i}=a(m-1,h-1)$$ for every $k=1,2,\dots,m$.  Indeed, the set $\mathbb{Z}^2 (3)_{1+}$ featured above consists of $a(1,2)=5$ points.

\vspace{.17in}

{\bf Exercises}

\begin{enumerate}

\item We have already evaluated the entries in the column of $H=\{h\}$ in the table on page \pageref{sumsizeformula}; here we verify the rest.   Prove each of the following.

\begin{enumerate}

\item $|\mathbb{Z}^m([0,s])|=a(m,s)$

\item $|\mathbb{N}_0^m([0,s])|={m+s \choose s}$

\item $|\hat{\mathbb{Z}}^m(\mathbb{N}_0)|=3^m$

\item $|\hat{\mathbb{N}}_0^m(\mathbb{N}_0)|=2^m$

\end{enumerate}

\item For each set below, first find the size of the set, then list all its elements.

\begin{enumerate}

\item $\mathbb{Z}^2([0,3])$

\item $\mathbb{N}_0^2([0,3])$

\item $\hat{\mathbb{Z}}^2([0,3])$

\item $\hat{\mathbb{N}}_0^2([0,3])$

\item $\mathbb{Z}^2(3)_{1+}$

\end{enumerate}

\end{enumerate}

\chapter{Group theory} \label{IntroGrTh}

A group is arguably the most important structure in abstract mathematics.  In essence, a {\em group} is any set of objects---for example, numbers, functions, vectors, etc.---combined with a binary operation---such as addition, multiplication, composition, etc.---satisfying certain fundamental properties.  More precisely, a set $G$ and an operation $\ast$ form a {\em group}, if each of the following four properties holds.
\begin{itemize}
\item {\em Closure property}: for any pair of elements $a$ and $b$ of $G$, $a\ast b$ is also in $G$.
\item {\em Associative property}: for any $a$, $b$, and $c$ in $G$, we have $(a\ast b) \ast c=a \ast(b \ast c)$.
\item {\em Identity property}: there is an element $z$ in $G$ so that $a \ast z=z \ast a=a$ holds for any element $a$ of $G$.
\item {\em Inverse property}: for any element $a$ of $G$, there exists an element $\overline{a}$, also in $G$, for which $a \ast \overline{a}=\overline{a} \ast a=z$.
\end{itemize} 
Here we discuss a special class of groups: abelian groups.  The group is said to be {\em abelian} if, in addition, the following holds.
\begin{itemize}
\item {\em Commutative property}: for any $a$ and $b$ in $G$, we have $a \ast b=b \ast a$.
\end{itemize}
One can show that, in a group, the identity element (denoted by $z$ above) is unique, and each element $a$ of the group has its unique inverse $\overline{a}$.

Abelian groups---named after the Norwegian mathematician Niels Abel (1802--1829)---play a central role in most branches of mathematics; our goal here is to investigate some of their fascinating number theoretic properties.  Since our focus in this book is on additive combinatorics, we restrict our attention to {\em additive groups}, where the operation $\ast$ is addition, denoted as $+$; we will also write $0$ for the identity $z$ and $-a$ for the inverse $\overline{a}$ of $a$.  

The additive groups most familiar to us are probably
\begin{itemize} 

\item the set of integers (nonnegative and negative whole numbers), denoted by $\mathbb{Z}$; 

\item the set of rational numbers (fractions of integers), denoted by $\mathbb{Q}$; and 

\item the set of real numbers (finite or infinite decimals), denoted by $\mathbb{R}$.

\end{itemize}  
Other well-studied abelian groups include the set of vectors in $n$-dimensional space, the set of $n$-by-$m$ real matrices, and the set of real polynomials.  Another example for an abelian group is the set of even integers (including positive and negative even integers and zero); however, the set of odd integers is not a group, since the closure and zero properties fail (neither $1+1$ nor $0$ is odd).  The sets of positive and nonnegative integers, denoted by $\mathbb{N}$ and $\mathbb{N}_0$, respectively, are also not groups: they both fail the negative inverse property (and $\mathbb{N}$ even fails the zero property).  

Our examples for groups above are all {\em infinite} groups; that is, they have infinitely many elements.  The groups we intend to study here, however, are {\em finite} groups: those that have only a finite number of elements.  Thus the title of this chapter really should be an introduction to the theory of finite abelian groups.

\section{Finite abelian groups} \label{0.1.1}

The simplest family of finite abelian groups are the cyclic groups: the {\em cyclic group} of size (or {\em order}) $n$, denoted by $\mathbb{Z}_n$, can be defined, as follows.  The elements of $\mathbb{Z}_n$ are the nonnegative integers up to $n-1$: $$\mathbb{Z}_n = \{0,1,2, \dots, n-1\}.$$  Addition is performed ``mod $n$;'' that is, for elements $a$ and $b$ of $\mathbb{Z}_n$, the sum $a+b$ is the remainder of $a+b$ when divided by $n$.  For example, in $\mathbb{Z}_{10}$ we have $9+4=3$, representing the fact that when we add two integers whose last decimal digits are 9 and 4, respectively, then their sum will have a last digit of 3.  In $\mathbb{Z}_{12}$, we have $9+4=1$: if our evening guests arrive at 9 p.m. and plan to stay for 4 hours, then they will leave at 1 a.m.  (Of course, in both $\mathbb{Z}_{10}$ and $\mathbb{Z}_{12}$ we still have $5+3=8$.)  

We call $\mathbb{Z}_n$ cyclic because if we add 1 repeatedly to itself, then within $n$ steps we run through each of the elements in the group, after which the cycle repeats itself.  The same holds for any other element $a$ that is relatively prime to $n$; for example, in $\mathbb{Z}_{10}$ we have
$$\{\lambda \cdot 3  \mid \lambda =0,1,2,\dots,9\} = \{0,3,6,9,2,5,8,1,4,7\}=\mathbb{Z}_{10}.$$  (The notation $\lambda \cdot 3$, for $\lambda \in \mathbb{N}_0$, stands for the sum of $\lambda$ terms with each term being 3.)
If $a$ and $n$ are not relatively prime, then we still get a cycle, it's just that the cycle will be shorter as we won't run through all $n$ elements; for example, the multiples of 2 in $\mathbb{Z}_{10}$ only yield the set $\{0,2,4,6,8\}$.  

The length of the cycle that an element $a$ of a group $G$ generates is called the {\em order} of $a$ in $G$, and is denoted by $\mathrm{ord}_G(a)$ or simply $\mathrm{ord}(a)$.  For example, in $\mathbb{Z}_{10}$ we have $\mathrm{ord}(1)=10$, $\mathrm{ord}(2)=5$, $\mathrm{ord}(3)=10$, and so on.  According to {\em Lagrange's Theorem}, \label{lagrange} in a group of order $n$, the order of any element is a divisor of $n$; for example, in $\mathbb{Z}_{10}$, only orders 1, 2, 5, and 10 are possible.

The set of elements in $G$ that have order $d$ is denoted by $\mathrm{Ord}(G,d)$.  For example, we have 
$$\mathrm{Ord}(\mathbb{Z}_{10},1) = \{0\}, $$
$$\mathrm{Ord}(\mathbb{Z}_{10},2) =\{5\}, $$
$$\mathrm{Ord}(\mathbb{Z}_{10},5)=\{2,4,6,8\},$$ and
$$\mathrm{Ord}(\mathbb{Z}_{10},10)=\{1,3,7,9\}.$$

From cyclic groups, we can build up other finite abelian groups using direct sums.  The {\em direct sum} (which is also called the {\em direct product}) of the groups $G_1$ and $G_2$, denoted by $G_1 \times G_2$, consists of all ordered pairs of the form $(a_1,a_2)$ where $a_1$ is any element of $G_1$ and $a_2$ is any element of $G_2$; formally, 
$$G_1 \times G_2 = \{ (a_1,a_2) \mbox{   } | \mbox{   } a_1 \in G_1, a_2 \in G_2 \}.$$  If $G_1$ and $G_2$ have orders $n_1$ and $n_2$, respectively, then the order of $G_1 \times G_2$ is $n_1n_2$.  For example, $\mathbb{Z}_2 \times \mathbb{Z}_5$ has ten elements:  
$$\mathbb{Z}_2 \times \mathbb{Z}_5 = \{ (0,0), (0,1), (0,2), (0,3), (0,4), (1,0), (1,1), (1,2), (1,3), (1,4)  \}.$$ 

We can also define the direct sum of more than two groups: the direct sum of finite abelian groups $G_1, G_2, \dots, G_r$ consists of the ordered $r$-tuples $(a_1,a_2,\dots,a_r)$ where $a_i$ is any element of $G_i$ (here $i=1,2,\dots,r$).  As a special case, if each component in the direct sum is the same group, then we often use exponential notation.  For example, $\mathbb{Z}_2 \times \mathbb{Z}_2 \times \mathbb{Z}_2$ is denoted by $\mathbb{Z}_2^3$, and consists of eight elements:
$$\mathbb{Z}_2^3 = \{ (0,0,0), (0,0,1), (0,1,0), (0,1,1), (1,0,0), (1,0,1), (1,1,0), (1,1,1)  \}.$$

We add and subtract the elements of the direct sum component-wise.  For example, in $\mathbb{Z}_2 \times \mathbb{Z}_5$ one has $(1,3)+(1,4)=(0,2)$, and  in $\mathbb{Z}_2^3$ we have $(1,0,1)-(1,1,0)=(0,1,1)$.

We must note that not all such direct sum compositions result in new types of groups---as we will see in the next section.

\vspace{.17in}

{\bf Exercises}

\begin{enumerate}

\item 

\begin{enumerate} \item For each positive integer $d$, find the set $\mathrm{Ord}(\mathbb{Z}_{18},d).$

\item For each positive integer $d$, find the set $\mathrm{Ord}(\mathbb{Z}_3 \times \mathbb{Z}_6,d)$.

\end{enumerate}

\item Exhibit the complete addition table of the groups $\mathbb{Z}_3 \times \mathbb{Z}_4$ and $\mathbb{Z}_2^3$.
 
\end{enumerate}

\section{Group isomorphisms} \label{0.1.2}

Let us examine again two of the groups mentioned above: $\mathbb{Z}_{10}$ and $\mathbb{Z}_2 \times \mathbb{Z}_5$.  Both of these groups have ten elements---furthermore, they form exactly the same structure as we now explain.  

Consider the following table.
$$\begin{array}{||c||c|c|c|c|c|c|c|c|c|c||} \hline \hline
\mathbb{Z}_{10} &  0 & 1  & 2 & 3 & 4  &  5 &  6 &  7 &  8 & 9     \\ \hline 
\mathbb{Z}_2 \times \mathbb{Z}_5 & (0,0) & (1,1) & (0,2) & (1,3) & (0,4) & (1,0) & (0,1) & (1,2) & (0,3) & (1,4)  \\ \hline \hline
\end{array}$$

The table exhibits a correspondence between the elements of the two groups with a very special property: if we add two elements in one group and then add the corresponding elements in the other group, then the two sums will also correspond to each other.  For example: 
$$\begin{array}{cccccc} 
\mathbb{Z}_{10}: & 9 & + & 4 & = & 3      \\
 & \updownarrow  & & \updownarrow &  & \updownarrow \\  
\mathbb{Z}_2 \times \mathbb{Z}_5: & (1,4)& + & (0,4) & = & (1,3)
\end{array}$$

We thus find that each entry in the addition table of $\mathbb{Z}_{10}$ (the ten by ten table that lists all possible pairwise sums of elements) corresponds to the appropriate entry in the addition table of $\mathbb{Z}_2 \times \mathbb{Z}_5$.  Therefore, the two groups are essentially the same; using standard terminology, we say that they are {\em isomorphic}---a fact that we denote as follows:
$$\mathbb{Z}_{10} \cong  \mathbb{Z}_{2} \times \mathbb{Z}_{5}.$$
More generally, one can prove that, if $n_1$ and $n_2$ are relatively prime integers, then $$\mathbb{Z}_{n_1} \times \mathbb{Z}_{n_2} \cong \mathbb{Z}_{n_1n_2}.$$  For example, the groups $\mathbb{Z}_{5} \times \mathbb{Z}_{12}$, $\mathbb{Z}_{4} \times \mathbb{Z}_{15}$, $\mathbb{Z}_{3} \times \mathbb{Z}_{20}$, and $\mathbb{Z}_{3} \times \mathbb{Z}_{4} \times \mathbb{Z}_{5}$ are all isomorphic to $\mathbb{Z}_{60}$.

One can also show that, if $n_1$ and $n_2$ are not relatively prime, then $\mathbb{Z}_{n_1} \times \mathbb{Z}_{n_2}$ and $\mathbb{Z}_{n_1n_2}$ are not isomorphic.  For example, $\mathbb{Z}_{6} \times \mathbb{Z}_{10}$ is an example of a group of order 60 that is not isomorphic to $\mathbb{Z}_{60}$, since 6 and 10 are not relatively prime.  (To see why there cannot possibly be an isomorphism between these two groups, note that $\mathbb{Z}_{60}$ is cyclic, but $\mathbb{Z}_{6} \times \mathbb{Z}_{10}$ is not cyclic as none of its elements has order 60---see Section \ref{0.1.3} below.)    

Group isomorphism is a very important and useful concept: when two groups are isomorphic, it suffices to study one of them (either one).  For example, if $G_1 \cong G_2$, then the number of elements of a certain order in $G_1$ will be the same as the number of elements of that order in $G_2$, and the same kind of property holds for other important group characteristics.  This  reduces our task of studying all finite abelian groups to that of studying only those that have different {\em isomorphism types} (are pairwise non-isomorphic).

\vspace{.17in}

{\bf Exercises}

\begin{enumerate}

\item Show that the groups $\mathbb{Z}_{3} \times \mathbb{Z}_{4}$ and $\mathbb{Z}_{12}$ are isomorphic by finding an explicit correspondence between the elements of the two groups.

\item Explain why the groups $\mathbb{Z}_{3} \times \mathbb{Z}_{6}$ and $\mathbb{Z}_{18}$ are not isomorphic by considering the orders of the various elements in the two groups.
 
\end{enumerate}

\section{The Fundamental Theorem of Finite Abelian Groups} \label{0.1.3}

Suppose that one wishes to study all abelian groups of order 60, that is, those that have exactly 60 elements.  First, we wish to determine how many possible isomorphism types there are for abelian groups of order 60.  We have already seen that the groups $$\mathbb{Z}_{3} \times \mathbb{Z}_{4} \times \mathbb{Z}_{5}, \; \; \mathbb{Z}_{5} \times \mathbb{Z}_{12}, \; \; \mathbb{Z}_{4} \times \mathbb{Z}_{15}, \; \; \mathbb{Z}_{3} \times \mathbb{Z}_{20}, \; \; \mbox{and} \; \; \mathbb{Z}_{60}$$ are all isomorphic to one another.  Considering all other possible factorizations of 60 into a product of positive integers, we see that there are five other direct sums to examine: 
$$\mathbb{Z}_{2} \times \mathbb{Z}_{2} \times \mathbb{Z}_{3} \times \mathbb{Z}_{5}, \; \; \mathbb{Z}_{2} \times \mathbb{Z}_{3} \times \mathbb{Z}_{10},\;  \; \mathbb{Z}_{2} \times \mathbb{Z}_{5} \times \mathbb{Z}_{6}, \; \; \mathbb{Z}_{2} \times \mathbb{Z}_{30}, \; \; \mbox{and} \; \; \mathbb{Z}_{6} \times \mathbb{Z}_{10}.$$  
(We should note that, since we only study abelian groups here, changing the order of the terms in a direct sum will not change its isomorphism type.)
 
It is not hard to determine that these five direct sums are also all isomorphic to each other.  For example, to see that $$\mathbb{Z}_{2} \times \mathbb{Z}_{30} \cong \mathbb{Z}_{6} \times \mathbb{Z}_{10},$$ note that 5 and 6 are relatively prime, and thus $\mathbb{Z}_{2} \times \mathbb{Z}_{30} \cong \mathbb{Z}_{2} \times \mathbb{Z}_{5} \times \mathbb{Z}_{6},$ from which, since 2 and 5 are also relatively prime, we get $\mathbb{Z}_{6} \times \mathbb{Z}_{10}$.  Thus, among the ten possible direct sums, we have two different isomorphism types. 

We say that the isomorphism relation is an {\em equivalence relation}: one can partition the collection of all groups into {\em equivalence classes} (sometimes called {\em isomorphism classes}) where all groups within a class are equivalent to each other (but not to those in other classes).  As we have just seen, the collection of abelian groups of order 60 form two isomorphism classes.

In general, we can find the number of different isomorphism classes among the groups of order $n$, as follows: if $$n=p_1^{\alpha_1} \cdots  p_k^{\alpha_k}$$ is the prime factorization of $n$, then the number of different isomorphism classes among the groups of order $n$ equals
$$p(\alpha_1)\cdot \cdots \cdot p(\alpha_k)$$ where $p$ is the partition function introduced in Section \ref{multnth}.  For example, since $$60=2^2 \cdot 3^1 \cdot 5^1,$$ there are $$p(2) \cdot p(1) \cdot p(1)=2$$ different isomorphism types among groups of order 60.

We even have a convenient way of listing representatives of the different isomorphism types for a given order.  The {\em Fundamental Theorem of Finite Abelian Groups} asserts that for any finite abelian group $G$ of order at least 2, there are positive integers $r$ and $n_1, \dots, n_r$ such that 
$$G \cong \mathbb{Z}_{n_1} \times \mathbb{Z}_{n_2} \times \cdots \times \mathbb{Z}_{n_r};$$ furthermore, we can assume that $n_{i+1}$ is divisible by $n_i$ for $i=1,2,\dots, r-1$ and that $n_1 \geq 2$ (and therefore $n_i \geq 2$ for all $i$).  The factorization of $G$ above, which is unique, is called the {\em invariant decomposition} of $G$.  We say that $r$ is the {\em rank} of $G$, and the largest invariant factor, $n_r$, is called the {\em exponent} of $G$; the exponent clearly equals the length of the longest cycle in $G$.  We should also note that our treatment here includes the possibility that $r=1$, in which case we simply have $G \cong Z_{n}$ (with $n=n_1$), and the group is cyclic.

Therefore, we get a list of all possible isomorphism types of abelian groups of order $n$ directly from the list of invariant factorizations of $n$ (see Section \ref{multnth}).   For example, there are two isomorphism classes of abelian groups of order 60, as shown by the two invariant decompositions $\mathbb{Z}_{60}$ and $\mathbb{Z}_{2} \times \mathbb{Z}_{30}$.  The group $\mathbb{Z}_{60}$ has rank 1 (as does any cyclic group), and the group $\mathbb{Z}_{2} \times \mathbb{Z}_{30}$ has rank 2; every cycle in the latter group has length at most 30.  From Section \ref{multnth} we also know that there are exactly 105 different isomorphism types among abelian groups of order 840000.

We should also mention two other notable members of the two isomorphism classes for order 60: $\mathbb{Z}_{3} \times \mathbb{Z}_{4} \times \mathbb{Z}_{5}$ (which is isomorphic to $\mathbb{Z}_{60}$) and $\mathbb{Z}_{2} \times \mathbb{Z}_{2} \times \mathbb{Z}_{3} \times \mathbb{Z}_{5}$ (which is isomorphic to $\mathbb{Z}_{2} \times \mathbb{Z}_{30}$).  These direct sums involve (prime and) prime power orders only, they are thus called {\em primary decompositions}.  Like the invariant decomposition of $G$, the primary decomposition is also unique (up to the order of terms, of course).  

We can easily convert any decomposition of a finite abelian group into its invariant decomposition or its primary decomposition.  Let us consider the example $$G=\mathbb{Z}_{40} \times \mathbb{Z}_{50} \times \mathbb{Z}_{60} \times \mathbb{Z}_{70}.$$  We start with the primary decomposition as that is easier: since $40=8 \cdot 5$, $50=2 \cdot 25$, $60=4 \cdot 3 \cdot 5$, and $70=2 \cdot 5 \cdot 7$, the primary decomposition of $G$ is $$G \cong \mathbb{Z}_{2} \times \mathbb{Z}_{2} \times \mathbb{Z}_{4} \times \mathbb{Z}_{8} \times 
\mathbb{Z}_{3} \times \mathbb{Z}_{5} \times \mathbb{Z}_{5} \times \mathbb{Z}_{5} \times \mathbb{Z}_{25} \times \mathbb{Z}_{7},$$ which we can condense as 
$$G \cong \mathbb{Z}_{2}^2 \times \mathbb{Z}_{4} \times \mathbb{Z}_{8} \times 
\mathbb{Z}_{3} \times \mathbb{Z}_{5}^3 \times \mathbb{Z}_{25} \times \mathbb{Z}_{7}.$$  From this, we find that the invariant decomposition is  
$$G \cong \mathbb{Z}_{10} \times \mathbb{Z}_{10} \times \mathbb{Z}_{20} \times \mathbb{Z}_{4200},$$ since $8 \cdot 3 \cdot 25 \cdot 7=4200$, $4 \cdot 5=20$, and $2 \cdot 5=10$ (recall from Section \ref{multnth} the procedure of turning a primary decomposition of a positive integer into an invariant decomposition).

\vspace{.17in}

{\bf Exercises}

\begin{enumerate}

\item 

\begin{enumerate}

\item How many different isomorphism types are there among abelian groups of order 72?  

\item For each isomorphism class of part (a), find the member of the class that is in invariant decomposition and the one that is in primary decomposition. 

\end{enumerate}

\item 

\begin{enumerate} 

\item Prove that every abelian group of prime order is cyclic.

\item Characterize all positive integers $n$ for which it is true that every abelian group of order $n$ is cyclic.

\end{enumerate}

\end{enumerate}

\section{Subgroups and cosets}  \label{0.1.4}

As we have seen before, every element $a$ in a group $G$ determines a cycle
$$\langle a \rangle = \{ \lambda \cdot a  \mid \lambda =0,1,\dots,d-1 \};$$ the length $d$ of the cycle is the order $\mathrm{ord}_G(a)$ of $a$.  We also learned that $\mathrm{ord}_G(a)$ must divide the order $|G|$ of $G$.

A key feature of the cycle $\langle a \rangle$ is that it is not only a subset, but a {\em subgroup} of $G$; that is, $\langle a \rangle$ itself forms a group (for the same operation).  In general, if a subset $H \subseteq G$ is itself a group for the same operation, then we call $H$ a subgroup; this fact is denoted by $H \leq G$.  For example, $G=\mathbb{Z}_{10}$ has four different subgroups: besides $\{0\}$ and $G$ itself (which are always subgroups), we have $$\langle 2 \rangle =\{0,2,4,6,8\} \leq  \mathbb{Z}_{10}$$ and $$\langle 5 \rangle =\{0,5\} \leq  \mathbb{Z}_{10}.$$  The group $\mathbb{Z}_{12}$ has more: $\{0\}$, $\{0,6\}$, $\{0,4,8\}$, $\{0,3,6,9\}$, $\{0,2,4,6,8,10\}$, and $\mathbb{Z}_{12}$ itself.  It is not hard to prove that the cyclic group $\mathbb{Z}_n$ has exactly $d(n)$ different subgroups (where $d(n)$ is the number of positive divisors of $n$): each subgroup $H$ of $\mathbb{Z}_n$ must have order $d$ for some divisor $d$ of $n$, and, conversely, for each divisor $d$ of $n$, there is a unique subgroup $H$ of $\mathbb{Z}_n$ with order $d$.   Therefore, groups of prime order $p$---which, by the exercise at the end of the previous section, must be cyclic and thus are isomorphic to $\mathbb{Z}_p$---only have the trivial subgroups: $\{0\}$ and $\mathbb{Z}_p$ itself. 

In noncyclic groups, the situation is considerably more complicated.  Clearly, if $G=G_1 \times G_2$, then for all subgroups $H_1$ of $G_1$ and $H_2$ of $G_2$, $H_1 \times H_2$ is a subgroup of $G$; these kinds of subgroups are called {\em subproducts}.  For example, we may take $H_1=\{0,1\} \leq \mathbb{Z}_2$ and $H_2=\{0,2\} \leq \mathbb{Z}_4$, with which the subset
$$H_1 \times H_2=\{(0,0), (1,0), (0,2), (1,2)\}$$ is a subgroup of $\mathbb{Z}_2 \times \mathbb{Z}_4$.  But the collection of subgroups of a noncyclic group is more varied: $G_1 \times G_2$ may have subgroups that are not even in the form $H_1 \times H_2$.  In $G=\mathbb{Z}_2 \times \mathbb{Z}_4$, for example, $$H=\{(0,0), (1,2)\}$$ is a subgroup of $G$ and is not of this form.  It still holds that any subgroup of a (noncyclic) group of order $n$ must have order $d$ for some divisor $d$ of $n$, and, for any positive divisor $d$ of $n$, a group of order $n$ has a subgroup of order $d$, but we may have more than one such subgroup.

The number of subgroups of finite abelian groups is not yet fully understood.  As we mentioned above---and can easily be seen---the number of subgroups of the cyclic group $\mathbb{Z}_n$ is $d(n)$.  Regarding groups of rank two, the group $\mathbb{Z}_{n_1} \times \mathbb{Z}_{n_2}$ obviously has $d(n_1) \cdot d(n_2)$ subproducts, but its total number of subgroups is given by 
$$\sum_{d_1 \in D(n_1), d_2 \in D(n_2)} \gcd (d_1,d_2).$$  This formula follows from an 1987 paper of Calhoun\index{Calhoun, W. C.} (cf.~\cite{Cal:1987a}) and was recently proved directly by Hampejs et al.~(see \cite{Ham:2013a}).\index{Hampejs, M.}\index{Holighaus, N.}\index{T\'oth, L.}\index{Wiesmeyr, C.} Thus, for example, the group $\mathbb{Z}_{2} \times \mathbb{Z}_{4}$ has eight subgroups; of these, six are subproducts and two are not.  Similarly simple formulae for groups of rank three or more are not yet known.

Next, we discuss another important concept: cosets of subgroups.  Given any subgroup $H \leq G$ and any element $a \in G$, the set $$a+H=\{a+h  \mid h \in H\}$$ is called a {\em coset} of $H$ in $G$.  For example, if $G=\mathbb{Z}_{10}$, then for $H=\{0,5\} \leq G$ and $3 \in G$  we have 
$$3+H=\{3,8\},$$ and if $G=\mathbb{Z}_2 \times\mathbb{Z}_{4}$, then for $H=\{(0,0),(1,2)\}$ and $(1,1) \in G$ we have
$$(1,1)+H=\{(1,1),(0,3)\}.$$  Observe also that cosets could be represented by any of their other elements.   So, for example, with $H=\{0,5\}$, the coset $\{3,8\}$ of $\mathbb{Z}_{10}$ can be written both as $3+H$ and $8+H$, since in $\mathbb{Z}_{10}$ we have
$$3+\{0,5\}=8+\{0,5\};$$
similarly, with $H=\{(0,0),(1,2)\}$, the coset $\{(1,1), (0,3)\}$ of $\mathbb{Z}_2 \times\mathbb{Z}_{4}$ can be written both as $(1,1)+H$ and $(0,3)+H$ as in  $\mathbb{Z}_2 \times\mathbb{Z}_{4}$ we have
$$(1,1)+\{(0,0),(1,2)\}=(0,3)+\{(0,0),(1,2)\}.$$ 

Clearly, if $H$ has order $d$, then every coset of $H$ has size $d$ as well.  It turns out that for any two elements $a$ and $b$ of $G$, the cosets $a+H$ and $b+H$ are either identical (as sets) or entirely disjoint.  Therefore, the collection of distinct cosets of $H$ partitions $G$; that is, if $H$ has order $d$, then one can find $n/d$ elements $a_1, \dots, a_{n/d}$ so that $$G=(a_1 +H) \cup \cdots \cup (a_{n/d}+H).$$ For example, with $H=\{0,5\}$ in $G=\mathbb{Z}_{10}$, we have 
\begin{eqnarray*}
\mathbb{Z}_{10} & =& (0 +H) \cup (1 +H) \cup (2 +H) \cup  (3 +H) \cup (4+H) \\
& = & \{0,5\} \cup \{1,6\} \cup \{2,7\} \cup \{3,8\} \cup \{4,9\},
\end{eqnarray*}
and with $H=\{(0,0),(1,2)\}$ in $G=\mathbb{Z}_{2} \times\mathbb{Z}_{4} $, we have 
\begin{eqnarray*}
\mathbb{Z}_{2} \times\mathbb{Z}_{4} &=& ((0,0) +H) \cup ((0,1) +H) \cup ((1,0) +H)  \cup ((1,1) +H) \\
&=& \{(0,0),(1,2)\} \cup \{(0,1),(1,3)\} \cup \{(1,0),(0,2)\}  \cup \{(1,1),(0,3)\}.
\end{eqnarray*}  As we pointed out above, the elements $a_1, \dots, a_{n/d}$ are not unique; for example, we could also write $$\mathbb{Z}_{10}=(5 +\{0,5\})\cup(6 +\{0,5\})\cup(7 +\{0,5\})\cup(8 +\{0,5\})\cup(9+\{0,5\}).$$

Furthermore, we also define the {\em sum of cosets} $a+H$ and $b+H$ as $$(a+H)+(b+H)=(a+b)+H.$$ It is easy to see that this addition operation is well-defined (does not depend on which representatives $a$ and $b$ we choose), is closed (the sum of two cosets is also a coset), is associative, has an identity (the coset $0+H=H$), and each coset has an additive inverse (the inverse of $a+H$ being $-a+H$).  Thus, with this operation, the collection of cosets is itself a group, called the {\em quotient group} of $H$ in $G$; it is denoted by $G/H$.  

For example, the quotient group $\mathbb{Z}_{10}/\{0,5\}$ consists of the five cosets 
$$0+\{0,5\}, \; 1+\{0,5\}, \; 2+\{0,5\}, \; 3+\{0,5\}, \; 4+\{0,5\};$$ since (for example) $1+\{0,5\}$ generates all five cosets, we have $$\mathbb{Z}_{10}/\{0,5\} \cong \mathbb{Z}_5.$$  Similarly, 
the quotient group $(\mathbb{Z}_{2} \times \mathbb{Z}_{4}) /\{(0,0),(1,2)\}$ consists of the four cosets
$$(0,0)+\{(0,0),(1,2)\} , \; (0,1)+\{(0,0),(1,2)\}, \; (1,0)+\{(0,0),(1,2)\},  \; (1,1)+\{(0,0),(1,2)\};$$ here $(0,1)+\{(0,0),(1,2)\}$ cycles through the four cosets and thus $$(\mathbb{Z}_{2} \times \mathbb{Z}_{4}) /\{(0,0),(1,2)\} \cong \mathbb{Z}_{4}.$$

\vspace{.17in}

{\bf Exercises}

\begin{enumerate}

\item \begin{enumerate}

\item Find the number of subgroups of $\mathbb{Z}_{100}$.  How many of them are cyclic?

\item Find the number of subgroups of $\mathbb{Z}_{10}^2$.  How many of them are subproducts? 
\item Using information presented above, prove that the number of subgroups of $\mathbb{Z}_{n_1} \times \mathbb{Z}_{n_2}$ that are not subproducts equals $$\sum_{d_1 \in D(n_1), d_2 \in D(n_2)} (\gcd (d_1,d_2)-1).$$
\item Use part (c) to prove that every finite abelian group of rank two has at least one subgroup that is not a subproduct.
\item Prove that $\mathbb{Z}_2^2$ is the only finite abelian group of rank two that has exactly one subgroup that is not a subproduct. 

\end{enumerate}

\item \begin{enumerate}

\item List all subgroups of $\mathbb{Z}_{18}$.  Find the isomorphism type of each subgroup. 

\item List all subgroups of $\mathbb{Z}_{3} \times \mathbb{Z}_{6}$.  Find the isomorphism type of each subgroup.

\end{enumerate}

\item \begin{enumerate}

\item Consider the subgroup $H=\{0,9\}$ of $G=\mathbb{Z}_{18}$.  Find each coset of $H$ in $G$, and find the isomorphism type of $G/H$. 

\item Consider the subgroup $H=\{(0,0),(0,3)\}$ of $G=\mathbb{Z}_{3} \times \mathbb{Z}_{6}$.  Find each coset of $H$ in $G$, and find the isomorphism type of $G/H$.

\end{enumerate}

\end{enumerate}

\section{Subgroups generated by subsets} \label{Subgroups generated by subsets}

Recall that, for any element $a$ of $G$, the set $\langle a \rangle$ consists of all $d$ multiples of $a$ where $d=\mathrm{ord}_G(a)$ is the order of $a$:
$$\langle a \rangle=\{\lambda \cdot a  \mid \lambda=0,1,\dots,d-1\}.$$ We should note that it makes no difference if we increase the range of $\lambda$: since $d \cdot a=0$, $(d+1) \cdot a = a$, $(d+2) \cdot a = 2 \cdot a$, etc., we may also write
$$\langle a \rangle=\{\lambda \cdot a  \mid \lambda \in \mathbb{N}_0\};$$ in fact, since $(-1) \cdot a=(d-1) \cdot a$, $(-2) \cdot a=(d-2) \cdot a$, and so on, we have
$$\langle a \rangle=\{\lambda \cdot a  \mid \lambda \in \mathbb{Z}\}.$$

Recall also that $\langle a \rangle$ is a subgroup of $G$; we call $\langle a \rangle$ the subgroup generated by $a$ as it is the smallest subgroup of $G$ that contains $a$.  More generally, we may look for the smallest subgroup of $G$ that contains all of the set $A=\{a_1,\dots,a_m\}$; this subgroup, denoted by $\langle A \rangle$, is called the {\em subgroup generated by $A$}.  It is not hard to see that, with $d_i= \mathrm{ord}(a_i)$ ($i=1,2,\dots,m$), we have
$$\langle A \rangle = \{ \lambda_1 \cdot a_1 + \cdots + \lambda_m \cdot a_m  \mid \lambda_i=0,1,\dots,d_i-1 \; \mbox{for} \; i=1,2,\dots,m\}.$$  For example, in $G=\mathbb{Z}_{15}$, the subset $A=\{6,10\}$ generates the subgroup
$$\langle A \rangle = \{ \lambda_1 \cdot 6 +\lambda_2 \cdot 10  \mid \lambda_1=0,1,2,3,4; \lambda_2=0,1,2\},$$ since 6 has order 5 and 10 has order 3 in $G$.  Evaluating the expressions above yields
$$\langle A \rangle = \{0,6,12,3,9, 10,1,7,13,4, 5,11,2,8,14\};$$ so in this case we have $\langle A \rangle=G$.
Again, it makes no difference if we let the range of coefficients increase:   
$$\langle A \rangle = \{ \lambda_1 \cdot a_1 + \cdots + \lambda_m \cdot a_m  \mid \lambda_i \in \mathbb{N}_0 \; \mbox{for} \; i=1,2,\dots,m\}$$ or
$$\langle A \rangle = \{ \lambda_1 \cdot a_1 + \cdots + \lambda_m \cdot a_m  \mid \lambda_i \in \mathbb{Z} \; \mbox{for} \; i=1,2,\dots,m\}.$$

In our example above, we have the 2-element set $A=\{6,10\}$ generating the cyclic group $G=\mathbb{Z}_{15}$; but, like any cyclic group, $G=\mathbb{Z}_{15}$ could be generated by a single element ($1$ or any other $a \in \mathbb{Z}_{15}$ that is relatively prime to 15).  The question may arise: how large does a subset of a group need to be in order to generate the entire group?  The general answer is provided by the following proposition. 

\begin{prop} \label{rankvsgen}
An abelian group of rank $r$ cannot be generated by fewer than $r$ elements; that is, if $A=\{a_1,\dots,a_m\}$ is a subset of $G$ for which $\langle A \rangle = G$, then $m \geq r$.
\end{prop}

For the proof of Proposition \ref{rankvsgen}, see page \pageref{proofofrankvsgen}.  We should note that a subset of size $r$ doesn't always generate a given group of rank $r$.  For example, $\{(1,2),(1,4)\}$ in the group $\mathbb{Z}_6^2$ generates only a subgroup of order 18 (only ordered pairs with an even second component get generated), and the subset $\{(1,2),(3,1)\}$ of $\mathbb{Z}^2_5$ only generates a subgroup of order 5 (each of the two elements alone, in fact, generates the other).  It is not easy to see in general how large of a subgroup a given subset $A$ of a group $G$ generates.

\vspace{.17in}

{\bf Exercises}

\begin{enumerate}

\item \begin{enumerate}

\item Above, we made the comment that for any $a \in \mathbb{Z}_{15}$ that is relatively prime to 15, $a$ generates the entire group.  Verify this statement.

\item Prove that, conversely, if an element $a \in \mathbb{Z}_{15}$ is not relatively prime to 15, then $a$ does not generate the entire group.

\end{enumerate}

\item  Consider the group $G=\mathbb{Z}_{15}^2$.  Find the subgroups of $G$ generated by each of the following subsets, and find the isomorphism type of each subgroup.
 
\begin{enumerate}

\item $\{(1,3),(1,5)\}$.

\item $\{(1,3),(2,3)\}$.

\item $\{(1,4),(4,1)\}$.

\end{enumerate}

\end{enumerate}

\section{Sumsets} \label{sumset}

As we have just seen, taking all integer {\em linear combinations}
$$\lambda_1 \cdot a_1 + \cdots + \lambda_m \cdot a_m$$
of a subset $A=\{a_1,\dots,a_m\}$ of some finite abelian group $G$ generates a subgroup  $\langle A \rangle$ of $G$; in fact, we arrive at $\langle A \rangle$ by taking only the above linear combinations with coefficients $$\lambda_i =0,1,\dots, \mathrm{ord}(a_i)-1$$ ($i = 1,2 \dots,m$). 

In many famous and still-investigated problems, we put limitations on the coefficients.  For example, while we saw in Section \ref{Subgroups generated by subsets} that $\{6,10\}$ generates all of $\mathbb{Z}_{15}$, to generate 7, we may use (for example) 6 twice and 10 once (thus, a total of three terms) to write
$$2 \cdot 6 + 1 \cdot 10,$$ but there is no way to get to 7 without using at least three terms.  The question then arises: which terms do we get with, say, at most two terms?  In other words, we want to find the linear combinations $$\lambda_1 \cdot 6 + \lambda_2 \cdot 10$$ with integer coefficients $\lambda_1$ and $\lambda_2$ satisfying 
$$|\lambda_1|+|\lambda_2| \leq 2.$$ 
The answer then is the set 
$$ \{0, \; \pm 1 \cdot 6, \; \pm 1 \cdot 10, \; \pm 2 \cdot 6, \; \pm 2 \cdot 10, \; \pm 1 \cdot 6 \pm 1 \cdot 10 \} = \{0, 1, 3, 4, 5, 6, 9, 10, 11, 12, 14\}.$$
These eleven elements correspond to the thirteen elements of the layer 
$$\mathbb{Z}^2(2)=\{ (\lambda_1,\lambda_2) \in \mathbb{Z}^2  \mid |\lambda_1|+|\lambda_2 | \leq  2 \}=\{(0,0),(\pm 1,0),(0, \pm 1),(\pm 2,0), (0, \pm 2), (\pm 1, \pm 1)\}$$ of the 2-dimensional integer lattice (see Section \ref{0.2.4}).  (Note that not all elements of $\mathbb{Z}^2(2)$ yield distinct group elements: $1 \cdot 10=(-2) \cdot 10$ and $2 \cdot 10= (-1) \cdot 10$.) 

More generally, given a subset $\Lambda \subseteq \mathbb{Z}$ and a subset $H \subseteq \mathbb{N}_0$, we consider {\em sumsets} of $A=\{a_1,\dots,a_m\}$ corresponding to
$$\Lambda^m (H) = \{ (\lambda_1,\lambda_2,\dots,\lambda_m) \in \Lambda^m  \mid |\lambda_1|+|\lambda_2|+\cdots+|\lambda_m| \in H \}$$
(see Section \ref{0.2.4}); namely, we consider the collection 
\begin{eqnarray*}
H_{\Lambda}A &=& \{\lambda_1a_1+\cdots +\lambda_m a_m \; |  \;  (\lambda_1,\dots, \lambda_m) \in \Lambda^m (H) \} \\ \\
&=& \{\lambda_1a_1+\cdots +\lambda_m a_m \; |  \;  \lambda_1,\dots, \lambda_m \in \Lambda, |\lambda_1|+\cdots+|\lambda_m| \in H  \}
\end{eqnarray*} in $G$.  In our example above we had $G=\mathbb{Z}_{15}$, $A=\{6,10\}$, $\Lambda = \mathbb{Z}$, and $H=[0,2]=\{0,1,2\}$.     

We see that the sumset of $A$ corresponding to $\Lambda $ and $H$ consists of sums of the elements of $A$ with the conditions that
\begin{itemize}

\item the repetition number of each element of $A$ in the sum must be from the set $\Lambda$; and

\item the total number of terms in each sum must be from the set $H$.

\end{itemize}
We need to point out that the order of the elements of $A$ is immaterial: the sumset remains unchanged when the elements $a_1,\dots,a_m$ are permuted.

There are a number of special cases of sumsets that are most often studied.  First, we introduce terms and notations for the case when sums in our sumset all have a fixed number of terms (that is, $H$ consists of a single nonnegative integer $h$) and when the coefficient-set is $\Lambda=\mathbb{N}_0, \mathbb{Z}, \{0,1\}$, or $\{-1,0,1\}$.  
More precisely, we introduce the following notations.  Suppose that $A=\{a_1,a_2,\dots,a_m\}$ is a subset of an abelian group $G$ (with $m \in \mathbb{N}$) and that $h$ is a nonnegative integer.  We then define:   
\begin{itemize}
\item the (ordinary) {\em $h$-fold sumset} $hA$ of $A$, consisting of sums of exactly $h$ (not necessarily distinct) terms of $A$:
$$hA = \{ \lambda_1a_1+\cdots +\lambda_m a_m \mbox{   } | \mbox{   }  \lambda_1, \dots, \lambda_m \in \mathbb{N}_0, \; \lambda_1+\cdots +\lambda_m =h \};$$
\item the {\em $h$-fold signed sumset} $h_{\pm}A$ of $A$, consisting of signed sums of exactly $h$ (not necessarily distinct) terms of $A$:
$$h_{\pm}A = \{ \lambda_1a_1+\cdots +\lambda_m a_m \mbox{   } | \mbox{   }  \lambda_1, \dots, \lambda_m \in \mathbb{Z}, \; |\lambda_1|+\cdots +|\lambda_m| =h \};$$
\item the {\em restricted $h$-fold sumset} $h \hat{\;} A$ of $A$, consisting of sums of exactly $h$ distinct terms of $A$:
$$h\hat{\;}A = \{ \lambda_1a_1+\cdots +\lambda_m a_m \mbox{   } | \mbox{   }  \lambda_1, \dots, \lambda_m \in \{0,1\}, \; \lambda_1+\cdots +\lambda_m =h \};$$ and
\item the {\em restricted $h$-fold signed sumset} $h \hat{_{\pm}}A$ of $A$, consisting of signed sums of exactly $h$ distinct terms of $A$:
$$h \hat{_{\pm}}A = \{ \lambda_1a_1+\cdots +\lambda_m a_m \mbox{   } | \mbox{   }  \lambda_1, \dots, \lambda_m \in \{-1,0,1\}, \; |\lambda_1|+\cdots +|\lambda_m| =h \}.$$
\end{itemize}   

A summary scheme of these four types of sumsets, indicating containment, can be given as follows.

\begin{tabular}{l|ccc}
& Terms must be & & Repetition of terms \\ 
& distinct & & is allowed \\ \hline
Terms can be & & & \\
added only & \; $h\hat{\;}A$ & $\subseteq$ & \; $hA$ \\ \\ 
& $| \bigcap$ & & $| \bigcap$  \\
Terms can be & & & \\ 
added or subtracted & \; $h \hat{_{\pm}}A$  & $\subseteq$ & \; $h_{\pm}A$ \\ \\
\end{tabular}

We should caution that the $h$-fold sumset $hA$ of a set $A$ is different from the so called {\em $h$-fold dilation} $$h \cdot A = \{h \cdot a_1, \dots, h \cdot a_m \}$$ of $A$, which only contains the multiples of the elements in $A$.  (However, for an element $g \in G$ and an integer $\lambda \in \mathbb{Z}$ we will use the notations $\lambda \cdot g$ and $\lambda g$ interchangeably.)  We clearly have $h \cdot A \subseteq hA \,$ but usually $hA$ is much larger than $h \cdot A$.  We will rarely use dilations in this book, except for the negative $-A$ of $A$, defined, of course, as $$-A=(-1) \cdot A=\{-a_1,\dots,-a_m\}.$$

Let us now illustrate the spectra of our four types of sumsets as $h$ increases, using the example of $G=\mathbb{Z}_{13}$ and $A=\{2,3\}$.   

The easiest to evaluate are the restricted sumsets; note that since $|A|=2$, for $h \geq 3$ we have $h \hat{\;} A=  h\hat{_{\pm}}A =\emptyset$.  The relevant values of 
$$\lambda_1 \cdot 2 + \lambda_3 \cdot 3$$
in $G$ are as follows:
$$\begin{array}{||l||c|c||} \hline \hline
\lambda_2=1   & 3 & 5     \\ \hline 
\lambda_2=0   &  0 & 2   \\ \hline \hline
  & \lambda_1=0 &  \lambda_1=1  \\ \hline \hline 
\end{array}$$ 
Therefore, we get
  $$\begin{array}{||c||r|r|r|r|r|r|r|r|r|r|r|r|r||} \hline \hline
0 \hat{\;} A &0&      &  &  &   &   &   &   &   &   &    &    &     \\ \hline
1 \hat{\;} A &&     & 2 & 3 &   &   &   &   &   &   &    &    &     \\ \hline 
2 \hat{\;} A &&     &   &   & & 5 & &   &   &   &    &    &     \\ \hline \hline
\end{array}$$

Similarly, for the restricted signed sumsets we have

$$\begin{array}{||l||c|c|c||} \hline \hline
\lambda_2=1  &1 & 3 & 5     \\ \hline 
\lambda_2=0   &11 & 0 & 2   \\ \hline 
\lambda_2=-1   &  8 &  10 & 12   \\ \hline \hline
 & \lambda_1=-1 & \lambda_1=0 &  \lambda_1=1  \\ \hline \hline 
\end{array}$$ 
and thus

$$\begin{array}{||c||r|r|r|r|r|r|r|r|r|r|r|r|r||} \hline \hline
 0\hat{_{\pm}}A & 0     &  &  &   &   &   &   &   &   &    &     &   &  \\ \hline 
 1\hat{_{\pm}}A &&      & 2 & 3 &   &   &   &   &   &   &  10  &    11 &     \\ \hline 
2\hat{_{\pm}}A && 1     &   &   & & 5 & &   &   8 &   &    &    &   12  \\ \hline \hline
\end{array}$$

Turning to unrestricted sumsets, we first observe that both  $hA$ and $h_{\pm}A$ equal the entire group $\mathbb{Z}_{13}$ for $h \geq 12$: indeed, we have
$$h_{\pm}A \supseteq hA=\{(h-i) \cdot 2+i \cdot 3  \mid i=0,1,\dots,h\} = \{2h+i  \mid i=0,1,\dots,h\},$$ which, if $h \geq 12$, gives the entire group $\mathbb{Z}_{13}$.  So we only need to exhibit the $h$-fold sumset and the $h$-fold signed sumset of $A$ for $h \leq 12$; a computation similar to the ones above yields the following results:

  $$\begin{array}{||c||r|r|r|r|r|r|r|r|r|r|r|r|r||} \hline \hline
0A & 0  &   &  &  &   &   &   &   &   &   &    &    &     \\ \hline
1A &   &   & 2 & 3 &   &   &   &   &   &   &    &    &     \\ \hline 
2A &   &   &   &   & 4 & 5 & 6 &   &   &   &    &    &     \\ \hline 
3A &   &   &   &   &   &   & 6 & 7 & 8 & 9 &    &    &     \\ \hline 
4A &   &   &   &   &   &   &   &   & 8 & 9 & 10 & 11 & 12  \\ \hline 
5A & 0 & 1 & 2 &   &   &   &   &   &   &   & 10 & 11 & 12  \\ \hline 
6A & 0 & 1 & 2 & 3 & 4 & 5 &   &   &   &   &    &    & 12  \\ \hline 
7A &   & 1 & 2 & 3 & 4 & 5 & 6 & 7 & 8 &   &    &    &     \\ \hline 
8A &   &   &   & 3 & 4 & 5 & 6 & 7 & 8 & 9 & 10 & 11 &     \\ \hline 
9A & 0 & 1 &   &   &   & 5 & 6 & 7 & 8 & 9 & 10 & 11 & 12  \\ \hline 
10A & 0 & 1 & 2 & 3 & 4 &   &   & 7 & 8 & 9 & 10 & 11 & 12 \\ \hline 
11A & 0 & 1 & 2 & 3 & 4 & 5 & 6 & 7 &   & 9 & 10 & 11 & 12 \\ \hline 
12A & 0 & 1 & 2 & 3 & 4 & 5 & 6 & 7 & 8 & 9 & 10 & 11 & 12 \\ \hline \hline

\end{array}$$

$$\begin{array}{||c||r|r|r|r|r|r|r|r|r|r|r|r|r||} \hline \hline
0_{\pm}A & 0  &   &  &  &   &   &   &   &   &   &  & &     \\ \hline
1_{\pm}A &   &   & 2 & 3 &   &   &   &   &   &   & 10 & 11 &     \\ \hline 
2_{\pm}A &   & 1 &   &   & 4 & 5 & 6 & 7 & 8 & 9 &    &    & 12  \\ \hline 
3_{\pm}A &   & 1 &   &   & 4 & 5 & 6 & 7 & 8 & 9 &    &    & 12  \\ \hline 
4_{\pm}A &   & 1 & 2 & 3 & 4 & 5 & 6 & 7 & 8 & 9 & 10 & 11 & 12  \\ \hline 
5_{\pm}A & 0 & 1 & 2 & 3 &   & 5 &   &   & 8 &   & 10 & 11 & 12  \\ \hline 
6_{\pm}A & 0 & 1 & 2 & 3 & 4 & 5 & 6 & 7 & 8 & 9 & 10 & 11 & 12  \\ \hline 
7_{\pm}A &   & 1 & 2 & 3 & 4 & 5 & 6 & 7 & 8 & 9 & 10 & 11 & 12  \\ \hline 
8_{\pm}A &   & 1 & 2 & 3 & 4 & 5 & 6 & 7 & 8 & 9 & 10 & 11 & 12  \\ \hline 
9_{\pm}A & 0 & 1 & 2 & 3 & 4 & 5 & 6 & 7 & 8 & 9 & 10 & 11 & 12  \\ \hline 
10_{\pm}A & 0 & 1 & 2 & 3 & 4 & 5 & 6 & 7 & 8 & 9 & 10 & 11 & 12 \\ \hline 
11_{\pm}A & 0 & 1 & 2 & 3 & 4 & 5 & 6 & 7 & 8 & 9 & 10 & 11 & 12 \\ \hline 
12_{\pm}A & 0 & 1 & 2 & 3 & 4 & 5 & 6 & 7 & 8 & 9 & 10 & 11 & 12 \\ \hline \hline
\end{array}$$

We can introduce similar notations and terminology for the cases when sums in our sumsets contain a limited number of terms (say $H=[0,s]=\{0,1,\dots,s\}$ for some $s \in \mathbb{N}_0$) or an arbitrary number of terms ($H=\mathbb{N}_0$).  For instance, with the set $A=\{2,3\}$ in $G=\mathbb{Z}_{13}$, we get
$$[0,3]A=\cup_{i=0}^3 hA = \{0\} \cup \{2,3\}, \{4,5,6\} \cup \{6,7,8,9\}=\{0,2,3,4,5,6,7,8,9\},$$
and
$$[0,3]_{\pm}A=\cup_{i=0}^3 h_{\pm} A = \{0\} \cup \{2,3,10,11\}, \{1,4,5,6,7,8,9,12\} \cup \{1,4,5,6,7,8,9,12\}=G.$$

Our notations and terminology are summarized in the following table.

$$\begin{array}{l||c|c|c|} \label{twelvesumsets}
  & H=\{h\}  & H=[0,s]  & H=\mathbb{N}_0 \\ \hline \hline

 & h A & [0,s] A & \langle A \rangle \\ 
\Lambda=\mathbb{N}_0  & & &  \\
   & \mbox{$h$-fold sumset} & \mbox{$[0,s]$-fold sumset} & \mbox{sumset} \\ \hline

  & h_{\pm}A & [0,s]_{\pm}A & \langle A \rangle \\
\Lambda=\mathbb{Z} & && \\ 
   & \mbox{$h$-fold signed sumset} & \mbox{$[0,s]$-fold signed sumset} & \mbox{signed sumset} \\ \hline

  & h\hat{\;}A & [0,s] \hat{\;} A & \Sigma A \\
\Lambda=\{0,1\}& & & \\ 
   & \mbox{restricted} & \mbox{restricted} & \mbox{restricted} \\
& \mbox{$h$-fold sumset} & \mbox{$[0,s]$-fold sumset} & \mbox{sumset} \\ \hline

  & h\hat{_{\pm}}A & [0,s]\hat{_{\pm}}A & \Sigma_{\pm}A \\
\Lambda=\{-1,0,1\} & & & \\ 
   & \mbox{restricted} & \mbox{restricted} & \mbox{restricted} \\
   & \mbox{$h$-fold signed sumset} & \mbox{$[0,s]$-fold signed sumset} & \mbox{signed sumset} \\ \hline

\end{array}$$

Recall that the sumset $\cup_{h=0}^{\infty} hA$ and the signed sumset $\cup_{h=0}^{\infty} h_{\pm}A$ of a subset $A$ both equal $\langle A \rangle$, the subgroup of $G$ generated by $A$.  For example, with our previous example of $A=\{2,3\}$ in $G=\mathbb{Z}_{13}$, we have $\langle A \rangle = \mathbb{Z}_{13}.$  (Since 13 is prime, $\mathbb{Z}_{13}$ has no subgroups other than $\{0\}$ and $\mathbb{Z}_{13}$.)

Next, we point out some obvious but useful identities.  

Suppose that $A \subseteq G$ and assume, as usual, that $A=\{a_1,\dots,a_m\}$.  Let $g$ be an arbitrary element of $G$.  We then define the set $A-g$ as $\{a_1-g,\dots,a_m-g\}$.  For a fixed $h \in \mathbb{N}_0$, we now examine the sumsets $h(A-g)$ and $h \hat{\;} (A-g)$.  By definition, we have
\begin{eqnarray*}
h(A-g) & = & h\{a_1-g,\dots,a_m-g\} \\
& = & \{\lambda_1(a_1-g)+ \cdots + \lambda_m (a_m-g)  \mid \lambda_1,\dots,\lambda_m \in \mathbb{N}_0, \lambda_1+\cdots+\lambda_m=h\} \\
& = & \{\lambda_1a_1+ \cdots + \lambda_m a_m-h \cdot g  \mid \lambda_1,\dots,\lambda_m \in \mathbb{N}_0, \lambda_1+\cdots+\lambda_m=h\} \\
& = & hA-h\cdot g. 
\end{eqnarray*}
Similarly,
\begin{eqnarray*}
h \hat{\;}  (A-g) & = & h \hat{\;} \{a_1-g,\dots,a_m-g\} \\
& = & \{\lambda_1(a_1-g)+ \cdots + \lambda_m (a_m-g)  \mid \lambda_1,\dots,\lambda_m \in \{0,1\}, \lambda_1+\cdots+\lambda_m=h\} \\
& = & \{\lambda_1a_1+ \cdots + \lambda_m a_m-h \cdot g  \mid \lambda_1,\dots,\lambda_m \in \{0,1\}, \lambda_1+\cdots+\lambda_m=h\} \\
& = & h \hat{\;} A-h\cdot g. 
\end{eqnarray*}
Therefore, we see that $h(A-g)$ has the same cardinality as $hA$, and $h \hat{\;} (A-g)$ has the same cardinality as $  h \hat{\;}A$.  
This is particularly useful when $g \in A$, since then $A-g$ is a subset of $G$ that contains zero.  
In summary, we have the following:
\begin{prop}  \label{shifting sets}
For any $G$, $A \subseteq G$, $g \in G$, and $h \in \mathbb{N}_0$, we have 
$$|h(A-g)|=|hA|$$ and $$|h \hat{\;} (A-g)|=|h \hat{\;} A|.$$
In particular, there is a subset $A_0$ of $G$ so that $0 \in A_0$, $|A_0|=|A|$, $|hA_0|=|hA|$, and $|h \hat{\;} A_0|=|h \hat{\;} A|$.
\end{prop}
We need to point out that the ``signed sumset'' versions of these identities (when subtraction of elements is allowed) do not necessarily hold.

The following identities are also obvious: 
\begin{prop}  \label{[0,s]A=}
For any $G$, $A \subseteq G$, and $s \in \mathbb{N}_0$, we have $$[0,s]A=s(A \cup \{0\})$$
and
$$[0,s]_{\pm}A=s_{\pm}(A \cup \{0\}).$$ 
\end{prop}
We should note that the restricted versions of these identities do not necessarily hold. 

We also find that the ($h$-fold, etc.) signed sumset of a subset $A$ is closely related to the ($h$-fold, etc.) sumset of $A \cup (-A)$.  Of course, for $h=0$, we have $$0_{\pm}A=0(A\cup (-A))=\{0\};$$ for $h=1$ we get
$$1_{\pm}A=1(A \cup (-A))=A \cup (-A).$$

For $h=2$, we can easily see that $2_{\pm}A$ and $2(A \cup (-A))$ are almost the same, except that $2(A \cup (-A))$ always contains 0 while $2_{\pm}A$ may not, so
$$2(A \cup (-A))=2_{\pm}A \cup \{0\}=2_{\pm}A \cup 0_{\pm}A.$$

More generally, we have the following identities.

\begin{prop} \label{sumsetidentities}
For every subset $A$ of $G$ and every nonnegative integer $h$ we have
$$h(A \cup (-A)) = h_{\pm}A \cup (h-2)_{\pm}A \cup (h-4)_{\pm}A \cup \cdots;$$ and therefore, for all $s \in \mathbb{N}$,
$$[0,s](A \cup (-A)) = [0,s]_{\pm}A $$ and
$$\langle A \cup (-A )\rangle = \langle A \rangle.$$
\end{prop}
We provide the proof on page \pageref{proofofsumsetidentities}.

Note that when $A$ is {\em symmetric}, that is, $A=-A$, then an even simpler relation holds:

\begin{prop}  \label{hA and h pm A with sym}
For every symmetric subset $A$ of $G$ and every nonnegative integer $h$ we have $h_{\pm}A=hA$.

\end{prop}
The easy proof can be found on page \pageref{proof of hA and h pm A with sym}.

Regarding restricted sumsets, we always have
$$h\hat{\;}(A \cup (-A)) \subseteq h\hat{_{\pm}}A \cup (h-2)\hat{_{\pm}}A \cup (h-4)\hat{_{\pm}}A \cup \cdots,$$ (see the last paragraph of the proof of Proposition \ref{sumsetidentities}), but an identity for our two types of restricted sumsets seems a lot more complicated.  Note that, for example, when $A$ contains distinct elements $a_1$ and $a_2$ for which $a_1=-a_2$, then $a_1-a_2$ is definitely an element of $2\hat{_{\pm}}A$, but not necessarily of $2\hat{\;}(A \cup (-A))$.  We pose the following vague problem.

\begin{prob}
Find identities similar to those in Proposition \ref{sumsetidentities} for restricted sumsets.

\end{prob}

\vspace{.17in}

{\bf Exercises}

\begin{enumerate}

\item Consider the subset $A=\{2,3\}$ in the cyclic group $G=\mathbb{Z}_{21}$.  List the elements of each of the following sumsets. 

\begin{enumerate}

\item $[0,3]_{\pm} A$

\item $[0,3] A$

\item $[0,3]\hat{_{\pm}} A$

\item $[0,3] \hat{\;} A$

\end{enumerate}

\item Consider the subset $A=\{(0,0,1),(0,1,1),(1,1,1)\}$ in the group $G=\mathbb{Z}_{3}^3$.  List the elements of each of the following sumsets. 

\begin{enumerate}

\item $2_{\pm} A$

\item $2A$

\item $2\hat{_{\pm}} A$

\item $2 \hat{\;} A$

\end{enumerate}

\end{enumerate}

\part{Appetizers}  \label{Appetizers}

The short articles in this part are meant to invite everyone to the main entrees of the menu, presented in the next part of the book. Our appetizers are carefully chosen so that they provide bite-size representative samples of the research projects, as well as make connections to other parts of mathematics that students might have encountered.  (The first article describes how this author's research in spherical geometry lead him to additive combinatorics.)

\newpage

\addcontentsline{toc}{section}{Spherical designs}  

\section*{Spherical designs}

This appetizer tells the (simplified) story of how the author's research in algebraic combinatorics and approximation theory lead him to additive combinatorics (cf.~\cite{Baj:1998a}, \cite{Baj:2000a}). 
 
Imagine that we want to scatter a certain number of points on a sphere---how can we do it in the most ``uniformly balanced'' way?  The answer might be obvious in certain cases, but far less clear in others.  Consider, for example, the sphere $S^2$ in 3-dimensional Euclidean space, and suppose that we need to place six points on it: the most balanced configuration for the points is undoubtedly at the six vertices of a regular octahedron (for example, at the points $(\pm 1, 0, 0), (0, \pm 1, 0), (0,0,\pm 1)$ in case the sphere is centered at the origin and has radius 1).  The answer is perhaps equally clear for four, eight, twelve, or twenty points: place them so that they form a regular tetrahedron, cube, icosahedron, and dodecahedron, respectively.  But how should we position five points?  Or how about seven or ten or a hundred points?   How can one do this in general?

The answer to our question depends, of course, on how we measure the degree to which our pointset is balanced.  For example, in the case of a {\em packing problem}, we want to place our points as far away from each other as possible; more precisely, we may want to maximize the minimum distance between any two of our points.  Or, when addressing a {\em covering problem}, we want to minimize the maximum distance that any place on the sphere has from the closest point of our pointset.    There are several other reasonable criteria, but the one that is perhaps the most well known and applicable is the one where we maximize the degree to which our pointset is in {\em momentum balance}.  We make this precise as follows.

As usual, we let $S^d$ denote the sphere in the $(d+1)$-dimensional Euclidean space $\mathbb{R}^{d+1}$; we also assume that $S^d$ is centered at the origin and has radius 1, in which case it can be described as the set of points that have distance 1 from the origin.
Given a finite pointset $P \subset S^d$ and a polynomial $f: S^d \rightarrow \mathbb{R}$, we define the average of $f$ over $P$ as 
$$\overline{f}_P = \frac{1}{|P|} \cdot \sum_{p \in P} f(p);$$ the average of $f$ over the entire sphere is given by
$$\overline{f}_{S^d} = \frac{1}{|S^d|} \cdot \int_{S^d} f(x) \mathrm{d}x$$  (here $|S^d|$ denotes the surface area of  $S^d$).  We then say that our finite pointset $P$ is a {\em spherical $t$-design on $S^d$}, if $$\overline{f}_P = \overline{f}_{S^d}$$ holds for all polynomials $f$ of degree up to $t$.

Spherical designs were introduced and first studied in 1977 by Delsarte, Goethals, and Seidel in \cite{DelGoeSei:1977a}.\index{Delsarte, P.}\index{Goethals, J. M.}\index{Seidel, J. J.} In this far-reaching paper they also established the following tight lower bound for the number of points needed to form a spherical $t$-design on $S^d$: 
$$
|P| \geq N_t^d={d+\lfloor t/2 \rfloor \choose d}+{d+\lfloor (t-1)/2 \rfloor  \choose d}. 
$$
We shall refer to the quantity $N_t^d$ as the {\em DGS bound}.  Spherical designs of this minimum size are called {\em tight}.  Bannai and Damerell\index{Bannai, Ei.}\index{Damerell, R. M.} (cf.~\cite{BanDam:1979a}, \cite{BanDam:1980a}) proved that tight spherical designs for $d \geq 2$ exist only for $t=1,2,3,4,5,7$, or $11$. All tight $t$-designs are known, except possibly for $t=4,5,$ or 7; in particular, there is a unique $11$-design (with $d$=23 and $|P|=196560$, coming from the {\em Leech lattice}).

Clearly, any finite pointset $P \subset S^d$ is a spherical 0-design.  Spherical 1-designs and 2-designs have basic interpretations in physics: it is easy to check that $P$ is a spherical 1-design exactly when $P$ is in mass balance (that is, its center of gravity is at the center of the sphere), and we can verify that $P$ is a spherical 2-design if, and only if, $P$ is in both mass balance and inertia balance.  Higher degrees of balance find their applications in a variety of areas, including crystallography, coding theory, astronomy, and viral morphology.  

According to the definition, to test whether a given pointset $P$ is a spherical $t$-design, one should check whether $\overline{f}_{P}$ agrees with $\overline{f}_{S^d}$ for every polynomial $f$ of degree at most $t$.  There is a well-known shortcut: it suffices to do this for non-constant homogeneous harmonic polynomials of degree at most $t$.  A polynomial is called {\em homogeneous} if all its terms have the same (total) degree; for example, $$f=x^3z+xz^3-6xy^2z$$ is a homogeneous polynomial (of three variables and of degree 4).  A polynomial is called {\em harmonic} if it satisfies the {\em Laplace equation}; that is, the sum of its unmixed second-order partial derivatives equals zero.  We see that the polynomial $f$ just mentioned is harmonic:
$$\frac{\partial^2 f}{\partial x^2} +  \frac{\partial^2 f}{\partial y^2} + \frac{\partial^2 f}{\partial z^2} = 6xz -12xz +6xz = 0.$$

The reduction to homogeneous harmonic polynomials has two advantages.  First, there are a lot fewer of them:  the set $\mathrm{Harm}_k(S^d)$ of homogeneous harmonic polynomials over $S^d$ of degree $k$ forms a vector space over $\mathbb{R}$ whose dimension is `only' \label{dim harm}
$$\dim \mathrm{Harm}_k(S^d) = {d+k \choose d} - {d+k-2 \choose d}.$$  Second, the average of non-constant harmonic polynomials over the sphere is zero.  Therefore, a pointset $P \subset S^d$ is a spherical $t$-design on $S^d$ if, and only if, $\overline{f}_P =0$ for every $f \in \mathrm{Harm}_k(S^d)$ and $1 \leq k \leq t$. 

So, how can we construct spherical $t$-designs?  We get an inspiration from the case $d=1$ (when our sphere is a circle).  Note that 
$$\dim \mathrm{Harm}_k(S^1) = {k+1 \choose 1} - {k-1 \choose 1}=2,$$ and we can verify that the real and imaginary parts of the polynomial $(x+iy)^k$ (with $i$ denoting the `imaginary' square root of $-1$) form a basis of $\mathrm{Harm}_k(S^1)$; for example, 
\begin{eqnarray*}
\mathrm{Harm}_1(S^1) & = &  \langle \mathrm{Re} (x+iy)^1, \mathrm{Im} (x+iy)^1 \rangle = \langle x, y \rangle, \\ \\ 
\mathrm{Harm}_2(S^1) & = &  \langle \mathrm{Re} (x+iy)^2, \mathrm{Im} (x+iy)^2 \rangle = \langle x^2-y^2, 2xy \rangle, \\ \\
\mathrm{Harm}_3(S^1) & = &  \langle \mathrm{Re} (x+iy)^3, \mathrm{Im} (x+iy)^3 \rangle = \langle x^3-3xy^2, 3x^2y-y^3 \rangle.
\end{eqnarray*}

As it is probably easy to guess and as we now verify: regular polygons yield spherical designs on the circle.  More precisely, we show that the set of vertices 
$$P=\left\{ \left( \cos \left( 2 \pi j/n \right) ,  \sin \left( 2 \pi j/n \right) \right) \; \mid \; j=0,1,2,\dots, n-1 \right\}$$
of the regular $n$-gon form a spherical $t$-design on $S^1$ for all $n \geq t+1$.  (Note that, by the DGS bound, a spherical $t$-design on $S^1$ must have at least $N_t^1=t+1$ points.)

Indeed, by identifying $S^1$ with complex numbers of norm 1, our pointset can be described as 
$$P=\{ z^j \; \mid \; j=0,1,2,\dots, n-1\},$$ where $z$ is the `first' $n$-th root of unity: $$z=\cos \left( 2 \pi/n \right) + i \sin \left( 2 \pi/n \right).$$  We see that, when $k$ is not a multiple of $n$, then $z^k \neq 1$, and thus
$$\sum_{j=0}^{n-1} \left( z^j \right)^k = \sum_{j=0}^{n-1} \left( z^k \right)^j = \frac{  \left( z^k \right)^n - 1}{z^k-1} = \frac{  \left( z^n \right)^k - 1}{z^k-1} = \frac{  1^k - 1}{z^k-1} = 0,$$ and thus the average of $\mathrm{Re} (x+iy)^k$ and $\mathrm{Im} (x+iy)^k$ over $P$ both equal zero.  (When $k$ is a multiple of $n$, then each term in the sum equals $1$, so only the average of $\mathrm{Im} (x+iy)^k$ equals zero over $P$.)  Therefore, $\overline{f}_P=0$ for all $f \in \mathrm{Harm}_k(S^1)$ and all $1 \leq k \leq n-1$, and thus $P$ is a spherical $t$-design for every $t \leq n-1$, as claimed.    

Let us see now how we can generalize this construction to higher dimensions.  We will assume that $d$ is odd; the case when $d$ is even can be reduced to this case.  Let $m=(d+1)/2$, and suppose that $A=\{a_1,\dots,a_m\}$ is a set of integers.  We construct $n$ points on $S^d$, as follows: for each $j=0,1,2,\dots,n-1$, we let
$$z_j= \left( \cos \left( 2 \pi j a_1/n \right) ,  \sin \left( 2 \pi j a_1/n \right), \dots ,  \cos \left( 2 \pi j a_m/n \right) ,  \sin \left( 2 \pi j a_m/n \right) \right);$$
we then set   
$$P(A)=\left\{1/\sqrt{m} \cdot z_j \; \mid \; j=0,1,2,\dots, n-1 \right\}$$ (the factor $1/\sqrt{m}$ is needed so that $P(A) \subset S^d$).
The question that we ask is then the following: for what values of $n$ and for which sets $A$ is $P(A)$ a spherical $t$-design on $S^d$?

We answer this question for $t=1$ first.  From our formula on page \pageref{dim harm}, we see that $\dim \mathrm{Harm}_1(S^d) = d+1$; clearly, the polynomials $x_i$ with $i=1,2,\dots,d+1$ form a basis for $\mathrm{Harm}_1(S^d)$.  Consequently, $P(A)$ is a spherical $1$-design if, and only if, 
$$\sum_{j=0}^{n-1} \cos \left( 2 \pi j a_i/n \right) =  \sum_{j=0}^{n-1} \sin \left( 2 \pi j a_i/n \right) = 0$$
for each $a_i \in A$; as we verified above, these equations hold whenever $a_i$ is not a multiple of $n$.  Therefore, with $$a_1=a_2=\cdots = a_m=1,$$ $P(A)$ is a spherical $1$-design for every $n \geq 2$.  (Note that the DGS bound requires that a spherical 1-design on $S^d$ has size at least $N_1^d=2$; obviously, a single point is never in mass balance on the sphere, but two or more points may be.)

Let us turn to $t=2$.  This time our formula from page \pageref{dim harm} yields
$$\dim \mathrm{Harm}_2(S^d) =  {d+2 \choose d} - {d \choose d} = {d+1 \choose 2} + d,$$ and we can verify that the set 
$$\{ x_{i_1} x_{i_2}  \; \mid \; 1 \leq i_1 < i_2 \leq d+1\} \cup \{ x_{i+1}^2 - x_i^2 \; \mid \; 1 \leq i \leq d\}$$ forms a basis for $\mathrm{Harm}_2(S^d)$.  In order to compute the average values of these functions over $P(A)$, we need to use the trigonometric identities
\begin{eqnarray*}
\sin \alpha \cdot \sin \beta &= &1/2 \cdot \left( \cos (\alpha - \beta) - \cos (\alpha + \beta) \right), \\
\cos \alpha \cdot \cos \beta &= &1/2 \cdot \left( \cos (\alpha - \beta) + \cos (\alpha + \beta) \right), \\
\sin \alpha \cdot \cos \beta &= &1/2 \cdot \left( \sin (\alpha - \beta) + \sin (\alpha + \beta) \right).
\end{eqnarray*}

We consider first $f=x_{i_1} x_{i_2}$ with $1 \leq i_1 < i_2 \leq d+1$.  Let $i_1'=\lceil i_1/2 \rceil$ and $i_2'=\lceil i_2/2 \rceil$, and note that $1 \leq i_1' \leq i_2' \leq m$.  Assume first that $i_1' < i_2'$.  With our notations, $f(z_j)$ equals the product of the cosine or sine of $2 \pi j a_{i_1'}/n$ and the cosine or sine of $2 \pi j a_{i_2'}/n$, and with our identities we can turn this into a linear combination of a cosine or sine of $2 \pi j (a_{i_1'} - a_{i_2'})/n$ and of $2 \pi j (a_{i_1'} + a_{i_2'})/n$.  As we have seen above, if neither $a_{i_1'}-a_{i_2'}$ nor 
$a_{i_1'}+a_{i_2'}$ is a multiple of $n$, then $$\sum_{j=0}^{n-1} \cos \left( 2 \pi j (a_{i_1'} - a_{i_2'})/n \right)=\sum_{j=0}^{n-1} \sin \left( 2 \pi j (a_{i_1'} - a_{i_2'})/n \right)=0$$ and 
$$\sum_{j=0}^{n-1} \cos \left( 2 \pi j (a_{i_1'} + a_{i_2'})/n \right)=\sum_{j=0}^{n-1} \sin \left( 2 \pi j (a_{i_1'} + a_{i_2'})/n \right)=0,$$ and thus the average of $f$ over $P(A)$ is zero.  Turning to the case when $i_1'=i_2'$:  note that we then must have $i_2=i_1+1$ with $i_1$ odd, so  
$$f(z_j)=\cos \left( 2 \pi j a_{i_1}/n \right) \cdot \sin \left( 2 \pi j a_{i_1}/n \right) = 1/2 \cdot \sin \left( 2 \pi j (2a_{i_1})/n \right)=
1/2 \cdot \mathrm{Im} z^{j},$$ where $z$ is the $2a_{i_1}$-th $n$-th root of unity.  Since, as noted above, $\sum_{j=0}^{n-1} \mathrm{Im} z^{j} $ equals zero for any $n$-th root of unity $z$, the average of $f$ over $P(A)$ equals zero in this case as well.

Now let $f=x_{i+1}^2-x_i^2$ with $1 \leq i \leq d$.  We then find that, for even values of $i$ we have
\begin{eqnarray*}
f(z_j) & = &   \cos^2 \left( 2 \pi j a_{i/2+1} /n \right) -  \sin^2 \left( 2 \pi j a_{i/2} /n \right) \\
&= & 1/2 \cdot \left( \cos \left( 2 \pi j (2a_{i/2+1}) /n \right)  +  \cos \left( 2 \pi j (2a_{i/2}) /n \right)  \right),
\end{eqnarray*}
so, if $2a$ is not a multiple of $n$ for any element $a \in A$, then $f$ has a zero average over $P(A)$; our computation is similar when $i$ is odd.  

Therefore, in summary, $P(A)$ is a spherical 2-design on $S^d$---that is, every polynomial in $\mathrm{Harm}_1(S^d)$ and $\mathrm{Harm}_2(S^d)$ has a zero average on $P(A)$---if, and only if, $A$ consists of $m$ distinct terms so that
\begin{itemize}
  \item no element of $A$,
  \item the difference of no two distinct elements of $A$, and
  \item the sum of no two (not-necessarily-distinct) elements of $A$
\end{itemize}  
is a multiple of $n$.  Recall that for a subset $A=\{a_1,\dots,a_m\}$ of an abelian group $G$ we defined $[1,2]_{\pm} A$ as 
$$[1,2]_{\pm} A=\{ \lambda_1 a_1+ \cdots + \lambda_m a_m \; \mid \; \lambda_1, \dots , \lambda_m \in \mathbb{Z}, \;  1 \leq |\lambda_1|+ \cdots | \lambda_m| \leq 2\},$$ so our three-part condition above can be summarized by saying that we want to assure that $A$ is an $m$-subset of $\mathbb{Z}_n$ for which $ 0 \not \in [1,2]_{\pm} A$; that is, using the terminology of Section \ref{5maxUSlimited}, $A$ is {\em 2-independent} in $\mathbb{Z}_n$. 
  
It is not hard to find specific 2-independent sets in $\mathbb{Z}_n$; for example, the set $$A=\{1,2,\dots,m\}$$ is clearly 2-independent in $\mathbb{Z}_n$ for every $n \geq 2m+1$.  Recall that we set $m=(d+1)/2$, so we have constructed explicit spherical $2$-designs on $S^d$ of size $n$ for every odd $d$ and $n \geq d+2$.  (Note that, by the DGS bound, a spherical $2$-design on $S^d$ must contain at least $d+2$ points.)  As we mentioned before, the case of even $d$ can be reduced to the odd case, but we only get spherical $2$-designs of sizes $n=d+2$ or $n \geq d+4$; it can be shown that when $d$ is even, spherical $2$-designs of size $d+3$ do not exist on $S^d$ (cf.~\cite{Mim:1990a}).\index{Mimura, Y.} So, on $S^2$, for example, we have 2-designs of size four (the regular tetrahedron) and any size $n \geq 6$, but there is no spherical 2-design formed by five points.

Moving to $t=3$, we define $A$ to be a {\em 3-independent set} in $\mathbb{Z}_n$ if $0 \not \in [1,3]_{\pm}A$, that is, in addition to our three-part condition above,
\begin{itemize}
  \item the sum of three (not-necessarily-distinct) elements of $A$ and
  \item the sum of two (not-necessarily-distinct) elements of $A$ minus any element of $A$
\end{itemize}
is never a multiple of $n$.  With similar techniques, we can prove that if  $A$ is 3-independent in $\mathbb{Z}_n$, then the set $P(A)$ defined above is a spherical 3-design on $S^d$.  

It is an interesting problem to construct 3-independent sets in $\mathbb{Z}_n$.  One quickly realizes that the set
$$\{1,3,5,\dots, 2m-1\}$$ works for every $n \geq 6m-2$, and this bound can be lowered to $n \geq 4m-1$ when $n$ is even.  We can do better yet when $n$ has a divisor $p$ with $p \equiv 5$ mod 6:
\begin{eqnarray*}
A & = & \{1,3,5, \dots, (p-2)/3 \} + \{0,p,2p, \dots, n-p\} \\ 
& = &\{2i+1+jp \; \mid \; i=0,1,\dots, (p-5)/6, \; j =0,1,\dots,  n/p-1\}
\end{eqnarray*} is then 3-independent in $\mathbb{Z}_n$, and this set has size $(p+1)/6 \cdot n/p$.
It turns out that we cannot do better: as Bajnok  and Ruzsa\index{Ruzsa, I.}\index{Bajnok, B.} proved in 2003 in \cite{BajRuz:2003a}, the minimum size $\tau_{\pm} (\mathbb{Z}_n,[1,3])$ of a 3-independent set in $\mathbb{Z}_n$ equals
$$\tau_{\pm}(\mathbb{Z}_n, [1,3]) = \left\{
\begin{array}{cl}
\left\lfloor \frac{n}{4} \right\rfloor & \mbox{if $n$ is even,}\\ \\
\left(1+\frac{1}{p}\right) \frac{n}{6} & \mbox{if $n$ is odd, has prime divisors congruent to 5 mod 6,} \\ & \mbox{and $p$ is the smallest such divisor,}\\ \\
\left\lfloor \frac{n}{6} \right\rfloor & \mbox{otherwise.}\\
\end{array}\right.$$
(The value of $\tau_{\pm}(G, [1,3])$ is not completely known in noncyclic groups $G$, and quite little is known about $\tau_{\pm}(G, [1,t])$ for $t>3$ even for cyclic groups---see Section \ref{5maxUSlimited}.)  As a consequence, we get spherical 3-designs on $S^d$ for odd $d$ and every size $n$ with
\begin{itemize}
  \item $n \geq 2d+2$ if $n$ is divisible by 4.
  \item $n \geq \left(1 - \frac{1}{p+1} \right) (3d+3)$ if $n$ has prime divisors congruent to 5 mod 6 and $p$ is the smallest such divisor, or
  \item $n \geq 3d+3$.
\end{itemize}

Explicit constructions of spherical $t$-designs get increasingly complicated as $t$ grows, the cases of `small' sizes $n$ are particularly difficult and largely open (see  \cite{Baj:1991a, Baj:1992a}).  For further information on this topic see the extensive survey paper \cite{BanBan:2009a} of Bannai and Bannai\index{Bannai, Ei.}\index{Bannai, Et.} and its nearly two hundred references.  The construction of spherical designs and, more generally, uniformly distributed pointsets on the sphere, was listed by Fields Medalist Steven Smale as \#7 on his list of the most important mathematical problems for the twenty-first century (see \cite{Sma:1998a}).\index{Smale, S.}

\newpage

\addcontentsline{toc}{section}{Caps, centroids, and the game SET}  \label{appetizer on SET}

\section*{Caps, centroids, and the game SET}

\label{appetizer on SET}

In this appetizer, we explore how one particular question in additive combinatorics connects such diverse topics as centroids of triangles, caps in affine geometry, and the card game SET.  We present each of these topics by some relevant puzzles.

Let us start with the popular and award-winning card game SET, published by SET Enterprises (cf.~\cite{Set}).  (The game was invented by the population geneticist Marsha Falco while studying epilepsy patterns in German Shepherds.)  The deck consists of 81 cards, each of which featuring four attributes: shapes (ovals, diamonds, or squiggles); number of shapes (one, two, or three); shadings (empty, striped, or solid); and colors (red, green, or purple).  Each of the 81 cards features a different combination of attributes.  The game centers on players identifying sets among the cards: three cards form a set if each of the four attributes, when considered individually on the three cards, is either all the same or all different.  For example, the three cards 
\begin{center}
\setlength{\unitlength}{.19in}
\begin{picture}(26,4)

\put(3,2){\oval(8,4)}

\put(2.75,1.3){\line(1,3){.5}}
\put(3.25,1.3){\line(1,3){.5}}

\put(3.25,1.3){\oval(1,1)[bl]}
\put(3.25,1.05){\oval(.5,.5)[r]}

\put(3.25,3.05){\oval(.5,.5)[l]}
\put(3.25,2.8){\oval(1,1)[tr]}

\put(13,2){\oval(8,4)}

\put(11.5,2){\line(2,5){.5}}
\put(11.5,2){\line(2,-5){.5}}
\put(14.5,2){\line(-2,-5){.5}}
\put(14.5,2){\line(-2,5){.5}}
\put(13.5,2){\line(2,5){.5}}
\put(13.5,2){\line(2,-5){.5}}
\put(12.5,2){\line(-2,-5){.5}}
\put(12.5,2){\line(-2,5){.5}}

\put(23,2){\oval(8,4)}

\put(23,2){\oval(.7,2.3)}
\put(21,2){\oval(.7,2.3)}
\put(25,2){\oval(.7,2.3)}

\end{picture}

\end{center}
feature all three shapes (squiggle, diamonds, and ovals), all three numbers (one, two, and three), but share the same shading (empty),  and they form a set either when they are of the same color (e.g.~all are green) or are of different colors (one green, one red, one purple).

In the most popular version of the game, twelve cards are dealt from a shuffled deck and shown face up on a table.  When a player sees a set, he or she calls ``set,'' and collects the three cards from the table.  The three cards then get replaced from the deck, and the game continues; when the deck runs out and no more sets are left, the game ends.  The player with the most sets collected wins.  Occasionally---as it turns out, with a probability of slightly more than 3 percent---the twelve cards that were dealt don't contain a set; in this case other cards are dealt from the deck one by one, until a set is present and claimed.  The question then arises: what is the maximum number of cards without containing any sets? 

We see that sets are quite ubiquitous; in fact, any pair of cards is contained in a (unique) set.  Indeed, for each of the four attributes, the two cards either share the attribute, in which case we extend it with a card that also shares that attribute, or they don't, in which case we need a card whose relevant attribute differs from that of both cards.  For example, if one card has two striped green diamonds and the other has one striped red squiggle, then we can extend them to form a set with the card containing three striped purple ovals.  Therefore, the probability that three cards don't form a set equals $78/79$, almost 99 percent.

A similar calculation shows that the probability that four cards don't form a set equals $75/79$, slightly under 95 percent.  The computations get increasingly difficult as the number of cards gets larger; the exact values have been calculated by Knuth\index{Knuth, D.} in \cite{Knu:2001a}.  It is clear that the probability that there are no sets among $k$ cards decreases rapidly as $k$ increases; as we have mentioned above, the probability that twelve cards don't contain any sets is just over 3 percent.  Our question above is equivalent to asking for the largest value of $k$ for which this probability is still not zero.  We will reveal the answer shortly.

Let us now move on to a (seemingly) very different topic: affine plane geometry.  We start by recalling two properties of the Euclidean plane: 
\begin{itemize}
  \item any two distinct points determine a unique line that contains both points, and
  \item any point $P$ and any line $l$ that does not contain $P$ determine a unique line that contains $P$ but contains no point of $l$.   
\end{itemize}
(The second property, which guarantees the existence of a line through $P$ that is {\em parallel} to $l$, is referred to as the {\em Parallel Postulate}.)  

Euclidean geometry has many other properties, of course, such as those involving distances.  But any structure on a given set $S$ and on a collection ${\cal L}$ of its subsets that satisfies the two axioms and is not trivial---that is, $S$ contains at least four points no three of which are {\em collinear} (are contained in the same element of ${\cal L}$)---is worth studying: we may simply identify $S$ with the set of points and ${\cal L}$ with the set of lines; if $S$ and ${\cal L}$ satisfy the properties we require, we say that they form an {\em affine plane}.  In particular, here we discuss finite affine planes: those that contain a finite number of points.

It turns out that the number of points in a finite affine plane cannot be arbitrary: it must be a square number.  In fact, not every square number is possible either: for example, it cannot be 36 or 100, but we still don't know if it can be 144.  One can prove that for every finite affine plane there is a positive integer $k$ so that 
\begin{itemize}
  \item there are exactly $k^2$ points and $k^2+k$ lines;
  \item each line contains $k$ points; and
  \item each point is contained in $k+1$ lines.
\end{itemize}     
The possible values of $k$ include all prime-powers, and it is one of the most famous open problems to prove that it does not include other values.  As we mentioned above, we know that $k$ cannot be 6 or 10, but we still don't know whether it can be 12.

The example that we are focusing on here is the one for $k=3$; the corresponding affine plane contains nine points and twelve lines.  Let us denote the points by $A, B, C, D, E, F, G, H,$ and $I$; we can then verify that by setting the twelve lines equal to the point-sets
$$\{A, B, C\} \; \; \; \{D,E,F\} \; \; \; \{G,H,I\} \; \; \; \{A, D, G\} \; \; \; \{B, E, H\} \; \; \; \{C, F, I\}$$
$$\{A, F, H\} \; \; \; \{B, D, I\} \; \; \; \{C, E, G\} \; \; \; \{A, E, I\} \; \; \; \{B, F, G\} \; \; \; \{C, D, H\}$$
our required properties hold.  For example: the points $B$ and $F$ determine the unique line $\{B, F, G\}$ that contains them both, and the point $A$ and the line $\{C, E, G\}$ determine the unique line $\{A, F, H\}$ that contains $A$ but none of $C, E,$ or $G$.  This affine geometry is denoted by AG(2,3) (with the 2 representing the fact that we are in the plane---affine geometries may be considered in higher dimensions as well), and we can visualize it by the following diagram:

\begin{center}
\setlength{\unitlength}{.45in}
\begin{picture}(8,8)

\put(4,2){\line(1,1){2}}
\put(3,3){\line(1,1){2}}
\put(2,4){\line(1,1){2}}
\put(2,4){\line(1,-1){2}}
\put(3,5){\line(1,-1){2}}
\put(4,6){\line(1,-1){2}}

\put(4,2){\line(1,0){2}}
\put(2,3){\line(1,0){3}}
\put(2,4){\line(1,0){4}}
\put(3,5){\line(1,0){3}}
\put(2,6){\line(1,0){2}}

\put(2,2){\line(0,1){2}}
\put(3,3){\line(0,1){3}}
\put(4,2){\line(0,1){4}}
\put(5,2){\line(0,1){3}}
\put(6,4){\line(0,1){2}}

\qbezier(3,6)(3,6.75)(3.6,7.2)
\qbezier(3.6,7.2)(4,7.5)(4.5,7.5)
\qbezier(4.5,7.5)(5.2,7.5)(5.4,7.2)
\qbezier(5.4,7.2)(6,6.75)(6,6)

\qbezier(2,2)(2,1.25)(2.6,.8)
\qbezier(2.6,.8)(3,.5)(3.5,.5)
\qbezier(3.5,.5)(4.2,.5)(4.4,.8)
\qbezier(4.4,.8)(5,1.25)(5,2)

\qbezier(2,3)(1.25,3)(.8,3.6)
\qbezier(.8,3.6)(.5,4)(.5,4.5)
\qbezier(.5,4.5)(.5,5.2)(.8,5.4)
\qbezier(.8,5.4)(1.25,6)(2,6)

\qbezier(6,2)(6.75,2)(7.2,2.6)
\qbezier(7.2,2.6)(7.5,3)(7.5,3.5)
\qbezier(7.5,3.5)(7.5,4.2)(7.2,4.4)
\qbezier(7.2,4.4)(6.75,5)(6,5)

\put(2,4){\circle*{.5}}
\put(3,3){\circle*{.5}}
\put(3,5){\circle*{.5}}
\put(4,2){\circle*{.5}}
\put(4,4){\circle*{.5}}
\put(4,6){\circle*{.5}}
\put(5,3){\circle*{.5}}
\put(5,5){\circle*{.5}}
\put(6,4){\circle*{.5}}

\put(1.7,4.2){A}
\put(2.8,2.6){D}
\put(2.7,5.2){B}
\put(4.1,1.6){G}
\put(3.8,3.6){E}
\put(3.7,6.2){C}
\put(5.1,2.6){H}
\put(4.8,4.6){F}
\put(6.1,3.6){I}

\end{picture}
\end{center}
(Note that eight of our lines are `straight' and four are ``curved.'')

Our question then is the following: What is the maximum number of points that one can find so that they do not contain all three points of any line?  Point-sets within the geometry that do not contain three collinear points are called {\em cap-sets}; we can thus rephrase our question to ask for the maximum size of a cap-set in AG$(2,3)$.  We defer answering this question until we first discuss yet another topic: centroids of triangles.

Recall that the {\em centroid} of a triangle is defined as the point where the three medians of the triangle (lines connecting a vertex with the midpoint of the opposite side)  intersect, as illustrated below: the centroid of the triangle with vertices A, B, and C is marked by M.  (It is a well-known fact that the three medians go through the same point.)  Informally, the centroid of the triangle is the point where the tip of a pin should be if one wants to balance the triangle (made of material with uniform density) on it.

\begin{center}
\setlength{\unitlength}{.5in}
\begin{picture}(5,5.2)

\put(0,0){\line(5,1){5}}
\put(0,0){\line(4,5){4}}
\put(4,5){\line(1,-4){1}}

\dashline{.2}(0,0)(4.5,3)
\dashline{.2}(2,2.5)(5,1)
\dashline{.2}(2.5,0.5)(4,5)

\put(-.5,0){A}
\put(4.2,5){B}
\put(5.2,0.7){C}
\put(3,1.5){M}

\put(0,0){\circle*{.2}}
\put(5,1){\circle*{.2}}
\put(4,5){\circle*{.2}}
\put(3,2){\circle*{.2}}

\end{picture}\end{center}

Given the coordinates of the vertices of the triangle in a coordinate system, one can determine the coordinates of the centroid by taking the arithmetic averages of the coordinates of the three vertices; for example, if A=$(a_1,a_2)$, B=$(b_1,b_2)$, and C=$(c_1,c_2)$, then $$\mathrm{M}=\left( \frac{a_1+b_1+c_1}{3},  \frac{a_2+b_2+c_2}{3} \right).$$  This formulation allows us to broaden our definition to include the degenerate case when the three points are collinear; for example, the centroid of the collinear points A=$(0,0)$, B=$(2,1)$, and C=$(16,8)$ is the point M=$(6,3)$:

\begin{center}
\setlength{\unitlength}{.2in}
\begin{picture}(16,8)

\put(0,0){\line(2,1){2}}
\put(2,1){\line(2,1){14}}

\put(-.9,-.1){A}
\put(1.6,1.4){B}
\put(16.6,7.8){C}
\put(5.4,3.4){M}

\put(0,0){\circle*{.5}}
\put(2,1){\circle*{.5}}
\put(16,8){\circle*{.5}}
\put(6,3){\circle*{.5}}

\end{picture}\end{center}

Our question about centroids is the following: What is the maximum number of lattice points that one can find in the integer lattice $\mathbb{Z}^2$ (points in the coordinate plane with integer coordinates) so that no three of them have their centroid at a lattice point? 

Well, let us reveal how the three questions we posed relate to one another and how they are special cases of a quantity we investigate in Section \ref{5maxRfixed} of this book.  Namely, our interest is in finding the maximum size of a {\em weakly zero-$3$-sum-free set} in the group $\mathbb{Z}_3^r$, defined as the maximum size---denoted by $\tau \hat{\;} (\mathbb{Z}_3^r, 3)$---of  a subset of $\mathbb{Z}_3^r$ without any three distinct elements adding to zero; that is, 
the quantity 
$$\tau \hat{\;} (\mathbb{Z}_3^r, 3)=\max \{|A| \; \mid \; A \subseteq \mathbb{Z}_3^r, 0 \not \in 3 \hat{\;} A\}.$$ (The word `weak' signifies that we only disallow the sum of three {\em distinct} elements to be zero.)  
The problem of finding $\tau \hat{\;} (\mathbb{Z}_3^r, 3)$ seems forbiddingly difficult at the present time.  Exact values are only known for $r \leq 6$:
$$\begin{array}{|c||r|r|r|c|c|c|} \hline 
r & 1 & 2 & 3 & 4 & 5 & 6 \\ \hline
\tau\hat{\;}(\mathbb{Z}_3^r,3) & 2 & 4 & 9 & 20 & 45 & 112 \\ \hline
\end{array}$$
(See \cite{GaoTha:2004a} by Gao and Thangadurai\index{Gao, W.}\index{Thangadurai, R.} and its references for the first five entries and \cite{Pot:2008a} by Potechin for the last.)
Here we only present the maximum-sized examples for $r \leq 4$: 

$\tau \hat{\;} (\mathbb{Z}_3, 3)=2:$
\begin{center}
\setlength{\unitlength}{.5in}
\begin{picture}(2,0)

\put(0,0){\line(1,0){2}}

\put(0,0){\circle*{.2}}
\put(1,0){\circle*{.2}}
\put(2,0){\circle{.2}}

\end{picture}\end{center}

$\tau \hat{\;} (\mathbb{Z}_3^2, 3)=4:$
\begin{center}
\setlength{\unitlength}{.5in}
\begin{picture}(2,2)

\put(0,0){\line(1,0){2}}
\put(0,1){\line(1,0){2}}
\put(0,2){\line(1,0){2}}
\put(0,0){\line(0,1){2}}
\put(1,0){\line(0,1){2}}
\put(2,0){\line(0,1){2}}

\put(0,0){\circle*{.2}}
\put(1,0){\circle*{.2}}
\put(0,1){\circle*{.2}}
\put(1,1){\circle*{.2}}
\put(2,0){\circle{.2}}
\put(2,1){\circle{.2}}
\put(2,2){\circle{.2}}
\put(1,2){\circle{.2}}
\put(0,2){\circle{.2}}

\end{picture}\end{center}

$\tau \hat{\;} (\mathbb{Z}_3^3, 3)=9:$
\begin{center}
\setlength{\unitlength}{.5in}
\begin{picture}(8,2)

\put(0,0){\line(1,0){2}}
\put(0,1){\line(1,0){2}}
\put(0,2){\line(1,0){2}}
\put(0,0){\line(0,1){2}}
\put(1,0){\line(0,1){2}}
\put(2,0){\line(0,1){2}}

\put(0,0){\circle*{.2}}
\put(1,0){\circle{.2}}
\put(0,1){\circle{.2}}
\put(1,1){\circle{.2}}
\put(2,0){\circle*{.2}}
\put(2,1){\circle{.2}}
\put(2,2){\circle*{.2}}
\put(1,2){\circle{.2}}
\put(0,2){\circle*{.2}}

\put(3,0){\line(1,0){2}}
\put(3,1){\line(1,0){2}}
\put(3,2){\line(1,0){2}}
\put(3,0){\line(0,1){2}}
\put(4,0){\line(0,1){2}}
\put(5,0){\line(0,1){2}}

\put(3,0){\circle{.2}}
\put(4,0){\circle{.2}}
\put(3,1){\circle{.2}}
\put(4,1){\circle*{.2}}
\put(5,0){\circle{.2}}
\put(5,1){\circle{.2}}
\put(5,2){\circle{.2}}
\put(4,2){\circle{.2}}
\put(3,2){\circle{.2}}

\put(6,0){\line(1,0){2}}
\put(6,1){\line(1,0){2}}
\put(6,2){\line(1,0){2}}
\put(6,0){\line(0,1){2}}
\put(7,0){\line(0,1){2}}
\put(8,0){\line(0,1){2}}

\put(6,0){\circle{.2}}
\put(7,0){\circle*{.2}}
\put(6,1){\circle*{.2}}
\put(7,1){\circle{.2}}
\put(8,0){\circle{.2}}
\put(8,1){\circle*{.2}}
\put(8,2){\circle{.2}}
\put(7,2){\circle*{.2}}
\put(6,2){\circle{.2}}

\end{picture}\end{center}

$\tau \hat{\;} (\mathbb{Z}_3^4, 3)=20:$

\begin{center}
\setlength{\unitlength}{.5in}
\begin{picture}(8,8)

\put(0,0){\line(1,0){2}}
\put(0,1){\line(1,0){2}}
\put(0,2){\line(1,0){2}}
\put(0,0){\line(0,1){2}}
\put(1,0){\line(0,1){2}}
\put(2,0){\line(0,1){2}}

\put(0,0){\circle*{.2}}
\put(1,0){\circle{.2}}
\put(0,1){\circle{.2}}
\put(1,1){\circle{.2}}
\put(2,0){\circle*{.2}}
\put(2,1){\circle{.2}}
\put(2,2){\circle*{.2}}
\put(1,2){\circle{.2}}
\put(0,2){\circle*{.2}}

\put(3,0){\line(1,0){2}}
\put(3,1){\line(1,0){2}}
\put(3,2){\line(1,0){2}}
\put(3,0){\line(0,1){2}}
\put(4,0){\line(0,1){2}}
\put(5,0){\line(0,1){2}}

\put(3,0){\circle{.2}}
\put(4,0){\circle{.2}}
\put(3,1){\circle{.2}}
\put(4,1){\circle*{.2}}
\put(5,0){\circle{.2}}
\put(5,1){\circle{.2}}
\put(5,2){\circle{.2}}
\put(4,2){\circle{.2}}
\put(3,2){\circle{.2}}

\put(6,0){\line(1,0){2}}
\put(6,1){\line(1,0){2}}
\put(6,2){\line(1,0){2}}
\put(6,0){\line(0,1){2}}
\put(7,0){\line(0,1){2}}
\put(8,0){\line(0,1){2}}

\put(6,0){\circle{.2}}
\put(7,0){\circle*{.2}}
\put(6,1){\circle*{.2}}
\put(7,1){\circle{.2}}
\put(8,0){\circle{.2}}
\put(8,1){\circle*{.2}}
\put(8,2){\circle{.2}}
\put(7,2){\circle*{.2}}
\put(6,2){\circle{.2}}

\put(0,3){\line(1,0){2}}
\put(0,4){\line(1,0){2}}
\put(0,5){\line(1,0){2}}
\put(0,3){\line(0,1){2}}
\put(1,3){\line(0,1){2}}
\put(2,3){\line(0,1){2}}

\put(0,3){\circle{.2}}
\put(1,3){\circle*{.2}}
\put(0,4){\circle*{.2}}
\put(1,4){\circle{.2}}
\put(2,3){\circle{.2}}
\put(2,4){\circle*{.2}}
\put(2,5){\circle{.2}}
\put(1,5){\circle*{.2}}
\put(0,5){\circle{.2}}

\put(3,3){\line(1,0){2}}
\put(3,4){\line(1,0){2}}
\put(3,5){\line(1,0){2}}
\put(3,3){\line(0,1){2}}
\put(4,3){\line(0,1){2}}
\put(5,3){\line(0,1){2}}

\put(3,3){\circle{.2}}
\put(4,3){\circle{.2}}
\put(3,4){\circle{.2}}
\put(4,4){\circle*{.2}}
\put(5,3){\circle{.2}}
\put(5,4){\circle{.2}}
\put(5,5){\circle{.2}}
\put(4,5){\circle{.2}}
\put(3,5){\circle{.2}}

\put(6,3){\line(1,0){2}}
\put(6,4){\line(1,0){2}}
\put(6,5){\line(1,0){2}}
\put(6,3){\line(0,1){2}}
\put(7,3){\line(0,1){2}}
\put(8,3){\line(0,1){2}}

\put(6,3){\circle*{.2}}
\put(7,3){\circle{.2}}
\put(6,4){\circle{.2}}
\put(7,4){\circle{.2}}
\put(8,3){\circle*{.2}}
\put(8,4){\circle{.2}}
\put(8,5){\circle*{.2}}
\put(7,5){\circle{.2}}
\put(6,5){\circle*{.2}}

\put(0,6){\line(1,0){2}}
\put(0,7){\line(1,0){2}}
\put(0,8){\line(1,0){2}}
\put(0,6){\line(0,1){2}}
\put(1,6){\line(0,1){2}}
\put(2,6){\line(0,1){2}}

\put(0,6){\circle{.2}}
\put(1,6){\circle{.2}}
\put(0,7){\circle{.2}}
\put(1,7){\circle*{.2}}
\put(2,6){\circle{.2}}
\put(2,7){\circle{.2}}
\put(2,8){\circle{.2}}
\put(1,8){\circle{.2}}
\put(0,8){\circle{.2}}

\put(3,6){\line(1,0){2}}
\put(3,7){\line(1,0){2}}
\put(3,8){\line(1,0){2}}
\put(3,6){\line(0,1){2}}
\put(4,6){\line(0,1){2}}
\put(5,6){\line(0,1){2}}

\put(3,6){\circle{.2}}
\put(4,6){\circle{.2}}
\put(3,7){\circle{.2}}
\put(4,7){\circle{.2}}
\put(5,6){\circle{.2}}
\put(5,7){\circle{.2}}
\put(5,8){\circle{.2}}
\put(4,8){\circle{.2}}
\put(3,8){\circle{.2}}

\put(6,6){\line(1,0){2}}
\put(6,7){\line(1,0){2}}
\put(6,8){\line(1,0){2}}
\put(6,6){\line(0,1){2}}
\put(7,6){\line(0,1){2}}
\put(8,6){\line(0,1){2}}

\put(6,6){\circle{.2}}
\put(7,6){\circle{.2}}
\put(6,7){\circle{.2}}
\put(7,7){\circle*{.2}}
\put(8,6){\circle{.2}}
\put(8,7){\circle{.2}}
\put(8,8){\circle{.2}}
\put(7,8){\circle{.2}}
\put(6,8){\circle{.2}}

\end{picture}\end{center}
Hopefully, these diagrams are self-explanatory: For example, for rank two, we show the four elements $(0,0), (0,1), (1,0),$ and $ (1,1)$ that form a weakly zero-3-sum-free subset in $\mathbb{Z}_3^2$, and within the three grids representing $\mathbb{Z}_3^3$, we feature the weakly zero-3-sum-free subset
$$\{(0,0,0), (0,0,2), (0,2,0), (0,2,2),(1,1,1),(2,0,1), (2,1,0), (2,1,2), (2,2,1)\}.$$  Besides verifying that the sets above are weakly zero-3-sum-free in their respective groups, one would also need to prove that they are of maximum size---these proofs get increasingly difficult as the rank of the group increases.  

So how do our three questions above relate to these values?  A key is the following observation: three elements of $\mathbb{Z}_3$ add to zero if, and only if, they are all distinct ($0+1+2=0$) or all the same (e.g. $1+1+1=0$).  If we assign coordinates to each attribute in the game SET (for example, the first coordinate represents shape: namely, 0, 1, and 2 denote ovals, diamonds, and squiggles, respectively; and similarly with the other three attributes), then each of our 81 cards corresponds to a unique element of $\mathbb{Z}_3^4$, with three cards forming a set exactly when their corresponding group elements add to $(0,0,0,0)$, the zero element of $\mathbb{Z}_3^4$.  Therefore, the maximum number of cards without having any sets among them is $\tau \hat{\;} (\mathbb{Z}_3^4, 3)=20.$    

Now let us turn to the question of finding the maximum size of a cap-set in the geometry AG$(2,3)$.  To start, we identify our nine points with the elements of $\mathbb{Z}_3^2$, as follows:
$A=(0,0), \; B=(0,1), \; C=(0,2), \; D=(1,0), \; E=(1,1), \; F=(1,2), \; G=(2,0), \; H=(2,1), \; I=(2,2).$    

Suppose that $X$, $Y$, and $Z$ are three of these points.  Recall that they are collinear whenever the vectors from $X$ to $Y$ and from $X$ to $Z$ are parallel, that is, when the vector $Y-X$ is a scalar multiple of the vector $Z-X$.  Since our points, and therefore our vectors, are considered here as elements of $\mathbb{Z}_3^2$, that scalar can only be 0, 1, or 2---but it cannot be 0 (since $X$ and $Y$ are distinct) and cannot be 1 (since $Y$ and $Z$ are distinct).  Therefore, $X$, $Y$, and $Z$ are collinear exactly when $Y-X=2\cdot (Z-X)$ in $\mathbb{Z}_3^2$ or, equivalently, $X+Y+Z=(0,0)$.  Indeed, it is easy to check that for each of our twelve lines, the sum of the three points within them equals zero; for example, two of our lines were $\{A,D,G\}$ and $\{B,F,G\}$, and indeed we have $A+D+G=(0,0)$ and $B+F+G=(0,0)$.  (In fact, the twelve lines exhibit all such 3-subsets: there are three of them where the first coordinates are equal and the second coordinates are all distinct, three where the second coordinates are equal and the first coordinates are all distinct, and six where both the first coordinates and second coordinates are pairwise distinct.)  Therefore, a set of points in AG$(2,3)$ is a cap-set exactly when the corresponding subset of $\mathbb{Z}_3^2$ is a weakly zero-3-sum-free set, and thus the maximum cap-set size in AG$(2,3)$ equals  $\tau \hat{\;} (\mathbb{Z}_3^2, 3)=4$.  In particular, the cap-set of size four that corresponds to our example of a weakly zero-3-sum-free set illustrated above consists of $A$, $B$, $D$, and $E$.     

Finally: the maximum size of a subset of $\mathbb{Z}^2$ without the centroid of any three of its points at a lattice point.  Note that the centroid of points $A=(a_1,a_2)$,  B=$(b_1,b_2)$, and C=$(c_1,c_2)$ is a lattice point if, and only if, $a_1+b_1+c_1$ and $a_2+b_2+c_2$ are both divisible by 3.  Now the sum of three integers is divisible by 3 exactly when they all leave the same remainder mod 3, or when they all leave different remainders mod 3.  

Let us map the lattice points of $\mathbb{Z}^2$ to elements of $\mathbb{Z}_3^2$ based on the reminders of their two coordinates.  By the previous paragraph, if three points of $\mathbb{Z}^2$ are mapped to the same element of $\mathbb{Z}_3^2$, then the centroid of the three points is a lattice point.  Furthermore, the same holds if the three points of $\mathbb{Z}^2$ are mapped to three distinct elements of $\mathbb{Z}_3^2$ whose sum equals $(0,0)$.  Therefore, if our point-set is such that no three of its points have their centroid at a lattice point, then they must be mapped to a subset of $\mathbb{Z}_3^2$ that is weakly zero-3-sum-free, and no more than two points can be mapped to the same element of $\mathbb{Z}_3^2$.  But then the maximum size of a subset of $\mathbb{Z}^2$ without the centroid of any three of its points at a lattice point equals $2 \cdot \tau \hat{\;} (\mathbb{Z}_3^2, 3)= 2 \cdot 4 =8$.  Our argument also shows that an example for such a set of eight points could be
$$\{(0,0), (0,1), (1,0), (1,1),  (3,3), (3,7), (4,9), (10,10)\}.$$

Hopefully, our set of examples provided an inviting appetizer for the entrees in this book.  One of these entrees is the question of finding $\tau \hat{\;} (\mathbb{Z}_3^r, 3)$, which (Fields Medalist and Breakthrough Prize Winner) Terence Tao\index{Tao, T.} calls ``perhaps [his] favorite open question'' (see \cite{Tao:2007a}).

\newpage

\addcontentsline{toc}{section}{How many elements does it take to span a group?}

\section*{How many elements does it take to span a group?}

Given a finite abelian group $G$ and a positive integer $h$, we are interested in finding the minimum possible size of a subset of $G$ so that each element of $G$ can be written as the sum of (exactly) $h$ elements of our set.  

Recall that the $h$-fold sumset of a subset $A$ of $G$, denoted by $hA$, consists of all possible $h$-term sums of (not necessarily distinct) elements of $A$; more formally, for $A=\{a_1,\dots,a_m\}$, we have
$$hA=\left \{\Sigma_{i=1}^m \lambda_i a_i \; \mid \; \lambda_i \in \mathbb{N}_0, \; \Sigma_{i=1}^m \lambda_i=h \right\}.$$  Thus our task is to find the minimum possible value of $m$ for which $G$ contains an $m$-subset $A$ so that $hA=G$.

Let us consider an example: suppose that $G$ is the cyclic group $\mathbb{Z}_{10}$ and $h=2$.  After some experimentation, one can find a variety of 5-subsets that yield the whole group; for example, with $A=\{0,1,2,5,7\}$, we get each element: $0=0+0$, $1=0+1$, $2=0+2$ (or $1+1$, or $5+7$), and so on, thus $2A=\mathbb{Z}_{10}$.  Can we do better?  That is, is there a smaller subset that generates the whole group?  

Well, it is easy to see that for $2A$ to contain all $10$ elements of our group, $A$ must have size at least four.  Indeed, more generally, even if all possible $h$-fold sums of an $m$-subset of $G$ yield distinct elements, the size of $hA$ is at most $ {m+h-1 \choose h}$, so to get all elements of the group, we must have
$$ {m+h-1 \choose h} \geq |G|.$$ Therefore, in our example, $m$ must be at least four.

Now we prove that, in fact, $m$ has to be at least five.  Suppose, indirectly, that there exists a set $A=\{a_1,a_2,a_3,a_4\}$ of size four for which $2A=\mathbb{Z}_{10}$.   We then have
$$2A=\{2a_1,2a_2,2a_3,2a_4, a_1+a_2,a_1+a_3,a_1+a_4,a_2+a_3,a_2+a_4,a_3+a_4\}.$$
Therefore, the ten elements listed must all be distinct, so are, in some order, equal to the elements of $\mathbb{Z}_{10}$; in particular, exactly five of them are even and five are odd.  But then exactly one of the last six elements listed in $2A$ is even; say it is $a_1+a_2$.  Now $$(a_1+a_2)+(a_3+a_4)$$ is the sum of an even and an odd number, thus odd.  However, this sum also equals $$(a_1+a_3)+(a_2+a_4),$$ the sum of two odd values, which is thus even.  This is a contradiction, so the minimum value of $m$ for which $\mathbb{Z}_{10}$ contains an $m$-subset $A$ with $2A=\mathbb{Z}_{10}$ equals five.

Staying with $\mathbb{Z}_{10}$, let us try to find the answer for $h=3$.  From our general inequality above, we see that no subset of size less than three will have a 3-fold sumset of size ten.  While there are numerous sets of size four that yield the entire group---for example, $\{0,1,4,7\}$ does---it seems impossible to find one of size three.  We can prove this, as follows.

Assume that there is a set $A=\{a_1,a_2,a_3\}$ of size three for which $3A=\mathbb{Z}_{10}$.  We see that $$3A=\{3a_1,3a_2,3a_3,2a_1+a_2,a_1+2a_2,2a_1+a_3,a_1+2a_3,2a_2+a_3,a_2+2a_3,a_1+a_2+a_3\}.$$  Again, the ten elements listed must all be distinct, so are, in some order, equal to the elements of $\mathbb{Z}_{10}$.  Thus, adding the elements in $3A$ and in $\mathbb{Z}_{10}$ should yield the same answer, so
$$10a_1+10a_2+10a_3=0+1+\cdots+9$$ or $0=5$, which is a contradiction.   

Letting $\phi (G,h)$ denote the minimum possible size of a subset of $G$ whose $h$-fold sumset contains each element of $G$, we just proved that $\phi (\mathbb{Z}_{10}, 2)=5$ and $\phi (\mathbb{Z}_{10}, 3)=4$.  


Let us now consider the variation where we are allowed not only to add terms of our subset but subtract them as well: in particular, we are looking for the smallest value of $m$ for which an $m$-subset of $G$ exists whose $h$-fold signed sumset is all of $G$; we let $\phi_{\pm} (G,h)$ denote this value.  Recall that the $h$-fold signed sumset of an $m$-subset $A=\{a_1,\dots,a_m\}$ of $G$ is the set 
$$h_{\pm}A=\left \{\Sigma_{i=1}^m \lambda_i a_i \; \mid \; \lambda_i \in \mathbb{Z}, \; \Sigma_{i=1}^m |\lambda_i|=h \right\}.$$     

Examining the cyclic group $\mathbb{Z}_{10}$, we run into a bit of uncertainty at $h=1$ already.  We would think that a set $A=\{a_1,a_2,a_3,a_4,a_5\}$ of size five is sufficient as its 1-fold signed sumset is $1_{\pm} A=\{\pm a_1, \pm a_2, \pm a_3, \pm a_4, \pm a_5\}$, which seems to contain ten elements.  However, we realize that, in $\mathbb{Z}_{10}$, every element and its inverse have the same parity; therefore, $1_{\pm} A$ cannot contain five even and five odd elements.  Since a set of size six with a 1-fold signed sumset of size ten can easily be found (for example, $A=\{0,1,2,3,4,5\}$), we conclude that $\phi_{\pm} (\mathbb{Z}_{10},1)=6$.  

It turns out that the case of $h=2$ is easy: a set $A=\{a_1,a_2\}$ of size two is clearly insufficient, as its 2-fold sumset, $$2_{\pm} A=\{\pm 2 a_1, \pm 2 a_2, \pm a_1 \pm a_2\},$$ will be of size at most eight; but a set of size three  whose 2-fold signed sumset is of size ten is not hard to find (e.g. $\{0,1,4\}$).  Therefore, $\phi_{\pm} (\mathbb{Z}_{10},2)=3$.  

Now let us attempt to find $\phi_{\pm} (\mathbb{Z}_{10},3)$.  Here we see that two elements may be enough: for $A=\{a_1,a_2\}$, we have 
$$3_{\pm} A=\{\pm 3 a_1, \pm 3 a_2, \pm 2 a_1 \pm a_2, \pm a_1 \pm 2a_2 \};$$ we see twelve (not necessarily different) elements listed.  We will prove, however, that this cannot yield all of $\mathbb{Z}_{10}$.  (Further analysis would prove that, in fact, $3_{\pm} A$ has size at most eight.)  

Note that, if $a_1$ and $a_2$ have the same parity (that is, are both even or both odd), then the twelve elements listed all share that parity; in particular, $3_{\pm} A$ has size at most five.  Therefore, $a_1$ and $a_2$ are of different parity; we will assume that $a_1$ is even and $a_2$ is odd, and therefore we see that the elements $\pm 3a_1$ and $\pm a_1 \pm 2a_2$ are even, and $\pm 3a_2$ and $\pm 2a_1 \pm a_2$ are odd.  We will need the fact that, if $a_1$ is even, then $5a_1=0$ in $\mathbb{Z}_{10}$.  

We also see that if $a_1 =0$, then $$3_{\pm} A=\{0,\pm a_2, \pm 2a_2, \pm 3a_2 \},$$ which is of size at most seven; so we assume that $a_1 \neq 0$. In that case, $\pm 3 a_1 \neq 0$ as well.  If 0 were to be an element of $3_{\pm} A$, then it would have to be one of $\pm a_1 \pm 2a_2$; let us assume here that $a_1+2a_2=0$ and thus $-a_1-2a_2=0$ as well.  (The case when $a_1-2a_2=0$ and $-a_1+2a_2=0$ can be examined the same way.)  Thus we find that the elements 2, 4, 6, and 8 of $\mathbb{Z}_{10}$ are, in some order, the elements $\pm 3a_1$ and $\pm (a_1-2a_2)$; in particular, these four elements must be distinct.  But from $-2a_2=a_1$ and $5a_1=0$  we get
$$a_1-2a_2=2a_1=2a_1-5a_1=-3a_1,$$ which is a clear contradiction with the four elements being distinct.  Therefore, 0 is not an element of $3_{\pm} A$, and thus $\phi_{\pm} (\mathbb{Z}_{10},3) \geq 3$.  Since a 3-subset of $\mathbb{Z}_{10}$ with a 3-fold signed sumset of size ten can be found easily ($\{0,1,2\}$ will do), we conclude that $\phi_{\pm} (\mathbb{Z}_{10},3)=3$.  

Using techniques similar to the ones we just saw, we can find the values of $\phi(\mathbb{Z}_{10},h)$ and $\phi_{\pm} (\mathbb{Z}_{10},h)$ for all values of $h$:
$$\begin{array}{||c||r|r|r|r|r|r|r|r|r||} \hline \hline
 & h=1 & h=2 & h=3 & h=4 & h=5 & h=6 & h=7 & h=8 & h \geq 9 \\ \hline \hline
\phi(\mathbb{Z}_{10},h) & 10 & 5 & 4 & 3 & 3 & 3 & 3 & 3 & 2   \\ \hline
\phi_{\pm} (\mathbb{Z}_{10},h) & 6 & 3 & 3 & 3 & 2 & 2 & 2 & 2 & 2  \\ \hline \hline
\end{array}$$

The problem of determining $\phi (G,h)$ and $\phi_{\pm} (G,h)$ for a general abelian group $G$ is widely open (not to mention other variations where restricted sums are considered or where the number of terms varies)---Chapter \ref{ChapterSpanning} discusses all that is known (to this author) on this subject. 

\newpage

\addcontentsline{toc}{section}{In pursuit of perfection}

\section*{In pursuit of perfection}

We are interested in the ``speed'' with which a given subset of a group generates the entire group.    As a (very nice) example, consider the set $A=\{3,4\}$ in the group $\mathbb{Z}_{25}$.  How ``fast'' does the set generate the group?  Let us explain what we mean by this question.  

Note that either element of $A$ alone generates $\mathbb{Z}_{25}$: since 3 and 4 are both relatively prime to 25, their multiples will yield all group elements, but this is not very ``fast.'' For example, to get the element 22, we need to add up 24 terms of 3s.  Even if subtraction is also allowed, the fastest way to get to 14 is to subtract 3 (from 0) twelve times.  The situation is not better with using the element 4 alone either.  However, if both 3s and 4s can be used, and we can both add and subtract them, then we can verify that each element of $\mathbb{Z}_{25}$ can be generated with three or fewer terms:
$$\begin{array}{|rcl|rcl|rcl|rcl|} \hline 
1& =& 4-3    & 7& =& 3+4    & 13& =& -4-4-4 & 19& =& -3-3    \\
2& =& 3+3-4  & 8& =& 4+4    & 14& =& -3-4-4 & 20& =& 3-4-4\\
3& =& 3      & 9& =& 3+3+3  & 15& =& -3-3-4 & 21& =& -4   \\
4& =& 4      & 10& =& 3+3+4 & 16& =& -3-3-3 & 22& =& -3  \\
5& =& 4+4-3  & 11& =& 3+4+4 & 17& =& -4-4   & 23& =& 4-3-3 \\
6& =& 3+3    & 12& =& 4+4+4 & 18& =& -3-4   & 24& =& 3-4 \\ \hline
\end{array}$$

As the table indicates, every element of $\mathbb{Z}_{25}$ can be generated by a signed sum of at most three terms of $A$.  (We consider $ 0$ to be generated trivially.)  We therefore call $A=\{3,4\}$ a \emph{3-spanning set} in $\mathbb{Z}_{25}$.

More generally, given a finite abelian group $G$, a subset $A=\{a_1,a_2,\dots, a_m\}$ of $G$, and a nonnegative integer $s$, we 
say that $A$ is an \emph{$s$-spanning set} in $G$, if every element of the group can be written as a linear combination $$\lambda_1a_1+\lambda_2a_2+\cdots +\lambda_m a_m$$ for some integers $\lambda_1, \lambda_2, \dots , \lambda_m$ with $$|\lambda_1|+|\lambda_2|+\cdots +|\lambda_m| \leq s.$$
Using our notations and terminology introduced previously, we can say that $A$ is an $s$-spanning set of $G$ if the $[0,s]$-fold signed sumset of $A$ is the entire group, that is, $$[0,s]_{\pm} A=\cup_{h=0}^s (h_{\pm} A)=G.$$In our example above, for $A=\{3,4\}$ in $\mathbb{Z}_{25}$, we find that
\begin{itemize}
\item $0_{\pm} A = \{ 0 \};$
\item $1_{\pm} A = \{ \pm 3, \; \pm 4 \}= \{ 3, 4, 21, 22 \};$
\item $2_{\pm} A = \{ \pm 2 \cdot 3, \; \pm 3 \pm 4, \; \pm 2 \cdot 4 \}= \{1, 6, 7, 8, 17, 18, 19, 24 \};$ and
\item $3_{\pm} A = \{ \pm 3 \cdot 3, \; \pm 2 \cdot 3 \pm 4, \; \pm 3 \pm 2 \cdot 4, \; \pm 3 \cdot 4 \}= \{2, 5, 9, 10, 11, 12, 13, 14, 15, 16, 20, 23 \}.$
\end{itemize}
Therefore, $$[0,3]_{\pm} A= (0_{\pm} A) \cup (1_{\pm} A ) \cup (2_{\pm} A) \cup (3_{\pm} A) =\mathbb{Z}_{25}.$$

In fact, the set $A=\{3,4\}$ has a remarkable property: the $$1+4+8+12=25$$ possible signed sums in $[0,3]_{\pm} A$ are all distinct elements of $\mathbb{Z}_{25}$; in other words, every element of the group can be written \emph{uniquely} as a signed sum of at most three elements of $A$.  We call such a set \emph{perfect}; perfect spanning sets generate the group most efficiently.  Unfortunately, they rarely exist!

Just how rare, we do not exactly know in general.  If all elements generated are distinct, then $[0,s]_{\pm} A$ must have size   
$$a(m,s)=\sum_{i \geq 0} {m \choose i} { s \choose i} 2^i$$ (see Section \ref{0.2.4} in Part I for a proof), and if each element of $G$ must be generated, then this quantity must equal the order of $G$.  Therefore, a necessary condition for $G$ to have a perfect $s$-spanning set of size $m$ is that $|G|=a(m,s)$.

The first few values of the function $a(m,s)$ are as follows:
    
\label{a(j,k)table} \begin{tabular}{|c||c|c|c|c|c|c|c|}  \hline
$a(m,s)$ & $s=0$  & $s=1$  & $s=2$  & $s=3$  & $s=4$  & $s=5$  & $s=6$ \\ \hline \hline
$m=0$ &1& 1& 1& 1& 1& 1& 1 \\ \hline
$m=1$ &1& 3& 5& 7& 9& 11& 13 \\ \hline
$m=2$ &1& 5& 13& 25& 41& 61& 85 \\ \hline
$m=3$ &1& 7& 25& 63& 129& 231& 377 \\ \hline
$m=4$ &1& 9& 41& 129& 321& 681& 1289 \\ \hline
$m=5$ &1& 11& 61& 231& 681& 1683& 3653 \\ \hline
$m=6$ &1& 13& 85& 377& 1289& 3653& 8989 \\ \hline
\end{tabular}

\noindent It may be helpful to observe that, besides the formula we gave for $a(m,s)$ above, the function can also be evaluated using the recurrence relation 
$$a(m,s)=a(m-1,s)+a(m-1,s-1)+a(m,s-1)$$ (together with the initial conditions $a(m,0)=a(0,s)=1$).  As we see, our function resembles Pascal's Triangle, except here we get to the entry in row $m$ and column $s$ by adding not only the entries directly above and to the left, but the entry in the ``above-left'' position as well.

We have already seen earlier that the set $A=\{3,4\}$ is a perfect 3-spanning set in $G=\mathbb{Z}_{25}$.  Here are some others:
\begin{itemize}

\item In any group of odd order, the set of nonzero elements can be partitioned into parts $K$ and $-K$, and both $K$ and $-K$ are perfect 1-spanning sets in the group.  For example, the set $\{1,2,\dots,m\}$ is a perfect 1-spanning set in $\mathbb{Z}_{2m+1}$.

\item The set $\{a\}$ is a perfect $s$-spanning set in $\mathbb{Z}_{2s+1}$ as long as $a$ and $2s+1$ are relatively prime.

\item The sets $\{1, 2s+1\}$ and $\{s,s+1\}$ are perfect $s$-spanning sets in $\mathbb{Z}_{2s^2+2s+1}$.

\end{itemize}

The first two statements are easy to see.  The claim regarding $\{s,s+1\}$ we just demonstrated for $s=3$ and $n=25$; we can illustrate the claim about $\{1,2s+1\}$ for the same parameters by exhibiting the values of
$$\lambda_1 \cdot 1 + \lambda_2 \cdot 7$$ in $\mathbb{Z}_{25}$ for all relevant coefficients $\lambda_1$ and $\lambda_2$:
$$\begin{array}{||l||c|c|c|c|c|c|c||} \hline \hline
  & \lambda_1=-3  & \lambda_1=-2  & \lambda_1=-1 & \lambda_1=0 &  \lambda_1=1 & \lambda_1=2  & \lambda_1=3 \\ \hline \hline 
\lambda_2=3 &  &  &  & 21 & &  & \\ \hline 
\lambda_2=2 &   &  & 13 & 14 & 15 & &  \\ \hline 
\lambda_2=1 &   & 5 & 6 & 7 & 8 & 9 &  \\ \hline 
\lambda_2=0 &  22 & 23 & 24  &  0 & 1 & 2 & 3 \\ \hline 
\lambda_2=-1 &   & 16 &  17 &  18 &  19 & 20 &  \\ \hline 
\lambda_2=-2 &   &  & 10 & 11 & 12 & &  \\ \hline 
\lambda_2=-3 &  &  &  & 4 & & &  \\ \hline \hline
\end{array}$$ 
The table shows a nice pattern for how the elements of the group arise.  (The general proofs for $\{1,2s+1\}$ and $\{s,s+1\}$ are provided on pages \pageref{proofofperfex} and \pageref{proofofp}, respectively.)  

Note that our three types of perfect spanning sets correspond to the cases 
\begin{itemize}

\item $s=1$ and $|G|=a(m,1)=2m+1$;

\item $m=1$ and $|G|=a(1,s)=2s+1$; and

\item $m=2$ and $|G|=a(2,s)=2s^2+2s+1$.

\end{itemize}
respectively.  Clearly, the cases of $s=1$ and $m=1$ are completely characterized by our list above:  the only possibilities are those listed.  For $m=2$, we may have other examples, though none are known.  We are also not aware of any perfect spanning sets for $s \geq 2$ and $m \geq 3$.  It might be an interesting problem to find and classify all perfect spanning sets (or to prove that no others exist besides the ones listed). 

As a modest attempt toward such a classification, we prove that there is neither a perfect 3-spanning set of size two, nor a perfect 2-spanning set of size three in $\mathbb{Z}_5^2$. \label{no perfects in Z25} (Note that $a(2,3)=a(3,2)=25$.  We have already seen that $\mathbb{Z}_{25}$ has perfect 3-spanning sets of size two; we know, via a computer search, that it has no perfect 2-spanning sets of size three.  Our proofs below show that, sometimes, it is easier to work with noncyclic groups than with cyclic ones of the same order.)

It is easy to rule out perfect 3-spanning sets of size 2 in $\mathbb{Z}_{5}^2$: indeed, if an element $a$ of the group were to be in a perfect 3-spanning set, then $2a$ and $-3a$ would need to be distinct elements; however, we have $5a=0$ for all $a \in \mathbb{Z}_{5}^2$, so $2a=-3a$, which we just ruled out.  

For our second claim, suppose, indirectly, that $A=\{a_1,a_2,a_3\}$ is a perfect 2-spanning set in $\mathbb{Z}_{5}^2$, in which case
$$\mathbb{Z}_5^2=\{0, \; \pm a_1, \; \pm a_2, \; \pm a_3, \; \pm2a_1, \; \pm2a_2, \; \pm2a_3, \; \pm a_1\pm a_2, \; \pm a_1\pm a_3,\; \pm a_2\pm a_3\}.$$  In particular, $a_1+a_2+a_3$ equals one of the elements listed.  Note that it cannot be any of the first seven; if it were, say, equal to $-a_1$, then we would get $a_2+a_3=-2a_1$, contradicting the fact that the 25 elements listed above are distinct.  Similarly, $a_1+a_2+a_3$ cannot equal any of the 12 elements among the last 18 where any of $a_1$, $a_2$, or $a_3$ appears with a positive sign, since cancelling would again result in a repetition.  This leaves only two possibilities: $a_1+a_2+a_3$ equals $-2a_i$ for some $1 \leq i \leq 3$, or it equals $-a_i-a_j$ for some $1 \leq i < j \leq 3$; without loss of generality, $$a_1+a_2+a_3=-2a_1$$ or $$a_1+a_2+a_3=-a_1-a_2.$$  But neither of these equations can occur: the first leads to
$$a_2+a_3=-3a_1=2a_1,$$ and the second yields $$a_3=-2a_1-2a_2,$$ from which we get $$2a_3=-4a_1-4a_2=a_1+a_2,$$ again contradicting the assumption that the 25 elements listed are distinct.

We will study perfect spanning sets later in more detail in Chapter \ref{ChapterSpanning}.

\newpage

\addcontentsline{toc}{section}{The declaration of independence}

\section*{The declaration of independence} 

Those familiar with linear algebra have undoubtedly heard of ``independent sets.''  A subset $A=\{a_1,\dots,a_m\}$ of elements in a vector space $V$ (over the set of real numbers) is called {\em independent}, if the zero element of $V$ cannot be expressed as a \emph{linear combination}  
$$\lambda_1a_1  + \cdots + \lambda_ma_m$$ with some (real number) coefficients $\lambda_1,  \dots, \lambda_m$ non-trivially, that is, without all coefficients $\lambda_1,  \dots, \lambda_m$ being zero.  There is an alternate definition: $A$ is independent if no element of $V$ can be expressed in the above form in two different ways.  It is a fundamental property of vector spaces that the two definitions are equivalent.

We immediately realize that the analogous concept in finite abelian groups behaves quite differently.  In fact, no (nonempty) set in a finite abelian group is independent: if $\lambda$ is the order of an element $a$ in $G$, then $\lambda a$ equals $0$, so even the 1-element set $\{a\}$ is not independent.   Thus to make our concept worth studying in finite abelian groups, we only require that a subset be independent ``to a certain degree'' rather than that it is ``completely'' independent.  More precisely, rather than considering all linear combinations of the elements of a subset, we limit our attention to those that use only a certain number of terms, that is, where the sum of the absolute values of the coefficients,
$$|\lambda_1| + \cdots + |\lambda_m|,$$ is at most some positive integer $t$.  We will refer to the linear combination above as a {\em signed sum} of $|\lambda_1| + \cdots + |\lambda_m|$ terms.

Given the two alternatives for defining independence mentioned above, we have two choices: declare a subset $A=\{a_1,\dots,a_m\}$ of an abelian group $G$ {\em $t$-independent} if
\begin{itemize}
  \item no nontrivial signed sum of at most $t$ terms equals zero (that is, the zero element of the group cannot be expressed as a linear combination 
$$\lambda_1a_1  + \cdots + \lambda_ma_m$$ using integer coefficients $\lambda_1,  \dots, \lambda_m$ with  $$|\lambda_1| + \cdots + |\lambda_m| \leq t,$$ unless $\lambda_i=0$ for every $i=1,2,\dots,m$); or
  \item no element of the group can be written as a signed sum of at most $t$ terms in two different ways (that is, we cannot have $$\lambda_1a_1+\lambda_2a_2+\cdots +\lambda_m a_m=\lambda_1'a_1+\lambda_2'a_2+\cdots +\lambda_m' a_m$$
for some integers $\lambda_1, \lambda_2, \dots , \lambda_m$ and $\lambda_1', \lambda_2', \dots , \lambda_m'$ with $$|\lambda_1|+|\lambda_2|+\cdots +|\lambda_m| \leq t \; \mbox{and} \; |\lambda_1'|+|\lambda_2'|+\cdots +|\lambda_m'| \leq t,$$ unless $\lambda_i=\lambda_i'$ for every $i=1,2,\dots,m$). 
\end{itemize}

It may come as a surprise, but the two possible definitions are not equivalent!  Consider, for example, the set $\{1,3\}$ in the cyclic group $\mathbb{Z}_{10}$: the set is 3-independent according to the first definition, since it is impossible to express zero as a non-trivial signed sum of at most three terms (try!), but the set is not 3-independent using the second definition, as, for example, $$2 \cdot 1 + (-1) \cdot 3 = 0 \cdot 1 + 3 \cdot 3.$$
So which definition should we declare to be the definition of independence?

It turns out that one definition is more powerful than the other; in fact, one of the two definitions is merely a special case of the other!  In particular, we can prove that a subset is $t$-independent by the second definition if, and only if, it is $2t$-independent by the first definition.  What this means is that studying $t$-independence for all possible $t$ using the second definition is equivalent to studying $t$-independence just for even $t$ values only using the first definition; therefore, the second definition is indeed superfluous.  

\label{declaration} To show that if a subset is $2t$-independent by the first definition then it is $t$-independent by the second definition, assume that two signed sums, each containing at most $t$ terms, equal one another.  We can then move all terms to the same side and thereby express zero as a signed sum of at most $2t$ terms.  Note that, as a result, some of the terms may get cancelled.  In fact, the assumption that our set is $2t$-independent by the first definition means that the signed sumset we created must be the trivial one, thus all terms got cancelled, which means that our two original signed sums were identical.

Conversely, if a signed sum of at most $2t$ terms equals zero, then we can rearrange this equation so that each side contains at most $t$ terms.  If we know that our set is $t$-independent by the second definition, then the two sides are identical, which means that our original signed sum was the trivial one.   

So we declare a set to be $t$-independent following the first definition; recalling our notations and terminology for the collection of signed sums of elements of a set $A$ with exactly $h$ terms being called the $h$-fold signed sumset of $A$ and denoted by $h_{\pm}A$, we can thus say that a subset $A$ is $t$-independent in a group $G$ whenever 
$$0 \not \in \cup_{h=1}^t (h_{\pm}A).$$ 
Thus, to verify that the subset $A=\{1,3\}$ of $\mathbb{Z}_{10}$ in our earlier example is 3-independent, we check that none of $1_{\pm} A$, $2_{\pm} A$, or $3_{\pm} A$ contains zero: this is obviously the case for $1_{\pm} A$ and $3_{\pm} A$ as they contain only odd elements, and $$2_{\pm} A=\{\pm (1+1), \pm (3+3), \pm 1 \pm 3\}=\{2,4,6,8\}$$ doesn't contain zero either.  Our set is not 4-independent, since $4_{\pm} A$ does contain zero (e.g. as $1+3+3+3$ or $1+1+1-3$).  It is also easy to check that $\mathbb{Z}_{10}$ has no 3-independent sets of size larger than two.

To see a bigger example, consider the set $A=\{1,4,6,9,11\}$ in $\mathbb{Z}_{25}$.  \label{Dec Indep} We find that:
\begin{itemize}
\item $1_{\pm} A = \{ 1, 4, 6, 9, 11, 14, 16, 19, 21, 24 \};$
\item $2_{\pm} A = \{ 2, 3, 5, 7, 8, 10, 12, 13, 15, 17, 18, 20, 22, 23 \};$ and
\item $3_{\pm} A = \{ 1, 2, 3, 4, 6, 7, 8, 9, 11, 12, 13, 14, 16, 17, 18, 19, 21, 22, 23, 24 \};$
\end{itemize}
and therefore $$ 0 \not \in \cup_{h=1}^3 h_{\pm} A,$$  implying that $A$ is 3-independent in $G$.  We also see, however, that $1+4+4-9$ (for example) equals zero, so
$0 \in 4_{\pm} A$, and therefore $A$ is not 4-independent.

Can we pack more than five elements in $\mathbb{Z}_{25}$ that are 3-independent?  Here is a quick argument that shows that we certainly cannot pack seven.  Suppose that $$A=\{a_1,a_2,a_3,a_4,a_5,a_6,a_7\} \subset \mathbb{Z}_{25}.$$  Now consider the following 28 signed sums:
$$0, \; \pm a_1, \; \pm a_2, \; \pm a_3, \; \pm a_4, \; \pm a_5, \; \pm a_6, \; \pm a_7,$$ $$a_1+a_1, \; a_1\pm a_2, \; a_1\pm a_3, \; a_1\pm a_4, \; a_1\pm a_5, \; a_1\pm a_6, \; a_1\pm a_7.$$  If $A$ were to be 3-independent, then all 28 expressions would be different.  (For example, if $-a_2$ were to equal $a_1-a_4$, then we would have $a_1 +a_2-a_4=0$, contradicting the fact that $A$ is 3-independent.)  Since $G$ only contains 25 elements, this cannot happen. 

It is also true that $\mathbb{Z}_{25}$ contains no six-element 3-independent sets either, but this is much harder to prove; it follows from the more general result by Bajnok and Ruzsa (see \cite{BajRuz:2003a})\index{Ruzsa, I.}\index{Bajnok, B.} that the maximum size of a $3$-independent set in the cyclic group $\mathbb{Z}_n$ equals
\begin{itemize}
  \item $\lfloor n/4 \rfloor$ when $n$ is even;
  \item $\left(1+ 1/p\right) n/6$ when $n$ is odd, has prime divisors congruent to 5 mod 6, and $p$ is the smallest such divisor; and
  \item $\lfloor n/6 \rfloor$ otherwise.
\end{itemize}
(The corresponding results for noncyclic groups or for $t>3$ are not yet known.)

Our concept of independence is one of the most important ones in this book, and is closely related to several other fundamental concepts, such as the zero-sum-free property, the sum-free property, and the Sidon property, as we now explain. 

Recall that we defined a subset $A$ of $G$ to be $t$-independent in $G$ if no signed sum of $t$ or fewer terms of $A$ (repetition of terms allowed), other than the all-zero sum, equals zero.  Observe that any equation expressing that a particular signed sum is zero can be rearranged so that on both sides of the equation the elements of $A$ appear with positive coefficients.  Therefore, an equivalent way of saying that $A$ is a $t$-independent set in $G$ is to say that, for all nonnegative integers $k$ and $l$ with $k+l \leq t$, the sum of $k$ (not necessarily distinct) elements of $A$ can only equal
the sum of $l$ (not necessarily distinct) elements of $A$ in a trivial way, that is, $k=l$ and
the two sums contain the same terms in some order.

Therefore, we can thus break up our definition of $A$ being $t$-independent into three conditions: 

\begin{itemize}

\item the {\em zero-sum-free property}: $$0 \not \in h A,$$ that is, the sum of $h$ elements of $A$ cannot equal zero---this needs to hold for $1 \leq h \leq t$; 

\item the {\em sum-free property}: $$(k  A) \cap (l  A) = \emptyset,$$ that is, the sum of $k$ elements of $A$ never equals the sum of $l$ elements of $A$---this needs to hold whenever $k$ and $l$ are distinct positive integers and $k+l \leq t$; and 

\item the {\em Sidon property}: $$|h  A | = {m+h-1 \choose h},$$ that is, two $h$-term sums of elements of $A$ can only equal if the sums contain the same terms---this needs to hold for $1 \leq h \leq \lfloor t/2 \rfloor$.   

\end{itemize}
(It is enough, in fact, to require these conditions for equations containing a total of $t$ or $t-1$ terms; therefore the total number of equations considered can be reduced to $2+(t-2)+1=t+1$.)

We can return to our earlier example of the set $A=\{1,4,6,9,11\}$ in $G=\mathbb{Z}_{25}$ to verify that $A$ is $3$-independent in $G$ using the equivalent three-part condition we just gave.  We find that
\begin{itemize}
\item $1A = \{ 1, 4, 6, 9, 11\},$
\item $2A = \{ 2, 5,7,8,10,12,13,15,17,18,20,22\},$ and
\item $3A = \{ 1, 2, 3, 4, 6, 8, 9, 11, 12, 13, 14, 16, 17, 18, 19, 21, 22, 23, 24\}.$
\end{itemize}
To conclude that $A$ is 3-independent in $G$, note that $A$ is zero-sum-free  for  $h=1,2,3$ as $0 \not \in hA$ then; sum-free  for  $(k,l)=(2,1)$ since $A$ and $2A$ are disjoint; and a Sidon set  for  $h=1$ (as are all sets).
The fact that $A$ is not 4-independent can be seen by the fact that it is not zero-$4$-sum-free ($1+4+9+11=0$), or that it is not $(3,1)$-sum-free ($1+4+4=9$), or that it is not a Sidon set for $h=2$ ($4+6=1+9$).  

We study zero-sum-free sets, sum-free sets, and Sidon sets in detail in Chapters \ref{ChapterZerosumfree}, \ref{ChapterSumfree}, and \ref{ChapterSidon}, respectively.

\part{Sides}  \label{Sides}

\setcounter{chapter}{4}
\setcounter{thm}{0}
\setcounter{section}{0}
\setcounter{subsection}{0}

Here we introduce and study some functions that prove valuable for several topics in this book.

\newpage

\addcontentsline{toc}{section}{The function $v_g(n,h)$}

\section*{The function $v_g(n,h)$}

Suppose that $h$ and $g$ are fixed positive integers.  Since we will only need the cases when $1 \leq g \leq h$, we make that assumption here.  Given a positive integer $n$, recall that $D(n)$ is the set of positive divisors of $n$.  We define  \label{vdef}
$$v_g(n,h)= \max \left\{ \left( \left \lfloor \frac{d-1-\mathrm{gcd} (d, g)}{h} \right \rfloor +1  \right) \cdot \frac{n}{d}  \mid d \in D(n) \right\}.$$ (Here we usually think of $v_g(n,h)$ as a function of $n$ for fixed values of $g$ and $h$.)  \label{0.3.5.3 page}

For example, we can compute $v_2(18,4)$ by evaluating, for each positive divisor $d$ of 18, the quantity $$\left( \left \lfloor \frac{d-1-\mathrm{gcd} (d, 2)}{4} \right \rfloor +1  \right) \cdot \frac{18}{d};$$ since the maximum occurs at $d=3$ and equals 6, we have $v_2(18,4)=6$.

The evaluation of $v_g(n,h)$ can be quite cumbersome in general.  For a prime $p$, however, $v_g(p,h)$ can be evaluated easily.  In this case, there are only two positive divisors to consider: $d=1$ and $d=p$.  For $d=1$, the expression yields 0; therefore we have 
$$v_g(p,h)= \max \left\{0,   \left \lfloor \frac{p-1-\mathrm{gcd} (p, g)}{h} \right \rfloor +1   \right\}.$$
The greatest common divisor of $p$ and $g$ is either $p$ or 1, depending on whether $p$ is a divisor of $g$ or not.  This yields the following.

\begin{prop} \label{nuforp}
If $p$ is a positive prime number, then
$$v_g(p,h) =\left\{
\begin{array}{cl}
0 & \mbox{if $p|g$,} \\ \\
\left \lfloor \frac{p-2}{h} \right \rfloor +1 & \mbox{otherwise.}\\
\end{array}\right.$$
\end{prop} 

With a bit more work, we can evaluate $v_g(n,h)$ for any prime power value of $n$, at least in the case when $g$ is not divisible by that prime; since we will use this later, we present the result here:

\begin{prop} \label{v for prime power}
If $p$ is a positive prime number that is not a divisor of $g$, then
$$v_g(p^r,h) =\left\{
\begin{array}{cl}
\frac{p^r-1}{h} & \mbox{if $p \equiv 1$ mod $h$,} \\ \\
\left(\left \lfloor \frac{p-2}{h} \right \rfloor +1 \right) \cdot p^{r-1} & \mbox{otherwise.}\\
\end{array}\right.$$
\end{prop}

We present the proof on page \pageref{proof of v for prime power}.  

Let us now turn to the evaluation of $v_g(n,h)$ when $n$ is arbitrary.  For $h=1$ and $h=2$, we can prove the following.

\begin{prop} \label{nusmall}  For all positive integers $n$ we have 
$$v_1(n,1)=n-1$$ and 
$$v_g(n,2) =\left\{
\begin{array}{ll}
\left\lfloor \frac{n}{2} \right\rfloor & \mbox{if $g=1$,} \\ \\
\left\lfloor \frac{n-1}{2} \right\rfloor & \mbox{if $g=2$.}\\
\end{array}\right.$$
\end{prop}
For a proof, see page \pageref{proofofnusmall}.

For higher values of $h$, the evaluation of $v_g(n,h)$ becomes more difficult using its definition above.  The following table lists the values of $v_1(n,3)$, $v_3(n,3)$, $v_1(n,4)$, $v_2(n,4)$, $v_4(n,4)$, $v_1(n,5)$, $v_3(n,5)$, $v_5(n,5)$ for $n \leq 40$.  (These sequences appear in \cite{OEIS} as A211316, A289435, A289436, A289437, A289438, A289439, A289440, and A289441, respectively.)

\newpage

$$\begin{array}{||c||cc||ccc||ccc||} \hline \hline
n & v_1(n,3) & v_3(n,3) & v_1(n,4) & v_2(n,4) & v_4(n,4) & v_1(n,5) & v_3(n,5) & v_5(n,5)  \\ \hline \hline
 2& 1& 1& 1& 0& 0& 1& 1& 1\\ \hline 
 3& 1& 0& 1& 1& 1& 1& 0& 1\\ \hline 
 4& 2& 2& 2& 1& 0& 2& 2& 2\\ \hline 
 5& 2& 2& 1& 1& 1& 1& 1& 0\\ \hline 
 6& 3& 3& 3& 2& 2& 3& 3& 3\\ \hline 
 7& 2& 2& 2& 2& 2& 2& 2& 2\\ \hline 
 8& 4& 4& 4& 2& 1& 4& 4& 4\\ \hline 
 9& 3& 2& 3& 3& 3& 3& 2& 3\\ \hline 
 10& 5& 5& 5& 2& 2& 5& 5& 5\\ \hline 
 11& 4& 4& 3& 3& 3& 2& 2& 2\\ \hline 
 12& 6& 6& 6& 4& 4& 6& 6& 6\\ \hline 
 13& 4& 4& 3& 3& 3& 3& 3& 3\\ \hline 
 14& 7& 7& 7& 4& 4& 7& 7& 7\\ \hline 
 15& 6& 6& 5& 5& 5& 5& 3& 5\\ \hline 
 16& 8& 8& 8& 4& 3& 8& 8& 8\\ \hline 
 17& 6& 6& 4& 4& 4& 4& 4& 4\\ \hline 
 18& 9& 9& 9& 6& 6& 9& 9& 9\\ \hline 
 19& 6& 6& 5& 5& 5& 4& 4& 4\\ \hline 
 20& 10& 10& 10& 5& 4& 10& 10& 10\\ \hline 
 21& 7& 6& 7& 7& 7& 7& 6& 7\\ \hline 
 22& 11& 11& 11& 6& 6& 11& 11& 11\\ \hline 
 23& 8& 8& 6& 6& 6& 5& 5& 5\\ \hline 
 24& 12& 12& 12& 8& 8& 12& 12& 12\\ \hline 
 25& 10& 10& 6& 6& 6& 5& 5& 4\\ \hline 
 26& 13& 13& 13& 6& 6& 13& 13& 13\\ \hline 
 27& 9& 8& 9& 9& 9& 9& 6& 9\\ \hline 
 28& 14& 14& 14& 8& 8& 14& 14& 14\\ \hline 
 29& 10& 10& 7& 7& 7& 6& 6& 6\\ \hline 
 30& 15& 15& 15& 10& 10& 15& 15& 15\\ \hline 
 31& 10& 10& 8& 8& 8& 6& 6& 6\\ \hline 
 32& 16& 16& 16& 8& 7& 16& 16& 16\\ \hline 
 33& 12& 12& 11& 11& 11& 11& 6& 11\\ \hline 
 34& 17& 17& 17& 8& 8& 17& 17& 17\\ \hline 
 35& 14& 14& 10& 10& 10& 10& 10& 10\\ \hline 
 36& 18& 18& 18& 12& 12& 18& 18& 18\\ \hline 
 37& 12& 12& 9& 9& 9& 8& 8& 8\\ \hline 
 38& 19& 19& 19& 10& 10& 19& 19& 19\\ \hline 
 39& 13& 12& 13& 13& 13& 13& 9& 13\\ \hline 
 40& 20& 20& 20& 10& 9& 20& 20& 20\\ \hline \hline
\end{array}$$

\newpage

The next theorem, based in part on work by Butterworth (cf.~\cite{But:2008a}),\index{Butterworth, J.} simplifies the evaluation of $v_g(n,h)$.

\begin{thm} [Bajnok; cf.~\cite{Baj:2017a}] \label{nu}\index{Bajnok, B.}
Suppose that $n$, $h$, and $g$ are positive integers and that $1 \leq g \leq h$.  Let $D(n)$ denote the set of positive divisors of $n$.  For $i=0,1,2,\dots,h-1$, let $$D_i(n)=\{ \; d \in D(n)  \mid d \equiv i \; (\mathrm{mod} \; h) \; \mbox{and} \; \mathrm{gcd}(d,g) <i \}.$$  We let $I$ denote those values of $i=0,1,2,\dots,h-1$ for which $D_i(n) \not = \emptyset$, and for each $i \in I$, we let $d_i$ be the smallest element of $D_i(n)$.

Then, the value of $v_g(n,h)$ is
$$v_g(n,h) =\left\{
\begin{array}{lll}
\frac{n}{h}  \cdot \mathrm{max} \left\{ 1+ \frac{h-i}{d_i}  \mid i \in I \right\} & \mbox{if} & I \not = \emptyset; \\ \\
\left\lfloor \frac{n}{h} \right\rfloor & \mbox{if} & I = \emptyset \; \mbox{and} \; g \not = h;\\
& & \\ 
\left\lfloor \frac{n-1}{h} \right\rfloor & \mbox{if} & I = \emptyset \; \mbox{and} \; g = h. \\
\end{array}\right.$$
 
\end{thm} 

The proof of Theorem \ref{nu} can be found on page \pageref{proofofnu}.  Observe that this result generalizes Proposition \ref{nusmall}.   Theorem \ref{nu} makes the computation of $v_g(n,h)$ considerably simpler; for instance, for the example of $v_2(18,4)$ above, we see that $I=\{3\}$, $d_3=3$, and thus $v_2(18,4)=18/4 \cdot (1+1/3)=6$.  

Let us now examine what explicit bounds we can deduce from Theorem \ref{nu}.  Clearly, $v_g(n,h) \geq \left\lfloor \frac{n-1}{h} \right\rfloor$.  To find an upper bound, we assume that $h \geq 2$.  First, note that for all $n$ we have $D_0(n)=D_1(n)=\emptyset$.  Furthermore, for $i \geq 2$ we have 
$$1+\frac{h-i}{d_i} \leq 1+ \frac{h-2}{2} = \frac{h}{2},$$ with equality if, and only if, $i=2$ and $d_i=2$.  This gives:

\begin{cor} \label{v function bounds}
For all integers $n$, $h \geq 2$, and $1 \leq g \leq h$, we have
$$\left\lfloor \frac{n-1}{h} \right\rfloor \leq v_g(n,h) \leq \frac{n}{2},$$ with $v_g(n,h) = \frac{n}{2}$ if, and only if, $n$ is even and $g$ is odd.
\end{cor}

We can use Theorem \ref{nu} to evaluate $v_g(n,h)$ explicitly when $h$ is relatively small.   Consider, for example, the case when $h=3$ and $g=1$.  Then $I =\{2\}$ or $I = \emptyset$ depending on whether $n$ has divisors that are congruent to 2 mod 3 or not.  Note that, when an integer has divisors that are congruent to 2 mod 3, then its smallest such divisor must be a prime, since the product of integers that are not congruent to 2 mod 3 will also not be congruent to 2 mod 3.  We thus see that 
\label{formula for v1(n,3)} $$v_1(n,3) =\left\{
\begin{array}{ll}
\left(1+\frac{1}{p}\right) \frac{n}{3} & \mbox{if $n$ has prime divisors congruent to 2 mod 3,} \\ & \mbox{and $p$ is the smallest such divisor,}\\ \\
\left\lfloor \frac{n}{3} \right\rfloor & \mbox{otherwise.}\\
\end{array}\right.$$   

We can use Theorem \ref{nu} to evaluate other cases similarly.  We find the following expressions.

$$v_2(n,3) =\left\{
\begin{array}{ll}
\left(1+\frac{1}{p}\right) \frac{n}{3} & \mbox{if $n$ has prime divisors congruent to 5 mod 6,} \\ & \mbox{and $p$ is the smallest such divisor,}\\ \\
\left\lfloor \frac{n}{3} \right\rfloor & \mbox{otherwise;}\\
\end{array}\right.$$

$$v_3(n,3) =\left\{
\begin{array}{ll}
\left(1+\frac{1}{p}\right) \frac{n}{3} & \mbox{if $n$ has prime divisors congruent to 2 mod 3,} \\ & \mbox{and $p$ is the smallest such divisor,}\\ \\
\left\lfloor \frac{n-1}{3} \right\rfloor & \mbox{otherwise;}\\
\end{array}\right.$$

\label{v_1(n,4)form} $$v_1(n,4) =\left\{
\begin{array}{ll}
\frac{n}{2}  & \mbox{if $n$ is even}\\ \\
\left(1+\frac{1}{p}\right) \frac{n}{4} & \mbox{if $n$ is odd and has prime divisors congruent to 3 mod 4,} \\ & \mbox{and $p$ is the smallest such divisor,}\\ \\
\left\lfloor \frac{n}{4} \right\rfloor & \mbox{otherwise;}\\
\end{array}\right.$$ 
\label{v_1(n,4)form}

\label{v2n4} $$v_2(n,4) =\left\{
\begin{array}{ll}
\left(1+\frac{1}{p}\right) \frac{n}{4} & \mbox{if $n$ has prime divisors congruent to 3 mod 4,} \\ & \mbox{and $p$ is the smallest such divisor,}\\ \\
\left\lfloor \frac{n}{4} \right\rfloor & \mbox{otherwise;}\\
\end{array}\right.$$

$$v_4(n,4) =\left\{
\begin{array}{ll}
\left(1+\frac{1}{p}\right) \frac{n}{4} & \mbox{if $n$ has prime divisors congruent to 3 mod 4,} \\ & \mbox{and $p$ is the smallest such divisor,}\\ \\
\left\lfloor \frac{n-1}{4} \right\rfloor & \mbox{otherwise;}\\
\end{array}\right.$$ 

$$v_1(n,5) =\left\{
\begin{array}{ll}
\mathrm{max} \left\{ \left(1+\frac{1}{p_4}\right) \frac{n}{5} , \left(1+\frac{2}{p_3}\right) \frac{n}{5}, \left(1+\frac{3}{p_2}\right) \frac{n}{5} \right\} & \mbox{if $n$ has prime divisors } \\ & \mbox{congruent to 2, 3, or 4 mod 5,} \\ & \mbox{and $p_2$, $p_3$, $p_4$ are the smallest} \\ & \mbox{such divisors, respectively,}\\ \\
\left\lfloor \frac{n}{5} \right\rfloor & \mbox{otherwise;}\\
\end{array}\right.$$ 

\label{v3n5} $$v_3(n,5) =\left\{
\begin{array}{ll}
\mathrm{max} \left\{ \left(1+\frac{1}{d_4}\right) \frac{n}{5} , \left(1+\frac{2}{p_3}\right) \frac{n}{5}, \left(1+\frac{3}{p_2}\right) \frac{n}{5} \right\} 
& \mbox{if $n$ has prime divisors congruent to} \\ 
& \mbox{2 mod 5 and $p_2$ is the smallest,} \\ 
& \mbox{or prime divisors congruent to 3 mod 5} \\
& \mbox{ other than 3 and $p_3$ is the smallest,} \\ 
& \mbox{or divisors $d_4$ congruent to 4 mod 5} \\ 
& \mbox{and $d_4$ is the smallest such divisor,}
\\ \\
\left\lfloor \frac{n}{5} \right\rfloor & \mbox{otherwise;}\\
\end{array}\right.$$

$$v_5(n,5) =\left\{
\begin{array}{ll}
\mathrm{max} \left\{ \left(1+\frac{1}{p_4}\right) \frac{n}{5} , \left(1+\frac{2}{p_3}\right) \frac{n}{5}, \left(1+\frac{3}{p_2}\right) \frac{n}{5} \right\} & \mbox{if $n$ has prime divisors } \\ & \mbox{congruent to 2, 3, or 4 mod 5,} \\ & \mbox{and $p_2$, $p_3$, $p_4$ are the smallest} \\ & \mbox{such divisors, respectively,}\\ \\
\left\lfloor \frac{n-1}{5} \right\rfloor & \mbox{otherwise.}\\
\end{array}\right.$$

Observe that in all the formulas above, with the exception of $v_3(n,5)$, the divisors playing a role can be assumed to be primes (a fact that can be proven easily).  In the case of $v_3(n,5)$, however, $d_4$ is not necessarily prime: for example, when $n=9$, we have $d_4=9$ (and no $p_2$ or $p_3$ exist), thus $v_3(9,5)=2$.  We should also point out that the maximum value in the formulae for $h=5$ above may occur with a prime that is not the smallest prime divisor of $n$ congruent to 2, 3, or 4 mod 5.  For example, for $n=437=19 \cdot 23$ we get $$v_5(437,5)=\mathrm{max} \left\{ \left(1+\frac{1}{19}\right) \frac{437}{5} , \left(1+\frac{2}{23}\right) \frac{437}{5} \right\} = \left(1+\frac{2}{23}\right) \frac{437}{5} = 95.$$

Similar expressions for $v_g(n,h)$ get more complicated for some other choices of $g$ and $h$.  It is not true, for example, that in Theorem \ref{nu} the minimal element $d_i$ of $D_i(n)$ is prime: $v_3(9,5)=2$ with $d_4=9$, and $v_6(16,11)=4$ with $d_4=4$.

The function $v_g(n,h)$ seems to possess many interesting properties; we offer the following rather vague problem.

\begin{prob} \label{nuother}
Investigate some of the properties of the function $v_g(n,h)$.
\end{prob}

\newpage

\addcontentsline{toc}{section}{The function $v_{\pm}(n,h)$}

\section*{The function $v_{\pm}(n,h)$} 

Here we introduce and investigate a close relative of the function $v_g(n,h)$ discussed above.

Unlike with $v_g(n,h)$, the function $v_{\pm}(n,h)$ we discuss here has no $g$: it depends only on positive integers $n$ and $h$.  (We can think of $v_{\pm}(n,h)$ as the ``plus--minus'' version of $v_1(n,h)$.)  We define  this function as \label{v pm def}
$$v_{\pm}(n,h)= \max \left\{ \left( 2 \cdot \left \lfloor \frac{d-2}{2h} \right \rfloor +1  \right) \cdot \frac{n}{d}  \mid d \in D(n) \right\}.$$   \label{0.3.5.3 pm page}

We can immediately see the close relationship between $v_{\pm}(n,h)$ and $v_1(n,h)$.  In particular, note that for every $d$ and $h$ we have 
$$2 \cdot \left \lfloor \frac{d-2}{2h} \right \rfloor  \leq \left \lfloor \frac{d-2}{h} \right \rfloor,$$ since if $q$ and $r$ are, respectively, the quotient and the remainder of $d-2$ when divided by $2h$, then $2 \cdot \left \lfloor (d-2)/(2h) \right \rfloor=2q$, while $\left \lfloor (d-2)/h \right \rfloor$ equals $2q$ or $2q+1$ (depending on whether $r$ is less than $h$ or not).  This yields:

\begin{prop} \label{v pm vs v}
For all positive integers $n$ and $h$ we have $v_{\pm}(n,h) \leq v_1(n,h)$.

\end{prop}

The analogue of Proposition \ref{nuforp}, this time even simpler as we have no $g$, is:

\begin{prop} \label{nu pm forp}
If $p$ is a positive prime number, then
$$v_{\pm}(p,h) =
2 \cdot \left \lfloor \frac{p-2}{2h} \right \rfloor +1 .$$
\end{prop} 

When $n$ is not (necessarily) prime, the evaluation of $v_{\pm}(n,h)$ is complicated even for small values of $h$.  For $h=1,2,3$, we can prove the following results.

\begin{prop} \label{nu pm small}  For all positive integers $n$ we have 
$$v_{\pm}(n,1)=
\left\{
\begin{array}{ll}
n-1 & \mbox{if \; $n$ is even,} \\ \\
n-2 & \mbox{if \; $n$ is odd.}\\
\end{array}\right.$$

$$v_{\pm}(n,2)= 
\left\{
\begin{array}{cl}
n/2 & \mbox{if \; $n$ is even,} \\ \\
(n-1)/2 & \mbox{if \; $n \equiv 3$  mod 4,} \\ \\
(n-3)/2 & \mbox{if \; $n \equiv 1$  mod 4.}\\
\end{array}\right.$$

$$v_{\pm}(n,3)= 
\left\{
\begin{array}{cl}
n/2 & \mbox{if \; $n$ is even,} \\ \\
n/3 & \mbox{if \; $n \equiv 3$  mod 6,} \\ \\
(n-2)/3 & \mbox{if \; $n \equiv 5$  mod 6,}\\ \\
(n-4)/3 & \mbox{if \; $n \equiv 1$  mod 6.}\\
\end{array}\right.$$
\end{prop}

For a proof, see page \pageref{proofofnu pm small}. 

We go one step further and evaluate  $v_{\pm}(n,4)$:

\begin{prop} \label{nu pm h=4}  For all positive integers $n$ we have 
$$v_{\pm}(n,4) =\left\{
\begin{array}{cl}
\frac{n}{2}  & \mbox{if $n$ is even,}\\ \\
\left(1+\frac{1}{d}\right) \frac{n}{4} & \mbox{if $n$ is odd and has divisors congruent to 3 mod 8,} \\ & \mbox{and $d$ is the smallest such divisor,}\\ \\
2 \cdot \left\lfloor \frac{n-2}{8}\right\rfloor +1 & \mbox{otherwise.}\\
\end{array}\right.$$

\end{prop}

The proof of Proposition \ref{nu pm h=4} appears on page \pageref{proofofnu pm h=4}.  We should note that the smallest divisor $d$ of $n$ that is congruent to 3 mod 8 is not necessarily a prime (consider, for example, $n=35$), unlike in the analogous formula for $v_1(n,4)$ involving the smallest divisor of $n$ that is congruent to 3 mod 4, which is always a prime.   

For higher values of $h$, the evaluation of $v_{\pm}(n,h)$ becomes increasingly difficult.  We offer:

\begin{prob}
Develop the analogue of Theorem \ref{nu} for $v_{\pm}(n,h)$. 

\end{prob}

While we don't yet have an analogue of Theorem \ref{nu} that helps us evaluate $v_{\pm}(n,h)$ exactly, we can still establish a tight upper bound for $v_{\pm}(n,h)$.  By Proposition \ref{v pm vs v} and Corollary \ref{v function bounds}, for $h \geq 2$ we get $$v_{\pm}(n,h) \leq v_1(n,h) \leq \frac{n}{2}.$$   When $n$ is even, we have $2 \in D(n)$, and thus with
$$g_d(n,h)=\left( 2 \cdot \left \lfloor \frac{d-2}{2h} \right \rfloor +1  \right) \cdot \frac{n}{d},$$ we see that
$v_{\pm}(n,h) \geq    g_2(n,2) =n/2;$ therefore, equality must hold.  When $n$ is odd and $h \geq 3$, we can show that $v_{\pm}(n,h) \leq n/3$, since $g_1(n,h)=-n$, and for odd $d \geq 3$, we have
$$g_d(n,h) \leq \left( 2 \cdot  \frac{d-3}{2h}  +1  \right) \cdot \frac{n}{d} = \left( \frac{h-3}{d} + 1  \right) \cdot \frac{n}{h} \leq \left( \frac{h-3}{3} + 1  \right) \cdot \frac{n}{h} =\frac{n}{3},$$ with equality holding if, and only if, $d=3$.  In summary, we get:

\begin{prop} \label{v pm upper bounds}
Let $n$ be a positive integer.
\begin{enumerate}
\item If $h \geq 2$, then
$v_{\pm}(n,h) \leq n/2$, with equality if, and only if, $2 | n$.
\item If $n$ is odd and $h \geq 3$, then
$v_{\pm}(n,h) \leq n/3$, with equality if, and only if, $3 | n$.
\end{enumerate}

\end{prop}
We should pint out that, in contrast, one can have $v_1(n,3)>n/3$ even for odd values of $n$ (see page \pageref{formula for v1(n,3)}).  

As with $v_g(n,h)$, we offer the following rather vague problem:

\begin{prob} \label{nu pm other}
Investigate some of the properties of the function $v_{\pm}(n,h)$.
\end{prob}

\newpage

\addcontentsline{toc}{section}{The function $u(n,m,h)$}

\section*{The function $u(n,m,h)$} \label{functionu(n,m,h)} \label{0.3.5.1}

The next function we discuss here is a variation on the famous Hopf--Stiefel function, introduced in the 1940s to study real division algebras.  Since then, the Hopf--Stiefel function has been studied in a variety of contexts, including topology, linear and bilinear algebra, and additive number theory.  

Suppose that $n$, $m$, and $h$ are fixed positive integers; we will also assume that $m \leq n$.  Recall that $D(n)$ is the set of all positive divisors of $n$.  For a fixed $d \in D(n)$, we set
$$f_d=f_d(m,h)=\left( h \cdot \left \lceil \frac{m}{d} \right \rceil -h +1  \right) \cdot d,$$
and then define \label{udef} $$u(n,m,h)=\min \left\{ f_d(m,h)  \mid d \in D(n) \right\}.$$ 

The evaluation of $u(n,m,h)$ can become quite cumbersome, particularly when $D(n)$ is relatively large.  Since for every positive integer $n$, 1 and $n$ are divisors of $n$, we have 
$$u(n,m,h) \leq  f_1(m,h)=\left( h \cdot \left \lceil \frac{m}{1} \right \rceil -h +1  \right) \cdot 1 = hm-h+1$$
and $$u(n,m,h)\leq f_n(m,h)=\left( h \cdot \left \lceil \frac{m}{n} \right \rceil -h +1  \right) \cdot n =n,$$ providing two upper bounds for $u(n,m,h)$.   
If $n$ is prime, then $D(n)=\{1,n\}$, and thus we have to only compare the two values above; thus we have the following proposition.

\begin{prop} \label{u(p,m,h)}
For a prime number $p$ we have
$$u(p,m,h)=\min \{ p, hm-h+1\}.$$

\end{prop}

We can also easily handle the cases of $m=1$ and $h=1$.  When $m=1$, we see that $\left \lceil \frac{m}{d} \right \rceil=1$ for every $d \in D(n)$, and thus $u(n,1,h)$ equals the least element of $D(n)$: 
$$u(n,1,h)=\min \left\{   d  \mid d \in D(n) \right\}=1.$$  Similarly, when $h=1$, we see that 
$$u(n,m,1)=\min \left\{ \left \lceil \frac{m}{d} \right \rceil  \cdot d  \mid d \in D(n) \right\};$$ since the expression $\left \lceil \frac{m}{d} \right \rceil  \cdot d$ attains its minimum when $d$ is a divisor of $m$, we can take $d=1$ and get
$$u(n,m,1)=m.$$

When $n$ is not prime, evaluating $u(n,m,h)$ is not easy in general.  For example, for $n=15$, we can tabulate the values of $u(15,m,h)$ for all $1 \leq m \leq n$ and for $h=2,3,4,5$: \label{u(15,m,h)}
$$\begin{array}{|c||c|c|c|c|c|c|c|c|c|c|c|c|c|c|c||} \hline
m & 1 & 2 & 3 & 4 & 5 & 6 & 7 & 8 & 9 & 10 & 11 & 12 & 13 & 14 & 15  \\ \hline \hline
u(15,m,2) & 1 & 3 & 3 & 5 & 5 & 9 & 13 & 15 & 15 & 15 & 15& 15 & 15 & 15 & 15 \\ \hline
u(15,m,3) & 1 & 3 & 3 & 5 & 5 & 12 & 15 & 15 & 15 & 15 & 15& 15 & 15 & 15 & 15 \\ \hline
u(15,m,4) & 1 & 3 & 3 & 5 & 5 & 15 & 15 & 15 & 15 & 15 & 15& 15 & 15 & 15 & 15 \\ \hline
u(15,m,5) & 1 & 3 & 3 & 5 & 5 & 15 & 15 & 15 & 15 & 15 & 15& 15 & 15 & 15 & 15 \\ \hline
\end{array}$$  

These values seem to indicate that we always have $m \leq u(n,m,h) \leq n.$  The upper bound we have already explained above; in fact, we have seen that $$u(n,m,h) \leq \min \{ n, hm-h+1\}.$$ We will prove that the lower bound holds, as follows.  We have seen that $u(n,m,1)=m$, so let us assume that $h \geq 2$.  Consider the expression 
$$f_d(m,h)= \left( h \cdot \left \lceil \frac{m}{d} \right \rceil -h +1  \right) \cdot d.$$ 
Note that $f_m(m,h)=m$.  What we will show is that $f_d(m,h)$ is greater than $m$ both when $d>m$ and when $d < m$.

If $d>m$, then $\left \lceil \frac{m}{d} \right \rceil =1$, and thus $f_d(m,h)=d$; since we assumed that $d >m$, we have $f_d(m,h) > m$.

If $d < m$, we see that
\begin{eqnarray*}
f_d(m,h) & = & \left( h \cdot \left \lceil \frac{m}{d} \right \rceil -h +1  \right) \cdot d \\
& \geq & \left( h  \cdot \frac{m}{d} -h +1  \right) \cdot d \\
& = & (h-1)(m-d) +m;
\end{eqnarray*}  
a value strictly greater than $m$ under the assumptions that $h \geq 2$ and $d < m$.

Therefore, we have the following result.

\begin{prop} \label{u(n,m,h)extreme}
For all positive integers $n$, $m$, and $h$, with $m \leq n$, we have
$$m \leq u(n,m,h) \leq \min \{ n, hm-h+1\};$$
furthermore, $u(n,m,h)=m$ if, and only if, $h=1$ or $n$ is divisible by $m$.

\end{prop}

In certain cases, it is useful to compare $u(n,m,h)$ to the smallest prime divisor $p$ of $n$.  
Suppose first that $m \leq p$.  In this case, we have 
$$u(n,m,h) = \min \{p, hm-h+1\},$$ since for $d=1$ we have $f_d(m,h)=hm-h+1$, and for all other divisors $d$ of $n$, we have $d \geq p \geq m$ and, therefore, 
$$f_d(m,h)=\left( h \cdot \left \lceil \frac{m}{d} \right \rceil -h +1  \right) \cdot d=d,$$ which attains its minimum when $d=p$.  On the other hand, when $p<m$, then $p<u(n,m,h)$ since $u(n,m,h) \geq m$, and, therefore,  
$$u(n,m,h) > \min \{p, hm-h+1\}.$$
In summary, we have shown the following.

\begin{prop} \label{u versus p}
Let $n$, $m$, and $h$ be positive integers, $m \leq n$, and let $p$ be the smallest prime divisor of $n$.  We then have 
$$u(n,m,h) \geq \min \{p, hm-h+1\},$$
with equality if, and only if, $m \leq p$.
\end{prop} 

Note that Proposition \ref{u versus p} is a generalization of Proposition \ref{u(p,m,h)}.

Proposition \ref{u(n,m,h)extreme} has a complete characterization for situations when $u(n,m,h)$ reaches its lower bound, $m$---but how about its upper bound?  In particular, when is $u(n,m,h)=n$?  From its definition, we can observe that the function $u(n,m,h)$ is a nondecreasing function of $m$ and also of $h$ (but not of $n$).  Furthermore, for fixed $n$ and $h \geq 1$, the values of the function as $m$ increases reach $n$; indeed, noting that if $d \in D(n)$ then $\left \lceil \frac{n}{d} \right \rceil=\frac{n}{d}$, we have
$$u(n,n,h)=\min \left\{ \left( h \cdot \left \lceil \frac{n}{d} \right \rceil -h +1  \right) \cdot d \mid d \in D(n) \right\}= 
\min \{  h n -d(h-1) \mid d \in D(n)\}= n.$$  

We call the minimum value of $m$ for which $u(n,m,h)=n$ the {\em $h$-critical number of $n$}; according to our computation we see that the $h$-critical number of $n$ is well-defined for every $h \geq 1$, and equals at most $n$.  This value is now known, as we can prove the following result:

\begin{thm}  \label{h-critical of n}
The $h$-critical number of $n$ equals $v_1(n,h)+1$.  

\end{thm}

The short proof of Theorem \ref{h-critical of n} can be found on page \pageref{proof of h-critical of n}.  (Campbell\index{Campbell, K.} in \cite{Cam:2014a} established Theorem \ref{h-critical of n} in some special cases, and provided lower bounds in other cases.)  According to Theorem \ref{h-critical of n}, the arithmetic function $v_1(n,h)$ is a certain inverse of the Hopf--Stiefel function.  In particular, the maximum value of $m$ for which $u(n,m,h)<n$ equals $v_1(n,h)$.

While the function $u(n,m,h)$ seems to be quite mysterious, it may be possible to find its values, at least in certain special cases, more easily, or to determine upper and lower bounds better than those above.  We pose the following rather vague problem.

\begin{prob}
Find some exact formulas or good upper or lower bounds, at least in certain special cases, for $u(n,m,h)$.
\end{prob}

\newpage

\addcontentsline{toc}{section}{The function $u \hat{\;}(n,m,h)$}

\section*{The function $u \hat{\;}(n,m,h)$} \label{functionuhat(n,m,h)} \label{0.3.5.2}

The function $u \hat{\;}(n,m,h)$ we introduce here is a relative of the Hopf--Stiefel function $u(n,m,h)$ discussed above.  (These two functions will express cardinalities of two related sets we study later.)

Suppose that $n$, $m$, and $h$ are fixed positive integers; we will also assume that $h< m \leq n$.  For a fixed $d \in D(n)$, we let $k$ and $r$ denote the positive remainder of $m$ mod $d$ and $h$ mod $d$, respectively.  That is, we write
$$m=cd+k \; \; \; \mbox{and} \; \; \; h=qd+r$$
with
$$1 \leq k \leq d \; \; \; \mbox{and} \; \; \; 1 \leq r \leq d.$$  Note that one can compute these quotients and remainders as 
$$c=\left \lceil \frac{m}{d} \right \rceil -1 \; \; \; \mbox{and} \; \; \; q= \left \lceil \frac{h}{d} \right \rceil -1;$$
$$k=m-d \left \lceil \frac{m}{d} \right \rceil +d \; \; \; \mbox{and} \; \; \; r=h-d \left \lceil \frac{h}{d} \right \rceil +d.$$

We then set
$$f\hat{_d}(n,m,h)=
\left\{\begin{array}{ll}
\min\{n, f_d, hm-h^2+1\}  & \mbox{if $h \leq \min\{k,d-1\}$,}\\ \\
\min\{n, hm-h^2+1 - \delta_d\} & \mbox{otherwise;}\\  
\end{array}\right.$$
where $f_d$ is the function $$f_d(m,h)=\left( h \cdot \left \lceil \frac{m}{d} \right \rceil -h +1  \right) \cdot d$$ defined in Section \ref{0.3.5.1}, and $\delta_d$ is a ``correction term'' defined as
$$\delta_d(n,m,h)=\left\{ 
\begin{array}{cl}
(k-r)r-(d-1)  & \mbox{if $r<k$,}\\ 
(d-r)(r-k)-(d-1)  & \mbox{if $k<r<d$,}\\
d-1  & \mbox{if $k=r=d$,}\\
0  & \mbox{otherwise.}  
\end{array}\right.$$
We then define  \label{uhatdef} $$u\hat{\;}(n,m,h)=\min \{ f\hat{_d}(n,m,h)  \mid d \in D(n) \}.$$ 

Evaluating $u\hat{\;}(n,m,h)$ is more complicated than evaluating $u(n,m,h)$ is; to see some specific values, we compute $u\hat{\;}(n,m,h)$ for $n=15$ and for all $h < m \leq n$ and $h=2,3,4,5$: \label{uhat(15,m,h)}
$$\begin{array}{|c||c|c|c|c|c|c|c|c|c|c|c|c|c|c|c||} \hline
m & 1 & 2 & 3 & 4 & 5 & 6 & 7 & 8 & 9 & 10 & 11 & 12 & 13 & 14 & 15  \\ \hline \hline
u\hat{\;}(15,m,2) &  &  & 3 & 5 & 5 & 9 & 11 & 13 & 15 & 15 & 15& 15 & 15 & 15 & 15 \\ \hline
u\hat{\;}(15,m,3) &  &  &  & 4 & 5 & 8 & 13 & 15 & 15 & 15 & 15& 15 & 15 & 15 & 15 \\ \hline
u\hat{\;}(15,m,4) &  &  &  &  & 5 & 9 & 13 & 15 & 15 & 15 & 15& 15 & 15 & 15 & 15 \\ \hline
u\hat{\;}(15,m,5) &  &  &  &  &  & 6 & 11 & 15 & 15 & 15 & 15& 15 & 15 & 15 & 15 \\ \hline
\end{array}$$

We pose the following rather vague problem.

\begin{prob} \label{redefineuhat}
Find a simpler way to define $u\hat{\;}(n,m,h)$.
\end{prob} 

Somewhat less ambitiously:

\begin{prob} \label{evaluateuhat}
Find a general expression for $u\hat{\;}(n,m,h)$, at least for some cases.
\end{prob}

As perhaps a step toward Problems \ref{redefineuhat} and \ref{evaluateuhat}, there may be a simplified expression for $\delta_d$.  Indeed, we can easily compute $\delta_d$ for small values of $d$.  We get
$$\delta_1 = 0 \; \; \mbox{for all $m$ and $h$}; $$
$$\delta_2 = \left\{\begin{array}{rl}
1  & \mbox{if $m$ and $h$ are both even,}\\
0 & \mbox{otherwise;} 
\end{array}\right.$$
and
$$\delta_3 = \left\{\begin{array}{rl}
2  & \mbox{if $m$ and $h$ are both divisible by 3,}\\
-1  & \mbox{if none of $m$, $h$, or $m-h$ is divisible by 3,}\\
0 & \mbox{otherwise.} 
\end{array}\right.$$
As these formulae suggest, one may be able to evaluate $\delta_d$ easily.

\begin{prob} \label{evaluatedelta}
Find a simpler general expression for $\delta_d(n,m,h)$.
\end{prob}

With a simpler general expression for $\delta_d$, we can hope for an easier way to evaluate $f\hat{_d}$.  For $d \leq 3$ we find that \label{fhatsmall}
$$f\hat{_1}=\min\{n,hm-h^2+1\};$$

$$f\hat{_2}= \left\{\begin{array}{ll}
\min\{n,hm-h^2\}  & \mbox{if $m$ and $h$ are both even,}\\ \\
\min\{n,hm-h^2+1\} & \mbox{otherwise;} 
\end{array}\right.$$
and

$$f\hat{_3}=\left\{\begin{array}{ll}
\min\{n,2m-1,3m-8\}  & \mbox{if $h=2$ and $3|m-2$,}\\ \\
\min\{n,2m-3,3m-8\}  & \mbox{if $h=2$ and $3|m$,}\\ \\
\min\{n,hm-h^2-1\}  & \mbox{if $3|m$ and $3|h$,}\\ \\
\min\{n,hm-h^2+2\}  & \mbox{if $h>1$, $3\not|m$, $3\not|h$, and $3 \not| m-h$},\\ \\
\min\{n,hm-h^2+1\} & \mbox{otherwise.} 
\end{array}\right.$$  

Having explicit formulae for $f\hat{_d}$ for all $n,m,h,$ and $d$ would enable us to evaluate $u\hat{\;}(n,m,h)$.  Regarding the case of $h=1$: we see that either $d \geq 2$ in which case 
$$f\hat{_d}(n,m,1)=\min\{n, f_d, hm-h^2+1\}=\min\{n, \left \lceil \frac{m}{d} \right \rceil \cdot d , m\}=m,$$
or $d=1$ in which case
$$f\hat{_1}(n,m,1)=\min\{n, hm-h^2+1\}=\min\{n, m\}=m;$$
and therefore 
$$u\hat{\;}(n,m,1)=m.$$

For $h=2$ and $h=3$ we have the following results.

\begin{prop}  \label{uhatforh=2}

For $h=2$, we have
$$u\hat{\;}(n,m,2)=\left\{\begin{array}{ll}
\min\{u(n,m,2),2m-4\}  & \mbox{if $n$ and $m$ are both even,}\\ \\
\min\{u(n,m,2),2m-3\} & \mbox{otherwise.} 
\end{array}\right.$$

\end{prop}

\begin{prop}  \label{uhatforh=3}

For $h=3$, we have
$$u\hat{\;}(n,m,3)=\left\{
\begin{array}{ll}
\min\{u(n,m,3), 3m-3-\gcd(n,m-1)\} & \mbox{if} \; \gcd(n,m-1) \geq 8; \\ \\
\min\{u(n,m,3), 3m-10\} & \mbox{if} \; \gcd(n,m-1) = 7, \; \mbox{or} \\
& \gcd(n,m-1) \leq 5, \; 3|n, \; \mbox{and} \; 3|m; \\ \\
\min\{u(n,m,3), 3m-9\} & \mbox{if} \; \gcd(n,m-1) = 6; \\ \\
\min\{u(n,m,3), 3m-8\} & \mbox{otherwise.} 
\end{array}
\right.$$

\end{prop}

The proof of Propositions \ref{uhatforh=2} and \ref{uhatforh=3} can be found on pages \pageref{proofofuhatforh=2} and \pageref{proofofuhatforh=3}, respectively.  As $h$ increases, the evaluation of $u\hat{\;}(n,m,h)$ gets more complicated; we offer the following problem.

\begin{prob}
Find an explicit formula for $u\hat{\;}(n,m,h)$ in terms of $u(n,m,h)$ for $h=4$, $h=5$, and $h=6$.
\end{prob}

We also know the value of $u\hat{\;}(n,m,h)$ when $n$ is prime; in that case $D(n)=\{1,n\}$.  We have already seen that
$$f\hat{_1}=\min\{n,hm-h^2+1\};$$ we can similarly verify that
$$f\hat{_n}=\min\{n,hm-h^2+1\}$$and thus we have the following proposition.

\begin{prop} \label{uhat(p,m,h)}
For a prime number $p$ we have
$$u\hat{\;}(p,m,h)=\min \{ p, hm-h^2+1\}.$$

\end{prop}

Other than the results above, we have no general formula for $u\hat{\;}(n,m,h)$.  Below we investigate some upper and lower bounds instead.  First, we compare $u\hat{\;}(n,m,h)$ and $u(n,m,h)$.

\begin{prop} \label{uhat lessthan u}
For all $h < m \leq n$ we have $u\hat{\;}(n,m,h) \leq u(n,m,h)$.

\end{prop}

The proof of Proposition \ref{uhat lessthan u} can be found on page \pageref{proofofuhat lessthan u}.  

Let us now see how much less $u\hat{\;}(n,m,h)$ can get below $u(n,m,h)$.   For $h=1$, we have already seen that $u(n,m,1)=u\hat{\;}(n,m,1)=m$.  For $h=2$, we use Proposition \ref{uhatforh=2} to see that, when $u\hat{\;}(n,m,2)<u(n,m,2)$ and $n$ and $m$ are both even, then
$$1 \leq u(n,m,2)-u\hat{\;}(n,m,2) = u(n,m,2)-(2m-4) \leq f_2(m,2) - (2m-4) = (2m -2)- (2m-4)=2,$$  and when $u\hat{\;}(n,m,2)<u(n,m,2)$ and $n$ or $m$ is odd, then
$$1 \leq u(n,m,2)-u\hat{\;}(n,m,2) = u(n,m,2)-(2m-3) \leq f_1(m,2) - (2m-3) = (2m -1)- (2m-3)=2.$$  

This proves the following.

\begin{prop} \label{uhat lessthan u h=2}
For all $m \leq n$ we have $$u(n,m,2) - 2 \leq u\hat{\;}(n,m,2) \leq u(n,m,2).$$

\end{prop}
According to Proposition \ref{uhat lessthan u h=2}, the only possible values of $u\hat{\;}(n,m,2)$ are $u(n,m,2)$, $u(n,m,2)-1$, and $u(n,m,2)-2$.

In contrast to Proposition \ref{uhat lessthan u h=2}, we can show that, when $h \geq 3$, then $u(n,m,h) - u\hat{\;}(n,m,h)$ can get arbitrarily large: larger than any given (positive) real number $C$.  As an example, take an arbitrary prime $p>h$ and a positive integer $t >2$ so that $$(h-2)p^{t-1} \geq C.$$
One can readily verify that for $n=p^t$ and $m=p^{t-1}+1$ we have 
$$u(n,m,h)=f_1=hp^{t-1}+1$$ and for $d=p^{t-1}$ we get 
$$u \hat{\;} (n,m,h) \leq f \hat{_{d}}=2p^{t-1}.$$  Thus, we have:

\begin{prop} \label{rho-hat very small}
For every $h \geq 3$ and for every positive real number $C$, one can find positive integers $n$ and $m$ for which $$u\hat{\;}(n,m,h)< u(n,m,h)- C.$$
\end{prop}

Recall that by Proposition \ref{u(n,m,h)extreme}, we have
$$m \leq u(n,m,h) \leq \min \{ n, hm-h+1\},$$ with $u(n,m,h)=m$ if, and only if, $h=1$ or $n$ is divisible by $m$.  For $u\hat{\;}(n,m,h)$ we have the following result.

\begin{prop} \label{uhat(n,m,h)extreme}
For all positive integers $n$, $m$, and $h$ with $h< m \leq n$ we have
$$m \leq u\hat{\;}(n,m,h) \leq  \min \{ u(n,m,h), hm-h^2+1\}\leq \min \{ n, hm-h^2+1\},$$
with $u\hat{\;}(n,m,h)=m$ if, and only if, (at least) one of the following holds:
\begin{enumerate}[(i)]
\item $h=1$,
\item $h=m-1$,
\item $n$ is divisible by $m$, or
\item $h=2$, $m=4$, and $n$ is even.
\end{enumerate}
\end{prop}

The upper bound in Proposition \ref{uhat(n,m,h)extreme} follows from the fact that 
$$u\hat{\;}(n,m,h) \leq f\hat{_1}(n,m,h)=\min\{n,hm-h^2+1\}$$ and Propositions \ref{uhat lessthan u} and \ref{u(n,m,h)extreme}.  The proof of the lower bound and the classification of its equality can be found on page \pageref{proofofuhat(n,m,h)extreme}.

We know considerably less about cases when $u\hat{\;}(n,m,h)$ reaches the upper bound $$\min \{ u(n,m,h), hm-h^2+1\}$$ in Proposition \ref{uhat(n,m,h)extreme}.  Equality clearly always holds when $h=1$, and from Proposition \ref{uhatforh=2} we see that, when $h=2$, equality also holds unless $n$ and $m$ are both even and $2m-3 \leq u(n,m,2)$.  When comparing the table on page \pageref{u(15,m,h)} to the table on page \pageref{uhat(15,m,h)}, we see that, for $n=15$ and $2 \leq h \leq 5$, there is only one case when the upper bound is not reached: $h=3$ and $m=6$.  We may observe that in both these cases we have $\mathrm{gcd}(n,m,h) >1$;  indeed, when $\mathrm{gcd}(n,m,h)=d>1$ and $hm-h^2+1 \leq u(n,m,h)$, then 
$$u\hat{\;}(n,m,h) \leq f\hat{_d}(n,m,h) =hm-h^2+1 -(d-1) < \mathrm{min}\{u(n,m,h), hm-h^2+1\}.$$     
It is not hard to find other scenarios when the upper bound is not reached, but we do not have a full understanding of all of them.  We offer the following challenging problem.

\begin{prob} \label{when uhat <}
Classify all situations when $u\hat{\;}(n,m,h) < \mathrm{min}\{u(n,m,h), hm-h^2+1\}.$

\end{prob}

As we did for $u(n,m,h)$, we now compare $u\hat{\;}(n,m,h)$ to the smallest prime divisor $p$ of $n$.  Clearly, if $p<m$ then, by Proposition \ref{uhat(n,m,h)extreme}, 
$$u\hat{\;}(n,m,h) \geq m >p \geq \min \{p, hm-h^2+1\}.$$
On the other hand, if $p \geq m$, then 
$$f\hat{_1}(n,m,h)=\min\{n,hm-h^2+1\} \geq \min\{p,hm-h^2+1\},$$
 and for $d \in D(n) \setminus \{1\}$, we have $d \geq p \geq m>h$ and, therefore, $r=h<m=k$, so $h \leq\min\{d-1,k\}$ and thus
\begin{eqnarray*}
f\hat{_d}(n,m,h)&=&\min\{n,\left( h \cdot \left \lceil \frac{m}{d} \right \rceil -h +1  \right) \cdot d,hm-h^2+1\} \\
&=&\min\{n,d,hm-h^2+1\} \\
&\geq& \min\{p,hm-h^2+1\},
\end{eqnarray*} with equality when (though not necessarily only when) $d=p$, and thus $$u\hat{\;}(n,m,h)=\min\{p,hm-h^2+1\}.$$
In summary, we have shown the following.

\begin{prop} \label{uhat versus p}
Let $n$, $m$, and $h$ be positive integers, $h< m \leq n$, and let $p$ be the smallest prime divisor of $n$.  We then have 
$$u\hat{\;}(n,m,h) \geq \min \{p, hm-h^2+1\},$$
with equality if, and only if, $m \leq p$.
\end{prop} 

Note that Proposition \ref{uhat versus p} is a generalization of Proposition \ref{uhat(p,m,h)}.

Recall that the function $u(n,m,h)$ is a nondecreasing function of both $m$ and $h$ (but not of $n$) and, for a fixed $h$, it reaches $n$ at a certain threshold  value of $m$, called the $h$-critical number of $n$.  In contrast, as the table on page \pageref{uhat(15,m,h)} indicates, $u\hat{\;}(n,m,h)$ is not a nondecreasing function of $h$.  With some effort, one can prove that $u\hat{\;}(n,m,h)$ is a nondecreasing function of $m$ and that it reaches $n$ eventually, so we are able to define the {\em restricted  $h$-critical number of $n$}.  We will only do this here for $h \leq 2$.

Since $u\hat{\;}(n,m,1)=m$ is nondecreasing with $m$ reaching $n$ at $m=n$, the restricted  $1$-critical number of $n$ is clearly well-defined and equals $n$.  For $h=2$, we use Proposition \ref{uhatforh=2} to show that $u\hat{\;}(n,m+1,2) \geq u\hat{\;}(n,m,2)$. Indeed, if $n$ is odd, then
$$u\hat{\;}(n,m+1,2)=\min\{u(n,m+1,2),2(m+1)-3\} \geq \min\{u(n,m,2),2m-3\}= u\hat{\;}(n,m,2);$$ if $n$ is even and $m$ is even, then
$$u\hat{\;}(n,m+1,2)=\min\{u(n,m+1,2),2(m+1)-3\} \geq \min\{u(n,m,2),2m-4\}= u\hat{\;}(n,m,2);$$
and, if $n$ is even and $m$ is odd, then
$$u\hat{\;}(n,m+1,2)=\min\{u(n,m+1,2),2(m+1)-4\} \geq \min\{u(n,m,2),2m-3\}= u\hat{\;}(n,m,2).$$  Therefore, $u\hat{\;}(n,m,2)$ is a nondecreasing function of $m$.  
Furthermore, using Theorem \ref{h-critical of n} and Propositions \ref{nusmall} and \ref{uhatforh=2}, we also see that 
$$u\hat{\;}(n,m,2)\left\{
\begin{array}{ll}
< n & \mbox{if $m \leq \left \lfloor \frac{n}{2} \right \rfloor +1$,} \\ 
= n & \mbox{if $m \geq \left \lfloor \frac{n}{2} \right \rfloor +2$,} 
\end{array}\right.$$
and thus we have the following.

\begin{thm} \label{restricted critical of n h=2}
The restricted  $2$-critical number of $n$ is $\left \lfloor \frac{n}{2} \right \rfloor +2$.
\end{thm}

\begin{prob}
For each $h \geq 3$, prove that the $h$-critical number of $n$ is well-defined and find its value.
\end{prob}

The arithmetic function $u\hat{\;}(n,m,h)$ seems quite interesting; we offer the following vague problem.

\begin{prob}
Find some exact formulas or good upper or lower bounds, at least in certain special cases, for $u\hat{\;}(n,m,h)$.
\end{prob}

\part{Entrees}  \label{Entrees}

\setcounter{chapter}{0}
\renewcommand{\thechapter}{\Alph{chapter}}

This is the main part of the book where our problems are posed and their backgrounds are explained.  At the present time, these problems are about maximum sumset size, spanning sets, Sidon sets, minimum sumset size,  critical numbers, zero-sum-free sets, and sum-free sets.  (We plan to include additional chapters in the near future.)  

Guide to section numbering: In Chapter X (where X $\in \{$A, B, C, \ldots$\}$), we study, for $A \subseteq G$, $H \subseteq \mathbb{N}_0$, $h \in \mathbb{N}_0$, and $s \in \mathbb{N}_0$: 

Section X.1: Unrestricted sumsets: $HA$ 

\hspace{.6in} Section X.1.1: Fixed number of terms: $hA$

\hspace{.6in} Section X.1.2: Limited number of terms: $[0,s]A$

\hspace{.6in} Section X.1.3: Arbitrary number of terms: $\langle A \rangle$

Section X.2: Unrestricted signed sumsets: $H_{\pm}A$

\hspace{.6in} Section X.2.1: Fixed number of terms: $h_{\pm}A$

\hspace{.6in} Section X.2.2: Limited number of terms: $[0,s]_{\pm}A$

\hspace{.6in} Section X.2.3: Arbitrary number of terms: $\langle A \rangle$

Section X.3: Restricted sumsets: $H\hat{\;}A$

\hspace{.6in} Section X.3.1: Fixed number of terms: $h\hat{\;}A$ 

\hspace{.6in} Section X.3.2: Limited number of terms: $[0,s] \hat{\;} A$ 

\hspace{.6in} Section X.3.3: Arbitrary number of terms: $\Sigma A$ 

Section X.4: Restricted signed sumsets: $H\hat{_{\pm}}A$

\hspace{.6in} Section X.4.1: Fixed number of terms: $h\hat{_{\pm}}A$

\hspace{.6in} Section X.4.2: Limited number of terms: $[0,s]\hat{_{\pm}}A$

\hspace{.6in} Section X.4.3: Arbitrary number of terms: $\Sigma_{\pm}A$

\noindent See the table on page \pageref{twelvesumsets} for definitions and terminology.

\chapter{Maximum sumset size}  \label{ChapterMaxsumsetsize}

Recall that for a given finite abelian group $G$, $m$-subset $A=\{a_1,\dots, a_m\}$ of $G$, $\Lambda \subseteq \mathbb{Z}$, and $H \subseteq \mathbb{N}_0$, we defined the sumset of $A$ corresponding to $\Lambda$ and $H$ as
$$H_{\Lambda}A = \{\lambda_1a_1+\cdots +\lambda_m a_m \mbox{    } |  \mbox{    }  (\lambda_1,\dots ,\lambda_m) \in \Lambda^m(H)  \}$$
where the index set $\Lambda^m(H)$ is defined as
$$\Lambda^m(H)=\{(\lambda_1,\dots ,\lambda_m) \in \Lambda^m \; |  \; |\lambda_1|+\cdots +|\lambda_m| \in H \}.$$ 

In this chapter we attempt to answer the following question: Given a finite abelian group $G$ and a positive integer $m$, how large can a sumset of an $m$-subset of $G$ be?  More precisely, our objective is to determine, for any $G$, $m$, $\Lambda$, and $H$ the quantity
$$\nu_{\Lambda}(G,m,H)=\mathrm{max} \{ |H_{\Lambda}A|  \mid A \subseteq G, |A|=m\}.$$

We can immediately see two upper bounds for $\nu_{\Lambda}(G,m,H)$: since $H_{\Lambda}A$ is a subset of $G$ and since each element of the index set $\Lambda^m(H)$ contributes a unique element toward $H_{\Lambda}A$ (which may or may not be distinct), we have the following obvious result.

\begin{prop} \label{nuupper}  We have
$$\nu_{\Lambda}(G,m,H) \leq \mathrm{min} \{|G|, |\Lambda^m(H)| \}.$$
\end{prop}

One of the most fascinating questions in additive combinatorics is to investigate the cases when equality occurs in Proposition \ref{nuupper}.  In particular, subsets $A$ of $G$ with $H_{\Lambda}A=G$ (that is, every element of $G$ can be written as an appropriate sum of the elements of $A$) are called spanning sets, and subsets with $|H_{\Lambda}A|=|\Lambda^m(H)|$ (i.e. all sums are distinct modulo the rearrangement of the terms) are called Sidon sets.  We investigate spanning sets in Chapter \ref{ChapterSpanning} and Sidon sets in Chapter \ref{ChapterSidon}.

In the following sections we consider $\nu_{\Lambda}(G,m,H)$ for special $\Lambda \subseteq \mathbb{Z}$ and $H \subseteq \mathbb{N}_0$.

\section{Unrestricted sumsets} \label{sectionmaxsumsetsizeU}

Our goal in this section is to investigate the quantity $$\nu (G,m,H) =\mathrm{max} \{ |HA|  \mid A \subseteq G, |A|=m\}$$ where $HA$ is the union of all $h$-fold sumsets $hA$ for $h \in H$.  We consider three special cases: when $H$ consists of a single nonnegative integer $h$, when $H$ consists of all nonnegative integers up to some value $s$, and when $H$ is the entire set of nonnegative integers. 

\subsection{Fixed number of terms} \label{sectionmaxsumsetsizeUfixed}

We first consider $$\nu (G,m,h) = \mathrm{max} \{ |hA|  \mid A \subseteq G, |A|=m\},$$ that is, the maximum size of an $h$-fold sumset of an $m$-element subset of $G$.

In general, we do not know the value of $\nu (G,m,h)$.  Obviously, for $h=0$ and $h=1$ we have $hA=\{0\}$ and $hA=A$, respectively, thus $$\nu (G,m,0)= 1$$ and $$\nu (G,m,1)=  m$$ for every $m \in \mathbb{N}$.  If $A$ consist of a single element $a$, then $hA=\{ha\}$ for every $h \in \mathbb{N}_0$, so $$\nu (G,1,h)= 1$$ for every $h \in \mathbb{N}_0$.  

Before turning to the case of $m=2$, recall that, according to Proposition \ref{shifting sets}, in order to find $\nu (G,m,h)$ (for any $G$, $m$, and $h$), we may assume that the subset $A$ of $G$ of size $m$ yielding $|hA|=\nu (G,m,h)$ will contain $0$.  Therefore, to find $\nu (G,2,h)$, we only need to look for an element $a \in G$ for which the size of $$h\{0,a\}=\{0,a,2a,\dots,ha\}$$ is maximal.  But the size of this set is clearly $\mathrm{min} \{|\langle a \rangle|,h+1\}$.  Here $|\langle a \rangle|$ is the order of $a$; its maximal value is the exponent of $G$, denoted by $\kappa$.  Therefore, $$\nu (G,2,h)= \mathrm{min} \{\kappa, h+1\}$$ for any $h \in \mathbb{N}_0$.

Summarizing our results thus far, we have the following.

\begin{prop} In any abelian group $G$ of order $n$ and exponent $\kappa$ we have
\begin{eqnarray*}
\nu(G,m,0) & = & 1, \\
\nu(G,m,1) & = & m, \\
\nu(G,1,h) & = & 1, \\
\nu(G,2,h) & = & \mathrm{min} \{\kappa, h+1\}.
\end{eqnarray*}
\end{prop} 

For values of $m \geq 3$ and $h \geq 2$, we have no exact values for $\nu (G,m,h)$ in general; as a consequence of Proposition \ref{nuupper}, we have, however, the following upper bound.

\begin{prop} \label{nuupperUfixed}  In any abelian group $G$ of order $n$ we have
$$\nu(G,m,h) \leq \mathrm{min} \left\{ n, {m+h-1 \choose h} \right\}.$$
\end{prop}
       
We pose the following very general problem.

\begin{prob}  \label{findnugeneralgroup}
Find the value of (or, at least, find good bounds for) $\nu(G,m,h)$ for noncyclic groups $G$ and integers $m$ and $h$.
\end{prob}

The case of cyclic groups is of special interest:

\begin{prob} \label{findnugeneralcyclic}
For positive integers $n$, $m$, and $h$, find the value of (or, at least, find good bounds for) $\nu(\mathbb{Z}_{n},m,h)$.
\end{prob}

The value of $\nu(\mathbb{Z}_{n},m,h)$ (or, more generally, the value of $\nu(G,m,h)$) tends to agree with the upper bound in Proposition \ref{nuupperUfixed}.  Perhaps the most intriguing aspect of Problem \ref{findnugeneralcyclic} (or, more generally, of Problem \ref{findnugeneralgroup}) is to analyze the exceptions to this predilection.  The number of exceptions seems to vary a great deal.  For example, Manandhar\index{Manandhar, D.} in \cite{Man:2012a} found (using a computer program) that, as $n$ ranges from 2 to 20, the number of such exceptions is as follows.
$$\begin{array}{|c||c|c|c|c|c|c|c|c|c|c|c|c|c|c|c|c|c|c|c|} \hline
 n & 2 & 3 & 4 & 5 & 6 & 7 & 8 & 9 & 10 & 11 & 12 & 13 & 14 & 15 & 16 & 17 & 18 & 19 & 20 \\ \hline \hline
\mbox{exceptions} & 0 & 0 & 0 & 0 & 1 & 0 & 0 & 1 & 2 & 2 & 3 & 1 & 2 & 2 & 2 & 4 & 5 & 4 & 6 \\ \hline 
\end{array}$$
The six exceptions for $n=20$ are listed in the following table.
$$\begin{array}{|c|c||c|c|} \hline
 h & m  & \mathrm{min} \left\{ 20, {m+h-1 \choose h} \right\}  & \nu(\mathbb{Z}_{20},m,h) \\ \hline \hline

2 & 5 & 15 & 14 \\ \hline
2 & 6 & 20 & 18 \\ \hline
3 & 4 & 20 & 18 \\ \hline
4 & 3 & 15 & 14 \\ \hline
5 & 3 & 20 & 17 \\ \hline
6 & 3 & 20 & 19 \\ \hline
\end{array}$$

These data suggest that Problem \ref{findnugeneralcyclic} might be difficult.  Therefore, we also pose the following special cases.

\begin{prob}
For positive integers $n$ and $h$, find $\nu(\mathbb{Z}_{n},3,h)$.
\end{prob}

\begin{prob}
For positive integers $n$  and $m$, find $\nu(\mathbb{Z}_{n},m,2)$.
\end{prob}

\begin{prob}
Find all (or, at least, infinitely many) positive integers $n$, $m$, and $h$, for which $$\nu(\mathbb{Z}_{n},m,h) < \mathrm{min} \left\{ n, {m+h-1 \choose h} \right\}.$$
\end{prob}

\subsection{Limited number of terms} \label{sectionmaxsumsetsizeUlimited}

Here we ought to consider, for a given group $G$, positive integer $m$ (with $m \leq n=|G|$), and nonnegative integer $s$, $$\nu (G,m,[0,s]) = \mathrm{max} \{ |[0,s]A|  \mid A \subseteq G, |A|=m\},$$ that is, the maximum size of $\cup_{h=0}^s hA$ for an $m$-element subset $A$ of $G$. 

Obviously, we have  $$\nu (G,n,[0,s]) =  |[0,s]G| =
\left\{ \begin{array}{cl}
1 & \mbox{if} \; s=0, \\
n & \mbox{if} \; s \geq 1;
\end{array}\right.
$$ hence we may restrict our attention to the cases when $m \leq n-1$.  However, we have the following result.    

\begin{prop} \label{rho(G,m,[0,s])}
For any group $G$, positive integer $m \leq n-1$, and nonnegative integer $s$ we have
$$\nu (G,m,[0,s]) = \nu (G,m+1,s).$$ 
\end{prop}

We can prove Proposition \ref{rho(G,m,[0,s])}, as follows.  Suppose first that $A$ is a subset of $G$ of size $m$ and that it has maximum-size $[0,s]$-fold sumset: $$|[0,s]A|=\nu (G,m,[0,s]).$$  Clearly, $|A \cup \{0\}| \leq m+1$, so $$|s(A \cup \{0\})| \leq \nu (G,m+1,s).$$  But $[0,s]A=s(A \cup \{0\})$  (see Proposition \ref{[0,s]A=}), and therefore $$\nu (G,m,[0,s]) \leq \nu (G,m+1,s).$$

For the other direction, choose a subset $A$ of $G$ of size $m+1$ for which $$|sA|=\nu (G,m+1,s).$$  By Proposition \ref{shifting sets}, we may assume that $0 \in A$; let $A'=A \setminus \{0\}$.  Then $A'$ has size $m$, and we have
$$\nu (G,m+1,s) = |sA| = |s (A' \cup \{0\})| = |[0,s]A'| \leq \nu (G,m,[0,s]).$$  Therefore, $$\nu (G,m,[0,s]) = \nu (G,m+1,s),$$ as claimed.

Proposition \ref{rho(G,m,[0,s])} makes this section superfluous (or makes Section \ref{sectionmaxsumsetsizeUfixed} superfluous via this section).

\subsection{Arbitrary number of terms} \label{sectionmaxsumsetsizeUarbitrary}

Here we consider, for a given group $G$ and positive integer $m$ (with $m \leq n$), $$\nu (G,m,\mathbb{N}_0) = \mathrm{max} \{ |\langle A \rangle|  \mid A \subseteq G, |A|=m\}.$$ Recall that $\langle A \rangle$ is the subgroup of $G$ generated by $A$.

We can easily determine that $$\nu (\mathbb{Z}_n,m,\mathbb{N}_0) = n$$ for all $m$ and $n$, since any set $A$ which contains an element of order $n$ will generate all of $\mathbb{Z}_n$.

More generally, consider the invariant decomposition of $G$,
$$G = \mathbb{Z}_{n_1} \times \mathbb{Z}_{n_2} \times \cdots \times \mathbb{Z}_{n_r},$$ where $r$ and $n_1, \dots, n_r$  are integers all at least 2 and $n_{i+1}$ is divisible by $n_i$ for $i=1,2,\dots, r-1$.  (Here $r$ is the rank of $G$ and $n_r=\kappa$ is the exponent of $G$.)  For each $i=1,2,\dots, r$, let $e_i$ denote the element $(0,\dots,0,1,0,\dots,0)$ of $G$ where the 1 occurs in position $i$.  When $m = r$, taking $A=\{e_1,\dots,e_m\}$ results in $\langle A \rangle = G$, and thus $\nu (G,m,\mathbb{N}_0) =n$; this holds when $m>r$ for any set $A$ containing $\{e_1,\dots,e_m\}$.   When $m < r$, we can take $A=\{e_{r-m+1},e_{r-m+2},\dots,e_r\}$, and this results in $$\langle A \rangle \cong \mathbb{Z}_{n_{_{r-m+1}}} \times \mathbb{Z}_{n_{_{r-m+2}}} \times \cdots \times \mathbb{Z}_{n_r}.$$  We conjecture that we cannot do better:     

\begin{conj} \label{rhoarbitraryconj}
Suppose that $G$ is given by its invariant decomposition, as above.  Prove that
$$\nu (G,m,\mathbb{N}_0) =\left\{
\begin{array}{cl}
n_{_{r-m+1}}  n_{_{r-m+2}} \cdots n_{_{r}} & \mbox{if $m \leq r$}, \\ \\
n & \mbox{if $m \geq r$}.
\end{array} \right.$$  

\end{conj}

\begin{prob}
Prove (or disprove) Conjecture \ref{rhoarbitraryconj}.
\end{prob}

\section{Unrestricted signed sumsets} \label{sectionmaxsumsetsizeUS}

In this section we investigate the quantity $$\nu_{\pm} (G,m,H) =\mathrm{max} \{ |H_{\pm}A|  \mid A \subseteq G, |A|=m\}$$  for various $H \subseteq \mathbb{N}_0$.  

\subsection{Fixed number of terms} \label{sectionmaxsumsetsizeUSfixed}

First, we consider $$\nu_{\pm} (G,m,h) = \mathrm{max} \{ |h_{\pm}A|  \mid A \subseteq G, |A|=m\},$$ that is, the maximum size of an $h$-fold signed sumset of an $m$-element subset of $G$.

In general, we do not know the value of $\nu_{\pm} (G,m,h)$, but we can evaluate it for $h=0$, $h=1$, and $m=1$, as follows.  Since $0_{\pm}A=\{0\}$ for every $A \subseteq G$, we have $$\nu_{\pm} (G,m,0)= 1.$$

Furthermore,  $1_{\pm}A=A \cup (-A)$, so to find $\nu_{\pm} (G,m,1)$, we need to find a subset $A$ of $G$ with $|A|=m$ for which $|A \cup (-A)|$ is maximal.  We can do this as follows.

Let $$L=\mathrm{Ord}(G,2) \cup \{0\}=\{g \in G  \mid 2g=0\}=\{g \in G  \mid g=-g\}.$$ Observe that the elements of $G \setminus L$ are distinct from their inverses, so we have a subset $K$ of $G \setminus L$ with which $G=L \cup K \cup (-K)$, and $L$, $K$, and $-K$ are pairwise disjoint.  We then see that to maximize the quantity $|A \cup (-A)|$, we can have $A \subseteq K$ when $m \leq |K|$; $K \subseteq A \subseteq K \cup L$  when $|K| \leq m \leq |K \cup L|$; and $K \cup L \subseteq A$ when $m \geq |K \cup L|$.  Thus, 
$$ \nu_{\pm} (G,m,1)= 
\left\{
\begin{array}{cl}
2m & \mbox{if \; $m \leq |K|$}\\ \\
m+|K| & \mbox{if \; $|K| \leq m \leq |K \cup L|$}\\ \\
n & \mbox{if \; $m \geq |K \cup L|;$}\\
\end{array}\right. 
$$
since $|K \cup L|=n-|K|$, this simplifies to
$$\nu_{\pm} (G,m,1) =\mathrm{min} \left\{2m, m+ |K|, n \right\} $$ or 
$$\nu_{\pm} (G,m,1) =\mathrm{min} \left\{n, 2m, m+ \frac{n-|\mathrm{Ord}(G,2)|-1}{2} \right\}.$$

Let's turn to the case of $m=1$.  If $A$ consists of a single element $a$, then $h_{\pm}A=\{ha,-ha\}$ for every $h \in \mathbb{N}_0$, so to find $\nu_{\pm} (G,1,h)$, we need to find, if possible, a one-element subset $A=\{a\}$ of $G$ for which the set $\langle A, h\rangle=\{ha,-ha\}$ has size 2.  This means that we only need to answer the following question: when does $G$ contain an element $a$ for which $ha \neq -ha$?   But, for any $a \in G$,  $ha \neq -ha$ is equivalent to $2ha \neq 0$, which is the same as saying that the order of $a$ does not divide $2h$.  Since the order of any element in $G$ is a divisor of the exponent $\kappa$ of $G$, if the order of $a$ does not divide $2h$, then $\kappa$ does not divide $2h$ either.  Conversely, if $\kappa$ does not divide $2h$, then choosing $a$ to be an element of order $\kappa$ will work.  Therefore, we get
$$ \nu_{\pm} (G,1,h) = 
\left\{
\begin{array}{ll}
1 & \mbox{if \; $\kappa | 2h$;}\\
2 & \mbox{otherwise.}\\
\end{array}\right. 
$$

In summary, we have the following.

\begin{prop} \label{nupmtrivial}
In any abelian group $G$ of order $n$ and exponent $\kappa$ we have
\begin{eqnarray*}
\nu_{\pm} (G,m,0) & = & 1; \\
\nu_{\pm} (G,m,1) & = & \mathrm{min} \left\{n, 2m, m+ \frac{n-|\mathrm{Ord}(G,2)|-1}{2} \right\}; \\
\nu_{\pm} (G,1,h) & = &  
\left\{
\begin{array}{ll}
1 & \mbox{if \; $\kappa | 2h$,}\\
2 & \mbox{otherwise.}\\
\end{array}\right.
\end{eqnarray*}

\end{prop}

For values of $h \geq 2$ and $m \geq 2$, we have no exact values for $\nu_{\pm} (G,m,h)$ in general; as a consequence of Proposition \ref{nuupper}, we have, however, the following upper bound.

\begin{prop} \label{nuupperUSfixed}  In any abelian group $G$ of order $n$ we have
$$\nu_{\pm}(G,m,h) \leq \mathrm{min} \left\{ n, c(h,m) \right\}.$$
\end{prop}
Recall that $$c(h,m)=\sum_{i \geq 0}  {m \choose i}{h-1 \choose i-1} 2^i .$$

Similarly to Section \ref{sectionmaxsumsetsizeUfixed}, we pose the following problems.

\begin{prob}
Find the value of (or, at least, find good bounds for)  $\nu_{\pm}(G,m,h)$ for noncyclic groups $G$ and integers $m$ and $h$.
\end{prob}

\begin{prob}
Find the value of (or, at least, find good bounds for)  $\nu_{\pm}(\mathbb{Z}_n,m,h)$ for all positive integers $n,m$, and $h$. 

\end{prob}

\begin{prob}
For positive integers $n$ and $h$, find $\nu_{\pm}(\mathbb{Z}_{n},2,h)$.
\end{prob}

\begin{prob}
For positive integers $n$  and $m$, find $\nu_{\pm}(\mathbb{Z}_{n},m,2)$.
\end{prob}

\begin{prob} \label{nu pm less}
Find all (or, at least, infinitely many) positive integers $n$, $m$, and $h$, with $m \geq 2$ and $h \geq 2$, for which $$\nu_{\pm}(\mathbb{Z}_{n},m,h) < \mathrm{min} \left\{ n, c(h,m) \right\}.$$
\end{prob}

Some of the cases when the inequality of Problem \ref{nu pm less} holds were found by Buell\index{Buell, M.} in \cite{Bue:2014a}; namely, when $2 \leq m \leq 5$, the set of $n$ values for which $\nu_{\pm}(\mathbb{Z}_{n},m,2)$ is less than $\mathrm{min} \left\{ n, c(2,m) \right\}$ are as follows:  
$$\begin{array}{||c|l||} \hline
 m & n \\ \hline \hline
2 & 6-8 \\ \hline
3 & 14-18, 20 \\ \hline
4 & 22-33, 36, 40 \\ \hline
5 & 34, 36-54, 56-58 \\ \hline \hline
\end{array}$$
From these data it appears that a complete solution to our problems above may be quite challenging.

\subsection{Limited number of terms} \label{sectionmaxsumsetsizeUSlimited}

Here we consider, for a given group $G$, positive integer $m$ (with $m \leq n=|G|$), and nonnegative integer $s$, $$\nu_{\pm} (G,m,[0,s]) = \mathrm{max} \{ |[0,s]_{\pm}A|  \mid A \subseteq G, |A|=m\},$$ that is, the maximum size of $[0,s]_{\pm}A$ for an $m$-element subset $A$ of $G$.

We note that we don't have a version of Proposition \ref{shifting sets} for signed sumsets, so we are not able to reduce this entire section to Section \ref{sectionmaxsumsetsizeUSfixed}.  (However, one may be able to apply similar techniques.) 

It is easy to see that, for every $m$, we have $$\nu_{\pm} (G,m,[0,0]) =1.$$  Furthermore,  $$[0,1]_{\pm}A=A \cup (-A) \cup \{0\},$$ whose maximum size we can find as we did for $A \cup (-A)$ in Section \ref{sectionmaxsumsetsizeUSfixed}, except here we write $G$ as the pairwise disjoint union of four (potentially empty) parts: $\{0\}$, $\mathrm{Ord}(G,2)$, $K$, and $-K$.  The computation this time yields 
$$\nu_{\pm} (G,m,[0,1]) =\mathrm{min} \left\{n, 2m+1, m+ \frac{n-|\mathrm{Ord}(G,2)|+1}{2} \right\}.$$

Regarding the case of $m=1$, we see that for $A=\{a\}$ we have $$[0,s]_{\pm}A =\{0,\pm a, \pm 2a,\dots,\pm sa\},$$
thus $$\nu_{\pm} (G,1,[0,s]) =\mathrm{min}\{\kappa,2s+1\}.$$ 

Summarizing our results thus far, we have the following.

\begin{prop} \label{nuUStrivial} In any abelian group $G$ of order $n$ and exponent $\kappa$ we have
\begin{eqnarray*}
\nu_{\pm} (G,m,[0,0]) & = & 1, \\
\nu_{\pm} (G,m,[0,1]) & = & \mathrm{min} \left\{n, 2m+1, m+ \frac{n-|\mathrm{Ord}(G,2)|+1}{2} \right\}, \\
\nu_{\pm} (G,1,[0,s]) & = & \mathrm{min}\{\kappa,2s+1\}.
\end{eqnarray*}
\end{prop} 

For values of $s \geq 2$ and $m \geq 2$, we have no exact values for $\nu_{\pm} (G,m,[0,s])$ in general; as a consequence of Proposition \ref{nuupper}, we have, however, the following upper bound.

\begin{prop} \label{nuupperUSlimited}  In any abelian group $G$ of order $n$ we have
$$\nu_{\pm}(G,m,[0,s]) \leq \mathrm{min} \left\{ n, a(m,s)\right\}.$$
\end{prop}
Recall that $$a(m,s)=\sum_{i \geq 0}  {m \choose i}{s \choose i} 2^i .$$

Similarly to previous sections, we pose the following problems.

\begin{prob}
Find the value of (or, at least, find good bounds for)  $\nu_{\pm}(G,m,[0,s])$ for noncyclic groups $G$ and integers $m$ and $s$.
\end{prob}

\begin{prob}
Find the value of (or, at least, find good bounds for)  $\nu_{\pm}(\mathbb{Z}_n,m,[0,s])$ for all integers $n$, $m$, and $s$. 

\end{prob}

\begin{prob}
For positive integers $n$ and $s$, find $\nu_{\pm}(\mathbb{Z}_{n},2,[0,s])$.
\end{prob}

\begin{prob}
For positive integers $n$  and $m$, find $\nu_{\pm}(\mathbb{Z}_{n},m,[0,2])$.
\end{prob}

\begin{prob} \label{nu pm less}
Find all (or, at least, infinitely many) positive integers $n$, $m$, and $s$ for which $$\nu_{\pm}(\mathbb{Z}_{n},m,[0,s]) < \mathrm{min} \left\{ n, a(m,s) \right\}.$$
\end{prob}

\subsection{Arbitrary number of terms} \label{sectionmaxsumsetsizeUSarbitrary}

This subsection is identical to Subsection \ref{sectionmaxsumsetsizeUarbitrary}.

\section{Restricted sumsets} \label{sectionmaxsumsetsizeR}

Our goal in this section is to investigate the quantity $$\nu\hat{\;} (G,m,H) =\mathrm{max} \{ |H\hat{\;}A|  \mid A \subseteq G, |A|=m\}$$ for various $H \subseteq \mathbb{N}_0$.  

\subsection{Fixed number of terms} \label{sectionmaxsumsetsizeRfixed}

Here we consider $$\nu\hat{\;} (G,m,h) = \mathrm{max} \{ |h\hat{\;}A|  \mid A \subseteq G, |A|=m\},$$ that is, the maximum size of a restricted $h$-fold sumset of an $m$-element subset of $G$.  

Observe first that for all $h>m$, $h\hat{\;}A =\emptyset$ and thus $\nu\hat{\;} (G,m,h)=0$. 
Furthermore, note that for every set $A=\{a_1,\dots,a_m\} \subseteq G$ of size $m$ and for every $h \in \mathbb{N}_0$, we have
$$(m-h)\hat{\;}A=(a_1+\cdots +a_m) - h\hat{\;}A;$$
in particular, $$|(m-h)\hat{\;}A|=|h\hat{\;}A|$$ and
$$\nu\hat{\;} (G,m,m-h) = \nu\hat{\;} (G,m,h).$$
Therefore, we can restrict our attention to the cases when $$h \leq \left\lfloor \frac{m}{2} \right\rfloor.$$

It is also quite obvious that we have $0\hat{\;}A=\{0\}$ and $1\hat{\;}A=A$; consequently, we have the following.

\begin{prop} \label{nuRtrivial} In any abelian group $G$ we have

$\nu\hat{\;} (G,m,0) = 1, $

$\nu\hat{\;}(G,m,1) = m, $

$\nu\hat{\;} (G,m,m-1) = m,$ 

$\nu\hat{\;} (G,m,m) = 1, $

$\nu\hat{\;} (G,m,h) = 0 \; \mbox{for} \; h>m.$

\end{prop}

For values of $2 \leq h \leq m-2$, we have no exact values for $\nu \hat{\;}(G,m,h)$ in general; as a consequence of Proposition \ref{nuupper}, we have, however, the following upper bound.

\begin{prop} \label{nuupperRfixed}  In any abelian group $G$ of order $n$ we have
$$\nu\hat{\;}(G,m,h) \leq \mathrm{min} \left\{ n, {m \choose h} \right\}.$$
\end{prop}

We pose the following problems.

\begin{prob}
Find the value of (or good bounds for) $\nu\hat{\;}(G,m,h)$ for noncyclic groups $G$ and integers $m$ and $h$.
\end{prob}

\begin{prob}
Find the value of (or good bounds for)  $\nu\hat{\;}(\mathbb{Z}_n,m,h)$ for all positive integers $n,m$, and $h$. 

\end{prob}

\begin{prob}
For positive integers $n$ and $m$, find $\nu\hat{\;}(\mathbb{Z}_n,m,2)$.
\end{prob}

\begin{prob} \label{nu pm less}
Find all (or, at least, infinitely many) positive integers $n$, $m$, and $h$ for which $$\nu\hat{\;}(\mathbb{Z}_n,m,h) <  \min \left\{ n, {m \choose h} \right\}.$$
\end{prob}

Manandhar\index{Manandhar, D.} (cf.~\cite{Man:2012a}) found (using a computer program) that for all $n \leq 20$, we have  
$$\nu\hat{\;}(\mathbb{Z}_{n},m,h)=\mathrm{min} \left\{ n, {m \choose h} \right\},$$ with the following exceptions:
$$\begin{array}{|c|c|c||c|c|} \hline
n & m & h & \mathrm{min} \left\{ n, {m \choose h} \right\} & \nu\hat{\;}(\mathbb{Z}_{n},m,h) \\ \hline \hline
10 & 5 & 2, 3 & 10 & 9 \\ \hline
14 & 6 & 2, 4 & 14 & 13 \\ \hline
15 & 6 & 2, 4 & 15 & 13 \\ \hline
16 & 6 & 2, 4 & 15 & 14 \\ \hline
17 & 6 & 2, 4 & 15 & 14 \\ \hline
18 & 6 & 2, 4 & 15 & 14 \\ \hline
18 & 7 & 2, 5 & 18 & 17 \\ \hline
19 & 7 & 2, 5 & 19 & 18 \\ \hline
20 & 7 & 2, 5 & 20 & 19 \\ \hline
\end{array}$$
It is interesting to note that all these exceptions occur with $h=2$ or $h=m-2$.

\subsection{Limited number of terms} \label{sectionmaxsumsetsizeRlimited}

Here we consider, for a given group $G$, positive integer $m$ (with $m \leq n=|G|$), and nonnegative integer $s$, $$\nu\hat{\;} (G,m,[0,s]) = \mathrm{max} \{ |[0,s]\hat{\;}A|  \mid A \subseteq G, |A|=m\},$$ that is, the maximum size of $[0,s]\hat{\;}A$ for an $m$-element subset $A$ of $G$.  We may assume that $s \leq m$, since for $s > m$ we have $$ \nu\hat{\;} (G,m,[0,s])=\nu\hat{\;} (G,m,[0,m]).$$

It is easy to see that, for every $m$, we have $$\nu\hat{\;} (G,m,[0,0]) =1$$ and $$\nu\hat{\;} (G,m,[0,1]) =\min \{n,m+1\}.$$  

For values of $2 \leq s \leq m$, we have no exact values for $\nu\hat{\;} (G,m,[0,s])$ in general; as a consequence of Proposition \ref{nuupper}, we have, however, the following upper bound.

\begin{prop} \label{nuupperRlimited}  In any abelian group $G$ of order $n$ we have
$$\nu\hat{\;}(G,m,[0,s]) \leq \mathrm{min} \left\{ n, \sum_{h=0}^s {m \choose h} \right\}.$$
\end{prop}

As in previous section, we have the following problems.

\begin{prob}
For positive integers $n$, $m$,  and $s$, find the value of (or, at least, good bounds for) $\nu\hat{\;}(\mathbb{Z}_{n},m,[0,s])$.  In particular, find $\nu\hat{\;}(\mathbb{Z}_{n},m,[0,2])$.
\end{prob}

\begin{prob}
Find all (or, at least, infinitely many) positive integers $n$, $m$, and $s$, for which $$\nu\hat{\;}(\mathbb{Z}_{n},m,[0,s]) < \mathrm{min} \left\{ n, \sum_{h=0}^s {m \choose h} \right\}.$$
\end{prob}

\begin{prob}
Find the value of (or, at least, good bounds for) $\nu\hat{\;}(G,m,[0,s])$ for noncyclic groups $G$ and integers $m$ and $s$.
\end{prob}

\subsection{Arbitrary number of terms} \label{sectionmaxsumsetsizeRarbitrary}

Here we consider, for a given group $G$ and positive integer $m$ (with $m \leq n$), $$\nu\hat{\;} (G,m,\mathbb{N}_0) = \mathrm{max} \{ |\Sigma  A|  \mid A \subseteq G, |A|=m\}.$$

Note that $$\nu\hat{\;} (G,m,\mathbb{N}_0) =\nu\hat{\;} (G,m,[0,s])$$ for every $s \geq m$, thus we could think of the problem of finding $\nu\hat{\;} (G,m,\mathbb{N}_0)$, the maximum size of a restricted sumset, as a special case of finding $\nu\hat{\;} (G,m,[0,s])$ when $s=m$ (or any $s \geq m$).  However, it may be worthwhile to separate this case as it is of special interest.

As a consequence of Proposition \ref{nuupper}, we have the following upper bound.

\begin{prop} \label{nuupperRarbitrary}  In any abelian group $G$ of order $n$ we have
$$\nu\hat{\;}(G,m,\mathbb{N}_0) \leq \mathrm{min} \left\{ n, 2^m \right\}.$$

\end{prop}

We can easily prove that, in the case of cyclic groups, equality holds in Proposition \ref{nuupperRarbitrary}.  We will consider two cases: when $n \geq 2^m$ and when $n < 2^m$.

If $n \geq 2^m$, then let $$A=\{1,2,2^2,\dots,2^{m-1}\}.$$  Clearly, $A$ has size $m$, and one can also see (recalling the base 2 representation of integers)  that $$|\Sigma A|=2^m=\mathrm{min} \left\{ n, 2^m \right\}.$$  

Suppose now that $n < 2^m$, and let $k=\lceil \log_2 n \rceil$; note that $2^{k-1} < n \leq 2^{k}$.  Therefore, for the set
$$A'=\{1,2, 2^2, \dots, 2^{k-1}\},$$ we get $\Sigma A'=\mathbb{Z}_n;$ and this implies that for any set $A$ that contains $A'$ we have $\Sigma A=\mathbb{Z}_n.$ But $n < 2^m$ implies that $m > \log_2 n$ and thus $m \geq k$, which means that, again we are able to find an $m$-subset $A$ of $\mathbb{Z}_n$ for which $$|\Sigma A|=n=\mathrm{min} \left\{ n, 2^m \right\}.$$       

These constructions, together with Proposition \ref{nuupperRarbitrary}, imply the following:

\begin{prop}

For positive integers $n$ and $m$ with $m \leq n$ we have
$$\nu\hat{\;}(\mathbb{Z}_n,m,\mathbb{N}_0) = \mathrm{min} \left\{ n, 2^m \right\}.$$

\end{prop}
As we have mentioned above, this also implies that 
$$\nu\hat{\;}(\mathbb{Z}_n,m,[0,s]) = \mathrm{min} \left\{ n, 2^m \right\}$$ holds for every integer $s \geq m$.

So, the case of cyclic groups has been settled, which leaves us with the following problem.

\begin{prob}
For every positive integer $m$ and noncyclic group $G$, find the value of $\nu\hat{\;}(G,m,\mathbb{N}_0)$.
\end{prob}

\section{Restricted signed sumsets} \label{sectionmaxsumsetsizeRS}

Our goal in this section is to investigate the quantity $$\nu\hat{_\pm} (G,m,H) =\mathrm{max} \{ |H\hat{_\pm} A|  \mid A \subseteq G, |A|=m\}$$  for various $H \subseteq \mathbb{N}_0$.  

\subsection{Fixed number of terms} \label{sectionmaxsumsetsizeRSfixed}

Here we consider $$\nu\hat{_\pm} (G,m,h) = \mathrm{max} \{ |h\hat{_\pm}A|  \mid A \subseteq G, |A|=m\},$$ that is, the maximum size of a restricted $h$-fold signed sumset of an $m$-element subset of $G$.  

Observe first that, for all $h>m$, $h\hat{_\pm}A =\emptyset$ and thus $\nu\hat{_\pm} (G,m,h)=0$, so only $h \leq m$ needs to be considered, as was the case in Section \ref{sectionmaxsumsetsizeRfixed}.  On the other hand, here we do not have the ``palindromic'' property of restricted sumsets that  $\nu\hat{\;} (G,m,m-h) = \nu\hat{\;} (G,m,h).$

As in Section \ref{sectionmaxsumsetsizeUSfixed}, we get 

\begin{prop} \label{nupmtrivial}
In any abelian group $G$ of order $n$ and for every $m \leq n$ we have
\begin{eqnarray*}
\nu\hat{_\pm} (G,m,0) & = & 1; \\
\nu\hat{_\pm} (G,m,1) & = & \mathrm{min} \left\{n, 2m, m+ \frac{n-|\mathrm{Ord}(G,2)|-1}{2} \right\}. \\
\end{eqnarray*}

\end{prop}

For values of $2 \leq h \leq m$, we have no exact values for $\nu\hat{_\pm} (G,m,h)$ in general; as a consequence of Proposition \ref{nuupper}, we have, however, the following upper bound.

\begin{prop} \label{nuupperRSfixed}  In any abelian group $G$ of order $n$ we have
$$\nu\hat{_\pm}(G,m,h) \leq \mathrm{min} \left\{ n, {m \choose h} \cdot 2^h \right\}.$$
\end{prop}

Not much else is known about the value of $\nu\hat{_\pm}(G,m,h)$ in general.  We thus pose the following.

\begin{prob}
Find the value of (or, at least, good bounds for) $\nu\hat{_\pm}(G,m,h)$ for noncyclic groups $G$ and integers $m$ and $2 \leq h \leq m$.
\end{prob}

At least for cyclic groups, we can find the exact answer when $h=m$.  As conjectured by Olans\index{Olans, T.} in \cite{Ola:2013a}, we have the following result:

\begin{prop} \label{nu m=h cyclic}  
For all positive integers $n$ and $m \leq n$ we have
$$\nu\hat{_\pm}(\mathbb{Z}_n,m,m) =
\left\{
\begin{array}{ll}
\min \left\{ n, 2^m \right\} & \mbox{if $n$ is odd} ;\\ \\
\min \left\{ n/2, 2^m \right\} & \mbox{if $n$ is even} .
\end{array}
\right.$$
\end{prop}

The easy proof is on page \pageref{proof of nu m=h cyclic}.  This result leaves us with the following open questions.

\begin{prob}  \label{problem cyclic nu hat plus minus}
Find the value of $\nu\hat{_\pm}(\mathbb{Z}_n,m,h)$ for all positive integers $n,m$, and $2 \leq h \leq m-1$. 

\end{prob}

The following two special cases of Problem \ref{problem cyclic nu hat plus minus} are worth mentioning separately.

\begin{prob}
For positive integers $n$  and $m$, find $\nu\hat{_\pm}(\mathbb{Z}_{n},m,2)$.
\end{prob}

\begin{prob}
Find all (or, at least, infinitely many) positive integers $n$, $m$, and $2 \leq h \leq m-1$, for which $$\nu\hat{_\pm}(\mathbb{Z}_{n},m,h) < \mathrm{min} \left\{ n, {m \choose h} \cdot 2^h \right\}.$$
\end{prob}

\subsection{Limited number of terms} \label{sectionmaxsumsetsizeRSlimited}

Here we consider, for a given group $G$, positive integer $m$ (with $m \leq n=|G|$), and nonnegative integer $s$, $$\nu\hat{_{\pm}} (G,m,[0,s]) = \mathrm{max} \{ |[0,s]\hat{_{\pm}}A|  \mid A \subseteq G, |A|=m\},$$ that is, the maximum size of $[0,s] \hat{_{\pm}}A$ for an $m$-element subset $A$ of $G$.  We may assume that $s \leq m$, since for $s > m$ we have $$ \nu\hat{_{\pm}} (G,m,[0,s])=\nu\hat{_{\pm}} (G,m,[0,m]).$$

As in Section \ref{sectionmaxsumsetsizeUSlimited}, we get 

\begin{prop} \label{nupmtrivial}
In any abelian group $G$ of order $n$ and for every $m \leq n$ we have
\begin{eqnarray*}
\nu\hat{_\pm} (G,m,[0,0]) & = & 1; \\
\nu\hat{_\pm} (G,m,[0,1]) & = & \mathrm{min} \left\{n, 2m+1, m+ \frac{n-|\mathrm{Ord}(G,2)|+1}{2} \right\}. \\
\end{eqnarray*}

\end{prop}

For other values of $s$, we have no exact values for $\nu\hat{_{\pm}} (G,m,[0,s])$ in general; as a consequence of Proposition  \ref{nuupper}, we have, however, the following upper bound.

\begin{prop} \label{nuupperRSlimited}  In any abelian group $G$ of order $n$ we have
$$\nu\hat{_{\pm}}(G,m,[0,s]) \leq \mathrm{min} \left\{ n, \sum_{h=0}^s {m \choose h} 2^h \right\}.$$
\end{prop}

The following problems are largely unsolved.

\begin{prob}
Find the value of (or, at least, good bounds for) $\nu\hat{_{\pm}}(\mathbb{Z}_n,m,[0,s])$ for all integers $n$, $m$, and $s$; in particular, find $\nu\hat{_{\pm}}(\mathbb{Z}_n,m,[0,2])$.

\end{prob}

\begin{prob}
Find all (or, at least, infinitely many) positive integers $n$, $m$, and $s$, for which $$\nu\hat{_\pm}(\mathbb{Z}_{n},m,[0,s]) < \mathrm{min} \left\{ n, \sum_{h=0}^s {m \choose h} 2^h \right\}.$$
\end{prob}

\begin{prob}
Find the value of (or, at least, good bounds for) $\nu\hat{_{\pm}}(G,m,[0,s])$ for noncyclic groups $G$ and integers $m$ and $s$.
\end{prob}

\subsection{Arbitrary number of terms} \label{sectionmaxsumsetsizeRSarbitrary}

Here we consider, for a given group $G$ and positive integer $m$ (with $m \leq n$), $$\nu\hat{_\pm} (G,m,\mathbb{N}_0) = \mathrm{max} \{ |\Sigma_{\pm} A|  \mid A \subseteq G, |A|=m\}.$$

Note that $$\nu\hat{_\pm} (G,m,\mathbb{N}_0) =\nu\hat{_\pm} (G,m,[0,s])$$ for every $s \geq m$, thus we could think of the problem of finding $\nu\hat{_\pm} (G,m,\mathbb{N}_0)$, the maximum size of a restricted sumset, as a special case of finding $\nu\hat{_\pm} (G,m,[0,s])$ when $s=m$ (or any $s \geq m$).  Nevertheless, we separate this subsection as it is of special interest.

As a consequence of Proposition \ref{nuupper}, we have the following upper bound.

\begin{prop} \label{nuupperRSarbitrary}  In any abelian group $G$ of order $n$ we have
$$\nu\hat{_\pm}(G,m,\mathbb{N}_0) \leq \mathrm{min} \left\{ n, 3^m \right\}.$$

\end{prop}

With an argument similar to the one in Section \ref{sectionmaxsumsetsizeRarbitrary}, we show that, when $G$ is cyclic, equality holds in Proposition \ref{nuupperRSarbitrary}.  
We will consider two cases: when $n \geq 3^m$ and when $n < 3^m$.

If $n \geq 3^m$, then let $$A=\{1,3,3^2,\dots,3^{m-1}\};$$  clearly, $A$ has size $m$.  Recall that every positive integer up to $3^m-1$ has a (unique) ternary representation of at most $m$ ternary digits.  Subtracting $$1+3+3^2+\cdots+3^{m-1}=\frac{3^{m}-1}{2},$$ every integer between $-(3^{m}-1)/2$ and $(3^{m}-1)/2$, inclusive, can be written (uniquely) as
$$r_0+r_1 \cdot 3 + \cdots + r_{m-1} \cdot 3^{m-1}$$ with $r_0, r_1, \dots, r_{m-1} \in \{-1,0,1\}$.  Therefore, $$|\Sigma_{\pm}A|=3^m=\mathrm{min} \left\{ n, 3^m \right\}.$$  

Suppose now that $n < 3^m$, and let $k=\lceil \log_3 n \rceil$; note that $3^{k-1} < n \leq 3^{k}$.  Therefore, for the set
$$A'=\{1,3, 3^2, \dots, 3^{k-1}\},$$ we get $\Sigma_{\pm}A'=\mathbb{Z}_n;$ and this implies that for any set $A$ that contains $A'$ we have $\Sigma_{\pm}A=\mathbb{Z}_n.$ But $n < 3^m$ implies that $m > \log_3 n$ and thus $m \geq k$, which means that, again we are able to find an $m$-subset $A$ of $\mathbb{Z}_n$ for which $$|\Sigma_{\pm}A|=n=\mathrm{min} \left\{ n, 3^m \right\}.$$       

These constructions, together with Proposition \ref{nuupperRSarbitrary}, imply the following:

\begin{prop}

For positive integers $n$ and $m$ with $m \leq n$ we have
$$\nu\hat{_\pm}(\mathbb{Z}_n,m,\mathbb{N}_0) = \mathrm{min} \left\{ n, 3^m \right\}.$$

\end{prop}
As we have mentioned above, this also implies that 
$$\nu\hat{_\pm}(\mathbb{Z}_n,m,[0,s]) = \mathrm{min} \left\{ n, 3^m \right\}$$ holds for every integer $s \geq m$.

So, the case of cyclic groups has been settled, which leaves us with the following problem.

\begin{prob}
For every positive integer $m$ and noncyclic group $G$, find the value of $\nu\hat{_\pm}(G,m,\mathbb{N}_0)$.
\end{prob}

\chapter{Spanning sets} \label{ChapterSpanning}

Recall that for a given finite abelian group $G$, $m$-subset $A=\{a_1,\dots, a_m\}$ of $G$, $\Lambda \subseteq \mathbb{Z}$, and $H \subseteq \mathbb{N}_0$, we defined the sumset of $A$ corresponding to $\Lambda$ and $H$ as
$$H_{\Lambda}A = \{\lambda_1a_1+\cdots +\lambda_m a_m \mbox{    } |  \mbox{    }  (\lambda_1,\dots ,\lambda_m) \in \Lambda^m(H)  \}$$
where the index set $\Lambda^m(H)$ is defined as
$$\Lambda^m(H)=\{(\lambda_1,\dots ,\lambda_m) \in \Lambda^m \; |  \; |\lambda_1|+\cdots +|\lambda_m| \in H \}.$$

The case when the sumset yields the entire group is of special interest; we say that $A$ is an {\em $H$-spanning set over $\Lambda$} if $H_{\Lambda}A=G$.  In this chapter we attempt to find the minimum possible size of an $H$-spanning set over $\Lambda$ in a given finite abelian group $G$.   Namely, our objective is to determine, for any $G$, $\Lambda \subseteq \mathbb{Z}$, and $H \subseteq \mathbb{N}_0$ the quantity
$$\phi_{\Lambda}(G,H)=\mathrm{min} \{ |A|  \mid A \subseteq G, H_{\Lambda}A=G\}.$$  If no $H$-spanning set exists, we put $\phi_{\Lambda}(G,H) = \infty.$

Note that we have a strong connection between the maximum possible size of sumsets, studied in Chapter \ref{ChapterMaxsumsetsize}, and the minimum size of spanning sets, studied here.  In particular, we have the following obvious proposition.

\begin{prop} \label{chapter3specialofchapter2}
For any group $G$ of size $n$, $\Lambda \subseteq \mathbb{Z}$, and $H \subseteq \mathbb{N}_0$, we have
$$\phi_{\Lambda}(G,H) = \mathrm{min} \{m  \mid \nu_{\Lambda}(G,m,H)=n\}. $$

\end{prop}
Therefore, in theory, the question of finding $\phi_{\Lambda}(G,H)$ is a special case of the more general problem of finding $\nu_{\Lambda}(G,m,H)$; however, there is enough interest in spanning sets alone to treat them separately.

We have the following obvious bound.

\begin{prop} \label{philower}
If $A$ is an $H$-spanning set over $\Lambda$ in a group $G$ of order $n$ and $|A|=m$, then
$$n \leq |\Lambda^m(H)|.$$
\end{prop}

Proposition \ref{philower} provides a lower bound for the size of $H$-spanning sets over $\Lambda$ in $G$.  The case of equality in Proposition \ref{philower} is of special interest: an $H$-spanning set over $\Lambda$ of size $m$ in a group $G$ is called {\em perfect}  when $n = |\Lambda^m(H)|$ holds.
 
In the following subsections we consider $\phi_{\Lambda}(G,H)$ for special coefficient sets $\Lambda$.

\section{Unrestricted sumsets} \label{3minU}

Our goal in this section is to investigate the quantity $$\phi(G,H) =\mathrm{min} \{ |A|  \mid A \subseteq G, H A =G\};$$ if no $H$-spanning set exists, we put $\phi(G,H) = \infty.$  Note that we have $$\phi(G,\{0\})=\infty$$ for all groups of order at least 2, but if $H$ contains at least one positive integer $h$, then $$\phi(G,H) \leq n$$ since we have $$G=\{(h-1) \cdot 0 + 1 \cdot g  \mid g \in G\}  \subseteq hG \subseteq HG$$ and thus $HG=G$.    

\subsection{Fixed number of terms} \label{3minUfixed}

In this subsection we ought to consider $$\phi (G,h) = \mathrm{min} \{ |A|  \mid A \subseteq G, hA =G\},$$ that is, the minimum value of $m$ for which $G$ contains a set $A$ of size $m$ with $hA=G.$  

However, using Propositions \ref{chapter3specialofchapter2} and  \ref{rho(G,m,[0,s])}, we can easily reduce the study of $\phi (G,h)$ to that of 
$$\phi (G,[0,h])=\mathrm{min} \{ |A|  \mid A \subseteq G, [0,h]A =G\}.$$  Namely, we have the following result.    

\begin{prop} \label{phi(G,[0,h])}
For any group $G$ and nonnegative integer $h$ we have
$$\phi (G,h) = \phi (G,[0,h])+1.$$ 
\end{prop}

According to Proposition \ref{phi(G,[0,h])}, it suffices to study only one of $\phi (G,h)$ or $\phi (G,[0,h])$.  We elect to study the latter---see Subsection \ref{3minUlimited} below.

\subsection{Limited number of terms} \label{3minUlimited}

A subset $A$ of $G$ for which $[0,s]A =G$ for some nonnegative integer $s$ is called an {\em $s$-basis} for $G$.  (The term is somewhat confusing as for a set to be a basis, it only needs to be spanning---no independence property is assumed.  Nevertheless, we keep this historical terminology.)

Here we investigate $$\phi (G,[0,s]) = \mathrm{min} \{ |A|  \mid A \subseteq G, [0,s]A =G\},$$ that is, the minimum size of an $s$-basis for $G$.  

As we pointed out above, we have $$\phi (G,[0,0]) =\infty$$ for all groups of order at least 2, but $$\phi(G,[0,s]) \leq n$$ for every $s \in \mathbb{N}$.  Furthermore, for every $G$ of order 2 or more we get $$\phi(G,[0,1])=n-1,$$ since for every $A \subseteq G$ we have $[0,1]A=\{0\} \cup A$.  But for $s \geq 2$, we do not have a formula for $\phi(G,[0,s])$.  Therefore we pose the following problem.

\begin{prob}
Find $\phi(G,[0,s])$ for every $G$ and $s \geq 2$.
\end{prob}

From Proposition \ref{philower}, we have the following general bound.

\begin{prop} \label{philowerUlimited}  If $A$ is an $s$-basis of size $m$ in $G$, then
$$n \leq {m+s \choose s}.$$
\end{prop}

We are particularly interested in the extremal cases of Proposition \ref{philowerUlimited}.  Namely, we would like to find {\em perfect $s$-bases}, that is, $s$-bases of size $m$ with $$n={m+s \choose s}.$$

It is easy to determine all perfect $s$-bases for $s=1$ or $m=1$:  

\begin{prop}
Let $G$ be a finite abelian group, $A \subseteq G$, and $a \in G$.

\begin{enumerate}
  
\item $A \subseteq G$ is a perfect $1$-basis in $G$ if, and only if, $A=G \setminus \{0\}$.
\item $\{a\} \subseteq G$ is a perfect $s$-basis in $G$ if, and only if, $G \cong \mathbb{Z}_{s+1}$ and $\gcd (a, s+1)=1$.
\end{enumerate}

\end{prop}

We are not aware of any other perfect bases; in fact, we believe that there are none:

\begin{conj} \label{no perfect bases}
There are no perfect $s$-bases of size $m$ in $G$, unless $s=1$ or $m=1$.

\end{conj}

The following result says that Conjecture \ref{no perfect bases} holds for $s=2$ and $s=3$:

\begin{thm} \label{no perfect bases s=2,3}
For $s \in \{2,3\}$, there are no perfect $s$-bases in $G$ of size $m \geq 2$.

\end{thm}
We present the proof starting on page \pageref{proof of no perfect bases s=2,3}.   Our proof there can probably be generalized, so we offer:

\begin{prob}
Prove Conjecture \ref{no perfect bases}.
\end{prob}

Let us return to the question of finding $\phi (G,[0,s])$.  As a consequence of Proposition \ref{philowerUlimited}, we have the following lower bound for $\phi (G,[0,s])$.  

\begin{prop} \label{philowergeneral}
For any abelian group $G$ of order $n$ and positive integer $s$ we have $$\phi (G,[0,s]) > \left \lceil \sqrt[s]{s! n} \right\rceil -s.$$

\end{prop}
  
To find some upper bounds for $\phi (G,[0,s])$, we exhibit an explicit $s$-basis, as follows.  Consider first the cyclic group $\mathbb{Z}_n$, and let $a=\lceil \sqrt[s]{n} \rceil$.  It is then easy to see that the set
$$A=\{ i \cdot a^j \mid i=1,2,\dots,a-1; j=0,1,\dots,s-1\}$$ is an $s$-basis for $\mathbb{Z}_n$ (this follows from the fact that all nonnegative integers up to $a^s-1$ have a base $a$ representation of at most $s$ digits, and that $a^s-1 \geq n-1$).  Since $|A| \leq (a-1) \cdot s$, we get the following result.

\begin{prop} \label{phiuppercyclic}
For positive integers $n$ and $s$ we have $$\phi (\mathbb{Z}_n,[0,s]) \leq s \cdot \left( \lceil \sqrt[s]{n} \rceil -1 \right) < s \cdot \sqrt[s]{n} .$$

\end{prop} 

\label{from cyclic to gen} Observe that if $A_1$ and $A_2$ are $s$-bases for groups $G_1$ and $G_2$, respectively, then $$A=\{(a_1,a_2) \mid a_1 \in A_1, a_2 \in A_2 \}$$ is an $s$-basis for $G_1 \times G_2$; therefore, $$\phi (G_1 \times G_2,[0,s]) \leq \phi (G_1,[0,s]) \cdot \phi (G_2,[0,s]).$$ So, from Proposition \ref{phiuppercyclic} we get the following corollary:

\begin{prop} \label{phiuppergeneral}
Let $G$ be an abelian group of rank $r$ and order $n$, and let $s$ be a positive integer.   We then have $$\phi (G,[0,s])  < s^r \cdot \sqrt[s]{n} .$$

\end{prop} 

A different---and in most cases much lower---upper bound was provided by Jia:\index{Jia, X.}

\begin{thm} [Jia; cf.~\cite{Jia:1990a}] \label{phiuppergeneral Jia}\index{Jia, X.}
Let $G$ be an abelian group of order $n$, and let $s$ be a positive integer.   We then have $$\phi (G,[0,s])  < s \cdot \left(1+ \frac{1}{\sqrt[s]{2}} \right)^{s-1}\cdot \sqrt[s]{n}.$$

\end{thm}

Let us now examine the case of cyclic groups in more detail.
Krasny\index{Krasny, O.} (see \cite{Kras:2010a}) developed a computer program that yields the following data. 

$ \phi (\mathbb{Z}_n,[0,2])= \left\{
\begin{array}{ll}
1 & \mbox{if $n=1, 2, 3;$}\\
2 & \mbox{if $n=4,5;$}\\
3 & \mbox{if $n=6, \dots, 9;$}\\
4 & \mbox{if $n=10, \dots, 13;$}\\
5 & \mbox{if $n=14, \dots, 17$, and $n=19;$}\\
6 & \mbox{if $n=18, 20, 21;$}\\
7 & \mbox{if $n=22, \dots, 27$, $n=29,30;$}\\
8 & \mbox{if $n=28$, $n=31, \dots, 35;$}\\
\end{array}\right.$

$ \phi (\mathbb{Z}_n,[0,3])= \left\{
\begin{array}{ll}
1 & \mbox{if $n=1, \dots, 4;$}\\
2 & \mbox{if $n=5, \dots, 8;$}\\
3 & \mbox{if $n=9, \dots, 16;$}\\
4 & \mbox{if $n=17, \dots, 25;$}\\
5 & \mbox{if $n=26, \dots, 40;$}\\
\end{array}\right.$

and 

$ \phi (\mathbb{Z}_n,[0,4]) = \left\{
\begin{array}{ll}
1 & \mbox{if $n=1, \dots, 5;$}\\
2 & \mbox{if $n=6, \dots, 11;$}\\
3 & \mbox{if $n=12, \dots, 27;$}\\
4 & \mbox{if $n=28, \dots, 49.$}\\
\end{array}\right.$

As these values suggest, $\phi (\mathbb{Z}_n,[0,s])$ behaves rather peculiarly; the following problem is wide open. 

\begin{prob}
Find $\phi (\mathbb{Z}_n,[0,s])$ for all $n$ and $s$.
\end{prob}

Observe that the coefficient of $\sqrt[s]{n}$ in the lower bound of Proposition \ref{philowergeneral} and the upper bound of Proposition \ref{phiuppercyclic} are $\sqrt[s]{s!}$ and $s$, respectively; therefore, even the following problems are of much interest.

\begin{prob} \label{c1forphi}
For a given $s \geq 2$, find a constant $c_1 > \sqrt[s]{s!}$ so that $$\phi (\mathbb{Z}_n,[0,s]) > c_1 \cdot \sqrt[s]{n}$$ holds for all $n$ (or, perhaps, all but finitely many $n$), or prove that such a constant cannot exist.
\end{prob}

\begin{prob} \label{c2forphi} 
For a given $s \geq 2$, find a constant $c_2 < s$ so that $$\phi (\mathbb{Z}_n,[0,s]) < c_2 \cdot \sqrt[s]{n}$$ holds for all $n$ (or, perhaps, all but finitely many $n$), or prove that such a constant cannot exist.
\end{prob}

For example, for $s=2$, Problems \ref{c1forphi} and \ref{c2forphi} ask for constants $c_1> \sqrt{2}$ and $c_2<2$, if they exist, for which 
$$c_1 \cdot \sqrt{n} < \phi (\mathbb{Z}_n,[0,2]) < c_2 \cdot \sqrt{n}$$ holds for all (but finitely many) $n$.   While there has been no progress toward finding such a $c_1$, there have been several results for $c_2$; after a series of papers by Fried\index{Fried, K.} in \cite{Fri:1988a}, Mrose\index{Mrose, A.} in \cite{Mro:1979a}, and Kohonen\index{Kohonen, J.} in \cite{Koh:2017a}, the best such result thus far was found by Jia and Shen\index{Jia, X.}\index{Shen, J.} who in \cite{JiaShe:2017a}  proved that one can take $c_2$ to be any real number larger than $\sqrt{3}$.  As a special case of Problem \ref{c2forphi}, we ask for an improvement of this construction:

\begin{prob} \label{c2forphiimp} 
Find a constant $c_2 < \sqrt{3}$ so that $$\phi (\mathbb{Z}_n,[0,2]) < c_2 \cdot \sqrt{n}$$ holds for all (but finitely many) $n$, or prove that such a constant cannot exist.
\end{prob}

It is also interesting to investigate $s$-bases from the opposite viewpoint: given positive integers $s$ and $m$, what are the possible groups $G$ for which $\phi (G,[0,s]) =m$?

For $m=1$ the answer is clear: 

\begin{prop}
Let $s$ be a positive integer and $G$ be an abelian group of order $n$.  Then $G$ contains an $s$-basis of size 1 if, and only if, $G$ is cyclic and $n \leq s+1$.

\end{prop}

For $m \geq 2$, the answer is not known, so we pose the following problem.

\begin{prob} \label{allGform1}
For each $s \geq 2$ and $m \geq 2$, find all groups $G$ for which $\phi(G,[0,s])=m$.
\end{prob}

As a special case of Problem \ref{allGform1}, we have

\begin{prob} \label{allnform}
For each $s \geq 2$ and $m \geq 2$, find all values of $n$ for which $\phi(\mathbb{Z}_{n},[0,s])=m$.
\end{prob}

For $m=2$, Maturo\index{Maturo, A.}  (see  \cite{Mat:2009a}) used a computer program to provide the following partial answer to Problem \ref{allnform}:
$$\begin{array}{|r|l|} \hline
s & \mbox{all} \; n \; \mbox{for which} \; \phi (\mathbb{Z}_n,[0,s])=2 \\ \hline
2 & 4, 5 \\
3 & 5, \dots, 8 \\
4 & 6, \dots, 11 \\
5 & 7, \dots, 16 \\
6 & 8, \dots, 19, 21 \\
7 & 9, \dots, 24, 26 \\
8 & 10, \dots, 31, 33 \\
9 & 11, \dots, 40 \\
10 & 12, \dots, 45, 47 \\
11 & 13, \dots, 52, 54,55,56 \\
12 & 14, \dots, 61, 63, 65 \\
13 & 15, \dots, 66, 68, \dots, 72,74 \\
14 & 16, \dots, 81, 84, 85 \\
15 & 17, \dots, 88, 90, 91, 92, 94, 95, 96 \\
16 & 18, \dots, 98, 100, 101, 102, 105, 107 \\
17 & 19, \dots, 120 \\
18 & 20, \dots, 121, 123, 124, 126, 127, 129, 131, 133 \\
19 & 21, \dots, 136, 138, 139, 140, 143, 144, 146 \\
20 & 22, \dots, 155, 157, 159,161 \\
21 & 23, \dots, 162, 164, \dots, 172, 174, 175, 176 \\
22 & 24, \dots, 180, 182, \dots, 185, 189, 191 \\
23 & 25, \dots, 196, 198, \dots, 208 \\ \hline
\end{array}$$

The data make Problem \ref{allnform} seem challenging even for $m=2$: 

\begin{prob} \label{allnform=2}
For each $s \in \mathbb{N}$, find all values of $n$ for which $\phi(\mathbb{Z}_{n},[0,s])=2$.
\end{prob}

We can formulate two sub-problems of Problem \ref{allnform}:

\begin{prob}  \label{allG cyclic form1}

For each $s \geq 2$ and $m \geq 2$, find the largest integer $f(m,s)$ for which $\phi(\mathbb{Z}_{n},[0,s]) \leq m$ holds for $n=f(m,s)$.
\end{prob}

\begin{prob}

For each $s \geq 2$  and $m \geq 2$, find the largest integer $g(m,s)$ for which $\phi(\mathbb{Z}_{n},[0,s]) \leq m$ holds for all $n \leq g(m,s)$.
\end{prob}

For example, from the table above we see that
$f(2,16)=107$, $g(2,16)=98$, and $f(2,17)=g(2,17)=120$.

Regarding $f(2,s)$, it is not hard to prove that $$f(2,s) \geq \left\lfloor \tfrac{s^2+4s+3}{3} \right\rfloor,$$  in particular, one can verify that the set $A=\{1, s+r\}$, where $r$ is the positive remainder of $s$ when divided by 3, is an $s$-basis in $\mathbb{Z}_n$ for $$n=\left\lfloor \tfrac{s^2+4s+3}{3} \right\rfloor.$$  For example, for $s=6$, the set $A=\{1,9\}$ is a 6-basis in $G=\mathbb{Z}_{21}$, since for nonnegative integer coefficients $\lambda_1$ and $\lambda_2$ satisfying $\lambda_1+\lambda_2 \leq 6$, the values of $$\lambda_1 \cdot 1 + \lambda_2 \cdot 9$$ yield all elements of the group:
$$\begin{array}{||l||c|c|c|c|c|c|c||} \hline \hline
  & \lambda_1=0  & \lambda_1=1  & \lambda_1=2 & \lambda_1=3 &  \lambda_1=4 & \lambda_1=5  & \lambda_1=6 \\ \hline \hline 
\lambda_2=0 &  0 & 1 & 2  &  3 & 4 & 5 & 6 \\ \hline 
\lambda_2=1 &  9 & 10 &  11 &  12 &  13 & 14 &  \\ \hline 
\lambda_2=2 &  18 & 19 & 20 & 0 & 1 & &  \\ \hline 
\lambda_2=3 & 6 & 7 & 8 & 9 & & &  \\ \hline 
\lambda_2=4 & 15 & 16 & 17 & & &  & \\ \hline 
\lambda_2=5 & 3  & 4 &  & & & &  \\ \hline 
\lambda_2=6 & 12  &  &  & & & &  \\ \hline \hline
\end{array}$$    

In fact, Morillo, Fiol, and F\`abrega\index{Fiol, M. A.}\index{Morillo, P.}\index{F\`abrega, J.} proved that $f(2,s)$ cannot be larger (this was independently re-proved by Hsu and Jia\index{Hsu, D. F.}\index{Jia, X.} in \cite{HsuJia:1994a} and conjectured by Maturo in \cite{Mat:2009a}):\index{Maturo, A.} 

\begin{thm} [Morillo, Fiol, and F\`abrega; cf.~\cite{MorFioFab:1985a}] \label{thmnform=2}\index{Fiol, M. A.}\index{Morillo, P.}\index{F\`abrega, J.}

Let $s$ be a positive integer.  The largest possible value $f(2,s)$ of $n$ for which $\phi(\mathbb{Z}_{n},[0,s]) \leq 2$ is $$f(2,s)=\left\lfloor \tfrac{s^2+4s+3}{3} \right\rfloor.$$

\end{thm}  

Exact values for $f(m,s)$ are not known in general when $m \geq 3$, though we have computational data for small $s$ and $m=3$ and for small $m$ and $s=2$.  Namely, in the paper \cite{HsuJia:1994a} of Hsu and Jia we see:\index{Hsu, D. F.}\index{Jia, X.}
$$\begin{array}{||c||c|c|c|c|c|c|c|c|c|c|c|c|c|c||} \hline \hline
s & 2 & 3 & 4 & 5 & 6 & 7 & 8 & 9 & 10 & 11 & 12 & 13 & 14 & 15 \\ \hline
f(3,s) & 9 & 16 & 27 & 40 & 57 & 78 & 111 & 138 & 176 & 217 & 273 & 340 & 395 & 462 \\ \hline \hline
\end{array}$$
In \cite{Haa:2004a}, Haanp\"a\"a presents:\index{Haanp\"a\"a, H.} 
$$\begin{array}{||c||c|c|c|c|c|c|c|c|c|c|c|c||} \hline \hline
m & 1& 2 & 3 & 4 & 5 & 6 & 7 & 8 & 9 & 10 & 11 & 12  \\ \hline
f(m,2) &3 & 5 & 9 & 13 & 19 & 21 & 30 & 35 & 43 & 51 & 63 & 67  \\ \hline \hline
\end{array}$$
(Some of these values had been determined earlier by Graham and Sloane in \cite{GraSlo:1980a}.)\index{Graham, R. L.}\index{Sloane, N. J. A.} 

More ambitiously than Problem \ref{allG cyclic form1}, but more modestly than Problem \ref{allGform1}, we may ask the following:

\begin{prob} \label{allGform11}
For each $s \geq 2$ and $m \geq 2$, find the largest integer $F(m,s)$ so that there is a group $G$ of order $n=F(m,s)$ for which $\phi(G,[0,s]) \leq m$.
\end{prob} 

Since obviously $F(m,s) \geq f(m,s)$, we already know from Theorem \ref{thmnform=2} that
$$F(2,s) \geq  \left\lfloor \tfrac{s^2+4s+3}{3} \right\rfloor.$$
In \cite{FioYebAleVal:1987a}, Fiol, Yebra, Alegre, and Valero\index{Fiol, M. A.}\index{Yebra, J. L. A.}\index{Alegre, I.}\index{Valero, M.}  proved that, when $s \not \equiv 1$ mod 3, then we cannot do better.  However, when $s \equiv 1$ mod 3, then the group $G=\mathbb{Z}_k \times \mathbb{Z}_{3k}$ has the $s$-basis $\{(0,1),(1,3k-1)\}$, where $k=(s+2)/3$.   Note that this group has order $\tfrac{s^2+4s+4}{3}$, and if $s \not \equiv 1$ mod 3, then $$\left\lfloor \tfrac{s^2+4s+3}{3} \right\rfloor=\left\lfloor \tfrac{s^2+4s+4}{3} \right\rfloor.$$  In summary, we have:

\begin{thm} [Fiol; cf.~\cite{Fio:2013a}]\index{Fiol, M. A.}
Let $s$ be any positive integer.  The largest value $F(2,s)$ for which there is a group $G$ of order $n=F(m,s)$ for which $\phi(G,[0,s]) \leq 2$ is
$$F(2,s)=\left \lfloor \tfrac{s^2+4s+4}{3} \right \rfloor.$$

\end{thm}

Exact values for $F(m,s)$ are not known in general when $m \geq 3$, but in \cite{Haa:2004a}, Haanp\"a\"a\index{Haanp\"a\"a, H.} computed the values of $F(m,2)$ for $m \leq 12$.  He found that $F(m,2)=f(m,2)$ for all $m \leq 12$, except as follows:
\begin{itemize}
\item $F(8,2)=36=f(8,2)+1$, as shown by the 2-basis 
$$\{(1, 0, 1), (0, 0, 1), (0, 0, 2), (1, 1, 0),
(1, 2, 0), (3, 0, 2), (3, 1, 0), (3, 2, 0)\}$$ in $G=\mathbb{Z}_4 \times \mathbb{Z}_3^2$;
\item $F(11,2)=64=f(11,2)+1$, as shown by the 2-basis
$$\{(0, 1), (0, 4), (1, 0), (1, 2), (2, 1),
(2, 2), (2, 6), (4, 5), (5, 0), (5, 2), (6, 5)\}$$ in $G=\mathbb{Z}_8^2$; and
\item $F(12,2)=72=f(12,2)+5$, as shown by the 2-basis
 $$\{(0, 1, 1), (0, 0, 2), (0, 2, 1), (0, 2, 4), (0, 2, 7), (0, 3, 1),$$ 
$$(1, 0, 3), (1, 0, 8), (1, 1, 1),
(1, 2, 5), (1, 2, 6), (1, 3, 1)\}$$ in $G=\mathbb{Z}_2 \times \mathbb{Z}_4 \times \mathbb{Z}_9$.
\end{itemize}

Turning now to $g(m,s)$, first note that, obviously, $g(m,s) \leq f(m,s)$.  We can find a lower bound for $g(m,s)$ as follows.  \label{lower bound for g}
Observe that the base $a$ representation of nonnegative integers guarantees that the set $$A=\{a^j \mid j=0,1,\dots,m-1\}$$ is an $s$-basis for $\mathbb{Z}_n$ for all $n \leq a^m$, as long as $s \geq m \cdot (a-1)$.  Thus, with $$a=\left \lfloor \frac{s}{m} \right \rfloor+1,$$ we get the following.

\begin{prop}
Let $m$ and $s$ be a positive integers.  The largest possible value $g(m,s)$ of $n$ for which $\phi(\mathbb{Z}_{n},[0,s]) \leq m$ holds for all $n \leq g(m,s)$ satisfies
$$g(m,s) \geq \left( \left \lfloor \frac{s}{m} \right \rfloor +1 \right)^m.$$

\end{prop}  

In fact, for $g(2,s)$ we can do slightly better:   

\begin{prop} \label{nform=2}  

Let $s$ be a positive integer.  The largest possible value $g(2,s)$ of $n$ for which $\phi(\mathbb{Z}_{n},[0,s]) \leq 2$ holds for all $n \leq g(2,s)$ satisfies
$$\left \lfloor \tfrac{s^2+6s+5}{4} \right \rfloor \leq g(2,s) \leq \left\lfloor \tfrac{s^2+4s+3}{3} \right\rfloor.$$

\end{prop}

The upper bound follows directly from Theorem \ref{thmnform=2} above.  To prove the lower bound, we can verify that for each $n$ up to this value, the set $$A=\left\{1, \left\lfloor \tfrac{s+3}{2} \right\rfloor \right\}$$ is an $s$-basis in $\mathbb{Z}_n$---the details can be found on page \pageref{proofofpropnform=2}.  Note that for $s \in \{2,3,4\}$, the values of $\lfloor \tfrac{s^2+6s+5}{4} \rfloor$ and $\lfloor \tfrac{s^2+4s+3}{3} \rfloor$ agree, thus, for these values, $g(2,s)$ is determined and Problem \ref{allnform=2} is answered.  However, we still have the following:

\begin{prob}

For each $s \geq 2$, find the largest integer $g(2,s)$ for which $\phi(\mathbb{Z}_{n},[0,s]) \leq 2$ holds for all $n \leq g(2,s)$.
\end{prob}

\subsection{Arbitrary number of terms} \label{3minUarbitrary}

\section{Unrestricted signed sumsets} \label{3minUS}

In this section we investigate the quantity $$\phi_{\pm}(G,H) =\mathrm{min} \{ |A|  \mid A \subseteq G, H_{\pm}A =G\};$$ if no $H$-spanning set exists, we put $\phi_{\pm}(G,H) = \infty.$  Note that we have $$\phi_{\pm}(G,\{0\})=\infty$$ for all groups of order at least 2, but if $H$ contains at least one positive integer $h$, then $$\phi_{\pm}(G,H) \leq n$$ since we have $$G=\{(h-1) \cdot 0 + 1 \cdot g  \mid g \in G\}  \subseteq h_{\pm}G \subseteq H_{\pm}G,$$ so $
H_{\pm}G=G$.

\subsection{Fixed number of terms} \label{3minUSfixed}

A subset $A$ of $G$ for which $h_{\pm}A =G$ for some nonnegative integer $h$ is called an {\em exact $h$-spanning set} for $G$.  Here we investigate $$\phi_{\pm} (G,h) = \mathrm{min} \{ |A|  \mid A \subseteq G, h_{\pm}A =G\},$$ that is, the minimum size of an exact $h$-spanning set for $G$.  

As we pointed out above, we have $$\phi_{\pm} (G,0) =\infty$$ for all groups of order at least 2, but $$\phi_{\pm}(G,h) \leq n$$ for every $h \in \mathbb{N}$.

For $h=1$, it is clear that $A$ is an exact $1$-spanning set if, and only if, for each $g \in G$, $A$ contains at least one of $g$ or $-g$; in particular, $A$ must contain 0, every element of order 2, and half of the elements of order more than 2.  Therefore, we have 

\begin{prop}  \label{phi pm US h=1}
For any finite abelian group $G$ we have $$\phi_{\pm} (G,1)=\frac{n+|\mathrm{Ord}(G,2)|+1}{2};$$
in particular,
$$   \phi_{\pm} (\mathbb{Z}_n,1) =\lfloor (n+2)/2 \rfloor.$$

\end{prop}

But for $h \geq 2$, we do not have a formula for $\phi_{\pm}(G,h)$.  Therefore we pose the following general problem.

\begin{prob} \label{Findphipm(G,h)}
Find $\phi_{\pm}(G,h)$ for every $G$ and $h \geq 2$.
\end{prob} 

To find a lower bound for $\phi_{\pm}(G,h)$, note that an $m$-subset of $G$ may have an $h$-fold signed sumset of size at most $$c(h,m)=\sum_{i \geq 0} {m \choose i} {h-1 \choose i-1} 2^i;$$ if the $h$-fold signed sumset is in fact $G$, then the elements that are their own inverses are generated twice, since replacing all coefficients in the linear combination by their negatives yields the same element.  Therefore, we have the following general bound: 

\begin{prop} \label{philowerUSfixed}  If $A$ is an exact $h$-spanning set of size $m$ in $G$, then
$$n \leq c(h,m)-|\mathrm{Ord}(G,2)|-1;$$ in particular,
for the cyclic group $\mathbb{Z}_n$ we must have
$$n \leq c(h,m)-1.$$
\end{prop}
(Note that $c(h,m)$ is even for all $n, h \in \mathbb{N}$.)

Proposition \ref{philowerUSfixed} provides us with a lower bound for $\phi_{\pm} (G,h)$.  For $h=1$, it yields $$\phi_{\pm} (G,1) \geq \frac{n+|\mathrm{Ord}(G,2)|+1}{2};$$ in fact, by Proposition \ref{phi pm US h=1}, we know that equality holds.  For $h=2$, we have:   

\begin{prop} \label{phipm lower h=2}

For all abelian groups $G$ of order $n$ we have $$\phi_{\pm} (G,2) \geq \sqrt{\frac{n+|\mathrm{Ord}(G,2)|+1}{2}};$$
in particular,
$$   \phi_{\pm} (\mathbb{Z}_n,2) \geq \sqrt{\lfloor (n+2)/2 \rfloor}.$$

\end{prop}
For $h \geq 3$ this bound is difficult to state explicitly.

To find some upper bounds for $\phi_{\pm} (G,h)$, we exhibit an explicit exact $h$-spanning set, as follows.  We will only consider the cyclic group $\mathbb{Z}_n$ here (for noncyclic groups we may use the method we discussed on page \pageref{from cyclic to gen}).  Let $a=\left \lceil \sqrt[h]{n+1} \right \rceil$.  We can readily verify that the set
$$A=\{0\} \cup \{ i \cdot a^j \mid i=1,2,\dots,\left \lfloor a/2 \right \rfloor; j=0,1,\dots,h-1\}$$ is an exact $h$-spanning set for $\mathbb{Z}_n$: indeed, all integers with absolute value at most $n/2$ are generated:
$$\left \lfloor \frac{a}{2} \right \rfloor \cdot (1+a+\cdots+a^{h-1})=\left \lfloor \frac{a}{2} \right \rfloor \cdot \frac{a^h-1}{a-1} \geq \frac{a-1}{2} \cdot \frac{a^h-1}{a-1}  = \frac{a^h-1}{2} \geq \frac{n}{2} $$ are generated.  Since 
$$|A|=1+\left \lfloor a/2 \right \rfloor \cdot h = 1+ h \cdot \left \lfloor \left \lceil \sqrt[h]{n+1} \right\rceil /2 \right \rfloor \leq 1+ h \cdot \left \lfloor ( \sqrt[h]{n}+1) /2 \right \rfloor \leq 1+ h \cdot (\sqrt[h]{n} +1)/2,$$ we get the following result.

\begin{prop} \label{phipmuppercyclic h}
For positive integers $n$ and $h$ we have $$\phi_{\pm}  (\mathbb{Z}_n,h) \leq h/2 \cdot \sqrt[h]{n} +(h+2)/2.$$

\end{prop} 

  The following problem is wide open.

\begin{prob}
For positive integers $n$ and $h \geq 2$, find $\phi_{\pm} (\mathbb{Z}_n,h)$.
\end{prob}

For $h=2$, from Propositions \ref{phipm lower h=2} and \ref{phipmuppercyclic h} we get 
$$ \sqrt{\lfloor (n+2)/2 \rfloor}  \leq \phi_{\pm} (\mathbb{Z}_n,2) \leq \sqrt{n} +2.$$
We offer the following interesting problem:

\begin{prob}
Find constants $c_1 > \frac{1}{\sqrt{2}}$ and $c_2<1$ so that for all (but perhaps finitely many) positive integers $n$ we have $$ c_1 \cdot  \sqrt{n} \leq \phi_{\pm} (\mathbb{Z}_n,2) \leq c_2 \cdot  \sqrt{n},$$ or prove that no such constants exist.
\end{prob}

It is also interesting to investigate exact $h$-spanning sets from the opposite viewpoint: given positive integers $h$ and $m$, what are the possible groups $G$ for which $\phi_{\pm} (G,h) =m$?  

The answer for $m=1$ follows immediately from Proposition \ref{philowerUSfixed}: since $c(h,1)=2$, the inequality, and thus $\phi_{\pm} (G,h) =1$, holds if, and only if, $n=1$.  For $m=2$, Reckner in \cite{Rec:2013a} has the following results:

\begin{prop} [Reckner; cf.~\cite{Rec:2013a}]  \label{Reckner m=2}\index{Reckner, T.}  Let $n$ and $h$ be positive integers.  
\begin{enumerate}
\item If $$2 \leq n \leq 2h+1,$$ then $\phi_{\pm} (\mathbb{Z}_n,h) =2$.  In particular, the set $\{0,1\}$ is an exact $h$-spanning set in $\mathbb{Z}_n$ for all $n$ in the given range.
\item If $n$ and $h$ are both odd and $$2h+3 \leq n \leq 3h,$$ then $\phi_{\pm} (\mathbb{Z}_n,h) =2$.    In particular, when $h$ is odd, the set $\{1,3\}$ is an exact $h$-spanning set in $\mathbb{Z}_n$ for all odd $n$ in the given range.
\end{enumerate}
\end{prop}

Proposition \ref{Reckner m=2} does not classify all $n$ and $h$ with $\phi_{\pm} (\mathbb{Z}_n,h) =2$, so we offer the following problem:

\begin{prob}   \label{other than Reckner m=2}
Find all other values of $n$ and $h$ for which $\phi_{\pm} (\mathbb{Z}_n,h) =2$.
\end{prob}

Note that, by Proposition \ref{philowerUSfixed}, it suffices to investigate Problem \ref{other than Reckner m=2} when $n \leq 4h-1$.

Moving on to $m=3$, we observe that, as an immediate consequence of Proposition \ref{p} in the next subsection, we have:

\begin{prop} \label{Reckner m=3}
Let $n$ and $h$ be positive integers.  If $$n \leq 2h^2+2h+1,$$ then $\phi_{\pm} (\mathbb{Z}_n,h) \leq 3$.  In particular, the set $\{0,h, h+1\}$ is an exact $h$-spanning set in $\mathbb{Z}_n$ for all $n \leq 2h^2+2h+1$.

\end{prop}

\begin{prob}   \label{other than Reckner m=3}
Find all other values of $n$ and $h$ for which  $\phi_{\pm} (\mathbb{Z}_n,h) \leq 3$.
\end{prob}

Note that, by Proposition \ref{philowerUSfixed}, it suffices to investigate Problem \ref{other than Reckner m=3} when $n \leq 4h^2+1$.

\begin{prob}
For each $h$ and $m \geq 4$, find all $n$ for which $\phi_{\pm} (\mathbb{Z}_n,h) = m$.
\end{prob}
And, more generally:
\begin{prob}
For all positive integers $h$ and $m$, find all groups $G$ for which $\phi_{\pm} (G,h) = m$.
\end{prob}

\subsection{Limited number of terms} \label{3minUSlimited}

A subset $A$ of $G$ for which $[0,s]_{\pm}A =G$ for some nonnegative integer $s$ is called an {\em $s$-spanning set} for $G$.  Here we investigate $$\phi_{\pm} (G,[0,s]) = \mathrm{min} \{ |A|  \mid A \subseteq G, [0,s]_{\pm}A =G\},$$ that is, the minimum size of an $s$-spanning set for $G$.  

As we pointed out above, we have $$\phi_{\pm} (G,[0,0]) =\infty$$ for all groups of order at least 2, but $$\phi_{\pm}(G,[0,s]) \leq n$$ for every $s \in \mathbb{N}$.

For $s=1$, it is clear that $A$ is $1$-spanning if, and only if, for each $g \in G$, $A$ contains at least one of $g$ or $-g$; in particular, $A$ must contain every element of order 2 and half of the elements of order more than 2.  Therefore, we have 

\begin{prop}  \label{1-spanning}
For any finite abelian group $G$ we have $$\phi_{\pm} (G,[0,1])=\frac{n+|\mathrm{Ord}(G,2)|-1}{2};$$
in particular,
$$   \phi_{\pm} (\mathbb{Z}_n,[0,1]) =\lfloor n/2 \rfloor.$$

\end{prop}

But for $s \geq 2$, we do not have a formula for $\phi_{\pm}(G,[0,s])$.  Therefore we pose the following general problem.

\begin{prob} \label{Findphipm(G,[0,s])}
Find $\phi_{\pm}(G,[0,s])$ for every $G$ and $s \geq 2$.
\end{prob} 

As a consequence of Proposition \ref{philower}, we have the following general bound.

\begin{prop} \label{philowerUSlimited}  If $A$ is an $s$-spanning set of size $m$ in $G$, then
$$n \leq a(m,s)=\sum_{i \geq 0} {m \choose i} {s \choose i} 2^i.$$
\end{prop}

The classification of the extremal cases of Proposition \ref{philowerUSlimited} is a particularly intriguing question.  Namely, we would like to find {\em perfect $s$-spanning sets}, that is, $s$-spanning sets of size $m$ with $n=a(m,s).$  The values of $a(m,s)$ can be tabulated for small values of $m$ and $s$ (see also Section \ref{0.2.3}):

$$\begin{array}{c||c|c|c|c|c|c|c|}
$a(m,s)$ & $s=0$ & $s=1$  & $s=2$  & $s=3$  & $s=4$  & $s=5$  & $s=6$ \\ \hline \hline
$m=1$ & {\bf 1} & {\bf 3}& {\bf 5}& {\bf 7}& {\bf 9}& {\bf 11}& {\bf 13} \\ \hline
$m=2$ & 1 & {\bf 5}& {\bf 13}& {\bf 25}& {\bf 41}& {\bf 61}& {\bf 85} \\ \hline
$m=3$ & 1 & {\bf 7}& 25& 63& 129& 231& 377 \\ \hline
$m=4$ & 1 & {\bf 9}& 41& 129& 321& 681& 1289 \\ \hline
$m=5$ & 1 & {\bf 11}& 61& 231& 681& 1683& 3653 \\ \hline
$m=6$ & 1 & {\bf 13}& 85& 377& 1289& 3653& 8989 \\ \hline
\end{array}$$

Cases where there exists a group of size $a(m,s)$ with a known perfect $s$-spanning set are marked with bold-face.  The following proposition exhibits perfect spanning sets for all known parameters.

\begin{prop}  \label{perfex}

Let $m$ be a positive integer and $s$ be a nonnegative integer, and let $G$ be an abelian group of order $n$.

\begin{enumerate}

\item If $n=2m+1$, then $G \setminus \{0\}$ can be partitioned into parts $K$ and $-K$, and both $K$ and $-K$ are perfect 1-spanning sets in $G$.  For example, the set $\{1,2,\dots,m\}$ is a perfect 1-spanning set in $\mathbb{Z}_n$.

\item If $n=2s+1$ and $\gcd(a,n)=1$, then the set $\{a\}$ is a perfect $s$-spanning set in $\mathbb{Z}_n$.

\item If $n=2s^2+2s+1$, then the sets $\{1, 2s+1\}$ and $\{s,s+1\}$ are perfect $s$-spanning sets in $\mathbb{Z}_n$.

\end{enumerate}

\end{prop}

The first two statements are obvious.  The fact that the set $\{1, 2s+1\}$ is a perfect $s$-spanning set in $\mathbb{Z}_n$ is provided on page \pageref{proofofperfex}, and the same claim for $\{s,s+1\}$ follows from Proposition \ref{p} below; note, however, that $\{s,s+1\}$ in Proposition \ref{p} cannot be replaced by $\{1,2s+1\}$.  (Both claims were  demonstrated for $s=3$ and $n=25$ in the Appetizer section ``In pursuit of perfection'')  While these two perfect spanning sets of size 2 are known, there has not been a characterization of any others.

\begin{prob}
Find all perfect $s$-spanning sets of size 2 in the cyclic group $\mathbb{Z}_{2s^2+2s+1}$.
\end{prob}   

We could not find perfect spanning sets for $s \geq 2$ and $m \geq 3$ for any $n$, and neither could we find any noncyclic groups with perfect spanning sets for $s \geq 2$ and $m=2$.  In particular, we definitely know that, other than the ones already mentioned, no perfect spanning sets exist in groups of order up to 100: Laza\index{Laza, N.}  (cf.~\cite{Laz:?}) has verified (using a computer program) that no perfect 2-spanning sets exist in $\mathbb{Z}_{25}$, $\mathbb{Z}_{41}$, $\mathbb{Z}_{61}$, or $\mathbb{Z}_{85}$ (note that $a(m,2)=25, 41, 61,$ and 85 for $m=3, 4, 5,$ and 6, respectively) and that no perfect 3-spanning set exists in $\mathbb{Z}_{63}$ (we have $a(3,3)=63$); on page \pageref{no perfects in Z25}, we presented simple arguments to prove that there is neither a perfect 3-spanning set of size two, nor a perfect 2-spanning set of size three in $\mathbb{Z}_5^2$; and Jankowski\index{Jankowski, T.} in \cite{Jan:2012a} proved that $\mathbb{Z}_{3} \times \mathbb{Z}_{21}$ cannot contain a perfect 3-spanning set of size 3 either. 

It might be an interesting problem to find and classify all perfect spanning sets.

\begin{prob}

Find perfect $s$-spanning sets in $G=\mathbb{Z}_{n}$ of size $m$ for some values $s \geq 2$ and $m \geq 3$, or prove that such perfect spanning sets do not exist.

\end{prob}

\begin{prob}

Find perfect $s$-spanning sets in noncyclic groups of size $m$ for some values $s \geq 2$ and $m \geq 2$, or prove that such perfect spanning sets do not exist.

\end{prob}

Let us return to the question of finding $\phi_{\pm} (G,[0,s])$.  Proposition \ref{philowerUSlimited} provides us with a lower bound for $\phi_{\pm} (G,[0,s])$, though for $s \geq 3$ this bound is difficult to state explicitly.  For $s=2$ we get:

\begin{prop} \label{phipm lower s=2}

For all abelian groups $G$ of order $n$ we have $$\phi_{\pm} (G,[0,2]) \geq \frac{\sqrt{2n-1}-1}{2}.$$

\end{prop}
  
As in the subsection \ref{3minUSfixed}, we can find some upper bounds for $\phi_{\pm} (\mathbb{Z}_n,[0,s])$ by letting $a=\lceil \sqrt[s]{n+1} \rceil$ and verifying that the set
$$A=\{ i \cdot a^j \mid i=1,2,\dots,\left \lfloor \frac{a}{2} \right \rfloor; j=0,1,\dots,s-1\}$$ is an $s$-spanning set for $\mathbb{Z}_n$.  Since
$$|A|=s \cdot \left \lfloor a/2 \right \rfloor  = s \cdot \left \lfloor \left \lceil \sqrt[s]{n+1} \right\rceil /2 \right \rfloor \leq s \cdot \left \lfloor ( \sqrt[s]{n}+1) /2 \right \rfloor \leq s \cdot (\sqrt[s]{n} +1)/2,$$ we get the following result.

\begin{prop} \label{phipmuppercyclic}
For positive integers $n$ and $s$ we have $$\phi_{\pm}  (\mathbb{Z}_n,[0,s]) \leq s/2 \cdot \sqrt[s]{n} +s/2.$$

\end{prop} 

The computational data of Laza\index{Laza, N.}  (see \cite{Laz:?}) shows that  

$ \phi_{\pm} (\mathbb{Z}_n,[0,2])= \left\{
\begin{array}{ll}
1 & \mbox{if $n=1, 2, 3, 4, {\bf 5};$}\\
2 & \mbox{if $n=6, 7, \dots, 12, {\bf 13};$}\\
3 & \mbox{if $n=14, 15, \dots, 21;$}\\
4 & \mbox{if $n=22, 23, \dots, 33$, and $n=35;$}\\
5 & \mbox{if $n=34$, $n=36, 37, \dots, 49$, and $n=51;$}\\
\end{array}\right.$

and 

$ \phi_{\pm} (\mathbb{Z}_n,[0,3]) = \left\{
\begin{array}{ll}
1 & \mbox{if $n=1,2, \dots, 6, {\bf 7};$}\\
2 & \mbox{if $n=8, 9 \dots, 24, {\bf 25};$}\\
3 & \mbox{if $n=26, 27, \dots, 50$, $n=52$, and $n=55;$}\\
4 & \mbox{if $n=51,53, 54$, $n=56, 57, \dots, 100$, and $n=104.$}\\
\end{array}\right.$

\noindent (Values marked in bold-face indicate perfect spanning sets, as explained above.)  The following problem is wide open.

\begin{prob}
For positive integers $n$ and $s \geq 2$, find $\phi_{\pm} (\mathbb{Z}_n,[0,s])$.
\end{prob}

For $s=2$, from Propositions \ref{phipm lower s=2} and \ref{phipmuppercyclic} we get 
$$ \frac{\sqrt{2n-1}-1}{2}  \leq \phi_{\pm} (\mathbb{Z}_n,[0,2]) \leq \sqrt{n} +1.$$
We offer the following interesting problem:

\begin{prob}
Find constants $c_1 > \frac{1}{\sqrt{2}}$ and $c_2<1$ so that for all (but perhaps finitely many) positive integers $n$ we have $$ c_1 \cdot  \sqrt{n} \leq \phi_{\pm} (\mathbb{Z}_n,[0,2]) \leq c_2 \cdot  \sqrt{n},$$ or prove that no such constants exist.
\end{prob}

It is also interesting to investigate $s$-spanning sets from the opposite viewpoint: given positive integers $s$ and $m$, what are the possible groups $G$ for which $\phi_{\pm} (G,[0,s]) =m$?   Here is what we know:

\begin{prop} [Bajnok; cf.~\cite{Baj:2004a}] \label{p}\index{Bajnok, B.} Suppose, as usual, that $G$ is an abelian group of order $n$.  Let $s \geq 1$ be an integer.

1.  We have $ \phi_{\pm} (G,[0,s]) =1$ if, and only if, $G$ is cyclic and  $$1 \leq n \leq 2s+1.$$  In particular, the set $\{1\}$ is $s$-spanning in $\mathbb{Z}_n$ for every $n \leq 2s+1$.

2.  We have $ \phi_{\pm} (\mathbb{Z}_n,[0,s]) =2$  if, and only if, $$2s+2 \leq n \leq 2s^2+2s+1.$$ 
In particular, the set $\{s,s+1\}$ is $s$-spanning in $\mathbb{Z}_n$ for every $n \leq 2s^2+2s+1$.

\end{prop}

The case of $m=1$ in Proposition \ref{p} is clear: for any given $s \in \mathbb{N}$, the only groups that contain an $s$-spanning set of size 1 are the cyclic groups $\mathbb{Z}_n$, and we have $\phi (\mathbb{Z}_n,[0,s]) =1$ if, and only if, $n \leq 2s+1$ (in which case the set $\{1\}$, for example, is $s$-spanning in $\mathbb{Z}_n$).  For $m = 2$, the ``only if'' part follows from the fact that for $n \leq 2s+1$ we have $\phi_{\pm} (\mathbb{Z}_n,[0,s]) =1$ (lower bound) and from Proposition \ref{philowerUSlimited} (upper bound).  For the rest, see page \pageref{proofofp}.

Note that the second part of Theorem \ref{p} treats only cyclic groups.  While it is clear from Proposition \ref{rankvsgen} that the only noncyclic groups $G$ possessing $s$-spanning sets of size two must have rank two, we do not have a characterization of those for which $ \phi_{\pm} (G,[0,s]) =2$.  We have the following obvious conditions:

\begin{prop} \label{nec and suff for m=2 s-spanning}
Let $G=\mathbb{Z}_{n_1} \times \mathbb{Z}_{n_2}$ be of rank two.   
\begin{enumerate}
  \item A necessary condition for  $ \phi_{\pm} (G,[0,s]) =2$ is that $$n_1 \cdot n_2 \leq 2s^2+2s+1.$$
  \item A sufficient condition for  $ \phi_{\pm} (G,[0,s]) =2$ is that  $$\lfloor n_1/2 \rfloor + \lfloor n_2/2 \rfloor \leq s.$$
\end{enumerate}
\end{prop}

The first claim is due to Proposition \ref{philowerUSlimited}, and the second claim follows from the fact that if the stated condition holds, then (for example) the set
$\{(0,1),(1,0)\}$ is $s$-spanning in $G$.

While the two conditions in Proposition \ref{nec and suff for m=2 s-spanning} are generally far away from each other, we can see what conclusions we have for small values of $s$.  For $s=1$,  by the necessary condition, the only possible group with a 1-spanning set of size two is $\mathbb{Z}_2^2$, but (e.g.~by Proposition \ref{1-spanning}) $ \phi_{\pm} (\mathbb{Z}_2^2,[0,1]) =3$.  For $s=2$, the only possible groups with a 2-spanning set of size two are $\mathbb{Z}_2^2$, $\mathbb{Z}_2 \times \mathbb{Z}_4$,  $\mathbb{Z}_3^2$, and $\mathbb{Z}_2 \times \mathbb{Z}_6$; of these, the first three have one (e.g.~the set $\{0,1),(1,1)\}$ works in each), but, as we can easily check, the fourth one does not.  For $s=3$, the necessary condition of Proposition \ref{nec and suff for m=2 s-spanning} yields six other possibilities besides the ones mentioned, which we can then be checked individually.  We summarize our findings as follows:

\begin{cor}
\begin{enumerate}
  \item There is  no group $G$ of rank two with  $ \phi_{\pm} (G,[0,1]) =2$.
\item There are exactly three groups $G$ of rank two with  $ \phi_{\pm} (G,[0,2]) =2$; namely, $\mathbb{Z}_2^2$, $\mathbb{Z}_3^2$, and $\mathbb{Z}_2 \times \mathbb{Z}_4$.
\item There are exactly six groups $G$ of rank two with  $ \phi_{\pm} (G,[0,3]) =2$; namely, $\mathbb{Z}_2 \times \mathbb{Z}_{2k}$ for $k \in \{1,2,3,4\}$, and  $\mathbb{Z}_3\times \mathbb{Z}_{3k}$ for $k \in \{1,2\}$.
\end{enumerate}
\end{cor}

For groups of rank two whose first component is of size two, we have the following more general result:

\begin{prop} \label{Z2 X Z2k s spanning}
Suppose that $k$ and $s$ are positive integers satisfying one of the following conditions:
\begin{itemize}
  \item $s=2$ and $k \in \{1,2\}$;
  \item $s=3$ and $k \in \{1,2,3,4\}$;
  \item $s=4$ and $k \in \{1,2,3,4,5,6,8\}$; or
  \item $s \geq 5$ and $k \leq 3s-4$.
\end{itemize}
Then $ \phi_{\pm} (\mathbb{Z}_2 \times \mathbb{Z}_{2k},[0,s]) =2$.
\end{prop}

The fact that $k=7$ is not listed for $s=4$ is no typo: $ \phi_{\pm} (\mathbb{Z}_2 \times \mathbb{Z}_{14},[0,4]) =3.$  The proof of Proposition \ref{Z2 X Z2k s spanning} starts on page \pageref{proof of Z2 X Z2k s spanning}.  

We also believe that the converse of Proposition \ref{Z2 X Z2k s spanning} holds as well:
\begin{conj} \label{Z2 X Z2k s spanning conj}
All values of $s$ and $k$ for which $ \phi_{\pm} (\mathbb{Z}_2 \times \mathbb{Z}_{2k},[0,s]) =2$ are listed in Proposition \ref{Z2 X Z2k s spanning}.
\end{conj}

We pose the following interesting problems.

\begin{prob}
Prove or disprove Conjecture \ref{Z2 X Z2k s spanning conj}.

\end{prob}

\begin{prob}
For each $s \geq 2$, find all values of $k$ and $l$ for which $ \phi_{\pm} (\mathbb{Z}_3 \times \mathbb{Z}_{3k},[0,s]) =2$ and $ \phi_{\pm} (\mathbb{Z}_4 \times \mathbb{Z}_{4l},[0,s]) =2$.

\end{prob}

\begin{prob}
For each $s \geq 2$, find all values of $k$ for which $ \phi_{\pm} (\mathbb{Z}_{k}^2,[0,s]) =2$.

\end{prob}

\begin{prob} \label{allGnformpm}
For each $s \geq 2$, find all groups $G$ of rank two for which $ \phi_{\pm} (G,[0,s]) =2$.
\end{prob}

\begin{prob} \label{allnform1pm}
For each $s \geq 2$ and $m \geq 3$, find all values of $n$ for which $\phi_{\pm}(\mathbb{Z}_n,[0,s])=m$.
\end{prob}

\begin{prob} \label{allGnformpm}
For each $s \geq 2$ and $m \geq 3$, find all groups $G$ for which $\phi_{\pm}(G,[0,s])=m$.
\end{prob}

We can formulate two sub-problems of Problem \ref{allnform1pm}:

\begin{prob}

For all positive integers $s$ and $m$, find the largest integer $f_{\pm}(m,s)$ for which $\phi_{\pm}(\mathbb{Z}_{n},[0,s]) \leq m$ holds for $n =f_{\pm}(m,s)$.
\end{prob}

\begin{prob}

For all positive integers $s$ and $m$, find the largest integer  $g_{\pm}(m,s)$ for which $\phi_{\pm}(\mathbb{Z}_{n},[0,s]) \leq m$ holds for all $n \leq g_{\pm}(m,s)$.
\end{prob}

Recall that Proposition \ref{philowerUSlimited} provides the following upper bound:
$$ g_{\pm}(m,s) \leq f_{\pm}(m,s) \leq a(m,s)=\sum_{i \geq 0} {m \choose i} {s \choose i} 2^i.$$
According to Proposition \ref{p}, we have
$$g_{\pm}(1,s)=f_{\pm}(1,s)=2s+1=a(1,s),$$ and
$$g_{\pm}(2,s)=f_{\pm}(2,s)=2s^2+2s+1=a(2,s).$$

We can find a general lower bound for $g_{\pm}(m,s)$, and thus for $f_{\pm}(m,s)$, using an analogous argument to the one on page \pageref{lower bound for g}; we get the following.

\begin{prop}
Let $m$ and $s$ be a positive integers.  The largest possible value $g_{\pm}(m,s)$ of $n$ for which $\phi_{\pm}(\mathbb{Z}_{n},[0,s])=m$ holds for all $n \leq g_{\pm}(m,s)$ satisfies
$$g_{\pm}(m,s) \geq \left(2  \left \lfloor \frac{s}{m} \right \rfloor +1 \right)^m.$$

\end{prop}

For example, for $m=3$ we get
$$\frac{8}{27}s^3 \sim \left(2  \left \lfloor \frac{s}{3} \right \rfloor +1 \right)^3 \leq g_{\pm}(3,s) \leq f_{\pm}(3,s) \leq a(3,s) \sim \frac{4}{3} s^3.$$

\begin{prob}
Find (if possible) better bounds for $g_{\pm}(3,s)$ and $ f_{\pm}(3,s) $ than those above.
\end{prob}

We mention that, as an attempt to find a lower bound for $g_{\pm}(3,s)$, Doskov, Pokhrel, and Singh\index{Doskov, N.}\index{Pokhrel, P.}\index{Singh, S.} in \cite{Dos:2003a} conjectured that the set
$$A(s)=\left \{  \left \lceil \frac{s^2}{2} \right \rceil -s+1, \; \frac{s^2-s}{2}+2, \; \frac{s^2-s}{2}+3 \right\}$$ is $s$-spanning in $\mathbb{Z}_n$ for all $n$ up to $s^3-s^2+6s+1$ and thus $g_{\pm}(3,s) \geq s^3-s^2+6s+1$; while this indeed holds for $s \leq 7$, $A(8)$ only generates 493 elements of $\mathbb{Z}_{497}$.

\subsection{Arbitrary number of terms} \label{3minUSarbitrary}

\section{Restricted sumsets} \label{3minR}

\subsection{Fixed number of terms} \label{3minRfixed}

\subsection{Limited number of terms} \label{3minRlimited}

A subset $A$ of $G$ for which $[0,s]\hat{\;}A =G$ for some nonnegative integer $s$ is called a {\em restricted $s$-basis} for $G$.  Here we investigate $$\phi \hat{\;}(G,[0,s]) = \mathrm{min} \{ |A|  \mid A \subseteq G, [0,s]\hat{\;}A =G\},$$ that is, the minimum size of a restricted $s$-basis for $G$.  

Clearly, we have $$\phi \hat{\;}(G,[0,0]) =\infty$$ for all groups of order at least 2, but, since $1 \hat{\;} G=G$, $$\phi \hat{\;} (G,[0,s]) \leq n$$ for every $s \in \mathbb{N}$.  Furthermore, for every $G$ of order 2 or more we get $$\phi\hat{\;}(G,[0,1])=n-1,$$ since for every $A \subseteq G$ we have $[0,1]\hat{\;} A=\{0\} \cup A$.  But for $s \geq 2$, we do not have a formula for $\phi\hat{\;}(G,[0,s])$.  Therefore we pose the following problem.

\begin{prob}
Find $\phi \hat{\;} (G,[0,s])$ for every $G$ and $s \geq 2$.
\end{prob}

From Proposition \ref{philower}, we have the following general bound.

\begin{prop} \label{philowerRlimited}  If $A$ is an $s$-basis of size $m$ in $G$, then
$$n \leq \sum_{h=0}^s {m \choose h}.$$
\end{prop}

We are particularly interested in the extremal cases of Proposition \ref{philowerRlimited}.  Namely, we would like to find {\em perfect restricted $s$-bases}, that is, restricted $s$-bases of size $m$ with $$n=\sum_{h=0}^s {m \choose h}.$$

It is easy to determine all perfect restricted $1$-bases:  

\begin{prop}
A subset $A$ of a finite abelian group is a perfect restricted $1$-basis in $G$ if, and only if, $A=G \setminus \{0\}$.

\end{prop}

Let us pursue a search for perfect restricted $2$-bases.  If $A$ is an $m$-subset of $G$ that is a perfect restricted $2$-basis, then we must have
$$n = \frac{m^2+m+2}{2}.$$  Thus, for $m=1$ we have $n=2$ and thus $G \cong \mathbb{Z}_2$; $\{1\}$ is indeed a perfect restricted $2$-basis in $\mathbb{Z}_2$.  For $m=2$ we have $n=4$ and thus $G \cong \mathbb{Z}_4$ or $G \cong \mathbb{Z}_2^2$; $\{1,2\}$ is a perfect restricted $2$-basis in $\mathbb{Z}_4$, and $\{(0,1),(1,0)\}$ is a perfect restricted $2$-basis in $\mathbb{Z}_2^2$.  Moving on to $m=3$, we have $n=7$ and thus $G \cong \mathbb{Z}_7$; $\{1,2,4\}$ (for example) is a perfect restricted $2$-basis in $\mathbb{Z}_7$.  However, as Krasny\index{Krasny, O.} in \cite{Kras:2010a} verified, there are no perfect restricted $2$-bases of size 4, 5, 6, 7, or 8 in the relevant cyclic groups.   In fact, we are not aware of any other perfect restricted $s$-bases.

\begin{prob}
Find restricted perfect $s$-bases for $s \geq 3$ or for $s=2$ and $m \geq 4$, or prove that they do not exist.
\end{prob}  

Let us return to the question of finding $\phi \hat{\;}(G,[0,s])$.  As a consequence of Proposition \ref{philowerRlimited}, we have a lower bound for $\phi \hat{\;} (G,[0,s])$, although this bound is difficult to exhibit exactly.  For $s=2$, though, we have the following. 

\begin{prop} \label{philowergeneralR}
For any abelian group $G$ of order $n$ we have $$\phi \hat{\;} (G,[0,2]) \geq \left \lceil \frac{\sqrt{8n-7}-1}{2} \right\rceil .$$

\end{prop}
  
Regarding the upper bound, note that our construction for Propositions \ref{phiuppercyclic} and \ref{phiuppergeneral} used distinct terms, thus we again have: 

\begin{prop} \label{phiuppercyclicR}
For positive integers $n$ and $s$ we have $$\phi \hat{\;} (\mathbb{Z}_n,[0,s]) \leq s \cdot \left( \lceil \sqrt[s]{n} \rceil -1 \right) < s \cdot \sqrt[s]{n} .$$

\end{prop} 

\begin{prop} \label{phiuppergeneralR}
Let $G$ be an abelian group of rank $r$ and order $n$, and let $s$ be a positive integer.   We then have $$\phi \hat{\;} (G,[0,s])  < s^r \cdot \sqrt[s]{n} .$$

\end{prop}

Krasny\index{Krasny, O.} (see \cite{Kras:2010a}) developed a computer program that yields the following data. 

$ \phi \hat{\;} (\mathbb{Z}_n,[0,2])= \left\{
\begin{array}{ll}
1 & \mbox{if $n=1, 2;$}\\
2 & \mbox{if $n=3,4;$}\\
3 & \mbox{if $n=5,6,7;$}\\
4 & \mbox{if $n=8,9,10;$}\\
5 & \mbox{if $n=11,12,13,14;$}\\
6 & \mbox{if $n=15,16,17,18,19,20;$}\\
7 & \mbox{if $n=21,22,23,24;$}\\
8 & \mbox{if $n=25,26,27,28,29,30;$}\\
9 & \mbox{if $n=31,32,33,34,35,36,37;$}\\
\end{array}\right.$

$ \phi \hat{\;}(\mathbb{Z}_n,[0,3])= \left\{
\begin{array}{ll}
1 & \mbox{if $n=1,2;$}\\
2 & \mbox{if $n=3,4;$}\\
3 & \mbox{if $n=5,6,7,8;$}\\
4 & \mbox{if $n=9,\dots,15;$}\\
5 & \mbox{if $n=16, \dots, 24;$}\\
6 & \mbox{if $n=25, \dots, 35;$}\\
7 & \mbox{if $n=36, \dots, 50;$}\\
\end{array}\right.$

and 

$ \phi \hat{\;} (\mathbb{Z}_n,[0,4]) = \left\{
\begin{array}{ll}
1 & \mbox{if $n=1, 2;$}\\
2 & \mbox{if $n=3,4;$}\\
3 & \mbox{if $n=5,6,7,8;$}\\
4 & \mbox{if $n=9,\dots,16;$}\\
5 & \mbox{if $n=17,\dots,31;$}\\
6 & \mbox{if $n=32,\dots,52.$}\\
\end{array}\right.$

As these values suggest, $\phi \hat{\;} (\mathbb{Z}_n,[0,s])$ behaves rather peculiarly; the following problem is wide open. 

\begin{prob}
Find $\phi \hat{\;} (\mathbb{Z}_n,[0,s])$ for all $n$ and $s$.
\end{prob}

In particular, it is worth considering the following special cases:

\begin{prob}
Find $\phi \hat{\;} (\mathbb{Z}_n,[0,2])$ for all $n$.
\end{prob}

\begin{prob}

For all positive integers $s$ and $m$, find the largest integer $f \hat{\;} (m,s)$ for which $\phi \hat{\;}(\mathbb{Z}_{n},[0,s]) \leq m$ holds for $n =f\hat{\;}(m,s)$.
\end{prob}

\begin{prob}

For all positive integers $s$ and $m$, find the largest integer  $g \hat{\;} (m,s)$ for which $\phi \hat{\;}(\mathbb{Z}_{n},[0,s]) \leq m$ holds for all $n \leq g \hat{\;}(m,s)$.
\end{prob}

\subsection{Arbitrary number of terms} \label{3minRarbitrary}

\section{Restricted signed sumsets} \label{3minRS}

\subsection{Fixed number of terms} \label{3minRSfixed}

\subsection{Limited number of terms} \label{3minRSlimited}

\subsection{Arbitrary number of terms} \label{3minRSarbitrary}

\chapter{Sidon sets} \label{ChapterSidon}

Recall that for a given finite abelian group $G$, $m$-subset $A=\{a_1,\dots, a_m\}$ of $G$, $\Lambda \subseteq \mathbb{Z}$, and $H \subseteq \mathbb{N}_0$, we defined the sumset of $A$ corresponding to $\Lambda$ and $H$ as
$$H_{\Lambda}A = \{\lambda_1a_1+\cdots +\lambda_m a_m \mbox{    } |  \mbox{    }  (\lambda_1,\dots ,\lambda_m) \in \Lambda^m(H)  \},$$
where the index set $\Lambda^m(H)$ is defined as
$$\Lambda^m(H)=\{(\lambda_1,\dots ,\lambda_m) \in \Lambda^m \; |  \; |\lambda_1|+\cdots +|\lambda_m| \in H \}.$$ 

The cases where the elements of the sumset corresponding to distinct elements of the index set are themselves distinct have generated much interest: we say that a subset $A$ of size $m$ is a {\em Sidon set over $\Lambda$} in $G$ if $$|H_{\Lambda}A|=|\Lambda^m(H)|.$$  Sidon sets are named after the Hungarian mathematician Simon Sidon who introduced them in the 1930s to study Fourier series.  (Traditionally, only the special case of $\Lambda=\mathbb{N}_0$ and $H=\{2\}$ is being referred to as a Sidon set.)

In this chapter we attempt to find the maximum possible size of a Sidon set over $\Lambda$ in a given finite abelian group $G$.   Namely, our objective is to determine, for any $G$, $\Lambda \subseteq \mathbb{Z}$, and $H \subseteq \mathbb{N}_0$, the quantity
$$\sigma_{\Lambda}(G,H)=\mathrm{max} \{ |A|  \mid A \subseteq G, |H_{\Lambda}A|=|{\Lambda}^{|A|}(H)|\}.$$  If no Sidon set exists, we put $\sigma_{\Lambda}(G,H) = 0.$

There are strong connections between the maximum possible size of Sidon sets and several other quantities studied in this book---we will point out these connections throughout.

We have the following obvious bound.

\begin{prop} \label{zetaupper}
If $A$ is a Sidon set over $\Lambda$ in a group $G$ of order $n$ and $|A|=m$, then
$$n \geq |{\Lambda}^{m}(H)|.$$
\end{prop}

Proposition \ref{zetaupper} provides an upper bound for the size of Sidon sets over $\Lambda$ in $G$.  The case of equality in Proposition \ref{zetaupper} is of special interest: a Sidon set for which equality holds coincides with perfect spanning sets (see Chapter \ref{ChapterSpanning}).
 
In the following sections we consider $\sigma_{\Lambda}(G,H)$ for special coefficient sets $\Lambda$.

\section{Unrestricted sumsets} \label{4maxU}

We first study, for any $G$ and $H \subseteq \mathbb{N}_0$, the quantity
$$\sigma(G,H)=\mathrm{max} \{ |A|  \mid A \subseteq G, |HA|=|\mathbb{N}_0^{|A|}(H)|\}.$$  If no such set exists, we put $\sigma(G,H) = 0;$  for example, we clearly have $\sigma(G,H) = 0$ whenever $H$ contains two elements $h_1$ and $h_2$ whose difference is divisible by $\kappa$ (the exponent of $G$).  Conversely, if all elements of $H$ leave a different remainder when divided by $\kappa$, then $\sigma(G,H)\geq 1$, since for any element $a$ of $G$ with order $\kappa$, at least the one-element set $\{a\}$ will be a Sidon set for $H$ over $\mathbb{N}_0$.

\begin{prop} \label{sigma0}
We have $\sigma(G,H) \geq 1$ if, and only if, the elements of $H$ are pairwise incongruent mod $\kappa$.   In particular, if $|H| > \kappa$, then $\sigma(G,H) =0$.

\end{prop}

\subsection{Fixed number of terms} \label{4maxUfixed}

In this section we investigate, for a given group $G$ and positive integer $h$ so-called {\em $B_h$-sets}; that is,  sets $A$ whose $h$-term sums are distinct up to the rearrangement of terms.  Thus, $A$ is a $B_h$-set of size $m$ if, and only if, the equality
$$|h  A | = {m+h-1 \choose h}$$ holds.

We are interested in finding the maximum size of a $B_h$-set in $G$---we denote this quantity by $\sigma(G,h)$.

Clearly, every subset of $G$ is a $B_1$ set, so $$\sigma(G,1)=n.$$  We can also see that, if $A$ is a $B_h$ set for some positive integer $h$, then it is also a $B_k$ set for every positive integer $k \leq h$; therefore, $\sigma(G,h)$ is a monotone nonincreasing function of $h$.  When $h \geq \kappa$ in a group of exponent $\kappa$, then the only $B_h$ sets in the group are its 1-subsets: Indeed, if $a_1$ and $a_2$ are two distinct elements of $G$, then
$$\kappa a_1 +(h - \kappa) a_2=h a_2,$$  so no $B_h$ set in the group can contain more than one element.  Therefore, by also noting Proposition \ref{sigma0}, we get:

\begin{prop} \label{sigma 1}
If $G$ is an abelian group with exponent $\kappa$ and $h \geq \kappa$, then $\sigma(G,h)=1$.
\end{prop}

According to Proposition \ref{sigma 1}, we may limit our investigation to $2 \leq h \leq \kappa -1$.   The case $h=2$ is worth special mentioning; traditionally, a $B_2$-set is called a {\em Sidon-set}.      

As a consequence of Proposition \ref{zetaupper}, we have the following bound.

\begin{prop} \label{zetaupperUfixed}  If $A$ is a $B_h$ set of size $m$ in $G$, then
$$n \geq {m+h-1 \choose h}.$$
\end{prop}

From Proposition \ref{zetaupperUfixed} we get an upper bound for $\sigma(G,h)$: 

\begin{cor} \label{corzetaupperUfixed}
The maximum size of a $B_h$ set in an abelian group $G$ of order $n$ satisfies
$$\sigma(G,h) \leq \left \lfloor \sqrt[h]{h!n} \right \rfloor.$$
\end{cor}

With a bit more work, Bravo, Ruiz, and Trujillo\index{Bravo, J. J.}\index{Ruiz, D. F.}\index{Trujillo, C. A.} improved this, as follows:

\begin{thm} [Bravo, Ruiz, and Trujillo; cf.~\cite{BraRuiTru:2012a}] \label{Bravo at al}\index{Bravo, J. J.}\index{Ruiz, D. F.}\index{Trujillo, C. A.}
The maximum size of a $B_h$ set in an abelian group $G$ of order $n$ satisfies
$$\sigma(G,h) \leq \left \lfloor \sqrt[h]{\lfloor h/2 \rfloor ! \lceil h/2 \rceil ! n} \right \rfloor + h-1.$$
\end{thm}

Note that Theorem \ref{Bravo at al} is indeed an improvement of Corollary \ref{corzetaupperUfixed}; for example, we get $\sigma(G,2) \leq \lfloor \sqrt{2n} \rfloor$ from Corollary \ref{corzetaupperUfixed}, but $\sigma(G,2) \leq \lfloor \sqrt{n} \rfloor +1$ from Theorem \ref{Bravo at al}. 

Below we demonstrate these results for $h=2$ and $h=3$; in fact, we slightly improve the bound.  Note that, if $A=\{a_1,a_2,\dots, a_m\}$ is a Sidon set in $G$, then the set $$\{a_i-a_j  \mid 1 \leq i \leq m, 1 \leq j \leq m, i \not = j \}$$ must have exactly $m(m-1)$ elements, since all the differences listed above must be distinct.  Furthermore, none of these elements are $0$.  Therefore, $m(m-1) +1 \leq n$, which yields the following:

\begin{prop} \label{sidon tight}
For every abelian group $G$ of order $n$ we have
$$\sigma(G,2) \leq \left \lfloor \frac{\sqrt{4n-3} +1}{2} \right \rfloor.$$
\end{prop}

Similarly, one can show, as we did in the proof of Theorem \ref{no perfect bases s=2,3} on page \pageref{proof of no perfect bases s=2,3}, that if $A$ is a $B_3$ set, then $2A-A$ has cardinality $$|2A-A|=m(m-1)(m-2)/2+m(m-1)+m,$$ $A-A$ has cardinality $$|A-A|=m(m-1)+1,$$ and that these two sets are disjoint.  Therefore, 
$$m(m-1)(m-2)/2+2m(m-1)+m+1 \leq n,$$ or $$m^3+m^2 \leq 2n-2,$$ and thus certainly $$m^3 \leq 2n-2.$$ This implies:

\begin{thm} \label{B3 better}
For every abelian group $G$ of order $n$ we have
$$\sigma(G,3) \leq \left \lfloor \sqrt[3]{2n-2} \right \rfloor.$$
\end{thm}

Let us now turn to constructions of $B_h$ sets.  We start with cyclic groups, for which we present three famous results that all guarantee Sidon sets of maximum possible cardinality.

\begin{thm} [Singer; cf.~\cite{Sin:1938a}] \label{Singer thm}\index{Singer, J.} 
Suppose that $q$ is a prime power and that $n=q^2+q+1$.  We then have  $$\sigma(\mathbb{Z}_n,2) = q+1.$$ 
\end{thm}

\begin{thm} [Bose; cf.~\cite{Bos:1942a}]  \label{Bose thm}\index{Bose, R. C.}
Suppose that $q$ is a prime power and that $n=q^2-1$.  We then have  $$\sigma(\mathbb{Z}_n,2) = q.$$ 
\end{thm}

\begin{thm} [Ruzsa; cf.~\cite{Ruz:1993a}] \label{Ruzsa thm}\index{Ruzsa, I.}
Suppose that $p$ is a prime and that $n=p^2-p$.  We then have  $$\sigma(\mathbb{Z}_n,2) = p-1.$$ 
\end{thm}

We can verify that in each of these three results, $$\sigma(\mathbb{Z}_n,2)=\left \lfloor \frac{\sqrt{4n-3} +1}{2} \right \rfloor,$$ and thus are best possible.  To illuminate this further, we can compute that for positive integers $m$ and $n$, $$\left \lfloor \frac{\sqrt{4n-3} +1}{2} \right \rfloor = m$$ holds if, and only if,
$$m^2-m+1 \leq n \leq m^2+m.$$   
\label{extreme singer}When $m=q+1$ for a prime power $q$, then $$m^2-m+1=q^2+q+1,$$ so Theorem \ref{Singer thm} treats the case of smallest possible $n$, while if $m=p-1$ for a prime $p$, then $$m^2+m=p^2-p,$$ so Theorem \ref{Ruzsa thm} addresses the case of largest possible $n$.  (Theorem \ref{Bose thm} handles an intermediate case.)  In particular, note that there is no `waste' in Singer's Theorem \ref{Singer thm},\index{Singer, J.} since for his Sidon set $A$ of size $q+1$, the $(q+1)q+1$ elements of $A-A$ cover all of $\mathbb{Z}_n$.  Such a set is called a {\em perfect difference set}.

Other cases of equality in Proposition \ref{sidon tight} also deserve special interest.

\begin{prob} \label{sidon tight equality}
Find all values of $n$ for which $$\sigma(\mathbb{Z}_n,2) = \left \lfloor \frac{\sqrt{4n-3} +1}{2} \right \rfloor;$$ that is, find all values of $m$ and $n$ with
$$m^2-m+1 \leq n \leq m^2+m,$$ for which a Sidon set of size $m$ exists in $\mathbb{Z}_n$.

\end{prob}  

For small values of $m$, we can provide the answer to Problem \ref{sidon tight equality}, as follows.  In the table below, for each $m \in \{1,2,3,4,5,6\}$, we list the values of $n$ that are in the relevant range, then those values of $n$ for which a Sidon set of size $m$ exists in $\mathbb{Z}_n$ with an example for such a set, and we also list  the values of $n$ for which such a Sidon set does not exist.

\begin{center}
\begin{tabular}{|c|c|c|c|c|} \hline 
$m$ & $n$ & yes & set & no \\ \hline \hline
1 & 1,2 & 1,2 & $\{0\}$ & -- \\ \hline
2 & 3,\dots,6 & 3,\dots,6 & $\{0,1\}$ & -- \\ \hline
3 & 7, \dots, 12 & 7, \dots, 12 & $\{0,1,3\}$ & -- \\ \hline
4 & 13, \dots, 20 & 13, \dots, 20 & $\{0,1,4,6\}$ & -- \\ \hline
5 & 21, \dots,30 & 21; 23, \dots, 30 & $\{0,2,7,8,11\}$ & 22 \\ \hline
6 & 31, \dots, 42 & 31; 35, \dots, 42 & $\{0,1,4,10,12,17\}$ & 32,33,34 \\ \hline
\end{tabular}
\end{center}
It is easy to verify that the sets given work for each relevant $n$; the fact that no such sets exist for the four cases listed was verified via the computer program \cite{Ili:2017a}.

Of course, we are interested in the value of $\sigma(\mathbb{Z}_n,2)$  even when it is not of the maximum size given by Proposition \ref{sidon tight}:

\begin{prob}

Find $\sigma(\mathbb{Z}_n,2)$ for all (or at least infinitely many) positive integers $n$.

\end{prob}

Very little is known about $\sigma(\mathbb{Z}_n,h)$ for $h \geq 3$.  As generalizations of Theorems \ref{Singer thm} and \ref{Bose thm}, we have the following results:

\begin{thm} [Bose and Chowla; cf.~\cite{BosCho:1962a}] \label{Singer thm Bose and Chowla}\index{Bose, R. C.}\index{Chowla, S.}
Let $h \geq 2$.  Suppose that $q$ is a prime power and that $n=q^h+q^{h-1}+\cdots+1$.  We then have  $$\sigma(\mathbb{Z}_n,h) \geq  q+1.$$ 
\end{thm}

\begin{thm} [Bose and Chowla; cf.~\cite{BosCho:1962a}]  \label{Bose thm Bose and Chowla}\index{Bose, R. C.}\index{Chowla, S.}
Let $h \geq 2$.  Suppose that $q$ is a prime power and that $n=q^h-1$.  We then have  $$\sigma(\mathbb{Z}_n,h) \geq  q.$$ 
\end{thm}

Regarding $\sigma(\mathbb{Z}_n,h)$ in general, we can exhibit a rather obvious lower bound in terms of $\sigma(\mathbb{Z}_{n_0},h)$ for some $n_0$ that is much smaller than $n$.  Suppose that $A$ is a set of nonnegative integers, so that $A$ is a $B_h$ set when viewed in $\mathbb{Z}_{n_0}$.  Note that the largest element of $hA$ is at most $h(n_0-1)$.  It is then easy to see that $A$ is also a Sidon set in $\mathbb{Z}_n$ for all $n \geq h(n_0-1)+1$.  This yields:

\begin{prop} \label{sidon big from small}
Suppose that $n$ and $n_0$ are positive integers so that $n \geq hn_0-h+1$.  We then have $\sigma(\mathbb{Z}_n,h) \geq \sigma(\mathbb{Z}_{n_0},h)$.

\end{prop}  

For example, we can apply Proposition \ref{sidon big from small} with Theorem \ref{Ruzsa thm}, as follows.  By the famous Chebyshev's Theorem (also known as Bertrand's Postulate), there is always a prime between a positive integer and its double; in particular, for every $n$, there is a prime $p$ so that $$\left \lfloor  \sqrt{n/8} \right \rfloor < p < 2 \left \lfloor  \sqrt{n/8} \right \rfloor.$$  For this prime $p$, by Theorem \ref{Ruzsa thm} we have $\sigma(\mathbb{Z}_{p^2-p},2) \geq p-1$, and since $n \geq 2(p^2-p)-1$, by Proposition \ref{sidon big from small} we get $\sigma(\mathbb{Z}_n,2) \geq p-1$.  Therefore:

\begin{prop}  \label{lower sidon gen}
For every positive integer $n$, we have
   $$\sigma(\mathbb{Z}_n,2) \geq \left \lfloor \sqrt{n/8} \right \rfloor.$$
\end{prop}
By using stronger results than Chebyshev's Theorem, stating that primes are more dense, one can greatly improve on Proposition \ref{lower sidon gen}.  

We seem to be quite far away from the exact answer:

\begin{prob}

Find (better lower bounds for) $\sigma(\mathbb{Z}_n,h)$ for all positive integers $n$ and $h \geq 2$.

\end{prob}

Let us now turn to noncyclic groups; in particular, groups of the form $\mathbb{Z}_k^r$.  The situation is quite trivial when $k=2$, since by Proposition \ref{sigma 1} we have

\begin{prop}  \label{sidon 2-groups}
For all $r \geq 1$ and $h \geq 2$, we have $\sigma(\mathbb{Z}_2^r,h) =1$.

\end{prop}

Moving to $k \geq 3$, we have the following exact results:

\begin{thm} [Babai and S\'os; cf.~\cite{BabSos:1985a}]  \label{babai sos thm}\index{Babai, L.}\index{sos@S\'os, V. T.}
Suppose that $p$ is an odd prime and that $r$ is even.  We then have $$\sigma(\mathbb{Z}_p^r,2) =p^{r/2}.$$   

\end{thm}  

\begin{thm} [Cilleruelo; cf.~\cite{Cil:2012a}]  \label{cilleruelo thm}\index{Cilleruelo, J.}
Suppose that $q$ is a prime power.  We then have $$\sigma(\mathbb{Z}_{q-1}^2,2) =q-1.$$   

\end{thm} 
By Proposition \ref{sidon tight}, we can verify that in both cases the results are maximum possible.  We are not aware of other such cases, so we offer:

\begin{prob}  \label{prob max sidon homocyclic}
Find all values of $k$ and $r$ for which $\sigma(\mathbb{Z}_k^r,2)$ equals the maximum value allowed by Proposition \ref{sidon tight}.  
\end{prob}

For example, when $k=3$, according to Theorem \ref{babai sos thm}, the answer to this problem includes all even values of $r$.  Since we have $\sigma(\mathbb{Z}_3,2)=2$, as shown by the Sidon set $\{0,1\}$, and $\sigma(\mathbb{Z}_3^3,2)=5$, as demonstrated by  
$$\{(0,0,0),(0,0,1),(0,1,0),(1,0,0),(1,1,1)\},$$ the answer also includes $r=1$ and $r=3$.  It would be nice to know whether $\mathbb{Z}_3^5$ has a Sidon set of size 16.

Note that when $r$ is even, the maximum value allowed by Proposition \ref{sidon tight} equals $k^{r/2}$, so Problem \ref{prob max sidon homocyclic} can be specialized as follows:
\begin{prob}  \label{prob max sidon homocyclic even}
Find all values of $k$ and even values of $r$ for which $$\sigma(\mathbb{Z}_k^r,2)=k^{r/2};$$  in particular, find all values of $k$ for which $$\sigma(\mathbb{Z}_k^2,2)=k.$$
\end{prob}

By Theorems \ref{babai sos thm} and \ref{cilleruelo thm}, for $r=2$, the answer to Problem \ref{prob max sidon homocyclic even} includes all prime values of $k$, as well as those where $k+1$ is a prime power.  The first $k$ for which neither condition applies is $k=9$.  With the computer program \cite{Ili:2017a}, we can determine that $\sigma(\mathbb{Z}_9^2,2) \geq 8$; an example for a Sidon set of size eight is
$$\{(0,0),(0,1),(0,3),(1,0),(1,4),(2,2), (3,7), (5,3)\}.$$

Generalizing Problem \ref{prob max sidon homocyclic}, we have:
\begin{prob}  \label{prob max sidon other}
Find all abelian groups $G$ for which $\sigma(G,2)$ equals the maximum value allowed by Proposition \ref{sidon tight}.  
\end{prob}

Let us now present some lower bounds for $\sigma(G,2)$.  First, an obvious observation:

\begin{prop} \label{sidon product triv}
For all finite abelian groups $G_1$ and $G_2$ and for all positive integers $h$, we have
$$\sigma(G_1 \times G_2, h) \geq \sigma(G_1, h).$$  

\end{prop}
Indeed: if $A$ is a $B_h$ set in $G_1$, then $A\times \{0\}$ is a $B_h$ set in $G_1 \times G_2$.

As an immediate corollary, from Theorem \ref{babai sos thm} we get:
 \begin{thm} [Babai and S\'os; cf.~\cite{BabSos:1985a}]  \label{babai sos thm odd r}\index{Babai, L.}\index{sos@S\'os, V. T.}
Suppose that $p$ is an odd prime and that $r$ is odd.  We then have $$\sigma(\mathbb{Z}_p^r,2)  \geq p^{(r-1)/2}.$$   

\end{thm} 

For a more general lower bound, one can generalize Proposition \ref{lower sidon gen}, as follows:    

\begin{thm} [Babai and S\'os; cf.~\cite{BabSos:1985a}]  \label{babai sos thm gen}\index{Babai, L.}\index{sos@S\'os, V. T.}
Suppose that $G$ is an abelian group of rank $r$, odd order $n$, and smallest invariant factor $n_1$.  We then have $$\sigma(G,2) \geq  \left \lfloor \sqrt{n_1/8^r} \right \rfloor.$$

\end{thm}
It is important to point out that Theorem \ref{babai sos thm gen}  is stated for all $n$ in \cite{BabSos:1985a}, but their concept of a Sidon set there is slightly different from ours and only coincides with ours when $n$ is odd. 
We should also add that, just like with Proposition \ref{lower sidon gen}, the constant 8 can be largely reduced.  In particular, we have:

\begin{cor} [Babai and S\'os; cf.~\cite{BabSos:1985a}]  \label{babai sos coroll odd r}\index{Babai, L.}\index{sos@S\'os, V. T.}
For all positive integers $r$ and $k$, with $k$ odd, we have $$\sigma(\mathbb{Z}_k^r,2)  \geq (k/8)^{r/2}.$$   

\end{cor} 

It is worth pointing out that when $k$ is an odd prime, then Corollary \ref{babai sos coroll odd r} is weaker than Theorem \ref{babai sos thm} when $r$ is even, but stronger than Theorem \ref{babai sos thm odd r} when $r$ is odd and $p$ is larger than $8^r$.

The main conjecture here is as follows:

\begin{conj} [Babai; cf.~\cite{Bab:2017a}] \label{babai conj}\index{Babai, L.}
There is a positive constant $C$ so that $\sigma(G,2)  \geq C \cdot \sqrt{n}$ holds for all abelian groups $G$ of odd order $n$.
\end{conj}
Note that, by Proposition \ref{sidon tight}, the constant $C$ cannot be more than 1.

\begin{prob}
Prove Conjecture \ref{babai conj}.
\end{prob}

Note that Conjecture \ref{babai conj} is false for groups of even order; see, for example, Proposition \ref{sidon 2-groups}.

Regarding values of $h \geq 3$, we are only aware of the following not-yet-published result, which generalizes Theorem \ref{babai sos thm}:

\begin{thm} [Ruiz and Trujillo; cf.~\cite{RuiTru:2016a}]\index{Ruiz, D. F.}\index{Trujillo, C. A.}
If $p$ is a prime so that $p >h$ and $r_0$ is a divisor of $r$ for which $r_0 \geq h$, then $$\sigma(\mathbb{Z}_p^r,h) \geq  p^{r/r_0}.$$

\end{thm}
The assumption that $p > h$ is necessary: see Proposition \ref{sigma 1}.  We should also note that \cite{RuiTru:2016a} only addresses the case when $r_0=h$, but our version immediately follows from that.

The most general question remains largely unsolved:

\begin{prob}

Find $\sigma(G,h)$ for abelian groups $G$ and $h \geq 2$.

\end{prob}

It is also interesting to approach these questions from another viewpoint.  Namely, we may fix positive integers $m$ and $h$, and ask for all groups $G$ that have a $B_h$ set of size $m$.  In other words:

\begin{prob} \label{sigma fixed m}
Let $m$ and $h$ be given positive integers.  Find all finite abelian groups $G$ for which $\sigma(G,h)\geq m$.

\end{prob}

We can answer Problem \ref{sigma fixed m} easily for $h=1$, $m=1$, and $m=2$, as follows.  Since $\sigma(G,1)=n$, we have $\sigma(G,1)\geq m$ if, and only if, $n \geq m$.  By Proposition \ref{sigma0}, we have $\sigma(G,h)\geq 1$ for all $G$.  Furthermore, we can see that $\sigma(G,h)\geq 2$ if, and only if, the exponent $\kappa$ of $G$ is at least $h+1$:  The `only if' part follows from Proposition \ref{sigma 1}, and for the `if' part, note that the set $\{0,a\}$ is a $B_h$ set for every $a \in G$ of order $\kappa$, since if we were to have distinct nonnegative integers $\lambda_1$ and $\lambda_2$ with $\lambda_1 \leq h$ and $\lambda_2 \leq h$ for which
$$\lambda_1 \cdot 0 + (h-\lambda_1) \cdot a = \lambda_2 \cdot 0 + (h-\lambda_2) \cdot a,$$ then $$(\lambda_1 - \lambda_2) \cdot a=0,$$ which is impossible since $$1 \leq |\lambda_1 - \lambda_2| \leq h < \kappa.$$  

Problem \ref{sigma fixed m} is unsolved for higher values of $m$ and $h$, and may already be challenging for $m=3$.  For example, we find (via a computer program) that
\label{nu(mathbbZ_n,3,h)}
$$\sigma(\mathbb{Z}_n,h)\geq 3 \; \mbox{holds for} \; 
\left\{ \begin{array}{ll}
n \geq 7 & \mbox{if} \; h=2, \\
n \geq 13 & \mbox{if} \; h=3, \\
n=19,  \geq 21 & \mbox{if} \; h=4, \\
n \geq 30 & \mbox{if} \; h=5, \\
n =37, 39, 40, 41, \geq 43 & \mbox{if} \; h=6. \\
\end{array}\right.
$$

We formulate the following two sub-problems of Problem \ref{sigma fixed m}:

\begin{prob}
Let $m$ and $h$ be given positive integers.  Find the smallest positive integer $f(m,h)$ for which $\sigma(\mathbb{Z}_n,h)\geq m$ holds for $n = f(m,h)$.

\end{prob}

\begin{prob}
Let $m$ and $h$ be given positive integers.  Find the smallest positive integer $g(m,h)$ for which $\sigma(\mathbb{Z}_n,h)\geq m$ holds for all $n \geq g(m,h)$.

\end{prob}

According to our considerations right above, we have $f(m,1)=g(m,1)=m$, $f(1,h)=g(1,h)=1$, and $f(2,h)=g(2,h)=h+1$, and as the data above indicate, we have $f(3,4)=19$, $g(3,4)=21$, and $f(3,5)=g(3,5)=30.$

For $h=2$, we have the following data:

\begin{center}
\begin{tabular}{||c||c|c|c|c|c|c|c|c|c|c|c|c|c|c||} \hline \hline
$m$ & 1 & 2 & 3 & 4 & 5 & 6 & 7 & 8 & 9 & 10 & 11 & 12 & 13 & 14 \\ \hline
$f(m,2)$ & 1 & 3 & 7 & 13 & 21 & 31 & 48 & 57 & 73 & 91 & 120 & 133 & 168 & 183
 \\ \hline \hline
\end{tabular}
\end{center}

Most of these entries follow directly from our comments on page \pageref{extreme singer}: If $m-1$ is a prime power, then Theorem \ref{Singer thm} provides the value of $f(m,2)$.  Regarding the remaining three entries: $f(7,2)$ was found by Graham and Sloane\index{Graham, R. L.}\index{Sloane, N. J. A.} in \cite{GraSlo:1980a}, and  $f(11,2)$ and $f(13,2)$ by Haanp\"a\"a, Huima, and \"Osterg\r{a}rd\index{Haanp\"a\"a, H.}\index{Huima, A.}\index{Oster@\"Osterg\r{a}rd, P.} in  \cite{HaaHuiOst:2004a};   note that these three values are given exactly by Theorem \ref{Bose thm}.

We can prove the existence and in fact find an upper bound for $g(m,h)$ (and thus for $f(m,h)$), by verifying that for every $m \geq 1$ and $h \geq 2$, the set
$$A=\{1,h,h^2, \dots, h^{m-1}\}$$ is a $B_h$ set in $\mathbb{Z}_n$ whenever $n \geq h^m$---we will do this on page \pageref{proof of rhoZ_nUfixed}.  This yields:

\begin{prop} \label{rhoZ_nUfixed}

For all positive integers $m \geq 1$ and $h \geq 2$, we have  $$f(m,h) \leq g(m,h)  \leq h^m.$$

\end{prop}

As our computations on page \pageref{nu(mathbbZ_n,3,h)} show,  the bound in Proposition \ref{rhoZ_nUfixed}  can be greatly lowered.  In particular, note that for $h=2$, from Proposition \ref{lower sidon gen} one gets $g(m,2) \leq 8m^2$ (approximately), substantially better than $2^m$.

\subsection{Limited number of terms} \label{4maxUlimited}

Sidon sets over $\mathbb{N}_0$ for the set $[0,s]$ are called {\em $B_{[0,s]}$-sets}; that is, $A \subseteq G$ is a $B_{[0,s]}$-set if, and only if, all ${m+s \choose m}$ linear combinations of at most $s$ terms of $A$ are distinct (ignoring the order of the terms).  

In this section we ought to investigate the quantity $\sigma(G,[0,s])$, denoting the maximum size of a $B_{[0,s]}$-set in $G$.  However, as we now prove, this section is superfluous as it can be reduced to Section \ref{4maxUfixed} above.   

Suppose first that $A$ is a $B_{[0,s]}$ set in $G$ for some positive integer $s$.  Then $0A=\{0\}$ and $1A=A$ must be disjoint, so $0 \not \in A$.  Let $B=A \cup \{0\}$.  Here $B$ must be a $B_s$ set in $G$, since if two different $s$-term sums were to be equal, then, after deleting all $0$s from these sums (if needed), we get two different sums with at most $s$ terms in each that are equal, contradicting the fact that $A$ is a $B_{[0,s]}$ set in $G$.  Therefore, $\sigma(G,s) \geq \sigma(G,[0,s])+1$.

Conversely, suppose that $A$ is a $B_s$ set in $G$.  Note that, for any $g \in G$, $$A-g=\{a-g \mid a \in A\}$$ is then also a $B_s$ set, since an $s$-term linear combination of the elements of $A$ is always exactly $sg$ more than the $s$-term linear combination of the corresponding elements of $A-g$.  In particular, for any $a \in A$, the set $B=A-a$ is then a $B_s$ set in $G$ for which $0 \in B$.  Let $C=B \setminus \{0\}$.  Then $C$ is a $B_{[0,s]}$ set in $G$, since if two different linear combinations of at most $s$ terms of $C$ were to be equal, then extending each by the necessary terms of $0$s we get two different linear combinations of exactly $s$ terms of $B$ that are equal, which is a contradiction.  Thus $\sigma(G,[0,s]) \geq \sigma(G,s)-1$.  This yields:

\begin{prop}  
For every group $G$ and all positive integers $s$ we have $$\sigma(G,[0,s]) = \sigma(G,s)-1.$$

\end{prop}

\subsection{Arbitrary number of terms} \label{4maxUarbitrary}

Here we ought to consider $\sigma (G,\mathbb{N}_0)$, but this quantity (e.g. by Proposition \ref{sigma0}) is clearly 0.

\section{Unrestricted signed sumsets} \label{4maxUS}

Here we study, for any $G$ and $H \subseteq \mathbb{N}_0$, the quantity
$$\sigma_{\pm}(G,H)=\mathrm{max} \{ |A|  \mid A \subseteq G, |H_{\pm}A|=|\mathbb{Z}^{|A|}(H)|\}.$$  If no such set exists, we put $\sigma_{\pm}(G,H) = 0.$  

We clearly have $\sigma_{\pm}(G,H) = 0$ whenever $H$ contains an element $h \in \mathbb{N}$ for which $2h$ is divisible by the exponent $\kappa$ of $G$, or distinct elements $h_1, h_2 \in \mathbb{N}_0$ whose sum or difference is divisible by $\kappa$.  Conversely, if there are no such elements in $H$, then at least the one-element set $\{a\}$ whose order equals $\kappa$ will be a Sidon set for $H$ over $\mathbb{Z}$.

\begin{prop} \label{sigma0pm}
We have $\sigma_{\pm}(G,H) \geq 1$ if, and only if, there is no $h \in H$ for which $2h$  is divisible by $\kappa$ and no distinct elements $h_1, h_2 \in H$ for which $h_1 \pm h_2$ is divisible by $\kappa$.   In particular, if $|H| > \lceil \kappa/2 \rceil$, then $\sigma_{\pm}(G,H) =0$.

\end{prop}

\subsection{Fixed number of terms} \label{4maxUSfixed}

In this section we investigate, for a given $G$ and positive integer $h$, the quantity
$$\sigma_{\pm}(G,h)=\mathrm{max} \{ |A|  \mid A \subseteq G, |h_{\pm}A|=|\mathbb{Z}^{|A|}(h)|\};$$ that is, the maximum size of a $B_h$ set over $\mathbb{Z}$.  Since, according to Section \ref{0.2.4}, $$ |\mathbb{Z}^{m}(h)|=c(h,m),$$ a $B_h$ set over $\mathbb{Z}$ that has size $m$ has unrestricted signed sumset size 
$$c(h,m)=\sum_{i \geq 0} {m \choose i} {h-1 \choose i-1} 2^i.$$  

The problem of finding $\sigma_{\pm}(G,h)$ is related to $\tau_{\pm}(G,2h)$---see Chapter \ref{5maxUSfixed}.

Clearly, a subset $A$ of $G$ is a $B_1$ set over $\mathbb{Z}$ if, and only if, $A$ and $-A$ are disjoint, from which we get: 

\begin{prop} \label{sigma pm h=1}
For any abelian group $G$ of order $n$ we have
$$\sigma_{\pm}(G,1)=\frac{n-1-|\mathrm{Ord}(G,2)|}{2};$$  in particular,
$$\sigma_{\pm}(\mathbb{Z}_n,1)=\left \lfloor (n-1)/2 \right \rfloor.$$

\end{prop}
We do not know the value of $\sigma_{\pm}(G,h)$ in general for $h \geq 2$.  

It is important to point out that, in contrast to Section \ref{4maxUfixed}, a $B_h$ set over $\mathbb{Z}$ is not necessarily a $B_k$ set over $\mathbb{Z}$ for values of $k \leq h$.  Suppose, for example, that $n$ is a positive integer that is divisible by 4, and consider the set $A=\{n/4\}$ in $\mathbb{Z}_n$.  We can then see that $A$ is a $B_h$ set over $\mathbb{Z}$ in $\mathbb{Z}_n$ if, and only if, $h$ is odd: Indeed, $$h \cdot n/4=-h \cdot n/4$$ in $\mathbb{Z}_n$ if, and only if, $h$ is even.   In particular, 
$$\sigma_{\pm}(\mathbb{Z}_4,h)=
\left\{
\begin{array}{ll}
1 & \mbox{if} \; h \; \mbox{is odd}, \\
0 & \mbox{if} \; h \; \mbox{is even}.
\end{array}
\right.$$

As a consequence of Proposition \ref{zetaupper}, we have the following bound.

\begin{prop} \label{zetaupperUSlimited}  If $A$ is a $B_{h}$ set over $\mathbb{Z}$ of size $m$ in $G$, then
$$n \geq c(h,m).$$
\end{prop}

From Proposition \ref{zetaupperUSlimited} we get an upper bound for $\sigma_{\pm}(G,h)$; in particular, since $c(2,m)=2m^2$, we have 

\begin{cor} \label{corzetaupperUSfixedh=2}
For any abelian group $G$ of order $n$ we have
$$\sigma_{\pm}(G,2) \leq \left \lfloor \sqrt{n/2} \right\rfloor.$$
\end{cor}
 
It would be particularly interesting to classify situations with equality in Corollary \ref{corzetaupperUSfixedh=2}:

\begin{prob}
For each $n \in \mathbb{N}$, find all groups $G$ of order $n$ for which $$\sigma_{\pm}(G,2) =\left \lfloor \sqrt{n/2} \right\rfloor.$$

\end{prob}

The following problems are wide open.

\begin{prob}

Find $\sigma_{\pm}(\mathbb{Z}_n,h)$ for all positive integers $n$ and $h \geq 2$.

\end{prob}

\begin{prob}

Find $\sigma_{\pm}(G,h)$ for all finite abelian groups $G$ and integers $h \geq 2$.

\end{prob}

It is also interesting to approach these questions from another viewpoint.  Namely, we may fix positive integers $m$ and $h$, and ask for all groups $G$ that have a $B_{h}$ set over $\mathbb{Z}$ of size $m$.  In other words:

\begin{prob} \label{sigma pm fixed m fixed}
Let $m$ and $h$ be given positive integers.  Find all finite abelian groups $G$ for which $\sigma_{\pm}(G,h)\geq m$.

\end{prob}

We can answer Problem \ref{sigma pm fixed m fixed} for $m=1$ and for $h=1$ easily; by Propositions \ref{sigma0pm} and  Proposition \ref{sigma pm h=1}, we get the following:

\begin{prop} \label{sigma >= m for m=1, h=1}

Suppose that $G$ is a finite abelian group with exponent $\kappa$; as usual, let $\mathrm{Ord}(G,2)$ denote the set of involutions (elements of order 2) in $G$.

\begin{enumerate}
  \item 
For a given positive integer $h$, we have $\sigma_{\pm}(G,h)\geq 1$ if, and only if, $2h$ is not divisible by $\kappa$.
\item
For a given positive integer $m$, we have $\sigma_{\pm}(G,1)\geq m$ if, and only if, $$|\mathrm{Ord}(G,2)| \leq n-2m-1.$$
\end{enumerate}
\end{prop}  
The answer to Problem \ref{sigma pm fixed m fixed} is not known when $h \geq 2$ and $m \geq 2$.

We formulate the following two sub-problems of Problem \ref{sigma pm fixed m fixed}:

\begin{prob}
Let $m$ and $h$ be given positive integers.  Find the smallest positive integer $f_{\pm}(m,h)$ for which $\sigma_{\pm}(\mathbb{Z}_n,h)\geq m$ holds for $n = f_{\pm}(m,h)$.

\end{prob}

\begin{prob}
Let $m$ and $h$ be given positive integers.  Find the smallest positive integer $g_{\pm}(m,h)$ for which $\sigma_{\pm}(\mathbb{Z}_n,h)\geq m$ holds for all $n \geq g_{\pm}(m,h)$.

\end{prob}

It is easy to see that Proposition \ref{sigma >= m for m=1, h=1} implies the following:

\begin{prop}  \label{f and g pm for small}
With $f_{\pm}(m,h)$ and $g_{\pm}(m,h)$ defined above, 
\begin{enumerate}
  \item 
$f_{\pm}(1,h)$ equals the smallest integer that is not a divisor of $2h$, and $g_{\pm}(1,h)=2h+1$; and 
\item
$f_{\pm}(m,1)=g_{\pm}(m,1)=2m+1$.
\end{enumerate}

\end{prop}
To verify the second statement, note that $\mathrm{Ord}(\mathbb{Z}_n,2)$ is empty when $n$ is odd, and contains exactly one element when $n$ is even.

We can prove the existence and in fact find an upper bound for $g_{\pm}(m,h)$ (and thus for $f_{\pm}(m,h)$), by verifying that for every $m \geq 1$ and $h \geq 2$, the set
$$A=\{1,2h,(2h)^2, \dots, (2h)^{m-1}\}$$ is a $B_h$ set in $\mathbb{Z}_n$ whenever $n \geq 2h^m+1$---we will do this on page \pageref{proof of rhoZ_nUSfixed}.  This yields:

\begin{prop} \label{rhoZ_nUSfixed}

For all positive integers $m$ and $h$, we have  $$f_{\pm}(m,h) \leq g_{\pm}(m,h)  \leq (2h)^m+1.$$

\end{prop}
Observe that, by Proposition \ref{f and g pm for small}, for $m=1$, equality holds  in Proposition \ref{rhoZ_nUSfixed}.

By Proposition \ref{rhoZ_nUSfixed}, we have $$g_{\pm}(2,h) \leq 4h^2+1.$$  However, Day in \cite{Day:2011a}\index{Day, R.} conjectured that the set $A=\{a_1,h\}$, with $a_1=2$ if $h$ is odd and $a_1=1$ if $h$ is even provides a signed sumset of size $c(h,2)=4h$ when $n>2h^2$.  Indeed, we can prove the following.

\begin{prop} \label{nupmmaxm=2}
For all positive integers $h$ we have
$$g_{\pm}(2,h) \leq 2h^2+1.$$
\end{prop}

For a proof of Proposition \ref{nupmmaxm=2}, see page \pageref{proofofnupmmaxm=2}.  Note that, by Proposition \ref{zetaupperUSlimited}, we have $$g_{\pm}(2,h) \geq 4h;$$ therefore, almost certainly, one can even improve the bound $2h^2+1$.

\begin{prob}
Find a better upper bound for $g_{\pm}(2,h).$

\end{prob}

For $m=3$, from Propositions \ref{zetaupperUSlimited} and \ref{rhoZ_nUSfixed} we have
$$4h^2+2 \leq g_{\pm}(3,h) \leq 8h^3+1.$$
Day in \cite{Day:2011a} conjectured that the set $A=\{a_1, h, a_3\}$, with $a_1=2$ and $a_3=h^2+h+4$ if $h$ is odd and $a_1=1$ and $a_3=h^2+h+1$ if $h$ is even provides a signed sumset of size $c(h,3)$ when $n>2ha_3$, from which we can conjecture the following.

\begin{conj} [Day; cf.~\cite{Day:2011a}] \label{Day conj}\index{Day, R.}
For all positive integers $h$ we have
$$g_{\pm}(3,h) \leq 2h^3+2h^2+8h+1.$$
\end{conj}

\begin{prob}
Prove (or disprove) Conjecture \ref{Day conj}, or---even better---find an improved upper bound for $g_{\pm}(3,h).$

\end{prob}

\subsection{Limited number of terms} \label{4maxUSlimited}

Here we ought to study, for a given $G$ and nonnegative integer $s$, the quantity
$$\sigma_{\pm}(G,[0,s])=\mathrm{max} \{ |A|  \mid A \subseteq G, |[0,s]_{\pm}A|=|\mathbb{Z}^{|A|}([0,s])|\},$$  that is, the size of the largest subset $A$ of $G$ for which all signed sums of at most $s$ terms of $A$ are distinct (ignoring the order of the terms).  However, as we proved on page \pageref{declaration}, this quantity agrees with the maximum size of a subset $A$ of $G$ for which no signed sum of at most $2s$ terms of $A$, other than the trivial one with zero terms, is equal to $0$---in Section \ref{5maxUSlimited}, we denote this latter quantity by $\tau_{\pm}(G,[1,2s])$.  Therefore:

\begin{prop} \label{sigma=zeta}

For every finite abelian group $G$ and every positive integer $s$, we have $$\sigma_{\pm}(G,[0,s])=\tau_{\pm}(G,[1,2s]).$$

\end{prop}

In addition, for $s=0$, we trivially have
$$\sigma_{\pm}(G,[0,0])=n;$$ thus, we reduced the question of finding $\sigma_{\pm}(G,[0,s])$ to Section \ref{5maxUSlimited}.

\subsection{Arbitrary number of terms} \label{4maxUSarbitrary}

Here we ought to consider $\sigma_{\pm} (G,\mathbb{N}_0)$, but this quantity (e.g. by Proposition \ref{sigma0pm}) is clearly 0.

\section{Restricted sumsets} \label{4maxR}

In this section we study, for any $G$ and $H \subseteq \mathbb{N}_0$, the quantity
$$\sigma\hat{\;}(G,H)=\mathrm{max} \{ |A|  \mid A \subseteq G, |H\hat{\;}A|=|\hat{\mathbb{N}}_0^{|A|}(H)|\};$$ that is, the size of the largest subset of $G$ whose restricted sumset (for a given $H$) has the same size as the corresponding index set. If no such set exists, we put $\sigma\hat{\;}(G,H) = 0.$

\subsection{Fixed number of terms} \label{4maxRfixed}

The analogue of a $B_h$ set for restricted addition is called a weak $B_h$ set.  More precisely, we call a subset $A$ of $G$ (with $|A|=m$) a {\em weak $B_h$ set} if $$|h\hat{\;}A| = |\hat{\mathbb{N}}_0^{m}(h)|={m \choose h}.$$  The maximum size of a weak $B_h$ set in $G$ is denoted by $\sigma\hat{\;} (G,h)$; if no such subset exists, we write $\sigma\hat{\;} (G,h)=0$.  

Note that, if $A$ is a $B_h$ set, then it is also a weak $B_h$ set, and therefore we have $$\sigma\hat{\;}(G,h) \geq \sigma(G,h).$$

We trivially have $\sigma\hat{\;}(G,1)=n$, since every subset of $G$ is a weak $B_h$ set.  We do not know the value of $\sigma\hat{\;}(G,h)$ in general for $h \geq 2$.

As a consequence of Proposition \ref{zetaupper}, we have the following bound.

\begin{prop} \label{zetaupperRfixed}  If $A$ is a weak $B_{h}$ set of size $m$ in $G$, then
$$n \geq {m \choose h}.$$
\end{prop}

From Proposition \ref{zetaupperRfixed} we get an upper bound for $\sigma\hat{\;}(G,h)$:

\begin{cor} \label{corzetaupperRfixed}
For any abelian group $G$ of order $n$ we have
$$\sigma\hat{\;}(G,h) \leq \left \lfloor \sqrt[h]{h!n} \right\rfloor+h-1.$$
\end{cor}
 
A weak $B_h$ set for $h=2$ is called a weak Sidon set.  Weak Sidon sets were introduced and studied by Ruzsa in \cite{Ruz:1993a}, though the same concept under the name ``well spread set'' was investigated both earlier (cf.~\cite{Kot:1972a}\index{Kotzig, A.}) and later (cf.~\cite{Kay:2004a}\index{Kayll, P. M.}, \cite{PhiWal:1999a}\index{Phillips, N. C. K.}\index{Wallis, W. D.}).  For weak Sidon sets a better upper bound is known, in terms of the number of order-two elements of $G$:

\begin{thm} [Haanp\"a\"a and \"Osterg\r{a}rd; cf.~\cite{HaaOst:2007a}] \label{Haanpaa abelian}\index{Haanp\"a\"a, H.}\index{Oster@\"Osterg\r{a}rd, P.} 
For any abelian group $G$ of order $n$ we have
$$\sigma\hat{\;}(G,2) \leq \left \lfloor \frac{\sqrt{4n+4|\mathrm{Ord}(G,2)|+5}+3}{2} \right\rfloor.$$
\end{thm}
In particular, since in a cyclic group $|\mathrm{Ord}(G,2)| \leq 1$, from Theorem \ref{Haanpaa abelian} we get the previously published result:

\begin{cor} [Haanp\"a\"a, Huima, and \"Osterg\r{a}rd; cf.~\cite{HaaHuiOst:2004a}] \label{Haanpaa cyclic}\index{Haanp\"a\"a, H.}\index{Huima, A.}\index{Oster@\"Osterg\r{a}rd, P.} 
For any positive integer $n$ we have
$$\sigma\hat{\;}(\mathbb{Z}_n,2) \leq \left \lfloor \frac{\sqrt{4n+9}+3}{2} \right\rfloor.$$
\end{cor}

As for a lower bound, we have the result of Babai and S\'os (which, unlike with Theorem \ref{babai sos thm gen} in Section \ref{4maxUfixed}, applies to both even and odd values of $n$):

\begin{thm} [Babai and S\'os; cf.~\cite{BabSos:1985a}]  \label{babai sos thm gen hat}\index{Babai, L.}\index{sos@S\'os, V. T.}
Suppose that $G$ is an abelian group of rank $r$, odd order $n$, and smallest invariant factor $n_1$.  We then have $$\sigma \hat{\;}(G,2)  \geq  \left \lfloor \sqrt{n_1/8^r} \right \rfloor.$$

\end{thm}
(We note that it is possible to improve the constant 8.)

The following problems are wide open.

\begin{prob}

Find $\sigma\hat{\;}(\mathbb{Z}_n,h)$ for positive integers $n$ and $h \geq 2$.

\end{prob}

\begin{prob}

Find $\sigma\hat{\;}(G,h)$ for finite abelian groups $G$ and integers $h \geq 2$.

\end{prob}

It is also interesting to approach these questions from another viewpoint.  Namely, we may fix positive integers $m$ and $h$, and ask for all groups $G$ that have a weak $B_{h}$ set of size $m$.  In other words:

\begin{prob} \label{sigma hat fixed m fixed}
Let $m$ and $h$ be given positive integers.  Find all finite abelian groups $G$ for which $\sigma\hat{\;}(G,h)\geq m$.

\end{prob}

We formulate the following two sub-problems of Problem \ref{sigma hat fixed m fixed}:

\begin{prob} \label{prob f hat sidon}
Let $m$ and $h$ be given positive integers.  Find the smallest positive integer $f\hat{\;}(m,h)$ for which $\sigma\hat{\;}(\mathbb{Z}_n,h)\geq m$ holds for $n = f\hat{\;}(m,h)$.

\end{prob}

\begin{prob} \label{prob g hat sidon}
Let $m$ and $h$ be given positive integers.  Find the smallest positive integer $g\hat{\;}(m,h)$ for which $\sigma\hat{\;}(\mathbb{Z}_n,h)\geq m$ holds for all $n \geq g\hat{\;}(m,h)$.

\end{prob}  

Problems \ref{prob f hat sidon} and \ref{prob g hat sidon} are trivial when $m < h$, and we have
$$f\hat{\;}(h,h)=g\hat{\;}(h,h)=h.$$  Furthermore, the answer is clear for $h=1$: since $\sigma \hat{\;}(G,1)=n$, we have $$f\hat{\;}(m,1)=g\hat{\;}(m,1)=m.$$

The general answers to Problems \ref{prob f hat sidon} and \ref{prob g hat sidon} are not known.  For $h=2$, Haanp\"a\"a, Huima, and \"Osterg\r{a}rd\index{Haanp\"a\"a, H.}\index{Huima, A.}\index{Oster@\"Osterg\r{a}rd, P.}   computed $f\hat{\;}(m,2)$ for $m \leq 15$ (see  \cite{HaaHuiOst:2004a}; some of these values had been determined earlier by Graham and Sloane\index{Graham, R. L.}\index{Sloane, N. J. A.} in \cite{GraSlo:1980a}), and Maturo and Yager-Elorriaga\index{Maturo, A.}\index{Yager-Elorriaga, D.} computed $g\hat{\;}(m,2)$ for $m \leq 9$
 (see \cite{MatYag:2008a}):

$$\begin{array}{||c||c|c|c|c|c|c|c|c|c|c|c|c|c|c||} \hline \hline
m & 2 & 3 & 4 & 5 & 6 & 7 & 8 & 9 & 10 & 11 & 12 & 13 & 14 & 15 \\ \hline
f\hat{\;}(m,2) & 2 & 3 & 6 & 11 & 19 & 28 & 40 & 56 & 72 & 96 & 114 & 147 & 178 & 183  \\ \hline 
g\hat{\;}(m,2) & 2 & 3 & 6 & 11 & 19 & 28 & 42 & 56 &  &  &  &  &  &  \\ \hline \hline
\end{array}$$

As we can see from the table, we have $f\hat{\;}(m,2)=g\hat{\;}(m,2)$ for $m \leq 9$, with one exception: $\mathbb{Z}_n$ contains a weak Sidon set of size 8 for $n=40$ and for every $n \geq 42$, but not for $n=41$.

We can easily find a general upper bound for both $f\hat{\;}(m,h)$ and $g\hat{\;}(m,h)$ (and thus prove their existence).  Note that, for every $h \in \mathbb{N}$, the set $$\{1,2,2^2,\dots,2^{m-1}\}$$ is a weak $B_h$ set in $\mathbb{Z}_n$ when $n$ is at least $2^m$.  Therefore:

\begin{prop} \label{rhoZ_nRfixed}
For all positive integers $m$ and $h$ we have
$$f\hat{\;}(m,h) \leq g\hat{\;}(m,h) \leq 2^m.$$
\end{prop}
Undoubtedly, this bound---which does not even depend on $h$---can be greatly reduced.

A bit more ambitiously than Problem \ref{prob f hat sidon}, but less ambitiously than Problem \ref{sigma hat fixed m fixed}, we may pose the following:

\begin{prob} \label{prob F hat sidon}
Let $m$ and $h$ be given positive integers.  Find the smallest positive integer $F\hat{\;}(m,h)$, for which there exists a group $G$ of order $F\hat{\;}(m,h)$ such that $\sigma\hat{\;}(G,h) = m$.

\end{prob}

For $h=2$, Haanp\"a\"a  and \"Osterg\r{a}rd\index{Haanp\"a\"a, H.}\index{Oster@\"Osterg\r{a}rd, P.} computed $F\hat{\;}(m,2)$ for $m \leq 15$ (see  \cite{HaaOst:2007a}):

$$\begin{array}{||c||c|c|c|c|c|c|c|c|c|c|c|c|c|c||} \hline \hline
m & 2 & 3 & 4 & 5 & 6 & 7 & 8 & 9 & 10 & 11 & 12 & 13 & 14 & 15 \\ \hline
F\hat{\;}(m,2) & 2 & 3 & 6 & 11 & 16 & 24 & 40 & 52 & 72 & 96 & 114 & 147 & 178 & 183  \\ \hline 
\hline
\end{array}$$

When comparing values of $F\hat{\;}(m,2)$ and $f\hat{\;}(m,2)$ above, we see that they agree for all $m \leq 15$, except as follows:
\begin{itemize}
\item $F\hat{\;}(6,2)=f\hat{\;}(6,2)-3$, as $\mathbb{Z}_2^2 \times \mathbb{Z}_4$ and $\mathbb{Z}_2^4$ both have weak Sidon sets of size 6;
\item $F\hat{\;}(7,2)=f\hat{\;}(7,2)-4$, as $\mathbb{Z}_2^3 \times \mathbb{Z}_3$ has a weak Sidon sets of size 7; and
\item $F\hat{\;}(9,2)=f\hat{\;}(9,2)-4$, as $\mathbb{Z}_2^2 \times \mathbb{Z}_{13}$ has a weak Sidon sets of size 9.
\end{itemize}

\subsection{Limited number of terms} \label{4maxRlimited}

We say that an $m$-subset $A$ of $G$ is a {\em weak $B_{[0,s]}$-set} if, and only if, the equality
$$|{[0,s]} \hat{\;} A | =|\hat{\mathbb{N}}_0^m([0,s])|=\sum_{h=0}^s {m \choose h}$$ holds.
We are interested in finding the maximum size of a weak $B_{[0,s]}$-set in $G$---we denote this quantity by $\sigma \hat{\;}(G,[0,s])$.

Clearly, every subset of $G \setminus \{0\}$ is a weak $B_{[0,1]}$ set, so $$\sigma \hat{\;}(G,[0,1])=n-1.$$  We do not know the value of $\sigma \hat{\;}(G,[0,s])$ in general for $s \geq 2$.    

As a consequence of Proposition \ref{zetaupper}, we have the following bound.

\begin{prop} \label{zetaupperRlimited}  If $A$ is a weak $B_{[0,s]}$ set of size $m$ in $G$, then
$$n \geq \sum_{h=0}^s {m \choose h}.$$
\end{prop}

Proposition \ref{zetaupperRlimited} yields an upper bound for $\sigma\hat{\;}(G,[0,s])$; for example, for $s=2$ we get: 

\begin{cor} \label{corzetaupperRlimited}
The maximum size of a $B_{[0,2]}$ set in an abelian group $G$ of order $n$ satisfies
$$\sigma\hat{\;} (G,[0,2]) \leq \left\lfloor \frac{\sqrt{8n-7}-1}{2} \right\rfloor.$$
\end{cor}
 
The following problems are wide open.

\begin{prob}

Find $\sigma\hat{\;}(\mathbb{Z}_n,[0,s])$ for positive integers $n$ and $s \geq 2$.

\end{prob}

\begin{prob}

Find $\sigma\hat{\;}(G,[0,s])$ for finite abelian groups $G$ and integers $s \geq 2$.

\end{prob}

It is also interesting to approach these questions from another viewpoint.  Namely, we may fix positive integers $m$ and $s$, and ask for all groups $G$ that have a weak $B_{[0,s]}$ set of size $m$.  In other words:

\begin{prob} \label{sigma restr fixed m limited}
Let $m$ and $s$ be given positive integers.  Find all finite abelian groups $G$ for which $\sigma\hat{\;}(G,[0,s])\geq m$.

\end{prob}

We formulate the following two sub-problems of Problem \ref{sigma restr fixed m limited}:

\begin{prob}
Let $m$ and $s$ be given positive integers.  Find the smallest positive integer $f\hat{\;}(m,[0,s])$ for which $\sigma\hat{\;}(\mathbb{Z}_n,[0,s])\geq m$ holds for $n = f\hat{\;}(m,[0,s])$.

\end{prob}

\begin{prob}
Let $m$ and $s$ be given positive integers.  Find the smallest positive integer $g\hat{\;}(m,[0,s])$ for which $\sigma\hat{\;}(\mathbb{Z}_n,[0,s])\geq m$ holds for all $n \geq g\hat{\;}(m,[0,s])$.

\end{prob}

We can easily find a general upper bound for both $f\hat{\;}(m,[0,s])$ and $g\hat{\;}(m,[0,s])$ (and thus prove their existence).  Note that, for every $s \in \mathbb{N}$, the set $$\{1,2,2^2,\dots,2^{m-1}\}$$ is a weak $B_{[0,s]}$ set in $\mathbb{Z}_n$ when $n$ is at least $2^m$.  Therefore: 

\begin{prop} \label{rhoZ_nRlimited}
For all positive integers $m$ and $s$ we have
$$g\hat{\;}(m,[0,s]) \leq 2^m.$$
\end{prop}

As Soma in \cite{Som:2011a}\index{Soma, T.} noticed, one can improve on the bound in Proposition \ref{rhoZ_nRlimited} slightly, as follows.  For a fixed positive integer $s$, consider the sequence $$\mathbf{a}(s)=(a_1(s),a_2(s),\dots),$$ defined as
$$
a_i(s)=
\left\{\begin{array}{cl}
2^{i-1} & \mbox{if} \; 1 \leq i \leq s, \\ \\
1+\sum_{j=i-s}^{i-1} a_j(s)& \mbox{if} \; i > s. \\
\end{array}\right.$$
For example, $$\mathbf{a}(2)=(1,2,4,7,12,20,33,54,\dots);$$ the $i$-th term we recognize as the $(i+1)$-st Fibonacci number minus 1.  Clearly, if $n \geq a_{m+1}(s)-1$, then $$A=\{a_1(s),a_2(s),\dots, a_m(s)\}$$ is a weak $B_{[0,s]}$ set in $\mathbb{Z}_n$.  While the bound on $n$ is hard to compute explicitly in general, for  $s=2$ we get $F_{m+2}-2$; that is: 

\begin{prop} [Soma; cf.~\cite{Som:2011a}] \label{Soma}\index{Soma, T.} 
For all positive integers $m$ we have
$$g\hat{\;}(m,[0,2]) \leq \left[ \tfrac{1}{\sqrt{5}} \left(\tfrac{1+\sqrt{5}}{2} \right)^{m+3}\right]-2.$$
\end{prop}
(As usual, $[x]$ denotes the closest integer to the real number $x$.)  This bound is still exponential in $m$ but its base is smaller than 2, thus it improves Proposition \ref{rhoZ_nRlimited}.   The measure of improvement is difficult to state  for $s>2$ and becomes less pronounced as $s$ gets closer to $m$.

\subsection{Arbitrary number of terms} \label{4maxRarbitrary}

Here we should consider, for a given group $G$, $$\sigma\hat{\;} (G,\mathbb{N}_0) = \mathrm{max} \{ m \mid A \subseteq G, \; |A|=m, \; |\Sigma  A| = |\hat{\mathbb{N}}_0^m(\mathbb{N}_0)|=2^m \};$$
that is, the maximum size of a subset of $G$ for which all restricted sums are distinct.

However, as we now show, this subsection is redundant in that the following holds:

\begin{prop} \label{C.3.3=F.4.3}
For any subset $A$ of $G$ we have $|\Sigma A|=2^m$ if, and only if, $0 \not \in \cup_{h=1}^{\infty} h \hat{_{\pm}} A$. 

\end{prop}

To see this, simply note that two restricted sums are equal if, and only if, their difference is zero; in particular, no two `different' restricted sums being equal is equivalent to no `nontrivial' signed sum equaling zero.

As a consequence of Proposition \ref{C.3.3=F.4.3}, we have
$$\sigma\hat{\;} (G,\mathbb{N}_0) = \tau \hat{_{\pm}} (G, \mathbb{N}),$$ so it suffices to study this quantity in Section \ref{5maxRSarbitrary}.

\section{Restricted signed sumsets} \label{4maxRS}

\subsection{Fixed number of terms} \label{4maxRSfixed}

\subsection{Limited number of terms} \label{4maxRSlimited}

\subsection{Arbitrary number of terms} \label{4maxRSarbitrary}

Here we consider, for a given group $G$, $$\sigma\hat{_{\pm}} (G,\mathbb{N}_0) = \mathrm{max} \{ m \mid A \subseteq G, \; |A|=m, \; |\Sigma_{\pm} A| = |\hat{\mathbb{Z}}^m(\mathbb{N}_0)|=3^m \}.$$

The case when $G$ is cyclic is easy.  Suppose that $G=\mathbb{Z}_n$, and let $$m=\lfloor \log_3 n \rfloor;$$ note that we then have
$$3^m \leq n < 3^{m+1}.$$
Consider the set $$A=\{1,3,\dots,3^{m-1}\}$$ in $\mathbb{Z}_n$.  Observe that $|A|=m$.  We can easily see (see Section \ref{sectionmaxsumsetsizeRSarbitrary}) that
$$\Sigma _{\pm}A=\left\{-\frac{3^{m}-1}{2},\dots,-1,0,1,\dots,\frac{3^{m}-1}{2}\right\},$$
so $$|\Sigma_{\pm} A|=3^m;$$  therefore, 
$$ \sigma\hat{_{\pm}}(\mathbb{Z}_n,\mathbb{N}_0) \geq m.$$  To show that we cannot do better, suppose, indirectly, that $A' \subseteq \mathbb{Z}_n$, $|A'| = m+1$, and $$|\Sigma_{\pm} A'| = 3^{m+1}.$$  But this implies that $$n \geq 3^{m+1},$$ a contradiction.  Thus we have proved the following.

\begin{prop}

For all positive integers $n$ we have
$$\sigma\hat{_{\pm}}(\mathbb{Z}_n,\mathbb{N}_0) = \lfloor \log_3 n \rfloor.$$

\end{prop}

So, the case of cyclic groups has been settled, which leaves us with the following problem.

\begin{prob}
Find the value of $\sigma\hat{_{\pm}}(G,\mathbb{N}_0)$ for noncyclic groups $G$.
\end{prob}

\chapter{Minimum sumset size}  \label{ChapterMinsumsetsize}

In this chapter we attempt to answer the following question: Given a finite abelian group $G$ and a positive integer $m$, how small can a sumset of an $m$-subset of $G$ be?  More precisely, our objective is to determine, for any $G$, $m$, $\Lambda$, and $H$ the quantity
$$\rho_{\Lambda}(G,m,H)=\mathrm{min} \{ |H_{\Lambda}A|  \mid A \subseteq G, |A|=m\}.$$

In the following sections we consider $\rho_{\Lambda}(G,m,H)$ for special $\Lambda \subseteq \mathbb{Z}$ and $H \subseteq \mathbb{N}_0$.

\section{Unrestricted sumsets} \label{sectionminsumsetsizeU}

Our goal in this section is to investigate the quantity $$\rho (G,m,H) =\mathrm{min} \{ |HA|  \mid A \subseteq G, |A|=m\}$$ where $HA$ is the union of all $h$-fold sumsets $hA$ for $h \in H$.  We consider three special cases: when $H$ consists of a single nonnegative integer $h$, when $H$ consists of all nonnegative integers up to some value $s$, and when $H$ is the entire set of nonnegative integers.

\subsection{Fixed number of terms} \label{sectionminsumsetsizeUfixed}

Here we consider $$\rho (G,m,h) = \mathrm{min} \{ |hA|  \mid A \subseteq G, |A|=m\},$$ that is, the minimum size of an $h$-fold sumset of an $m$-element subset of $G$.

Clearly, for all $A$ we have $0A=\{0\}$ and $1A=A$, so $$\rho (G,m,0) = 1$$ and 
$$\rho (G,m,1) = m.$$
It is also easy to see that for all $A=\{a_1,\dots,a_m\}$ and $h \geq 1$ we have $$\{(h-1)a_1+a_i  \mid i=1,2,\dots,m\} \subseteq hA,$$ and thus we have the following obvious bound.

\begin{prop} \label{rho (G,m,h)>=m}
For all $h \geq 1$ we have $\rho (G,m,h) \geq m$.

\end{prop}

As we have said above, equality always holds when $h=1$; in Proposition \ref{|hA|=m} below, we classify all other cases for equality in Proposition \ref{rho (G,m,h)>=m}.   But first we provide a general construction that, as it turns out, provides the minimum possible size of an $h$-fold sumset that an $m$-element subset of $G$ can have.

How can one find $m$-subsets $A$ in a group $G$ that have small $h$-fold sumsets $hA$?  Two ideas come to mind.  First, observe that if $A$ is a subset of a subgroup $H$ of $G$, then $hA$ will be a subset of $H$ as well; a bit more generally, if $A$ is a subset of any coset $g+H$ of $H$, then $hA$ will be a subset of the coset $hg+H$.  For example, with $G=\mathbb{Z}_{15}$, $m=4$, and $h=2$, we may choose $A$ to be any 4-subset of $H=\{0,3,6,9,12\}$; we then have $2A \subseteq H$ and thus can conclude that  $\rho (\mathbb{Z}_{15},4,2) \leq 5$.  (More generally, $A$ may be a 4-subset of any coset of $H$.)

Our second idea is based on the observation that, when $A$ is an arithmetic progression
$$A=\{a,a+g,a+2g,\dots,a+(m-1)g\}$$ for some $a, g \in G$, then many of the $h$-fold sums coincide; in particular, we have
$$hA=\{ha,ha+g,ha+2g,\dots,ha+(hm-h)g\}.$$  For example, with $G=\mathbb{Z}_{15}$, $m=4$, and $h=2$, we may choose $A$ to be the set $\{0,1,2,3\}$, in which case $2A=\{0,1,2,\dots,6\}$ and thus $\rho (\mathbb{Z}_{15},4,2) \leq 7$.  While this result is worse than that of the one obtained above, in other instances the second construction may be better; for example, with $G=\mathbb{Z}_{15}$, $m=7$, and $h=2$, not having a subgroup $H$ of $G$ of size at least 7 other than $G$ itself, the first construction gives $\rho (\mathbb{Z}_{15},7,2) \leq 15$, while the second construction yields the better bound $\rho (\mathbb{Z}_{15},7,2) \leq 13$.  (We will soon see that, in fact, $\rho (\mathbb{Z}_{15},4,2) = 5$ and $\rho (\mathbb{Z}_{15},7,2) =13$.)

The general construction we are about to present is based on the combination of these two ideas: we choose a subgroup $H$ of $G$, and then select an $m$-subset $A$ so that its elements are in as few cosets of $H$ as possible; furthermore, we want these cosets to form an arithmetic progression.   

More explicitly, let us discuss this for the cyclic group $G=\mathbb{Z}_n$.  Fixing a divisor $d$ of $n$, we consider the (unique) subgroup $H$ of order $d$ of $\mathbb{Z}_n$,  $$H=\left \{j \cdot \frac{n}{d}  \mid j=0,1,2,\dots,d-1 \right\}.$$      
With $|A|=m$ and $|H|=d$, the number of cosets of $H$ that we need is $\left \lceil \frac{m}{d} \right \rceil$. (Note that $m \leq n$ assures that $\left \lceil \frac{m}{d} \right \rceil \leq \frac{n}{d}$ and thus $A$ has $m$ distinct elements.)  
\label{Ad(n,m) defined}
Let $k$ be the positive remainder of $m$ when divided by $d$; that is, write $m=cd+k$ with integers $c$ and $k$ satisfying $1 \leq k \leq d$.  We note, in passing, that $c$ and $k$ can be computed as 
$$c=\left \lceil \frac{m}{d} \right \rceil-1$$ and
$$k=m-cd=m-\left( \left \lceil \frac{m}{d} \right \rceil-1 \right) \cdot d,$$ respectively.  Now we choose our cosets to be
$i +H$ with $i=0,1,2,\dots, c$, and set
$$A_d(n,m)=\bigcup_{i=0}^{c-1} (i+H) \cup \left\{c + j \cdot \frac{n}{d}  \mid j=0,1,2,\dots,k-1 \right\}.$$ 
Then $A_d(n,m)$ has size $m$.  Note that, when $c=0$ (that is, when $m \leq d$), then $$\bigcup_{i=0}^{c-1} (i+H)=\emptyset,$$ so $A$ lies entirely within a single coset and forms the arithmetic progression
$$\left \{j \cdot \frac{n}{d}  \mid j=0,1,2,\dots,m-1 \right\}.$$  On the other hand, $k \geq 1$, so $$\left\{c + j \cdot \frac{n}{d}  \mid j=0,1,2,\dots,k-1 \right\}$$ is never the empty-set.

It is easy to see that we have
$$hA_d(n,m)=\bigcup_{i=0}^{hc-1} (i+H) \cup \left\{hc + j \cdot \frac{n}{d}  \mid j=0,1,2,\dots,h(k-1) \right\},$$
and thus
\begin{eqnarray*}
|hA_d(n,m)|&=&\min\{n, \; hcd+\min\{d, \; h(k-1)+1\}\} \\ \\
&=&\min\{n, \; (hc+1)d, \; hc d + h(k-1)+1\} \\ \\
&=&\min\{n, \; \left(h\left \lceil \frac{m}{d} \right \rceil-h +1 \right) \cdot d, \; hm-h+1\}.
\end{eqnarray*}
Recalling our notation $$f_d=f_d(m,h)=\left(h\left \lceil \frac{m}{d} \right \rceil-h +1 \right) \cdot d$$ from Section \ref{0.3.5.1}, we get the following:
\begin{prop} \label{|hA_d(n,m)|}
For integers $n$, $m$, $h$, a given divisor $d$ of $n$, and $A_d(n,m)$ defined as above, we have
$$|hA_d(n,m)|=\min\{n, \; f_d, \; hm-h+1\}.$$
\end{prop}

Recalling from Section \ref{0.3.5.1} that 
$$u(n,m,h)=\min \{f_d  \mid d \in D(n)\},$$
and noting that $f_n=n$ and $f_1=hm-h+1$, we can 
 summarize our findings to say that

$$\rho(\mathbb{Z}_n,m,h) \leq \min\{|hA_d(n,m)|  \mid d \in D(n) \}=u(n,m,h).$$

It is somewhat surprising that we get the same bound by taking a potentially larger set that is easier to work with.  Namely, completing the part of $A_d(n,m)$ that falls into the last coset, we get
$$\overline{A}_d(n,m)=\bigcup_{i=0}^{c} (i+H),$$ a set that is the union of $c+1=\left \lceil \frac{m}{d} \right \rceil$ cosets and thus is of size
$$|\overline{A}_d(n,m)| = \left \lceil \frac{m}{d} \right \rceil \cdot d \geq m.$$
We now find that
\begin{eqnarray*}
|h\overline{A}_d(n,m)|&=&\min\{n, \; (hc+1)d\} \\ \\
&=&\min\{n, \; f_d\} ,
\end{eqnarray*}
and therefore
$$\rho(\mathbb{Z}_n,m,h) \leq \min\{|h \overline{A} _d(n,m)|  \mid d \in D(n) \}=u(n,m,h).$$

Next, we ask whether it is possible for the $h$-fold sumset of an $m$-subset of $\mathbb{Z}_n$ to have size less than $u(n,m,h)$.  The question of whether one can improve on the bound $$\rho (G,m,h) \leq u(n,m,h)$$ or not for $\mathbb{Z}_n$ or other abelian groups has a long history.  Cauchy's result\index{Cauchy, A--L.} from 1813 in \cite{Cau:1813a} implies that for $h=2$ and for cyclic groups of prime order we have equality above; since for a prime $p$ we have 
$$u(p,m,2)=\min \{ p, 2m-1\}$$ 
(see Proposition \ref{u(p,m,h)}), we thus have 
$$\rho (\mathbb{Z}_p,m,2)=\min \{ p, 2m-1\}.$$
In 1935, Davenport in \cite{Dav:1935a}  rediscovered Cauchy's result, which is now known as the {\em Cauchy--Davenport Theorem}.\index{Cauchy, A--L.}\index{Davenport, H.}  (See page \pageref{CDthm} for the general statement of the theorem.  Davenport\index{Davenport, H.}  was unaware of Cauchy's result\index{Cauchy, A--L.} until twelve years later; see \cite{Dav:1947a}.)  Finally, after various partial results by several researchers, a sequence of papers at the beginning of the twenty-first century, including \cite{EliKer:1998a} by Eliahou and Kervaire; \cite{Pla:2003a}\index{Eliahou, S.}\index{Kervaire, M.} by 
Plagne; \cite{EliKerPla:2003a} by Eliahou, Kervaire, and Plagne;\index{Eliahou, S.}\index{Kervaire, M.} and \cite{Pla:2006a} by Plagne\index{Plagne, A.} established that, indeed, equality holds for all finite abelian groups and all $h$:

\begin{thm} [Plagne; cf.~\cite{Pla:2006a}] \label{rho (G,m,h)}\index{Plagne, A.} 
For any finite abelian group $G$ of order $n$ we have $$\rho (G,m,h)=u(n,m,h).$$
\end{thm}

It may be worthwhile to state the following more explicit consequence of Theorem \ref{rho (G,m,h)} in light of Proposition \ref{u versus p}.

\begin{cor} \label{rho vs p}
Let $G$ be any finite abelian group of order $n \geq m$, and let $p$ be the smallest prime divisor of $n$.  We then have
$$\rho (G,m,h) \geq \min \{p, hm-h+1\},$$
with equality if, and only if, $m \leq p$.  In particular, 
$$\rho (\mathbb{Z}_p,m,h) = \min \{p, hm-h+1\}$$ for all $m \leq p$.
\end{cor}

Note that Corollary \ref{rho vs p} is a generalization of the Cauchy--Davenport Theorem.\index{Cauchy, A--L.}\index{Davenport, H.}  

According to Theorem \ref{rho (G,m,h)}, the question of finding the minimum size $\rho (G,m,h)$ of the $h$-fold sumset of an $m$-subset of $G$ is solved.  We then may ask for a classification of all subsets $A$ of $G$ for which $|hA|=\rho (G,m,h)$.  Since, as we have seen, we have $$|0A|=\rho (G,m,0)=1$$ and $$|1A|=\rho (G,m,1)=m$$ for all $G$ and $A$ (with $|A|=m \leq n$), we may assume that $h \geq 2$. 

\begin{prob} \label{classifyrho(G,m,h)}
For each $G$, $m$, and $h \geq 2$, find a characterization of all $m$-subsets of $G$ with $h$-fold sumset of size $\rho (G,m,h)$.
\end{prob}  

For example, let us consider $G=\mathbb{Z}_{15}$, $m=6$, and $h=2$, for which we have $\rho (\mathbb{Z}_{15},6,2)=9$ (see page \pageref{u(15,m,h)}).  As we can verify, a 6-subset $A$ of $G$ has 2-fold sumset of size 9 if, and only if, it is of the form $$A=(a_1+H) \cup (a_2+H)$$ where $H=\{0,5,10\}$ is the subgroup of $G$ of order 3.  Thus, in this situation, the only subsets with minimum sumset size are those constructed above.

We get a more mixed picture when we consider  $G=\mathbb{Z}_{15}$, $m=7$, and $h=2$.  This time, $\rho (\mathbb{Z}_{15},7,2)=13$, and we can determine that 7-subsets of $G$ with 2-fold sumset of size 13 come in a variety of forms: as an arithmetic progression
$$\{a,a+g,a+2g, \dots, a+6g\}$$ for some $a$ and $g$ of $G$; as 
$$(a_1+H) \cup \{a_2,a_3\}$$ with $H$ being the subgroup of $G$ of order 5 and for certain $a_1,a_2,a_3 \in G$; or in the form
$$(a_1+H) \cup (a_2+H) \cup \{a_3\},$$ where $H$ is the subgroup of $G$ of order 3, and $a_1,a_2,a_3 \in G$.  For example, $$\{0,1,2,3,4,5,6\},$$ $$\{0,3,6,9,12\} \cup \{1,4\},$$ and $$\{0,5,10\} \cup \{1,6,11\} \cup \{2\}$$ all yield minimum sumset size.   As these examples indicate, Problem \ref{classifyrho(G,m,h)} may be quite difficult in general.  

We are able, however, to say more about the cases when $u(n,m,h)$ achieves its extreme values.  Recall that, by Proposition \ref{u(n,m,h)extreme}, we have
$$m \leq u(n,m,h) \leq \min \{ n, hm-h+1\},$$
with $u(n,m,h)=m$ if, and only if, $h=1$ or $n$ is divisible by $m$. 

In the case when  $u(n,m,h)=m$ and $h \geq 2$, we must have $m \in D(n)$, and therefore a subgroup of $G$ of order $m$ exists.  Taking $A$ to be a coset of this subgroup, $hA$ will also be a coset of the same subgroup, and thus $|hA|=m$.  It turns out that the converse of this holds as well, as we have the following.

\begin{prop} \label{|hA|=m}
Let $A$ be an $m$-subset of $G$ and $h \geq 1$.  Then $|hA| \geq m$ with equality if, and only if, $h=1$ or $A$ is a coset of a subgroup of $G$.
\end{prop}     
For the proof of Proposition \ref{|hA|=m}, see page \pageref{proofof|hA|=m}.  In light of Proposition \ref{|hA|=m},  Problem \ref{classifyrho(G,m,h)} is solved in the case when $u(n,m,h)=m$.

Let us now turn to the upper bounds of $u(n,m,h)$ in Proposition \ref{u(n,m,h)extreme}.  The case when $u(n,m,h)=n$ translates to the situation when every $m$-subset $A$ of $G$ has $hA=G$; sets like these are called {\em $h$-fold bases} (cf.~Chapter \ref{ChapterSpanning}), and the minimum value of $m$ for which every $m$-subset of $G$ is an $h$-fold basis of $G$ is called the {\em $h$-critical number of $G$}.  We allocate a separate chapter to critical numbers; in particular, we will study the $h$-critical number of groups in Section \ref{CritUfixed}.

The case when $u(n,m,h)=hm-h+1<n$ is also quite interesting, and we offer this question as a special case of Problem \ref{classifyrho(G,m,h)}.

\begin{prob} \label{classifyrho(G,m,h)=hm-h+1}
Find all $G$, $m$, and $h \geq 2$ for which $\rho (G,m,h)=hm-h+1<n$, and for each such instance, find a characterization of all $m$-subsets of $G$ with $h$-fold sumset of size $hm-h+1$.
\end{prob} 

Problem \ref{classifyrho(G,m,h)=hm-h+1}, while only a special case of Problem \ref{classifyrho(G,m,h)}, is still quite elusive in general, as the example of $G=\mathbb{Z}_{15}$, $m=7$, and $h=2$ above indicates.  We do have a complete answer, however, when $G$ is of prime order.

\begin{thm} \label{|hA|=hm-h+1<p}
Let $p$ be a prime, $h \geq 2$, and suppose that $hm-h+1<p$.  Then an $m$-subset $A$ of $\mathbb{Z}_p$ satisfies $|hA|=hm-h+1$ if, and only if, $A$ is an arithmetic progression; that is, $$A=\{a,a+g,a+2g, \dots, a+(m-1)g\}$$ for some $a,g \in \mathbb{Z}_p$.
\end{thm}
Theorem \ref{|hA|=hm-h+1<p} can be reduced to  Vosper's Theorem\index{Vosper, A. G.} (cf.~\cite{Vos:1956a}); see page \pageref{proofof|hA|=hm-h+1<p}.  As we pointed out above, the ``only if'' part of Theorem \ref{|hA|=hm-h+1<p} doesn't hold when $h=1$, and it is also false when $hm-h+1=p$.  For example, the set $A=\{0,1,2,4\}$ is not an arithmetic progression in $\mathbb{Z}_7$ (something that can be easily verified), but  $|2A|=7=2 \cdot 4-2+1$. 

Let us see what we can say when $n$ is not prime and the lower bound of Corollary \ref{rho vs p} is achieved so that 
$$\rho (G,m,h) = \min \{p, hm-h+1\},$$ where $p$ is the smallest prime divisor of $n$.  According to Corollary \ref{rho vs p}, we then have
$p \geq m$.  We examine three cases.

When $m \leq p <  hm-h+1$, we can find $m$-subsets $A$ of $G$ for which $$|hA|=p=\min \{p, hm-h+1\},$$ as follows.  Let $H$ be any subgroup of $G$ with $|H|=p$.  Then for any $g \in G$, the coset $g+H$ contains $p$ elements, and its $h$-fold sumset is another coset of $H$.  Thus, any $m$-subset $A$ of $g+H$ has $|hA|=p$.  It turns out that there are no others:

\begin{thm} \label{thm Gryn |hA|=p}
Let $p$ be the smallest prime divisor of $n$, $A$ be an $m$-subset of $G$, and assume that $m \leq p <  hm-h+1$.  Then $|hA|=p$ if, and only if, $A$ is contained in a coset of some subgroup $H$ of $G$ with $|H|=p$.   
\end{thm}
As Grynkiewicz\index{Grynkiewicz, D. J.} pointed out in \cite{Gry:2017a}, Theorem \ref{thm Gryn |hA|=p} follows easily from Kneser's Theorem\index{Kneser, M.} \cite{Kne:1953a}; we explain this on page \pageref{proof of thm Gryn |hA|=p}. 

Assume now that $p>hm-h+1$.  This time, we can find $m$-subsets $A$ of $G$ for which $$|hA|=hm-h+1=\min \{p, hm-h+1\},$$ as follows.  Let $g \in G$ be of order $p$.  If $A$ is the arithmetic progression $$A=\{a,a+g,a+2g, \dots, a+(m-1)g\},$$ then we clearly have $$hA=\{ha,ha+g,ha+2g, \dots, ha+h(m-1)g\};$$ since $p>h(m-1)+1$, these $hm-h+1$ elements are all distinct.

As a generalization of Theorem \ref{|hA|=hm-h+1<p} above, we can prove that the converse holds as well:

\begin{thm} \label{thm |hA|=hm-h+1}
Let $h \geq 2$, $p$ be the smallest prime divisor of $n$, $A$ be an $m$-subset of $G$, and assume that $p >  hm-h+1$.  Then $|hA|=hm-h+1$ if, and only if, $A$ is an arithmetic progression.   
\end{thm}

The proof of Theorem \ref{thm |hA|=hm-h+1} follows from Kemperman's\index{Kemperman, J. H. B.}  famous results in \cite{Kem:1960a}---see page \pageref{proof of thm |hA|=hm-h+1}.  

Finally, the third case, when $p=hm-h+1$: this case seems more complicated.  Not only do we have the option of arithmetic progressions of length $m$ and cosets of a subgroup of order $p$, but there are other possibilities as well.  For example, with $p=7$, $m=4$, and $h=2$ (as is the case, for example, with $\rho (\mathbb{Z}_{49}, 4, 2)=7$) the subset $$A=\{0,a,(n-a)/2,(n+a)/2\}$$ works as well, as we get
$$2A=\{0, a, 2a, (n-a)/2, (n+a)/2, (n+3a)/2, n-a\}$$ (since $p=7$ implies that $n$ must be odd, we need to assume that $a$ is odd).  The following problem seems particularly intriguing:

\begin{prob}
Suppose that $n$, $mm$, and $h$ are positive integers, $p$ is the smallest prime divisor of $n$, and $p=hm-h+1$.  Classify all $m$-subsets $A$ of $G$ for which $|hA|=p$.
\end{prob}

\subsection{Limited number of terms} \label{2minUlimited}

Here we ought to consider, for a given group $G$, positive integer $m$ (with $m \leq n=|G|$), and nonnegative integer $s$, $$\rho (G,m,[0,s]) = \mathrm{min} \{ |[0,s]A|  \mid A \subseteq G, |A|=m\},$$ that is, the minimum size of $\cup_{h=0}^s hA$ for an $m$-element subset $A$ of $G$. 

However, we have the following result.    

\begin{prop} \label{rho  (G,m,[0,s])}
For any group $G$, positive integer $m \leq n$, and nonnegative integer $s$ we have
$$\rho (G,m,[0,s]) = \rho (G,m,s).$$ 
\end{prop}

We can prove Proposition \ref{rho  (G,m,[0,s])}, as follows.  Suppose first that $A$ is a subset of $G$ of size $m$ and that it has minimum-size $[0,s]$-fold sumset: $$|[0,s]A|=\rho (G,m,[0,s]).$$  Clearly, $sA \subseteq [0,s]A$, so  $$\rho (G,m,[0,s])=|[0,s]A| \geq |sA| \geq \rho (G,m,s).$$

For the other direction, choose a subset $A$ of $G$ of size $m$ for which $$|sA|=\rho (G,m,s).$$  By Proposition \ref{shifting sets}, we may assume that $0 \in A$, so $[0,s]A=sA$, which implies that
$$\rho (G,m,[0,s]) \leq |[0,s]A| = |sA| = \rho (G,m,s).$$ 

Proposition \ref{rho  (G,m,[0,s])} makes this subsection superfluous.

\subsection{Arbitrary number of terms} \label{2minUarbitrary}

Here we consider, for a given group $G$ and positive integer $m$ (with $m \leq n=|G|$) the quantity $$\rho (G,m,\mathbb{N}_0) = \mathrm{min} \{ |\langle A \rangle|  \mid A \subseteq G, |A|=m\}.$$  Recall that $\langle A \rangle$ is the subgroup of $G$ generated by $A$; this immediately gives the following result:  

\begin{prop} \label{rho sigma}
For any group $G$ and positive integer $m \leq n$, we have
$$\rho (G,m,\mathbb{N}_0) = \min\{d \in D(n) \mid d \geq m\}.$$ 
\end{prop}

\section{Unrestricted signed sumsets} \label{2minUS}

Our goal in this section is to investigate the quantity $$\rho_{\pm} (G,m,H) =\mathrm{min} \{ |H_{\pm} A|  \mid A \subseteq G, |A|=m\}$$ where $H_{\pm}A$ is the union of all $h$-fold signed sumsets $h_{\pm} A$ for $h \in H$.  We consider three special cases: when $H$ consists of a single nonnegative integer $h$, when $H$ consists of all nonnegative integers up to some value $s$, and when $H$ is the entire set of nonnegative integers.

\subsection{Fixed number of terms} \label{2minUSfixed}

Here we consider $$\rho_{\pm} (G,m,h) = \mathrm{min} \{ |h_{\pm}A|  \mid A \subseteq G, |A|=m\},$$ that is, the minimum size of an $h$-fold signed sumset of an $m$-element subset of $G$.

It is easy to see that $$\rho_{\pm} (G,1,h)=1,$$ $$\rho_{\pm} (G,m,0)=1,$$ and $$\rho_{\pm} (G,m,1)=m.$$ (To see the last equality, it suffices to verify that one can always find a {\em symmetric} subset of size $m$ in $G$, that is, an $m$-subset $A$ of $G$ for which $A=-A$.)   Therefore, for the rest of this subsection, we assume that $m \geq 2$ and $h \geq 2$.

Perhaps surprisingly, we find that, while the $h$-fold signed sumset of a given set is generally much larger than its sumset, $\rho_{\pm} (G, m, h)$ often agrees with $\rho (G, m, h)$; in particular, this is always the case when $G$ is cyclic:  

\begin{thm} [Bajnok and Matzke; cf.~\cite{BajMat:2014a}] \label{cyclic}\index{Bajnok, B.}\index{Matzke, R.} For all positive integers $n$, $m$, and $h$, we have 
$$\rho_{\pm} (\mathbb{Z}_n, m, h)= \rho (\mathbb{Z}_n, m, h).$$
\end{thm}

The situation seems considerably more complicated for noncyclic groups: in contrast to $\rho (G, m,h)$, the value of $\rho_{\pm} (G, m,h)$ depends on the structure of $G$ rather than just the order $n$ of $G$.  We do not have the answer to the following general problem:

\begin{prob}
Find the value of $\rho_{\pm} (G, m, h)$ for all noncyclic groups $G$.

\end{prob}  

We do, however, have some tight bounds.  Observe that by Theorem \ref{rho (G,m,h)}, we have the lower bound
$$\rho_{\pm} (G, m,h) \geq u(n,m,h)=\min \{f_d (m,h) \; : \; d \in D(n)\}.$$ In \cite{BajMat:2014a}, we proved that with a certain subset $D(G,m)$ of $D(n)$, we have
$$\rho_{\pm} (G, m,h) \leq u_{\pm} (G,m,h)=\min \{f_d (m,h) \; : \; d \in D(G,m)\};$$ here $D(G,m)$ is defined in terms of the {\em type} $(n_1,\dots,n_r)$ of $G$, that is, via  integers $n_1,\dots,n_r$ such that $n_1 \geq 2$, $n_i$ is a divisor of $n_{i+1}$ for each $i \in \{1,\dots, r-1\}$, and for which $G$ is isomorphic to  the invariant product  
$$\mathbb{Z}_{n_1} \times \cdots \times \mathbb{Z}_{n_r}.$$
Namely, we proved the following result:

\begin{thm} [Bajnok and Matzke; cf.~\cite{BajMat:2014a}]  \label{u pm with f}\index{Bajnok, B.}\index{Matzke, R.} The minimum size of the $h$-fold signed sumset of an $m$-subset of a group $G$ of type $(n_1,\dots,n_r)$ satisfies
$$\rho_{\pm} (G, m,h) \leq  u_{\pm} (G,m,h),$$
where $$u_{\pm} (G,m,h)=\min \{f_d (m,h) \; : \; d \in D(G,m) \}$$ with $$D(G,m)=\{d \in D(n) \; : \; d= d_1 \cdots d_r, d_1 \in D(n_1), \dots, d_r \in D(n_r), dn_r \geq d_rm \}.$$
\end{thm}
Observe that, for cyclic groups of order $n$, $D(G,m)$ is simply $D(n)$. 

We are aware of only one type of scenario where $u_{\pm} (G,m,h)$ does not yield the actual value of $\rho_{\pm} (G, m,h)$; it only occurs when $h=2$.  

As a prime example, consider the group $\mathbb{Z}_p^2$ for an odd prime $p$, let $m=(p^2-1)/2$, and $h=2$.  We find that 
$$D(\mathbb{Z}_p^2, (p^2-1)/2) =\{p,p^2\};$$
and we have $$f_p=f_{p^2}=p^2,$$hence 
$$u_{\pm} (\mathbb{Z}_p^2, (p^2-1)/2,2)=p^2.$$  However, observe that each nonzero element of $\mathbb{Z}_p^2$ has order $p \geq 3$, thus we can partition $\mathbb{Z}_p^2$ as $$\{0\} \cup A \cup (-A)$$  for some subset $A$ of $\mathbb{Z}_p^2$ of size $(p^2-1)/2$.  Since clearly $0 \not \in 2_{\pm}A$, we have $$\rho_{\pm} (\mathbb{Z}_p^2, (p^2-1)/2,2)\leq p^2-1.$$

A bit more generally, if $d$ is an odd element of $D(n)$ so that $d \geq 2m+1$, then the same argument yields 
$$\rho_{\pm} \left(G, m, 2   \right) \leq d-1,$$
and therefore we have the following:

\begin{prop} [Bajnok and Matzke; cf.~\cite{BajMat:2014a}]\index{Bajnok, B.}\index{Matzke, R.}
Suppose that $G$ is an abelian group of order $n$ and type $(n_1, \dots, n_r)$.  Let $m \leq n$, and let $d_m$ be the smallest odd element of $D(n)$ that is at least $2m+1$; if no such element exists, set $d_m=\infty$.  We then have
$$\rho_{\pm} \left(G, m, 2   \right) \leq \min\{u_{\pm}(G,m,2), d_m-1\}.$$

\end{prop}

We make the following conjecture:

\begin{conj} [Bajnok and Matzke; cf.~\cite{BajMat:2014a}] \label{conj for rho pm}\index{Bajnok, B.}\index{Matzke, R.}
Suppose that $G$ is an abelian  group of order $n$ and type $(n_1, \dots, n_r)$.

If $h \geq 3$, then $$\rho_{\pm} \left(G, m, h   \right) =u_{\pm}(G,m,h).$$

If each odd divisor of $n$ is less than $2m$, then $$\rho_{\pm} \left(G, m, 2  \right) =u_{\pm}(G,m,2).$$

If there are odd divisors of $n$ greater than $2m$, let $d_m$ be the smallest one.  We then have
$$\rho_{\pm} \left(G, m, 2   \right) = \min\{u_{\pm}(G,m,2), d_m-1\}.$$
\end{conj}

\begin{prob}  Prove or disprove Conjecture \ref{conj for rho pm}.
\end{prob}

We are able to say more about minimum sumset size in elementary abelian groups; in particular, we wish to study
$\rho_{\pm} (\mathbb{Z}_p^r, m, h)$,
where $p$ denotes a positive prime and $r$ is a positive integer.  By Theorem \ref{cyclic}, we assume that $r \geq 2$, and, since obviously   
$$\rho_{\pm} (\mathbb{Z}_2^r, m, h) = \rho (\mathbb{Z}_2^r, m, h)$$ for all $m$, $h$, and $r$, we will also assume that $p \geq 3$.  

Let us first exhibit a sufficient condition for $\rho_{\pm} (\mathbb{Z}_p^r, m, h)$ to equal $\rho(\mathbb{Z}_p^r, m, h)$.  When $p \leq h$, our result is easy to state:

\begin{thm} [Bajnok and Matzke; cf.~\cite{BajMat:2014b}] \label{p leq h}\index{Bajnok, B.}\index{Matzke, R.}

If $p \leq h$, then for all values of $1 \leq m \leq p^r$ we have $$\rho_{\pm} (\mathbb{Z}_p^r, m, h) = \rho (\mathbb{Z}_p^r, m, h).$$

\end{thm}

The case $h \leq p-1$ is more complicated and delicate.  In order to state our results, we will need to introduce some notations.  Suppose that $m \geq 2$ is a given positive integer.  First, we let $k$ be the maximal integer for which
$$p^k  +\delta \leq hm-h+1,$$  where 
$\delta$ equals 0 or 1, depending on whether $p-1$ is divisible by $h$ or not.  Second, we let $c$ be the maximal integer for which
$$(hc+1) \cdot p^k + \delta \leq hm-h+1.$$ Note that $k$ and $c$ are nonnegative integers and $c \leq p-1,$ since for $c \geq p$ we would have
$$(hc+1) \cdot p^k \geq p^{k+1} +\delta > hm-h+1.$$  It is also worth noting that $$f_1(m,h)=hm-h+1.$$

Our sufficient condition can now be stated as follows:

\begin{thm} [Bajnok and Matzke; cf.~\cite{BajMat:2014b}]  \label{sufficient for =}\index{Bajnok, B.}\index{Matzke, R.}

Suppose that $2 \leq h \leq p-1$, and let $k$ and $c$ be the unique nonnegative integers defined above.  If 
$$m \leq (c+1) \cdot p^k,$$then  
$$\rho_{\pm} (\mathbb{Z}_p^r, m, h) = \rho (\mathbb{Z}_p^r, m, h).$$

\end{thm}

In fact, we believe that this condition is also necessary:

\begin{conj}[Bajnok and Matzke; cf.~\cite{BajMat:2014b}] \label{conj p square}\index{Bajnok, B.}\index{Matzke, R.}

The converse of Theorem \ref{sufficient for =} is true as well; that is, if $2 \leq h \leq p-1$, $k$ and $c$ are the unique nonnegative integers defined above, and 
$$m > (c+1) \cdot p^k,$$then  
$$\rho_{\pm} (\mathbb{Z}_p^r, m, h) > \rho (\mathbb{Z}_p^r, m, h).$$
\end{conj}

\begin{prob}  Prove or disprove Conjecture \ref{conj p square}.
\end{prob}

We are able to prove that Conjecture \ref{conj p square} holds in the case of $\rho_{\pm} (\mathbb{Z}_p^2, m, 2)$:

\begin{thm} [Bajnok and Matzke; cf.~\cite{BajMat:2014b}]  \label{thm p square}\index{Bajnok, B.}\index{Matzke, R.}
Let $p$ be an odd prime and $m \leq p^2$ be a positive integer.  Then  
$$\rho_{\pm} (\mathbb{Z}_p^2, m, 2) = \rho (\mathbb{Z}_p^2, m, 2),$$ if, and only if, one of the following holds:
\begin{itemize}
  \item $m \leq p$,
  \item $m \geq (p^2+1)/2$, or
  \item there is a positive integer $c \leq (p-1)/2$ for which
  $$c \cdot p+(p+1)/2 \leq m \leq (c+1) \cdot p.$$
\end{itemize}

\end{thm}

Theorem \ref{thm p square} does not tell us the value of $\rho_{\pm} (\mathbb{Z}_p^2, m, 2)$ when it is more than $\rho (\mathbb{Z}_p^2, m, 2)$; for these cases, we have the following recent result of Lee:

\begin{thm} [Lee; cf.~\cite{Lee:2015a}]  \label{mitchell lee}\index{Lee, M.} 
Let $p$ be an odd prime.

\begin{enumerate}
  \item If $m=c \cdot p +v$ with $1 \leq c \leq (p-3)/2$ and $1 \leq v \leq (p-1)/2$, then $$\rho_{\pm} (\mathbb{Z}_p^2, m, 2)=(2c+1)p.$$
  \item If $m=c \cdot p +v$ with $c= (p-1)/2$ and $1 \leq v \leq (p-1)/2$, then $$\rho_{\pm} (\mathbb{Z}_p^2, m, 2)=p^2-1.$$
\end{enumerate}  
\end{thm}

Combining Theorems \ref{thm p square} and \ref{mitchell lee}, we get:

\begin{cor} [Bajnok and Matzke; cf.~\cite{BajMat:2014b} and Lee;\index{Lee, M.} cf.~\cite{Lee:2015a}] \label{prime square exact}\index{Bajnok, B.}\index{Matzke, R.}

Let us write $m$ as $m=cp+v$ with $$0 \leq c \leq p-1 \; \; \; \mbox{and} \; \; \; 1 \leq v \leq p.$$  We then have:
$$\begin{array}{cc|ccccc} 
c & v & \rho (\mathbb{Z}_p^2,m,2) && \rho_{\pm} (\mathbb{Z}_p^2,m,2) && u_{\pm}(\mathbb{Z}_p^2,m,2) \\ \hline 
& & & & \\
 & v \leq (p-1)/2 & 2m-1 &=& 2m-1 &=& 2m-1 \\ 
\raisebox{1.5ex}[0pt]{$0$} & v \geq (p+1)/2 & p & =&p & =&p \\ 
& & & & \\
 & v \leq (p-1)/2 & 2m-1 & <& (2c+1)p  & =&(2c+1)p \\ 
\raisebox{1.5ex}[0pt]{$1 \leq c \leq (p-3)/2$} & v \geq (p+1)/2 & (2c+1)p & =&(2c+1)p & =&(2c+1)p \\  
& & & & \\
 & v \leq (p-1)/2 & 2m-1 & <& p^2-1  & <& p^2 \\ 
\raisebox{1.5ex}[0pt]{$c=(p-1)/2$} & v \geq (p+1)/2 & p^2 & = & p^2 & =&p^2 \\ 
& & & & \\
c \geq (p+1)/2 & \mbox{any } v & p^2 & =&p^2 & =&p^2 \\ 
\end{array}$$

\end{cor}

We can also see that we have settled Conjecture \ref{conj for rho pm} for the group $\mathbb{Z}_p^2$.  Indeed, from Corollary \ref{prime square exact}, we see that $$\rho_{\pm} (\mathbb{Z}_p^2,m,2)=u_{\pm}(\mathbb{Z}_p^2,m,2)$$ in all cases, except when $m=cp+v$ with $c = (p-1)/2$ and $1 \leq v \leq (p-1)/2$, this case coincides exactly with $$d_m-1 < u_{\pm}(\mathbb{Z}_p^2,m,2)$$ for $$d_m = \min \{ d \in D(n) \mid d \mbox{ odd and } d >2m\}.$$  Therefore:

\begin{cor} [Bajnok and Matzke; cf.~\cite{BajMat:2014b} and Lee;\index{Lee, M.} cf.~\cite{Lee:2015a}]\index{Bajnok, B.}\index{Matzke, R.}
Let $p$ be an odd prime.  Conjecture \ref{conj for rho pm} holds for elementary abelian $p$-groups of rank two.
\end{cor}

According to Theorem \ref{thm p square}, for a given $p$, there are exactly $(p-1)^2/4$ values of $m$ for which $\rho_{\pm} (\mathbb{Z}_p^2, m, 2)$ and $\rho(\mathbb{Z}_p^2, m, 2)$ disagree---fewer than $1/4$ of all possible values.  We have not been able to find any groups where this proportion is higher than $1/4$, so make the following conjecture:

\begin{conj} \label{fewer than 1/4}
For any abelian group $G$ of order $n$, we have fewer than $n/4$ values of $m$ for which $\rho_{\pm} (G, m, 2)$ and $\rho(G, m, 2)$ disagree.

\end{conj} 

\begin{prob}
Prove (or disprove) Conjecture \ref{fewer than 1/4}.

\end{prob}

We do not know the generalization of Theorem \ref{thm p square} for higher rank:

\begin{prob}
Generalize Theorem \ref{thm p square} for rank $r \geq 3$.

\end{prob}

We also have the following ``inverse type'' result from \cite{BajMat:2014a} regarding subsets that achieve $\rho_{\pm} \left(G, m, h   \right)$.  Given a group $G$ and a positive integer $m \leq |G|$, we define a certain collection ${\cal A}(G,m)$ of $m$-subsets of $G$.  
We let 
\begin{itemize}
  \item $\mathrm{Sym}(G,m)$ be the collection of {\em symmetric} $m$-subsets of $G$, that is, $m$-subsets $A$ of $G$ for which $A=-A$;
  \item $\mathrm{Nsym}(G,m)$ be the collection of {\em near-symmetric} $m$-subsets of $G$, that is, $m$-subsets $A$ of $G$ that are not symmetric, but for which $A\setminus \{a\}$ is symmetric for some $a \in A$; 
  \item $\mathrm{Asym}(G,m)$ be the collection of {\em asymmetric} $m$-subsets of $G$, that is, $m$-subsets $A$ of $G$ for which $A \cap (-A)=\emptyset$.  
  \end{itemize}  We then let 
$${\cal A}(G,m)=\mathrm{Sym}(G,m) \cup \mathrm{Nsym}(G,m)\cup \mathrm{Asym}(G,m).$$  In other words, ${\cal A}(G,m)$ consists of those $m$-subsets of $G$ that have exactly $m$, $m-1$, or $0$ elements whose inverse is also in the set. 

\begin{thm} [Bajnok and Matzke; cf.~\cite{BajMat:2014a}] \label{symmetry thm}\index{Bajnok, B.}\index{Matzke, R.}
For every $G$, $m$, and $h$, we have
$$\rho_{\pm} (G,m,h)= \min \{|h_{\pm} A| \; : \; A \in {\cal A}(G,m)\}.$$

\end{thm}
We should add that each of the three types of sets are essential as can be seen by examples (cf.~\cite{BajMat:2014a}).  However, we do not know the answers to the following problems:

\begin{prob}
Classify each situation where the minimum signed sumset size is achieved by symmetric sets.
\end{prob}

\begin{prob}
Classify each situation where the minimum signed sumset size is achieved by near-symmetric sets.
\end{prob}

\begin{prob}
Classify each situation where the minimum signed sumset size is achieved by asymmetric sets.
\end{prob}

\begin{prob}
Classify each situation where the minimum signed sumset size is achieved by sets that are not in ${\cal A}(G,m)$.
\end{prob}

\subsection{Limited number of terms} \label{2minUSlimited}

Here we consider, for a given group $G$, positive integer $m$ (with $m \leq n=|G|$), and nonnegative integer $s$, $$\rho_{\pm} (G,m,[0,s]) = \mathrm{min} \{ |[0,s]_{\pm}A|  \mid A \subseteq G, |A|=m\},$$ that is, the minimum size of $[0,s]_{\pm}A$ for an $m$-element subset $A$ of $G$.

We note that we don't have a version of Proposition \ref{shifting sets} for signed sumsets, so we are not able to reduce this entire section to Section \ref{2minUSfixed}.  (However, one may be able to apply similar techniques.) 

It is easy to see that, for every $m$, we have $$\rho_{\pm} (G,m,[0,0]) =1,$$  and for every $s$, we have $$\rho_{\pm} (G,1,[0,s]) =1.$$

Furthermore, we can evaluate $\rho_{\pm} (G,m,[0,1])$, as follows.  Note that $$[0,1]_{\pm}A=A \cup (-A) \cup \{0\},$$ so $|[0,1]_{\pm}A| \geq m$ with equality if, and only if, $A$ is symmetric (that is, $A=-A$) and $0 \in A$.  Let us see the conditions that allow for such a set $A$.  We can partition $G$ into the pairwise disjoint union of four (potentially empty) parts: $\{0\}$, $\mathrm{Ord}(G,2)$, $K$, and $-K$.  Therefore, if 
$\mathrm{Ord}(G,2) \neq \emptyset$, we can take $A$ to be the set containing 0, together with the right number of pairs of the form $\pm k$ with $k \in K$, and some elements from $\mathrm{Ord}(G,2)$.  Likewise, if $\mathrm{Ord}(G,2) = \emptyset$ but $m$ is odd, we can let $A$ be 0 together with $(m-1)/2$ elements of $K$ together with their inverses.

This leaves only the case when $\mathrm{Ord}(G,2) = \emptyset$ (that is, $n$ is odd), and $m$ is even.  It is easy to find sets $A$ with $|[0,1]_{\pm}A| =m+1$, and we can see that that is the best we can do.  Indeed, we will have either $0 \not \in A$ (but $0 \in [0,1]_{\pm}A$), or an element $a \in A$ for which $-a \not \in A$ (but $-a \in [0,1]_{\pm}A$).

Summarizing our results thus far, we have the following.

\begin{prop} \label{rhoUStrivial} Let $G$ be an abelian group $G$ of order $n$, and let $m$ and $s$ positive integers.
\begin{itemize}
\item For all $m$ we have $$\rho_{\pm} (G,m,[0,0]) = 1.$$
\item For all $s$ we have $$\rho_{\pm} (G,1,[0,s]) =1.$$
\item For all $m$ we have $$\rho_{\pm} (G,m,[0,1]) = 
\left \{
\begin{array}{cll}
m & \mbox{if} & n \; \mbox{is even or $m$ is odd,} \\ \\
m+1 & \mbox{if} & n \; \mbox{is odd and $m$ is even.}
\end{array} \right.$$
\end{itemize}
\end{prop} 

For values of $s \geq 2$ and $m \geq 2$, we have no exact values for $\rho_{\pm} (G,m,[0,s])$ in general.  

However, Matzke\index{Matzke, R.} in \cite{Mat:2013a} provided the following upper bound for $\rho_{\pm} (G,m,[0,s])$ for the case when $G$ is cyclic.  First, we introduce some notations.  For a positive integer $k$ we let $v(k)$ denote the highest power of 2 that is a divisor of $k$.  For given positive integers $n$ and $m $ we then define 
$$D_1(n) =\{ d \in D(n) \mid v(n) \geq v(d  \lceil m/d \rceil)\}$$ and
 $$D_2(n) =\{ d \in D(n) \mid v(n) < v(d  \lceil m/d \rceil)\}.$$ Finally, we set
$$u_{\pm}(n,m,[0,s])=\min\{\min \{ f_d(m,s) \mid d \in D_1(n)\}, \; \min \{ f_d(m+d,s) \mid d \in D_2(n)\}\}.$$
(Note that for $d \in D_2(n)$ we have $d  \lceil m/d \rceil < n$ and thus $m+d \leq n$.)

\begin{thm} [Matzke; cf.~\cite{Mat:2013a}] \label{matzke upper for [0,s]}\index{Matzke, R.}
With the notations just introduced, we have
$$\rho_{\pm} (\mathbb{Z}_n,m,[0,s]) \leq u_{\pm}(n,m,[0,s]).$$

\end{thm}

Furthermore, Matzke\index{Matzke, R.} believes that equality holds in Theorem \ref{matzke upper for [0,s]}:

\begin{conj} [Matzke; cf.~\cite{Mat:2013a}] \label{Matzke limited conj}\index{Matzke, R.}
With the notations introduced above, we have
$$\rho_{\pm} (\mathbb{Z}_n,m,[0,s]) = u_{\pm}(n,m,[0,s]).$$
\end{conj}

\begin{prob}
Prove or disprove Conjecture \ref{Matzke limited conj}.  

\end{prob}

As Grynkiewicz\index{Grynkiewicz, D. J.} pointed out (cf.~\cite{Gry:2017a}), we can prove Conjecture \ref{Matzke limited conj}  for prime values of $n$:

\begin{thm} \label{Matzke limited prime}
For odd prime values of $p$ we have $$\rho_{\pm} (\mathbb{Z}_p,m,[0,s]) = u_{\pm}(p,m,[0,s])=\min\{p,2s  \lfloor m/2 \rfloor +1\}.$$

\end{thm}  

We present the short proof on page \pageref{proof of Matzke limited prime}.

We do not know much about the value of $\rho_{\pm}(G,m,[0,s])$ for noncyclic groups:

\begin{prob}
Find the value of  (or, at least, find good bounds for) $\rho_{\pm}(G,m,[0,s])$ for noncyclic groups $G$ and integers $m$ and $s$.
\end{prob}

\subsection{Arbitrary number of terms} \label{2minUSarbitrary}

This subsection is identical to Subsection \ref{2minUarbitrary}.

\section{Restricted sumsets} \label{2minR}

Our goal in this section is to investigate the quantity $$\rho\hat{\;} (G,m,H) =\mathrm{min} \{ |H\hat{\;}A|  \mid A \subseteq G, |A|=m\}$$ where $H\hat{\;}A$ is the union of all $h$-fold restricted sumsets $h\hat{\;}A$ for $h \in H$.  We consider three special cases: when $H$ consists of a single nonnegative integer $h$, when $H$ consists of all nonnegative integers up to some value $s$, and when $H$ is the entire set of nonnegative integers.

\subsection{Fixed number of terms} \label{2minRfixed}

Here we consider $$\rho\hat{\;} (G,m,h) = \mathrm{min} \{ |h\hat{\;}A|  \mid A \subseteq G, |A|=m\},$$ that is, the minimum size of an $h$-fold restricted sumset of an $m$-element subset of $G$.

Note that, when $h>m$, we obviously have $h\hat{\;}A = \emptyset$ for every $m$-subset $A$ of $G$, thus we may assume that $h \leq m$.  Clearly, for all $A$ we have $0\hat{\;}A=\{0\}$  and $1\hat{\;}A=A$, so $$\rho\hat{\;} (G,m,0) = 1$$ and 
$$\rho\hat{\;} (G,m,1) = m.$$

Furthermore, as we saw in Section \ref{sectionmaxsumsetsizeRfixed}, for all $m$, $h$, and $m$-subsets $A$ of $G$ we have
$$|(m-h)\hat{\;}A|=|h\hat{\;}A|,$$ and thus
$$\rho\hat{\;} (G,m,m-h)=\rho\hat{\;} (G,m,h).$$  Therefore, it suffices to study cases when $$2 \leq h \leq \left\lfloor \frac{m}{2} \right \rfloor.$$

By these considerations, we see that the only case with $m \leq 4$ in which $\rho\hat{\;} (G,m,h)$ is not immediate is $m=4$ and $h=2$, for which we prove the following result.

\begin{prop} \label{rhohat m=4}
Let $G$ be an abelian group with $\mathrm{Ord}(G,2)$ as the set of its elements of order 2.  We have
$$\rho\hat{\;} (G,4,2)=
\left\{\begin{array}{ll}
3  & \mbox{if $\; |\mathrm{Ord}(G,2)| \geq 2$,}\\ 
4  & \mbox{if $\; |\mathrm{Ord}(G,2)| =1$,}\\ 
5  & \mbox{if $\; |\mathrm{Ord}(G,2)| =0$.}  
\end{array}\right.$$

\end{prop}

The proof of Proposition \ref{rhohat m=4} can be found on page \pageref{proof of rhohat m=4}; the proposition on page \pageref{h=2 m=4 size} also provides a complete characterization of 4-subsets $A$ of $G$ with $|2 \hat{\;} A|$ attaining all possible values.  We can rephrase the conditions in Proposition \ref{rhohat m=4} by using the invariant factorization of $G$.  Let
$$G \cong \mathbb{Z}_{n_1} \times \cdots \times \mathbb{Z}_{n_r}$$ with $2 \leq n_1$, $n_1|n_2| \cdots | n_r$.  We then have
$$\rho\hat{\;} (G,4,2)=
\left\{\begin{array}{ll}
3  & \mbox{if} \; r \geq 2 \; \mbox{and} \; n_{r-1} \; \mbox{is even},\\ 
4  & \mbox{if} \; r \geq 2 \; \mbox{and} \; n_{r-1} \; \mbox{is odd but} \; n_r \; \mbox{is even, or} \; r=1 \; \mbox{and} \; n_r \; \mbox{is even},\\ 
5  & \mbox{if} \; n_{r} \; \mbox{is odd}.  
\end{array}\right.$$
(Note that $n_r$ is even if, and only if, $n=|G|$ is even.)

As Proposition \ref{rhohat m=4} demonstrates, the value of $\rho\hat{\;} (G,4,h)$ may be less than $4$. (This is in contrast to $\rho (G,m,h) \geq m$; see Proposition \ref{rho (G,m,h)>=m}.)  However, we can easily prove that  $$\rho\hat{\;} (G,5,2) \geq 5,$$ as follows.  Let $A=\{a_1,\dots,a_5\}$ (with $|A|=5$); we then have
$$2\hat{\;}A=\{a_i+a_j  \mid 1 \leq i < j \leq 5\}.$$  It is easy to see that $$B=\{a_1+a_i  \mid i=2,3,4,5\}$$ is a 4-subset of $2\hat{\;}A$, thus $|2\hat{\;}A| \geq 4$; furthermore, for $2\hat{\;}A$ to have size 4, we must have (among other things) $a_2+a_3 \in B$ and $a_4+a_5 \in B$.  Since $$a_2+a_3 \not \in \{a_1+a_2,a_1+a_3\},$$ we have, w.l.o.g., $$a_2+a_3=a_1+a_4;$$ similarly, $$a_4+a_5=a_1+a_2.$$  Adding these two equations and cancelling yields $$a_3+a_5=2a_1,$$ and, therefore, $2a_1 \in B$, but that cannot happen.  Therefore, $|2\hat{\;}A| \geq 5$.

In fact, Carrick\index{Carrick, E.} evaluated $\rho\hat{\;} (G,5,2)$ for every group $G$:

\begin{thm} [Carrick; cf.~\cite{Car:2015a}] \label{carrick}\index{Carrick, E.}

For an abelian group $G$ of order $n$ we have
$$\rho\hat{\;} (G,5,2)=
\left\{\begin{array}{ll}
5  & \mbox{if} \; n \; \mbox{is divisible by 5},\\ 
6  & \mbox{if} \; n \; \mbox{is not divisible by 5 but is divisible by 6},\\ 
7  & \mbox{otherwise}.  
\end{array}\right.$$

\end{thm}

The evaluation of $\rho\hat{\;} (G,m, 2)$ gets more difficult as $m$ increases.  We still offer:

\begin{prob}  \label{rho hat m=6}
Evaluate $\rho\hat{\;} (G,6,2)$ for every finite abelian group $G$.
\end{prob}

The values of $\rho\hat{\;} (G,m, h)$ are largely unknown in general; we first attempt to evaluate them in the case when $G$ is cyclic.  
The general question regarding $\rho\hat{\;} (\mathbb{Z}_n,m,h)$ remains open:

\begin{prob}
Find the exact value of $\rho\hat{\;} (\mathbb{Z}_n,m,h)$ for all $n$, $m$, and $h$.

\end{prob}

Below we summarize what we know about $\rho\hat{\;} (\mathbb{Z}_n,m,h)$.

We start by wondering what the minimum possible value of $\rho\hat{\;} (\mathbb{Z}_n,m,h)$ is.  Suppose that $1 \leq h \leq m-1$.  Let $A=\{a_1,\dots,a_m\}$, but assume that the elements (as integers) are between 0 and $n-1$ and that they are in increasing order.  Now consider the sets  
$$A_1=\{a_1+a_2+\cdots+a_{h-1}+a_i  \mid i=h,h+1,\dots,m\}$$ and 
$$A_2=\{a_1+a_2+\cdots+a_h-a_j+a_m  \mid j=1,2,\dots,h-1\}.$$  Note that $A_1 \subseteq h\hat{\;}A$ and (since $h \leq m-1$) $A_2 \subseteq h\hat{\;}A$.  Clearly, $|A_1|+|A_2|=m$; furthermore, $A_1$ and $A_2$ are disjoint, since (as integers), we have
$$a_1+a_2+\cdots+a_h < \cdots < a_1+\cdots+a_{h-1}+a_m < a_1+\cdots + a_{h-2}+a_h+a_m < \cdots < a_2+\cdots+a_h+a_m,$$ and the smallest and largest sums differ by $a_m-a_1<n$, and thus the $m$ elements are distinct in $\mathbb{Z}_n$.  Therefore, we have the following.

\begin{prop} \label{rhohat(G,m,h)>=mcyclic}
For all $1 \leq h \leq m-1$ we have $\rho\hat{\;} (\mathbb{Z}_n,m,h) \geq m$.

\end{prop}
  
Clearly, equality holds for $h=1$ and $h=m-1$.  Another obvious example for $\rho\hat{\;} (\mathbb{Z}_n,m,h) = m$ is the case when $m$ is a divisor of $n$: the (unique) subgroup $A$ of $\mathbb{Z}_n$ of order $m$ has $|h\hat{\;}A|=|A|=m$; more generally, $A$ can be any coset of this subgroup.  As it turns out (see below), there are no other cases where equality holds for the cyclic group.

We now turn to the question of finding upper bounds for $\rho\hat{\;} (\mathbb{Z}_n,m,h)$.  An obvious upper bound is that
$$\rho\hat{\;} (\mathbb{Z}_n,m,h) \leq \rho (\mathbb{Z}_n,m,h),$$ where, according to Theorem \ref{rho (G,m,h)}, 
$$\rho (\mathbb{Z}_n,m,h)=u(n,m,h)=\min \left \{ \left( h \cdot \left \lceil \frac{m}{d} \right \rceil -h+1 \right) \cdot d  \mid d \in D(n) \right \}.$$

Recall that we provided, for each $d \in D(n)$, an $m$-subset $A_d(n,m)$ of $\mathbb{Z}_n$ with $h$-fold sumset size 
$$|hA_d(n,m)|=\min\{n,f_d,hm-h+1\},$$ where
$$f_d=\left( h \cdot \left \lceil \frac{m}{d} \right \rceil -h+1 \right) \cdot d$$  (see Proposition \ref{|hA_d(n,m)|}).  To get an upper bound for $\rho\hat{\;} (\mathbb{Z}_n,m,h)$, we now compute the size of the restricted $h$-fold sumset of $A_d(n,m)$.

  We recall the construction from Section \ref{sectionminsumsetsizeUfixed}, as follows.  For a divisor $d$ of $n$, consider the (unique) subgroup $H$ of order $d$ of $\mathbb{Z}_n$,  namely, $$H=\left \{j \cdot \frac{n}{d}  \mid j=0,1,2,\dots,d-1 \right\},$$ and then let $A_d(n,m)$ be a certain subset of the ``first'' $\left \lceil \frac{m}{d} \right \rceil$ cosets of $H$.  Namely, we set 
$$A_d(n,m)=\bigcup_{i=0}^{c-1} (i+H) \cup \left\{\left \lceil \frac{m}{d} \right \rceil-1 + j \cdot \frac{n}{d}  \mid j=0,1,2,\dots,k-1 \right\} $$ where
$$c=\left \lceil \frac{m}{d} \right \rceil-1$$ and $$k=m-d c$$ assures that $A=A_d(n,m)$ has size $m$.  We also see that $1 \leq k \leq d$ and thus we take at least 1 but at most $d$ elements of the coset $c+H$;  on the other hand, when $c=0$ (that is, when $m \leq d$), we see that $$\bigcup_{i=0}^{c-1} (i+H)=\emptyset,$$ so $A$ lies entirely within a single coset and forms the arithmetic progression
$$\left \{j \cdot \frac{n}{d}  \mid j=0,1,2,\dots,m-1 \right\}.$$

Recall also the function $u\hat{\;}(n,m,h)$ from Section \ref{0.3.5.2}.  For that, we introduced the notations  $h=qd+r$ with
$$q=\left \lceil \frac{h}{d} \right \rceil-1$$ and $$r=h-d q,$$
then set
$$f\hat{_d}(n,m,h)=
\left\{\begin{array}{ll}
\min\{n, f_d, hm-h^2+1\}  & \mbox{if $h \leq \min\{k,d-1\}$,}\\ \\
\min\{n, hm-h^2+1 - \delta_d\} & \mbox{otherwise;}\\  
\end{array}\right.$$
where $\delta_d$ is a ``correction term'' defined as
$$\delta_d(n,m,h)=\left\{ 
\begin{array}{cl}
(k-r)r-(d-1)  & \mbox{if $r<k$,}\\ 
(d-r)(r-k)-(d-1)  & \mbox{if $k<r<d$,}\\
d-1  & \mbox{if $k=r=d$,}\\
0  & \mbox{otherwise.}  
\end{array}\right.$$

We can then prove the following result. 

\begin{thm} [Bajnok; cf.~\cite{Baj:2013a}] \label{uhattheorem}\index{Bajnok, B.}
Suppose that $1 \leq h < m \leq n$.  We then have $$|h\hat{\;}A_d(n,m)| = f\hat{_d}(n,m,h).$$
\end{thm}

The proof of Proposition \ref{uhattheorem} can be found on page \pageref{proofofuhattheorem}.

Obviously, $|h\hat{\;}A_d(n,m)|$ provides an upper bound for $\rho \hat{\;} (\mathbb{Z}_n,m,h)$ for every $d \in D(n)$.  In Section \ref{0.3.5.2} we also defined  $$u\hat{\;}(n,m,h)=\min \{ f\hat{_d}(n,m,h)  \mid d \in D(n) \},$$ with which Theorem \ref{uhattheorem} implies:

\begin{cor} \label{uhatcorollary}
Suppose that $1 \leq h < m \leq n$.  For the cyclic group of order $n$ we have $$\rho \hat{\;} (\mathbb{Z}_n,m,h) \leq u\hat{\;}(n,m,h).$$
\end{cor}

Manandhar\index{Manandhar, D.} in \cite{Man:2012a} has computed the value of $\rho \hat{\;} (\mathbb{Z}_n,m,h)$ for all $n \leq 25$ and $m \leq 12$; Malec\index{Malec, M.} in \cite{Mal:2013a} extended this search to the range $n \leq 40$ and all $m$.  They found that in all cases, $\rho \hat{\;} (\mathbb{Z}_n,m,h)$ agrees with $u\hat{\;}(n,m,h)$, with the following exceptions:
$$\begin{array}{||c|c|c||c|c||}  \hline \hline 
n & m & h &   u\hat{\;}(n,m,h) & \rho \hat{\;} (\mathbb{Z}_n,m,h) \\ \hline \hline 
10 & 6 & 3 & 10 & 9 \\ \hline
12 & 7 & 2 & 11 & 10 \\ \hline
15 & 8 & 3 & 15 & 14 \\ \hline
20 & 6 & 3 & 10 & 9 \\ \hline 
20 & 11 & 2 & 19 & 18 \\ \hline
24 & 13 & 2 & 23 & 22 \\ \hline
28 & 15 & 2 & 27 & 26 \\ \hline
35 & 12 & 5 & 35 & 34 \\ \hline
40 & 11 & 2 & 19 & 18 \\ \hline \hline
\end{array}$$
(Only those with $h \leq \left \lfloor \frac{m}{2} \right \rfloor$ are listed.)

From these data it seems that $u\hat{\;}(n,m,h)$ is a remarkably good upper bound for $\rho \hat{\;}(\mathbb{Z}_n, m, h)$: indeed, $u \hat{\;}(n,m,h)$ agrees with $\rho \hat{\;}(\mathbb{Z}_n, m, h)$ in the overwhelming majority (over 99\%) of cases, and when it does not, it differs by only 1.  In light of this we pose the following: 

\begin{prob} \label{problem rho hat < u hat}
Classify all situations when $$\rho \hat{\;}(\mathbb{Z}_n, m, h)< u \hat{\;}(n,m,h);$$ in particular, how much smaller than $u \hat{\;}(n,m,h)$ can $\rho \hat{\;}(\mathbb{Z}_n, m, h)$ be?  

\end{prob}
As we just mentioned, in all cases that we are aware of, we have $\rho \hat{\;}(\mathbb{Z}_n, m, h)= u \hat{\;}(n,m,h)$ or $\rho \hat{\;}(\mathbb{Z}_n, m, h)= u \hat{\;}(n,m,h)-1$.
 
Let us mention three---as it turns out, quite predicative---examples for the case when $\rho \hat{\;}(\mathbb{Z}_n, m, h)= u \hat{\;}(n,m,h)-1$.
\begin{itemize}
\item $u \hat{\;}(12,7,2)=11$, but $\rho \hat{\;}(\mathbb{Z}_{12}, 7, 2)=10$ as shown by the set $$C_1=\{0,4\}\cup\{1,5,9\}\cup\{6,10\}.$$ \label{example c_1} 
\item $u \hat{\;}(10,6,3)=10$, but $\rho \hat{\;}(\mathbb{Z}_{10}, 6, 3)=9$ as shown by the set $$C_2=\{0,2,4,6\}\cup\{7,9\}.$$ 
\item $u \hat{\;}(15,8,3)=15$, but $\rho \hat{\;}(\mathbb{Z}_{15}, 8, 3)=14$ as shown by the set $$C_3=\{0,3,6,9\}\cup\{10,13,1,4\}.$$ 
\end{itemize}
Needless to say, we chose to represent our sets in the particular formats above for a reason.  In fact, all known situations where $\rho \hat{\;}(\mathbb{Z}_n, m, h)< u \hat{\;}(n,m,h)$ can be understood by a particular modification of our sets $A_d(n,m)$, as we now describe.

Observe that the $m$ elements of $A_d(n,m)$ are within $\lceil m/d \rceil$ cosets of the order $d$ subgroup $H$ of $\mathbb{Z}_n$, and at most one of these cosets does not lie entirely in $A_d(n,m)$.  We now consider the situation when the $m$ elements are still within $\lceil m/d \rceil$ cosets of $H$, but exactly two of these cosets don't lie entirely in our set.  In order to do so, we write $m$ in the form
$$m=k_1+(c-1)d+k_2$$ for some positive integers $c$, $k_1$, and $k_2$; we assume that $k_1<d$, $k_2<d$, but $k_1+k_2>d$.  We then are considering $m$-subsets $B$ of $\mathbb{Z}_n$ of the form \label{def of B_d sets}
$$B=B_d(n,m;k_1,k_2,g,j_0)=B' \cup \bigcup_{i=1}^{c-1} (ig+H) \cup B'',$$
where $H$ is the subgroup of $\mathbb{Z}_n$ with order $d$, $g$ is an element of $\mathbb{Z}_n$, $B'$ is a proper subset of $H$ given by 
$$B'=\left\{ j \cdot \frac{n}{d} \; : \;  j=0,1,\dots,k_1-1 \right\},$$ and $B''$ is a proper subset of $cg+H$ of the form 
$$B''=\left\{ cg+(j_0+j) \cdot \frac{n}{d} \; : \;  j=0,1,\dots,k_2-1 \right\}$$ for some integer $j_0$ with $0 \leq j_0 \leq d-1$.

It turns out that our set $B$ (under some additional assumptions to be made precise in \cite{Baj:2013a})  has the potential to have a restricted $h$-fold sumset of size less than $u \hat{\;} (n,m,h)$ in only three cases:
\begin{itemize}
\item $h=2$, $m-1$ is not a power of 2, and $n$ is divisible by $2m-2$; 
\item $h=3$, $m=6$, and $n$ is divisible by 10; or
\item $h$ is odd, $m+2$ is divisible by $h+2$, and $n$ is divisible by $hm-h^2$.
\end{itemize}
In particular, we have the following.

\begin{thm} [Bajnok; cf.~\cite{Baj:2013a}] \label{special sets listed}\index{Bajnok, B.}
Let $B_d(n,m;k_1,k_2,g,j_0)$ be the $m$-subset of $\mathbb{Z}_n$ defined above, and let $h$ be a positive integer with $h \leq m-1$.  
\begin{itemize}
\item If $h=2$, $m-1$ is not a power of 2, $n$ is divisible by $2m-2$, and $d$ is an odd divisor of $m-1$ with $d>1$, then 
$$B_d \left(n,m;\tfrac{d+1}{2}, \tfrac{d+1}{2}, \tfrac{n}{2m-2},\tfrac{d-1}{2} \right)$$  has a restricted 2-fold sumset of size $2m-4$.
\item If $h=3$, $m=6$, and $n$ is divisible by 10, then
$$B_5 \left(n,6;4,2, \tfrac{n}{10},3 \right)$$  has a restricted 3-fold sumset of size $3m-9=9$.

\item If $h$ is odd, $m+2$ is divisible by $h+2$, and $n$ is divisible by $hm-h^2$, then
$$B_{h+2} \left(n,m;h+1,h+1, \tfrac{n}{hm-h^2},\tfrac{h+3}{2} \right)$$  has a restricted $h$-fold sumset of size $hm-h^2-1$.
\end{itemize}

\end{thm}
Our examples above demonstrate the three cases of Theorem \ref{special sets listed} in order: we have
$$C_1=B_3(12,7;2,2,1,1),$$
$$C_2=B_5(10,6;4,2,1,3),$$ and
$$C_3=B_5(15,8;4,4,1,3).$$
The proof of Theorem \ref{special sets listed} is an easy exercise; it also follows from \cite{Baj:2013a} where we verify that, in a certain sense that we make precise, there are no other such sets.

While we seem far away from knowing the value of $\rho \hat{\;} (\mathbb{Z}_n,m,h)$ in general, we do have the answer in the case when $n$ is prime.  
Recall that, by Proposition \ref{uhat(p,m,h)}, for a prime $p$ we have  
$$u\hat{\;}(p,m,h)=\min \{ p, hm-h^2+1\},$$ thus Corollary \ref{uhatcorollary}, when $n$ equals a prime number $p$, simplifies to 
$$\rho \hat{\;} (\mathbb{Z}_p,m,h)  \leq \min \{ p, hm-h^2+1\}.$$  The conjecture that equality holds here 
 has been known since the 1960s as the {\em Erd\H{o}s--Heilbronn Conjecture}\index{Erd\H{o}s, P.}\index{Heilbronn, H.}  (not mentioned in \cite{ErdHei:1964a} but in \cite{ErdGra:1980a}).  Three decades later, Dias Da Silva and Hamidoune \cite{DiaHam:1994a}\index{Dias Da Silva, J. A.}\index{Hamidoune, Y. O.} succeeded in proving the Erd\H{o}s--Heilbronn Conjecture and thus we have

\begin{thm} [Dias Da Silva and Hamidoune] \label{Dias Da Silva and Hamidoune}\index{Dias Da Silva, J. A.}\index{Hamidoune, Y. O.}
For a prime $p$ and integers $1 \leq h \leq m \leq p$ we have $$\rho\hat{\;} (\mathbb{Z}_p,m,h) = \mathrm{min} \{p, hm-h^2+1\}.$$

\end{thm}  
(The result was reestablished and extended, using different methods, by Alon,\index{Alon, N.} Nathanson,\index{Nathanson, M. B.} and Ruzsa;\index{Ruzsa, I.} see \cite{AloNatRuz:1995a}, \cite{AloNatRuz:1996a}, and \cite{Nat:1996a}.)

For composite values of $n$, the evaluation of $\rho\hat{\;} (\mathbb{Z}_n,m,h)$ seems considerably more difficult; we attempt to summarize here what is known for $h=2$ and $h=3$.

For $h=2$, recall Proposition \ref{uhatforh=2}:
$$u\hat{\;}(n,m,2)=\left\{\begin{array}{ll}
\min\{u(n,m,2),2m-4\}  & \mbox{if $n$ and $m$ are both even,}\\ \\
\min\{u(n,m,2),2m-3\} & \mbox{otherwise.} 
\end{array}\right.$$
We can then use Theorem \ref{special sets listed} to determine the values of $n$ and $m$ that allow for an improvement over $u\hat{\;} (n,m,2)$.  Only one case applies: when $m-1$ is not a power of 2 and $n$ is divisible by $2m-2$, in which a special set $B$ exists with $|2\hat{\;} B|=2m-4$.  Therefore, using Theorem \ref{rho (G,m,h)} as well, we have the following result.  
\begin{cor} \label{cor h=2}
For all positive integers $n$ and $m$ with $3 \leq m \leq n$ we have
$$\rho\hat{\;} (\mathbb{Z}_n,m,2) \leq \left\{
\begin{array}{ll}
\min\{\rho(\mathbb{Z}_n,m,2), 2m-4\} & \mbox{if} \; 2|n \; \mbox{and} \; 2|m, \; \mbox{or} \\ 
& (2m-2)|n \; \mbox{and} \; \not \exists k \in \mathbb{N}, m=2^k+1; \\ \\
\min\{\rho(\mathbb{Z}_n,m,2), 2m-3\} & \mbox{otherwise.} 
\end{array}
\right.$$

\end{cor}

We have performed a computer search for all $m$-subsets of $\mathbb{Z}_n$ with $3 \leq m \leq n \leq 40$, and in each case we found that equality holds in Corollary \ref{cor h=2}.

\begin{conj} \label{conj rhohatforh=2}
For all $n$ and $m$, we have equality in Corollary \ref{cor h=2}.

\end{conj}

\begin{prob}
Prove or disprove Conjecture \ref{conj rhohatforh=2}.

\end{prob}

We can carry out a similar analysis for the case of $h=3$, by first recalling from Proposition \ref{uhatforh=3}: 
$$u\hat{\;} (n,m,3)=\left\{
\begin{array}{ll}
\min\{u(n,m,3), 3m-3-\gcd(n,m-1)\} & \mbox{if} \; \gcd(n,m-1) \geq 8; \\ \\
\min\{u(n,m,3), 3m-10\} & \mbox{if} \; \gcd(n,m-1) = 7, \; \mbox{or} \\
& \gcd(n,m-1) \leq 5, \; 3|n, \; \mbox{and} \; 3|m; \\ \\
\min\{u(n,m,3), 3m-9\} & \mbox{if} \; \gcd(n,m-1) = 6; \\ \\
\min\{u(n,m,3), 3m-8\} & \mbox{otherwise.} 
\end{array}
\right.$$

We then examine Theorem \ref{special sets listed} for the case $h=3$ to see if we can do better.  We find two such instances: when $n$ is divisible by 10 and $m=6$, and when $n$ is divisible by $3m-9$ and $m-3$ is divisible by 5.  Observe that, in the latter case, we may assume that $m$ is even, since otherwise $d_0=(3m-9)/2 \in D(n)$ and thus
$$u\hat{\;} (n,m,3) \leq u(n,m,3) \leq f_{d_0} =d_0 \leq 3m-10.$$  Therefore, we have the following.

\begin{cor} \label{cor h=3}
Let $n$ and $m$ be positive integers, $4 \leq m \leq n$, and set $d_0=\gcd(n,m-1)$.  We then have
$$\rho\hat{\;} (\mathbb{Z}_n,m,3) \leq \left\{
\begin{array}{ll}
\min\{u(n,m,3), 3m-3-d_0\} & \mbox{if} \; d_0 \geq 8; \\ \\
\min\{u(n,m,3), 3m-10\} & \mbox{if} \; d_0 = 7, \; \mbox{or} \\
&  d_0 \leq 5, \; 3|n, \; \mbox{and} \; 3|m, \; \mbox{or} \\
&  d_0 \leq 5, \; (3m-9)|n, \; \mbox{and} \; 5|(m-3); \\ \\
\min\{u(n,m,3), 3m-9\} & \mbox{if} \; d_0 = 6, \; \mbox{or} \\
& m=6 \; \mbox{and} \; 10|n \; \mbox{but} \; 3 \not | n; \\ \\
\min\{u(n,m,3), 3m-8\} & \mbox{otherwise.} 
\end{array}
\right.$$

\end{cor}

Motivated by computational data mentioned above, we have:

\begin{conj} \label{conj rhohatforh=3}
For all $n$ and $m$, we have equality in Corollary \ref{cor h=3}.

\end{conj}

\begin{prob}
Prove or disprove Conjecture \ref{conj rhohatforh=3}.

\end{prob}

Few other results or even conjectures are known for the exact value of $\rho\hat{\;} (\mathbb{Z}_n,m,h)$ or (especially) for $\rho\hat{\;} (G,m,h)$ in general.  One such result is the following:

\begin{prop} \label{upper bound for restricted 2-critical} 

Let $G$ be an abelian group of order $n$.  We have $\rho \hat{\;} (G,m,2)=n$ if, and only if,
$$m \geq \frac{n+|\mathrm{Ord}(G,2)|+3}{2}.$$ In particular, 
\begin{itemize}
  \item $\rho \hat{\;} (G,n,2)=n$ if, and only if, $G$ is not isomorphic to the elementary abelian 2-group; and 
  \item $\rho \hat{\;} (\mathbb{Z}_n,m,2)=n$ if, and only if,
$m \geq \left \lfloor n/2 \right \rfloor+2.$
\end{itemize}

\end{prop}
Note that $n+|\mathrm{Ord}(G,2)|+3$ is always even.  We should also mention that Proposition \ref{upper bound for restricted 2-critical} appeared several times in the literature: cf.~\cite{RotLem:1992a} by\index{Roth, R. M.}\index{Lempel, A.} Roth and Lempel; \cite{Chi:1993a} (for the ``if'' part only) and \cite{Chi:2002a}\index{Chiaselotti, G.} by Chiaselotti; and \cite{Baj:2014a}\index{Bajnok, B.} by Bajnok.     

According to Proposition \ref{upper bound for restricted 2-critical}, the largest value of $m$ for which $\rho \hat{\;} (G,m,2)$ is less than $n$ equals $$m_0=\frac{n+|\mathrm{Ord}(G,2)|+1}{2}.$$  It is then an interesting question to find $\rho \hat{\;} (G,m_0,2)$.  By Proposition \ref{upper bound for restricted 2-critical}, $$\rho \hat{\;} (G,m_0,2) \leq n-1.$$  We prove the following:

\begin{prop}  \label{rho hat for chi hat -1}
For $$m_0=\frac{n+|\mathrm{Ord}(G,2)|+1}{2}$$ we have 
$$\rho \hat{\;} (G,m_0,2) \leq \left\{
\begin{array}{ll}
n-1 & \mbox{if} \; \exists k \in \mathbb{N}, n=2^k;\\ \\
n-2 & \mbox{otherwise.}
\end{array}
\right.$$ 
In particular,
$$\rho \hat{\;} \left(\mathbb{Z}_n, \left \lfloor \frac{n}{2} \right \rfloor+1, 2 \right) \leq \left\{
\begin{array}{ll}
n-1 & \mbox{if} \; \exists k \in \mathbb{N}, n=2^k;\\ \\
n-2 & \mbox{otherwise.}
\end{array}
\right.$$
\end{prop}
As we mentioned, the first case of our result holds by Proposition \ref{upper bound for restricted 2-critical}; the proof of the second case can be found on page \pageref{proof of rho hat for chi hat -1}.

In 2002, Gallardo, Grekos,\index{Gallardo, L.}\index{Grekos, G.}\index{Habsieger, L.}\index{Hennecart, F.}\index{Landreau, B.}\index{Plagne, A.}  et al.~proved that, when $G$ is cyclic, Proposition \ref{rho hat for chi hat -1} holds with equality:

\begin{thm} [Gallardo, Grekos, et al.; cf.~\cite{GalGre:2002a}] \label{thm g G et al}\index{Gallardo, L.}\index{Grekos, G.}\index{Habsieger, L.}\index{Hennecart, F.}\index{Landreau, B.}\index{Plagne, A.} 
For every positive integer $n \geq 2$ we have
$$\rho \hat{\;} \left(\mathbb{Z}_n, \left \lfloor \frac{n}{2} \right \rfloor+1, 2 \right) = \left\{
\begin{array}{ll}
n-1 & \mbox{if} \; \exists k \in \mathbb{N}, n=2^k;\\ \\
n-2 & \mbox{otherwise.}
\end{array}
\right.$$
\end{thm}
 Note that Conjecture \ref{conj rhohatforh=2}, once established, would generalize Theorem  \ref{thm g G et al}.

We do not know if Proposition \ref{rho hat for chi hat -1} holds with equality when $G$ is not cyclic:

\begin{prob}  \label{prob rho hat for chi hat -1 h}
Decide whether equality holds in Proposition \ref{rho hat for chi hat -1} when $G$ is noncyclic.

\end{prob}

We can generalize Problem \ref{prob rho hat for chi hat -1 h}, as follows.  In Chapter \ref{ChapterCriticalnumber}, we define and investigate the {\em critical numbers} of groups; in particular, in Section \ref{CritRfixed} we discuss the {\em restricted $h$-critical number} $\chi \hat{\;} (G,h)$ of $G$, defined as the smallest positive integer $m$ (if it exists) for which $\rho \hat{\;} (G,m,h)=n$ or, equivalently, the smallest positive integer $m$ for which $h \hat{\;}A=G$ holds for every $m$-subset $A$ of $G$.  For example, according to Proposition \ref{upper bound for restricted 2-critical}, for every group $G$ we have
$$ \chi \hat{\;} (G,2)   = \frac{n+|\mathrm{Ord}(G,2)|+3}{2}.$$
The value of $\chi \hat{\;} (G,h)$ is not known in general---in Section \ref{CritRfixed} we summarize what we know.  

With this notation, we can restate Theorem \ref{thm g G et al} in the form
$$\rho \hat{\;} (\mathbb{Z}_n, \chi \hat{\;} (\mathbb{Z}_n,2)-1 ,2)= \left\{
\begin{array}{ll}
n-1 & \mbox{if} \; \exists k \in \mathbb{N}, n=2^k;\\ \\
n-2 & \mbox{otherwise.}
\end{array}
\right.$$
We can then ask for the following:

\begin{prob}
For each $G$ and $h$, find $\rho \hat{\;} (G, \chi \hat{\;} (G,h)-1 ,h)$.

\end{prob}

As an example, for $G=\mathbb{Z}_{15}$ we have the following values:

$$\begin{array}{||c||c|c||} \hline \hline
h & \chi \hat{\;} (G,h) & \rho \hat{\;} (G, \chi \hat{\;} (G,h)-1 ,h) \\ \hline \hline
1 & 15 & 14 \\ \hline
2 & 9 & 13 \\ \hline
3 & 9 & 14 \\ \hline
4 & 8 & 13 \\ \hline
5 & 9 & 14 \\ \hline
6 & 9 & 13 \\ \hline
7-13 & h+2 & h+1 \\ \hline
14 & 15 & 1 \\ \hline \hline
\end{array}$$

The most general problem, of course, is:

\begin{prob}
Find the exact value of $\rho\hat{\;} (G,m,h)$ for all groups $G$ and positive integers $m$ and $h$.

\end{prob}

We should point out that, while $\rho (G,m,h)$ depends only on the order $n$ of $G$ and not on the structure of $G$, this is definitely not the case for $\rho\hat{\;} (G,m,h)$ (see, for example, Proposition \ref{rhohat m=4} or further results below).

Very few general results are known for $\rho\hat{\;} (G,m,h)$ when $G$ is not cyclic; indeed, most exact values thus far have been for $\rho\hat{\;} (\mathbb{Z}_p^r,m,2)$ when $p$ is prime.  

The case $p=2$ is easy.  Consider any $m$-subset $A=\{a_1,\dots,a_m\}$ of $G$.  We then have
$$2\hat{\;}A=\{a_i+a_j \; | \; 1 \leq i < j \leq m\}.$$  Since
$$2 A=\{a_i+a_j \; | \; 1 \leq i \leq j \leq m\}=2\hat{\;}A \cup \{0\}$$  and $0 \not \in 2\hat{\;}A$, we have the following result.

\begin{prop} \label{rho hat Z_2^r}
For all positive integers $r$ and $m$, we have
$$\rho \hat{\;} (\mathbb{Z}_2^r,m,2)=\rho (\mathbb{Z}_2^r,m,2)-1=u(2^r,m,2)-1.$$

\end{prop}

The corresponding answer for $h \geq 3$ is not known.

\begin{prob}
Find $\rho \hat{\;} (\mathbb{Z}_2^r,m,h)$ for $h \geq 3$.

\end{prob}

For $p \geq 3$ we do not have a complete answer even for $h=2$.  Here is what we do know:  

\begin{thm} [Eliahou and Kervaire; cf.~\cite{EliKer:1998a}] \label{rho hat Z_p^r}\index{Eliahou, S.}\index{Kervaire, M.}

Suppose that $p$ is an odd prime, and $r$ and $m$ are positive integers with $2 \leq m \leq p^r$.

\begin{enumerate}

\item The following statements are equivalent:

\begin{enumerate}
  \item ${2m-2 \choose m-1}$ is divisible by $p$
  \item The base $p$ representation of $m-1$ contains a digit that is at least $\frac{p+1}{2}$
  \item $\rho (\mathbb{Z}_p^r,m,2)=u(p^r,m,2)<2m-1$
  \item $\rho \hat{\;} (\mathbb{Z}_p^r,m,2)=u(p^r,m,2)$
\end{enumerate}

\item The following statements are equivalent:

\begin{enumerate}
  \item ${2m-2 \choose m-1}$ is not divisible by $p$
  \item Each digit in the base $p$ representation of $m-1$ is at most $\frac{p-1}{2}$
  \item $\rho (\mathbb{Z}_p^r,m,2)=u(p^r,m,2)=2m-1$
  \item $\rho \hat{\;} (\mathbb{Z}_p^r,m,2)<u(p^r,m,2)$ and $\rho \hat{\;} (\mathbb{Z}_p^r,m,2) \in \{2m-3,2m-2\}$
\end{enumerate}

\end{enumerate} 

\end{thm}

The fact that statements (a), (b), (c) are equivalent---which, obviously, only need to be proved in 1 or 2---can be established by some elementary considerations (recall that $u(n,m,2) \leq f_1(m,2)=2m-1$ holds for all $n$ and $m$).  Observe that the only unsettled case in Theorem \ref{rho hat Z_p^r} occurs in case 2: we don't know whether the value of $\rho \hat{\;} (\mathbb{Z}_p^r,m,2)$ equals $2m-2$ or $2m-3$.  In \cite{EliKer:2001a} the authors write ``The problem of deciding this alternative is unsolved and looks amazingly difficult in general.''  They have the following additional partial results.

\begin{thm} [Eliahou and Kervaire; cf.~\cite{EliKer:1998a, EliKer:2001a, EliKer:2001b}] \label{rho hat Z_p^r m-1}\index{Eliahou, S.}\index{Kervaire, M.}
Keeping the notations of Theorem \ref{rho hat Z_p^r}, suppose that we are in case 2.

\begin{enumerate}
  \item 
If $m-1$ is not divisible by $p$, then
$\rho \hat{\;} (\mathbb{Z}_p^r,m,2)=2m-3$.

\item 

For $m=p+1$ we have $\rho \hat{\;} (\mathbb{Z}_p^r,m,2)=2m-3$ when $p=3$ and $\rho \hat{\;} (\mathbb{Z}_p^r,m,2)=2m-2$ when $p \geq 5$.

\item 

Let $m=9k+4$ for some nonnegative integer $k$; we then have $\rho \hat{\;} (\mathbb{Z}_3^r,m,2)=2m-3$.

\item 

Let $m=9 \cdot 3^k+1$ for some nonnegative integer $k$; we then have $\rho \hat{\;} (\mathbb{Z}_3^r,m,2)=2m-2$.

\end{enumerate} 
\end{thm}

We should point out that in the second statement of Theorem \ref{rho hat Z_p^r m-1} the case $p=3$ follows immediately from Proposition \ref{rhohat m=4}; this simple fact, in turn, implies the third statement of the theorem.  After Theorems \ref{rho hat Z_p^r} and \ref{rho hat Z_p^r m-1} we see that $\rho \hat{\;} (\mathbb{Z}_3^r,m,2)$ has been settled for all cases except when $m \equiv 1$ mod 27 or when $m \equiv 10$ mod 27.  For $p \geq 5$, the simplest unsolved case seems to be when $m=2p+1$.

\begin{prob}
Find the value of $\rho \hat{\;} (\mathbb{Z}_3^r,m,2)$ for each $r \geq 2$ and $m \equiv 1$ mod $27$ or $m \equiv 10$ mod $27$.

\end{prob}

\begin{prob}
Find the value of $\rho \hat{\;} (\mathbb{Z}_p^r,2p+1,2)$ for each $r \geq 2$ and each prime $p \geq 5$.

\end{prob}  

Furthermore:

\begin{prob}
Find the value of $\rho \hat{\;} (\mathbb{Z}_p^r,m,h)$ for each prime $p \geq 3$, $r \geq 2$, $m \geq 5$, and $h \geq 3$.

\end{prob} 

We should also mention an ``inverse'' result of the same authors:

\begin{thm} [Eliahou and Kervaire; cf.~\cite{EliKer:2001a}] \label{inverse rho hat Z_p^r p+1}\index{Eliahou, S.}\index{Kervaire, M.} 
Suppose that $p \geq 5$ is prime, $A \subset \mathbb{Z}_p^r$ with $|A|=p+1$, and that $|2 \hat{\;} A|=\rho \hat{\;} (\mathbb{Z}_p^r,p+1,2)=2p.$  Then $A$ is of the form $$A=(a_1+H) \cup \{a_2\}$$ for some $a_1 , a_2 \in G$ and $H \leq G$ with $|H|=p$.

\end{thm}
It is easy to verify that the set $A$ above has a restricted 2-fold sumset of size $2p$; according to Theorem \ref{inverse rho hat Z_p^r p+1}, all sets with minimum sumset size are of this form.

Turning now to general abelian groups, we have the following (unsurprising) result for the case of $m=n$:

\begin{thm} \label{rho hat m=n}
For positive integers $r$, we have $2\hat{\;} \mathbb{Z}_2^r=\mathbb{Z}_2^r \setminus \{0\},$ but for all other finite abelian groups $G$ and $1 \leq h \leq \lfloor n/2 \rfloor$, we have $h\hat{\;} G=G.$  Consequently, 
$$\rho\hat{\;} (G,n,h)=
\left \{
\begin{array}{cl}
n-1 & \mbox{if} \; G \cong \mathbb{Z}_2^r \; \mbox{and} \; h=2; \\ \\
n & \mbox{otherwise.}
\end{array}
\right.$$

\end{thm}
The proof of Theorem \ref{rho hat m=n} can be found on page \pageref{proof of rho hat m=n} (for $h=2$, see Proposition \ref{upper bound for restricted 2-critical} as well).

With exact results for $\rho\hat{\;} (G,m,h)$ few and far between, we turn to lower and upper bounds.

We start with some lower bounds on $\rho\hat{\;} (G,m,h)$.  Repeating an argument from earlier, we see that for all $A=\{a_1,\dots,a_m\} \subseteq G$ and $1 \leq h \leq m$, we have $$\{a_1+a_2+\cdots+a_{h-1}+a_i  \mid i=h,h+1,\dots,m\} \subseteq h\hat{\;}A,$$ and these $m-h+1$ elements are distinct, thus $|h\hat{\;}A| \geq m-h+1$.  The dual inequality (switching $h$ and $m-h$) gives $|h\hat{\;}A| \geq h+1$, and therefore we have the following obvious lower bound. 

\begin{prop} \label{rhohat(G,m,h)>=m-h+1}
For all $1 \leq h \leq m-1$ we have $$\rho\hat{\;} (G,m,h) \geq \mathrm{max} \{m-h+1,h+1\};$$ consequently, for  $1 \leq h \leq \left\lfloor \frac{m}{2} \right \rfloor$ we have $$\rho\hat{\;} (G,m,h) \geq m-h+1.$$

\end{prop}

Naturally, we would like to know about the cases where equality occurs in Proposition \ref{rhohat(G,m,h)>=m-h+1}--- this has been very recently answered by Girard, Griffiths, and Hamidoune in \cite{GirGriHam:2012a}.\index{Hamidoune, Y. O.}\index{Girard, B.}\index{Griffiths, S.} For $h=1$, the answer, of course, is obvious: equality holds for all $G$ and $m$.  For $h=2$, we have already seen that equality holds if $m=4$ and $G$ has more than one element of order 2.  More generally, if $A$ is any coset $g+H$ of some subgroup $H=\{0,h_2,\dots,h_m\}$ of $G_2=\mathrm{Ord}(G,2)\cup \{0\}$, then, since $h_i+h_j \in H \setminus \{0\}$ for any $2 \leq i < j \leq m$, we have 
$$2\hat{\;}A=\{2g+h_2, \dots, 2g+h_m\},$$ and thus $|2\hat{\;}A|=m-1$.  
It turns out that these are the only cases when equality holds in Proposition \ref{rhohat(G,m,h)>=m-h+1}:

\begin{thm} [Girard, Griffiths, and Hamidoune; cf.~\cite{GirGriHam:2012a}] \label{rhohat(G,m,h)=m-h+1}\index{Hamidoune, Y. O.}\index{Girard, B.}\index{Griffiths, S.} 

Suppose that $A$ is an $m$-subset of $G$ and that $1 \leq h \leq \left\lfloor \frac{m}{2} \right \rfloor$.   Then $$|h\hat{\;}A|=m-h+1$$ if, and only if, $h=1$ (and $A$ arbitrary) or $h=2$ and $A$ is a coset of some subgroup of $\mathrm{Ord}(G,2)\cup \{0\}$.

\end{thm}   

As a consequence of Theorem \ref{rhohat(G,m,h)=m-h+1}, we get that $$\rho\hat{\;} (G,m,h) \geq m-1,$$ with equality if, and only if, $h \in \{2,m-2\}$ and $\mathrm{Ord}(G,2)\cup \{0\}$ possesses a subgroup of order $m$ (and thus $m=2^k$ for some $2 \leq k \leq e$ where $e$ is the number of even orders in the invariant factorization of $G$).  

As it turns out, we also have a full characterization of cases where $$\rho\hat{\;} (G,m,h) = m.$$  This equality clearly holds for $h=1$, $h=m-1$, or when $m$ is a divisor of $n$ (in which case $G$ has a subgroup of order $m$).  

Let us demonstrate another example.
Let $H=\{0,h_2,\dots,h_m,h_{m+1}\}$ be a subgroup of $G_2=\mathrm{Ord}(G,2)\cup \{0\}$ in a (noncyclic) group $G$, and let $A$ be a coset of $H$ in $G$ with one element removed:
$$A=\{g,g+h_2,\dots,g+h_m\}.$$  As above, we see that  $$2\hat{\;}A=\{2g+h_2, \dots, 2g+h_m\},$$ and thus $|2\hat{\;}A|=m$.     
The question of lower bounds for $\rho\hat{\;} (G,m,h)$ is settled for general $G$ by the following results.

\begin{thm} [Girard, Griffiths, and Hamidoune; cf.~\cite{GirGriHam:2012a}] \label{rhohat(G,m,h)>=m}\index{Hamidoune, Y. O.}\index{Girard, B.}\index{Griffiths, S.} 

Suppose that $m \geq 5$ and $1 \leq h \leq \left\lfloor \frac{m}{2} \right \rfloor$, and let $A$ be an $m$-subset of $G$ that is not a coset of some subgroup of $\mathrm{Ord}(G,2)\cup \{0\}$.   Then $|h\hat{\;}A| \geq m$.  Furthermore, $|h\hat{\;}A| = m$ if, and only if, (at least) one of the following holds:
\begin{enumerate}[(i)]
  \item $h=1$ (and $A$ arbitrary), 
  \item $A$ is a coset of some subgroup of $G$, or
  \item $h=2$  and $A$ is a coset of some subgroup of $\mathrm{Ord}(G,2)\cup \{0\}$ minus one element.
\end{enumerate}
\end{thm}  
Note that the assumption that $m \geq 5$ is permissible as we have the corresponding results for $m \leq 4$ in  Proposition \ref{rhohat m=4}  above.

It is worthwhile to compare Theorem \ref{rhohat(G,m,h)>=m} to Proposition \ref{|hA|=m}.  We state the following corollary explicitly.

\begin{cor}  \label{cor to GGH}
Suppose that $m \geq 5$ and $1 \leq h \leq m-1$.   Let $e$ be the number of even orders in the invariant factorization of $G$.  We then have 
$$\rho\hat{\;} (G,m,h) \left\{\begin{array}{ll}
=m-1  & \mbox{if $h \in \{2,m-2\}$ and $m=2^k$ for some $2 \leq k \leq e$;}\\ \\
=m & \mbox{if $h \in \{1,m-1\}$, or} \\
& \mbox{$m|n$, or} \\
& \mbox{$h \in \{2,m-2\}$ and $m=2^k-1$ for some $2 \leq k \leq e$;}\\ \\
\geq m+1 & \mbox{otherwise}.
\end{array}\right.$$

\end{cor}

We can attempt to find a different type of lower bound for $\rho\hat{\;} (G,m,h)$ by considering Proposition \ref{uhat versus p}.  We make the following conjecture:

\begin{conj} \label{rhohat versus p}
Let $G$ be a group of order $n$, $p$ be the smallest prime divisor of $n$, and assume that $1 \leq h< m \leq n$.  We then have 
$$\rho\hat{\;} (G,m,h) \geq \min \{p, hm-h^2+1\},$$
with equality if, and only if, $p \geq m$.
\end{conj} 

Note that Conjecture \ref{rhohat versus p} is a generalization of Theorem \ref{Dias Da Silva and Hamidoune}.

K\'arolyi  succeeded in providing a proof for the case $h=2$:

\begin{thm} [K\'arolyi; cf.~\cite{Kar:2003a, Kar:2004a}] \label{Karolyi1}\index{K\'arolyi, Gy.}
For any abelian group $G$ of order $n$, we have
$$\rho\hat{\;} (G,m,2) \geq \min \{p, 2m-3\},$$
where $p$ is the smallest prime divisor of $n$.  Equality may occur if, and only if, $p \geq m$.
\end{thm}

A very difficult problem is the following.

\begin{prob}
Prove Conjecture \ref{rhohat versus p} for $h \geq 3$.
\end{prob}

Let us attempt to characterize situations where the lower bound of Conjecture \ref{rhohat versus p} is achieved and 
$$\rho\hat{\;} (G,m,h) = \min \{p, hm-h^2+1\},$$ where $p$ is the smallest prime divisor of $n$.  We examine four cases.

First, let us assume that $p < m$, in which case $\min \{p, hm-h^2+1\}=p$ for all $h$.  If $p$ is odd, then $n$ is odd, and thus, by Theorem \ref{rhohat(G,m,h)>=m}, we have $|h\hat{\;}A| \geq m$, so we cannot have $|h\hat{\;}A| =p$.  If $p=2<m$, then $|h\hat{\;}A| =p$ cannot happen either, since by Proposition \ref{rhohat(G,m,h)>=m-h+1} this would imply that 
$$2 \geq \mathrm{max} \{m-h+1,h+1\},$$ which cannot hold for any $m$ and $h$ with $2 < m$.  So this case yields no examples, and we, in fact, verified the ``only if'' part of the last claim in Conjecture \ref{rhohat versus p}; we are about to prove the ``if'' part as well.

Suppose now that $m \leq p <  hm-h^2+1$.  Let $H$ be any subgroup of $G$ with $|H|=p$, and let $A$ be any $m$-subset of $H$.  By Theorem \ref{Dias Da Silva and Hamidoune}, $$|h\hat{\;}A| = \mathrm{min} \{p, hm-h^2+1\}=p.$$  Therefore, any $m$-subset of $H$ and, more generally, any $m$-subset of any coset of $H$, is an example for a set for which the lower bound of Conjecture \ref{rhohat versus p} is achieved.  We then question whether there are any other such sets:

\begin{conj} \label{conj |h hat A|=p}
Let $p$ be the smallest prime divisor of $n$, $A$ be an $m$-subset of $G$, and assume that $m \leq p <  hm-h^2+1$.  Then $|h\hat{\;}A|=p$ if, and only if, $A$ is contained in a coset of some subgroup $H$ of $G$ with $|H|=p$.   
\end{conj}

\begin{prob}
Prove (or disprove) Conjecture \ref{conj |h hat A|=p}.
\end{prob}

Assume now that $p>hm-h^2+1$.  This time, we can find $m$-subsets $A$ of $G$ for which $$|h\hat{\;}A|=hm-h^2+1=\min \{p, hm-h^2+1\},$$ as follows.  First of all, any $m$-subset $A$ of $G$ will do, if $h=1$ or $h=m-1$.  Second, by our proposition on page \pageref{h=2 m=4 size}, we see that when $h=2$ and $m=4$, 4-subsets of $G$ of the form $$A=\{a,a+d_1,a+d_2,a+d_1+d_2\}$$ work for any $a, d_1, d_2 \in G$.  (Note that, since $p >hm-h^2+1=5$ in this case, $n$ is odd, and thus $G$ has no elements of order 2.  Therefore, $|2 \hat{\;}A|=5.$)

We can find further examples by using arithmetic progressions.  Let $g \in G$ be of order $p$.  If $A$ is the arithmetic progression $$A=\{a,a+g,a+2g, \dots, a+(m-1)g\},$$ then, as we have seen, $$h\hat{\;}A=\left\{ha+\frac{h^2-h}{2}g,\dots, ha+\left(h(m-1)-\frac{h^2-h}{2} \right)g \right\};$$ since $p>hm-h^2+1$, these $hm-h^2+1$ elements are all distinct.  We believe that there are no other examples.

\begin{conj} \label{conj |hA|=hm-h^2+1}
Let $p$ be the smallest prime divisor of $n$, $A$ be an $m$-subset of $G$, and assume that $p >  hm-h^2+1$.  Then $|h\hat{\;}A|=hm-h^2+1$ if, and only if, (at least) one of the following holds:
\begin{enumerate}[(i)]
\item $h=1$  (and $A$ arbitrary);
\item $h=m-1$  (and $A$ arbitrary);
\item $h=2$, $m=4$, and $A$ is of the form $A=\{a,a+d_1,a+d_2,a+d_1+d_2\}$ (for arbitrary $a, d_1, d_2 \in G$);
\item $A$ is an arithmetic progression (of size $m$).
\end{enumerate}   
\end{conj}

By a result of K\'arolyi, we have:

\begin{thm} [K\'arolyi; cf.~\cite{Kar:2005a}] \label{Karolyi2}\index{K\'arolyi, Gy.}
Conjecture \ref{conj |hA|=hm-h^2+1} holds for $h =2$.
\end{thm}

A difficult problem is the following.

\begin{prob}
Prove (or disprove) Conjecture \ref{conj |hA|=hm-h^2+1} for $h \geq 3$.
\end{prob}

Finally, suppose that $p=hm-h^2+1$.  This case (mirroring the analogous situation for unrestricted sumsets) is more complicated.  For example, with $p=7$, $m=5$, and $h=2$ (as is the case, for example, with $\rho \hat{\;} (\mathbb{Z}_{49},5,2)=7$) we see that the set
$$A=\{0,a,2a,(n+a)/2, (n+3a)/2 \}$$ provides an example, as we have
$$2 \hat{\;} A=\{a, 2a, 3a, (n+a)/2, (n+3a)/2, (n+5a)/2, (n+7a)/2\}.$$
We have the following intriguing problem:

\begin{prob}
Suppose that $n$, $m$, and $h$ are positive integers, $p$ is the smallest prime divisor of $n$, and $p=hm-h^2+1$.  Classify all $m$-subsets $A$ of $G$ for which $\rho \hat{\;} (G,m,h)=p$.

\end{prob}

Of course, these last type of lower bounds on $\rho\hat{\;} (G,m,h)$ are not meaningful when $p$ (the smallest prime divisor of $|G|$) is small.  For example, when a subset of $G$ is within a subgroup of small size (or within a coset of that subgroup), then its restricted sumset---indeed, any sumset---will not be larger than that subgroup (a coset of that subgroup).   Therefore, when there is a divisor $d$ of $n$ with $d \geq m$, then  $\rho \hat{\;} (G,m,h)$ cannot be more than $d$.  Thus, to form more meaningful lower bounds, one may wonder what happens if we assume that no such divisor $d$ exists.  
The following theorem provides one answer to this question.

\begin{thm} [Hamidoune, Llad\'o, and Serra; cf.~\cite{HamLlaSer:2000a}] \label{Hamidoune Llado, Serra}\index{Hamidoune, Y. O.}\index{Llad\'o, A. S.}\index{Serra, O.}  
Suppose that $G$ is an abelian group that is either cyclic or is of odd order, and let $m$ be an integer with $m \geq 33$ and $m \geq 21$ in these two cases, respectively.  Furthermore, suppose that the only divisor of $n$ that is greater than or equal to $m$ is $n$ itself.  Then
$$\rho\hat{\;} (G,m,2) \geq \min\{n, 3m/2\}.$$
\end{thm}

We can prove that Conjecture \ref{conj rhohatforh=2}, with the additional assumption that the only divisor of $n$ that is greater than or equal to $m$ is $n$ itself, implies Theorem \ref{Hamidoune Llado, Serra} when $G$ is cyclic, and we only need to know that $m \geq 8$.  (The claim is false for $m=7$ as example $C_1$ on page \pageref{example c_1} demonstrates.)  Indeed, we will show that
$$\min\{u(n,m,2),2m-4\} \geq \min\{n, 3m/2\}.$$
Note that $m \geq 8$ implies that $2m-4 \geq 3m/2$; therefore, our claim clearly holds when $$u(n,m,2) \geq \min\{n, 3m/2\}.$$ We can show that, in fact, 
if the only divisor of $n$ that is greater than or equal to $m$ is $n$ itself, then $$u(n,m,h) \geq \min\{n, (h+1)m/2\}$$ holds for all positive integers $h$.

To see this, let $d \in D(n)$ be such that $$u(n,m,h)=f_d(m,h)=\left( h \left \lceil \frac{m}{d} \right \rceil-h+1 \right)d.$$
Note that when $d< m \leq 2d$, we have
$$\left( h \left \lceil \frac{m}{d} \right \rceil-h+1 \right)d =(h+1)d \geq (h+1)m/2,$$ and in the case when $m > 2d$, we have
$$\left( h \left \lceil \frac{m}{d} \right \rceil-h+1 \right)d  \geq \left( h \frac{m}{d}-h+1 \right)d=hm-(h-1)d > hm - (h-1)m/2= (h+1)m/2.$$  
The remaining case to consider is when $m \leq d$, but then, by assumption, $d=n$, for which we have $f_n=n$, and thus our claim is established.

To see that Theorem \ref{Hamidoune Llado, Serra} is not true for an arbitrary finite abelian group, consider $G=\mathbb{Z}_2^r \times \mathbb{Z}_3$, and let $A=\mathbb{Z}_2^r \times \{0,1\}$.  Note that $n=3 \cdot 2^r$, $m=2 \cdot 2^r$, so the only divisor of $n$ that is greater than or equal to $m$ is $n$ itself.  We can see that $$2 \hat{\;} A=G \setminus \{(0,0),(0,2)\},$$ thus $$|2 \hat{\;} A|=n-2=3m/2-2.$$  
This example prompts us to wonder the following:

\begin{prob}
Is there a constant $C$ (independent of $n$ and $m$) for which for every $n$ and $m$ with the property that the only divisor of $n$ that is greater than or equal to $m$ is $n$ itself, the inequality
$$\rho\hat{\;} (G,m,2) \geq \min\{n, 3m/2\}-C$$
holds for every abelian group $G$ of order $n$? 
\end{prob}
As our example above demonstrates, if there is such a constant, it must be at least 2.

A similar argument shows that, for $h=3$, Conjecture \ref{conj rhohatforh=3} implies that, if $m \geq 8$ and the only divisor of $n$ that is greater than or equal to $m$ is $n$ itself, then
$$\rho\hat{\;} (\mathbb{Z}_n,m,3) \geq \min\{n, 2m-2\}.$$
Indeed, according to Conjecture \ref{conj rhohatforh=3}, we have 
$$\rho\hat{\;} (\mathbb{Z}_n,m,3) \geq \min\{u(n,m,3), 3m-10,3m-3-\gcd(n,m-1)\}.$$  The fact that our assumptions imply that $$u(n,m,3) \geq \min\{n, 2m\} \geq \min\{n, 2m-2\}$$ was shown above, and for $m \geq 8$ we obviously have $3m-10 \geq 2m-2$ and $3m-3-\gcd(n,m-1) \geq 2m-2$.

Therefore, we pose the following:

\begin{prob}
Is there a constant $C$ (independent of $n$ and $m$) for which for every $n$ and $m$ with the property that the only divisor of $n$ that is greater than or equal to $m$ is $n$ itself, the inequality
$$\rho\hat{\;} (G,m,3) \geq \min\{n, 2m\}-C$$
holds for every abelian group $G$ of order $n$? 
\end{prob}
Since $\rho\hat{\;}(\mathbb{Z}_{15},6,3)=u \hat{\;}(15,6,3)=8$, $C$ would have to be at least 4.

More generally:

\begin{prob}
Let $h$ be a given positive integer.  Is there a constant $C(h)$ (independent of $n$ and $m$ but dependent on $h$) for which for every $n$ and $m$ with the property that the only divisor of $n$ that is greater than or equal to $m$ is $n$ itself, the inequality
$$\rho\hat{\;} (G,m,h) \geq \min\{n, (h+1)m/2\}-C(h)$$
holds for every abelian group $G$ of order $n$? 
\end{prob}

We can propose some other lower bounds for $\rho\hat{\;} (G,m,h)$ (at least for $h=2$) in terms of $\rho (G,m,h)$.  One such conjecture is the following:

\begin{conj} [Plagne; cf.~\cite{Pla:2006b}] \label{conj rhohat lower}\index{Plagne, A.} 
For any finite abelian group $G$ and positive integer $m \geq 2$  we have $$\rho \hat{\;} (G,m,2) \geq \rho (G,m,2)-2.$$
\end{conj}
If true, Conjecture \ref{conj rhohat lower} limits the value of $\rho \hat{\;} (G,m,2)$ to one of only three possibilities: $\rho (G,m,2)$, $\rho (G,m,2)-1$, or $\rho (G,m,2)-2$.  We can verify that each exact result mentioned in this section satisfies this conjecture.

\begin{prob}

Prove or disprove Conjecture \ref{conj rhohat lower}.
\end{prob}

We should point out that Conjecture \ref{conj rhohatforh=2} implies Conjecture \ref{conj rhohat lower} for cyclic groups.  Indeed, we have $$\rho (\mathbb{Z}_n,m,2) \leq f_1=2m-1;$$ furthermore, when $n$ is even, we have $$\rho (\mathbb{Z}_n,m,2) \leq f_2=2m-2,$$ and when $n$ is divisible by $2m-2$, we find that $$\rho (\mathbb{Z}_n,m,2) \leq f_{2m-2}=2m-2.$$

Another conjecture of this type is the following:

\begin{conj} [Lev; cf.~\cite{Lev:2000a}] \label{conj rhohat lower with ord}\index{Lev, V. F.} 
For any $G$ and  $m$, we have $$\rho \hat{\;} (G,m,2) \geq \min\{\rho (G,m,2), 2m-3-|\mathrm{Ord}(G,2)|\}.$$
\end{conj}

\begin{prob}

Prove or disprove Conjecture \ref{conj rhohat lower with ord}.
\end{prob}

Lev\index{Lev, V. F.} proved the somewhat weaker result that 
$$\rho \hat{\;} (G,m,2) \geq \min\{\rho (G,m,h), \theta m-3-|\mathrm{Ord}(G,2)|\}$$ where $\theta=\frac{1+\sqrt{5}}{2} \approx 1.6$ is the golden ratio.  For the case when $G$ is cyclic, $\mathrm{Ord}(G,2)$ can have at most one element (namely, $\frac{n}{2}$ if $n$ is even), and thus Conjecture \ref{conj rhohatforh=2} implies Conjecture \ref{conj rhohat lower with ord} for cyclic groups.

We can further illuminate Conjectures \ref{conj rhohat lower} and \ref{conj rhohat lower with ord} by the following example for a set $A$ for which we have
both $$|2\hat{\;}A|=|2A|-2$$ and $$|2\hat{\;}A|=2|A|-3-|\mathrm{Ord}(G,2)|.$$  Let $G$ be any group of even order that has at least one element of order  more than 2 (and thus $G$ is not isomorphic to $\mathbb{Z}_2^r$).  

Recall that, by the Fundamental Theorem of Finite Abelian Groups, we have
$$G \cong \mathbb{Z}_2^{\alpha_1} \times \mathbb{Z}_4^{\alpha_2} \times \cdots \times \mathbb{Z}_{2^k}^{\alpha_k} \times H,$$
where $k \in \mathbb{N}$ (since $n$ is even), $\alpha_i \in \mathbb{N}_0$ for each $i=1,2,\dots,k$ with at least one being positive, and $H$ is a group of odd order (perhaps $|H|=1$).   One can also easily verify that
$$G_2=\mathrm{Ord}(G,2) \cup \{0\}=\{g \in G \; | \; 2g=0\}$$ is a subgroup of $G$.  Let $a$ be an element of $G$ of order at least 3, and set $$A=G_2 \cup \{a+g_2 \; | \; g_2 \in G_2\}.$$  We can then check that $$2A=G_2 \cup \{a+g_2 \; | \; g_2 \in G_2\} \cup \{2a+g_2 \; | \; g_2 \in G_2\}$$ and $$2\hat{\;}A=2A \setminus \{0,2a\}.$$  (Note that the sum of two distinct elements of order 2 in $G$ is again an element of order 2, and the sum of 0 and an element of order 2 is, of course, an element of order 2.)  Therefore, we have $$|2\hat{\;}A|=|2A|-2=2|A|-3-|\mathrm{Ord}(G,2)|.$$

By Proposition \ref{rho-hat very small}, \label{further than C} for every $h \geq 3$ and for every positive real number $C$, one can find positive integers $n$ and $m$ for which $$\rho \hat{\;} (\mathbb{Z}_n,m,h) \leq   u\hat{\;}(n,m,h)< u(n,m,h)- C=\rho (\mathbb{Z}_n,m,h)-C.$$  Therefore, for $h \geq 3$, one cannot expect a claim similar to Conjecture \ref{conj rhohat lower}.  However, we offer:

\begin{prob}

Generalize Conjecture \ref{conj rhohat lower with ord} for $h \geq 3$.
\end{prob}

Let us now turn to upper bounds for $\rho \hat{\;} (G,m,h)$.  Of course, most trivially, we have
$$\rho \hat{\;} (G,m,h) \leq \rho (G,m,h).$$  More meaningful bounds were given by Plagne for $h=2$:\index{Plagne, A.} 

\begin{thm} [Plagne; cf.~\cite{Pla:2006b}] \label{thm rhohat upper}\index{Plagne, A.} 
For any $G$ and  $m$, we have $$\rho \hat{\;} (G,m,2) \leq \min\{\rho (G,m,2), 2m-2\}.$$  Furthermore, if $n$ has a prime divisor $p$ that does not divide $m-1$, then $$\rho \hat{\;} (G,m,2) \leq \min\{\rho (G,m,2), 2m-3\}.$$
\end{thm}

We are not aware of similar bounds for $h \geq 3$ and thus we pose:

\begin{prob}

Find a function $f(h)$ (as ``small'' as possible) for which $$\rho \hat{\;} (G,m,h) \leq \min\{\rho (G,m,h), hm-h^2+1+f(h)\}$$ holds for any $G$, $m$, and  $h$.

\end{prob}
Recall that, by Proposition \ref{uhat(n,m,h)extreme} we have $$u \hat{\;} (n,m,h) \leq \min\{u(n,m,h), hm-h^2+1\}$$ for any $n$, $m$, and $h $.

A particularly intriguing question is the following:

\begin{prob}
Find a set analogous to $A_d(n,m)$ (described earlier) that yields a good upper bound for $\rho \hat{\;} (G,m,h)$ when $G$ has rank 2 (or more).
\end{prob}

Before closing this section, we define two new quantities related to restricted addition.  For positive integers $m$ and $h$, we let

$$\rho \hat{\;} (m,h)_{\min}=\min \{ \rho \hat{\;}(G,m,h) \; \mid \; |G| \geq m\},$$ and
$$\rho \hat{\;} (m,h)_{\max}=\max \{ \rho \hat{\;}(G,m,h) \; \mid \; |G| \geq m\}.$$
(The condition that the order of $G$ is at least $m$ is necessary to avoid trivialities.)

Clearly, $$\rho \hat{\;} (m,1)_{\min}=\rho \hat{\;} (m,1)_{\max}=m,$$ and, as before, we may assume $2 \leq h \leq \lfloor m/2 \rfloor$.  By these considerations, the first non-trivial example is $(m,h)=(4,2)$: According to Proposition \ref{rhohat m=4}, we have $\rho \hat{\;} (4,2)_{\min}=3$ and $\rho \hat{\;} (4,2)_{\max}=5$.  Thus we may assume that $m \geq 5$.

The value of $\rho \hat{\;} (m,h)_{\min}$ is already known for all $(m,h)$: By Corollary \ref{cor to GGH}, we have:

\begin{thm}
Suppose that $m$ and $h$ are positive integers with $m \geq 5$ and $h \leq \lfloor m/2 \rfloor$.  Then
$$\rho \hat{\;} (m,h)_{\min}= \left\{\begin{array}{cl}
m-1  & \mbox{if $h=2$ and $m=2^k$ for some $k \in \mathbb{N}$;}\\ \\
m & \mbox{otherwise}.
\end{array}\right.$$

\end{thm}

The value of $\rho \hat{\;} (m,h)_{\max}$ is not known in general but, from previous results, we can find lower and upper bounds.  For a lower bound, recall that by Theorem \ref{Dias Da Silva and Hamidoune}, when $p$ is a prime with $$p \geq hm-h^2+1,$$ then $$\rho \hat{\;} (\mathbb{Z}_p,m,h)=hm-h^2+1.$$  This yields the lower bound $hm-h^2+1$. 

For an upper bound, recall that by Theorem \ref{rho (G,m,h)}, we have
$$\rho (G,m,h)=u(n,m,h),$$ and if $$n \geq hm-h+1,$$ then by Proposition \ref{u(n,m,h)extreme}, $$u(n,m,h) \leq hm-h+1,$$   which then is an upper bound for $\rho \hat{\;} (m,h)_{\max}$.  In summary:

\begin{thm}
Suppose that $m$ and $h$ are positive integers with $m \geq 5$ and $h \leq \lfloor m/2 \rfloor$.  Then
$$hm-h^2+1 \leq \rho \hat{\;} (m,h)_{\max} \leq hm-h+1.$$

\end{thm}

For $h=2$, we know a bit more: by Theorem \ref{thm rhohat upper}, we see that $$\rho \hat{\;} (m,2)_{\max} \leq 2m-2.$$  Therefore:

\begin{thm}
For all positive integers $m$, 
$\rho \hat{\;} (m,2)_{\max}$ is either $2m-3$ or $2m-2.$

\end{thm}

We offer the following intriguing problem:

\begin{prob}
Find the value of $\rho \hat{\;} (m,2)_{\max}$ for all positive integers $m$.

\end{prob}

We note that, by Theorem \ref{rho hat Z_p^r m-1}, we have $$\rho \hat{\;} (m,2)_{\max} = 2m-2$$ when $m-1 \geq 5$ is a prime, and also when $m=9 \cdot 3^k+1$ for some $k \in \mathbb{N}_0$.  As a consequence of Theorem \ref{carrick}, for $m=5$ we have:

\begin{thm} [Carrick; cf.~\cite{Car:2015a}]\index{Carrick, E.}
We have $\rho \hat{\;} (5,2)_{\max}=7$.

\end{thm} 

We know very little about $\rho \hat{\;} (m,h)_{\max}$ for $h \geq 3$:

\begin{prob}
Find the value of $\rho \hat{\;} (m,h)_{\max}$ for all positive integers $m$ and $3 \leq h \leq \lfloor m/2 \rfloor$.

\end{prob}

\subsection{Limited number of terms} \label{2minRlimited}

Here we consider, for a given group $G$, positive integer $m$ (with $m \leq n=|G|$), and nonnegative integer $s$, $$\rho \hat{\;} (G,m,[0,s]) = \mathrm{min} \{ |[0,s] \hat{\;} A|  \mid A \subseteq G, |A|=m\},$$ that is, the minimum size of $\cup_{h=0}^s h \hat{\;} A$ for an $m$-element subset $A$ of $G$.  Clearly, we may restrict our attention to $0 \leq s \leq m$; in fact, when $s \geq m$, $$\rho \hat{\;} (G,m,[0,s]) =\rho \hat{\;} (G,m,\mathbb{N}_0)$$ (cf.~Subsection \ref{2minRarbitrary}).

Note that, although we always have   
$$\rho (G,m,[0,s]) = \rho  (G,m,s)$$
(cf.~Section \ref{2minUlimited}), the quantities $\rho \hat{\;} (G,m,[0,s])$ and $\rho \hat{\;} (G,m,s)$ are not necessarily equal.  (Of course, the inequality
$$\rho \hat{\;} (G,m,[0,s]) \geq \rho \hat{\;} (G,m,s)$$ 
holds.)

\begin{prob} 
Find $\rho \hat{\;} (G,m,[0,s])$ for any group $G$, positive integer $m \leq n$, and nonnegative integer $s$.
\end{prob}

Some authors have investigated the minimum size of $[0,s] \hat{\;} A$ among all $m$-subsets of $G$ that possess some additional properties.  One such pursuit is the attempt to find
$$\rho_A \hat{\;} (G,m,[0,s]) = \mathrm{min} \{ |[0,s] \hat{\;} A|  \mid A \subseteq G, |A|=m, A \cap -A =\emptyset \},$$ that is, the minimum size of $[0,s] \hat{\;} A$ among {\em asymmetric} $m$-subsets of $G$.
Since no element of $\mathrm{Ord}(G,2) \cup \{0\}$ may be in an asymmetric set in $G$, for an asymmetric set of size $m$ to exist, we need to assume that $$m \leq \frac{n-|\mathrm{Ord}(G,2)|-1}{2}.$$

We can get a lower bound in the case of the cyclic group $\mathbb{Z}_n$ by considering the set 
$$A=\{1,2,\dots,m\};$$ this results in 
$$[0,s] \hat{\;} A=\{0,1,2,\dots, ms-(s^2-s)/2\}$$ and thus we get:

\begin{prop} \label{asym lower [0,s] hat}
For positive integers $n$, $m$, and $s$ with $m \leq \lfloor (n-1)/2 \rfloor$ and $s \leq m$, we have $$\rho_A \hat{\;} (\mathbb{Z}_n,m,[0,s]) \leq \min \{n, ms-(s^2-s)/2+1 \}.$$ 

\end{prop}

As was recently proved by Balandraud\index{Balandraud, \'E.} in \cite{Bal:2017a}, equality occurs in Proposition \ref{asym lower [0,s] hat} when $n$ is prime:

\begin{thm} [Balandraud; cf.~\cite{Bal:2017a}]  \label{Balandraud asym [0,s]}\index{Balandraud, \'E.}
If $p$ is prime and $m \leq \lfloor (p-1)/2 \rfloor$, we have $$\rho_A \hat{\;} (\mathbb{Z}_p,m,[0,s]) = \min \{p, ms-(s^2-s)/2+1 \}.$$ 
\end{thm}

The value of $\rho_A \hat{\;} (G,m,[0,s])$ is not known for groups with composite order.

\begin{prob}
Find $\rho_A \hat{\;} (\mathbb{Z}_n,m,[0,s])$ for composite values of $n$.
\end{prob} 

\begin{prob}
Find $\rho_A \hat{\;} (G,m,[0,s])$ for noncyclic groups $G$.
\end{prob}

\subsection{Arbitrary number of terms} \label{2minRarbitrary}

In this subsection we discuss what we know about 
$$\rho \hat{\;} (G,m,\mathbb{N}_0) = \mathrm{min} \{ |\Sigma A|  \mid A \subseteq G, |A|=m\},$$ that is, the minimum size of $$\Sigma A=\cup_{h=0}^m h \hat{\;} A$$ among all $m$-subsets $A$ of $G$.  

Suppose that $G$ is cyclic.  We can develop an upper bound for $\rho \hat{\;} (\mathbb{Z}_n,m,\mathbb{N}_0)$, as follows.  Let $d$ be an arbitrary positive divisor of $n$; we write $m$ as 
$$m=cd+k$$ where $1 \leq k \leq d$ and thus $$c=\lceil m/d \rceil-1.$$  We also let $H$ be the subgroup of $\mathbb{Z}_n$ of size $d$:
$$H=\{j \cdot n/d \mid j=0,1,\dots,d-1\}.$$  We separate two cases depending on the parity of $c$.  

When $c$ is odd, we let 
$$A=\bigcup_{i=-(c-1)/2}^{(c-1)/2} (i+H) \bigcup \left\{(c+1)/2+j \cdot n/d \mid j=0,1,\dots,k-1\right\}.$$  Then $|A|=m$.  We also see that
$$\Sigma A=\bigcup_{i=-\frac{c^2-1}{8} \cdot d}^{\frac{c^2-1}{8} \cdot d+\frac{c+1}{2} \cdot k} (i+H),$$
and thus
\begin{eqnarray*}
|\Sigma A|&=&\min \left \{ n, \; \left( \frac{c^2-1}{4} \cdot d+\frac{c+1}{2} \cdot k+1 \right) \cdot d \right \} \\ \\
&=&\min \left \{ n, \; \left( \frac{c^2-1}{4} \cdot d+\frac{c+1}{2} \cdot (m-cd)+1 \right) \cdot d \right \} \\ \\
&=&\min \left \{ n, \; \left( \frac{c+1}{2} \cdot m - \frac{(c+1)^2}{4} \cdot d+1 \right) \cdot d \right \} .
\end{eqnarray*}

The case when $c$ is even is similar; this time we let 
$$A=\bigcup_{i=-c/2}^{c/2-1} (i+H) \bigcup \left\{c/2+j \cdot n/d \mid j=0,1,\dots,k-1\right\}.$$  Again $|A|=m$, and we have
$$\Sigma A=\bigcup_{i=-\frac{c^2+2c}{8} \cdot d}^{\frac{c^2-2c}{8} \cdot d+\frac{c}{2} \cdot k} (i+H),$$
and thus
\begin{eqnarray*}
|\Sigma A|&=&\min \left \{ n, \; \left( \frac{c^2}{4} \cdot d+\frac{c}{2} \cdot k+1 \right) \cdot d \right \} \\ \\
&=&\min \left \{ n, \; \left( \frac{c^2}{4} \cdot d+\frac{c}{2} \cdot (m-cd)+1 \right) \cdot d \right \} \\ \\
&=&\min \left \{ n, \; \left( \frac{c}{2} \cdot m - \frac{c^2}{4} \cdot d+1 \right) \cdot d \right \}.
\end{eqnarray*}

We can combine our two cases and say that
\begin{eqnarray*}
|\Sigma A|&=& \min \left \{ n, \; \left( \left \lceil \frac{c}{2} \right \rceil  \cdot m - \left \lceil   \frac{c}{2}  \right \rceil ^2 \cdot d+1 \right) \cdot d \right \} \\ \\
&=& \min \left \{ n, \; \left( \left \lceil \frac{\lceil m/d -1 \rceil }{2} \right \rceil  \cdot m - \left \lceil   \frac{\lceil m/d -1 \rceil}{2}  \right \rceil ^2 \cdot d+1 \right) \cdot d \right \} \\ \\
&=& \min \left \{ n, \; \left( \left \lceil \frac{ m/d -1  }{2} \right \rceil  \cdot m - \left \lceil   \frac{ m/d -1 }{2}  \right \rceil ^2 \cdot d+1 \right) \cdot d \right \}.
\end{eqnarray*}
Observe also that for $d=n$, we have $\left \lceil (m/n -1)/2 \right \rceil =0$, so
$$\left( \left \lceil \frac{ m/n -1  }{2} \right \rceil  \cdot m - \left \lceil   \frac{ m/n -1 }{2}  \right \rceil ^2 \cdot n+1 \right) \cdot n=n.$$ 
This provides the following upper bound:

\begin{prop} \label{upper for |Sigma|}
With $D(n)$ denoting the set of positive divisors of $n$, we have
$$\rho \hat{\;} (\mathbb{Z}_n,m,\mathbb{N}_0) \leq \min \left \{ \left( \left \lceil \frac{ m/d -1  }{2} \right \rceil  \cdot m - \left \lceil   \frac{ m/d -1 }{2}  \right \rceil ^2 \cdot d+1 \right) \cdot d \mid d \in D(n) \right \}.$$

\end{prop}

We risk the following:

\begin{conj}  \label{conj upper for |Sigma|}
For all positive integers $n$ and $m$ with $m \leq n$, we have equality in Proposition \ref{upper for |Sigma|}.

\end{conj}

\begin{prob} 

Prove or disprove Conjecture \ref{conj upper for |Sigma|}.

\end{prob}

Evaluating our expression in Proposition \ref{upper for |Sigma|} for $d=1$ yields
$$\left \lceil (m -1  )/2 \right \rceil  \cdot m - \left \lceil   (m -1  )/2  \right \rceil ^2 +1 = \lfloor m^2/4 \rfloor +1,$$ 
so we get the following: 

\begin{cor} \label{cor upper for |Sigma|}
For every positive integer $n$ we have
$$\rho \hat{\;} (\mathbb{Z}_n,m,\mathbb{N}_0) \leq \min \left \{ n, \lfloor m^2/4 \rfloor + 1 \right \}.$$

\end{cor}

We can, in fact, prove that, for cyclic groups of prime order $p$, equality holds in Corollary \ref{cor upper for |Sigma|}, since clearly
$$\rho \hat{\;} (\mathbb{Z}_p,m,\mathbb{N}_0) \geq \rho \hat{\;} (\mathbb{Z}_p,m,\lfloor m/2 \rfloor),$$ and, by Theorem \ref{Dias Da Silva and Hamidoune}, we have 
$$\rho \hat{\;} (\mathbb{Z}_p,m,\lfloor m/2 \rfloor)= \min \left \{ p, \lfloor m/2 \rfloor \cdot m- \lfloor m/2 \rfloor^2 + 1 \right \}= \min \left \{ p, \lfloor m^2/4 \rfloor + 1 \right \}.$$

\begin{prop}  \label{prop for |Sigma| prime}
For all positive primes $p$ and positive integers $m \leq p$, $$\rho \hat{\;} (\mathbb{Z}_p,m,\mathbb{N}_0) = \min \left \{ p, \lfloor m^2/4 \rfloor + 1 \right \}.$$

\end{prop}

We do not know the value of $\rho \hat{\;} (G,m,\mathbb{N}_0)$ for noncyclic groups:

\begin{prob}

Find the value of $\rho \hat{\;} (G,m,\mathbb{N}_0)$ for all positive integers $m \leq n$ and noncyclic groups $G$. 

\end{prob}

Some authors have investigated the minimum size of $\Sigma A$ among all $m$-subsets of $G$ that possess some additional properties.  One such pursuit is the attempt to find
$$\rho_A \hat{\;} (G,m,\mathbb{N}_0) = \mathrm{min} \{ |\Sigma A|  \mid A \subseteq G, |A|=m, A \cap -A =\emptyset \},$$ that is, the minimum size of $\Sigma A$ among {\em asymmetric} $m$-subsets of $G$.
Since no element of $\mathrm{Ord}(G,2) \cup \{0\}$ may be in an asymmetric set in $G$, for an asymmetric set of size $m$ to exist, we need to assume that $$m \leq \frac{n-|\mathrm{Ord}(G,2)|-1}{2}.$$

We can get a lower bound in the case of the cyclic group $\mathbb{Z}_n$ by considering the set 
$$A=\{1,2,\dots,m\};$$ this results in 
$$\Sigma A=\{0,1,2,\dots, m(m+1)/2\}$$ and thus we get:

\begin{prop} \label{asym lower sigma}
For positive integers $n$ and $m$ with $m \leq \lfloor (n-1)/2 \rfloor$, we have $$\rho_A \hat{\;} (\mathbb{Z}_n,m,\mathbb{N}_0) \leq \min \{n, (m^2+m+2)/2 \}.$$ 

\end{prop}

As a consequence of Balandraud's Theorem \ref{Balandraud asym [0,s]}\index{Balandraud, \'E.} with $s=m$ (and as proved by him in \cite{Bal:2012a}), equality occurs in Proposition \ref{asym lower sigma} when $n$ is prime:

\begin{thm} [Balandraud; cf.~\cite{Bal:2012a}, \cite{Bal:2012b}]  \label{Balandraud asym}\index{Balandraud, \'E.}
If $p$ is prime and $m \leq \lfloor (p-1)/2 \rfloor$, we have $$\rho_A \hat{\;} (\mathbb{Z}_p,m,\mathbb{N}_0) = \min \{p, (m^2+m+2)/2 \}.$$ 
\end{thm}

The value of $\rho_A \hat{\;} (G,m,\mathbb{N}_0)$ is not known for groups with composite order.

\begin{prob}
Find $\rho_A \hat{\;} (\mathbb{Z}_n,m,\mathbb{N}_0)$ for composite values of $n$.
\end{prob} 

\begin{prob}
Find $\rho_A \hat{\;} (G,m,\mathbb{N}_0)$ for noncyclic groups $G$.
\end{prob} 

As a variation, we may consider
$$\rho_A \hat{\;} (G,m,\mathbb{N}) = \mathrm{min} \{ |\Sigma^*  A|  \mid A \subseteq G, |A|=m, A \cap -A =\emptyset \},$$ that is, the minimum size of $$\Sigma^*   A =\cup_{h=1}^{\infty} h \hat{\;} A$$ among asymmetric $m$-subsets of $G$.

The subset $$A=\{1,2,\dots,m\}$$ of $\mathbb{Z}_n$ provides a lower bound again:
\begin{prop} \label{asym lower sigma star}
For positive integers $n$ and $m$ with $m \leq \lfloor (n-1)/2 \rfloor$, we have $$\rho_A \hat{\;} (\mathbb{Z}_n,m,\mathbb{N}_0) \leq \min \{n, (m^2+m)/2 \}.$$ 

\end{prop}

And we have:

\begin{thm} [Balandraud; cf.~\cite{Bal:2012a}, \cite{Bal:2012b}] \label{Balandraud}\index{Balandraud, \'E.}
If $p$ is prime and $m \leq \lfloor (p-1)/2 \rfloor$, we have $$\rho_A \hat{\;} (\mathbb{Z}_p,m,\mathbb{N}) = \min \{p, (m^2+m)/2 \}.$$ 
\end{thm}

\begin{prob}
Find $\rho_A \hat{\;} (\mathbb{Z}_n,m,\mathbb{N})$ for composite values of $n$.
\end{prob} 

\begin{prob}
Find $\rho_A \hat{\;} (G,m,\mathbb{N})$ for noncyclic groups $G$.
\end{prob}

Another related quantity of interest is
$$\widehat{\rho_*} \hat{\;} (G,m,\mathbb{N}_0) = \mathrm{min} \{ |\Sigma A|  \mid A \subseteq G \setminus \{0\}, |A|=m, \langle A \rangle =G \},$$ that is, the minimum size of $\Sigma A$ among those $m$-subsets $A$ of $G$ that do not contain zero (hence the $_*$) but generate the entire group (hence the $\widehat{\;}$).

A relevant result is the following:

\begin{thm} [Hamidoune; cf.~\cite{Ham:1998a}] \label{Ham min sumset size}\index{Hamidoune, Y. O.}
Let $S \subseteq G \setminus \{0\}$ be such that $|S| \geq 3$ and $\langle S \rangle=G$; furthermore, if $|S|=3$, then suppose that $S$ is not of the form $\{\pm a, 2a\}$ for any $a \in G$, and that if $|S|=4$, then $S$ is not of the form $\{\pm a, \pm 2a\}$ for any $a \in G$.  Then
$$|\Sigma S| \geq \min \{n-1, 2|S|\}.$$

\end{thm}

Combining Theorem \ref{Ham min sumset size} with the fact that if $n \geq 10$, then $\chi \hat{\;} (G^*,\mathbb{N}) \leq n/2$ (see Section \ref{CritRarbitrary}), we get the following:

\begin{cor} \label{Hamidoune got it wrong}
Suppose that $n \geq 10$, and let $S \subseteq G \setminus \{0\}$ be such that $|S| \geq 5$ and $\langle S \rangle=G$.  Then
$$|\Sigma S| \geq \min \{n, 2|S|\}.$$

\end{cor}
Indeed, if $|S| \leq (n-1)/2$, then the result follows from Theorem \ref{Ham min sumset size}; if $|S| \geq n/2$, then it is due to the fact that $\chi \hat{\;} (G^*,\mathbb{N}) \leq n/2$.  (We mention that Corollary \ref{Hamidoune got it wrong} was stated in \cite{Ham:1998a} only for cyclic groups $G$; its proof there, however, relied on Lemma 3.3 in \cite{Ham:1998a}, which is actually false for all groups of even order; cf.~\cite{Baj:2016a}.)

Furthermore, as it was noted by Hamidoune\index{Hamidoune, Y. O.} in \cite{Ham:1998a}, if $n$ is divisible by 3, and $$S=(H \setminus \{0\}) \cup \{s\}$$ for an index 3 subgroup $H$ of $G$ and $s \in G \setminus H$, then $|\Sigma S|=2n/3$.  Consequently,

\begin{cor}
If $n$ is divisible by 3 and $n \geq 15$, then $$\widehat{\rho_*} \hat{\;} (G,n/3,\mathbb{N}_0)  = 2n/3.$$

\end{cor}

\begin{prob}
Evaluate $\widehat{\rho_*} \hat{\;} (G,m,\mathbb{N}_0)$ for other $G$ and $m$.
\end{prob}

Yet another variation, introduced by Eggleton and Erd\H{o}s in \cite{EggErd:1972a},\index{Eggleton, R. B.} is to find the following quantity:
$$\rho_Z \hat{\;} (G,m,\mathbb{N}) = \mathrm{min} \{ |\Sigma^*  A|  \mid A \subseteq G, |A|=m, 0 \not \in \Sigma^*  A \},$$ that is, the minimum size of $\Sigma^*  A$ among all weakly zero-sum-free $m$-subsets $A$ of $G$; if no such set $A$ exists, we set $\rho_Z \hat{\;} (G,m,\mathbb{N})=\infty$.

For $m=1$ the answer is obvious, since $|\Sigma^*  A|=1$ for all 1-subsets of $G$.  The case of $m=2$ is not much harder: if $A=\{a,b\}$ with $a \neq b$, $a \neq 0$, and $b \neq 0$, then $\Sigma^*  A=\{a,b,a+b\}$ has size 3; we just need to make sure that $a+b \neq 0$.  We can summarize:

\begin{prop}  \label{m=1,2 zero-free min size}
We have: 
$$\rho_Z \hat{\;} (G,1,\mathbb{N})=\left\{
\begin{array}{cl}
\infty & \mbox{if $n =1$}, \\
1 & \mbox{if $n \geq 2$};
\end{array}
\right.$$
and
$$\rho_Z \hat{\;} (G,2,\mathbb{N})=\left\{
\begin{array}{cl}
\infty & \mbox{if $n \leq 3$}, \\
3 & \mbox{if $n \geq 4$}.
\end{array}
\right.$$

\end{prop}

For $m=3$, we have the following result:

\begin{prop}  \label{m=3 zero-free min size}
We have: 
$$\rho_Z \hat{\;} (G,3,\mathbb{N})=\left\{
\begin{array}{cl}
\infty & \mbox{if $n \leq 5$}; \\
5 & \mbox{if $n \geq 6$, $n$ is even, and $G \not \cong \mathbb{Z}_2^r$}; \\
6 & \mbox{if $n \geq 7$ and $n$ is odd}; \\
7 & \mbox{if $n \geq 8$ and $G \cong \mathbb{Z}_2^r$}.
\end{array}
\right.$$

\end{prop}
The proof of Proposition \ref{m=3 zero-free min size} is on page \pageref{proof of m=3 zero-free min size}.

\begin{prob}
Evaluate $\rho_Z \hat{\;} (G,4,\mathbb{N})$, $\rho_Z \hat{\;} (G,5,\mathbb{N})$, etc.~for all abelian groups $G$.
\end{prob}

A considerably easier, but still largely unknown, special case is the evaluation of $\rho_Z \hat{\;} (G,m,\mathbb{N})$ for cyclic groups $G$; we review what is currently known.  

From Section \ref{5maxRarbitrary}, we recall that the set $A=\{1,2,\dots,m\}$ is weakly zero-sum-free  in $\mathbb{Z}_n$ when $n \geq (m^2+m+2)/2$; since $$\Sigma^*  A=\{1,2,\dots,(m^2+m)/2\},$$ this implies:

\begin{prop} \label{rho Z easy}
If $n \geq (m^2+m+2)/2$, then $$\rho_Z \hat{\;} (\mathbb{Z}_n,m,\mathbb{N}) \leq (m^2+m)/2.$$ 

\end{prop}

The two constructions of Selfridge mentioned on page \pageref{Selfridge constructions} yield slightly better bounds for certain values of $n$.  The first, based on the fact that $$A=\{1,-2,3,4,\dots,m\}$$ is weakly zero-sum-free in $\mathbb{Z}_n$ when $n \geq (m^2+m-2)/2$ and $n \geq 6$, has an advantage over Proposition \ref{rho Z easy} only when $n=(m^2+m-2)/2$ or $n=(m^2+m)/2$, since we have
$$\Sigma^*  A=\{-1,-2,1,2,\dots,(m^2+m-4)/2\},$$  which has the same size as the set yielding that theorem when $n \geq (m^2+m+2)/2$.  However, for $n=(m^2+m-2)/2$, $$-2=(m^2+m-6)/2$$ and $$-1=(m^2+m-4)/2;$$ and for $n=(m^2+m)/2$, $$-2=(m^2+m-4)/2.$$  We thus get that if $A$ is the $m$-subset of the cyclic group of order $(m^2+m-2)/2$ described above, then $A$ is weakly zero-sum-free with $$|\Sigma^*  A|=(m^2+m-4)/2;$$ and if $A$ is the $m$-subset of the cyclic group of order $(m^2+m)/2 $ described above, then $A$ is weakly zero-sum-free with $$|\Sigma^*  A|=(m^2+m-2)/2.$$  We can take this further by noting that whenever $A$ is weakly zero-sum-free in $\mathbb{Z}_{n}$, then $d \cdot A$ is weakly zero-sum-free in $\mathbb{Z}_{dn}$; this holds for every $d \geq 1$ and neither the size of the subset nor the size of its sumset changes.  We thus arrive at the following result:

\begin{prop} \label{rho Z less easy}
If $n  \geq 6$ and $n$ is divisible by $(m^2+m-2)/2 $, then $$\rho_Z \hat{\;} (\mathbb{Z}_n,m,\mathbb{N}) \leq (m^2+m-4)/2,$$ and if $n \geq 6$ and $n$ is divisible by  $(m^2+m)/2$, then $$\rho_Z \hat{\;} (\mathbb{Z}_n,m,\mathbb{N}) \leq (m^2+m-2)/2.$$

\end{prop}

The second construction assumes that $n$ is even, and considers the set
$$A=\{1,2,\dots, \left \lfloor (m-1)/2 \right \rfloor\} \cup \{n/2, n/2+1,n/2+2,\dots, n/2+\left \lfloor m/2 \right \rfloor\},$$
for which
$$\Sigma^*  A=\{1,2,\dots, \left \lfloor m^2/4 \right \rfloor\} \cup \{n/2, n/2+1,n/2+2,\dots, n/2+\left \lfloor m^2/4 \right \rfloor\}.$$
Therefore, if $n/2 \geq \left \lfloor m^2/4 \right \rfloor+1$, that is, if $n \geq \left \lfloor m^2/2 \right \rfloor+2$, then $A$ is weakly zero-sum-free.  We get:

\begin{prop}  \label{rho Z less less easy}
If $n$ is even and $n  \geq \left \lfloor m^2/2 \right \rfloor+2$, then 
$$\rho_Z \hat{\;} (\mathbb{Z}_n,m,\mathbb{N}) \leq \left \lfloor m^2/2 \right \rfloor+1.$$
\end{prop}

Let us see now what our propositions tell us for small values of $m$.  The cases of $m \leq 3$ follow from Propositions \ref{m=1,2 zero-free min size} and \ref{m=3 zero-free min size}.  For $m=4$,
\begin{enumerate}[(i)]
\item Proposition \ref{rho Z easy} yields that if $n \geq 11$, then $\rho_Z \hat{\;} (\mathbb{Z}_n,m,\mathbb{N}) \leq 10$;
\item Proposition \ref{rho Z less easy} yields that if $9|n$, then $\rho_Z \hat{\;} (\mathbb{Z}_n,m,\mathbb{N}) \leq 8$;
\item Proposition \ref{rho Z less easy} also yields that if $10|n$, then $\rho_Z \hat{\;} (\mathbb{Z}_n,m,\mathbb{N}) \leq 9$; and
\item Proposition \ref{rho Z less less easy} yields that if $n \geq 10$ and $2|n$, then $\rho_Z \hat{\;} (\mathbb{Z}_n,m,\mathbb{N}) \leq 9$.
\end{enumerate}
Note that statement (iii) follows from (iv) and thus is unnecessary.  In \cite{BhoHalSch:2011a}, Bhowmik, Halupczok, and Schlage-Puchta\index{Bhowmik, G.}\index{Halupczok, I.}\index{Schlage-Puchta, J-C.} presented one additional construction: if $n \geq 12$ and $3|n$, then the set $$A=\{1, n/3, n/3+1, 2n/3+1\}$$ is weakly zero-sum-free in $\mathbb{Z}_n$, and $\Sigma^*  A$ has size 9.  Furthermore, relying on a computer program, they proved that one can never do better:

\begin{thm} [Bhowmik, Halupczok, and Schlage-Puchta; cf.~\cite{BhoHalSch:2011a}]\index{Bhowmik, G.}\index{Halupczok, I.}\index{Schlage-Puchta, J-C.}
We have: 
$$\rho_Z \hat{\;} (\mathbb{Z}_n,4,\mathbb{N})=\left\{
\begin{array}{cl}
\infty & \mbox{if $n \leq 8$}; \\
8 & \mbox{if $9|n$}; \\
9 & \mbox{if $n \geq 10$ and $9 \not | n$ but ($2|n$ or $3|n$)}; \\
10 & \mbox{otherwise}.
\end{array}
\right.$$

\end{thm}

For $m=5$ and $m=6$, the same authors proved that our Propositions \ref{rho Z easy}, \ref{rho Z less easy}, and \ref{rho Z less less easy} provide the right values:

\begin{thm} [Bhowmik, Halupczok, and Schlage-Puchta; cf.~\cite{BhoHalSch:2011a}]\index{Bhowmik, G.}\index{Halupczok, I.}\index{Schlage-Puchta, J-C.}
We have: 
$$\rho_Z \hat{\;} (\mathbb{Z}_n,5,\mathbb{N})=\left\{
\begin{array}{cl}
\infty & \mbox{if $n \leq 13$}; \\
13 & \mbox{if $n \geq 14$ and $2|n$}; \\
14 & \mbox{if $15|n$}; \\
15 & \mbox{otherwise};
\end{array}
\right.$$
and 
$$\rho_Z \hat{\;} (\mathbb{Z}_n,6,\mathbb{N})=\left\{
\begin{array}{cl}
\infty & \mbox{if $n \leq 19$}; \\
19 & \mbox{if $n \geq 20$ and $2|n$}; \\
20 & \mbox{if $21|n$}; \\
21 & \mbox{otherwise}.
\end{array}
\right.$$

\end{thm}
The same authors also exhibited the precise (but quite a bit more complicated) formula for $\rho_Z \hat{\;} (\mathbb{Z}_n,7,\mathbb{N})$ (see \cite{BhoHalSch:2011a}).

\begin{prob}
Evaluate $\rho_Z \hat{\;} (\mathbb{Z}_n,8,\mathbb{N})$, $\rho_Z \hat{\;} (\mathbb{Z}_n,9,\mathbb{N})$, etc.~for all values of $n$.

\end{prob} 

There is more known about $\rho_Z \hat{\;} (\mathbb{Z}_n,m,\mathbb{N})$ for prime values of $n$.  First, recall that, according to Balandraud's result\index{Balandraud, \'E.} from \cite{Bal:2012a} (cf.~Theorem \ref{Balandraud coroll}), if $$1+2+ \cdots +m \geq p,$$ then $\mathbb{Z}_p$ has no weakly zero-sum-free subsets of size $m$.  We can restate this as follows:

\begin{thm} [Balandraud; cf.~\cite{Bal:2012a}, \cite{Bal:2012b}] \label{Balandraud coroll rho hat Z}\index{Balandraud, \'E.}
If $p$ is prime for which $p \leq (m^2+m)/2$, then $\rho_Z \hat{\;} (\mathbb{Z}_p,m,\mathbb{N})=\infty$.
\end{thm}

On the other hand, from Proposition \ref{rho Z easy} we get that if $p \geq (m^2+m+2)/2$, then $$\rho_Z \hat{\;} (\mathbb{Z}_p,m,\mathbb{N}) \leq (m^2+m)/2.$$
Olson  proved that when $p$ is large, equality holds:

\begin{thm} [Olson; cf.~Theorem 2 in ~\cite{Ols:1968a}] \label{Olson rho hat Z p}\index{Olson, J. E.} 
Let $p$ be a prime for which $p \geq m^2+m-1$ when $m$ is even, and 
$p \geq m^2+(3m-5)/2$ when $m$ is odd.  
Then 
$$\rho_Z \hat{\;} (\mathbb{Z}_p,m,\mathbb{N}) = (m^2+m)/2.$$
\end{thm}

This raises the following:

\begin{prob} \label{prob p rho hat Z}
Decide whether $$\rho_Z \hat{\;} (\mathbb{Z}_p,m,\mathbb{N}) = (m^2+m)/2$$ holds for all primes $p$ with $p \geq (m^2+m+2)/2.$

\end{prob}
According to Theorem \ref{Olson rho hat Z p}, only finitely many primes need to be considered to answer Problem \ref{prob p rho hat Z} for a given value of $m$.  We should observe that the answer to Problem \ref{prob p rho hat Z} is affirmative for $m=1$ and $m=2$ (by Proposition \ref{m=1,2 zero-free min size}), for $m=3$ (by Proposition \ref{m=3 zero-free min size}), and for $m \in \{4,5,6,7\}$ (by \cite{BhoHalSch:2011a}).

Rather than finding $\rho_Z \hat{\;} (G,m,\mathbb{N})$ for all $G$ and $m$, Eggleton and Erd\H{o}s in \cite{EggErd:1972a}\index{Eggleton, R. B.}\index{Erd\H{o}s, P.}  proposed the potentially easier problem of evaluating $f(m)$, defined as the minimum possible value of $\rho_Z \hat{\;} (G,m,\mathbb{N})$ for any group $G$.  According to Propositions \ref{m=1,2 zero-free min size} and \ref{m=3 zero-free min size}, we have $f(1)=1$, $f(2)=3$, and $f(3)=5$.  We also have:
\begin{itemize}
\item $f(4)=8$ (Eggleton and Erd\H{o}s; cf.~\cite{EggErd:1972a});\index{Eggleton, R. B.}\index{Erd\H{o}s, P.}
\item $f(5)=13$ (Gao, et al.; cf.~\cite{GaoEtAl:2008a});\index{Gao, W.}\index{Li, Y.}\index{Peng, J.}\index{Sun, F.}
\item $f(6)=19$ (Gao, et al.; cf.~\cite{GaoEtAl:2008a});\index{Gao, W.}\index{Li, Y.}\index{Peng, J.}\index{Sun, F.}
\item $f(7)=24$ (Yuan and Zeng; cf.~\cite{YuaZen:2010a}).\index{Yuan, P.}\index{Zeng, X.} 
\end{itemize}
The proof of $f(6)=19$ is long with a very large number of cases, and the proof of $f(7)=24$ relies on a computer program.

\begin{prob}
Evaluate (perhaps relying on a computer program) $f(8)$, $f(9)$, etc.
\end{prob}

We can also observe that for $m \leq 7$, $f(m)$ agrees with the smallest possible value of $\rho_Z \hat{\;} (G,m,\mathbb{N})$ for cyclic groups $G$, and 
Eggleton and Erd\H{o}s believed that this is always the case:\index{Erd\H{o}s, P.}

\begin{conj} [Eggleton and Erd\H{o}s; cf.~\cite{EggErd:1972a}] \label{Erdos Eggleton conj}\index{Eggleton, R. B.}\index{Erd\H{o}s, P.}
For every $m \in \mathbb{N}$ there is an $n \in \mathbb{N}$ for which $f(m)=\rho_Z \hat{\;} (\mathbb{Z}_n,m,\mathbb{N}).$

\end{conj}

\begin{prob}
Prove Conjecture \ref{Erdos Eggleton conj}.

\end{prob}

We should mention that by Proposition \ref{rho Z less less easy}, we get:

\begin{cor} \label{cor f(m)}
For every $m \in \mathbb{N}$, we have $f(m) \leq \left \lfloor m^2/2 \right \rfloor+1.$

\end{cor}
From our stated values above, we have equality in Corollary \ref{cor f(m)} for $m \in \{1,2,3,5,6\}$ but not for $m \in \{4,7\}$. 

A particularly intriguing problem is the following:

\begin{prob}
Find infinitely many values of $m$ for which $f(m) \leq \left \lfloor m^2/2 \right \rfloor$ or prove that this is not possible.

\end{prob}

We also have a lower bound for $f(m)$:

\begin{thm} [Olson; cf.~\cite{Ols:1975a}]\index{Olson, J. E.} 
For every $m \in \mathbb{N}$, we have $f(m) \geq \left \lceil m^2/9 \right \rceil.$

\end{thm}

\begin{prob}
Find a real number $c > 1/9$ so that $f(m) \geq \left \lceil c \cdot m^2 \right \rceil$ holds for all (but finitely many ) $m \in \mathbb{N}$.
\end{prob}

\section{Restricted signed sumsets} \label{2minRS}

\subsection{Fixed number of terms} \label{2minRSfixed}

\subsection{Limited number of terms} \label{2minRSlimited}

\subsection{Arbitrary number of terms} \label{2minRSarbitrary}

\chapter{The critical number} \label{ChapterCriticalnumber}

Recall that in Chapter \ref{ChapterMinsumsetsize} we investigated, for given $\Lambda$ and $H$, the minimum sumset size of an $m$-subset of $G$:
$$\rho_{\Lambda}(G,m,H)=\mathrm{min} \{ |H_{\Lambda}A|  \mid A \subseteq G, |A|=m\}.$$
As a special case, here we are interested in the minimum value of $m$ for which $$\rho_{\Lambda}(G,m,H)=n;$$ that is, the minimum value of $m$ for which every $m$-subset of $G$ spans all of $G$.
This value, if exists, is called the {\em $(\Lambda,H)$-critical number} of $G$ and is denoted by $\chi_{\Lambda}(G,H)$.

In the following sections we consider $\chi_{\Lambda}(G,H)$ for special $\Lambda \subseteq \mathbb{Z}$ and $H \subseteq \mathbb{N}_0$.

\section{Unrestricted sumsets} \label{critU}

Our goal in this section is to investigate $\chi (G,H)$, the minimum value of $m$ for which $$HA=G$$ holds for every $m$-subset of $G$.  (Recall that $HA$ is the union of all $h$-fold sumsets $hA$ for $h \in H$.)  Since $0A=\{0\}$ for every subset $A$ of $G$ but $hG=G$ for every positive integer $h$, we see that $\chi (G,H)$ does not exist for $n \geq 2$ when $H=\{0\}$, but $\chi (G,H)$ does exist and is at most $n$ when $H$ contains at least one positive integer.  

 We consider three special cases: when $H$ consists of a single nonnegative integer $h$, when $H$ consists of all nonnegative integers up to some value $s$, and when $H$ is the entire set of nonnegative integers. 

\subsection{Fixed number of terms} \label{CritUfixed}

Here we ought to consider, for fixed $G$ and positive integer $h$, the quantity $\chi (G,h)$, that is, the minimum value of $m$ for which the $h$-fold sumset of every $m$-element subset of $G$ is $G$ itself.  However, according to Theorem  \ref{rho (G,m,h)}, the $h$-critical number of a group of order $n$ is the minimum value of $m$ for which $u(n,m,h)=n$, and this value was determined in Section \ref{0.3.5.1} by Theorem \ref{h-critical of n} to be $v_1(n,h)+1$ where
$$v_1(n,h)= \max \left\{ \left( \left \lfloor \frac{d-2}{h} \right \rfloor +1  \right) \cdot \frac{n}{d}   \mid d \in D(n) \right\}.$$
Thus we have:

\begin{thm}  \label{h crit numb}
For all finite abelian groups $G$ of order $n$ and all positive integers $h$ we have $$\chi (G,h)=v_1(n,h)+1.$$
\end{thm}

Having found the value of $\chi (G,h)$, we are now interested in the inverse problem of classifying all $m$-subsets $A$ of $G$ with $$m=\chi (G,h)-1=v_1(n,h)$$ for which $hA \neq G$.  The problem being trivial for $h=1$, we let $h \geq 2$.  

We consider $h=2$ first.  (As we explain below, the case of $h=2$ seems more complicated than the case of $h \geq 3$.)  

When $n$ is even, we have $v_1(n,2)=n/2$.  Recall that, by Theorem \ref{rho (G,m,h)}, we have
$$\rho(G,m,h)=u(n,m,h)=\min \{f_d(m,h) \mid d \in D(n) \},$$ where $D(n)$ is the set of positive divisors of $n$ and $$f_d(m,h)=\left( h \cdot \lceil m/d \rceil -h+1 \right) \cdot d.$$  In particular, when $f_d(n/2,2) < n$ for some $d \in D(n)$, then we are guaranteed to find subsets $A_d$ of $G$ with $|A_d|=n/2$ and $|2A_d|<n$.

One can easily determine that
  $$\left( 2 \cdot \lceil n/(2d) \rceil -1 \right) \cdot d = \left \{
\begin{array}{cl}
n-d & \mbox{if $n/d$ is even}, \\ \\
n & \mbox{if $n/d$ is odd}.
\end{array}
\right.$$
When $n$ and $n/d$ are both even for some $d \in D(n)$, we can, in fact, find explicit subsets $A_d$ of $G$ of size $n/2$ whose two-fold sumset has size $n-d$; we will explain this here for the case when $G$ is cyclic.  Recall the set $A_d(n,m)$ from page \pageref{Ad(n,m) defined}.  In particular, for $m=n/2$ and when $n/d$ is even, we have
$$A_d(n,n/2)=\cup_{i=0}^{n/(2d)-1} (i+H),$$ where $H$ is the subgroup of $\mathbb{Z}_n$ with order $d$.
We see that $|A_d(n,n/2)|=n/2$, and $$2A_d(n,n/2)=\cup_{i=0}^{n/d-2} (i+H),$$ so $|2A_d(n,n/2)|=n-d.$
We thus have:

\begin{prop}
Suppose $n$ is even and that it has a divisor $d$ for which $n/d$ is also even.  Then $\mathbb{Z}_n$ has a subset $A$ of size $n/2$ for which $2A$ has size $n-d$.  
\end{prop}     

For example, $\mathbb{Z}_{20}$ has subsets $A_d(20,10)$ for $d \in \{1,2,5,10\}$, each of size ten, so that $2A_d(20,10)$ has size $20-d$.  Using the computer program \cite{Ili:2017a}, we checked that for $d \in \{2,5,10\}$, there are essentially (ignoring equivalences) no other 10-subsets $A$ with $2A \neq \mathbb{Z}_{20}$ besides $A_d(20,10)$.  However, there are many 10-subsets whose 2-fold sumset has size 19 other than the (arithmetic progression) $A_1(20,10)$ constructed above: for example,
$$\{0,1,2,3,4,5,6,7\} \cup C$$ with   
$C=\{8,10\}$, $C=\{9,11\}$, $C=\{16,17\}$, etc.   

We pose the following questions, in increasing order of difficulty:

\begin{prob}
For each even value of $n$, classify all subsets of $\mathbb{Z}_n$ of size $n/2$ whose two-fold sumset is not $\mathbb{Z}_n$.

\end{prob}

\begin{prob}
For each abelian group $G$ of even order $n$, classify all subsets of size $n/2$ whose two-fold sumset is not $G$.

\end{prob}

\begin{prob}
For each abelian group $G$ of odd order $n$, classify all subsets of size $(n-1)/2$ whose two-fold sumset is not $G$.

\end{prob}

Let us now turn to the case of $h \geq 3$; we again assume that $n$ is even, in which case $v_1(n,h)=n/2$ by Corollary \ref{v function bounds}.  For divisors $d \in D(n)$, we can compute $f_d(n/2,h)$, as follows.  We see that $f_n(n/2,h)=n$, $f_{n/2}(n/2,h)=n/2$, and (in the case when $n$ is divisible by 3) $f_{n/3}(n/2,h)=(h+1) \cdot n/3 >n$.  For any other $d$ we get
$$f_d=\left( h \cdot \lceil n/(2d) \rceil -h+1 \right) \cdot d \geq \left( h \cdot  n/(2d) -h+1 \right) \cdot d  \geq h \cdot n/2 - (h-1) \cdot n/4 \geq n.$$ So $f_d <n$ holds only for $d=n/2$, in which case $f_d=n/2$.   
This suggests that, if $A$ is a subset of $G$ of size $n/2$ for which $hA \neq G$, then $A=H$ or $A=G \setminus H$ where $H \leq G$ has order $n/2$ (and in both cases $|hA|=n/2$).

\begin{conj} \label{conj inverse h>=3, incomplete}
Suppose that $G$ is an abelian group of even order $n$, $h$ is an integer with $h \geq 3$, and $A$ is a subset of $G$ with $|A|=n/2$ and $|hA| < n$.  Then $G$ has a subgroup $H$ of order $n/2$ for which $A=H$ or $A=G \setminus H$.
\end{conj}

We know that this claim holds for cyclic groups:

\begin{thm} [Navarro; cf.~\cite{Nav:2016a}]\index{Navarro, A.} 
Conjecture \ref{conj inverse h>=3, incomplete} holds when $G$ is cyclic.

\end{thm}

We have the following open questions:

\begin{prob}
Prove Conjecture \ref{conj inverse h>=3, incomplete} for noncyclic groups $G$.

\end{prob}

\begin{prob}
For each $n$ and $h$ with $n$ odd and $h \geq 3$, characterize all subsets $A$ of $\mathbb{Z}_n$ of size $v_1(n,h)$ for which $hA \neq \mathbb{Z}_n$. 

\end{prob} 

\begin{prob}  \label{odd sets w size v1 but not h-span}
For each $G$ of odd order $n$ and for each $h \geq 3$, characterize all subsets $A$ of $G$ of size $v_1(n,h)$ for which $hA \neq G$. 

\end{prob} 

We have the following result of Lev\index{Lev, V. F.} that not only answers Problem \ref{odd sets w size v1 but not h-span} in a special case, but accomplishes more:

\begin{thm} [Lev; cf.~\cite{Lev:2016a}]\index{Lev, V. F.} 
Suppose that $A$ is an $m$-subset of $\mathbb{Z}_5^r$ with $$(3 \cdot 5^{r-1}-1)/2 \leq m \leq 2 \cdot 5^{r-1} = v_1(5^r,3)$$ for which $3A \neq \mathbb{Z}_5^r$.  Then $A$ is contained in a union of two cosets of a subgroup of index 5.

\end{thm}

Note that if $A$ is contained in a union of two cosets of a subgroup of index 5, then indeed $3A \neq \mathbb{Z}_5^r$; furthermore, Lev\index{Lev, V. F.} constructed an example in \cite{Lev:2016a} that shows that the lower bound on $m$ is tight.

We should mention that, as an analogue of $\widehat{\chi}  (G,[0,s])$ discussed in Section \ref{CritUlimited} below, one may define 
the variation where only generating subsets of $G$ are considered:
$$\widehat{\chi}  (G,h) = \min \{m \mid A \subseteq G, \langle A \rangle =G, |A| \geq m \Rightarrow hA=G \}.$$ 
We have the following---somewhat surprising---result:

\begin{thm} [Bajnok; cf.~\cite{Baj:2016b}]\index{Bajnok, B.}
For all $G$ and $h$, we have
$$\widehat{\chi}  (G,h) = \chi (G,h)=v_1(n,h)+1.$$
\end{thm}

It is likely an interesting question to find inverse results as well:

\begin{prob}

For each $G$ and $h$, classify all generating subsets $A$ of $G$ for which $|A|=\widehat{\chi}  (G,h)-1$ and $hA \neq G$.
\end{prob}

\subsection{Limited number of terms} \label{CritUlimited}

For given groups $G$ and positive integer $s$, here we consider $\chi (G,[0,s])$, that is, the minimum value of $m$ for which the $[0,s]$-fold sumset of every $m$-subset of $G$ is $G$ itself.  By  Proposition \ref{rho  (G,m,[0,s])} and Theorem \ref{h crit numb}, we have
\begin{thm}  \label{[0,s] crit numb}
For all finite abelian groups $G$ of order $n$ and all positive integers $s$ we have $$\chi (G,[0,s])=v_1(n,s)+1.$$
\end{thm}

While $$\chi (G,[0,s])= \min \{m \mid A \subseteq G,  |A| \geq m \Rightarrow [0,s]A=G \}$$ has thus been evaluated for each $G$ and $s$, there is a variation that has been considered in the literature for which much less is known.  Namely, Klopsch and Lev\index{Klopsch, B.}\index{Lev, V. F.} in \cite{KloLev:2009a} have investigated the quantity $ \widehat{\chi}  (G,[0,s])$: the minimum value of $m$ for which the $[0,s]$-fold sumset of every $m$-subset of $G$ that generates $G$ is $G$ itself, that is,
$$\widehat{\chi}  (G,[0,s]) = \min \{m \mid A \subseteq G, \langle A \rangle =G, |A| \geq m \Rightarrow [0,s]A=G \}.$$

It turns out that the following relative of the arithmetic function 
$v_1(n,s)$ of page \pageref{0.3.5.3 page} will be useful: we define
$$\widehat{v}(n,s)  = \max \left\{ \left( \left \lfloor \frac{d-2}{s} \right \rfloor +1  \right) \cdot \frac{n}{d}   \mid d \in D(n), d \geq s+2 \right\};$$ we adhere to the convention that the maximum element of the empty-set equals zero, and thus if $n \leq s+1$, we have $\widehat{v}(n,s)=0$.  (Note also that we omitted the index 1 as unnecessary here.)   

From the paper \cite{KloLev:2009a} of Klopsch and Lev\index{Klopsch, B.}\index{Lev, V. F.} we are able to deduce the following upper bound for $\widehat{\chi}  (G,[0,s])$:

\begin{thm}  \label{From Klopsch Lev}
For every $G$ and $s$ we have $$\widehat{\chi}  (G,[0,s]) \leq \widehat{v}(n,s)+1.$$

\end{thm}

Since this result was not stated in \cite{KloLev:2009a}, we provide a proof---see page \pageref{proof of From Klopsch Lev}.

We can easily see that equality holds in Theorem \ref{From Klopsch Lev}, when $G$ is cyclic.  This is obvious when $n \leq s+1$, since then $$\widehat{\chi}  (\mathbb{Z}_n,[0,s])=\widehat{v}(n,s)+1=1.$$

Assume now that $n \geq s+2$.  Recall from Section \ref{sectionminsumsetsizeUfixed} (see page \pageref{Ad(n,m) defined}) that, for $1 \leq m \leq n$ and a divisor $d$ of $n$, we defined the set $A_d(n,m)$.  Suppose now that $d \in D(n)$, $d \geq s+2$ (possible since $n \geq s+2$), and $$m=\widehat{v}(n,s)=\left( \left \lfloor \frac{d-2}{s} \right \rfloor +1  \right) \cdot \frac{n}{d}.$$  Then, with $H$ as the order $n/d$ subgroup of $\mathbb{Z}_n$ and $c=\lfloor (d-2)/s \rfloor$, we get
$$A_{n/d} (n,m)=\cup_{i=0}^{c} (i+H).$$ Now $A_{n/d} (n,m)$ has size $m$, and 
$$[0,s]A_{n/d} (n,m)=\cup_{i=0}^{sc} (i+H)$$ has size $$(sc+1) \cdot \frac{n}{d}=\left(s \left \lfloor \frac{d-2}{s} \right \rfloor +1 \right) \cdot \frac{n}{d} \leq (d-1) \cdot \frac{n}{d} <n,$$ so $[0,s]A \neq \mathbb{Z}_n$.  Furthermore, since $d \geq s+2$, we have $c \geq 1$ and thus $1 \in A$, which implies that $A$ generates $\mathbb{Z}_n$.  This yields:

\begin{prop} \label{chi gen cyclic}
For all positive integers $n$ and $s$ we have
$$\widehat{\chi}  (\mathbb{Z}_n,[0,s]) \geq \widehat{v}(n,s) +1.$$
\end{prop}

Combining Proposition \ref{chi gen cyclic} with Theorem \ref{From Klopsch Lev}, we get:

\begin{thm} [Klopsch and Lev; cf.~\cite{KloLev:2009a}] \label{Klo Lev cyclic}\index{Klopsch, B.}\index{Lev, V. F.} 
For all positive integers  $n$ and $s$ we have
$$\widehat{\chi}  (\mathbb{Z}_n,[0,s])=
\widehat{v}(n,s) +1 .
$$

\end{thm} 

Considerably less is known about $\widehat{\chi}  (G,[0,s])$ for noncyclic $G$; in fact, in contrast to $\chi(G,[0,s])$, which only depends on the order $n$ of $G$ (see Theorem \ref{[0,s] crit numb} above), $\widehat{\chi}  (G,[0,s])$ is generally greatly dependent on the structure of $G$ itself.  

We have the following general lower bound:

\begin{prop} [Bajnok; cf.~\cite{Baj:2016b}]  \label{chi cap lower [0,s]}\index{Bajnok, B.}
Let $G$ be an abelian group of order $n$, and let $H$ be a subgroup of $G$ of index $d>1$ for which $G/H$ is of type $(d_1,\dots,d_t)$.  For each $i=1,\dots, t$, let $c_i$ be a positive integer with $c_i \leq d_i-1$, and suppose that
$$\Sigma_{i=1}^t \left \lceil (d_i-1)/c_i \right \rceil \geq s+1.$$  Then we have
$$\widehat{\chi}  (G,[0,s]) \geq \left ( 1+ \Sigma_{i=1}^t c_i \right) \cdot n/d +1.$$

\end{prop}

As an application, we consider $\mathbb{Z}_2^r$, the elementary abelian 2-group of rank $r$.  Trivially, when $r \leq s$, we have $$\widehat{\chi}  (\mathbb{Z}_2^r, [0,s]) =1,$$
so assume that $s+1 \leq r$, and let $t$ be an integer with $$s+1 \leq t \leq r.$$  Then choosing $H=\mathbb{Z}_2^t$ and $c_i=1$ for all $i \in \{1,\dots,t\}$, 
Proposition \ref{chi cap lower [0,s]} implies that 
$$\widehat{\chi}  (\mathbb{Z}_2^r, [0,s]) \geq  (t+1) \cdot 2^{r-t}+1;$$ in particular,
we have
$$\widehat{\chi}  (\mathbb{Z}_2^r, [0,s]) \geq  (s+2) \cdot 2^{r-s-1}+1.$$ 

It turns out that equality holds:  

\begin{thm} [Lev; cf.~\cite{Lev:2003a}]  \label{Lev 2003a}\index{Lev, V. F.} 
Let $r$ and $s$ be positive integers, $s \geq 2$.  If $r \leq s$, then $\widehat{\chi}  (\mathbb{Z}_2^r,[0,s])=1$; otherwise we have
$$\widehat{\chi}  (\mathbb{Z}_2^r,[0,s])=(s+2) \cdot 2^{r-s-1}+1.$$
\end{thm}   
We thus find that $$\widehat{\chi}  (\mathbb{Z}_2^r,[0,2])=\widehat{\chi}  (\mathbb{Z}_{2^r},[0,2])$$ for all $r \geq 3$, but
$$\widehat{\chi}  (\mathbb{Z}_2^r,[0,s])<\widehat{\chi}  (\mathbb{Z}_{2^r},[0,s])$$ for all $r \geq s+1 \geq 4$.

\begin{prob}
Find other applications of Proposition \ref{chi cap lower [0,s]}.

\end{prob}

The following additional results are known:
\begin{thm} [Klopsch and Lev; cf.~\cite{KloLev:2009a}]  \label{thm Klopsch Lev extremes}\index{Klopsch, B.}\index{Lev, V. F.} 
Suppose that $G$ is a finite abelian group of order $n \geq 2$, rank $r$, and invariant factorization $\mathbb{Z}_{n_1} \times \cdots \times \mathbb{Z}_{n_r};$  we set $$D=n_1+\cdots+n_r-r.$$ ($D$ is called the {\em positive diameter} of $G$.)  
\begin{enumerate}
  \item If $G \not \cong \mathbb{Z}_2$, then $\widehat{\chi}  (G,[0,1])=n$.
  \item If $G \not \cong \mathbb{Z}_2, \mathbb{Z}_2^2$, then $\widehat{\chi}  (G,[0,2])=\lfloor n/2 \rfloor +1$.
  \item If $G \not \cong \mathbb{Z}_2^r$, then 
$$\widehat{\chi}  (G,[0,3])=\left\{
  \begin{array}{cl}
\left( 1+ \frac{1}{d} \right) \cdot \frac{n}{3}+1 & \mbox{if} \; \mbox{$G$ has a subgroup whose order is congruent to 2 mod 3}\\
& \mbox{and which is not isomorphic to an elementary abelian 2-group,} \\
& \mbox{and $d$ is the minimum size of such a subgroup}; \\ \\
\left \lfloor \frac{n}{3}\right \rfloor +1& otherwise.
\end{array}
\right.$$

  \item If $G \not \cong \mathbb{Z}_2$, then $\widehat{\chi}  (G,[0,D-1])=r+2$.
  \item If $s \geq D$, then $\widehat{\chi}  (G,[0,s])=1$.
\end{enumerate}

\end{thm} 

We can observe that, according to Theorem \ref{thm Klopsch Lev extremes}, for $s \in \{1,2\}$ we have $$\widehat{\chi}  (G,[0,s])=\chi  (G,[0,s])=v_1(n,s)+1.$$  

To assess the case of $s=3$, recall from page \pageref{formula for v1(n,3)} that 
$$v_1(n,3) =\left\{
\begin{array}{ll}
\left(1+\frac{1}{p}\right) \frac{n}{3} & \mbox{if $n$ has prime divisors congruent to 2 mod 3,} \\ & \mbox{and $p$ is the smallest such divisor,}\\ \\
\left\lfloor \frac{n}{3} \right\rfloor & \mbox{otherwise.}\\
\end{array}\right.$$  
Similarly, we get
$$\widehat{v}_1(n,3) =\left\{
\begin{array}{ll}
\left(1+\frac{1}{d}\right) \frac{n}{3} & \mbox{if $n$ has divisors congruent to 2 mod 3 that are greater than 2,} \\ & \mbox{and $d$ is the smallest such divisor,}\\ \\
\left\lfloor \frac{n}{3} \right\rfloor & \mbox{otherwise.}\\
\end{array}\right.$$
(Note that $d$ need not be prime.)  
Therefore, we see that
$$\widehat{\chi}  (G,[0,3]) \leq \widehat{v}_1(n,3)+1 \leq \chi  (G,[0,3])=v_1(n,3)+1.$$ While equality holds throughout for cyclic groups, this may not be the case for noncyclic groups; for example, for $G=\mathbb{Z}_2^2 \times \mathbb{Z}_6$, we get
$\widehat{\chi}  (G,[0,3])=9$, $\widehat{v}_1(n,3)+1=10$, and $\chi  (G,[0,3])=v_1(n,3)+1=13$.

Combining all results above, we see that $\widehat{\chi}  (G,[0,s])$ has been determined for all $G$ and $s$, except for noncyclic groups of exponent more than two and for $4 \leq s \leq D-2$.  

\begin{prob}
Find $\widehat{\chi}  (G,[0,s])$ for every noncyclic group $G$ and every $4 \leq s \leq D-2$.
\end{prob}

Since the general problem is probably difficult, we offer the following special cases:

\begin{prob}
Find $\widehat{\chi}  (G,[0,4])$ for every noncyclic group $G$.
\end{prob}

\begin{prob}
Find $\widehat{\chi}  (G,[0,D-2])$ for every noncyclic group $G$.
\end{prob}

\begin{prob}
Find $\widehat{\chi}  (\mathbb{Z}_k^2,[0,s])$ for all $k \geq 4$ and $s \geq 4$.
\end{prob}

\begin{prob}
Find $\widehat{\chi}  (G,[0,s])$ for every $s \geq 4$ and every noncyclic group $G$ of odd order.
\end{prob}

While exact values of $\widehat{\chi}  (G,[0,s])$ might be difficult to get in general, we can find a tight upper bound for it.  Recall that by Theorem \ref{From Klopsch Lev}, we have
$$\widehat{\chi}  (G,[0,s]) \leq  \widehat{v}(n,s) +1.$$ 
Observe that, when $d \in D(n)$ for which $d \geq s+2$, then
\begin{eqnarray*}
\left( \left \lfloor \frac{d-2}{s} \right \rfloor +1  \right) \cdot \frac{n}{d}  & \leq & \left(  \frac{d-2}{s} +1  \right) \cdot \frac{n}{d} \\ \\
& = & \left(  \frac{s-2}{d} +1  \right) \cdot \frac{n}{s} \\ \\ 
& \leq & \left(  \frac{s-2}{s+2} +1  \right) \cdot \frac{n}{s} \\ \\ 
& = & \frac{2n}{s+2}.
\end{eqnarray*}
Therefore, we have
$$\widehat{\chi}  (G,[0,s]) \leq \frac{2n}{s+2} +1,$$ with equality if, and only if, $n$ is divisible by $s+2$.

We have the following extension of this result:   

\begin{thm} [Klopsch and Lev; cf.~\cite{KloLev:2009a}]  \label{thm Klopsch Lev upper}\index{Klopsch, B.}\index{Lev, V. F.} 
For every $G$  abelian group of order $n$ and integer $s$ we have 
$$\widehat{\chi}  (G,[0,s]) \leq \frac{2n}{s+2} +1;$$ furthermore, when $s \geq 3$, equality holds if, and only if, there is a subgroup $H$ of order $n/(s+2)$ in $G$ for which $G/H$ is cyclic. 
\end{thm}
Indeed, when $H$ is a subgroup of order $n/(s+2)$ in $G$ for which $G/H$ is cyclic with $g+H$ as a generator, then the set $$A=H \cup (g+H)$$ has size $|A|=2n/(s+2)$, and we have $\langle A \rangle =G$, but $$[0,s]A=\cup_{i=0}^s (ig+H) \neq G.$$

Alternately, as did Margotta\index{Margotta, M. T.} in \cite{Mar:2010a}, instead of $\widehat{\chi}  (G,[0,s])$, we may study the quantity $ \widehat{s}  (G,m)$: the minimum value of $s$ for which the $[0,s]$-fold sumset of every $m$-subset of $G$ that generates $G$ is $G$ itself, that is,
$$\widehat{s}  (G,m) = \min \{s \mid A \subseteq G, \langle A \rangle =G, |A| \geq m \Rightarrow [0,s]A=G \}.$$ While, in theory, it suffices to study only one of $\widehat{s}  (G,m)$ or $ \widehat{\chi}  (G,[0,s]) $, it may be possible to gain different results via the two different perspectives.  

To establish a lower bound for $\widehat{s}  (\mathbb{Z}_n,m) $ for all $n \geq m \geq 2$, we can observe that for the subset
$$A=\{0,1,2,\dots,m-1\}$$ of $\mathbb{Z}_n$  we get
$$[0,s]A=\{0,1,2,\dots,s(m-1)\},$$ which immediately implies the following:

\begin{prop} \label{s hat bound}
For every positive integer $n$ and $m$ with $n \geq m \geq 2$, we have $$\widehat{s}  (\mathbb{Z}_n,m) \geq \left \lfloor \frac{n+m-3}{m-1} \right\rfloor.$$

\end{prop}

It turns out that for prime values of $n$, equality holds in Proposition \ref{s hat bound}: Indeed, by Corollary \ref{rho vs p} (the generalization of the Cauchy--Davenport\index{Davenport, H.}\index{Cauchy, A--L.} Inequality), with
$$s=\left \lfloor \frac{p+m-3}{m-1} \right\rfloor,$$ for all $m$-subsets $A$ of $\mathbb{Z}_p$ we get
\begin{eqnarray*}
|[0,s]A| & \geq & |sA| \\
& \geq & \min\{p, sm-s+1\} \\
& = & \min \left\{p, \; \left \lfloor \frac{p+m-3}{m-1} \right\rfloor \cdot (m-1)+1 \right\} \\
& \geq & \min  \left\{p, \; \frac{p-1}{m-1} \cdot (m-1)+1 \right\} \\
& = & p.
\end{eqnarray*}
Therefore:

\begin{prop}
For every prime $p$ and positive integer $m$ with $p \geq m \geq 2$, we have $$\widehat{s}  (\mathbb{Z}_p,m) = \left \lfloor \frac{p+m-3}{m-1} \right\rfloor.$$

\end{prop}

While determining $\widehat{s}  (G,m)$ was easy when the order of $G$ is prime, this seems not to be the case when $G$ has subgroups of many different sizes.  As a case in point, observe that for every divisor $d$ of $n$, by Proposition \ref{chi gen cyclic} above we have
$$\widehat{\chi}  (\mathbb{Z}_n,[0,n/d-2]) \geq \widehat{v}(n,n/d-2) +1 \geq \left( \left \lfloor \frac{n/d-2}{n/d-2} \right \rfloor +1 \right) \cdot \frac{n}{n/d} +1=2d+1.$$
Therefore, if $d \geq m/2$, then there must be an $m$-subset $A$ of $\mathbb{Z}_n$ that generates $\mathbb{Z}_n$ but for which $$[0,n/d-2]A \neq \mathbb{Z}_n;$$ this yields the following lower bound:

\begin{prop}
Suppose that $n$ and $m$ are positive integers, and $d \in D(n)$ with $d \geq m/2$.  Then
$$\widehat{s}  (\mathbb{Z}_n,m) \geq n/d-1.$$

\end{prop}

It turns out that for $m \leq 5$ our two lower bounds above actually determine $\widehat{s}  (\mathbb{Z}_n,m)$; it can be shown that (for $n \geq m$):
\begin{eqnarray*}
\widehat{s}  (\mathbb{Z}_n,2) &=& n-1; \\ \\
\widehat{s}  (\mathbb{Z}_n,3) &=&\lfloor n/2 \rfloor; \\ \\
\widehat{s}  (\mathbb{Z}_n,4) &=& \left\{
\begin{array}{cl}
n/2-1 & \mbox{if $n$ is even}; \\ \\
\lfloor (n+1)/3 \rfloor & \mbox{if $n$ is odd};
\end{array}
\right. \\ \\
\widehat{s}  (\mathbb{Z}_n,5) &=& \left\{
\begin{array}{cl}
n/3-1 & \mbox{if $n$ is divisible by 3}; \\ \\
\lfloor (n+2)/4 \rfloor & \mbox{if $n$ is not divisible by 3}.
\end{array}
\right.
\end{eqnarray*}
(Margotta\index{Margotta, M. T.} in \cite{Mar:2010a} conjectured the first three formulae.)  However, as $m$ increases, the result becomes less transparent; for example, for $n \geq 7$ we get
$$ \widehat{s}  (\mathbb{Z}_n,6) = \left\{
\begin{array}{cl}
n/3-1 & \mbox{if $n$ is divisible by 3}; \\ \\
\lfloor n/4 \rfloor & \mbox{$n$ is even but not divisible by 3;} \\ \\
\lfloor (n+3)/5 \rfloor & \mbox{otherwise}.
\end{array}
\right.$$

We pose the following (potentially difficult) problems:

\begin{prob}
Find a concise formula for $ \widehat{s}  (\mathbb{Z}_n,m)$ for all $m \leq n$.

\end{prob}

\begin{prob}
Evaluate, or find bounds for $ \widehat{s}  (G,m)$ for arbitrary $m$ and noncyclic group $G$.  

\end{prob}

\subsection{Arbitrary number of terms} \label{CritUarbitrary}

Here we consider, for a given group $G$, the quantity $$\chi (G,\mathbb{N}_0) = \min \{m \mid A \subseteq G,  |A| \geq m \Rightarrow \langle A \rangle =G \},$$  where $\langle A \rangle$ is the subgroup of $G$ generated by $A$.   By Proposition \ref{rho sigma}, we have:  

\begin{prop} 
Let $G$ be any abelian group of order $n \geq 2$, and let $p$ be the smallest prime divisor of $n$.  Then
$$\chi (G,\mathbb{N}_0) = n/p+1.$$ 
\end{prop}

We also note that the variation 
$$\widehat{\chi}  (G,\mathbb{N}_0) = \min \{m \mid A \subseteq G, \langle A \rangle =G, |A| \geq m \Rightarrow \langle A \rangle =G \}$$
that would correspond to the analogous $\widehat{\chi}  (G,[0,s])$ of Subsection \ref{CritUlimited} is trivial and thus of no interest.

\section{Unrestricted signed sumsets} \label{critUS}

Our goal in this section is to investigate $\chi_{\pm} (G,H)$, the minimum value of $m$ for which $$H_{\pm} A=G$$ holds for every $m$-subset of $G$.  (Recall that $H_{\pm}A$ is the union of all $h$-fold signed sumsets $h_{\pm} A$ for $h \in H$.)  Since $0_{\pm}A=\{0\}$ for every subset $A$ of $G$ but $h_{\pm}G=G$ for every positive integer $h$, we see that $\chi_{\pm} (G,H)$ does not exist for $n \geq 2$ when $H=\{0\}$, but $\chi_{\pm} (G,H)$ does exist and is at most $n$ when $H$ contains at least one positive integer.  

 We consider three special cases: when $H$ consists of a single nonnegative integer $h$, when $H$ consists of all nonnegative integers up to some value $s$, and when $H$ is the entire set of nonnegative integers.

\subsection{Fixed number of terms} \label{CritUSfixed}

Our goal here is to find, for a given group $G$ and positive integer $h$, the quantity $$\chi_{\pm} (G,h)=\min \{m \mid A \subseteq G,  |A| \geq m \Rightarrow h_{\pm}A=G \}.$$  

It is easy to see that $$\chi_{\pm} (G,1)=n$$ for each group $G$: indeed, $1_{\pm}G=G$, but $1_{\pm}(G \setminus \{0\})=G \setminus \{0\}$.  

We can also evaluate $\chi_{\pm} (G,2)$.  First, observe that 
$$\chi_{\pm} (G,2) \leq \chi (G,2) = \lfloor n/2 \rfloor +1.$$  Clearly, if $n$ is even, then for a subgroup $H$ of order $n/2$ we have $2_{\pm}H=H$, so $\chi_{\pm} (G,2)$ cannot be $n/2$ or less.  When $n$ is odd, $G$ can be partitioned as
$$G=\{0\} \cup K \cup (-K);$$ here $0 \not \in 2_{\pm}K$, so  $\chi_{\pm} (G,2)$ cannot be $(n-1)/2$ or less.  In summary, we have:

\begin{prop}
For all groups $G$ of order $n$, we have
$$\chi_{\pm} (G,1)=n$$ and 
$$\chi_{\pm} (G,2)=\lfloor n/2 \rfloor +1.$$

\end{prop}

Furthermore, as an immediate consequence of Theorems \ref{cyclic} and \ref{h crit numb}, we get:

\begin{thm}  \label{chi pm cyclic h}
For all $n$ and $h$ we have
$$\chi_{\pm} (\mathbb{Z}_n,h)=v_1(n,h)+1.$$
\end{thm} 

This leaves us with the following problem:

\begin{prob}
Evaluate $\chi_{\pm} (G,h)$ for noncyclic groups $G$ and integers $h \geq 3$.
\end{prob}

\subsection{Limited number of terms} \label{CritUSlimited}

For given groups $G$ and positive integer $s$, here we consider $\chi_{\pm} (G,[0,s])$, that is, the minimum value of $m$ for which the $[0,s]$-fold signed sumset of every $m$-subset of $G$ is $G$ itself.  

By  Proposition \ref{rhoUStrivial}, we have:
\begin{prop}  \label{[0,1] pm crit numb}
For all finite abelian groups $G$ of order $n \geq 3$, we have $$\chi_{\pm}(G,[0,1])=
\left\{
\begin{array}{cll}
n-1 & \mbox{if} & \mbox{$n$ is odd,} \\ \\
n & \mbox{if} & \mbox{$n$ is even.} 
\end{array} \right.$$
\end{prop}

For $s \geq 2$, we do not know the value of $\chi_{\pm} (G,[0,s])$ in general, but we have the following obvious upper bound:

\begin{prop} \label{chi US upper triv}
For every $G$ and $s$ we have $$\chi_{\pm}(G,[0,s]) \leq \chi_{\pm}(G,s).$$

\end{prop}

We can also establish a lower bound in the case when $G$ is cyclic, as follows.  Let $d$ be any positive divisor of $n$, and let $H$ be a subgroup of order $n/d$ in $G$.  Consider the set
$$A=\bigcup_{i=- \lfloor (d-2)/(2s) \rfloor }^{ \lfloor (d-2)/(2s) \rfloor} (i+H).$$
We then see that $A$ has size 
$$\left( 2 \cdot  \left \lfloor \frac{d-2}{2s} \right \rfloor +1  \right) \cdot \frac{n}{d}.$$  Furthermore,
$$[0,s]_{\pm} A = \bigcup_{i=- s\lfloor (d-2)/(2s) \rfloor }^{s \lfloor (d-2)/(2s) \rfloor} (i+H),$$
so $[0,s]_{\pm} A$ has size
$$ \left( 2s \cdot  \left \lfloor \frac{d-2}{2s} \right \rfloor +1  \right)  \cdot \frac{n}{d} \leq (d-1) \cdot \frac{n}{d}<n,$$
and thus $[0,s]_{\pm} A \neq \mathbb{Z}_n$.
Recalling the function 
$$v_{\pm}(n,h)= \max \left\{ \left( 2 \cdot \left \lfloor \frac{d-2}{2h} \right \rfloor +1  \right) \cdot \frac{n}{d}  \mid d \in D(n) \right\}$$
from page \pageref{0.3.5.3 pm page}, 
we get:

\begin{prop}  \label{rho pm cyclic lower}
For all positive integers $n$ and $s$ we have $$\chi_{\pm}(\mathbb{Z}_n,[0,s]) \geq v_{\pm}(n,s)+1.$$

\end{prop}

Combining Propositions \ref{chi pm cyclic h}, \ref{chi US upper triv}, and \ref{rho pm cyclic lower}, we get:

\begin{prop} \label{rho pm cyclic upper lower}
For all positive integers $n$ and $s$ we have $$v_{\pm}(n,s)+1 \leq \chi_{\pm}(\mathbb{Z}_n,[0,s]) \leq v_1(n,s)+1.$$

\end{prop}

We believe that the lower bound is exact:

\begin{conj} \label{conj chi pm cyclic}
For all positive integers $n$ and $s$ we have $$\chi_{\pm}(\mathbb{Z}_n,[0,s])= v_{\pm}(n,s)+1.$$

\end{conj}

\begin{prob}
Prove (or disprove) Conjecture \ref{conj chi pm cyclic}.

\end{prob}

From Propositions \ref{nu pm small} and \ref{[0,1] pm crit numb}, we see that Conjecture \ref{conj chi pm cyclic} holds for $s=1$.  By Proposition \ref{rho pm cyclic upper lower}, Conjecture \ref{conj chi pm cyclic} also holds whenever $v_{\pm}(n,s)=v_1(n,s)$.  In particular, from Propositions \ref{v function bounds} and \ref{v pm upper bounds} we get:

\begin{prop}
When $n$ is even and $s \geq 2$, we have $$\chi_{\pm}(\mathbb{Z}_n,[0,s])= v_{\pm}(n,s)+1=v_1(n,s)+1=n/2+1.$$

\end{prop}

For the case when $n$ is odd and divisible by 3, by Proposition \ref{v pm upper bounds}, Conjecture \ref{conj chi pm cyclic} becomes:

\begin{conj}  \label{conj chi pm cyclic 3}
When $n$ is odd and divisible by 3 and $s \geq 3$, we have $$\chi_{\pm}(\mathbb{Z}_n,[0,s])= v_{\pm}(n,s)+1=n/3+1.$$

\end{conj}

As a modest step towards Conjecture \ref{conj chi pm cyclic}, we offer:

\begin{prob}
Prove (or disprove) Conjecture \ref{conj chi pm cyclic 3}.

\end{prob}

We can also prove Conjecture \ref{conj chi pm cyclic} for groups of prime order $p$.  Recall that, by Theorem \ref{Matzke limited prime}, we have  
$$\rho_{\pm} (\mathbb{Z}_p,m,[0,s]) =\min\{p,2s  \lfloor m/2 \rfloor +1\}.$$
Since for 
$$m \geq 2 \lfloor (p-2)/(2s) \rfloor +2$$
we have
$$2s  \lfloor m/2 \rfloor +1 \geq 2s \lfloor (p-2)/(2s) \rfloor +2s+1 \geq 2s \left( \frac{p-2-(2s-1)}{2s} \right) +2s+1 =p,$$ but for 
$$m \leq 2 \lfloor (p-2)/(2s) \rfloor +1$$
we have
$$2s  \lfloor m/2 \rfloor +1 \leq 2s \lfloor (p-2)/(2s) \rfloor +1 \leq p-1,$$
we get 
$$\chi_{\pm}  (\mathbb{Z}_p,[0,s]) = 2 \lfloor (p-2)/(2s) \rfloor +2.$$
Therefore, recalling Proposition \ref{nu pm forp}, this yields:

\begin{thm} \label{crit prime from KloLev}
For every prime $p$ and positive integer $s$, we have $$\chi_{\pm}  (\mathbb{Z}_p,[0,s]) = 
v_{\pm}(p,s)+1=2 \lfloor (p-2)/(2s) \rfloor +2.$$

\end{thm}

We know little about $\chi_{\pm} (G,[0,s])$ for noncyclic groups:

\begin{prob}
Find $\chi_{\pm} (G,[0,2])$ for all noncyclic groups $G$.

\end{prob}

\begin{prob}
Find $\chi_{\pm} (G,[0,s])$ for all noncyclic groups $G$ and positive integers $s \geq 3$.

\end{prob}

Analogously to $\widehat{\chi}  (G,[0,s])$ of Subsection \ref{CritUlimited} above, Klopsch and Lev\index{Klopsch, B.}\index{Lev, V. F.} in \cite{KloLev:2003a} investigated the quantity $ \widehat{\chi}_{\pm}  (G,[0,s])$: the minimum value of $m$ for which the $[0,s]$-fold signed sumset of every $m$-subset of $G$ that generates $G$ is $G$ itself, that is,
$$\widehat{\chi}_{\pm}  (G,[0,s]) = \min \{m \mid A \subseteq G, \langle A \rangle =G, |A| \geq m \Rightarrow [0,s]_{\pm}A=G \}.$$

It turns out that the following relative of the function 
$v_1(n,s)$ of page \pageref{0.3.5.3 page} will be useful: if $n \geq 2s+2$, we define
$$\widehat{v}_{\pm}(n,s)  = \max \left\{ \left(2 \left  \lfloor \frac{d-2}{2s} \right \rfloor +1  \right) \cdot \frac{n}{d}   \mid d \in D(n), d \geq 2s+2 \right\};$$ we adhere to the convention that the maximum element of the empty-set equals zero, and thus if $n \leq 2s+1$, we have $\widehat{v}(n,s)=0$.  (Note also that we omitted the index 1 as unnecessary here.)   

We can easily see that for all $n$ and $s$ we have
$$\widehat{\chi}_{\pm}  (\mathbb{Z}_n,[0,s]) \geq \widehat{v}_{\pm}(n,s) +1.$$  This is trivial if $s \geq \lfloor n/2 \rfloor$, so assume $s \leq \lfloor n/2 \rfloor -1$.   
Suppose that $d \in D(n)$, $d \geq 2s+2$, and $$m=\widehat{v}_{\pm}(n,s)=\left( 2\left \lfloor \frac{d-2}{2s} \right \rfloor +1  \right) \cdot \frac{n}{d}.$$  Then, with $H$ as the order $n/d$ subgroup of $\mathbb{Z}_n$ and $c=\lfloor (d-2)/(2s) \rfloor$, we define
$$A=\cup_{i=-c}^{c} (i+H).$$ Now $A$ has size $m$, and 
$$[0,s]_{\pm}A=\cup_{i=-sc}^{sc} (i+H)$$ has size $$(2sc+1) \cdot \frac{n}{d}=\left(2s \left \lfloor \frac{d-2}{2s} \right \rfloor +1 \right) \cdot \frac{n}{d} \leq (d-1) \cdot \frac{n}{d} <n,$$ so $[0,s]_{\pm}A \neq \mathbb{Z}_n$.  Furthermore, since $d \geq 2s+2$, we have $c \geq 1$ and thus $1 \in A$, which implies that $A$ generates $\mathbb{Z}_n$.  This yields:

\begin{prop} \label{chi pm gen cyclic}
For positive integers $n$ and $s$, we have
$$\widehat{\chi}_{\pm}  (\mathbb{Z}_n,[0,s]) \geq \widehat{v}_{\pm}(n,s) +1.$$
\end{prop}

As it turns out, we have equality in Proposition \ref{chi pm gen cyclic}:

\begin{thm} [Klopsch and Lev; cf.~\cite{KloLev:2003a}]\index{Klopsch, B.}\index{Lev, V. F.} 
Let $n$ and $s$ be positive integers; $n \geq 2$.  We have
$$\widehat{\chi}_{\pm}  (\mathbb{Z}_n,[0,s]) = \widehat{v}_{\pm}(n,s) +1.$$

\end{thm} 

As a special case, we have
\begin{cor} \label{cor prime from KloLev}
For every prime $p$ and positive integer $s$, we have $$\widehat{\chi}_{\pm}  (\mathbb{Z}_p,[0,s]) = 
\left\{
\begin{array}{cll}
2 \lfloor (p-2)/(2s) \rfloor +2 & \mbox{if} & s \leq (p-3)/2; \\ \\
1 & \mbox{if} & s \geq (p-1)/2.
\end{array}
\right.$$

\end{cor}

Let us turn to $\widehat{\chi}_{\pm}  (G,[0,s])$ for noncyclic $G$.  First, we recall that, by Theorem \ref{Lev 2003a}, we have 
$$\widehat{\chi}_{\pm}  (\mathbb{Z}_2^r,[0,s])=\widehat{\chi}  (\mathbb{Z}_2^r,[0,s])=(s+2) \cdot 2^{r-s-1}+1.$$

The following additional results are known:
\begin{thm} [Klopsch and Lev; cf.~\cite{KloLev:2003a}]  \label{thm Klopsch Lev extremes pm}\index{Klopsch, B.}\index{Lev, V. F.} 
Suppose that $G$ is a finite abelian group of order $n \geq 2$, rank $r$, and invariant factorization $\mathbb{Z}_{n_1} \times \cdots \times \mathbb{Z}_{n_r};$  we set $$D_{\pm}=\lfloor n_1/2 \rfloor +\cdots+  \lfloor n_r/2 \rfloor.$$ ($D_{\pm}$ is called the {\em diameter} of $G$.)  
\begin{enumerate}
  \item If $D_{\pm} \geq 2$, 
then $$\widehat{\chi}_{\pm}  (G,[0,1])=\left\{
\begin{array}{cll}
n-1 & \mbox{if} & \mbox{$n$ is odd,} \\ \\
n & \mbox{if} & \mbox{$n$ is even.} 
\end{array} \right.$$

  \item If $D_{\pm} \geq 3$, then $$\widehat{\chi}_{\pm}  (G,[0,2])=\left\{
\begin{array}{ll}
\frac{n_r-1}{2} \cdot \frac{n}{n_r} & \mbox{if} \; n_r \equiv 1 \;(4), \\ \\
\frac{n_r-1}{2} \cdot \frac{n}{n_r}+1 & \mbox{if} \; n_r \equiv 3 \;(4), \\ \\
\frac{n}{2} & \mbox{if $G \cong \mathbb{Z}_{2^k}$ for some $k \in \mathbb{N}$,} \\ \\
\frac{n}{2}+1  & \mbox{otherwise.}
\end{array} \right.$$

  \item If $D_{\pm} \geq 2$, then $$\widehat{\chi}_{\pm}  (G,[0,D_{\pm}-1])=2r+2-n(\mathbb{Z}_2)+2 \cdot \lfloor n(\mathbb{Z}_3)/2 \rfloor,$$
 where $n(\mathbb{Z}_2)$ and $n(\mathbb{Z}_3)$ are the number of $\mathbb{Z}_2$ and $\mathbb{Z}_3$ factors in the invariant factorization of $G$, respectively. 
  \item If $s \geq D_{\pm}$, then $\widehat{\chi}_{\pm}  (G,[0,s])=1$.
\end{enumerate}

\end{thm} 

Thus we see that $\widehat{\chi}_{\pm}  (G,[0,s])$ has been determined for all $G$ and $s$, except for noncyclic groups of exponent more than two and for $3 \leq s \leq D_{\pm}-2$.  

\begin{prob}
Find $\widehat{\chi}_{\pm}  (G,[0,s])$ for every noncyclic group $G$ and every $3 \leq s \leq D_{\pm}-2$.
\end{prob}

Since the general problem is probably difficult, we offer the following special cases:

\begin{prob}
Find $\widehat{\chi}_{\pm}  (G,[0,3])$ for every noncyclic group $G$.
\end{prob}

\begin{prob}
Find $\widehat{\chi}_{\pm}  (G,[0,D_{\pm}-2])$ for every noncyclic group $G$.
\end{prob}

While exact values of $\widehat{\chi}_{\pm}  (G,[0,s])$ might be difficult to get in general, there is a tight upper bound for it.  Recall that by Proposition \ref{chi pm gen cyclic}, for all $n$ and $s$  we have
$$\widehat{\chi}_{\pm}  (\mathbb{Z}_n,[0,s]) \geq \widehat{v}_{\pm}(n,s) +1;$$ in the case when $n$ happens to be divisible by $2s+2$, we can further see from our formula for $\widehat{v}_{\pm}(n,s)$ that
$$\widehat{\chi}_{\pm}  (\mathbb{Z}_n,[0,s]) \geq \frac{3n}{2s+2} +1.$$  As  Klopsch and Lev\index{Klopsch, B.}\index{Lev, V. F.} proved in \cite{KloLev:2003a}, this is as good a bound as one can get:   

\begin{thm} [Klopsch and Lev; cf.~\cite{KloLev:2003a}]  \label{thm Klopsch Lev upper pm}\index{Klopsch, B.}\index{Lev, V. F.} 
For every $G$  abelian group of order $n$ and integer $s \geq 3$ we have 
$$\widehat{\chi}_{\pm}  (G,[0,s]) \leq \frac{3n}{2s+2} +1;$$ furthermore, equality holds if, and only if, there is a subgroup $H$ of order $n/(2s+2)$ in $G$ for which $G/H$ is cyclic. 
\end{thm}

\subsection{Arbitrary number of terms} \label{CritUSarbitrary}

This subsection is identical to Subsection \ref{CritUarbitrary}.

\section{Restricted sumsets} \label{critR}

Our goal in this section is to investigate $\chi \hat{\;} (G,H)$, the minimum value of $m$ for which $$H \hat{\;} A=G$$ holds for every $m$-subset of $G$.  (Recall that $H \hat{\;} A$ is the union of all restricted $h$-fold sumsets $h \hat{\;} A$ for $h \in H$.)  In contrast to $\chi  (G,H)$ that exists for every group $G$ when $H$ contains at least one positive element, there are some less trivial situations for which $\chi \hat{\;} (G,H)$ does not exist.   

In the subsections below, we consider three special cases: when $H$ consists of a single nonnegative integer $h$, when $H$ consists of all nonnegative integers up to some value $s$, and when $H$ is the entire set of nonnegative integers.  

\subsection{Fixed number of terms} \label{CritRfixed}

Analogously to the $h$-critical number of a group $G$, we define the {\em restricted $h$-critical number of $G$} to be the minimum value of $m$ for which every $m$-subset $A$ of $G$ has $h \hat{\;} A=G$; this quantity, if exists, is denoted by $\chi \hat{\;} (G,h)$. 

Let us make some initial observations.  First, note that $0 \hat{\;} A=\{0\}$ for every subset $A$ of $G$, thus $\chi \hat{\;} (G,0)$ only exists when $n=1$ (in which case it obviously equals 1).  Second, since $1 \hat{\;} A=A$ for every $A \subseteq G$, the restricted $1$-critical number of any $G$ is clearly just $n$.  Third, if $h > n$, then $h \hat{\;}G=\emptyset$, and when $h=n$, then $h \hat{\;}G$ consists of exactly one element.  Furthermore, for $h=n-1$, we have $h \hat{\;}G=G$: to see this, note that 
$$(n-1) \hat{\;} G = \left\{- g+\Sigma_{g \in G} g  \mid g \in G \right\}=G,$$ so $\chi \hat{\;} (G,n-1)$ exists for all $G$ and is at most $n$.  Since we have $|(n-1) \hat{\;} A|=1$  for every subset $A$ of $G$ that has size $n-1$, we also see that $\chi \hat{\;} (G,n-1)=n.$
          
We summarize our findings, as follows:

\begin{prop} \label{restricted h-critical easy}
Let $G$ be an abelian group of order $n$.
\begin{enumerate}
  \item If $n=1$, then $\chi \hat{\;} (G,0)=1$; if $n \geq 2$, then $\chi \hat{\;} (G,0)$ does not exist.
  \item We have $\chi \hat{\;} (G,1)=n$.
 \item We have $\chi \hat{\;} (G,n-1)=n$.
  \item If $n=1$, then $\chi \hat{\;} (G,n)=1$; if $n \geq 2$, then $\chi \hat{\;} (G,n)$ does not exist.
  \item If $h>n$, then $\chi \hat{\;} (G,h)$ does not exist.
\end{enumerate}

\end{prop}

By Proposition \ref{restricted h-critical easy}, it suffices to investigate $\chi \hat{\;} (G,h)$ for $2 \leq h \leq n-2$.  The question then arises: when does $\chi \hat{\;} (G,h)$ exist?  This question is clearly equivalent to deciding when $h \hat{\;} G=G$ holds, for which the answer is provided by Theorem \ref{rho hat m=n}.  Therefore:

\begin{thm} \label{crit exists}
The restricted $h$-critical number $\chi \hat{\;} (G,h)$ of an abelian group $G$ of order $n$ exists for all $G$ and $1 \leq h \leq n-1$, except for the elementary abelian 2-group for $h=2$ or $h=n-2$. 

\end{thm}

For $h=2$ we recall Proposition \ref{upper bound for restricted 2-critical}, which can be rephrased as follows:
\begin{prop}   \label{2-chrom G} 
Suppose that $G$ is of order $n \geq 3$ and is not isomorphic to the elementary abelian 2-group.  Then $$\chi \hat{\;} (G,2) = \frac{n+|\mathrm{Ord}(G,2)|+3}{2}.$$  In particular, $$\chi \hat{\;} (\mathbb{Z}_n,2) = \lfloor n/2 \rfloor +2.$$

\end{prop}

For $h \geq 3$, we know the value of $\chi \hat{\;} (G,h)$ for all $h$ when $n$ is even:

\begin{thm} [Roth and Lempel; cf.~\cite{RotLem:1992a}] \label{Roth and Lempel even}\index{Roth, R. M.}\index{Lempel, A.}

Suppose that $G$ is an abelian group of even order $n \geq 12$.  

If $G \in \{\mathbb{Z}_2^r, \mathbb{Z}_2^{r-1} \times \mathbb{Z}_4\}$ and $h \in \{3, n/2-2\}$, then 
$\chi \hat{\;} (G,h) =n/2+2$.  

In all other cases:
$$\chi \hat{\;} (G,h) = \left \{
\begin{array}{cl}
n/2+1 & \mbox{if} \; 3 \leq h \leq n/2-2; \\ \\
h+3 & \mbox{if} \; n/2-1 \leq h \leq (n+|\mathrm{Ord}(G,2)|-3)/2; \\ \\
h+2 & \mbox{if} \; (n+|\mathrm{Ord}(G,2)|-1)/2 \leq h \leq n-2.
\end{array}
\right.$$

\end{thm}

Observe that the assumption that $n \geq 12$ is necessary: we have $\chi \hat{\;} (\mathbb{Z}_{10},3)=7$ as shown by the subset $A=\{1,2,4,6,8,9\}$ for which $0 \not \in 3 \hat{\;} A$.  We also note that Theorem \ref{Roth and Lempel even} was proved independently for cyclic groups by Bajnok; cf.~\cite{Baj:2014a}.\index{Bajnok, B.}

Let us now turn to the case of odd values of $n$, for which we know much less.  

First, we find a lower bound, as follows.  Assume that $A$ is an $(h+1)$-subset of $G$.  Then $$|h \hat{\;} A| =h+1 \leq n-1.$$ Therefore: 

\begin{prop} \label{triv lower for chi hat}
For all abelian groups $G$ of order $n$ and positive integers $h \leq n-2$ we have
$\chi \hat{\;} (G,h) \geq h+2$.  

\end{prop}

Next, we show that, in fact, equality holds in Proposition \ref{triv lower for chi hat} for every group $G$ of odd order and for all $$(n-1)/2 \leq h \leq n-2.$$  Let $A$ be an $(h+2)$-subset of $G$.  Then, by symmetry, $|h \hat{\;} A| =|2 \hat{\;} A|;$ since 
$$|A|=h+2 \geq  (n+3)/2,$$  
by Proposition \ref{2-chrom G} we have $$|h \hat{\;} A| =n.$$ This proves the following:

\begin{prop} \label{chi hat big h}
Let $G$ be an abelian group of odd order $n$, and suppose that $h$ is a positive integer with $$(n-1)/2 \leq h \leq n-2.$$  Then $\chi \hat{\;} (G,h) = h+2.$

\end{prop}

This leaves us with the following problem:

\begin{prob} \label{restricted critical number of odd}
For each abelian group $G$ of odd order $n$  and each $h$ with $$3 \leq h \leq (n-3)/2,$$ find the restricted $h$-critical number of $G$.

\end{prob}

We now summarize what we know about cyclic groups.  First, a lower bound.  Consider the set
$$A=\{1,2,\dots,\lfloor (n-2)/h \rfloor +h \}$$ in $\mathbb{Z}_n$.  We can easily see that
$$h \hat{\;} A=\{h(h+1)/2, h(h+1)/2+1, \dots,h \lfloor (n-2)/h \rfloor +h(h+1)/2 \};$$ in particular,
$$h(h+1)/2-1 \not \in h \hat{\;} A.$$
Therefore:

\begin{prop}  \label{chi hat lower cyclic}
For all positive integers $n$ and $h$ with $h \leq n-1$ we have
$$\chi \hat{\;} (\mathbb{Z}_n,h) \geq \lfloor (n-2)/h \rfloor +h+1.$$

\end{prop}

For cyclic groups of prime order, Theorem \ref{Dias Da Silva and Hamidoune} implies that equality holds in Proposition \ref{chi hat lower cyclic}:

\begin{thm}  \label{prime rest crit}
For any positive integer $h$ and prime $p$ with $h \leq p-1$ we have
$$\chi \hat{\;} (\mathbb{Z}_p,h)=\left \lfloor (p-2)/h\right \rfloor +h+1.$$

\end{thm}

For cyclic groups, this leaves us with the following question:

\begin{prob} \label{restricted critical number of Z_n}
Find the restricted $h$-critical number of $\mathbb{Z}_n$ for each odd composite value of $n$ and for $3 \leq h \leq (n-3)/2.$

\end{prob}

As an example, we mention that for $n=15$ we find the following values: \label{n=15 crit}
$$\chi \hat{\;} (\mathbb{Z}_{15},h)=\left\{
\begin{array}{cl}
15 & h=1 \\
9 & h=2, 3 \\
8 & h=4 \\
9 & h=5,6 \\
h+2 & h=7, 8, \dots, 13 \\
15 & h=14.
\end{array} \right.$$ 
The answers for $h \leq 2$ and for $h \geq 7$ follow from our results above; the rest were determined by a computer program (of course, $\chi \hat{\;} (\mathbb{Z}_{15},h)$ does not exist for $h=0$ or $h \geq 15$).  As these values indicate, Problem \ref{restricted critical number of Z_n} may be challenging in general.

Relying on Corollary \ref{cor h=3}, we have additional results for $h=3$:

\begin{prop} [Bajnok; cf.~\cite{Baj:2014a}] \label{h=3 rest crit}\index{Bajnok, B.}
Let $n$ be an arbitrary integer with $n \geq 16$.  
\begin{enumerate}
  \item If $n$ has prime divisors congruent to $2$ mod $3$ and $p$ is the smallest such divisor, then
  $$\chi \hat{\;} (\mathbb{Z}_n,3) \geq 
\left\{
\begin{array}{ll}
\left(1+\frac{1}{p} \right) \frac{n}{3} +2 & \mbox{if} \; n=3p;  \\ \\
\left(1+\frac{1}{p} \right) \frac{n}{3} +1 & \mbox{otherwise}. 
\end{array} \right.$$
  \item If $n$ has no prime divisors congruent to $2$ mod $3$, then
  $$\chi \hat{\;} (\mathbb{Z}_n,3) \geq 
\left\{
\begin{array}{ll}
\left \lfloor \frac{n}{3} \right \rfloor +4 & \mbox{if $n$ is divisible by $9$};  \\ \\
\left \lfloor \frac{n}{3} \right \rfloor +3 & \mbox{otherwise}. 
\end{array} \right.$$
\end{enumerate}

\end{prop}

Observe that the case when $n$ is even follows from Theorem \ref{Roth and Lempel even}, since $$\left(1+\frac{1}{2} \right) \frac{n}{3} +1=\frac{n}{2}+1;$$  and the case when $n$ is prime follows from Theorem \ref{prime rest crit} since
$$\left \lfloor \frac{p-2}{3} \right \rfloor +3+1 =
\left\{
\begin{array}{ll}
\left(1+\frac{1}{p} \right) \frac{p}{3} +3 & \mbox{if} \; p \equiv 2 \; \mbox{mod} \; 3;  \\ \\
\left \lfloor \frac{p}{3} \right \rfloor +3 & \mbox{otherwise}. 
\end{array} \right.
$$  

We make the following conjecture:

\begin{conj} \label{conj for chi hat}
For all values of $n \geq 16$, equality holds in Proposition \ref{h=3 rest crit}.

\end{conj}

We have verified that Conjecture \ref{conj for chi hat} holds for all values of $n \leq 50$, and by Theorems \ref{prime rest crit} and \ref{Roth and Lempel even}, it holds when $n$ is prime or even.  As additional support, we have the following:

\begin{thm} [Bajnok; cf.~\cite{Baj:2014a}] \label{conj rhohatforh=3 implies conj for chi hat}\index{Bajnok, B.}
Conjecture \ref{conj rhohatforh=3} implies Conjecture \ref{conj for chi hat}.

\end{thm}

It is worth mentioning the following special case of Conjecture \ref{conj for chi hat}:
\begin{conj}  \label{conj on 3-fold odd}
If $n \geq 31$ is an odd integer, then $$\chi \hat{\;} (\mathbb{Z}_n,3) \leq \tfrac{2}{5} n+1.$$

\end{conj}
(The additive constant could be adjusted to include odd integers less than 31.)  This conjecture was made by Gallardo, Grekos, et al.\index{Gallardo, L.}\index{Grekos, G.}\index{Habsieger, L.}\index{Hennecart, F.}\index{Landreau, B.}\index{Plagne, A.}in \cite{GalGre:2002a}, and (for large $n$) proved by Lev\index{Lev, V. F.} via the following more general result:

\begin{thm} [Lev; cf.~\cite{Lev:2002a}] \label{Lev on 3-fold} \index{Lev, V. F.} 
Let $G$ be an abelian group of order $n$ with $$n \geq 312 \cdot |\mathrm{Ord}(G,2)|+1235.$$  Then for any subset $A$ of $G$, at least one of the following possibilities holds:
\begin{itemize}
  \item $|A| \leq \tfrac{5}{13}n$; 
  \item $A$ is contained in a coset of an index-two subgroup of $G$;
  \item $A$ is contained in a union of two cosets of an index-five subgroup of $G$; or
  \item $3 \hat{\;} A=G$.
\end{itemize}

\end{thm}

So, in particular, if $n$ is odd, is at least 1235, and a subset $A$ of $G$ has size more than $2n/5$, then the last possibility must hold, so we get:

\begin{cor} [Lev; cf.~\cite{Lev:2002a}] \label{Lev on 3-fold cor} \index{Lev, V. F.} 
If $n \geq 1235$ is an odd integer, then $$\chi \hat{\;} (\mathbb{Z}_n,3) \leq \tfrac{2}{5} n+1.$$
\end{cor}
The bound on $n$ in Corollary \ref{Lev on 3-fold cor} can hopefully be reduced to the one in Conjecture \ref{conj on 3-fold odd}:

\begin{prob}
Prove that $$\chi \hat{\;} (\mathbb{Z}_n,3) \leq \tfrac{2}{5} n+1$$ holds for odd integer values of $n$ between 31 and 1235 (inclusive).

\end{prob}

As another special case of Conjecture \ref{conj for chi hat}, we have

\begin{conj}  \label{conj chi 3 5}
If $n \geq 83$ is odd and not divisible by five, then $$\chi \hat{\;} (\mathbb{Z}_n,3) \leq \tfrac{4}{11} n+1.$$ 

\end{conj}
Theorem \ref{Lev on 3-fold} does not quite yield Conjecture \ref{conj chi 3 5}: while a careful read of \cite{Lev:2002a} enables us to reduce the coefficient $5/13$ to $(3-\sqrt{5})/2$ (at least for large enough $n$), this is still higher than $4/11$.  Hence we pose:

\begin{prob}
Prove Conjecture \ref{conj chi 3 5}.

\end{prob}

Combining Theorem \ref{h crit numb} with Conjecture \ref{conj for chi hat}, we claim that, when $n \geq 11$, we have
$$\chi  (\mathbb{Z}_n,3)  \leq  \chi \hat{\;} (\mathbb{Z}_n,3) \leq \chi  (\mathbb{Z}_n,3) + 3.$$

Before closing this subsection, we should mention that, unlike in Subsection \ref{CritRarbitrary} below, there is no point in considering the quantity
$$\chi \hat{\;} (G^*,h)=\min \{m \mid A \subseteq G \setminus \{0\},  |A| \geq m \Rightarrow h \hat{\;} A=G \}.$$ (The study of critical numbers originated with the paper \cite{ErdHei:1964a} of Erd\H{o}s and Heilbronn,\index{Erd\H{o}s, P.}\index{Heilbronn, H.}   where they studied only subsets of $G \setminus \{0\}$.)  Indeed, we have the following easy result:

\begin{prop} \label{fixed h crit nonzero}
Let $G$ be a finite abelian group of order $n \geq 6$, and let $h$ be an integer with $2 \leq h \leq n-2$.  Suppose that ${\chi} \hat{\;} (G,h)$ exists (that is, if $h \in \{2,n-2\}$ then $G$ is not an elementary abelian 2-group; cf.~Proposition \ref{crit exists}).  Then
$$\chi \hat{\;} (G^*,h)={\chi} \hat{\;} (G,h).$$

\end{prop}
The short proof can be found on page \pageref{proof of fixed h crit nonzero}.

\subsection{Limited number of terms} \label{CritRlimited}

For a finite abelian group $G$ and a nonnegative integer $s$, we define the {\em restricted $[0,s]$-critical number of $G$} to be the minimum value of $m$ for which every $m$-subset $A$ of $G$ has $[0,s] \hat{\;} A=G$; this quantity, if exists, is denoted by $\chi \hat{\;} (G,[0,s])$. 

Let us make some initial observations.  First, note that $0 \hat{\;} A=\{0\}$ for every subset $A$ of $G$, thus $\chi \hat{\;} (G,[0,0])$ only exists when $n=1$ (in which case it obviously equals 1).  Second, since $1 \hat{\;} A=A$ for every $A \subseteq G$, the restricted $[0,1]$-critical number of any $G$ is clearly just $n$.  Third, since for each $s \geq 1$, $$G = 1 \hat{\;} G \subseteq [0,s] \hat{\;} G,$$ $\chi \hat{\;} (G,[0,s])$ exists (and is at most $n$) for all $G$ and $s \geq 1$.  
          
We summarize our observations, as follows:

\begin{prop} \label{restricted [0,s]-critical easy}
Let $G$ be an abelian group of order $n$, and let $s$ be a nonnegative integer.
\begin{enumerate}
  \item If $n=1$, then $\chi \hat{\;} (G,[0,0])=1$; if $n \geq 2$, then $\chi \hat{\;} (G,[0,0])$ does not exist.
  \item For all $s \geq 1$, $\chi \hat{\;} (G,[0,s])$ exists and is at most $n$.
  \item We have $\chi \hat{\;} (G,[0,1])=n$.
  
  \end{enumerate}

\end{prop}

Furthermore, $\chi \hat{\;} (G,s)$ is clearly an upper bound for $\chi \hat{\;} (G,[0,s])$, and $\chi (G,[0,s])$ is a lower bound for it, so by Theorem \ref{[0,s] crit numb}, we have:  

\begin{prop}

For all $G$ and $s \geq 1$, $$v_1(n,s)+1= \chi (G,[0,s]) \leq \chi \hat{\;} (G,[0,s]) \leq \chi \hat{\;} (G,s).$$

\end{prop}

For $s=2$ we see that by Proposition \ref{2-chrom G}, we have
$$\chi \hat{\;} (\mathbb{Z}_n,[0,2]) \leq  \chi \hat{\;} (\mathbb{Z}_n,2) = \lfloor n/2 \rfloor +2;$$ we can show that equality holds by finding a subset $A$ of $\mathbb{Z}_n$   for which $|A|=\lfloor n/2 \rfloor +1$ but $[0,2]\hat{\;} A \neq \mathbb{Z}_n$.  Indeed, we see that when $n$ is odd and
$$A=\{0,1,\dots,(n-1)/2\},$$ then $|A|=\lfloor n/2 \rfloor +1=(n+1)/2$ and $n-1 \not \in [0,2] \hat{\;} A;$ when $n$ is divisible by 4, then with 
$$A=\{0,1,\dots,n/4\} \cup \{n/2+1, n/2+2, \dots, 3n/4\},$$ $|A|=n/2+1$ and $n/2 \not \in [0,2] \hat{\;} A;$ and when $n-2$ is divisible by 4, then with 
$$A=\{0,1,\dots,(n-2)/4\} \cup \{n/2, n/2+1, \dots, (3n-2)/4\},$$ $|A|=n/2+1$ and $n/2-1 \not \in [0,2] \hat{\;} A.$
Therefore:

\begin{prop} \label{[0,2]-chrom G}
For all integers $n \geq 3$, $$\chi \hat{\;} (\mathbb{Z}_n,[0,2]) = \lfloor n/2 \rfloor +2.$$

\end{prop}
(We mention that Lemma 3.3 in \cite{Ham:1998a} says that if $A \subseteq \mathbb{Z}_n \setminus \{0\}$ and $|A| \geq n/2$, then $[1,2] \hat{\;} A = \mathbb{Z}_n$, but, as we have just seen, this is always false when $n$ is even---see \cite{Baj:2016a} for more information.)

Recall that, by Theorem \ref{Roth and Lempel even}, for even values of $n \geq 12$  we have
$$\chi \hat{\;} (\mathbb{Z}_n,3)  =n/2+1.$$
But for even $n$, we must have $$\chi \hat{\;} (\mathbb{Z}_n,[0,s])  \geq n/2+1$$ for all positive integers $s$, so for $s \geq 3$ we have
$$n/2+1 \leq \chi \hat{\;} (\mathbb{Z}_n,[0,s]) \leq \chi \hat{\;} (\mathbb{Z}_n,[0,3]) \leq \chi \hat{\;} (\mathbb{Z}_n,3)  =n/2+1.$$  Therefore:

\begin{thm}
For all even values of $n \geq 12$ and every $s \geq 3$, we have
$$\chi \hat{\;} (\mathbb{Z}_n,[0,s])=n/2+1.$$ 

\end{thm}

This leaves us with the following open question:

\begin{prob}  \label{prob chi hat [0,s] n odd}
Find $\chi \hat{\;} (\mathbb{Z}_n,[0,s])$ for all $s \geq 3$ and for odd values of $n$.

\end{prob}

We can get a lower bound for $\chi \hat{\;} (\mathbb{Z}_n,[0,s])$, as follows.  Suppose that $n \geq s^2-s+2$, and consider
$$A= \{0,1,\dots,\lfloor (n-2)/s) \rfloor \} \cup \{n-s+1,n-s+2,\dots,n-1\} ;$$ we may consider this set as the interval
$$A= \{-(s-1),-(s-2), \dots, -1, 0, 1, \dots, \lfloor (n-2)/s) \rfloor \}.$$
For the size of $A$ we have
$$|A|= \lfloor (n-2)/s) \rfloor +s$$ (note that this value is less than $n$).

By the assumption that $n \geq s^2-s+2$, we have $$ \lfloor (n-2)/s) \rfloor \geq s-1,$$ and thus 
for the $[0,s]$-fold restricted sumset of $A$ we get
$$[0,s] \hat{\;} A = \{-s(s-1)/2, -s(s-1)/2+1, \dots, s  \lfloor (n-2)/s) \rfloor - s(s-1)/2\},$$ which we can rewrite as
$$ [0,s] \hat{\;} A = \{0,1, \dots, s  \lfloor (n-2)/s) \rfloor - s(s-1)/2\}  \cup \{n-s(s-1)/2, n-s(s-1)/2+1, \dots, n-1\}.$$
Here $$0 \leq s  \lfloor (n-2)/s) \rfloor - s(s-1)/2 < n - s(s-1)/2 -1 < n-s(s-1)/2 <n;$$ so
$$n - s(s-1)/2 -1 \not \in  [0,s] \hat{\;} A.$$ 
We just proved the following:

\begin{prop} \label{lower for chi hat [0,s]}
For all positive integers $s$ and $n$ with $n \geq s^2-s+2$, we have
$$\chi \hat{\;} (\mathbb{Z}_n,[0,s]) \geq \lfloor (n-2)/s) \rfloor +s+1.$$
\end{prop}

We should mention that Proposition \ref{lower for chi hat [0,s]} is not tight in that the condition $n \geq s^2-s+2$ is not necessary for the conclusion to hold, and the bound does not always give the value of $\chi \hat{\;} (\mathbb{Z}_n,[0,s])$.  So Problem \ref{prob chi hat [0,s] n odd} is very much still open.

We know little about noncyclic groups:

\begin{prob}
Find $\chi \hat{\;} (G,[0,2])$ for all noncyclic groups $G$.

\end{prob}

\begin{prob}
Find $\chi \hat{\;} (G,[0,s])$ for all noncyclic groups $G$ and for all $s \geq 3$.

\end{prob}

\subsection{Arbitrary number of terms} \label{CritRarbitrary}

In this subsection we determine the {\em restricted critical number} of $G$, which we define as 
$$\chi \hat{\;} (G,\mathbb{N}_0) = \min \{m \mid A \subseteq G,  |A| \geq m \Rightarrow \Sigma  A=G \}$$  where, for $A=\{a_1,\dots,a_m\} \subseteq G$, 
$$\Sigma  A =\cup_{h=0}^{\infty} h \hat{\;} A = \{ \lambda_1 a_1 + \cdots + \lambda_m a_m \mid \lambda_1, \dots, \lambda_m \in \{0,1\} \}. $$  

Before doing so, we introduce three variations: 
\begin{eqnarray*}
\chi \hat{\;} (G^*,\mathbb{N}_0) & = &  \min \{m \mid A \subseteq G \setminus \{0\},  |A| \geq m \Rightarrow \Sigma  A=G \}, \\ \\
{\chi} \hat{\;} (G,\mathbb{N}) & = &  \min \{m \mid A \subseteq G,  |A| \geq m \Rightarrow \Sigma^*   A=G \}, \\ \\
\chi \hat{\;} (G^*,\mathbb{N}) & = &  \min \{m \mid A \subseteq G \setminus \{0\},  |A| \geq m \Rightarrow \Sigma^*   A=G \}, 
\end{eqnarray*}
where
$$\Sigma^*   A =\cup_{h=1}^{\infty} h \hat{\;} A = \{ \lambda_1 a_1 + \cdots + \lambda_m a_m \mid \lambda_1, \dots, \lambda_m \in \{0,1\}, \lambda_1 + \cdots + \lambda_m \geq 1 \}. $$

We can determine if these four quantities are well-defined, as follows.  Since $$G=1 \hat{\;} G \subseteq \Sigma^*   G \subseteq \Sigma  G,$$ we see that $\chi \hat{\;} (G,\mathbb{N}_0)$ and $\chi \hat{\;} (G,\mathbb{N})$ are well-defined (and are at most $n$) for any group $G$.  Similarly, $$G=\{0\} \cup 1 \hat{\;} (G \setminus \{0\})  \subseteq   \Sigma(G \setminus \{0\}),$$ so $\chi \hat{\;} (G^*,\mathbb{N}_0)$ is also well-defined (and at most $n-1$) for any group $G$ of order $n \geq 2$.  The same way, we see that $\chi \hat{\;} (G^*,\mathbb{N})$ is well-defined (and is at most $n-1$) if, and only if, $0 \in \Sigma^*   (G \setminus \{0\})$.  This is clearly the case if $G$ has an element of order three or more (the element and its inverse are distinct and add to zero); if $G$ is an elementary abelian 2-group of rank two or more, then, for example, $e_1=1000\dots$, $e_2=0100\dots$, and $e_1+e_2$ are three distinct elements that add to zero.  That leaves us with the group of order two, but ${\chi} \hat{\;} (\mathbb{Z}_2^*,\mathbb{N})$ cannot exist.  In summary:

\begin{prop}
The quantities $\chi \hat{\;} (G,\mathbb{N}_0)$ and ${\chi} \hat{\;} (G,\mathbb{N})$ are well-defined for every group $G$; $\chi \hat{\;} (G^*,\mathbb{N}_0)$ is well-defined for every group $G$ of order at least two; and $\chi \hat{\;} (G^*,\mathbb{N})$ is well-defined for every group $G$ of order at least three.

\end{prop}

The four quantities are strongly related---see Theorem \ref{four quantities} below.  The last quantity, $\chi \hat{\;} (G^*,\mathbb{N})$, is the one that has been studied most; it is in fact the one that has been coined the {\em critical number} of $G$.  Thus we begin our investigation with $\chi \hat{\;} (G^*,\mathbb{N})$.

First, following a construction of Erd\H{o}s and Heilbronn in \cite{ErdHei:1964a}\index{Erd\H{o}s, P.}\index{Heilbronn, H.}  that was improved by Griggs\index{Griggs, J. R.}  in \cite{Gri:2001a}, we show that $${\chi} \hat{\;} (\mathbb{Z}_n^*,\mathbb{N}) \geq \lfloor 2 \sqrt{n-2} \rfloor.$$

Assume that $n \geq 3$ (as noted above, ${\chi} \hat{\;} (\mathbb{Z}_2^*,\mathbb{N})$ does not exist).  Letting $k=\lfloor 2 \sqrt{n-2} \rfloor$, we set \label{set leading to Erdos Griggs}
$$A=\left \{
\begin{array}{ll}
\{\pm 1, \pm 2, \dots, \pm (k-1)/2 \} & \mbox{if $k$ is odd;} \\ \\
\{\pm 1, \pm 2, \dots, \pm (k-2)/2, k/2 \} & \mbox{if $k$ is even.}
\end{array}
\right.$$
Then $|A|=k-1$, and $\Sigma^*   A$ is an interval consisting of all integers  between $-(k^2-1)/8$ and $(k^2-1)/8$ (inclusive) when $k$ is odd, and  between $-(k^2-2k)/8$ and $(k^2+2k)/8$ (inclusive) when $k$ is even.  This yields 
$$|\Sigma^*    A|=\left \{
\begin{array}{ll}
\frac{k^2-1}{8}+1+\frac{k^2-1}{8} =\frac{k^2+3}{4} & \mbox{if $k$ is odd;} \\ \\
\frac{k^2-2k}{8}+1+\frac{k^2+2k}{8} =\frac{k^2+4}{4} & \mbox{if $k$ is even.}
\end{array}
\right.$$
Therefore, 
$$|\Sigma^*   A| \leq \frac{4(n-2)+4}{4} =n-1,$$
and we get

\begin{prop} [Cf.~\cite{ErdHei:1964a}, \cite{Gri:2001a}]  \label{sqrt lower for prime}
With $n \geq 3$, we have
$${\chi} \hat{\;} (\mathbb{Z}_n^*,\mathbb{N}) \geq \lfloor 2 \sqrt{n-2} \rfloor.$$

\end{prop}

It turns out that, when $n$ is prime, the lower bound of Proposition \ref{sqrt lower for prime} is sharp.  We will need the following inequality:

\begin{lem} \label{lemma with prime sqrt}
For an odd integer $n \geq 3$, $k=\lfloor 2 \sqrt{n-2} \rfloor$, and $h=\lfloor (k+1)/2 \rfloor$ we have
$$\lfloor (n-2)/h \rfloor +h \leq k.$$

\end{lem}

The short and easy proof of Lemma \ref{lemma with prime sqrt} is on page \pageref{proof of lemma with prime sqrt}.

We now show how Theorem \ref{prime rest crit} (via Lemma \ref{lemma with prime sqrt}) implies that, for any odd prime $p$, 
$${\chi} \hat{\;} (\mathbb{Z}_p^*,\mathbb{N}) \leq \lfloor 2 \sqrt{p-2} \rfloor.  $$   We follow  the proof of Dias Da Silva and Hamidoune in \cite{DiaHam:1994a}.\index{Dias Da Silva, J. A.}\index{Hamidoune, Y. O.}

The claim is obvious for $p=3$, so we assume that $p \geq 5$.  Consider any subset $A$ of $\mathbb{Z}_p \setminus \{0\}$ of size $k=\lfloor 2 \sqrt{p-2} \rfloor$.  Then $B=A \cup \{0\}$ has size $k+1$, so by Lemma \ref{lemma with prime sqrt}, for $h=\lfloor (k+1)/2 \rfloor$, we have
$$|B| \geq \lfloor (p-2)/h \rfloor +h +1.$$  
Note also that $h \leq p-1$.  Therefore, Theorem \ref{prime rest crit} implies that
$h \hat{\;} B = \mathbb{Z}_p$.  But  
$$h \hat{\;} B =h \hat{\;} (A \cup \{0\}) = h \hat{\;} A \cup (h-1) \hat{\;} A,$$
so (since $p \geq 5$ implies that $h \geq 2$) $\Sigma^*  A=\mathbb{Z}_p$, proving our claim.    

Combining this upper bound with the lower bound of Proposition \ref{sqrt lower for prime}, we get:

\begin{thm}  [Dias Da Silva and Hamidoune; cf.~\cite{DiaHam:1994a} and Griggs; cf.~\cite{Gri:2001a}]  \label{critical for prime}\index{Dias Da Silva, J. A.}\index{Hamidoune, Y. O.}\index{Griggs, J. R.}  If $p$ is an odd prime, then
$${\chi} \hat{\;} (\mathbb{Z}_p^*,\mathbb{N}) = \lfloor 2 \sqrt{p-2} \rfloor.  $$

\end{thm}

Let us now turn to groups of composite order.  We can easily find a lower bound for $\chi \hat{\;} (G^*,\mathbb{N})$, as follows.  

Following Diderrich's construction\index{Diderrich, G. T.} in \cite{Did:1975a}, we let $p$ denote the smallest prime divisor of $n$, and consider the set
$$A=(H \setminus \{0\}) \cup (g+K)$$ where  $H$ is a subgroup of $G$ with index $p$, $K$ is a subset of $H$ of size $p-2$, and $g$ is any element of $G \setminus H$.  (This is possible as $n$ being composite implies that $p-2 \leq n/p$.)  Since $(p-1) \cdot g$ (and, in fact, every element of $(p-1) \cdot g +H$) is outside of $\Sigma^*   A$, we get the following lower bound:

\begin{prop} [Cf.~\cite{Did:1975a}]  \label{lower for composite}
Let $n$ be a composite integer with smallest prime divisor $p$.  Then for every abelian group $G$ of order $n$, we have
$$\chi \hat{\;} (G^*,\mathbb{N}) \geq n/p+p-2.$$
\end{prop}

It took about 35 years after Diderrich's lower\index{Diderrich, G. T.} bound to determine the value of $\chi \hat{\;} (G^*,\mathbb{N})$ for all groups $G$.  As it turns out, in most cases, the lower bound above is sharp.  In particular, we have the following results:

\begin{thm} [Diderrich and Mann; cf.~\cite{DidMan:1973a}] \label{Diderrich and Mann}\index{Diderrich, G. T.}\index{Mann, H. B.}
If $G$ is of even order $n \geq 4$, then  
$$\chi \hat{\;} (G^*,\mathbb{N}) =\left\{
\begin{array}{ll}
n/2+1 & \mbox{if $G \cong \mathbb{Z}_4, \mathbb{Z}_6, \mathbb{Z}_8, \mathbb{Z}_2^2$, or $\mathbb{Z}_2 \times \mathbb{Z}_4$}; \\ \\
n/2 & \mbox{otherwise}.
\end{array}
\right.$$

\end{thm}

\begin{thm} [Mann and Wou; cf.~\cite{ManWou:1986a}] \label{Mann and Wou}\index{Mann, H. B.}\index{Wou, Y. F.} 
Let $p$ be an odd prime.  We have  
$${\chi} \hat{\;} ( (\mathbb{Z}_p^{2})^*,\mathbb{N}) =\left\{
\begin{array}{ll}
2p-1 & \mbox{if $p=3$}; \\ \\
2p-2 & \mbox{otherwise}.
\end{array}
\right.$$

\end{thm}

\begin{thm} [Gao and Hamidoune; cf.~\cite{GaoHam:1999a}]  \label{Gao and Hamidoune}\index{Hamidoune, Y. O.}\index{Gao, W.}
Suppose that $G$ is an abelian group of odd order $n$.  Let $p$ be the smallest prime divisor of $n$.  If $n/p$ is a composite number, then
$$\chi \hat{\;} (G^*,\mathbb{N}) = n/p+p-2.$$
\end{thm}

These results leave us with the case of cyclic groups whose order $n$ is the product of two (not necessarily different) odd primes.  When the two primes are far from one another, we have the following:

\begin{thm} [Diderrich; cf.~\cite{Did:1975a}]  \label{Diderrich}\index{Diderrich, G. T.}
Suppose that $p$ and $q$ are odd primes.  If $q \geq 2p+1$, then 
$${\chi} \hat{\;} (\mathbb{Z}_{pq}^*,\mathbb{N}) = p+q-2.$$
\end{thm}

In the same paper, Diderrich also proved that\index{Diderrich, G. T.}
 $${\chi} \hat{\;} (\mathbb{Z}_{pq}^*,\mathbb{N}) \leq p+q-1$$ holds for all odd primes $p$ and $q$, and thus, by Proposition \ref{lower for composite}, the critical number of $\mathbb{Z}_{pq}$ is either $p+q-2$ or $p+q-1$.   We can observe that, by Proposition \ref{sqrt lower for prime}, we also have
$${\chi} \hat{\;} (\mathbb{Z}_{pq}^*,\mathbb{N}) \geq \lfloor 2 \sqrt{pq-2} \rfloor.$$

We claim that for odd integers $p$ and $q$ with $3 \leq p \leq q$, we have 
$$\lfloor 2 \sqrt{pq-2} \rfloor  \leq p+q-1,$$ with equality if, and only if, $$q \leq p+\lfloor 2 \sqrt{p-2} \rfloor +1.$$
Indeed, the first claim follows from the fact that $$2 \sqrt{pq-2}< p+q$$ is  equivalent to $$(q-p)^2+8>0.$$  
To prove the second claim, we can square and rearrange the inequality $$ 2 \sqrt{pq-2} \geq p+q-1$$ to get
$$(q-p-1)^2 \leq 4p-8,$$ from which the claim follows.  

Therefore, if $p$ and $q$ are odd integers with $$3 \leq p \leq q \leq p+\lfloor 2 \sqrt{p-2} \rfloor +1,$$ then 
$${\chi} \hat{\;} (\mathbb{Z}_{pq}^*,\mathbb{N}) \geq p+q-1.$$
Consequently, we have the following result:

\begin{thm} [Cf.~Diderrich; cf.~\cite{Did:1975a} and Griggs; cf.~\cite{Gri:2001a}] \label{griggs and easier}\index{Diderrich, G. T.}\index{Griggs, J. R.} 
Suppose that $p$ and $q$ are odd primes.  If $$p \leq q \leq p+\lfloor 2 \sqrt{p-2} \rfloor +1,$$ then 
$${\chi} \hat{\;} (\mathbb{Z}_{pq}^*,\mathbb{N}) = p+q-1=\lfloor 2 \sqrt{pq-2} \rfloor.$$
\end{thm}
(Griggs in \cite{Gri:2001a}\index{Griggs, J. R.} provided a more constructive proof for the lower bound than our argument above: he showed that, given the condition for $p$ and $q$ in Theorem \ref{griggs and easier}, the set 
 $$A=\{\pm 1, \pm 2, \dots, \pm (p+q-2)/2\}$$ does not generate the elements $\pm (pq-1)/2$ in $\mathbb{Z}_{pq}$.)

Finally, nearly a half century after Erd\H{o}s and Heilbronn\index{Erd\H{o}s, P.}\index{Heilbronn, H.}  posed the original problem of finding the critical number of an abelian group, the remaining case was decided in \cite{FreGaoGer:2009a} (see also \cite{FreGaoGer:2015a} for a correction):

\begin{thm} [Freeze, Gao, and Geroldinger; cf.~\cite{FreGaoGer:2009a}, \cite{FreGaoGer:2015a}] \label{Freeze et al}\index{Freeze, M.}\index{Gao, W.}\index{Geroldinger, A.}
Suppose that $p$ and $q$ are odd primes.  If $$p+\lfloor 2 \sqrt{p-2} \rfloor +1 < q < 2p+1,$$ then 
$${\chi} \hat{\;} (\mathbb{Z}_{pq}^*,\mathbb{N}) = p+q-2.$$
\end{thm}
(We note that $$p+\lfloor 2 \sqrt{p-2} \rfloor +1 <  2p+1$$ holds for all odd primes $p$.)

We can summarize these results as follows:

\begin{thm} \label{thm combined}
Suppose that $n \geq 3$ is an integer, let $p$ be the smallest prime divisor of $n$, and set $k=\lfloor 2 \sqrt{p-2} \rfloor$.  Then
$$\chi \hat{\;} (G^*,\mathbb{N})=\left\{
\begin{array}{cl}
k & \mbox{if $n=p$}; \\ \\
n/p+p-1 & \mbox{if $G \cong \mathbb{Z}_4, \mathbb{Z}_6, \mathbb{Z}_8, \mathbb{Z}_2 \times \mathbb{Z}_4$, or $\mathbb{Z}_3^2$}, \\
& \mbox{or $G$ is cyclic, $n/p$ is prime, and $3 \leq p \leq n/p \leq p+k+1$}; \\ \\ 
n/p+p-2 & \mbox{otherwise}.
\end{array}
\right.
$$
\end{thm}

As we noted above, for odd integers $p$ and $n/p$ with $$3 \leq p \leq n/p\leq p+\lfloor 2 \sqrt{p-2} \rfloor +1,$$ we have 
$$n/p+p  -1 = \lfloor 2 \sqrt{n-2} \rfloor .$$ Thus for $n \geq 10$ this allows for the second line in Theorem \ref{thm combined} to be combined with the first: 

\begin{cor} \label{combined simpler}

Suppose that $n \geq 10$, and let $p$ be the smallest prime divisor of $n$.  Then
$$\chi \hat{\;} (G^*,\mathbb{N})=\left\{
\begin{array}{ll}
\lfloor 2 \sqrt{n-2} \rfloor  & \mbox{if $G$ is cyclic of order $n=p$ or $n=pq$ where} \\
& \mbox{$q$ is prime and $3 \leq p \leq q \leq p+\lfloor 2 \sqrt{p-2} \rfloor+1$}, \\ \\ 
n/p+p-2 & \mbox{otherwise}.
\end{array}
\right.
$$
\end{cor}

With $\chi \hat{\;} (G^*,\mathbb{N})$ thus determined for any finite abelian group $G$, let us now turn to our other three quantities: $\chi \hat{\;} (G^*,\mathbb{N}_0)$, ${\chi} \hat{\;} (G,\mathbb{N}_0)$, and ${\chi} \hat{\;} (G,\mathbb{N})$.

First, we note the obvious facts that 
$$\chi \hat{\;} (G^*,\mathbb{N}_0) \leq \chi \hat{\;} (G^*,\mathbb{N})$$ and 
$${\chi} \hat{\;} (G,\mathbb{N}_0) \leq {\chi} \hat{\;} (G,\mathbb{N}).$$

Next, we show that 
$${\chi} \hat{\;} (G,\mathbb{N}) \leq \chi \hat{\;} (G^*,\mathbb{N})+1$$ and 
$${\chi} \hat{\;} (G,\mathbb{N}_0) \leq \chi \hat{\;} (G^*,\mathbb{N}_0)+1.$$
Indeed, if $A$ is a subset of $G$ of size ${\chi} \hat{\;} (G,\mathbb{N})-1$ so that $\Sigma^*   A \neq G$, then $A \setminus \{0\}$ is a subset of $G \setminus \{0\}$ of size at least ${\chi} \hat{\;} (G,\mathbb{N})-2$ and $\Sigma^*    (A \setminus \{0\}) \neq G$.  This implies our first inequality; the second can be shown similarly.

We now prove that for every group $G$ of order at least ten, we have $${\chi} \hat{\;} (G,\mathbb{N}_0) \geq \chi \hat{\;} (G^*,\mathbb{N})+1.$$  (We here ignore the exceptional cases of Theorems  \ref{Diderrich and Mann} and \ref{Mann and Wou} that may occur when $n \leq 9$.) 

Our strategy is to point to a(n already-mentioned) subset $A$ of $G$ for which (i) $0 \not \in A$, (ii) $0 \in  \Sigma^*  A$, (iii) $\Sigma^*  A \neq G$, and (iv) $|A|=\chi \hat{\;} (G^*,\mathbb{N})-1$.  Then, by (i) and (iv), $B=A \cup \{0\}$ has size $|B|=\chi \hat{\;} (G^*,\mathbb{N})$; and by (ii) and (iii), $\Sigma B = \Sigma^*  A \neq G$.  Therefore, 
$${\chi} \hat{\;} (G,\mathbb{N}_0) \geq |B|+1= \chi \hat{\;} (G^*,\mathbb{N})+1$$ follows.

To find a set $A$ in $G$ satisfying properties (i)--(iv), recall that, when $G$ is cyclic of order $n$, the set $A$ we exhibited on page \pageref{set leading to Erdos Griggs} has size $$|A|=k-1=\lfloor 2 \sqrt{n-2} \rfloor-1,$$  and satisfies (i), (ii), and (iii) above; furthermore, when $n$ is prime or a product of odd primes $p$ and $q$ with $$p \leq q \leq p+k+1,$$ then, by Theorem \ref{combined simpler}, $|A|=\chi \hat{\;} (G^*,\mathbb{N})-1$.  Additionally, with $p$ denoting the smallest prime divisor of the composite number $n$, the set leading to Proposition \ref{lower for composite} above also satisfies properties (i), (ii), and (iii); and by Theorem \ref{combined simpler}, in all remaining cases (that is, when $G$ is not cyclic, or $n$ is even, or $n/p$ is composite, or when $n/p$ equals a prime $q$ that is greater than $p+k+1$), its size $n/p+p-3$ equals $\chi \hat{\;} (G^*,\mathbb{N})-1$ as well.  Therefore, $${\chi} \hat{\;} (G,\mathbb{N}_0) \geq \chi \hat{\;} (G^*,\mathbb{N})+1,$$ as claimed.

The combination of our five inequalities enables us to evaluate each of $\chi \hat{\;} (G^*,\mathbb{N}_0)$, ${\chi} \hat{\;} (G,\mathbb{N}_0)$, and ${\chi} \hat{\;} (G,\mathbb{N})$ in terms of the already-determined value of $\chi \hat{\;} (G^*,\mathbb{N})$:

\begin{thm} \label{four quantities}
For any abelian group $G$ of order at least ten we have
$$\chi \hat{\;} (G^*,\mathbb{N}_0) = \chi \hat{\;} (G^*,\mathbb{N}) \lessdot  {\chi} \hat{\;} (G,\mathbb{N}_0) = {\chi} \hat{\;} (G,\mathbb{N}),$$ where $x \lessdot y$ means that $x$ is exactly one less than $y$.

\end{thm}
In particular, for ${\chi} \hat{\;} (G,\mathbb{N}_0)$ we have:
\begin{cor}  \label{corollary combined and simple}
Suppose that $n \geq 10$ is an integer, and let $p$ be the smallest prime divisor of $n$.  Then
$${\chi} \hat{\;} (G,\mathbb{N}_0)=\left\{
\begin{array}{ll}
\lfloor 2 \sqrt{n-2} \rfloor +1 & \mbox{if $G$ is cyclic of order $n=p$ or $n=pq$ where} \\
& \mbox{$q$ is prime and $3 \leq p \leq q \leq p+\lfloor 2 \sqrt{p-2} \rfloor+1$}, \\ \\ 
n/p+p-1 & \mbox{otherwise}.
\end{array}
\right.
$$
\end{cor}

Now that we have determined the restricted critical number of all finite groups, we move on to the inverse problem of classifying the extremal sets.   Given that we have four different versions for the critical number, we have four corresponding inverse problems:

\begin{enumerate}[{P}1]

\item \label{prob chi extreme}
Classify all subsets $A$ of $G$ with size ${\chi} \hat{\;} (G,\mathbb{N}_0)-1$ for which $\Sigma A \neq G$.  (We shall refer to these sets as P1-sets.)

\item \label{prob chi extreme}
Classify all subsets $A$ of $G \setminus \{0\}$ with  size ${\chi} \hat{\;} (G^*,\mathbb{N}_0)-1$ for which $\Sigma A \neq G$.  (We shall refer to these sets as P2-sets.)

\item \label{prob chi extreme}
Classify all subsets $A$ of $G$ with size ${\chi} \hat{\;} (G,\mathbb{N})-1$ for which $\Sigma^* A \neq G$.  (We shall refer to these sets as P3-sets.)

\item \label{prob chi extreme}
Classify all subsets $A$ of $G \setminus \{0\}$ with size ${\chi} \hat{\;} (G^*,\mathbb{N})-1$ for which $\Sigma^* A \neq G$.  (We shall refer to these sets as P4-sets.)

\end{enumerate}

We can prove that (for $n \geq 10$) the first three problems are equivalent and that the fourth problem incorporates the first three:

\begin{thm} \label{four inverse probs}
Let $G$ be an abelian group of order at least ten.  The following are equivalent:
\begin{itemize}
  \item $A$ is a P1-set in $G$;
  \item $0 \in A$ and $A \setminus \{0\}$ is a P2-set in $G$;
  \item $A$ is a P3-set in $G$.
\end{itemize}
Furthermore, each of the above implies that $0 \in A$ and $A \setminus \{0\}$ is a P4-set in $G$.
\end{thm}

The short proof is on page \pageref{proof of four inverse probs}.  We believe that P4 is also equivalent to the other three problems:

\begin{prob}
Prove that for groups $G$ of order ten or more, every P4-set in $G$ is also a P2-set in $G$.

\end{prob}
Note that having a P4-set in $G$ that is not a P2-set would mean that there is a subset $A$ of $G \setminus \{0\}$ of size ${\chi} \hat{\;} (G^*,\mathbb{N})-1$ for which $\Sigma^* A=G \setminus \{0\}$. 

Assuming that the four problems are indeed equivalent, we focus on P1:

\begin{prob} \label{prob chi extreme}
For each finite abelian groups $G$,
classify all subsets $A$ of size ${\chi} \hat{\;} (G,\mathbb{N}_0)-1$ for which $\Sigma A \neq G$.
\end{prob}

Problem \ref{prob chi extreme} is quite extensive; we discuss next what we already know. 

Let us first consider the case when $G$ is cyclic of order $n$.  It is helpful to view the elements of $\mathbb{Z}_n$ as integers in the interval $(-n/2, \; n/2]$.  Suppose that $A \subseteq \mathbb{Z}_n$, and let $$A=\{a_1,\dots,a_t,-a_{t+1},\dots,-a_m\}$$ with $a_1,\dots,a_t$ nonnegative and $-a_{t+1},\dots,-a_m$ negative (of course, we may have $t=0$ or $t=m$).  The {\em norm} of $A$, denoted by $||A||$, is the sum $a_1+\cdots +a_m$; for example, for the set $$A=\{0,2,5,8\}=\{0,2,5,-2\}$$ in $\mathbb{Z}_{10}$, we have $$||A||=0+2+5+2=9.$$

We will show that if $A$ is a subset of $\mathbb{Z}_n$ with norm $||A|| \leq n-2,$ then $\Sigma A \neq \mathbb{Z}_n$; in fact, we can easily see that, using our notations from above, $$a_1+\cdots+a_t +1 \not \in \Sigma A.$$  Indeed, this follows right away from our assumption, since it is equivalent to
$$a_1+\cdots+a_t +1<n-(a_{t+1}+ \cdots + a_m).$$  
Thus we have
\begin{prop} \label{prop p norm}
Let $A \subseteq \mathbb{Z}_n$.  If $||A|| \leq n-2$, then  $\Sigma A \neq \mathbb{Z}_n$. 

\end{prop}
Note also that if $\Sigma A \neq \mathbb{Z}_n$ for some $A \subseteq \mathbb{Z}_n$, then for any $b \in \mathbb{Z}_n$, $\Sigma (b \cdot A) \neq \mathbb{Z}_n$ as well, where $$b \cdot A=\{b \cdot a \mid a \in A\}$$ is a {\em dilate} of $A$; indeed, the size of $\Sigma (b \cdot A)$ cannot be more than the size of $\Sigma A$.

We have the following question:

\begin{prob}  \label{p for which norm small}
Find all odd prime values $p$ for which whenever $\mathbb{Z}_p$ contains a subset $A$ of size $$|A|={\chi} \hat{\;} (\mathbb{Z}_p,\mathbb{N}_0)-1=\lfloor 2 \sqrt{p-2} \rfloor$$ so that $\Sigma A \neq \mathbb{Z}_p$, then there is an element $b \in \{1,\dots,p-1\}$ for which $||b \cdot A|| \leq p-2$.

\end{prob}

We have verified that all primes less than 30 satisfy the requirements of Problem \ref{p for which norm small}, except for $p=17$: the subset $$A=\{0,1,3,4,5,12,14\}$$ of $\mathbb{Z}_{17}$ (for example) has size $|A|={\chi} \hat{\;} (\mathbb{Z}_{17},\mathbb{N}_0)-1=7$ and norm 21, does not generate the element 11, and has no dilate with a norm less than 21.  

Related to Problem \ref{p for which norm small}, we mention two results:

\begin{thm} [Nguyen, Szemer\'edi, and Vu; cf.~\cite{NguSzeVu:2008a}]  \label{Nguyen, Szemeredi, and Vu}\index{Nguyen, N. H.}\index{Vu, V. H.}\index{Szemer\'edi, E.} 
There is a positive constant $C$, so that whenever $A$ is a subset of $\mathbb{Z}_p$ for an odd prime $p$, if $|A| \geq 1.99 \sqrt{p}$  and $\Sigma A \neq \mathbb{Z}_p$, then there is an element $b \in \{1,\dots,p-1\}$ for which $$||b \cdot A|| \leq p+ C \sqrt{p}.$$

\end{thm}

\begin{thm} [Nguyen  and Vu; cf.~\cite{NguVu:2009a}] \label{Nguyen  and Vu}\index{Nguyen, N. H.}\index{Vu, V. H.} 
There is a positive constant $c$, so that whenever $A$ is a subset of $\mathbb{Z}_p$ for an odd prime $p$, if $\Sigma A \neq \mathbb{Z}_p$, then there is a subset $A'$ of $A$ of size at most $c p^{6/13} \log p$ and an element $b \in \{1,\dots,p-1\}$, for which $$||b \cdot (A \setminus A')|| < p.$$

\end{thm}

Theorem \ref{Nguyen, Szemeredi, and Vu} says that, for large values of $p$, if $A$ has size near ${\chi} \hat{\;} (\mathbb{Z}_p,\mathbb{N}_0)-1$ but $\Sigma A \neq \mathbb{Z}_p$, then $A$ has a dilate with norm not much larger than $p$; Theorem \ref{Nguyen and Vu} tells us that, again for large values of $p$, if $\Sigma A \neq \mathbb{Z}_p$, then after removing a relatively small subset from $A$, the remaining part has a dilate with norm less than $p$.  

The inverse problem for groups of prime order is thus not fully settled:

\begin{prob}

For each odd prime $p$, classify all subsets $A$ of $\mathbb{Z}_p$ of size ${\chi} \hat{\;} (\mathbb{Z}_p,\mathbb{N}_0)-1$ for which $\Sigma A \neq \mathbb{Z}_p$.

\end{prob}

Moving on to groups of composite order $n$, let us first consider the case when $n$ is even.  According to Corollary \ref{corollary combined and simple}, for $n \geq 10$ we have
$${\chi} \hat{\;} (G,\mathbb{N}_0)=n/2+1.$$  Clearly, if $A$ is a subgroup of $G$ of order $n/2$, then $\Sigma A \neq G$; it turns out that (for groups of order at least 16) the converse of this is true as well:

\begin{thm} \label{size n/2 not spanning}
Suppose that $G$ is an abelian group of order $n \geq 16$ and that $n$ is even.  Let $A$ be a subset of $G$ with $|A|=n/2$.  Then $\Sigma A \neq G$ if, and only if, $A$ is a subgroup of $G$.

\end{thm}  
We will show how results from the paper \cite{GaoHamLlaSer:2003a} by Gao, Hamidoune, Llad\'o, and Serra  imply\index{Hamidoune, Y. O.}\index{Gao, W.}\index{Llad\'o, A. S.}\index{Serra, O.} Theorem \ref{size n/2 not spanning}; see page \pageref{proof of size n/2 not spanning}.   We should point out that the lower bound of 16 on $n$ cannot be reduced: in $\mathbb{Z}_{14}$, for example, the set
$$\{0,\pm 1, \pm 2, \pm 3\}$$ is not a subgroup, yet it does not generate 7. 

Suppose now that $n$ is odd, divisible by 3, and $n/3$ is composite.  We then have the following result:

\begin{thm} \label{size n/3 not spanning}
Suppose that $G$ is an abelian group of order $n$ where $n$ is odd, divisible by 3, $n/3$ is composite, and $n \neq 27, 45$.  Let $A$ be a subset of $G$ with $|A|=n/3+1$.  Then $\Sigma A \neq G$ if, and only if, there is a subgroup $H$ of $G$ of size $n/3$ and an element $a \in A$ so that $$A=H \cup \{a\}.$$

\end{thm}  
Like for Theorem \ref{size n/2 not spanning}, our proof for Theorem \ref{size n/3 not spanning} will come from a careful reexamination of a corresponding result by Gao, Hamidoune, Llad\'o, and Serra\index{Hamidoune, Y. O.}\index{Gao, W.}\index{Llad\'o, A. S.}\index{Serra, O.}  in \cite{GaoHamLlaSer:2003a}; see page \pageref{proof of size n/3 not spanning}.  We can settle the case of $n=27$: The claim is false in $\mathbb{Z}_{27}$ as seen by the set $$\{0,\pm 1, \pm 2, \pm 3, \pm 4, 5\},$$ but true in the other two groups of order 27, as verified by the computer program \cite{Ili:2017a}.  The case of $n=45$ is still open and poses the following:

\begin{prob}
Decide if the claim of Theorem \ref{size n/3 not spanning} above is true for the two groups of order 45.

\end{prob}

In a similar fashion:

\begin{thm} \label{size n/p not spanning}
Suppose that $G$ is an abelian group of order $n$ where $n$ is odd, has smallest prime divisor $p \geq 5$, $n/p$ is composite, and $n \neq 125$.  Let $A$ be a subset of $G$ with $|A|=n/p+p-2$.  Then $\Sigma A \neq G$ if, and only if, 
$$A=H \cup A_1 \cup A_2,$$ where $H$ is a subgroup of $G$ of size $n/p$, and there is an element $a \in A$ so that $A_1 \subset a+H$ and $A_2 \subset -a+H$. 

\end{thm}
This result is essentially Theorem 4.1 in \cite{GaoHamLlaSer:2003a} by Gao, Hamidoune, Llad\'o, and Serra.\index{Hamidoune, Y. O.}\index{Gao, W.}\index{Llad\'o, A. S.}\index{Serra, O.}  Note that their assumption that $n \geq 7p^2+7p$ can be reduced to $n \geq 7p^2$ in our case, which, considering that $p \geq 5$, excludes only $n=125$.   

\begin{prob}
Decide if the claim of Theorem \ref{size n/p not spanning} above is true for the three groups of order 125.

\end{prob}

This leaves us with the case when $n$ is the product of two odd primes.  One such case was resolved by Qu, Wang, Wang, and Guo:

\begin{thm} [Qu, Wang, Wang, and Guo; cf.~\cite{Qu:2014a}]
Suppose that $p$ and $q$ are odd primes so that $q \geq 2p+3.$  Let $A$ be a subset of $\mathbb{Z}_{pq}$ of size $p+q-2$ for which $\Sigma A \neq \mathbb{Z}_{pq}$.  Then 
$$A=H \cup A_1 \cup A_2,$$ where $H$ is a subgroup of $G$ of size $q$, and there is an element $a \in A$ so that $A_1 \subset a+H$ and $A_2 \subset -a+H$. 

\end{thm}
(We should caution that the paper \cite{Qu:2014a} contains some inaccuracies: In Theorem A, the critical numbers of $\mathbb{Z}_2^2$ and cyclic groups of order $p^2$ with prime $p$ are given incorrectly; the other main result in the paper, regarding groups of even order, had been done previously (cf.~Theorem \ref{Diderrich and Mann}); and Example 4.1 is the wrong example for the point that the authors make.)

We separate the remaining cases of Problem \ref{prob chi extreme} into three parts:

\begin{prob} \label{prob extr p squared}
For each prime $p \geq 5$,
classify all subsets $A$ of size $2p-2$ of $\mathbb{Z}_p^2$ for which $\Sigma A \neq \mathbb{Z}_p^2$.
\end{prob}

\begin{prob} \label{prob extr p q small}
For all primes $p$ and $q$ with $$3 \leq p \leq q \leq p+ \lfloor 2 \sqrt{p-2} \rfloor +1,$$
classify all subsets $A$ of size $p+q-1$ of $\mathbb{Z}_{pq}$ for which $\Sigma A \neq \mathbb{Z}_{pq}$.
\end{prob}

\begin{prob} \label{prob extr p q large}
For all primes $p$ and $q$ with $$p+ \lfloor 2 \sqrt{p-2} \rfloor +2 \leq q \leq 2p+2,$$
classify all subsets $A$ of size $p+q-2$ of $\mathbb{Z}_{pq}$ for which $\Sigma A \neq \mathbb{Z}_{pq}$.
\end{prob}

Recall that our examples for sets in Problems \ref{prob extr p squared} and \ref{prob extr p q large} were generated by cosets of subgroups, but for Problem \ref{prob extr p q small}, they came from sets with small norm.  

\section{Restricted signed sumsets} \label{critRS}

\subsection{Fixed number of terms} \label{CritRSfixed}

\subsection{Limited number of terms} \label{CritRSlimited}

\subsection{Arbitrary number of terms} \label{CritRSarbitrary}

\chapter{Zero-sum-free sets} \label{ChapterZerosumfree}

Recall that for a given finite abelian group $G$, $m$-subset $A=\{a_1,\dots, a_m\}$ of $G$, $\Lambda \subseteq \mathbb{Z}$, and $H \subseteq \mathbb{N}_0$, we defined the sumset of $A$ corresponding to $\Lambda$ and $H$ as
$$H_{\Lambda}A = \{\lambda_1a_1+\cdots +\lambda_m a_m \mbox{    } |  \mbox{    }  (\lambda_1,\dots ,\lambda_m) \in \Lambda^m(H)  \}$$
where the index set $\Lambda^m(H)$ is defined as
$$\Lambda^m(H)=\{(\lambda_1,\dots ,\lambda_m) \in \Lambda^m \; |  \; |\lambda_1|+\cdots +|\lambda_m| \in H \}.$$ 

In this chapter we investigate the maximum possible size of a zero-sum-free set over $\Lambda$ in a given finite abelian group $G$.   Namely, our objective is to determine, for any $G$, $\Lambda \subseteq \mathbb{Z}$, and $H \subseteq \mathbb{N}_0$ the quantity
$$\tau_{\Lambda}(G,H)=\mathrm{max} \{ |A|  \mid A \subseteq G, 0 \not \in H_{\Lambda}A\}.$$  If no zero-sum-free set exists, we put $\tau_{\Lambda}(G,H) = 0.$  

Since, by definition, for every subset $A$ of size $\chi_{\Lambda} (G,H)$ or more $H_{\Lambda}A$ equals the entire group, we have the following:

\begin{prop} \label{zersumfree less than crit}
Let $\Lambda \subseteq \mathbb{Z}$ and $H \subseteq \mathbb{N}_0$.  If $\chi_{\Lambda} (G,H)$ exists, then
$$\tau_{\Lambda}(G,H) \leq \chi_{\Lambda} (G,H)-1.$$

\end{prop}

In the following sections we consider $\tau_{\Lambda}(G,H)$ for special coefficient sets $\Lambda$.

\section{Unrestricted sumsets} \label{5maxU}

Our goal in this section is to investigate the maximum possible size of a zero-$H$-sum-free set, that is, the quantity $$\tau (G,H) =\mathrm{max} \{ |A|  \mid A \subseteq G, 0 \not \in HA \}.$$  Clearly, we always have $\tau(G,H)=0$ when $0 \in H$; in fact, $\tau(G,H)=0$ whenever $H$ contains any multiple of the exponent of $G$.  However, when $H$ contains no multiples of the exponent $\kappa$, then $\tau  (G,H) \geq 1$: for any $a \in G$ with order $\kappa$, at least for the one-element set $A=\{a\}$ we have $0 \not \in H A$.

It is often useful to consider $G$ of the form $G_1 \times G_2$.  (We may do so even when $G$ is cyclic if its order has at least two different prime divisors.)  It is not hard to see that, if $A_1 \subseteq G_1$ is zero-$H$-sum-free in $G_1$, then $$A=\{(a,g)  \mid a \in A_1, g \in G_2\}$$ is zero-$H$-sum-free in $G$.  Indeed, if $hA$ were to contain $(0,0)$ (the zero-element of $G=G_1 \times G_2$) for some $h \in H$, then we would have $h$ (not necessarily distinct) elements of $A$ adding to $(0,0)$; this would then mean that $h$ (not necessarily distinct) elements of $A_1$ would add to $0$ in $G_1$, a contradiction.  Thus, we have the following.    

\begin{prop} \label{zetafordirectsum}
For all finite abelian groups $G_1$ and $G_2$ and for all $H \subseteq \mathbb{N}_0$ we have
$$\tau(G_1 \times G_2,H) \geq \tau(G_1,H) \cdot |G_2|.$$

\end{prop}

Below we consider two special cases: when $H$ consists of a single positive integer $h$ and when $H$ consists of all positive integers up to some value $t$.  The cases when $H = \mathbb{N}_0$ or $H = \mathbb{N}$, as we just mentioned, yield no zero-sum-free sets.

\subsection{Fixed number of terms} \label{5maxUfixed}

In this section we investigate, for a given group $G$ and positive integer $h$, the quantity $$\tau (G,h) = \mathrm{max} \{ |A|  \mid A \subseteq G, 0 \not \in hA\},$$ that is, the maximum size of a zero-$h$-sum-free subset of $G$.

For $h=1$, we see that a subset of $G$ is zero-1-sum-free if, and only if, it does not contain 0, hence the unique maximal zero-1-sum-free set in $G$ is $G \setminus \{0\}$.

We can also easily determine the value of $\tau(G,2)$.  First, note that a zero-2-sum-free set $A$ cannot contain any element of $\{0\} \cup \mathrm{Ord}(G,2)$ (the elements of order at most 2), and neither can it contain any element with its negative; to get a maximum zero-2-sum-free set in $G$, take exactly one of each element or its negative in $G \setminus \mathrm{Ord}(G,2)\setminus \{0\}$. 
To summarize:
\begin{prop} \label{z for h=1 and h=2}
We have $$\tau(G,1)=n-1$$ and $$\tau(G,2)=\frac{n-|\mathrm{Ord}(G,2)|-1}{2};$$ in particular, $$\tau(\mathbb{Z}_n,2)=\lfloor (n-1)/2 \rfloor.$$

\end{prop}

Suppose now that $h=3$.  For the cyclic group $\mathbb{Z}_n$, we can find explicit zero-3-sum-free sets as follows.  For every $n$, the positive integers that are less than $n/3$ form a zero-3-sum-free set; that is, the set $$A=\left\{1, 2, \dots, \left\lfloor (n-1)/3 \right\rfloor \right\}$$ is zero-3-sum-free in $\mathbb{Z}_n$, since $$3A=\left\{3,4,\dots, 3 \cdot\left\lfloor (n-1)/3 \right\rfloor \right\}$$ does not contain 0. 

We can do better in some cases.  For example, when $n$ is even, we may take the larger set $$\{1,3,5,\dots,n-1\},$$ which is zero-3-sum-free since no three odd numbers can add to a number divisible by $n$ when $n$ is even.  More generally, suppose that $n$ has a prime divisor $p$ which is congruent to 2 mod 3.  It is not hard to see that the set 
$$\left\{(p+1)/3+ i+pj \mbox{    } | \mbox{    } i=0,1,\dots,(p-2)/3, \; j=0,1,\dots,n/p-1 \right\}$$
is zero-3-sum-free.  Indeed, if $$k=(p+1)+(i_1+i_2+i_3) +p(j_1+j_2+j_3)=0$$ in $\mathbb{Z}_n$ for some $$i_1,i_2,i_3 \in \left\{0,1,\dots,(p-2)/3\right\}$$ and $$j_1,j_2,j_3 \in \left\{0,1,\dots,n/p-1 \right\},$$ then the integer $k$ is divisible by $n$ and thus by $p$.  Therefore, $1+(i_1+i_2+i_3)$ would have to be divisible by $p$, but this is not possible as 
$$1 \leq 1+(i_1+i_2+i_3) \leq p-1.$$

Recalling the function $v_g(n,h)$ from page \pageref{0.3.5.3 page}, we see that we have established the following.

\begin{prop} \label{zn3}  We have
$$\tau(\mathbb{Z}_n,3) \geq v_3(n,3) = \left\{
\begin{array}{ll}
\left(1+\frac{1}{p}\right) \frac{n}{3} & \mbox{if $n$ has prime divisors congruent to 2 mod 3,} \\ & \mbox{and $p$ is the smallest such divisor;}\\ \\
\left\lfloor \frac{n-1}{3} \right\rfloor & \mbox{otherwise.}\\
\end{array}\right.$$

\end{prop}

We believe that equality holds in Proposition \ref{zn3}, but no proof of this has been discovered.  More generally, we make the following conjecture.

\begin{conj} \label{zconj}
For all positive integers $h$ and $n$, the maximum size of a zero-$h$-sum-free set in the cyclic group $\mathbb{Z}_n$ is given by
$$\tau(\mathbb{Z}_n,h) = v_h(n,h).$$
\end{conj} 

We have already seen that Conjecture \ref{zconj} holds when $h=1$ or $h=2$, and we have the following result.

\begin{thm} \label{thmz}
For all positive integers $h$ and $n$, we have
$$v_h(n,h) \leq \tau(\mathbb{Z}_n,h) \leq v_1(n,h).$$

\end{thm} 

The upper bound here follows directly from Proposition \ref{zersumfree less than crit} and Theorem \ref{h crit numb}; the proof of the lower bound can be found on page \pageref{proof of thm thmz}.

The lower and upper bounds in Theorem \ref{thmz} often agree---in which case the value of $\tau(\mathbb{Z}_n,h)$ is determined.  For example, for $h=3$, we have $v_3(n,3)=v_1(n,3)$ if, and only if, $n$ has at least one prime divisor which is congruent to 2 mod 3 or if all its prime divisors are congruent to 1 mod 3.  So, for $h=3$, the first few values of $n$ for which $v_3(n,3) \not =v_1(n,3)$ are $$n=3, 9, 21, 27, 39, 57, 63, \dots.$$  A potentially difficult problem is the following.

\begin{prob}

Prove or disprove Conjecture \ref{zconj}; in particular, settle the case for $h=3$.

\end{prob} 

It is worth mentioning that, using Theorem \ref{thmz},  we can determine the size of the largest zero-$h$-sum-free set in cyclic groups of prime order.  According to Proposition \ref{nuforp}, for a prime $p$ and an arbitrary positive integer $h$ we have 
$$v_1(p,h)=\left \lfloor \frac{p-2}{h} \right \rfloor +1$$
and  
$$v_h(p,h) =\left\{
\begin{array}{ll}
0 & \mbox{if $p|h$,} \\ \\
\left \lfloor \frac{p-2}{h} \right \rfloor +1& \mbox{otherwise.}\\
\end{array}\right.$$
Therefore, if $h$ is not divisible by $p$, then $v_1(p,h)=v_h(p,h)$, and by Theorem \ref{thmz} this yields $\tau(\mathbb{Z}_p,h)$.  On the other hand, if $h$ is divisible by $p$, then clearly no (nonempty) subset of $\mathbb{Z}_p$ is zero-$h$-sum-free since for every element $a$ of the group we have $ha=0$.  In summary, we have the following.

\begin{thm} \label{zforp}
The size of the largest zero-$h$-sum-free set in  cyclic group of prime order $p$ is 
$$\tau(\mathbb{Z}_p,h) =v_h(p,h)=\left\{
\begin{array}{ll}
0 & \mbox{if $p|h$,} \\ \\
\left \lfloor \frac{p-2}{h} \right \rfloor +1& \mbox{otherwise.}\\
\end{array}\right.$$ 
\end{thm}

At this point we know considerably less about the value of $\tau(G,h)$ for $h \geq 3$ when $G$ is not cyclic.  An immediate consequence of Proposition \ref{zetafordirectsum} and Theorem \ref{thmz} is the following.

\begin{cor} \label{cor noncyclic tau bounds}
Suppose that $G$ is an abelian group of order $n$ and exponent $\kappa$.  Then, for all positive integers $h$ we have
$$v_h(\kappa,h) \cdot \frac{n}{\kappa} \leq \tau(G,h) \leq v_1(n,h).$$
\end{cor}

When the two bounds above coincide, then of course we have the exact value of $\tau(G,h)$.  Two such instances are worth mentioning: one coming from Corollary \ref{v function bounds}, and the other from our formulas for $v_1(n,3)$ and $v_3(n,3)$ on page \pageref{formula for v1(n,3)}:

\begin{cor}
If $n$ is even and $h \geq 3$ is odd, then $$\tau(G,h) = v_1(n,h)=n/2.$$
\end{cor}
\begin{cor}
If $n$ is divisible by a prime $p$ with $p \equiv 2$ mod 3 and $p$ is the smallest such prime, then 
$$\tau(G,3) = v_1(n,3)=\left(1+\frac{1}{p}\right) \frac{n}{3}.$$

\end{cor}

We also know the value of $\tau(G,h)$ when $h$ is relatively prime to $n$:

\begin{thm}  \label{tau=chi -1 when h and n are coprime}
If $h$ is relatively prime to $n$, then $$\tau(G,h)=\chi (G,n)-1=v_1(n,h).$$

\end{thm}
We verify Theorem \ref{tau=chi -1 when h and n are coprime} following an idea of Lemma 2.1 in the paper \cite{RotLem:1992a} of Roth\index{Roth, R. M.}\index{Lempel, A.} and Lempel. Note that it suffices to find a set $A$ in $G$ of size $\chi (G,n)-1$ for which $0 \not \in hA$.  Since $h$ and $n$ are relatively prime, we have a positive integer $h'$ for which $hh' \equiv 1 $ mod $n$.  Now let $B$ be a subset of $G$ of size $\chi (G,n)-1$ for which $hB \neq G$, choose an element $g \in G \setminus hB$, and set $A=-h'g+B$.  As $$hA=-hh'g+hB=-g+hB,$$ we indeed find that $0 \not \in hA$.

Next, for a positive prime $p$ and $r \geq 1$, we consider $\mathbb{Z}_p^r$, the elementary abelian $p$-group of rank $r$.  Let $h \geq 2$.  When $h$ is divisible by $p$, then, as we have mentioned above, $\tau (\mathbb{Z}_p^r, h)=0$.  If $h$ is not divisible by $p$, then $p$ cannot be divisible by $h$ either (since that would imply $p=h$ and thus $p|h$).  If $p$ leaves a remainder of at least 2 mod $h$, then, by Corollary \ref{cor noncyclic tau bounds} and by Propositions \ref{nuforp} and \ref{v for prime power}, $$\tau (\mathbb{Z}_p^r, h)=v_1(p^r,h)=p^{r-1} \cdot \left( 1+ \lfloor p/h \rfloor \right).$$
This leaves one case: when $p \equiv 1$ mod $h$, which can be treated using Finn's\index{Finn, C.} construction in \cite{Fin:2013a} that we explain (in a generalized form) next.   

Suppose that $G_i$ is a finite abelian group for $i=1,\dots,r$; $A_i$ is a zero-$[1,h]$-sum-free set in $G_i$ for $i=1,\dots,r-1$; and  $A_r$ is a zero-$h$-sum-free set in $G_r$.  (We say that $A$ is zero-$[1,h]$-sum-free in $G$ if $0 \not \in \cup_{j=1}^h jA$---see Section \ref{5maxUlimited}.)  Now set
$$G=G_1 \times \cdots \times G_r$$ and
$$A=(A_1 \times G_2 \times \cdots \times G_r) \cup (\{0\} \times A_2 \times G_3 \times \cdots \times G_r) \cup \dots \cup (\{0\} \times \cdots \times \{0\} \times A_r).$$      
It is easy to verify that $A$ is zero-$h$-sum-free in $G$, and therefore we get the following result:

\begin{prop} [Finn; cf.~\cite{Fin:2013a}] \label{Finn's constr}\index{Finn, C.}
For all abelian groups $G_1, \dots, G_r$ we have
$$
\tau(G_1 \times \cdots \times G_r,h) \geq  \sum_{i=1}^r \tau (G_i, [1,h]) \cdot \Pi_{j=i+1}^r |G_j|
.$$

\end{prop}   
(Of course, $\Pi_{j=r}^r |G_j| = |G_r|$ and $\Pi_{j=r+1}^r |G_j| = 1$; $\tau (G,[1,h])$ denotes the maximum size of a zero-$[1,h]$-sum-free set in $G$.)

We can apply Proposition \ref{Finn's constr} to the case when each $G_i$ is cyclic by observing that the set
$$\{1,2,\dots, \lfloor (n-1)/h \rfloor\}$$ is zero-$[1,h]$-free (and thus also zero-$h$-sum-free) in $\mathbb{Z}_n$.  In particular, when $p$ is a prime with $p \equiv 1$ mod $h$, then $$\tau (\mathbb{Z}_p,[1,h]) \geq (p-1)/h,$$
so in this case
$$\tau(\mathbb{Z}_p^r, h) \geq \frac{p-1}{h} \cdot (p^{r-1}+ \cdots + p +1) = \frac{p^r-1}{h}=v_1(p^r,h).$$

Therefore, we have the following:

\begin{thm} [Finn; cf.~\cite{Fin:2013a}] \label{tau elementary}\index{Finn, C.}
Let $p$ be a positive prime and $r \geq 1$.  

If $h$ is divisible by $p$, then $\tau (\mathbb{Z}_p^r, h)=0$.
 
If $h$ is not divisible by $p$, then $$\tau (\mathbb{Z}_p^r, h)=v_1(p^r,h)= \left \{
\begin{array}{cll}
(p^r-1)/h & \mbox{if} & p \equiv 1 \; \mbox{mod} \; h; \\ \\
p^{r-1} \cdot \left( 1+ \lfloor p/h \rfloor \right) & \mbox{if} & p \not \equiv 1 \; \mbox{mod} \; h.
\end{array}
\right.$$

\end{thm}
Note that Theorem \ref{tau elementary} is a generalization of Theorem \ref{zforp}.

The general problem of finding $\tau(G,h)$ for noncyclic groups $G$ is quite intriguing, even for $h=3$.

\begin{prob}

Determine the value of $\tau(G,h)$ for noncyclic $G$.

\end{prob}

\subsection{Limited number of terms} \label{5maxUlimited}

Here we investigate, for a given group $G$ and positive integer $t$, the quantity $$\tau (G,[1,t]) = \mathrm{max} \{ |A|  \mid A \subseteq G, 0 \not \in [1,t]A\},$$ that is, the maximum size of a zero-$[1,t]$-sum-free subset of $G$.

For $t=1$ and $t=2$ our considerations of Section \ref{5maxUfixed} above yield the answers here as well:

\begin{prop}
We have $$\tau(G,[1,1])=n-1$$ and $$\tau(G,[1,2])=\frac{n-|\mathrm{Ord}(G,2)|-1}{2};$$ in particular, $$\tau(\mathbb{Z}_n,[1,2])=\lfloor (n-1)/2 \rfloor.$$

\end{prop}

The following general upper bound is easy to establish (see Corollary 2.3 in \cite{Alo:1987a}):

\begin{thm} [Alon; cf.~\cite{Alo:1987a}] \label{Alon bound}\index{Alon, N.}
In any group of order $n$ we have
$$\tau(G,[1,t]) \leq  \lfloor (n-1)/t \rfloor.$$

\end{thm}

Note that the set
$$\{1,2,\dots, \lfloor (n-1)/t \rfloor\}$$ is zero-$[1,t]$-free in $\mathbb{Z}_n$, hence we get:

\begin{cor} 
 For all positive integers $t$ and $n$, we have
$$\tau (\mathbb{Z}_n,[1,t]) =\lfloor (n-1)/t \rfloor.$$
\end{cor}

The case of noncyclic groups remains open:

\begin{prob}
Find the value of $\tau (G, [1,t])$ for noncyclic groups $G$ for $t \geq 3$.

\end{prob}

\subsection{Arbitrary number of terms} \label{5maxUarbitrary}

Here we ought to consider $$\tau (G,H) =\mathrm{max} \{ |A|  \mid A \subseteq G, 0 \not \in HA \}$$ for the case when $H$ is the set of all nonnegative or all positive integers.  However, as we have already mentioned, we have $\tau (G,H)=0$ whenever $H$ contains a multiple of the exponent of the group (including $0$).  Thus, there are no zero-sum-free sets when the addition of an arbitrary number of terms is allowed.

\section{Unrestricted signed sumsets} \label{5maxUS}

Our goal in this section is to investigate the maximum possible size of a zero-sum-free set over the set of all integers, that is, the quantity $$\tau_{\pm} (G,H) =\mathrm{max} \{ |A|  \mid A \subseteq G, 0 \not \in H_{\pm}A \}.$$   

Clearly, we have $\tau_{\pm}  (G,H)=0$ whenever $H$ contains a multiple of the exponent of the group (including $0$).  However, when $H$ contains no multiples of the exponent $\kappa$, then $\tau_{\pm}  (G,H) \geq 1$: for any $a \in G$ with order $\kappa$, at least the one-element set $\{a\}$ will be zero-sum-free for $H$.

It is important to note \label{over Z not same} that Proposition \ref{zetafordirectsum} does not carry through to zero-sum-free sets over the set of all integers.  
For example, the subset $A_1=\{1\}$ of $\mathbb{Z}_{10}$ is clearly zero-4-sum-free over the integers, since $0 \not \in 4_{\pm}A$, but $A_1 \times \mathbb{Z}_{10}$ is not zero-4-sum-free over the integers in $\mathbb{Z}_{10}^2$, since (for example) 
$$(1,1)+(1,6)-(1,3)-(1,4)=(0,0).$$

We consider two special cases: when $H$ consists of a single positive integer $h$ and when $H$ consists of all positive integers up to some value $t$.  The cases when $H = \mathbb{N}_0$ or $H = \mathbb{N}$ are trivial as we then have $\tau_{\pm} (G, H)=0$ since $0 \in H_{\pm}A$ for any nonempty subset $A$ of $G$.

\subsection{Fixed number of terms} \label{5maxUSfixed}

In this section we investigate, for a given group $G$ and positive integer $h$, the quantity $$\tau_{\pm} (G,h) = \mathrm{max} \{ |A|  \mid A \subseteq G, 0 \not \in h_{\pm}A\},$$ that is, the maximum size of a zero-$h$-sum-free set over $\mathbb{Z}$. 

Since both for $h=1$ and for $h=2$, a subset of $G$ is zero-$h$-sum-free over $\mathbb{Z}$ if, and only if, it is zero-$h$-sum-free (over $\mathbb{N}_0$), Proposition \ref{z for h=1 and h=2} implies: 

\begin{prop} \label{zpm for h=1 and h=2}
We have $$\tau_{\pm}(G,1)=n-1$$ and $$\tau_{\pm}(G,2)=\frac{n-|\mathrm{Ord}(G,2)|-1}{2};$$ in particular, $$\tau_{\pm}(\mathbb{Z}_n,2)=\lfloor (n-1)/2 \rfloor.$$

\end{prop}

Now consider the case $h=3$.  For a subset $A$ of $G$, we have $0 \not \in 3_{\pm}A$ if, and only if, $A$ is both zero-3-sum-free and sum-free, that is, $0 \not \in 3A$ and $0 \not \in 2A-A$.  Let us consider $G=\mathbb{Z}_n$.  In Section \ref{5maxUfixed}, we showed that  
$$\tau(\mathbb{Z}_n,3) \geq v_3(n,3)=\left\{
\begin{array}{ll}
\left(1+\frac{1}{p}\right) \frac{n}{3} & \mbox{if $n$ has prime divisors congruent to 2 mod 3,} \\ & \mbox{and $p$ is the smallest such divisor,}\\ \\
\left\lfloor \frac{n-1}{3} \right\rfloor & \mbox{otherwise.}\\
\end{array}\right.$$
Here we show that $\tau_{\pm}(\mathbb{Z}_n,3) \geq v_3(n,3)$ as well.  As before, when $n$ has a divisor $p$ which is congruent to 2 mod 3, we take 
$$A=\left\{(p+1)/3+ i+pj \mbox{    } | \mbox{    } i=0,1,\dots,(p-2)/3, \; j=0,1,\dots,n/p-1 \right\}.$$  We have already seen that $A$ is zero-3-sum-free, it is also easy to see that it is sum-free as well.  Indeed, if $$k=(p+1)/3+(i_1+i_2-i_3) +p(j_1+j_2-j_3)=0$$ in $\mathbb{Z}_n$ for some $$i_1,i_2,i_3 \in \left\{0,1,\dots,(p-2)/3\right\}$$ and $$j_1,j_2,j_3 \in \left\{0,1,\dots,n/p-1 \right\},$$ then the integer $k$ is divisible by $n$ and thus by $p$.  Therefore, $(p+1)/3+(i_1+i_2-i_3)$ would have to be divisible by $p$, but this is not possible as 
$$1 \leq (p+1)/3+(i_1+i_2-i_3) \leq p-1.$$

Now suppose that $n$ has no such divisor---in this case, $n$ itself cannot be congruent to 2 mod 3.  Then $n$ is either divisible by 3, in which case the set 
$$\left\{ \frac{n}{3}  +1,  \frac{n}{3} +2, \dots, 2 \cdot \frac{n}{3} -1  \right\}$$ works, or it is congruent to 1 mod 3, in which case we can take $$\left\{ \frac{n-1}{3}  +1,  \frac{n-1}{3} +2, \dots, 2 \cdot \frac{n-1}{3}  \right\}.$$  Both of these sets have size $\left \lfloor \frac{n-1}{3} \right \rfloor$, completing the proof of the inequality above.

Therefore, we have the following.

\begin{prop} \label{zf}
For all positive integers $n$, we have $$\tau_{\pm}(\mathbb{Z}_n,3) \geq v_3(n,3)=\left\{
\begin{array}{ll}
\left(1+\frac{1}{p}\right) \frac{n}{3} & \mbox{if $n$ has prime divisors congruent to 2 mod 3,} \\ & \mbox{and $p$ is the smallest such divisor,}\\ \\
\left\lfloor \frac{n-1}{3} \right\rfloor & \mbox{otherwise.}\\
\end{array}\right.$$
\end{prop}

Our conjecture is that, in Proposition \ref{zf}, the lower bound gives the actual value:

\begin{conj} \label{zfconj}
For all positive integers $n$, we have
$$\tau_{\pm}(\mathbb{Z}_n,3) =v_3(n,3).$$
\end{conj}  
We pose the following problem.

\begin{prob}

Prove Conjecture \ref{zfconj}.

\end{prob} 

We should mention that Conjecture \ref{zfconj} obviously holds when $v_3(n,3)=v_1(n,3)$, as is the case, for example, if $n$ is even (more generally, has a prime divisor congruent to 2 mod 3) or is prime. 

Let us now turn to noncyclic groups.  As we noted on page \pageref{over Z not same}, Proposition \ref{zetafordirectsum} does not carry through to zero-sum-free sets over the integers.  

However, if $h \geq 3$ is odd, and for each odd integer $k \leq h$, $A_1$ is zero-$k$-sum-free over the integers in $G_1$, then we still see that $A_1 \times G_2$ is zero-$h$-sum-free over the integers in $G_1 \times G_2$: Indeed, if we were to have 
$$(a_1,g_1)+\cdots+(a_i,g_i)-(a_{i+1},g_{i+1})- \cdots - (a_h,g_h)=(0,0)$$ with $a_i \in A_1$ and $g_i \in G_2$ for some $i=1,2,\dots,h$ (ignoring subtractions when $i=h$),  then, looking at the first components and cancelling identical terms that are both added and subtracted (if any), we arrive at a signed sum of $k$ elements of $A_1$ for some odd integer $k \leq h$, a contradiction.  In fact, assuming that $A_1$ is zero-1-sum-free is not necessary as that is implied by it being zero-$h$-sum-free for $h$.  This results in the following:

\begin{prop} \label{h odd direct prod}
Suppose that $G_1$ and $G_2$ are finite abelian groups, $h \geq 3$ is odd, and that $A_1$ is zero-$k$-sum-free over the integers in $G_1$ for each $k \in \{3,5, \dots, h\}$.  Then $A_1 \times G_2$ is zero-$h$-sum-free over the integers in $G_1 \times G_2$.   

\end{prop}  
(A special case of this result was discovered by Matys\index{Matys, E.} in \cite{Mat:2014a}.)

As a consequence of Propositions \ref{cor noncyclic tau bounds},  \ref{zf}, and  \ref{h odd direct prod}, we get:

\begin{cor} \label{cor noncyclic tau Z bounds}
If $G$ is an abelian group of order $n$ and exponent $\kappa$, then
$$v_3(\kappa,3) \cdot \frac{n}{\kappa} \leq \tau_{\pm} (G,3) \leq v_1(n,3).$$
\end{cor}

We do not have a conjecture for $\tau_{\pm} (G,3)$ when $G$ is not cyclic.

\begin{prob}
Evaluate $\tau_{\pm} (G,3)$ for noncyclic groups $G$.

\end{prob}

Turning to the case of $h = 4$, we see a radical difference: as we are about to see, we cannot expect a formula for $\tau_{\pm}(\mathbb{Z}_n,4)$ that is a linear function of $n$.  Indeed, we have the following results.

\begin{prop} \label{zer2hisBh}

Let $A \subseteq G$, and suppose that $h \in \mathbb{N}$.  
If $A$ is a $B_h$ set over $\mathbb{Z}$, then $A$ is also a zero-$2h$-sum-free set over $\mathbb{Z}$, and therefore
$$\sigma_{\pm} (G, h) \leq \tau_{\pm}(G,2h).$$

\end{prop} 

\begin{prop} \label{Bhiszero4}

Let $A \subseteq G$, and suppose that $h$ is a positive integer that is divisible by 4.  
If $A$ is a zero-$h$-sum-free set over $\mathbb{Z}$, then $A$ is also a $B_2$ set over $\mathbb{Z}$, and therefore
$$\tau_{\pm}(G,h) \leq \sigma_{\pm} (G, 2).$$

\end{prop} 

For the proofs of Propositions \ref{zer2hisBh} and \ref{Bhiszero4}, see pages  \pageref{proofofzer2hisBh} and \pageref{proofofBhiszero4}, respectively. 

By these results:

\begin{cor}
A set is zero-$4$-sum-free set over $\mathbb{Z}$ in $G$ if, and only if, it is a $B_2$ set over $\mathbb{Z}$ in $G$.

\end{cor}

Therefore, by Proposition \ref{sidon tight}, we must have $$\tau_{\pm}(G,4) = \sigma_{\pm} (G, 2) \leq \sigma (G, 2) \leq \lfloor \sqrt{n} \rfloor+ 1.$$
We do not know $\tau_{\pm}(\mathbb{Z}_n,4)$.

\begin{prob}

Find the value of $\tau_{\pm}(\mathbb{Z}_n,4)$ for all $n$.

\end{prob} 

A bit more modestly:

\begin{prob} \label{explicit zero-4}

Find a positive constant $c$ and an explicit zero-4-sum-free set over $\mathbb{Z}$ in $\mathbb{Z}_n$ which has at least $c \cdot \sqrt{n}$ elements for every large enough $n$.

\end{prob}

The following result comes short of the demands of Problem \ref{explicit zero-4}: it provides a set of (asymptotically) smaller size, and the set is not given explicitly.

\begin{prop} \label{greedy zero-even}
For each $h$ there exists a positive constant $c_h$ for which $$\tau_{\pm}(\mathbb{Z}_n,h) >c_h \cdot n^{1/(h-1)}.$$  In particular, if $$n > \frac{4}{3}m^3 + \frac{8}{3}m,$$ then there is a zero-4-sum-free set over $\mathbb{Z}$ in $\mathbb{Z}_n$ that has size $m$.  

\end{prop} 
The proof of Proposition \ref{greedy zero-even} is on page \pageref{proof of greedy zero-even}.  (The result for $h=4$ was found by Phillips\index{Phillips, K.} in \cite{Phi:2011a}.)  

We also offer the following problem.

\begin{prob}

Find the value of $\tau_{\pm}(\mathbb{Z}_n,h)$ for even values of $h \geq 6$.

\end{prob} 

The situation is very different for odd values of $h$.  We have already discussed the case $h=3$ above; here we see what we can say about the general case when $h \geq 5$ is odd.  As Matys\index{Matys, E.} in \cite{Mat:2014a} pointed out, we cannot hope for a result similar to Proposition \ref{cor noncyclic tau Z bounds}: 
we see that $v_5(9,5)=3$, but $\tau_{\pm}(\mathbb{Z}_9,5)=2$.  (We can verify this last assertion, as follows.  Suppose that $A$ is a zero-5-sum-free set over the integers in $\mathbb{Z}_9$.  Clearly, $0 \not \in A$.  Also, if $a \in A$, then $4a \not \in A$; since the sequence $(1,4,16,64)$ becomes $(1,4,7,1)$ mod 9, at most one of 1, 4, or 7 can be in $A$. Similarly, at most one of 2, 8, or 5 can be in $A$.  Furthermore, neither 3 nor 6 can be in $A$ together with any of 1, 4, 7, 2, 8, or 5: for example, if $3 \in A$ and $a \in A$ for some $a \equiv 1$ mod 3, then 
$$3+3+a+a+a=3(a+2) \equiv 0 \; \mbox{mod} \; 9,$$ and if $3 \in A$ and $a \in A$ for some $a \equiv 2$ mod 3, then 
$$-3-3+a+a+a=3(a-2) \equiv 0 \; \mbox{mod} \; 9.$$  This proves that $\tau_{\pm}(\mathbb{Z}_9,5) \leq 2$; the set $\{4,5\}$ shows that equality holds.)  

Here is what we know about $\tau_{\pm}(G,h)$ for odd $h \geq 5$ when $G$ is cyclic.  As Matys\index{Matys, E.} in \cite{Mat:2014a} observed, when $n$ is even, then the odd elements of $\mathbb{Z}_n$ form a zero-$h$-sum-free set over the integers, hence $\tau_{\pm}(\mathbb{Z}_n,h) \geq n/2$ then.  Since $\tau(\mathbb{Z}_n,h)$, and thus $v_1(n,h)$ provides an upper bound and, by Corollary \ref{v function bounds}, $v_1(n,h)=n/2$ in this case, we get:

\begin{prop} [Matys; cf.~\cite{Mat:2014a}]\index{Matys, E.} 
If $n$ is even and $h \geq 3$ is odd, then $\tau_{\pm}(\mathbb{Z}_n,h) = n/2$.

\end{prop}

When $n$ and $h$ are both odd, we can find a linear lower bound for $\tau_{\pm}(\mathbb{Z}_n,h)$ as follows.  Let $A$ be the set of integers (viewed as elements of $G = \mathbb{Z}_n$) that are strictly between $\frac{h-1}{2}\frac{n}{h}$ and $\frac{h+1}{2}\frac{n}{h}$.  We can then prove that $0 \not \in h_{\pm}A$, and evaluating the size of $A$ results in the following proposition.

\begin{prop} [Matys; cf.~\cite{Mat:2014a}] \label{zetapmhodd}\index{Matys, E.} 
For all odd positive integers $n$ and $h$,
$$\tau_{\pm}(\mathbb{Z}_n,h) \geq 2  \left \lfloor \frac{n+h-2}{2h} \right \rfloor.$$
\end{prop}

The proof of Proposition \ref{zetapmhodd} (different from the one given by Matys\index{Matys, E.} in \cite{Mat:2014a}) is on page \pageref{proofofzetapmhodd}.

The following problem seems rather intriguing:

\begin{prob}

Find $\tau_{\pm}(\mathbb{Z}_n,h)$ for odd values of $h \geq 5$.

\end{prob}

As always, we are also interested in noncyclic groups.
\begin{prob}

Find the value of $\tau_{\pm}(G,h)$ for noncyclic groups and for $h \geq 4$.

\end{prob}

\subsection{Limited number of terms} \label{5maxUSlimited}

A subset $A$ of $G$ for which $$0 \not \in [1,t]_{\pm}A =\cup_{h=1}^t h_{\pm}A$$ for some positive integer $t$ is called a {\em $t$-independent set} in $G$.  Here we investigate $$\tau_{\pm} (G,[1,t]) = \mathrm{max} \{ |A|  \mid A \subseteq G, 0 \not \in [1,t]_{\pm}A\},$$ that is, the maximum size of a $t$-independent set in $G$.  

As we pointed out above, with $\kappa$ denoting the exponent of the group, adding any element of the group to itself $\kappa$ times results in $0$.  Therefore, if $t \geq \kappa$, then no element of $G$ can belong to a $t$-independent set and we have $\tau_{\pm} (G,[1,t]) =0.$  However, if $t < \kappa$, then $\tau_{\pm} (G,[1,t]) \geq 1$: at the least the one-element set $\{a\}$, where $a$ is any element of order $\kappa$, will be $t$-independent.

The cases of $t=1$ and $t=2$ are easy to handle.  Clearly, a set $A$ is 1-independent if, and only if, $0 \not \in A$.  Regarding $t=2$, first note that a 2-independent set cannot contain any element of $\{0\} \cup \mathrm{Ord}(G,2)$ (the elements of order at most 2); to get a maximum 2-independent set in $G$, take exactly one of each element or its negative in $G \setminus \mathrm{Ord}(G,2)\setminus \{0\}$.
In summary:

\begin{prop}  \label{tau pm t=1,2}
For all groups $G$ we have $$\tau_{\pm} (G,[1,1])=n-1$$ and $$\tau_{\pm} (G,[1,2])=\frac{n-|\mathrm{Ord}(G,2)|-1}{2};$$  
in particular, 
 $$\tau_{\pm} (\mathbb{Z}_n,[1,2])=\lfloor (n-1)/2 \rfloor.$$

\end{prop}

Let us now consider $t=3$.  Before we proceed, it is helpful to state that a subset $A$ is 3-independent in $G$ if, and only if, none of the equations
$$x=0, \; x+y=0, \; x+y+z=0, \mbox{and} \; x+y=z$$
have a solution in $A$.

We first find explicit 3-independent sets in the cyclic group $\mathbb{Z}_n$, as follows.  For every $n$, the odd integers that are less than $n/3$ form a 3-independent set; that is, the set $$\left\{2i+1 \mbox{    } | \mbox{    } i = 0, 1, \dots, \left\lfloor n/6 \right\rfloor -1 \right\}$$ is 3-independent in $\mathbb{Z}_n$.  Furthermore, if $n$ is even, we can go up to (but not including) $n/2$ as then the sum of two odd integers cannot equal $n$; so, when $n$ is even, the set $$\left\{2i+1 \mbox{    } | \mbox{    } i = 0, 1, \dots, \left\lfloor n/4 \right\rfloor -1 \right\}$$ is 3-independent in $\mathbb{Z}_n$.  

We can do better in one special case when $n$ is odd; namely, when $n$ has a prime divisor $p$ that is congruent to 5 mod 6, one can show that the set 
$$\left\{ pi_1+2i_2+1 \mbox{    } | \mbox{    } i_1=0,1,\dots,n/p-1, \mbox{  } i_2=0,1,\dots, (p-5)/6 \right\}$$
is 3-independent.  For example, when $n=25$, we may take $p=5$, with which we get $$\{ 5i_1+2i_2+1 \mbox{    } | \mbox{    } i_1=0,1,\dots,4, \mbox{  } i_2=0\} = \{1, 6, 11, 16, 21\}.$$ This set is 3-independent in $G=\mathbb{Z}_{25}$.  Note that, when determining the independence number of a subset, an element and its negative play the same role, thus the set $\{1, 6, 11, 16, 21\}$ above is essentially the same as the set $\{1,4,6,9,11\}$ of the example on page \pageref{Dec Indep}.

In summary, we have
$$\tau_{\pm}(\mathbb{Z}_n, [1,3]) \geq \left\{
\begin{array}{ll}
\left\lfloor \frac{n}{4} \right\rfloor & \mbox{if $n$ is even,}\\ \\
\left(1+\frac{1}{p}\right) \frac{n}{6} & \mbox{if $n$ is odd, has prime divisors congruent to 5 mod 6,} \\ & \mbox{and $p$ is the smallest such divisor,}\\ \\
\left\lfloor \frac{n}{6} \right\rfloor & \mbox{otherwise.}\\
\end{array}\right.$$

Let us now turn to noncyclic groups.  First, an observation.  Suppose that $G$ is of the form $G_1 \times G_2$.  (Note that every group can be written in the form $G_1 \times G_2$ with $G_2=\mathbb{Z}_{\kappa}$ where $\kappa$ is the exponent of $G$.)  It is not hard to see that, if $A_2 \subseteq G_2$ is 3-independent in $G_2$, then $$A=\{(g,a)  \mid g \in G_1, a \in A_2\}$$ is 3-independent in $G$, so we have the following:    

\begin{prop} \label{zeta3fordirectsum}
Let $G_1$ and $G_2$ be finite abelian groups, $G=G_1 \times G_2$, and suppose that $A_2\subseteq G_2$ is 3-independent in $G_2$.  Then $$A=\{(g,a)  \mid g \in G_1, a \in A_2\}$$ is 3-independent in $G$; in particular,
$$\tau_{\pm}(G_1 \times G_2,[1,3]) \geq |G_1| \cdot \tau_{\pm}(G_2,[1,3]).$$

\end{prop}
(This proposition does not hold for $t$-independent sets for $t \geq 4$.)

Combining Proposition \ref{zeta3fordirectsum} with our results for cyclic groups above, we see that, for a group of exponent $\kappa$ we get:
$$\tau_{\pm}(G, [1,3]) \geq \left\{
\begin{array}{ll}
\frac{n}{4}  & \mbox{if $\kappa$ is divisible by 4,}\\ \\
\frac{n}{\kappa} \cdot \frac{\kappa-2}{4} & \mbox{if $\kappa$ is even but not divisible by 4,}\\ \\
\left(1+\frac{1}{p}\right) \frac{n}{6} & \mbox{if $\kappa$ is odd, has prime divisors congruent to 5 mod 6,} \\ & \mbox{and $p$ is the smallest such divisor,}\\ \\
\left\lfloor \frac{\kappa}{6} \right\rfloor \cdot \frac{n}{\kappa} & \mbox{otherwise.}\\
\end{array}\right.$$

Next, we show that we can do slightly better when $\kappa$ is even but not divisible by 4.  Let $G=G_1 \times \mathbb{Z}_{\kappa}$.  Suppose that $A_1$ is a 2-independent set in $G_1$, and $A_2$ is the 3-independent set $$A_2=\left\{2i+1  \mid i=0,1,\dots,(\kappa -6)/4 \right\}$$ in $\mathbb{Z}_{\kappa}$.  We can then verify that the set
$$A=\{(g,a)  \mid g \in G_1, a \in A_2\} \cup \{(a, \kappa/2)  \mid a \in A_1\}$$
is 3-independent in $G$.  (When checking that the equations $x+y+z=0$ and $x+y=z$ have no solutions in $A$, note that $\kappa$ is even, but $\kappa/2$ and all elements of $A_2$ are odd.)

Suppose further that $A_1$ is of maximum size; that is, $$|A_1|=\frac{|G_1|-|\mathrm{Ord}(G_1,2)|-1}{2};$$   since
$$|\mathrm{Ord}(G_1,2)|=\frac{|\mathrm{Ord}(G,2)|-1}{2},$$
we have
$$|A_1|=\frac{n}{2\kappa}-\frac{|\mathrm{Ord}(G,2)|+1}{4},$$ 
with which $$|A|=\frac{\kappa -2}{4} \cdot \frac{n}{\kappa}+\frac{n}{2\kappa}-\frac{|\mathrm{Ord}(G,2)|+1}{4} = \frac{n-|\mathrm{Ord}(G,2)|-1}{4}.$$

Therefore, we have
$$\tau_{\pm}(G, [1,3]) \geq \left\{
\begin{array}{ll}
\frac{n}{4}  & \mbox{if $\kappa$ is divisible by 4,}\\ \\
\frac{n-|\mathrm{Ord}(G,2)|-1}{4} & \mbox{if $\kappa$ is even but not divisible by 4,}\\ \\
\left(1+\frac{1}{p}\right) \frac{n}{6} & \mbox{if $\kappa$ is odd, has prime divisors congruent to 5 mod 6,} \\ & \mbox{and $p$ is the smallest such divisor,}\\ \\
\left\lfloor \frac{\kappa}{6} \right\rfloor \cdot \frac{n}{\kappa} & \mbox{otherwise.}\\
\end{array}\right.$$

In \cite{BajRuz:2003a}, Bajnok and Ruzsa\index{Ruzsa, I.}\index{Bajnok, B.} showed that, in the first three cases, equality holds; for the last case, they only proved that $\tau_{\pm}(G, [1,3])$ cannot be more than $\left\lfloor \frac{n}{6} \right\rfloor.$  Namely, we have the following result.

\begin{thm}[Bajnok and Ruzsa; cf.~\cite{BajRuz:2003a}] \label{3freeG}\index{Ruzsa, I.}\index{Bajnok, B.}
As usual, let $\kappa$ be the exponent of $G$.  We have
$$\tau_{\pm}(G, [1,3]) = \left\{
\begin{array}{ll}
\frac{n}{4}  & \mbox{if $\kappa$ is divisible by 4,}\\ \\
\frac{n-|\mathrm{Ord}(G,2)|-1}{4} & \mbox{if $\kappa$ is even but not divisible by 4,}\\ \\
\left(1+\frac{1}{p}\right) \frac{n}{6} & \mbox{if $\kappa$ is odd, has prime divisors congruent to 5 mod 6,} \\ & \mbox{and $p$ is the smallest such divisor;}\\
\end{array}\right.$$ 
furthermore, if $\kappa$ (iff $n$) is odd and has no prime divisors congruent to 5 mod 6, then $$\left\lfloor \frac{\kappa}{6} \right\rfloor \frac{n}{\kappa} \leq \tau_{\pm}(G, [1,3]) \leq \frac{n}{6}.$$ 

\end{thm}

As a consequence, we see that Theorem \ref{3freeG} has settled the problem of finding the maximum size of a 3-independent set in cyclic groups: 

\begin{thm}[Bajnok and Ruzsa; cf.~\cite{BajRuz:2003a}] \label{3free}\index{Ruzsa, I.}\index{Bajnok, B.}

For the cyclic group $G=\mathbb{Z}_n$ we have
$$\tau_{\pm}(\mathbb{Z}_n, [1,3]) = \left\{
\begin{array}{ll}
\left\lfloor \frac{n}{4} \right\rfloor & \mbox{if $n$ is even,}\\ \\
\left(1+\frac{1}{p}\right) \frac{n}{6} & \mbox{if $n$ is odd, has prime divisors congruent to 5 mod 6,} \\ & \mbox{and $p$ is the smallest such divisor,}\\ \\
\left\lfloor \frac{n}{6} \right\rfloor & \mbox{otherwise.}\\
\end{array}\right.$$

\end{thm}

As a corollary to Theorem \ref{3free}, we get the following.

$\tau_{\pm}(\mathbb{Z}_n, [1,3]) = \left\{
\begin{array}{ll}
1 & \mbox{if $n={\bf 4}, 5, 6, 7, 9;$}\\
2 & \mbox{if $n={\bf 8}, 10, 11, 13;$}\\
3 & \mbox{if $n={\bf 12}, 14, 15, 17, 19, 21;$}\\
4 & \mbox{if $n={\bf 16}, 18, 23, 27;$}\\
5 & \mbox{if $n={\bf 20}, 22, 25, 29, 31;$}\\
6 & \mbox{if $n={\bf 24}, 26, 33, 37, 39.$}\\
\end{array}\right.$

(Entries in bold-face mark tight 3-independent sets, as we will explain shortly.)

Since the publication of \cite{BajRuz:2003a}, Green and Ruzsa\index{Ruzsa, I.}\index{Green, B.} succeeded in determining the maximum size of a sum-free set in any $G$: 
$$\mu(G, \{1,2\})= v_{1}(\kappa,3) \cdot \frac{n}{\kappa}$$
(cf.~Theorem \ref{(2,1)all}).  We can use this result to determine $\tau_{\pm}(G, [1,3])$ when $\kappa$ is only divisible by primes congruent to 1 mod 6.  Indeed, since a 3-independent set $A$ must be asymmetric: $A \cap -A$ must be empty.  Furthermore, the set $A \cup -A$ has to be sum-free, so we get
$$2|A|=|A \cup -A| \leq v_{1}(\kappa,3) \cdot \frac{n}{\kappa}= \left \lfloor \frac{\kappa }{3} \right \rfloor \cdot \frac{n}{\kappa}= \frac{\kappa -1}{3}\cdot \frac{n}{\kappa} ,$$ and thus $$|A| \leq \frac{\kappa -1}{6}\cdot \frac{n}{\kappa} $$ (note that $(\kappa-1)/3$ is an even integer).  Since in our case $$\left\lfloor \frac{\kappa}{6} \right\rfloor \frac{n}{\kappa}=\frac{\kappa -1}{6}\cdot \frac{n}{\kappa}$$ as well, from the last case of Theorem \ref{3freeG} we get the following:

\begin{cor}
Suppose that $G$ has exponent $\kappa$ and that every prime divisor of $\kappa$ is congruent to 1 mod 6.  Then  
$$\tau_{\pm}(G, [1,3])=\frac{\kappa -1}{6}\cdot \frac{n}{\kappa}.$$
\end{cor}

This leaves only one case open:

\begin{prob} \label{q(G,3)}

Suppose that $G$ is not cyclic and that its exponent is the product of a positive integer power of 3 and perhaps some primes that are congruent to 1 mod 6.  Find $\tau_{\pm}(G, [1,3])$.
\end{prob}

For example, by Theorem \ref{3freeG} we see that $\tau_{\pm}(\mathbb{Z}_9^2, [1,3])$ is at least 9 and at most 13; according to Laza\index{Laza, N.}  (see \cite{Laz:?}), we have $\tau_{\pm}(\mathbb{Z}_9^2, [1,3]) \geq 11$ as the set
$$A=\{(0,1), (1,0), (1,3), (1,6), (2,1), (2,4), (2,7), (4,1), (4,4), (4,7), (6,1)\}$$
is 3-independent in $\mathbb{Z}_9^2$.  (It is worth noting that $[1,3]_{\pm} A$ contains all nonzero elements of the group.)

For $t \geq 4$, exact results seem more difficult.  
With the help of a computer, Laza\index{Laza, N.}  (see \cite{Laz:?}) generated the following values:

$\tau_{\pm}(\mathbb{Z}_n, [1,4]) = \left\{
\begin{array}{ll}
1 & \mbox{if $n={\bf 5}, 6, \dots, 12;$}\\
2 & \mbox{if $n={\bf 13}, 14, \dots, 26;$}\\
3 & \mbox{if $n=27, 28, \dots, 45$, and $n=47;$}\\
4 & \mbox{if $n=46$, $n=48, 49, \dots, 68$, and $n=72, 73;$}\\
5 & \mbox{if $n=69, 70, 71$, and $n=74, 75, \dots, 102;$}\\
\end{array}\right.$

$\tau_{\pm}(\mathbb{Z}_n, [1,5]) = \left\{
\begin{array}{ll}
1 & \mbox{if $n={\bf 6}, 7, \dots, 17$, and $n=19, 20$;}\\
2 & \mbox{if $n={\bf 18}$, $n=21, 22, \dots, 37$, $n=39, 40, 41$, $n=43, 44, 45, 47$;}\\
3 & \mbox{if $n={\bf 38}, 42, 46$, $n=48, 49, \dots, 69$}, \\
 &                         $n=71, 72, 73, 75, 76, 77, 79, 81, 83, 85, 87;$\\
\end{array}\right.$

$\tau_{\pm}(\mathbb{Z}_n, [1,6]) = \left\{
\begin{array}{ll}
1 & \mbox{if $n={\bf 7}, 8, 9, \dots, 24$;}\\
2 & \mbox{if $n={\bf 25}, 26, 27, \dots, 69$;}\\
3 & \mbox{if $n=70, 71, \dots, 151$, and $n=153, 154, 155, 158, 159, 160.$}\\
\end{array}\right.$

(Values marked in bold-face will be discussed shortly.)

The following two problems seem difficult, but even partial answers (those for certain special values of $n$) would be very interesting:

\begin{prob} 

Find $\tau_{\pm}(\mathbb{Z}_n, [1,4])$ (at least for some infinite family of $n$ values).
\end{prob}

\begin{prob} 

Find $\tau_{\pm}(\mathbb{Z}_n, [1,5])$ (at least for some infinite family of $n$ values).
\end{prob}

It is also interesting to investigate $t$-independent sets from the opposite viewpoint: given a nonnegative integer $m$ and positive integer $t$, what are the possible groups $G$ for which $\tau_{\pm} (G,[1,t]) =m$?

The answer for $m=0$ is clear and has already been discussed.  For $m=1$ we can see from Laza's\index{Laza, N.} work \cite{Laz:?} that, in the cyclic group $\mathbb{Z}_n$, we have
 
$\tau_{\pm}(\mathbb{Z}_n, [1,t]) =1 \Leftrightarrow \left\{
\begin{array}{ll}
n=3-4 & \mbox{if $t=2$;}\\
n=4-7, 9 & \mbox{if $t=3$;}\\
n=5-12 & \mbox{if $t=4$;}\\
n=6-17, 19-20 & \mbox{if $t=5$;}\\
n=7-24 & \mbox{if $t=6$;}\\
n=8-31, 33, 35 & \mbox{if $t=7$;}\\
n=9-40 & \mbox{if $t=8$;}\\
n=10-49, 51-53 & \mbox{if $t=9$;}\\
n=11-60 & \mbox{if $t=10$.}\\
\end{array}\right.$

We can prove the following:

\begin{prop}  \label{q} As usual, let $G$ be an abelian group of order $n$ and exponent $\kappa$.  Let $t \geq 2$ be an integer. 

\begin{enumerate}

\item  We have  $\tau_{\pm} (G,[1,t]) =0$ if, and only if, $\kappa \leq t$.

\item

  Set $b_t=\left \lfloor t^2/2 \right \rfloor +t$.
\begin{enumerate}
\item Suppose that $t$ is even.  Then $\tau_{\pm} (\mathbb{Z}_n,[1,t]) =1$ if, and only if, $t+1 \leq n \leq b_t.$  In particular, the set $\{1\}$ is $t$-independent in $\mathbb{Z}_n$ for $n \geq t+1$, and the set $\{t/2, t/2 +1\}$ is $t$-independent in $\mathbb{Z}_n$ for $n \geq b_t+1$.

\item  Suppose that $t$ is odd.

\begin{enumerate}
\item If $t+1 \leq n \leq b_t,$ then $\tau_{\pm} (\mathbb{Z}_n,[1,t]) =1$; in particular, the set $\{1\}$ is $t$-independent in $\mathbb{Z}_n$.

\item  {\bf (Miller; cf.~\cite{Mil:2013a})}\index{Miller, Z.} If $\tau_{\pm} (\mathbb{Z}_n,[1,t]) =1$, then $t+1 \leq n \leq b_t+(t+1)/2;$ 
in particular, the set $\{(t+1)/2,(t+3)/2\}$ is $t$-independent in $\mathbb{Z}_n$  for $n \geq b_t+(t+3)/2$.
    
\item If $n  = b_t+1$, then $\tau_{\pm} (\mathbb{Z}_n,[1,t]) =2$; in particular, the set $\{1,t\}$ is $t$-independent in $\mathbb{Z}_n$.
\item {\bf (Miller; cf.~\cite{Mil:2013a})}\index{Miller, Z.} If $t$ is congruent to 3 mod 4, $n$ is even, and $n \geq b_t+1,$ then $\tau_{\pm} (\mathbb{Z}_n,[1,t]) \geq 2$; in particular, the set $\{(t-1)/2,(t+3)/2\}$ is $t$-independent in $\mathbb{Z}_n$.

\item {\bf (Miller; cf.~\cite{Mil:2013a})}\index{Miller, Z.} If $t$ is congruent to 1 mod 4, $n$ is congruent to 2 mod 4, and $n \geq b_t+1,$ then $\tau_{\pm} (\mathbb{Z}_n,[1,t]) \geq 2$; in particular, the set $\{(t-3)/2,(t+5)/2\}$ is $t$-independent in $\mathbb{Z}_n$.
\end{enumerate}

\end{enumerate}
\end{enumerate}
\end{prop}

Some of these statements we have already seen: the set $\{1\}$ is clearly $t$-independent in $\mathbb{Z}_n$ for $n \geq t+1$, and 
Proposition \ref{indepbound} implies that
\begin{itemize}
\item if there is a $t$-independent set of size 1 in $\mathbb{Z}_n$, then $n \geq t+1$, and
\item if there is a $t$-independent set of size 2 in $\mathbb{Z}_n$, then $n \geq \left \lfloor t^2/2 \right \rfloor +t+1$.
\end{itemize}
To complete the proof, we need to verify that the given two-element sets are indeed $t$-independent in their corresponding groups---for that, see page \pageref{proofofq}.  We can observe that, when $t \equiv 1$ mod 4, then $b_t+1 \equiv 2$ mod 4, so statements (iv) and (v) together imply that if $n  = b_t+1$, then $\tau_{\pm} (\mathbb{Z}_n,[1,t]) =2$---cf.~statement (iii).

We should point out that Proposition \ref{q} does not completely determine all $n$ for which $\tau_{\pm} (\mathbb{Z}_n,[1,t]) =1$, although for $t \leq 10$ it rules out all $n$ values other than those listed above.  For example, for $t=9$, $n \leq 9$ and $n \geq 55$ are ruled out by statement 2 (b) (ii), $n=50$ is ruled out by both 2 (b) (iii) and (v), and $n=54$ is ruled out by 2 (b) (v).   

We offer the following open problems.

\begin{prob} \label{m=1,todd}

Examine the remaining cases of statement 2 (b) of Proposition \ref{q} to determine, for each odd value of $t$, all values of $n$ for which $\tau_{\pm} (\mathbb{Z}_n,[1,t]) =1$.
\end{prob}

\begin{prob} \label{m=2,teven}

For each even value of $t$, find all values of $n$ for which $\tau_{\pm} (\mathbb{Z}_n,[1,t]) =2$.
\end{prob}

\begin{prob} \label{m=2,todd}

For each odd value of $t$, find all values of $n$ for which $\tau_{\pm} (\mathbb{Z}_n,[1,t]) =2$.
\end{prob}

The general problem for cyclic groups is as follows:

\begin{prob} \label{mgeneral,tgeneral}

For each value of $t$ and each $m \geq 3$, find all values of $n$ for which $\tau_{\pm} (\mathbb{Z}_n,[1,t]) =m$.
\end{prob}

For noncyclic groups, we offer

\begin{prob} \label{m=2,G}

For each value of $t$, find all noncyclic groups for which $\tau_{\pm} (G,[1,t]) =1$.
\end{prob}

With exact values few and far between, we are also interested in some good bounds for the value of $\tau_{\pm} (G,[1,t])$.  We can derive an upper bound for $\tau_{\pm} (G,[1,t])$, as follows.  Since the case $t=1$ is trivial, we assume that $t \geq 2$.

When $t$ is even, we already pointed out in Section \ref{4maxUSlimited} that $A=\{a_1,\dots,a_m\}$ (with $|A|=m$) being $t$-independent is equivalent to $$|[0,\tfrac{t}{2}]_{\pm}A|=|\mathbb{Z}^m([0,\tfrac{t}{2}])|$$ where 
$$\mathbb{Z}^m([0,\tfrac{t}{2}]))=\{(\lambda_1,\dots,\lambda_m) \in \mathbb{Z}^m  \mid |\lambda_1|+\cdots+|\lambda_m| \leq \tfrac{t}{2}\}.$$
Therefore, for even values of $t$, if a set $A$ of size $m$ is $t$-independent in $G$, then 
$$|\mathbb{Z}^m([0,\tfrac{t}{2}])| \leq n$$ must hold.
Recalling from Section \ref{0.2.4} that
$$|\mathbb{Z}^m([0,\tfrac{t}{2}])| =a(m, t/2)= \sum_{i\geq 0} {m \choose i} {t/2 \choose i} 2^i,$$
we thus must have
$$a(m, t/2)= \sum_{i\geq 0} {m \choose i} {t/2 \choose i} 2^i \leq n.$$

If $t$ is odd, we similarly must have $$|[0,\tfrac{t-1}{2}]_{\pm}A|=|\mathbb{Z}^m([0,\tfrac{t-1}{2}])|;$$ in addition, signed sums
$$\lambda_1a_1+\cdots+\lambda_ma_m$$ corresponding to the index set 
$$\mathbb{Z}^m(\tfrac{t+1}{2})_{1+}=\{(\lambda_1,\dots,\lambda_m) \in \mathbb{Z}^m  \mid \lambda_1 \geq 1, \lambda_1+|\lambda_2|+\cdots+|\lambda_m| = \tfrac{t+1}{2}\}$$
must be pairwise distinct and distinct from $[0,\tfrac{t}{2}]_{\pm}A$ as well.  Therefore, 
for odd values of $t$, if a set $A$ of size $m$ is $t$-independent in $G$, then 
$$|\mathbb{Z}^m([0,\tfrac{t-1}{2}])|+ |\mathbb{Z}^m(\tfrac{t+1}{2})_{1+}| \leq n$$ must hold.

Of course, there is nothing special about using $\lambda_1$ in the argument above---we could have replaced $\lambda_1$ by any one of the indices.  It is important to note that, while our condition above for the case when $t$ is even is equivalent for $A$ to be $t$-independent, for odd $t$ we only have a necessary condition: we could, in theory, still have, say, $\tfrac{t+1}{2} \cdot a_2$ equal one of the elements of $[0,\tfrac{t-1}{2}]_{\pm}A$, preventing $A$ from being $t$-independent. 

Recalling from Section \ref{0.2.4} that
$$|\mathbb{Z}^m(\tfrac{t+1}{2})_{1+}| =a(m-1, (t-1)/2)= \sum_{i\geq 0} {m-1 \choose i} {(t-1)/2 \choose i} 2^i,$$
we must have
$$a(m,(t-1)/2) +a(m-1, (t-1)/2)= \sum_{i\geq 0} \left[ {m \choose i} + {m-1 \choose i} \right] {(t-1)/2 \choose i} 2^i \leq n.$$
Note also that by Proposition \ref{functionsac}, we can rewrite the left-hand side as
$$a(m,(t-1)/2) +a(m-1, (t-1)/2)=c(m,(t+1)/2)=\sum_{i\geq 0} {m-1 \choose i-1} {(t+1)/2 \choose i} 2^i.$$

In summary, we have the following:

\begin{prop} \label{indepbound}
Suppose that $t \geq 2$ and that $A$ is a $t$-independent set of size $m$ in a group $G$ of order $n$.  Then
$$n \geq \left\{
\begin{array}{cl}
a(m,t/2) = \sum_{i\geq 0} {m \choose i} {t/2 \choose i} 2^i & \mbox{if} \; $t$\;  \mbox{is even}, \\ \\
c(m,(t+1)/2)= \sum_{i\geq 0} {m-1 \choose i-1} {(t+1)/2 \choose i} 2^i & \mbox{if} \; $t$\;  \mbox{is odd}.
\end{array}
\right.$$
\end{prop}

Proposition \ref{indepbound} gives a lower bound for $\tau_{\pm}(G,[1,t])$.  For example, we have
\begin{eqnarray*}
\tau_{\pm}(G,[1,2])&\leq& \left\lfloor \frac{n-1}{2} \right\rfloor, \\ \\
\tau_{\pm}(G,[1,3]) &\leq &\left\lfloor \frac{n}{4} \right\rfloor, \\ \\
\tau_{\pm}(G,[1,4]) &\leq &\left\lfloor \frac{\sqrt{2n-1}-1}{2} \right\rfloor, \\ \\
\tau_{\pm}(G,[1,5])& \leq &\left\lfloor \frac{\sqrt{n-2}}{2} \right\rfloor;
\end{eqnarray*}
the bound is less explicit for $t \geq 6$, but we certainly have
$$\tau_{\pm}(G,[1,t]) \leq \left\lfloor \frac{1}{2} \sqrt[s]{s! \cdot  n}+s \right\rfloor,$$ where $s=\lfloor t/2 \rfloor$.  In particular:

\begin{cor} \label{upper for ind}
For every $t \geq 2$ there is a positive constant $C_t$ so that 
$$\tau_{\pm}(G,[1,t]) \leq  C_t \cdot  \sqrt[s]{n},$$ where $s=\lfloor t/2 \rfloor$.

\end{cor}

For a lower bound, Bajnok and Ruzsa proved the following:\index{Ruzsa, I.}\index{Bajnok, B.}

\begin{thm} [Bajnok and Ruzsa; cf.~\cite{BajRuz:2003a}] \label{BajRuz lower}\index{Ruzsa, I.}\index{Bajnok, B.}
For every $t \geq 2$ there is a positive constant $c_t$ so that 
$$\tau_{\pm}(\mathbb{Z}_n,[1,t]) \geq  c_t \cdot  \sqrt[s]{n},$$ where $s=\lfloor t/2 \rfloor$.

\end{thm}

Letting $$c \cdot g_1(n) \lesssim f(n) \lesssim C \cdot  g_2(n)$$
mean that for every real number $\epsilon >0$, if $n$ is large enough (depending on $\epsilon$), then  
$$(c- \epsilon) \cdot g_1(n) \leq f(n) \le (C+\epsilon) \cdot  g_2(n),$$
we can combine Corollary \ref{upper for ind} and Theorem \ref{BajRuz lower} to write
$$c_t \cdot \sqrt[s]{n} \lesssim \tau_{\pm}(\mathbb{Z}_n,[1,t]) \lesssim C_t \cdot  \sqrt[s]{n}.$$

What can we say about the values of $c_t$ and $C_t$?  From Proposition \ref{tau pm t=1,2} we see that we may take $c_2=C_2=1/2$, and therefore
$$\lim \frac{\tau_{\pm}(\mathbb{Z}_n,[1,2])}{n} =1/2.$$  From Theorem \ref{3free}, we can have
$c_3=1/6$ and $C_3 =1/4$; we also see that both of these are tight as $\tau_{\pm}(\mathbb{Z}_n,[1,3])$ approaches both $n/6$ and $n/4$ infinitely often and thus
$$\lim \frac{\tau_{\pm}(\mathbb{Z}_n,[1,3])}{n}$$
does not exist.  

For $t=4$ and $t=5$ we have:
\begin{thm} [Bajnok and Ruzsa; cf.~\cite{BajRuz:2003a}] \label{Baj and Ruz t=4,5}\index{Ruzsa, I.}\index{Bajnok, B.}
We have
$$1/\sqrt{8} \cdot \sqrt{n} \lesssim \tau_{\pm}(\mathbb{Z}_n,[1,4]) \lesssim 1/\sqrt{2} \cdot \sqrt{n}$$
and
$$1/\sqrt{15} \cdot \sqrt{n} \lesssim \tau_{\pm}(\mathbb{Z}_n,[1,5]) \lesssim 1/\sqrt{2} \cdot \sqrt{n}.$$

\end{thm}

The following problems seem intriguing:

\begin{prob} \label{t=4}

Find (if possible) a value higher than $1/\sqrt{8}$ for $c_4$ and one lower than $1/\sqrt{2}$ for $C_4$.
\end{prob}

\begin{prob} \label{t=5}
Find (if possible) a value higher than $1/\sqrt{15}$ for $c_5$ and one lower than $1/\sqrt{2}$ for $C_5$.
\end{prob}

In \cite{BajRuz:2003a} we find the following conjectures:

\begin{conj} [Bajnok and Ruzsa; cf.~\cite{BajRuz:2003a}]\index{Ruzsa, I.}\index{Bajnok, B.}
We have $$\lim \frac{\tau_{\pm}(\mathbb{Z}_n,[1,4])}{\sqrt{n}} = 1/\sqrt{3},$$ but $$\lim \frac{\tau_{\pm}(\mathbb{Z}_n,[1,5])}{\sqrt{n}}$$ does not exist.

\end{conj}
Note that this conjecture, if true, would imply that one can take $c_4=C_4$ but not $c_5=C_5$. 

While we see that for a given $t \geq 3$, $\tau_{\pm}(\mathbb{Z}_n, [1,t])$ is not a monotone function of $n$, we may find that the subsequence of even or odd $n$ values is monotonic:

\begin{prob}
Decide whether for each even value of $t$ (but at least for $t=4$), the sequence $\tau_{\pm}(\mathbb{Z}_n, [1,t])$ is monotone for odd values of $n$. 

\end{prob}

\begin{prob}
Decide whether for each odd value of $t$ (but at least for $t=5$), the sequence $\tau_{\pm}(\mathbb{Z}_n, [1,t])$ is monotone for even values of $n$. 

\end{prob}

We are also interested in 4-independent and 5-independent sets in noncyclic groups:

\begin{prob} \label{r=2,t=4}

Determine the value of $\tau_{\pm}(G,[1,4])$ for noncyclic groups $G$.
\end{prob}

\begin{prob} \label{r=2,t=5}

Determine the value of $\tau_{\pm}(G,[1,5])$ for noncyclic groups $G$.
\end{prob}

We close this section by a particularly intriguing question: when can the size of $t$-independent sets achieve the maximum possible value allowed by the upper bound of Proposition \ref{indepbound}?  (We mention that for $t=1$ the answer is trivial: the set $G \setminus \{0\}$ is the maximum size 1-independent set in $G$.)  For $t \geq 2$, we are thus interested in classifying all {\em tight $t$-independent sets}, that is, $t$-independent sets of size $m$ with $n=b(m,t)$ where
$$b(m,t)= \left\{
\begin{array}{cl}
a(m, t/2)= \sum_{i\geq 0} {m \choose i} {t/2 \choose i} 2^i & \mbox{if} \; $t$\;  \mbox{is even}, \\ \\
c(m, (t+1)/2)=\sum_{i\geq 0}  {m-1 \choose i-1}  {(t+1)/2 \choose i} 2^i & \mbox{if} \; $t$\;  \mbox{is odd}.
\end{array}
\right.$$

The values of $b(m,t)$ can be tabulated for small values of $m$ and $t$:

$$\begin{array}{c||c|c|c|c|c|c|c|c|}
$b(m,t)$   &  $t=2$  & $t=3$  & $t=4$  & $t=5$  & $t=6$ & $t=7$ & $t=8$  \\ \hline \hline
$m=1$ &  {\bf 3}& {\bf 4}& {\bf 5}& {\bf 6}& {\bf 7}& {\bf 8}& {\bf 9} \\ \hline
$m=2$ &  {\bf 5}& {\bf 8}& {\bf 13}& {\bf 18}& {\bf 25}& {\bf 32}& {\bf 41} \\ \hline
$m=3$ &  {\bf 7}& {\bf 12}& 25& {\bf 38}& 63& 88& 129 \\ \hline
$m=4$ &  {\bf 9}& {\bf 16}& 41& 66& 129& 192& 321 \\ \hline
$m=5$ &  {\bf 11}& {\bf 20}& 61& 102& 231& 360& 681 \\ \hline
$m=6$ &  {\bf 13}& {\bf 24}& 85& 146& 377& 608& 1289 \\ \hline
\end{array}$$

Cases where there exists a group of size $b(m,t)$ with a known tight $t$-independent set are marked with bold-face.  The following proposition exhibits tight $t$-independent sets for all known parameters.

\begin{prop}  \label{tightex}

Let $m$ and $t$ be positive integers, $t \geq 2$, and let $G$ be an abelian group of order $n$ and exponent $\kappa$.

\begin{enumerate}

\item If $n=2m+1$, then $G \setminus \{0\}$ can be partitioned into parts $K$ and $-K$, and both $K$ and $-K$ are tight 2-independent sets in $G$.  For example, the set $\{1,2,\dots,m\}$ is a tight 2-independent set in $\mathbb{Z}_n$.

\item If $n=4m$ and $\kappa$ is divisible by 4, then the set
$$\{(g,2i+1)  \mid g \in G_1, i=0,1,\dots, \kappa/4 -1\}$$ is a tight 3-independent set in $G = G_1 \times \mathbb{Z}_{\kappa}$. 
For example, the set $\{1,3,\dots,2m-1\}$ is a tight 3-independent set in $\mathbb{Z}_n$.

\item If $n=t+1$, then the set $\{1\}$ is a perfect $t$-independent set in $\mathbb{Z}_n$.

\item Let $n=\lfloor t^2/2 \rfloor +t+1$.  
\begin{itemize} \item If $t$ is even, then the set $\{t/2 ,t/2+1\}$ is a tight $t$-independent set in $\mathbb{Z}_n$; \item if $t$ is odd, then the set $\{1,t\}$ is a tight $t$-independent set in $\mathbb{Z}_n$;
\item the set $\{(1,1),(1,3)\}$ is a tight $5$-independent set in $\mathbb{Z}_{3} \times \mathbb{Z}_{6}$.
\end{itemize}

\item The set $\{1, 7, 11\}$ is a tight $5$-independent set in $\mathbb{Z}_{38}$.

\end{enumerate}

\end{prop}

Each part of Proposition \ref{tightex}, other than the two specific examples (which can be easily verified by hand or on a computer), is either obvious or follows from previously discussed results.  Note that the sets given in Proposition \ref{tightex} are not unique.  

We could not find tight independent sets for $t \geq 4$ and $m \geq 3$ for any $n$ other than the seemingly sporadic example listed last.  It might be an interesting problem to find and classify all tight independent sets.

\begin{prob}
For each positive integer $n$, find all tight $t$-independent sets in $\mathbb{Z}_n$ of size 2.

\end{prob}

\begin{prob}
For any given noncyclic group $G$, find all tight $t$-independent sets in $G$ of size 2.

\end{prob}

\begin{prob}

Find tight $t$-independent sets in $G=\mathbb{Z}_{n}$ of size $m$ for some values $t \geq 4$ and $m \geq 3$ other than for $t=5$ and $n=38$, or prove that such tight independent sets do not exist.

\end{prob}

\begin{prob}

Find tight $t$-independent sets in noncyclic groups of size $m$ for values $t \geq 4$ and $m \geq 2$ (cf.~the example for $\mathbb{Z}_{3} \times \mathbb{Z}_{6}$ above).

\end{prob}

Finally, a more general problem.

\begin{prob} \label{ind}

For each given value of $n$, find all values of $m$ and $t$ for which the group $\mathbb{Z}_n$ has a $t$-independent set of size $m$.
\end{prob}   
Proposition \ref{indepbound} puts a necessary condition on $n$, $m$, and $t$, but that inequality is not necessarily sufficient.

\subsection{Arbitrary number of terms} \label{5maxUSarbitrary}

Here we ought to consider $$\tau_{\pm} (G,H) =\mathrm{max} \{ |A|  \mid A \subseteq G, 0 \not \in H_{\pm}A \}$$ for the case when $H$ is the set of all nonnegative or all positive integers.  However, as we have already mentioned, we have $\tau_{\pm} (G,H)=0$ whenever $H$ contains a multiple of the exponent of the group (including $0$).  Thus, there are no ``infinitely independent'' sets in a finite group.

\section{Restricted sumsets} \label{5maxR}

In this section we investigate the maximum possible size of a weakly zero-$H$-sum-free set, that is, the quantity $$\tau \hat{\;} (G,H) =\mathrm{max} \{ |A|  \mid A \subseteq G, 0 \not \in H \hat{\;} A \}$$  (if there is no subset $A$ for which $ 0 \not \in H \hat{\;} A$, we let $\tau \hat{\;} (G,H) =0$).  Clearly, we always have $\tau \hat{\;} (G,H)=0$ when $0 \in H$; however, when $ 0 \not \in H$ and $n \geq 2$, then $\tau \hat{\;} (G,H) \geq 1$: for any $a \in G \setminus \{0\}$, for the one-element set $A=\{a\}$ we obviously have $0 \not \in H \hat{\;} A$.

It is important to note that Proposition \ref{zetafordirectsum} does not carry through for $\tau \hat{\;} (G,H)$: for example, $\{1\}$ is trivially weakly zero-$\mathbb{N}$-sum-free in $\mathbb{Z}_5$, but $\{1\} \times \mathbb{Z}_5$ is not weakly zero-$\mathbb{N}$-sum-free in $\mathbb{Z}_5^2$ since $$(1,0)+(1,1)+(1,2)+(1,3)+(1,4)=(0,0).$$  

We consider three special cases: when $H$ consists of a single positive integer $h$,  when $H$ consists of all positive integers up to some value $t$, and when $H = \mathbb{N}$.  As noted above, we have $\tau \hat{\;} (G,H)=0$ whenever $0 \in H$, so the cases when $H =[0,t]$ or $H = \mathbb{N}_0$ yield no weakly zero-sum-free sets.

\subsection{Fixed number of terms} \label{5maxRfixed}

The analogue of a zero-$h$-sum-free set for restricted addition is called a weak or weakly zero-$h$-sum-free set.  In particular, in this section we investigate, for a given group $G$ and positive integer $h$, the quantity $$\tau\hat{\;} (G,h) = \mathrm{max} \{ |A|  \mid A \subseteq G, 0 \not \in h\hat{\;}A\},$$ that is, the maximum size of a weak zero-$h$-sum-free subset of $G$.

Since we trivially have $\tau\hat{\;}(G,h)=n$ for every $h \geq n+1$, we assume below that $h \leq n$.

We start with the following obvious bounds:

\begin{prop} \label{tau hat bounds}
When $\chi \hat{\;}(G,h)$ exists, then we have
$$\tau (G,h) \leq  \tau\hat{\;}(G,h) \leq \chi \hat{\;}(G,h) -1.$$ 
\end{prop}

As before, we can easily verify the following:

\begin{prop} \label{tau hat h=1,2}  In a group of order $n$, we have
$$\tau\hat{\;}(G,1)=n-1$$
and
$$\tau\hat{\;}(G,2)=\frac{n+|\mathrm{Ord}(G,2)|+1}{2};$$  
as a special case, for the cyclic group of order $n$ we have 
 $$\tau\hat{\;}(\mathbb{Z}_n,2)=\left\lfloor \frac{n+2}{2} \right\rfloor.$$

\end{prop}

Next, we consider  the other end of the spectrum: $h=n$.  As Blyler\index{Blyler, N.} in (\cite{Bly:2010a}) observed, in $\mathbb{Z}_n$ we have $$0+1+\cdots+(n-1) =\frac{(n-1) \cdot n }{2}  \left\{
\begin{array}{cl}
\neq 0 & \mbox{if} \; n \; \mbox{even}; \\ \\
= 0 & \mbox{if} \; n \; \mbox{odd}, \\ 
\end{array}\right.$$
and thus 
$$\tau\hat{\;}(\mathbb{Z}_n,n)=\left\{
\begin{array}{cl}
n & \mbox{if} \; n \; \mbox{even}; \\ \\
n-1 & \mbox{if} \; n \; \mbox{odd}. \\ 
\end{array}\right.$$  

More generally, suppose that $G$ has exponent $\kappa$ and $G \cong G_1 \times \mathbb{Z}_{\kappa}$.   The elements of $G_1 \times \mathbb{Z}_{\kappa}$ sum to \label{tau hat in direct prod} 
 $$ \sum_{g \in G_1} \left(\kappa \cdot g, \sum_{i=0}^{\kappa-1} i \right) = \sum_{g \in G_1} \left(0, \sum_{i=0}^{\kappa-1} i \right) = \left(0, |G_1| \cdot \sum_{i=0}^{\kappa-1} i \right) = \left(0, |G_1| \cdot \frac{(\kappa-1)  \cdot \kappa }{2} \right).$$ Here the second component is zero if, and only if, $|G_1|= n/ \kappa$ is even or if $\kappa$ is odd.  Therefore:

\begin{prop} \label{tau hat h=n}
Suppose that $G$ has order $n$ and exponent $\kappa$.  Then  
$$\tau\hat{\;}(G,n)=\left\{
\begin{array}{cl}
n & \mbox{if}  \; \kappa \; \mbox{is even and} \; n/\kappa \; \mbox{is odd}; \\ \\
n-1 & \mbox{otherwise}. \\ 
\end{array}\right.$$
\end{prop}

We can also evaluate $\tau\hat{\;}(G,n-1)$.  Let $g \in G$ denote the sum of all $n$ elements of $G$.  Then $G \setminus \{g\}$ is not a weak zero-$(n-1)$-sum-free set, since its elements add to zero, so $\tau\hat{\;}(G,n-1)  \leq n-1$.  On the other hand, for every $g' \in G \setminus \{g\}$, the set $G \setminus \{g'\}$ is a weak zero-$(n-1)$-sum-free set, since its elements add to $g-g' \neq 0$, so $\tau\hat{\;}(G,n-1) \geq n-1$.  Thus:

\begin{prop}
For every group $G$ of order $n$ we have $\tau\hat{\;}(G,n-1)=n-1$.
\end{prop}

Furthermore, Bajnok and Edwards\index{Bajnok, B.}\index{Edwards, S.} determined the value of $\tau\hat{\;}(G,h)$ for every $G$ in a wide range of $h$ values:

\begin{thm} [Bajnok and Edwards; cf.~\cite{BajEdw:2016a}]  \label{Sam tau hat large h}\index{Bajnok, B.}\index{Edwards, S.}
For every abelian group $G$ of order $n$ and all $$\frac{n+|\mathrm{Ord}(G,2)|-1}{2} \leq h \leq n-2,$$ we have $\tau\hat{\;}(G,h)=h+1$, with the following exceptions:
\begin{itemize}
  \item $\tau\hat{\;}(G,n-3)=n-3$ when $\kappa=3$ (that is, $G$ is isomorphic to an elementary abelian 3-group);
  \item $\tau\hat{\;}(G,n-2)=n-2$ when $|\mathrm{Ord}(G,2)|=1$ and $\kappa \equiv 2$ mod 4.
\end{itemize}  
\end{thm}

Note that Theorem \ref{Sam tau hat large h} does not apply to elementary abelian 2-groups; for these groups we have the value of $\tau\hat{\;}(G,h)$ for all $h \geq n/2-1$:

\begin{thm} [Bajnok and Edwards; cf.~\cite{BajEdw:2016a}]  \label{Sam tau hat large h 2-group}\index{Bajnok, B.}\index{Edwards, S.}
For every positive integer $r$ we have $$\tau\hat{\;}(\mathbb{Z}_2^r,h)=\left \{
\begin{array}{cl}
h+2 & \mbox{if $2^{r-1}-1 \leq h \leq 2^r-5$ or $h=2^r-3$}; \\  \\
h & \mbox{if $h=2^r-4$}.
\end{array}
\right.$$
\end{thm}

For the elementary abelian 2-group, we also have the value $\tau\hat{\;}(\mathbb{Z}_2^r,3)$.  
 We readily have two different weakly zero-3-sum-free sets of size $2^{r-1}+1$: 
\begin{itemize}
  \item the zero element together with all elements with last component 1; and
  \item the zero element together with all elements with an odd number of 1 components.
   \end{itemize}
Therefore, $\tau\hat{\;}(\mathbb{Z}_2^r,3)$ is at least $2^{r-1}+1$.
On the other hand, by Proposition \ref{tau hat bounds} and Theorem \ref{Roth and Lempel even}, $\tau\hat{\;}(\mathbb{Z}_2^r,3)$ is at most $2^{r-1}+1$.  Therefore:

\begin{thm}  \label{thm tau hat 2-group 3}
For all positive integers $r$ we have $\tau\hat{\;}(\mathbb{Z}_2^r,3) = 2^{r-1}+1$.

\end{thm}
(Edwards\index{Edwards, S.} proved Theorem \ref{thm tau hat 2-group 3} independently in \cite{Edw:2015a}.)

Comparing the results thus far in this section to the corresponding results in Section \ref{CritRfixed}, we may get the impression that $\tau\hat{\;}(G,h)$ and $\chi\hat{\;}(G,h)$ are always very close to each other, in fact, tend to only differ by one (cf.~Proposition \ref{tau hat bounds}).  Our next example shows that, in fact, $\tau\hat{\;}(G,h)$ and $\chi\hat{\;}(G,h)$ can be arbitrarily far from one another.

Suppose that $A$ is a weak zero-4-sum-free set in $\mathbb{Z}_2^r$; we will show that it is then a weak Sidon set as well.  Indeed, if $$a_1+a_2=a_3+a_4$$ for some elements $a_1,a_2,a_3,a_4 \in \mathbb{Z}_2^r$ with $a_1 \neq a_2$ and $a_3 \neq a_4$, then we also have $$a_1+a_2+a_3+a_4=0,$$ which implies that at least two of the elements are equal, but that leads to $\{a_1,a_2\}=\{a_3,a_4\}$.  We can show the same way that the converse holds as well, and therefore:

\begin{prop}
A subset of the elementary abelian 2-group $\mathbb{Z}_2^r$ is a weak zero-4-sum-free set if, and only if, it is a weak Sidon set.  Consequently, $$\tau\hat{\;}(\mathbb{Z}_2^r,4)=\sigma \hat{\;}(\mathbb{Z}_2^r,2).$$

\end{prop} 

Now recall that, by Proposition \ref{zetaupperRfixed}, if an $m$-subset of $\mathbb{Z}_2^r$ is a weak Sidon set, then $$2^r \geq {m \choose 2};$$ consequently, $$\tau\hat{\;}(\mathbb{Z}_2^r,4) \leq 2^{(r+1)/2}+1.$$  On the other hand, we have
$$\chi \hat{\;} (\mathbb{Z}_2^r,4) \geq \chi (\mathbb{Z}_2^r,4)=v_1(2^r,4)+1=2^{r-1}+1.$$  Therefore:

\begin{prop}
We have
$$\lim_{r \rightarrow \infty} \left( \chi \hat{\;} (\mathbb{Z}_2^r,4) - \tau \hat{\;} (\mathbb{Z}_2^r,4) \right) = \infty.$$
\end{prop}

We have the following open questions:

\begin{prob}
For positive integer $r$ and each $h$ with 
$3 \leq h \leq 2^{r-1}-2,$ find $\tau\hat{\;}(\mathbb{Z}_2^r,h)$.

\end{prob}

\begin{prob}
For each abelian group $G$ of order $n$ and each $h$ with 

$$3 \leq h \leq \frac{n+|\mathrm{Ord}(G,2)|-3}{2},$$ find $\tau\hat{\;}(G,h)$.

\end{prob}

We can develop some useful lower bounds for $\tau\hat{\;}(G,h)$ by constructing explicit weak zero-$h$-some-free sets, as follows.  We start by considering cyclic groups.  The most general result we have thus far is as follows:

\begin{prop} \label{Zlower}
For positive integers $n$ and $h \leq n$, let $r$ be the nonnegative remainder of $(h^2-h-2)/2$ when divided by $\gcd(n,h)$.   We have
$$\tau\hat{\;}(\mathbb{Z}_n,h) \geq  \left\lfloor \frac{n+h^2- r-2}{h} \right\rfloor .$$
\end{prop}
We should note that $(h^2-h-2)/2$ is an integer for all $h$.  
For the proof, see page \pageref{proofofZlower}.

Since $r$ is at most $\gcd(n,h)-1$, we have the following immediate consequence: 

\begin{cor} \label{cor Zlower 1}
For all positive integers $n$ and $h \leq n$, 
$$\tau\hat{\;}(\mathbb{Z}_n,h) \geq  \left\lfloor \frac{n+h^2- \mathrm{gcd} (n,h) -1}{h} \right\rfloor .$$
\end{cor}

Furthermore, when $n$ is divisible by $h$, then $\gcd(n,h)=h$; if $h$ is even, then 
$$\frac{h^2-h-2}{2}=\left( \frac{h}{2} - 1 \right) \cdot h + \frac{h-2}{2},$$ so $r=(h-2)/2$, and we get the following: 

\begin{cor} \label{cor Zlower 2}
If $h$ is even and $n$ is divisible by $h$, then
$$\tau\hat{\;}(\mathbb{Z}_n,h) \geq \frac{n}{h}+h-1.$$
\end{cor}

Note that when $n$ is divisible by $h$, then $$\left\lfloor \frac{n+h^2- \mathrm{gcd} (n,h) -1}{h} \right\rfloor =\frac{n}{h}+h-2,$$ so Corollary \ref{cor Zlower 2} is stronger than Corollary \ref{cor Zlower 1} in this case.

We can do very slightly better when $n$ is divisible by $h^2$ and $h>2$.  This was first observed and proved by Yager-Elorriaga\index{Yager-Elorriaga, D.} in \cite{Yag:2008a} for $h=3$ and $h=4$; we provide a general proof.

\begin{prop} \label{Zlower2}
For positive integers $n$ and $h \geq 3$ for which $n$ is divisible by $h^2$ we have
$$\tau\hat{\;}(\mathbb{Z}_n,h) \geq \frac{n}{h}+h.$$
\end{prop}
For the proof, see page \pageref{proofofZlower2}.

There are some additional instances where we can do better---we present these (rather easy) results below.

First, we observe that the set
$$A=\left\{\pm 1, \pm 2, \dots, \pm (h+1)/2\right\}$$ is weakly zero-$h$-sum-free when $h$ is odd and $n \geq h+2$; indeed, we have $h\hat{\;}A=A$.  The similar, but slightly less obvious set
$$A=\left\{0,1,\pm 2, \pm 3, \dots, \pm h/2, h/2+1 \right\}$$ works when $h$ is even and $n \geq h+3$, since in this case we have $$h\hat{\;}A = \left\{1,2, \dots, h/2+2, h/2+4,  h/2+5, h+2 \right\}.$$  This gives the following.

\begin{prop} \label{Zlower11}
Suppose that $n \geq h+2$ when $h$ is odd and $n \geq h+3$ when $h$ is even.  Then $\tau\hat{\;}(\mathbb{Z}_n,h) \geq h+1$.

\end{prop} 

Another simple construction applies when $n$ is even and $h$ is odd: clearly, the set
$$\{1,3,5,\dots,n-1\}$$ is then weakly zero-$h$-some-free.  As a generalization, assume that $n$ is divisible by some positive integer $d$, and that, for some nonnegative integer $k$, none of $h, h+1,\dots,h+k$ is divisible by $d$.  Consider then the set $A=A_1 \cup A_2$ where
$$A_1=\left\{1,d+1,2d+1,\dots,\left(n/d-1 \right)d+1\right\}$$ and
$$A_2=\left\{2,d+2,2d+2,\dots,(k-1)d+2\right\}.$$

To prove that $A$ is weakly zero-$h$-some-free, consider an element $x$ of $h \hat{\;} A$ where $h_1$ terms come from $A_1$ and $h_2$ terms come from $A_2$.  This sum is then of the form $$x=(q_1d+h_1)+(q_2d+2h_2)=(q_1+q_2)d+h+h_2$$ for some integers $q_1$ and $q_2$.  Since $h \leq h+h_2 \leq h+k$, by assumption, $h+h_2$ is not divisible by $d$, and therefore $x$ is not divisible by $d$.  But then it cannot be divisible by $n$ either, proving that $0 \not \in h \hat{\;} A$.  We just proved the following:

\begin{prop} \label{Zlower1}
Let  $n$ and $h \leq n$ be positive integers, and let $d$ be a positive divisor of $n$ for which none of $h, h+1,\dots,h+k$ is divisible by $d$.  We then have 
$$\tau\hat{\;}(\mathbb{Z}_n,h) \geq \frac{n}{d}+k.$$
\end{prop}

The following result is not difficult to establish either.  For fixed positive integers $n_1$, $n_2$, and $h$, consider the set $$A=\left\{0,1,2,\dots, \left\lfloor n_1/h \right\rfloor \right\}$$ in $\mathbb{Z}_{n_1}$ and any weakly zero-$h$-sum-free set $B$ in $\mathbb{Z}_{n_2}$.  We show that $$A \times B=\{(a,b)  \mid a \in A, b \in B\}$$ is weakly zero-$h$-sum-free in $\mathbb{Z}_{n_1} \times \mathbb{Z}_{n_2}$.  

Indeed, if the sum of $h$ elements $(a_1,b_1), \dots, (a_h,b_h)$ of $A \times B$ were to equal $(0,0)$, then, since $$0 \leq a_1+\cdots + a_h \leq h \cdot \left\lfloor n_1/h \right\rfloor \leq n_1,$$ either $a_1=\cdots =a_h=0$ or ($n_1$ is divisible by $h$ and) $a_1=\cdots=a_h=n_1/h$.  But if all first coordinates are equal, then all second coordinates must be distinct, but that would mean that we have $h$ distinct elements in $B$ that sum to 0 in $\mathbb{Z}_{n_2}$, which is a contradiction.  Thus, we proved the following result.

\begin{prop} \label{Zlower3}  
For positive integers $n_1$, $n_2$, and $h$ we have
$$\tau\hat{\;}(\mathbb{Z}_{n_1} \times \mathbb{Z}_{n_2} ,h) \geq \left( \left\lfloor n_1/h \right\rfloor +1 \right) \cdot \tau\hat{\;}(\mathbb{Z}_{n_2} ,h).$$
\end{prop}
Note that Proposition \ref{Zlower3} applies even to a cyclic group, as long as its order has more than one prime divisor.

Relying on a computer program, Yager-Elorriaga\index{Yager-Elorriaga, D.} in \cite{Yag:2008a} exhibited the values of $\tau\hat{\;}(\mathbb{Z}_n,h)$ for $h\in \{3,4\}$ and all $n \leq 30$; Blyler\index{Blyler, N.} in \cite{Bly:2010a} extended this to $h \leq 10$---we provide these data in the table below.    Most entries are in agreement with the highest of the values that Propositions  \ref{Zlower}, \ref{Zlower2},  \ref{Zlower11},  \ref{Zlower1}, or \ref{Zlower3} yield.  (None of these propositions is superfluous.)  For some entries this is not the case; we marked these entries in the table by a * and listed them here with an exemplary set.

$$\begin{array}{ll}
\tau\hat{\;}(\mathbb{Z}_{10},3)=6 & \{1,2,4,6,8,9\} \\
\tau\hat{\;}(\mathbb{Z}_{15},5)=8 & \{1,2,4,5,7,10,11,14\} \\
\tau\hat{\;}(\mathbb{Z}_{15},6)=8 & \{0,1,2,3,4,6,12,13\} \\
\tau\hat{\;}(\mathbb{Z}_{27},6)=10 & \{0,2,4,7,9,11,13,18,20,22\} \\
\tau\hat{\;}(\mathbb{Z}_{21},7)=9 & \{0,1,2,3,4,5,7,8,9\} \\
\tau\hat{\;}(\mathbb{Z}_{20},8)=10 & \{0,1,2,3,4,5,6,8,9,18\} \\
\tau\hat{\;}(\mathbb{Z}_{28},8)=11 & \{0,1,2,3,4,5,6,8,24,25,26\} \\
\tau\hat{\;}(\mathbb{Z}_{21},9)=11 & \{0,1,2,3,4,5,6,7,10,12,13\} \\
\tau\hat{\;}(\mathbb{Z}_{27},9)=12 & \{0,1,2,3,9,10,11,12,18,19,20,21\} \\
\tau\hat{\;}(\mathbb{Z}_{25},10)=12 & \{0,1,2,3,4,5,6,7,8,10,13,15\}
\end{array}$$


$$\begin{array}{||c||c|c|c|c|c|c|c|c||} \hline \hline
n & \tau\hat{\;}(\mathbb{Z}_n,3) & \tau\hat{\;}(\mathbb{Z}_n,4) & \tau\hat{\;}(\mathbb{Z}_n,5) & \tau\hat{\;}(\mathbb{Z}_n,6) & \tau\hat{\;}(\mathbb{Z}_n,7)  & \tau\hat{\;}(\mathbb{Z}_n,8) & \tau\hat{\;}(\mathbb{Z}_n,9) & \tau\hat{\;}(\mathbb{Z}_n,10) \\ \hline \hline
 1& 1& 1 & 1 & 1 & 1 &   1&    1 &  1             \\ \hline 
 2& 2& 2 &  2 & 2 & 2 &  2 &   2 &   2              \\ \hline
 3& 2& 3 &  3 & 3 & 3 &  3 &  3 & 3                    \\ \hline  
 4& 3& 4 &  4 & 4 & 4 & 4 & 4 &  4                        \\ \hline
 5& 4& 4&  4 & 5 & 5 &  5 &  5 &   5                   \\ \hline
 6& 4& 4 & 5 & 6 &  6 & 6 & 6 &  6                       \\ \hline
7& 4& 5 & 6 &  6 & 6 & 7 &  7 &  7                        \\ \hline
8& 5& 5 &  6 & 7  & 7 &  8 & 8 & 8                       \\ \hline
9& 6& 5 &  6 & 7  & 8 &  8 &  8 &   9                \\ \hline
10& 6 *& 5 & 6  & 7  & 8 &  8 & 9 &  10                   \\ \hline
11& 6& 6 & 6 &  7 & 8 & 9 &  10 &   10                       \\ \hline
 12& 6& 6 & 7 & 7  & 8 &  9 & 10  &  10                           \\ \hline
13& 6&  6 & 7 & 7 & 8 & 9 &  10 &    11                           \\ \hline
 14& 7& 6 &  7 & 8  & 8  &   9 &  10 &  11                            \\ \hline
15& 8& 7 & 8* &  8 * & 8 & 9 &  10 & 11                               \\ \hline
 16& 8& 8 & 8 & 8 & 9 &  9  & 10 &  11                               \\ \hline
17& 8& 7 &  8 & 8 & 9 &  9 &  10 &   11                          \\ \hline
 18& 9& 7 & 9 &  8  & 9 & 9 &  10  &   11                                       \\ \hline
19& 8& 8 &  8 & 8 & 9 &   10 & 10 &  11                                          \\ \hline
 20& 10& 8 &  10 & 9  & 10 &  10 * &  11 & 11                                 \\ \hline
21& 8& 8 &  8 & 8 & 9 * &  10 &   11 * &  11                                        \\ \hline
 22& 11& 8 &  11 & 9 & 11 & 10 & 11 &   12                                         \\ \hline
23& 10& 9 &  9 & 9 & 10 & 10 &   11 &   12                                 \\ \hline
 24& 12& 9 &  12 & 9  & 12 &  10  &  12 &  12                                \\ \hline
25& 10& 9 &  10 & 9 & 10 &   10 & 11 &   12 *                                   \\ \hline
 26& 13& 9 &  13 & 10  & 13 &  10 &  13 &  12                          \\ \hline
27& 12& 10 &  10 & 10 * & 10 & 11 &  12*  &  12                              \\ \hline
 28& 14& 10 & 14 &  10 & 14 & 11 * &  14 &  12                         \\ \hline
29& 12&  10 & 10 & 10 &  10 &   11 & 12 &  12                           \\ \hline
 30& 15& 11 & 15 & 10 & 15 &  11 &    15 &  12                       \\ \hline
\hline
\end{array}$$

After these lower bounds, we mention some exact results.  First, the value of $\tau\hat{\;}(\mathbb{Z}_p,h)$ for $p$ prime immediately follows from combining Propositions \ref{tau hat bounds} and \ref{Zlower} and Theorem \ref{prime rest crit}:

\begin{thm}  \label{Zforp} 

Suppose that $p$ is a positive prime and $1 \leq h \leq p-1$.  We then have
$$\tau\hat{\;}(\mathbb{Z}_p,h) =
\left \lfloor \frac{p-2}{h} \right \rfloor + h.$$
\end{thm}

Next, we have:

\begin{thm}  \label{Zfor n even h odd} 

If $n \geq 12$ is even and $3 \leq h \leq n-1$ is odd, then 
$$\tau\hat{\;}(\mathbb{Z}_n,h) = \chi \hat{\;} (\mathbb{Z}_{n},h) -1 = \left \{
\begin{array}{cl}
n-1 & \mbox{if} \; h=1; \\ \\
n/2 & \mbox{if} \; 3 \leq h \leq n/2-2; \\ \\
n/2+1 & \mbox{if} \; h=n/2-1; \\ \\
h+1 & \mbox{if} \; n/2 \leq h \leq  n-2; \\ \\
n-1 & \mbox{if} \; h=n-1.
\end{array}
\right.$$
\end{thm}

The proof can be found on page \pageref{proof of Zfor n even h odd}.

We pose the following problems.

\begin{prob}

Find the exact formulas for $\tau\hat{\;}(\mathbb{Z}_n,3)$ and $\tau\hat{\;}(\mathbb{Z}_n,5)$ for all odd $n$.

\end{prob}

\begin{prob}

Find $\tau\hat{\;}(\mathbb{Z}_n,4)$.

\end{prob}

Regarding noncyclic groups, we can use the same argument as the one provided for Theorem \ref{tau=chi -1 when h and n are coprime} to prove the following:

\begin{thm}
If $h$ is relatively prime to $n$, then $$\tau\hat{\;}(G,h)=\chi \hat{\;}(G,h)-1.$$

\end{thm}

Let us now turn to elementary abelian groups.  We have already considered the elementary abelian 2-group  $\mathbb{Z}_2^r$, so we move on to the elementary abelian 3-groups.  Pursuing the determination of $\tau\hat{\;}(\mathbb{Z}_3^r,3)$ is attracting much attention as it is related to finite geometry.   
For example:
\begin{itemize}
  \item $\tau\hat{\;}(\mathbb{Z}_3^r,3)$ is the maximum size of a cap (a collection of points without any three being collinear) in the affine space $\mathrm{AG}(r,3)$; and
  \item $2 \cdot \tau\hat{\;}(\mathbb{Z}_3^r,3)$ is the maximum number of points in the integer lattice $\mathbb{Z}^r$ so that the centroid of no three of them is a lattice point. 
  \end{itemize} 
We elaborated more on these questions in an `appetizer' section, see page \pageref{appetizer on SET}.
 
The following table summarizes all values of $\tau\hat{\;}(\mathbb{Z}_3^r,3)$ that are known (see \cite{GaoTha:2004a} by Gao and Thangadurai\index{Gao, W.}\index{Thangadurai, R.}  and its references for the first five entries and \cite{Pot:2008a} by Potechin for the last): \label{tau hat 3 table}\index{Potechin, A.} 
$$\begin{array}{|c||r|r|r|c|c|c|} \hline 
r & 1 & 2 & 3 & 4 & 5 & 6 \\ \hline
\tau\hat{\;}(\mathbb{Z}_3^r,3) & 2 & 4 & 9 & 20 & 45 & 112 \\ \hline
\end{array}$$
We also know that the sets of maximum size for $r \leq 6$ are essentially unique.  

Even good bounds are difficult to achieve for $\tau\hat{\;}(\mathbb{Z}_3^r,3)$---we summarize what is currently known.  Starting with lower bounds, observe that the set
$\{0,1\}^r$ is clearly weakly zero-3-sum-free in  $\mathbb{Z}_3^r$, hence 
$$ \tau\hat{\;}(\mathbb{Z}_3^r,3) \geq 2^r.$$

A better and very nice lower bound can be developed, as follows.  \label{bui'sidea} For $i=0,1,\dots,r-1$, consider the collection of elements $A_i$ of $\mathbb{Z}_3^r$ that contain exactly $i$ 0-components, and whose remaining $r-i$ components are all 1, and let $-A_i$ denote the negatives of the elements of $A_i$ (which contain exactly $i$ 0-components, and whose remaining $r-i$ components are all 2).  (We will not use $A_r$.)  Note that $$|A_i \cup -A_i|=2 \cdot {r \choose i}$$ for each $i=0,1,\dots,r-1$.  Bui (see \cite{Bui:2009a})\index{Bui, C.} proved the following interesting result: If $I$ is a subset of $\{0,1,\dots,r-1\}$ for which the equation $$i_1+i_2-i_3=r$$ has no solution with $i_1,i_2,i_3 \in I$, then $A=\cup_{i \in I} (A_i \cup -A_i)$ is weakly zero-3-sum-free.  (Note that we had to exclude $r$ from $I$ as $i_1=i_2=i_3=r$ does yield a solution to the equation.)  For example, for $r=3$, the index set $I=\{0,1\}$ yields no solution for $$i_1+i_2-i_3=3,$$ and, correspondingly, the set $A=\pm A_0 \cup \pm A_1$, consisting of the eight elements
$$(1,1,1), (2,2,2), (0,1,1), (0,2,2), (1,0,1), (2,0,2), (1,1,0), (2,2,0),$$
is, as can be easily verified, weakly zero-3-sum-free in  $\mathbb{Z}_3^3$: no three of them add to $(0,0,0)$.  (The index set $I=\{0,2\}$ would work as well.)  

We can select the index set $I$ to maximize $|A|$ by positioning it in the ``middle'' of $\{0,1,\dots,r-1\}$.  More precisely, if $r=3q+s$ with a remainder $s=0,1,2$, then letting $$I=\{q+s-1,q+s,q+s+1,\dots,2q+s-1\}$$ works, since for any $i,_1,i_2,i_3 \in I$, we have $i_1+i_2-i_3$ at most equal to $$(2q+s-1)+(2q+s-1)-(q+s-1)=3q+s-1=r-1<r,$$ and it yields an optimal $|A|$.  This results in the following lower bound:

\begin{thm} [Bui; cf.~\cite{Bui:2009a}] \label{bui}\index{Bui, C.}
Let $r$  be any positive integer, and write $r=3q+s$ with $s=0,1$, or 2.  Then we have
$$\tau\hat{\;}(\mathbb{Z}_3^r,3) \geq 2 \cdot \sum_{i=q+s-1}^{2q+s-1} {r \choose i}.$$
\end{thm} 

Theorem \ref{bui} gives a remarkably good lower bound for the values of $\tau\hat{\;}(\mathbb{Z}_3^r,3)$; for example, for $r=4$, 5, and 6, we get that $\tau\hat{\;}(\mathbb{Z}_3^r,3)$ is at least 20 (the actual value, see above), 40, and 82, respectively.

It would be very interesting to see if Bui's idea\index{Bui, C.} generalizes to other settings; we pose the following problems.

\begin{prob}

Find a lower bound for $\tau\hat{\;}(\mathbb{Z}_k^r,3)$ for $k \geq 4$.

\end{prob} 

\begin{prob}

Find a lower bound for $\tau\hat{\;}(\mathbb{Z}_3^r,h)$ for $h \geq 4$.

\end{prob}

Beyond the lower bound of Theorem \ref{bui}, Edel in \cite{Ede:2004a},\index{Edel, Y.} improving on a result of Frankl, Graham, and R\"odl\index{Frankl, P.}\index{Graham, R. L.}\index{Rodl@R\"odl, V.} in \cite{FraGraRod:1987a}, has shown that there are infinitely many values of $r$ for which 
$$\tau\hat{\;}(\mathbb{Z}_3^r,3) > 2.2^r.$$  (Note that the lower bound of Theorem \ref{bui} is less than $2^{r+1}$.)  

An easy upper bound can be derived from the fact that if $A$ is weakly zero-3-sum-free in $\mathbb{Z}_3^r$, then it is disjoint from $$B=-a_1-(A \setminus \{a_1\})$$ for any fixed $a_1 \in A$.  Indeed, if $-a_1-a_i =a_j$ for some 
$a_i \in A \setminus \{a_1\}$ and $a_j \in A$, then $a_1+a_i+a_j=0$, so the fact that $A$ is weakly zero-3-sum-free implies that these three elements of $A$ cannot be all distinct, but in $\mathbb{Z}_3^r$ this actually implies that all three are equal, contradicting $a_i \in A \setminus \{a_1\}$.  Since $|B|=|A|-1$, we get
$$ \tau\hat{\;}(\mathbb{Z}_3^r,3) \leq (3^r+1)/2.$$ 
This bound was improved to
$$ \tau\hat{\;}(\mathbb{Z}_3^r,3) \leq (r+1) \cdot 3^r /r^2$$
by Bierbrauer and Edel in \cite{BieEde:2002a},\index{Bierbrauer, J.}\index{Edel, Y.} where the upper bound is of approximate size $3^r/r$ for large $r$.  More recently, Bateman and Katz\index{Bateman, M.}\index{Katz, N. H.} proved in their lengthy paper \cite{BatKat:2012a} that we actually have $$\lim_{r \rightarrow \infty} \tau\hat{\;}(\mathbb{Z}_3^r,3)/(3^r/r) =0.$$
We offer the following very difficult problem:

\begin{prob}

Find better lower and upper bounds for $\tau\hat{\;}(\mathbb{Z}_3^r,3)$.

\end{prob}

The case when $h$ equals the exponent $\kappa$ of the group, as with the elementary abelian 3-group and $h=3$ above, has attracted much attention: the value of $\tau\hat{\;}(G, \kappa)+1$ (the smallest integer such that each subset of $G$ of that cardinality has a subset size $\kappa$ whose elements sum to 0) is called the {\em Harborth constant} of $G$.  Below we summarize what is known about $\tau\hat{\;}(G, \kappa)$.

It is obvious that for every $G$ of order $n$ and exponent $\kappa$ we have 
$$\kappa -1 \leq \tau\hat{\;}(G, \kappa) \leq n;$$ we can easily classify cases when either inequality becomes an equality, as follows.  For the elementary abelian 2-group, we have $\tau\hat{\;}(\mathbb{Z}_2^r, 2)=2^r$ by Proposition \ref{tau hat h=1,2}.  Furthermore, for the cyclic group, by Proposition \ref{tau hat h=n}, we have $\tau\hat{\;}(\mathbb{Z}_n, n)=n$ if $n$ is even and $\tau\hat{\;}(\mathbb{Z}_n, n)=n-1$ if $n$ is odd.  Next, we show that in all other cases,  
$$\kappa  \leq \tau\hat{\;}(G, \kappa) \leq n-1.$$
Suppose that $G \cong G_1 \times \mathbb{Z}_{\kappa}$, where $\kappa \geq 3$ (so $G$ is not isomorphic to an elementary abelian 2-group) and $|G_1| \geq 2$ (so $G$ is not cyclic).  To prove our claim, we need to find a subset of size $\kappa$ whose elements don't sum to zero, and another subset of size $\kappa$ whose elements do sum to zero.  Note that the elements in $\{0\} \times \mathbb{Z}_{\kappa}$ sum to zero if $\kappa$ is odd, and don't sum to zero if $\kappa$ is even.  Furthermore, for any $g \in G_1 \setminus \{0\}$, if $\kappa$ is odd, then the elements in $\{0\} \times (\mathbb{Z}_{\kappa}\setminus \{0\}) \cup \{(g,0)\}$ don't sum to zero, and if $\kappa$ is even, then the elements in $\{0\} \times (\mathbb{Z}_{\kappa}\setminus \{0, \kappa/2\}) \cup \{(g,1), (-g,\kappa-1)\}$ sum to zero.    We thus proved:

\begin{prop} \label{Harborth extreme}
For all groups $G$ of order $n$ and exponent $\kappa$ we have 
$$\kappa -1 \leq \tau\hat{\;}(G, \kappa) \leq n.$$  Furthermore:
\begin{itemize}
  \item the lower bound holds if, and only if, $G$ is cyclic of odd order, and
  \item the upper bound holds if, and only if, $G$ is cyclic of even order or $G$ is an elementary abelian 2-group.
\end{itemize}
\end{prop} 
We should note that the claim regarding the upper bound appeared by Gao and Geroldinger in \cite{GaoGer:2006a} (see Lemma 10.1).\index{Gao, W.}\index{Geroldinger, A.}

Regarding exact values, we first consider $\tau\hat{\;}(\mathbb{Z}_k^r,k)$.  As Kemnitz\index{Kemnitz, A.} observed in \cite{Kem:1983a}, the set
$$A=\{0,1\}^{r-1} \times \{0, 1, \dots, k-2\},$$ that is, the collection of elements whose first $r-1$ components are 0 or 1 and whose last component is not $k-1$, is zero-$k$-sum-free: indeed, the sum of $k$ distinct elements in $A$ can only be the zero element of $\mathbb{Z}_k^r$ if they agree in each of their first $r-1$ components, but then their last components must be distinct, which is impossible since we only have $k-1$ choices.  We can do slightly better when $k$ is even: $k$ distinct elements of $$A=\{0,1\}^{r-1} \times \mathbb{Z}_k$$ again must share their first $r-1$ components, but then their last components must be distinct and thus add to 
$$0+1+\cdots+k-1=k/2 \neq 0$$  in $\mathbb{Z}_k$.  Thus we get

\begin{prop} [Kemnitz; cf.~\cite{Kem:1983a}] \label{Kemnitz bounds}\index{Kemnitz, A.} 
For each $k \geq 2$ and $r \geq 1$ we have
$$\tau\hat{\;}(\mathbb{Z}_k^r,k) \geq \left\{
\begin{array}{ccl}
(k-1) \cdot 2^{r-1} & \mbox{if} & k \; \mbox{odd}; \\ \\
k \cdot 2^{r-1} & \mbox{if} & k \; \mbox{even}.
\end{array}
\right.$$
\end{prop}

Note that, by Proposition \ref{tau hat h=n}, equality holds for $r=1$, and it has been conjectured that equality holds for $r=2$ as well:

\begin{conj} [Gao and Thangadurai; cf.~\cite{GaoTha:2004a}] \label{conj GaoTha}\index{Gao, W.}\index{Thangadurai, R.} 
For each $k \geq 2$ we have 
$$\tau\hat{\;}(\mathbb{Z}_k^2,k) = \left\{
\begin{array}{ccl}
2k-2 & \mbox{if} & k \; \mbox{odd}; \\ \\
2k & \mbox{if} & k \; \mbox{even}.
\end{array}
\right.$$
\end{conj}

Kemnitz\index{Kemnitz, A.} in \cite{Kem:1983a} proved that Conjecture \ref{conj GaoTha} holds for $k \in \{2, 3, 5,7\}$, and in \cite{GaoTha:2004a} Gao and Thangadurai\index{Gao, W.}\index{Thangadurai, R.}  established Conjecture \ref{conj GaoTha} for prime values of $k$ with $k \geq 67$; this has been improved somewhat by Gao, Geroldinger, and Schmid:\index{Gao, W.}\index{Geroldinger, A.}\index{Schmid, W.} 

\begin{thm} [Gao, Geroldinger, and Schmid; cf.~\cite{GaoGerSch:2007a}] \label{thm GaoGerSch}\index{Gao, W.}\index{Geroldinger, A.}\index{Schmid, W.} 
If $p \geq 47$ is a prime, then 
$$\tau\hat{\;}(\mathbb{Z}_p^2,p) =2p-2.$$
\end{thm}

Furthermore, in \cite{GaoGerSch:2007a} the authors also determine all weakly zero-$p$-sum-free subsets of $\mathbb{Z}_p$ of size $2p-2$ (for $p \geq 47$).  

Additionally, Gao and Thangadurai\index{Gao, W.}\index{Thangadurai, R.}  showed in \cite{GaoTha:2004a} that $\tau\hat{\;}(\mathbb{Z}_4^2,4) =8$, and Schmid\index{Schmid, W.} and co-authors verified that $\tau\hat{\;}(\mathbb{Z}_6^2,6) =12$; cf.~\cite{Sch:2017a}.

We pose the following problems:

\begin{prob}
Prove that Conjecture \ref{conj GaoTha} holds for prime numbers $11 \leq p \leq 43$.

\end{prob}

\begin{prob}
Prove that Conjecture \ref{conj GaoTha} holds for small composite numbers $k \geq 8$.

\end{prob}

However, we know of at least one case when the inequality is strict in Proposition \ref{Kemnitz bounds}:\index{Kemnitz, A.} We have $\tau \hat{\;} (\mathbb{Z}_3^3,3) \geq 9$, since
$$\{(0,0,0),(1,0,0),(0,1,0),(0,0,1),(1,1,0),(1,1,1),(1,2,1),(2,0,1),(0,1,2)\}$$ is weakly zero-3-sum-free in $\mathbb{Z}_3^3$.

Let us now examine general noncyclic groups, in particular groups of rank 2.  Let $n_1$ and $n_2$ be positive integers with $n_1 \geq 2$ and $n_1 |n_2$, and consider $\mathbb{Z}_{n_1} \times \mathbb{Z}_{n_2}$.  It is easy to see that the set $$(\{0\} \times A_2)  \cup (\{1\} \times A_2')$$ is weakly zero-$n_2$-sum-free for all $A_2$ and $A_2'$ of $\mathbb{Z}_{n_2}$ when $|A_2|=n_2-1$ and $|A_2'|=n_1-1$.  And, like we mentioned above, when $n_2$ is even, we may replace $A_2$ by all of $\mathbb{Z}_{n_2}$.  This yields the following:

\begin{prop}  \label{tau hat rank two lower triv}
Let $n_1$ and $n_2$ be integers greater than 1 so that $n_2$ is divisible by $n_1$.  We then have $$\tau\hat{\;}(\mathbb{Z}_{n_1} \times \mathbb{Z}_{n_2},n_2) \geq n_1+n_2-2;$$ if $n_2$ is even, we have $$\tau\hat{\;}(\mathbb{Z}_{n_1} \times \mathbb{Z}_{n_2},n_2) \geq n_1+n_2-1.$$

\end{prop} 

The question then arises whether we can do better than Proposition \ref{tau hat rank two lower triv}; we exhibit two cases when we can: one for $n_1=2$, the other for $n_1=3$.  

For the case when $n_1=2$, we follow the construction of Marchan, Ordaz, Ramos, and Schmid; cf.~\cite{MarOrdRamSch:2013a}.\index{Schmid, W.}\index{Marchan, L. E.} \index{Ordaz, O.}\index{Ramos, D.} 
Suppose that $n_2=2k$ with $k$ odd, in which case $\mathbb{Z}_{2} \times \mathbb{Z}_{n_2}$ is isomorphic to $\mathbb{Z}_{2} \times \mathbb{Z}_{2} \times \mathbb{Z}_{k}$.  We define the subsets $B_1$ and $B_2$ of $\mathbb{Z}_k$ as
$$B_1=\{0, 1,2,\dots,(k-1)/2\}$$ and  
$$B_2=\{(k+1)/2, (k+1)/2+1, \dots, k-1, 0\},$$ and set
$$A=(\{0\} \times \{0\} \times B_1) \cup (\{0\} \times \{1\} \times B_2) \cup (\{1\} \times \{0\} \times B_2) \cup (\{1\} \times \{1\} \times B_1).$$
Then $A$ is a subset of $\mathbb{Z}_{2} \times \mathbb{Z}_{2} \times \mathbb{Z}_{k}$ of size $2k+2$, and it is easy to verify that the sum of these $2k+2$ elements equals $(0,0,0)$.  Therefore, having $k$ distinct elements of $A$ adding to $(0,0,0)$ is equivalent to having two of them add to $(0,0,0)$, but this is clearly not the case.  Therefore, when $n_1=2$ and $n_2 \equiv 2$ mod 4, we have  
$$\tau\hat{\;}(\mathbb{Z}_{n_1} \times \mathbb{Z}_{n_2},n_2) \geq n_1+n_2.$$
It turns out that we cannot do better:

\begin{thm} [Marchan, Ordaz, Ramos, and Schmid; cf.~\cite{MarOrdRamSch:2013a}] \label{tau hat Z2 by Z4}\index{Schmid, W.}\index{Marchan, L. E.}\index{Ordaz, O.}\index{Ramos, D.} 
For each positive integer $k$, 
$$\tau\hat{\;}(\mathbb{Z}_2 \times \mathbb{Z}_{2k}, 2k) = \left\{
\begin{array}{ccl}
2k+2 & \mbox{if} & k \; \mbox{odd}; \\ \\
2k+1 & \mbox{if} & k \; \mbox{even}.
\end{array}
\right.$$
\end{thm}

Let us turn now to $n_1=3$.  We already know from Proposition \ref{tau hat rank two lower triv} that
$$\tau\hat{\;}(\mathbb{Z}_3 \times \mathbb{Z}_{3k}, 3k) \geq 3k+2$$ when $k$ is even; we now show that the same bound holds when $k$ is odd and $k \geq 2$ (for $k=1$ we have $\tau\hat{\;}(\mathbb{Z}_3^2, 3) =4$; cf.~Conjecture \ref{conj GaoTha} above).  
Kiefer in \cite{Kie:2016a} presents the following subset of $\mathbb{Z}_3 \times \mathbb{Z}_{3k}$:
$$A=(\{0\} \times  (\mathbb{Z}_{3k} \setminus \{3k-2,3k-1\}) ) \cup \{(1,0),(1,1),(1,3k-6),(2,0)\}.$$   
Then $A$ has size $3k+2$ (note that $k \geq 3$), and these elements add to $(2,3k-2)$.  As with the previous example, for $k$ distinct elements to add to $(0,0)$, we must have two distinct elements that add to $(2,3k-2)$, but one can quickly see that that is impossible.     This yields:

\begin{prop} [Kiefer; cf.~\cite{Kie:2016a}]\index{Kiefer, C.} \label{Kiefer bound}
For all values of $k >1$ we have $$\tau\hat{\;}(\mathbb{Z}_3 \times \mathbb{Z}_{3k}, 3k) \geq 3k+2.$$

\end{prop}

It turns out that one can do a bit better for $k=3$.  Namely, consider the set 
$$A=\{(0,0),(0,1),(0,3),(0,4),(0,6),(0,7),(1,0),(1,1),(1,3),(1,4),(1,6),(1,7)\}$$ in $\mathbb{Z}_3 \times \mathbb{Z}_{9}$; note that $A$ can be rewritten as $$A=\{0,1\} \times (\{0,3,6\}+\{0,1\}).$$  Since the twelve elements add to $(0,6)$, $A$ is weakly zero-9-sum-free if, and only if, $(0,6) \not \in 3 \hat{\;} A$.  Indeed, for three distinct elements of $A$ to add to $(0,6)$, the three first components would need to agree (all $0$s or all $1$s), but then the second components must be distinct.  There are only two possibilities for the second components to add to a number divisible by 3: $0+3+6$ or $1+4+7$, but neither sum equals 6 in $\mathbb{Z}_9$.  Therefore,
$$\tau\hat{\;}(\mathbb{Z}_3 \times \mathbb{Z}_{9}, 9) \geq 12.$$

Schmid\index{Schmid, W.} writes in \cite{Sch:2017a} that he and co-authors have proved that the lower bound in Proposition \ref{Kiefer bound} holds for large enough primes, but the general question remains open:

\begin{prob}

Evaluate $\tau\hat{\;}(\mathbb{Z}_3 \times \mathbb{Z}_{3k}, 3k)$ for each positive integer $k$.

\end{prob}

And, more generally:

\begin{prob}

Evaluate $\tau\hat{\;}(\mathbb{Z}_{n_1} \times \mathbb{Z}_{n_2}, n_2)$ for positive integers $n_1 \geq 4$, $n_1|n_2$,  and $n_2 > n_1$.

\end{prob}

More generally still:

\begin{prob}

Evaluate $\tau\hat{\;}(G, \kappa)$ for each noncyclic group $G$ of exponent $\kappa$.

\end{prob}

\subsection{Limited number of terms} \label{5maxRlimited}

In this section we investigate, for a given group $G$ and positive integer $t$, the quantity $$\tau\hat{\;} (G,[1,t]) = \mathrm{max} \{ |A|  \mid A \subseteq G, 0 \not \in [1,t]\hat{\;}A\},$$ that is, the maximum size of a weak zero-$[1,t]$-sum-free subset of $G$.

As before, we can easily verify the following:

\begin{prop} \label{tau hat t=1,2}  In a group of order $n$, we have
$$\tau\hat{\;}(G,[1,1])=n-1$$
and
$$\tau\hat{\;}(G,[1,2])=\frac{n+|\mathrm{Ord}(G,2)|-1}{2};$$  
as a special case, for the cyclic group of order $n$ we have 
 $$\tau\hat{\;}(\mathbb{Z}_n,[1,2])=\left\lfloor \frac{n}{2} \right\rfloor.$$

\end{prop}

We have the following general upper bound for $\tau\hat{\;} (G,[1,t])$:

\begin{prop}  \label{tau hat t vs [1,t]}
For every $G$ and $t \geq 2$, we have
$$\tau\hat{\;} (G,[1,t]) \leq \tau\hat{\;} (G,t)-1.$$
\end{prop}

Indeed, if $A$ is a weak zero-$[1,t]$-sum-free in $G$ of size $\tau\hat{\;} (G,[1,t])$, then $0 \not \in A$, so $A \cup \{0\}$ has size $\tau\hat{\;} (G,[1,t])+1$.  Furthermore, 
$0 \not \in (t-1) \hat{\;} A$ and $0 \not \in t \hat{\;} A$, which implies that $ 0 \not \in t \hat{\;} (A \cup \{0\})$.  This proves that $\tau\hat{\;} (G,t)$ is at least one more than $\tau\hat{\;} (G,[1,t])$, as claimed.  As Propositions \ref{tau hat h=1,2} and \ref{tau hat t=1,2} show, we have equality in Proposition \ref{tau hat t vs [1,t]} for $t=2$.

Let us now consider $t \geq 3$.  Clearly, if $m$ is a positive integer for which $$m+(m-1)+\cdots+(m-t+1) \leq n-1,$$ then the set
$$A=\{1,2,\dots,m\}$$ is a weak zero-$[1,t]$-sum-free subset of the cyclic group $\mathbb{Z}_n$; this gives:

\begin{prop} \label{lower for tau hat [1,t] cyclic}
For all positive integers $n$ and $t$ we have $$\tau\hat{\;}(\mathbb{Z}_n,[1,t]) \geq \left \lfloor \frac{n-1}{t}+\frac{t-1}{2} \right \rfloor.$$

\end{prop} 
(Observe that, by Proposition \ref{tau hat t=1,2},  equality holds for $t=1$ and $t=2$.)  

We are able to do better in certain cases.  Consider, for example, the case when $n-2$ is divisible by 6, and define the subset $A$ of $\mathbb{Z}_n$ as $$A=\{1,2,\dots, (n-2)/6\} \cup \{n/2, n/2+1, \dots, (2n+2)/3\}.$$ (Note that when $n-2$ is divisible by 6, then $n$ is even and $2n+2$ is divisible by 6.)  An easy computation shows that $0 \not \in [1,3] \hat{\;}A$; furthermore, $A$ contains $(n+4)/3$ elements.  Therefore, we get:

\begin{prop} \label{two part set for tau hat [1,3]}
If $n-2$ is divisible by 6, then $$\tau\hat{\;}(\mathbb{Z}_n,[1,3]) \geq (n+4)/3.$$

\end{prop}

In a similar manner, Yin proved that the set $$A=\{1,2,\dots, n/6\} \cup \{n/2, n/2+1, \dots, 2n/3\}$$ is a weak zero-$[1,3]$-sum-free set in $\mathbb{Z}_n$ when $n$ is divisible by 6, and that the set
$$\{1,2,\dots,\left \lfloor (n+2)/8 \right \rfloor\} \cup \{ n/2, n/2+1, \dots, n/2+ \left \lfloor (n+6)/8 \right \rfloor\}$$ is a weak zero-$[1,4]$-sum-free set in $\mathbb{Z}_n$ when $n-2$ is divisible by 4.  These results yield:

\begin{prop} [Yin; cf.~\cite{Yin:2016a}] \label{Yin two part set for tau hat [1,3]}\index{Yin, Y.} 
\begin{enumerate} \item If $n$ is divisible by 6, then $$\tau\hat{\;}(\mathbb{Z}_n,[1,3]) \geq n/3+1.$$
\item If $n-2$ is divisible by 4, then $$\tau\hat{\;}(\mathbb{Z}_n,[1,4]) \geq (n+6)/4.$$
\end{enumerate}
\end{prop}

We should note that Proposition \ref{two part set for tau hat [1,3]} and both parts of Proposition \ref{Yin two part set for tau hat [1,3]} yield  values one larger than Proposition \ref{lower for tau hat [1,t] cyclic} does.

For an upper bound, we have the following conjecture of Hamidoune (see Conjecture 1.1 in \cite{Ham:1998a}):\index{Hamidoune, Y. O.} 
\begin{conj} [Hamidoune; cf.~\cite{Ham:1998a}] \label{Hamidoune conj tau hat [1,t]}\index{Hamidoune, Y. O.} 
For all integers $n \geq t$ we have $$\tau\hat{\;}(\mathbb{Z}_n,[1,t]) \leq (n-2)/t + t-1.$$
\end{conj}
In particular, for $t=3$ we believe that $$\tau\hat{\;}(\mathbb{Z}_n,[1,3]) \leq (n+4)/3,$$ and Proposition \ref{two part set for tau hat [1,3]} above shows that this inequality, if true, is sharp.

It is a very interesting question to see if Conjecture \ref{Hamidoune conj tau hat [1,t]} is sharp in general:

\begin{prob}
Find all pairs of integers $t$ and $n$ for which a weak zero-$[1,t]$-sum-free subset of size  $\lfloor (n-2)/t \rfloor + t-1 $  exists in
the cyclic group $\mathbb{Z}_n$.

\end{prob}

A less exact upper bound was established by Alon in \cite{Alo:1987a}:\index{Alon, N.} if $n$ is large enough (depending on $t$ and any fixed  positive real number $\epsilon$), then   
$$\tau\hat{\;}(\mathbb{Z}_n,[1,t]) \leq (1/t + \epsilon) \cdot n.$$  Lev\index{Lev, V. F.} generalized this to arbitrary finite abelian groups as follows:

\begin{thm} [Lev; cf.~\cite{Lev:2004a}] \index{Lev, V. F.} 
For any $t$ and positive real number $\epsilon$, there is a constant $n_0(t, \epsilon)$ so that $$\tau\hat{\;}(G,[1,t]) \leq (1/t + \epsilon) \cdot n$$ holds for every $n \geq n_0(t, \epsilon) \cdot |\mathrm{Ord}(G,2)|$.

\end{thm}

We have very few exact values known for $\tau\hat{\;}(G,[1,t])$, particularly when $G$ is not cyclic.  We mention the following result:

\begin{prop} [Yin; cf.~\cite{Yin:2016a}] \label{Yin Z2}\index{Yin, Y.} 
For all positive integers $r$, we have $\tau\hat{\;}(\mathbb{Z}_2^r,[1,3])=2^{r-1}.$

\end{prop}
We present the short proof on page \pageref{proof of Yin Z2}.

For the group $\mathbb{Z}_3^3$, Bhowmik and Schlage-Puchta computed\index{Bhowmik, G.}\index{Schlage-Puchta, J-C.} in  \cite{BhoSch:2007a} three values: 
$$\tau\hat{\;}(\mathbb{Z}_3^3,[1,3])=8, \; \tau\hat{\;}(\mathbb{Z}_3^3,[1,4])=7, \; \tau\hat{\;}(\mathbb{Z}_3^3,[1,5])=7.$$
In \cite{BhoSch:2009a}, the same authors prove that any weak zero-$[1,3]$-free subset of size 8 in $\mathbb{Z}_3^3$ is of the form
$$\{a_1, a_2, a_3, a_1+a_2, a_1+a_2+a_3, a_1+2a_2+a_3, 2a_1+a_3, a_2+2a_3\}.$$

The general problem of finding the exact value of $\tau\hat{\;}(G,[1,t])$ is still wide open:

\begin{prob}
Find $\tau\hat{\;}(G,[1,t])$ for all groups $G$ and integers $t \geq 3$.
\end{prob}

\subsection{Arbitrary number of terms} \label{5maxRarbitrary}

Here we investigate the maximum value of $m$ for which there exists an $m$-subset $A$ of $G$ that is {\em zero-sum-free} in $G$; that is, for which its sumset $$ \Sigma^*   A=\cup_{h=1}^{\infty} h \hat{\;} A$$ does not contain zero: 
$$\tau \hat{\;}(G, \mathbb{N}) =\mathrm{max} \{ |A|  \mid A \subseteq G, 0 \not \in \Sigma^*   A \}.$$  The quantity $\tau \hat{\;}(G, \mathbb{N})+1$, that is, the smallest integer for which for any subset $A$ of $G$ of that size (or more) one has $0 \in  \Sigma^*   A$, is called the {\em Olson's constant} of $G$, and has been the subject of much attention since Erd\H{o}s and Heilbronn\index{Erd\H{o}s, P.}\index{Heilbronn, H.} first discussed it in \cite{ErdHei:1964a} in 1964.  (The term  ``Olson's constant'' was coined at a 1994 meeting in Venezuela to honor Olson\index{Olson, J. E.} who contributed tremendously to this and other questions in additive combinatorics.) 

Let us consider the cyclic group first.  It is easy to see that if $1+2+\cdots+m <n$, then the set $\{1,2,\dots,m\}$ is zero-sum-free in $\mathbb{Z}_n$, and therefore we get:
\begin{prop} \label{zero-sum-free simple}
For all positive integers $n$ we have
$$\tau \hat{\;}(\mathbb{Z}_n, \mathbb{N}) \geq \left \lfloor (\sqrt{8n-7} -1)/2 \right \rfloor.$$
\end{prop}

Selfridge (see Problem C.15 in \cite{Guy:2004a}; cf.~also page 95 in \cite{ErdGra:1980a}) offered two other constructions that are sometimes better than the one above. \label{Selfridge constructions}  First, observe that there is no harm replacing the element 2 with $-2$, so when $1+3+4+\cdots+m<n$, then    
the set $\{1,-2,3,\dots,m\}$ is zero-sum-free in $\mathbb{Z}_n$.  (For $-2 \not \in \{1,3,4,\dots,m\}$ we must have $n \geq 6$.)  We thus get:
\begin{prop} [Selfridge; cf.~Problem C.15 in \cite{Guy:2004a}] \label{zero-sum-free simple -2}\index{Guy, R.} 
For all positive integers $n \geq 6$ we have
$$\tau \hat{\;}(\mathbb{Z}_n, \mathbb{N}) \geq \left \lfloor (\sqrt{8n+9} -1)/2 \right \rfloor.$$
\end{prop}

The other construction of Selfridge assumes that $n$ is even, and considers the set 
$$A=\left\{
\begin{array}{cl}
\{1,2,\dots,(m-1)/2\} \cup \{n/2, n/2+1, \dots, n/2+(m-1)/2\} & \mbox{if $m$ is odd}; \\ \\
\{1,2,\dots,m/2-1\} \cup \{n/2, n/2+1, \dots, n/2+m/2\} & \mbox{if $m$ is even}.
\end{array} \right.$$  
(When $m=1$, we simply take $A=\{n/2\}$.)

Clearly, $A$ is zero-sum-free in $\mathbb{Z}_n$ when $m$ is odd and $$2 \cdot (1+2+\cdots+(m-1)/2) < n/2,$$ or when $m$ is even and 
$$2 \cdot (1+2+\cdots+(m/2-1))+m/2 < n/2.$$  Therefore, we can say that 
$$\tau \hat{\;}(\mathbb{Z}_n, \mathbb{N}) \geq
\left\{
\begin{array}{cl}
\left \lfloor \sqrt{2n-3} \right \rfloor & \mbox{if this value is odd}; \\ \\
\left \lfloor \sqrt{2n-4} \right \rfloor & \mbox{if this value is even}.
\end{array} \right.$$
(Note that these two conditions are not mutually exclusive, but, as the next few lines demonstrate, they do cover all cases.)  
Observe that $\left \lfloor \sqrt{2n-3} \right \rfloor$ and $\left \lfloor \sqrt{2n-4} \right \rfloor$ are only unequal when $2n-3$ happens to be a square number, but when that is the case, then it must be odd (and thus $\left \lfloor \sqrt{2n-3} \right \rfloor$ is odd as well).  Therefore, we can conclude:
\begin{prop} [Selfridge; cf.~Problem C.15 in \cite{Guy:2004a}] \label{zero-sum-free simple n even}\index{Guy, R.} 
For all positive even integers $n$ we have
$$\tau \hat{\;}(\mathbb{Z}_n, \mathbb{N}) \geq \left \lfloor \sqrt{2n-3} \right \rfloor.$$
\end{prop}

We can also verify that Proposition \ref{zero-sum-free simple -2} never trumps Proposition \ref{zero-sum-free simple n even}; indeed, for all $n \geq 14$, we have
$$\sqrt{2n-3}  \geq (\sqrt{8n+9} -1)/2.$$  (For $n \in \{6,8,10,12\}$ the two floors are equal.)  We are not aware of any even $n$ for which Proposition \ref{zero-sum-free simple n even} is not exact, and we have:
\begin{conj} [Selfridge; cf.~Problem C.15 in \cite{Guy:2004a}]  \label{conj zero-sum-free simple n even}\index{Guy, R.}
For all positive even integers $n$ we have
$$\tau \hat{\;}(\mathbb{Z}_n, \mathbb{N}) = \left \lfloor \sqrt{2n-3} \right \rfloor.$$
\end{conj}
Subocz in \cite{Sub:2000a}\index{Subocz G., J. C.} verified Conjecture \ref{conj zero-sum-free simple n even} for all even $n \leq 64$.  

\begin{prob}
Prove (or disprove) Conjecture \ref{conj zero-sum-free simple n even}.

\end{prob}

The case when $n$ is odd seems more complicated.  From Subocz's paper \cite{Sub:2000a}\index{Subocz G., J. C.}  we see that Proposition \ref{zero-sum-free simple -2} gives the correct value of $\tau \hat{\;}(\mathbb{Z}_n, \mathbb{N})$ for all odd $n$ with $7 \leq n \leq 63$, except for $n=25$; that exception is demonstrated by the fact that
$$\{1,6,11,16,21\} \cup \{5,10\}$$ is zero-sum-free in $\mathbb{Z}_{25}$.  

After a variety of partial results, the case when $n$ is prime is now settled: As Balandraud\index{Balandraud, \'E.} pointed out in \cite{Bal:2012a} (cf.~\cite{Bal:2012b}) and as we now explain, for a prime $p$, the value of $\tau \hat{\;}(\mathbb{Z}_p, \mathbb{N})$ easily follows from his Theorem \ref{Balandraud}.  The cases of $p \in \{2,3,5\}$ can be easily evaluated, so we assume that $p \geq 7$.  Given Proposition \ref{zero-sum-free simple -2}, we only need to prove that if $A \subseteq \mathbb{Z}_p$ has size $$|A| \geq \left \lfloor (\sqrt{8p+9} -1)/2 \right \rfloor + 1 > (\sqrt{8p+9} -1)/2,$$ then $A$ is not zero-sum-free in $\mathbb{Z}_p$.

Clearly, if $A$ is not asymmetric, that is, if $A$ and $-A$ are not disjoint, then $0 \in 2 \hat{\;} A$, and $A$ is not zero-sum-free.  On the other hand, if $A$ is asymmetric, then $|A| \leq (p-1)/2$, and by Theorem \ref{Balandraud}, $\Sigma^*   A$ has size at least $$\min \left \{p, |A| \cdot (|A|+1)/2 \right\} \geq \min \left \{p, (\sqrt{8p+9} -1) \cdot (\sqrt{8p+9} +1)/8 \right\}= \min \{p, p+1\}=p.$$ But this means that $\Sigma^*   A=\mathbb{Z}_p$, thus $A$ is not zero-sum-free. 

\begin{thm} [Balandraud; cf.~\cite{Bal:2012a}, \cite{Bal:2012b}]\index{Balandraud, \'E.}
We have $\tau \hat{\;}(\mathbb{Z}_2, \mathbb{N})=1$, $\tau \hat{\;}(\mathbb{Z}_3, \mathbb{N})=1$, and $\tau \hat{\;}(\mathbb{Z}_5, \mathbb{N})=2$; furthermore, for every prime $p \geq 7$, we have
$$\tau \hat{\;}(\mathbb{Z}_p, \mathbb{N})=\left \lfloor (\sqrt{8p+9} -1)/2 \right \rfloor.$$

\end{thm}
Perhaps it is worth pointing out that the $+9$ in the expression above can be omitted, since for every prime $p \geq 7$,  
$$\left \lfloor (\sqrt{8p+9} -1)/2 \right \rfloor= \left \lfloor \sqrt{2p} -1/2 \right \rfloor.$$  To see this, note that for $(\sqrt{8p+c} -1)/2$ to equal an integer $k$, one must have $c \equiv 1$ mod 8, and if $(\sqrt{8p+9} -1)/2$  or $(\sqrt{8p+1} -1)/2$ equals $k$, then $p$ equals $(k+2)(k-1)/2$ or $k(k+1)/2$, but neither of these quantities can equal a prime $p \geq 7$.  Dropping the $+9$ from our formula also has the advantage that the cases of $p \in \{2,3,5\}$ do not need to be separated, thus we get:

\begin{cor} [Balandraud; cf.~\cite{Bal:2012a}, \cite{Bal:2012b}] \label{Balandraud coroll}\index{Balandraud, \'E.}
For every prime $p $, we have
$$\tau \hat{\;}(\mathbb{Z}_p, \mathbb{N})=\left \lfloor \sqrt{2p} -1/2 \right \rfloor.$$

\end{cor}

The question of finding the exact value of $\tau \hat{\;}(\mathbb{Z}_n, \mathbb{N})$ for composite $n$ remains open:

\begin{prob}
Evaluate $\tau \hat{\;}(\mathbb{Z}_n, \mathbb{N})$ for odd composite $n$.  In particular, find all instances when 
$$\tau \hat{\;}(\mathbb{Z}_n, \mathbb{N}) > \left \lfloor (\sqrt{8n+9} -1)/2 \right \rfloor.$$
\end{prob}

Let us now move on to noncyclic groups.  First, following Gao, Ruzsa, and Thangadurai (see \cite{GaoRuzTha:2004a}),\index{Ruzsa, I.}\index{Gao, W.}\index{Thangadurai, R.}  we establish a recursive lower bound for $\tau \hat{\;}(G, \mathbb{N})$.  Suppose that $G$ has exponent $\kappa$ and $G \cong G_1 \times \mathbb{Z}_{\kappa}$.  Let $A_1$ be a zero-sum-free set in $G_1$ and, if possible, let $A_2$ be any subset of $G_1$ of size $\kappa -1$; if $|G_1| < \kappa -1$, then let $A_2=G_1$.  It is then easy to see that
$$A=\left(A_1  \times \{0\} \right) \cup \left( A_2 \times  \{1\} \right)$$ is zero-sum-free in $G_1 \times \mathbb{Z}_{\kappa}$.  Therefore:

\begin{prop}
For every finite abelian group $G_1$ and positive integer $\kappa$ we have 
$$\tau \hat{\;}(G_1 \times \mathbb{Z}_{\kappa}, \mathbb{N}) \geq \tau \hat{\;}(G_1, \mathbb{N})+ \min\{|G_1|, \kappa -1\}.$$
\end{prop}

Consequently, we have the following lower bound for groups of the form $\mathbb{Z}_k^r$:

\begin{cor} [Gao, Ruzsa, and Thangadurai; cf.~\cite{GaoRuzTha:2004a}] \label{Gao, Ruzsa, and Thangadurai recursive}\index{Ruzsa, I.}\index{Gao, W.}\index{Thangadurai, R.} 
For all $k \geq 2$ and $r \geq 2$, we have $$\tau \hat{\;}(\mathbb{Z}_k^r, \mathbb{N}) \geq \tau \hat{\;}(\mathbb{Z}_k^{r-1}, \mathbb{N}) +k-1.$$
\end{cor}

Furthermore, Gao, Ruzsa,\index{Ruzsa, I.}\index{Gao, W.} and Thangadurai\index{Thangadurai, R.}  in \cite{GaoRuzTha:2004a} proved that when $k$ is a very large prime and $r=2$, then equality holds in Corollary \ref{Gao, Ruzsa, and Thangadurai recursive}; the requirement on the prime being large was then greatly reduced by Bhowmik and Schlage-Puchta in \cite{BhoSch:2010a},\index{Bhowmik, G.} though it is still formidable:

\begin{thm} [Bhowmik and Schlage-Puchta; cf.~\cite{BhoSch:2010a}] \label{Bhowmik and Sclage-Puchta p^2}\index{Bhowmik, G.}
For every prime $p > 6000$, we have $$\tau \hat{\;}(\mathbb{Z}_p^2, \mathbb{N}) = \tau \hat{\;}(\mathbb{Z}_p, \mathbb{N}) +p-1.$$
\end{thm} 

Combining this with Corollary \ref{Balandraud coroll} yields:

\begin{cor} [Bhowmik and Schlage-Puchta; cf.~\cite{BhoSch:2010a}] \label{Bhowmik and Sclage-Puchta p^2 cor}\index{Bhowmik, G.}
For every prime $p > 6000$, we have $$\tau \hat{\;}(\mathbb{Z}_p^2, \mathbb{N}) = p+\left \lfloor \sqrt{2p} -1/2 \right \rfloor-1.$$

\end{cor}

The obvious question here is whether the bound on $p$ can be reduced:

\begin{prob}
Prove that the conclusion of Corollary \ref{Bhowmik and Sclage-Puchta p^2 cor} (or, equivalently, of Theorem \ref{Bhowmik and Sclage-Puchta p^2}) holds for primes $p <6000$.
\end{prob}

In fact, we believe that the same equality holds for composite values of $k$ as well:

\begin{conj} [Gao, Ruzsa, and Thangadurai; cf.~\cite{GaoRuzTha:2004a}] \label{Gao, Ruzsa, and Thangadurai conj}\index{Ruzsa, I.}\index{Gao, W.}\index{Thangadurai, R.} 
For all values of $k \geq 2$, we have $$\tau \hat{\;}(\mathbb{Z}_k^2, \mathbb{N}) = \tau \hat{\;}(\mathbb{Z}_k, \mathbb{N}) +k-1.$$

\end{conj}

\begin{prob}
Prove (or disprove) Conjecture \ref{Gao, Ruzsa, and Thangadurai conj}.\index{Ruzsa, I.}\index{Gao, W.}\index{Thangadurai, R.} 
\end{prob}

The paper \cite{GaoRuzTha:2004a} actually had a stronger version of Conjecture \ref{Gao, Ruzsa, and Thangadurai conj}:\index{Ruzsa, I.}\index{Gao, W.}\index{Thangadurai, R.} the authors conjectured that equality always holds in Corollary \ref{Gao, Ruzsa, and Thangadurai recursive}.\index{Ruzsa, I.}\index{Gao, W.}\index{Thangadurai, R.} However, as it was first pointed out by Gao and Geroldinger\index{Gao, W.}\index{Geroldinger, A.} in \cite{GaoGer:2006a}, for every prime power $k \geq 3$ there exists an $r \geq 2$ for which the claim fails.  Furthermore, as Ordaz et al.~proved\index{Ordaz, O.} in \cite{OrdEtAl:2011a} and as we review below, even for $r=3$, equality in Corollary \ref{Gao, Ruzsa, and Thangadurai recursive}\index{Ruzsa, I.}\index{Gao, W.}\index{Thangadurai, R.} can only hold for at most finitely many prime values of $k$ (possibly only for $k=2$).

Ordaz at al. proved lower bounds for $\tau \hat{\;}(\mathbb{Z}_k^r, \mathbb{N})$ for all $r \geq 3$:

\begin{thm} [Ordaz, Phillipp, Santos, and Schmid; cf.~\cite{OrdEtAl:2011a}] \label{Ordaz et al r=3} \index{Ordaz, O.}\index{Phillipp, A.},\index{Santos, I.}\index{Schmid, W.} 
For all integers $k \geq 2$ we have
$$\tau \hat{\;}(\mathbb{Z}_k^3, \mathbb{N}) \geq 
\left\{
\begin{array}{cl}
3(k-1) & \mbox{if $k \leq 3$}; \\ \\
2k+ \left \lfloor \left(\sqrt{8k-31} -1 \right)/2 \right \rfloor & \mbox{if $k \geq 4$}.
\end{array}
\right.$$
\end{thm}

\begin{thm} [Ordaz, Phillipp, Santos, and Schmid; cf.~\cite{OrdEtAl:2011a}] \label{Ordaz et al r>=4} \index{Ordaz, O.}\index{Phillipp, A.},\index{Santos, I.}\index{Schmid, W.} 
For all integers $k \geq 2$ and $r \geq 4$ we have
$$\tau \hat{\;}(\mathbb{Z}_k^r, \mathbb{N}) \geq 
\left\{
\begin{array}{cl}
r(k-1) & \mbox{if $k \leq r+1$}; \\ \\
(r-1)k + \left \lfloor \left(\sqrt{8(k-r)+1} -1 \right) \left/2 \right. \right \rfloor & \mbox{if $k \geq r+2$}.
\end{array}
\right.$$

\end{thm}
(Observe that there are slight discrepancies between the formulae as well as the conditions of Theorems \ref{Ordaz et al r=3} and \ref{Ordaz et al r>=4}.)  

Combining Theorem \ref{Ordaz et al r=3} and Corollary \ref{Bhowmik and Sclage-Puchta p^2 cor}, for a prime $p > 6000$ we thus get:
$$\tau \hat{\;}(\mathbb{Z}_p^3, \mathbb{N}) \geq 2p+ \left \lfloor \left(\sqrt{8p-31} -1 \right)/2 \right \rfloor > 2p+ \left \lfloor \sqrt{2p} -1/2 \right \rfloor-2 =\tau \hat{\;}(\mathbb{Z}_p^2, \mathbb{N})+p-1,$$ which shows that equality in Corollary \ref{Gao, Ruzsa, and Thangadurai recursive}\index{Ruzsa, I.}\index{Gao, W.}\index{Thangadurai, R.} does not hold for primes $p > 6000$.  This disproves the conjecture in \cite{GaoRuzTha:2004a} that we mentioned above.

The authors of \cite{OrdEtAl:2011a} also believe (though shy away from a conjecture) that equality holds everywhere in Theorems \ref{Ordaz et al r=3} and \ref{Ordaz et al r>=4}; we thus have the following interesting question:

\begin{prob}
Decide whether equality always holds in Theorems  \ref{Ordaz et al r=3} and \ref{Ordaz et al r>=4}.
\end{prob}

We do know that equality holds in Theorems \ref{Ordaz et al r=3} and \ref{Ordaz et al r>=4} for $k \in \{2,3,4,5\}$ and in Theorem \ref{Ordaz et al r=3} when $k \in \{6,7\}$: 

\begin{thm} [Subocz; cf.~\cite{Sub:2000a}] \label{tau hat elementary k=2}\index{Subocz G., J. C.} 
For every $r \geq 3$ we have $\tau \hat{\;}(\mathbb{Z}_2^r, \mathbb{N})=r$ and $\tau \hat{\;}(\mathbb{Z}_3^r, \mathbb{N})=2r$.
\end{thm}

\begin{thm} [Ordaz, Phillipp, Santos, and Schmid; cf.~\cite{OrdEtAl:2011a}] \label{Ordaz et al z2 by z4} \index{Ordaz, O.}\index{Phillipp, A.},\index{Santos, I.}\index{Schmid, W.}  
For every $r \geq 4$ we have $\tau \hat{\;}(\mathbb{Z}_4^r, \mathbb{N})=3r$ and $\tau \hat{\;}(\mathbb{Z}_5^r, \mathbb{N})=4r$.  Furthermore, $\tau \hat{\;}(\mathbb{Z}_4^3, \mathbb{N})=8$, $\tau \hat{\;}(\mathbb{Z}_5^3, \mathbb{N})=11$, $\tau \hat{\;}(\mathbb{Z}_6^3, \mathbb{N})=13$, and $\tau \hat{\;}(\mathbb{Z}_7^3, \mathbb{N})=16$.
\end{thm}

We also have the following upper bound:

\begin{thm} [Ordaz, Phillipp, Santos, and Schmid; cf.~\cite{OrdEtAl:2011a}] \label{prime power Ordaz} \index{Ordaz, O.}\index{Phillipp, A.},\index{Santos, I.}\index{Schmid, W.} 
For every $r \geq 1$ and prime power $q$, $$\tau \hat{\;}(\mathbb{Z}_q^r, \mathbb{N}) \leq r(q-1).$$

\end{thm}

Combining Theorems \ref{prime power Ordaz} and \ref{Ordaz et al r>=4} yields:

\begin{thm} [Ordaz, Phillipp, Santos, and Schmid; cf.~\cite{OrdEtAl:2011a}] \label{Ordaz et al z2 by z4} \index{Ordaz, O.}\index{Phillipp, A.},\index{Santos, I.}\index{Schmid, W.}
For every $r \geq 4$ and prime power $q$ for which $q \leq r+1$, $$\tau \hat{\;}(\mathbb{Z}_q^r, \mathbb{N}) = r(q-1).$$

\end{thm}

We mention some other exact results:

\begin{thm} [Ordaz, Phillipp, Santos, and Schmid; cf.~\cite{OrdEtAl:2011a}] \label{Ordaz et al z2 by z4} \index{Ordaz, O.}\index{Phillipp, A.},\index{Santos, I.}\index{Schmid, W.} 
For every $r_1, r_2 \geq 1$, we have $$\tau \hat{\;}(\mathbb{Z}_2^{r_1} \times \mathbb{Z}_4^{r_2}, \mathbb{N}) = r_1+3r_2.$$

\end{thm}

\begin{thm} [Ordaz, Phillipp, Santos, and Schmid; cf.~\cite{OrdEtAl:2011a}] \label{Ordaz et al p group} \index{Ordaz, O.}\index{Phillipp, A.},\index{Santos, I.}\index{Schmid, W.} 
Suppose that $G$ has invariant factorization
$$\mathbb{Z}_{p}^{r_1} \times \mathbb{Z}_{p^2}^{r_2} \times \cdots \times \mathbb{Z}_{p^k}^{r_k}$$ where $p$ is prime, $r_1, \dots, r_k \geq 0$, $r_k >0$, and $r_1+ \cdots + r_k \geq p^k$.  Then we have $$\tau \hat{\;}(G, \mathbb{N}) = r_1 (p-1) + r_2 (p^{2}-1) + \cdots + r_k (p^{k}-1).$$

\end{thm}
(This last result on $\tau \hat{\;}(G, \mathbb{N})$ for $p$-groups, under stronger assumptions on the rank of $G$, was given by Gao and Geroldinger\index{Gao, W.}\index{Geroldinger, A.} as Corollary 7.4 in \cite{GaoGer:1999a}.)

Note that by Theorem \ref{Ordaz et al p group}, it suffices to determine $\tau \hat{\;}(\mathbb{Z}_3^{r_1} \times \mathbb{Z}_9^{r_2}, \mathbb{N})$ for $r_1+r_2 \leq 8$:

\begin{prob}
Evaluate $\tau \hat{\;}(\mathbb{Z}_3^{r_1} \times \mathbb{Z}_9^{r_2}, \mathbb{N})$ for all $r_1, r_2 \geq 1$ with $r_1+r_2 \leq 8$.

\end{prob}

\begin{prob}
Evaluate $\tau \hat{\;}(\mathbb{Z}_6^r, \mathbb{N})$ for all $r \geq 4$.

\end{prob}

Turning now to general finite abelian groups, we first pose a 1973 conjecture of Erd\H{o}s:\index{Erd\H{o}s, P.}

\begin{conj} [Erd\H{o}s; cf.~\cite{Erd:1973a}] \label{Erdos on tau hat}\index{Erd\H{o}s, P.}
For every finite abelian group of order $n$, we have $$\tau \hat{\;}(G, \mathbb{N}) < \sqrt{2n} .$$

\end{conj}

According to Corollary \ref{Bhowmik and Sclage-Puchta p^2 cor}, Conjecture \ref{Erdos on tau hat} holds for all prime values of $n > 6000$; furthermore, we see from Subocz's\index{Subocz G., J. C.} tables in \cite{Sub:2000a} that the conjecture holds for all groups of order $n \leq 50$ and all cyclic groups of order $n \leq 64$.

In 1975, Olson proved the following:

\begin{thm} [Olson; cf.~\cite{Ols:1975a}] \label{Olson 1975 bound}\index{Olson, J. E.} 
For every $G$ we have  
$$\tau \hat{\;}(G, \mathbb{N}) < 3 \sqrt{n} .$$
\end{thm}

The best current general result is the following:

\begin{thm} [Hamidoune and Z\'emor; cf.~\cite{HamZem:1996a}] \label{Ham Zem bound}\index{Hamidoune, Y. O.}\index{Zemor@Z\'emor, G.}
There exists a positive real number $C$ for which 
$$\tau \hat{\;}(G, \mathbb{N}) < \sqrt{2n}  + C \cdot \sqrt[3]{n}\ln n$$ holds for every finite abelian group of order $n$.
\end{thm}

We also mention the following interesting conjecture:

\begin{conj} [Subocz; cf.~\cite{Sub:2000a}] \label{Subocz on tau hat}\index{Subocz G., J. C.} 
For every finite abelian group of order $n$, we have $$\tau \hat{\;}(G, \mathbb{N}) \leq \tau \hat{\;}(\mathbb{Z}_n, \mathbb{N}).$$ 

\end{conj}

We should note that the difference between $\tau \hat{\;}(G, \mathbb{N})$ (for $|G|=n$) and $\tau \hat{\;}(\mathbb{Z}_n, \mathbb{N})$ can be arbitrarily large: For example, we have $\tau \hat{\;}(\mathbb{Z}_2^r, \mathbb{N})=r$ (see Theorem \ref{tau hat elementary k=2}) but $\tau \hat{\;}(\mathbb{Z}_{2^r}, \mathbb{N})  \sim 2^{(r+1)/2}$ (see Proposition \ref{zero-sum-free simple n even}).

In closing this section, we mention that results regarding the structure of zero-sum-free subsets $A$ of $\mathbb{Z}_p$ (with $p$ prime) for which $|A|$ is close to $\tau \hat{\;}(\mathbb{Z}_p, \mathbb{N})$ are discussed by Deshouillers and Prakash\index{Prakash, G.}\index{Deshouillers, J--M.} in \cite{DesPra:2011a}; by Nguyen, Szemer\'edi, and Vu\index{Nguyen, N. H.}\index{Vu, V. H.}\index{Szemer\'edi, E.} in \cite{NguSzeVu:2008a}; and by Nguyen and Vu in \cite{NguVu:2009a}.  For example, the analogue of Theorem \ref{Nguyen  and Vu} is:

\begin{thm} [Nguyen  and Vu; cf.~\cite{NguVu:2009a}] \label{Nguyen  and Vu zero-free}\index{Nguyen, N. H.}\index{Vu, V. H.} 
There is a positive constant $c$, so that whenever $A$ is a zero-sum-free subset of $\mathbb{Z}_p$ for an odd prime $p$, then there is a subset $A'$ of $A$ of size at most $c p^{6/13} \log p$ and an element $b \in \{1,\dots,p-1\}$, for which $$||b \cdot (A \setminus A')|| < p.$$

\end{thm}

Zero-sum-free subsets of $\mathbb{Z}_p^2$ (with $p>6000$ prime) that have size exactly $\tau \hat{\;}(\mathbb{Z}_p^2, \mathbb{N})$ are classified by Bhowmik and Schlage-Puchta in \cite{BhoSch:2010a} and \cite{BhoSch:2011a}.\index{Bhowmik, G.}\index{Schlage-Puchta, J-C.}

We offer the following difficult problems:

\begin{prob}
Prove Conjecture \ref{Erdos on tau hat}, or at least improve on the bound in Theorem \ref{Ham Zem bound}.

\end{prob}

\begin{prob}
for each group $G$, classify all zero-sum-free sets in $G$ of size exactly $\tau \hat{\;}(G, \mathbb{N})$.

\end{prob}

In particular, one may be able to improve on Theorem \ref{Nguyen  and Vu zero-free}:

\begin{prob}
For any odd prime $p$, classify all zero-sum-free sets in $\mathbb{Z}_p$ of size exactly $\tau \hat{\;}(\mathbb{Z}_p, \mathbb{N})=\left \lfloor \sqrt{2p} -1/2 \right \rfloor$.

\end{prob}

\section{Restricted signed sumsets} \label{5maxRS}

In this section we investigate the quantity $$\tau \hat{_{\pm}}(G,H) =\mathrm{max} \{ |A|  \mid A \subseteq G, 0 \not \in H \hat{_{\pm}} A \}$$  (if there is no subset $A$ for which $ 0 \not \in H \hat{_{\pm}} A$, we let $\tau \hat{_{\pm}} (G,H) =0$).  Clearly, we always have $\tau \hat{_{\pm}} (G,H)=0$ when $0 \in H$; however, when $ 0 \not \in H$ and $n \geq 2$, then $\tau \hat{_{\pm}} (G,H) \geq 1$: for any $a \in G \setminus \{0\}$, for the one-element set $A=\{a\}$ we obviously have $0 \not \in H \hat{_{\pm}} A$.

We consider three special cases: when $H$ consists of a single positive integer $h$,  when $H$ consists of all positive integers up to some value $t$, and when $H = \mathbb{N}$.

\subsection{Fixed number of terms} \label{5maxRSfixed}

Here we investigate, for a given group $G$ and positive integer $h$, the quantity $$\tau \hat{_{\pm}} (G,h) = \mathrm{max} \{ |A|  \mid A \subseteq G, 0 \not \in h \hat{_{\pm}}A\}.$$

As before, we can easily verify the following:

\begin{prop} \label{tau hat pm h=1,2}  In a group of order $n$, we have
$$\tau \hat{_{\pm}}(G,1)=n-1;$$
$$\tau \hat{_{\pm}}(G,2)=\frac{n+|\mathrm{Ord}(G,2)|+1}{2};$$ and  
$$\tau \hat{_{\pm}}(G,h)=n$$ for every $h \geq n+1$. 
\end{prop} 

Next, we evaluate $\tau \hat{_{\pm}}(G,n)$, as follows.  Note that, trivially, $\tau \hat{_{\pm}}(G,n) \geq n-1$.  Recall also that on page \pageref{tau hat in direct prod} we found that when $G$ has exponent $\kappa$ and $G \cong G_1 \times \mathbb{Z}_{\kappa}$, then the elements of $G_1 \times \mathbb{Z}_{\kappa}$ sum to $(0,0)$ if $|G_1|= n/ \kappa$ is even or if $\kappa$ is odd, and to $(0, \kappa/2)$ when $n/ \kappa$ is odd and $\kappa$ is even.  We can thus also observe that when $\kappa$ is divisible by 4 and $n/ \kappa$ is odd, then adding all elements of $G_1 \times \mathbb{Z}_{\kappa}$, except for the element $(0, \kappa/4)$ that we instead subtract, we get $(0,0)$.  That leaves us with the case when $\kappa \equiv 2$ mod 4 and $n/\kappa$ is odd.  But in this case, no signed sum of the $n$ elements results in $(0,0)$, since the second component of the signed sum will differ from the second component of the sum (which is odd) in an even number (if we subtract an element instead of adding it, the sum gets reduced by an even value).  Therefore, we proved:

\begin{prop} \label{tau hat pm h=n}
Suppose that $G$ has order $n$ and exponent $\kappa$.  Then  
$$\tau \hat{_{\pm}}(G,n)=\left\{
\begin{array}{cl}
n & \mbox{if}  \; \kappa \equiv 2 \; \mbox{mod 4 and} \; n/\kappa \; \mbox{is odd}; \\ \\
n-1 & \mbox{otherwise}. \\ 
\end{array}\right.$$
\end{prop}

For $3 \leq h \leq n-1$, the value of $\tau \hat{_{\pm}} (G,h)$ is not known in general.  

The values of $\tau \hat{_{\pm}} (G,h)$ behave quite differently for even and odd values of $h$.  For example, let us consider $h=4$ first.  Note that, if for a set $A \subseteq G$, we have $0 \not \in 4\hat{_{\pm}}A$, then $A$ must be a weak Sidon set in $G$: Indeed, if we had $$a_1+a_2=a_3+a_4$$ for some $a_1,a_2,a_3,a_4 \in A$ with $a_1 \neq a_2$ and $a_3 \neq a_4$, then the four elements cannot be all distinct as then we would have $0 \in 4\hat{_{\pm}}A$.  Therefore, we must have, for example, $a_1 \in \{a_3,a_4\}$, which proves that $A$ is a weak Sidon set in $G$.  Therefore, by Corollary \ref{corzetaupperRlimited}, we get the following upper bound:

\begin{prop}
For every $G$ of order $n$, we have $$\tau \hat{_{\pm}}(G,4) \leq \lfloor \sqrt{2n} \rfloor -1.$$ 

\end{prop}

For odd values of $h$, however, we may have lower bounds for $\tau \hat{_{\pm}} (G,h)$ that are linear in $n$.  Consider first the cyclic group $\mathbb{Z}_n$ with $n$ even.  Note that, when $h$ is odd, for the set
$$A=\{1,3,5, \dots, n-1\}$$ we have $h \hat{_{\pm}}A \subseteq A$; in particular, $0 \not \in h \hat{_{\pm}}A$.  This shows that $$\tau \hat{_{\pm}}(\mathbb{Z}_n,h) \geq n/2.$$  On the other hand, 
$$\tau \hat{_{\pm}} (G,h) \leq \tau \hat{\;} (G,h),$$ so by Theorem \ref{Zfor n even h odd}, when $n \geq 12$ and $3 \leq h \leq n/2-2$, then 
$$\tau \hat{_{\pm}}(\mathbb{Z}_n,h) \leq n/2.$$
Therefore:

\begin{thm}
When $n$ is even, $h$ is odd, and $3 \leq h \leq n/2-2$, then $$\tau \hat{_{\pm}}(\mathbb{Z}_n,h) = n/2.$$

\end{thm}

For the case when $n$ and $h$ are both odd, Collins has the following lower bound:

\begin{prop} [Collins; cf.~\cite{Col:2014a}] \label{collins}\index{Collins, K.} 
When $n$ and $h \leq n$ are both odd, we have $$\tau \hat{_{\pm}}(\mathbb{Z}_n,h) \geq \left \lfloor \frac{n}{h}+ \frac{h^2-3}{2h} \right \rfloor.$$

\end{prop}

We present a proof for Proposition \ref{collins} on page \pageref{proof of collins}.  

We believe that Proposition \ref{collins} is quite accurate.  For example, when $p$ is prime, we know from Theorem \ref{Zforp} that
$$\tau \hat{_{\pm}}(\mathbb{Z}_p,h) \leq \tau \hat{\;}(\mathbb{Z}_p,h) = \left \lfloor \frac{p-2}{h} \right \rfloor + h.$$  In particular, for $h=3$ we get the following:

\begin{cor}
For every prime $p$, we have $$(p+1)/3 \leq \tau \hat{_{\pm}}(\mathbb{Z}_p,3) \leq (p+7)/3.$$

\end{cor} 

We offer the following intriguing problems:

\begin{prob}
Find $\tau \hat{_{\pm}}(\mathbb{Z}_p,3)$ for each prime $p$.

\end{prob}
 
\begin{prob}
Find $\tau \hat{_{\pm}}(\mathbb{Z}_p,h)$ for each prime $p$ and odd $h \leq p$.

\end{prob}

\begin{prob}
For each odd $n$ and odd $h$ with $3 \leq h \leq n-1$, evaluate $\tau \hat{_{\pm}}(\mathbb{Z}_n,h)$.

\end{prob}

\begin{prob}
For each even $n$ and odd $h$ with $n/2 - 1 \leq h \leq n-1$, evaluate $\tau \hat{_{\pm}}(\mathbb{Z}_n,h)$.

\end{prob}

\begin{prob}
For given $n$ and even $h$, evaluate, or at least find a good lower bound for $\tau \hat{_{\pm}}(\mathbb{Z}_n,h)$, in particular for $\tau \hat{_{\pm}}(\mathbb{Z}_n,4)$.

\end{prob}

Let us now turn to noncyclic groups.  The general question, of course, is as follows:

\begin{prob}
For each $G$ and $h$, evaluate, or at least find a good lower bound for $\tau \hat{_{\pm}}(G,h)$, in particular for $\tau \hat{_{\pm}}(G,3)$.

\end{prob}

The case of $h=\kappa$ attracts special interest.  Recall that the value $\tau\hat{\;}(G, \kappa)+1$ (the smallest integer such that each subset of $G$ of that cardinality has a subset size $\kappa$ whose elements sum to 0) is called the Harborth constant of $G$; analogously, here we define the {\em signed Harborth constant of $G$} as the value of $\tau \hat{_{\pm}}(G, \kappa)+1$.  This value is determined in Proposition \ref{tau hat pm h=n} above for cyclic groups; we have the following additional result: 
  
\begin{thm} [Marchan, Ordaz, Ramos, and Schmid; cf.~\cite{MarOrdRamSch:2013a}] \label{tau hat pm Z2 by Z4}\index{Schmid, W.}\index{Marchan, L. E.}\index{Ordaz, O.}\index{Ramos, D.} 
For a positive integer $k$, we have
$$\tau \hat{_{\pm}}(\mathbb{Z}_2 \times \mathbb{Z}_{2k}, 2k) = \left\{
\begin{array}{ccl}
4 & \mbox{if} & k \in \{1,2\}; \\ \\
2k+1 & \mbox{if} & k \geq 3.
\end{array}
\right.$$
\end{thm}
Furthermore, in \cite{MarOrdRamSch:2015a}  the same authors determined all subsets $A$ of size $2k+1$ in $\mathbb{Z}_2 \times \mathbb{Z}_{2k}$ for which $0 \not \in (2k) \hat{_{\pm}} A$.  

We pose the following problem:

\begin{prob}
Find $\tau \hat{_{\pm}}(G, \kappa)$ for each noncyclic group $G$ of exponent $\kappa \geq 3$.

\end{prob}
If $\kappa=2$, then the answer is given by Proposition \ref{tau hat pm h=1,2}: $$\tau \hat{_{\pm}}(\mathbb{Z}_2^r,2)=2^r$$ for every $r \in \mathbb{N}$.

\subsection{Limited number of terms} \label{5maxRSlimited}

The analogue of a $t$-independent set for restricted addition is called a {\em weak $t$-independent set}.  In particular, in this section we investigate, for a given group $G$ and positive integer $t$, the quantity $$\tau\hat{_{\pm}} (G,[1,t]) = \mathrm{max} \{ |A|  \mid A \subseteq G, 0 \not \in [1,t] \hat{_{\pm}} A\},$$ that is, the maximum size of a weak $t$-independent subset of $G$.

For $t=1$, we see that $A$ is weakly 1-independent if, and only if, $0 \not \in A$.  We can easily determine the value of $\tau\hat{_{\pm}}(G,[1,2])$ as well.  Obviously, a weak 2-independent set $A$ cannot contain 0; furthermore, for each element $g$ of $G$, $A$ cannot contain both $g$ and $-g$, unless $g=-g$.  So, to get a maximum weak 2-independent set in $G$, take exactly one of each element or its negative in $G \setminus \mathrm{Ord}(G,2)\setminus \{0\}$, and take all elements of $\mathrm{Ord}(G,2)$.  We can summarize these finding as follows:

\begin{prop}
For all finite abelian groups $G$ of order $n$ we have
$$\tau\hat{_{\pm}}(G,[1,1])=n-1$$ and
$$\tau\hat{_{\pm}}(G,[1,2])=\frac{n+|\mathrm{Ord}(G,2)|-1}{2}.$$  
As a special case, for the cyclic group of order $n$ we have 
 $$\tau\hat{_{\pm}}(\mathbb{Z}_n,[1,2])=\lfloor n/2 \rfloor.$$
\end{prop}

For $t \geq 3$, the value of $\tau\hat{_{\pm}}(G,[1,t])$ is not known.   In \cite{BajRuz:2003a}, Bajnok and Ruzsa\index{Ruzsa, I.}\index{Bajnok, B.} proved the following general bounds:

\begin{thm} [Bajnok and Ruzsa; cf.~\cite{BajRuz:2003a}] \label{Baj Ruz tau hat pm}\index{Ruzsa, I.}\index{Bajnok, B.}
For every $G$ and $t \geq 2$ we have
$$\left( t! / 2^t \cdot n  \right) ^{1/t} -t/2 < \tau\hat{_{\pm}}(G,[1,t]) < \left( \lfloor t/2 \rfloor ! \cdot n  \right) ^{1/\lfloor t/2 \rfloor} +t/2.$$

\end{thm}
In particular, for $t \geq 4$, $\tau\hat{_{\pm}}(G,[1,t])$ is not a linear function of $n$.

For $t=3$, the upper bound in Theorem \ref{Baj Ruz tau hat pm} is trivial; we can do better, as follows.  Observe that if $A$ is a weak $3$-independent set in $G$, then the sets $\{0\}$, $A$, and $2 \hat{\;} A$ must be pairwise disjoint: indeed, if we were to have, say, $a_1+a_2=a_3$ for some $a_1,a_2,a_3\in A$ with $a_1 \neq a_2$, then we would have to have $a_1=a_3$ or $a_2=a_3$, but that would imply that one of the elements is 0, a contradiction.  Note also that for each $a_0 \in A$, $$\{a_0 +a \mid a \in A, a \neq a_0\} \subseteq 2 \hat{\;} A,$$ so $2 \hat{\;} A$ has size at least $|A|-1$.  Therefore, $$n \geq |\{0\} \cup A \cup 2 \hat{\;} A| \geq 1+|A|+|A|-1=2|A|,$$ which yields the following upper bound:

\begin{prop} \label{tau hat pm  <= n/2}
For all $G$ of order $n$ we have $$\tau\hat{_{\pm}}(G,[1,3]) \leq n/2.$$

\end{prop} 
 A better upper bound for $\tau\hat{_{\pm}}(G,[1,3])$ can be derived in the case when $n$ is odd. Note that, in this case, if $A$ is weakly 3-independent, then $A$ and $-A$ must be disjoint.  (This may be false if $n$ is even: the elements of order 2 may belong to both $A$ and $-A$.)  Furthermore, the set $A \cup (-A)$ is weakly zero-3-sum-free: the sum of three distinct elements of $A  \cup (-A)$ cannot be zero as it is either an element of $A \cup (-A)$ or the signed sum of three distinct elements of $A$.  Therefore, $A \cup (-A)$ has size at most equal to $\tau \hat{\;} (G,3)$, which yields:

\begin{prop} \label{tau hat pm with tau hat}
For all $G$ of odd order $n$ we have $$\tau\hat{_{\pm}}(G,[1,3]) \leq \tau \hat{\;} (G,3)/2.$$ 

\end{prop} 

For $t=3$ we can also derive a much better lower bound for $\tau\hat{_{\pm}}(G,[1,3])$ than what Theorem \ref{Baj Ruz tau hat pm} gives based on the idea that if $A$ is a (regular) 3-independent set in a group $G_2$, then $G_1 \times A$ is a weak 3-independent set in $G_1 \times G_2$ for any group $G_1$:

\begin{prop} \label{weak 3-indep direct prod}
For all groups $G_1$ and $G_2$ we have $$\tau\hat{_{\pm}}(G_1 \times G_2,[1,3]) \geq  |G_1| \cdot \tau_{\pm} (G_2,[1,3]).$$

\end{prop}    
Since any group $G$ of order $n$ and exponent $\kappa$ is isomorphic to $G_1 \times \mathbb{Z}_{\kappa}$ for some group $G_1$ of order $n/\kappa$, Proposition \ref{weak 3-indep direct prod} implies:
\begin{cor} \label{tau hat pm 3 with tau pm 3}
For every group $G$ of order $n$ and exponent $\kappa$ we have $$\tau\hat{_{\pm}}(G,[1,3]) \geq  n/\kappa \cdot \tau_{\pm} (\mathbb{Z}_{\kappa},[1,3]).$$

\end{cor}
Recall that the value of $\tau_{\pm} (\mathbb{Z}_{\kappa},[1,3])$ was explicitly evaluated in Theorem \ref{3free}.

While we have the various bounds just discussed, exact values for $\tau\hat{_{\pm}}(G,[1,3])$ are not known even for cyclic groups:

\begin{prob}

Find the exact values of $\tau\hat{_{\pm}}(\mathbb{Z}_n,[1,3])$.

\end{prob}

Let us now consider noncyclic groups, in particular, groups of the form $\mathbb{Z}^r_k$.  

First, we consider $\mathbb{Z}_2^r$.  It is not hard to see that the set $\mathbb{Z}_2^{r-1} \times \{1\}$ is weakly 3-independent in $\mathbb{Z}^2_r$ (note that no two distinct elements add to zero), and therefore $$\tau\hat{_{\pm}}(\mathbb{Z}_2^r,[1,3]) \geq 2^{r-1}.$$ Together with Proposition \ref{tau hat pm  <= n/2}, this implies:

\begin{prop}
For every $r \geq 1$ we have $$\tau\hat{_{\pm}}(\mathbb{Z}_2^r,[1,3])=2^{r-1}.$$

\end{prop}

We can obtain a large weakly 3-independent set in $\mathbb{Z}_3^r$ by using Bui's idea\index{Bui, C.} described on page \pageref{bui'sidea}: (using notations described there) the set $\cup_{i \in I} A_i$ is weakly 3-independent in $\mathbb{Z}^r_3$, yielding the following lower bound.

\begin{thm} \label{tau hat pm Biu}
Let $r$  be any positive integer, and write $r=3q+s$ with $s=0,1$, or 2.  Then we have
$$\tau\hat{_{\pm}}(\mathbb{Z}_3^r,[1,3]) \geq \sum_{i=q+s-1}^{2q+s-1} {r \choose i}.$$
\end{thm}

The lower bounds provided by Theorem \ref{tau hat pm Biu} are quite good.  Indeed, by Proposition \ref{tau hat pm with tau hat} above and considering the values for $\tau\hat{\;}(\mathbb{Z}_3^r,3)$ on page \pageref{tau hat 3 table}, we see that they give the exact values for $r \leq 4$:
$$\tau\hat{_{\pm}}(\mathbb{Z}_3,[1,3])=1, \; \tau\hat{_{\pm}}(\mathbb{Z}_3^2,[1,3])=2, \; \tau\hat{_{\pm}}(\mathbb{Z}_3^3,[1,3])=4, \; \tau\hat{_{\pm}}(\mathbb{Z}_3^4,[1,3])=10.$$

For $k \geq 4$, we know no better lower bounds for $\tau\hat{_{\pm}}(\mathbb{Z}^r_k,[1,3])$  than what Corollary \ref{tau hat pm 3 with tau pm 3} above provides.

We pose the following problems:

\begin{prob}

Find the exact values of $\tau\hat{_{\pm}}(\mathbb{Z}^r_3,[1,3])$ for $r \geq 5$.

\end{prob}

\begin{prob}

Find the exact values of $\tau\hat{_{\pm}}(\mathbb{Z}^r_4,[1,3])$, $\tau\hat{_{\pm}}(\mathbb{Z}^r_5,[1,3])$, etc.

\end{prob}

\begin{prob}

Find $\tau\hat{_{\pm}}(G,[1,3])$ for other noncyclic groups $G$.

\end{prob}

Of course, the ultimate goal is to answer this question:

\begin{prob} \label{prob tau hat pm G, 1,t}
Find $\tau\hat{_{\pm}}(G,[1,t])$ for all groups $G$ and integers $t \geq 3$.
\end{prob}

As a modest step toward Problem \ref{prob tau hat pm G, 1,t}, we establish:

\begin{prop} \label{prob tau hat pm G, 1,t n small}
If $n$ and $t$ are positive integers so that $2^{t-1} \leq n < 2^t$, then $$\tau\hat{_{\pm}}(\mathbb{Z}_n,[1,t])=t-1.$$  

\end{prop}
The short proof of Proposition \ref{prob tau hat pm G, 1,t n small} is on page \pageref{proof of prob tau hat pm G, 1,t n small}.

\subsection{Arbitrary number of terms} \label{5maxRSarbitrary}

In this subsection we are trying to evaluate
$$\tau\hat{_{\pm}}(G,\mathbb{N})=\max \{ |A|  \mid A \subseteq G, 0 \not \in \cup_{h=1}^{\infty} h \hat{_{\pm}} A\}.$$  A subset $A$ of $G$ for which $$0 \not \in \cup_{h=1}^{\infty} h \hat{_{\pm}} A$$ holds is called a {\em dissociated} subset of $G$.  Recall that, by Proposition \ref{C.3.3=F.4.3},  our condition is equivalent to $$|\Sigma A|=2^{|A|}.$$  Indeed, we have
$$\tau\hat{_{\pm}}(G,\mathbb{N})=\sigma\hat{_{\pm}}(G,\mathbb{N}_0).$$

The case when $G$ is cyclic is easy.  Suppose that $G=\mathbb{Z}_n$, and let $$m=\lfloor \log_2 n \rfloor;$$ note that we then have
$$2^m \leq n < 2^{m+1}.$$
According to Proposition \ref{prob tau hat pm G, 1,t n small}, we have $$\tau\hat{_{\pm}}(\mathbb{Z}_n,[1,m+1])=m.$$  In particular, there is an $m$-subset $A$ of $\mathbb{Z}_n$ (namely, as the proof of Proposition \ref{prob tau hat pm G, 1,t n small} shows, $$A=\{1,2,\dots,2^{m-1}\}$$ works) for which $0 \not \in [1,m+1] \hat{_{\pm}} A$.  But, since $|A|=m$, $$[1,m+1] \hat{_{\pm}} A = \cup_{h=1}^{\infty} h \hat{_{\pm}} A,$$ which proves that 
$\tau\hat{_{\pm}}(\mathbb{Z}_n,\mathbb{N}) \geq m.$

On the other hand, clearly,
$$ \tau\hat{_{\pm}}(\mathbb{Z}_n,\mathbb{N}) \leq \tau\hat{_{\pm}}(\mathbb{Z}_n,[1,m+1])=m.$$ Thus we have:

\begin{prop}  \label{tau hat pm for cyclic groups}
For every positive integer $n$, $$ \tau\hat{_{\pm}}(\mathbb{Z}_n,\mathbb{N}) =\lfloor \log_2 n \rfloor.$$

\end{prop} 

We do not quite have the value of $\tau \hat{_{\pm}}(G,\mathbb{N})$ for noncyclic groups $G$, but we can find lower and upper bounds, as follows.  For an upper bound, note that if $A$ is dissociated in $G$, then $|\Sigma A| = 2^m$, so $\tau \hat{_{\pm}}(G,\mathbb{N}) \leq \lfloor \log_2 n \rfloor.$  

For a lower bound, let us suppose that $G$ is of type $(n_1,\dots,n_r);$ that is, 
we have integers $n_1,\dots,n_r$ so that $2 \leq n_1$ and $n_i$ is a divisor of $n_{i+1}$ for each $i=1,2,\dots,r-1$, and for which
$$G \cong \mathbb{Z}_{n_1} \times \cdots \mathbb{Z}_{n_r}.$$  Let $m_i = \lfloor \log_2 n_i \rfloor$.  It is easy to see that the set 
$$\left( \{1,2,\dots,2^{m_1-1}\} \times \{0\}^{r-1} \right) \cup \cdots \cup \left( \{0\}^{r-1}  \times \{1,2,\dots,2^{m_r-1}\}\right)$$ is dissociated in $G$, and thus 
$$\tau \hat{_{\pm}}(G,\mathbb{N}) \geq m_1+\cdots+m_r.$$  

Therefore, we get:

\begin{prop} \label{bounds for dissoc in G}
Suppose that $G$ is an abelian group of type $(n_1,\dots,n_r)$ and order $n$.  We then have
$$\lfloor \log_2 n_1 \rfloor + \cdots + \lfloor \log_2 n_r \rfloor \leq \tau \hat{_{\pm}}(G,\mathbb{N}) \leq \lfloor \log_2 n \rfloor.$$

\end{prop}

Considering that $$ \log_2 n_1  + \cdots + \log_2 n_r = \log_2 (n_1 \cdots n_r) = \log_2 n ,$$ we can say that the lower and upper bounds in Proposition \ref{bounds for dissoc in G} are close.  However, they are not equal: For example, for the groups $\mathbb{Z}_6^2$ or $\mathbb{Z}_7^2$, the lower bound equals 4 and the upper bound equals 5, and, indeed, the set
$$\{(0,1),(1,0),(1,2),(1,4),(3,2)\}$$ (for example) is dissociated in both groups.

Hence the general question remains open:

\begin{prob}
Find the value of $\tau \hat{_{\pm}}(G,\mathbb{N})$ for noncyclic groups $G$.
\end{prob}

Let us now turn to the inverse problem of classifying all dissociated subsets of $G$ of maximum size.  We consider cyclic groups first; in particular, the cyclic group of order $n=2^k$ where $k \in \mathbb{N}$---note that, by Proposition \ref{tau hat pm for cyclic groups}, we have $ \tau\hat{_{\pm}}(\mathbb{Z}_{2^k},\mathbb{N}) =k$.  

We provide the following recursive construction for a collection ${\cal A}_k$ of $k$-subsets of $\mathbb{Z}_{2^k}$: 
\begin{itemize}
  \item We let ${\cal A}_1$ consist of the single subset $\{1\}$ of $\mathbb{Z}_2$.
  \item Suppose that  ${\cal A}_j$ is already constructed for some positive integer $j$.  For a given member $A_j=\{a_1,\dots,a_j\}$ of ${\cal A}_j$ and for a given element $\epsilon = (\epsilon_1, \dots, \epsilon_j)$ of $\mathbb{Z}_2^j$,  we define the set $$A_{j+1}(A_j,\epsilon) = \{2^j\} \cup \{a_1+\epsilon_1 \cdot 2^j , \dots, a_j+\epsilon_j \cdot 2^j \};$$ and then set $${\cal A}_{j+1} = \{A_{j+1}(A_j,\epsilon) \mid A_j \in {\cal A}_j \; \mbox{and} \; \epsilon \in \mathbb{Z}_2^j\}.$$
  \end{itemize}
So, for $k=2$, we get
$${\cal A}_{2} = \{ \{2\} \cup \{1+ 0 \cdot 2\}, \{2\} \cup \{1+ 1 \cdot 2\}=\{\{1,2\},\{2,3\}\};$$
and for $k=3$, we have
\begin{eqnarray*}
{\cal A}_{3} & = & \{ \{4 \} \cup \{1+ 0 \cdot 4, 2+ 0 \cdot 4\}, \{4 \} \cup \{1+ 0 \cdot 4, 2+ 1 \cdot 4\}, \{4 \} \cup \{1+ 1 \cdot 4, 2+ 0 \cdot 4\}, \\
& & \{4 \} \cup \{1+ 1 \cdot 4, 2+ 1 \cdot 4\}, \{4 \} \cup \{2+ 0 \cdot 4, 3+ 0 \cdot 4\}, \{4 \} \cup \{2+ 0 \cdot 4, 3+ 1 \cdot 4\}, \\
& & \{4 \} \cup \{2+ 1 \cdot 4, 3+ 0 \cdot 4\}, \{4 \} \cup \{2+ 1 \cdot 4, 3+ 1 \cdot 4\} \} \\
& = & \{ \{1,2,4\}, \{1,4,6\}, \{2,4,5\}, \{4,5,6\} , \{2,3,4\}, \{2,4,7\}, \{3,4,6\}, \{4,6,7\}\}.
\end{eqnarray*}

We have the following conjecture for dissociated subsets of $\mathbb{Z}_{2^k}$ of maximum size $$\tau \hat{_{\pm}}(\mathbb{Z}_{2^k},\mathbb{N})=k.$$

\begin{conj} \label{conj on inverse tau hat pm 2 power}
A $k$-subset $A$ of $\mathbb{Z}_{2^k}$ is dissociated if, and only if $A \in {\cal A}_k$. 

\end{conj}

Using the computer program \cite{Ili:2017a}, we have verified Conjecture \ref{conj on inverse tau hat pm 2 power} for $k \in \{1,2,3,4\}$.   Note that Conjecture \ref{conj on inverse tau hat pm 2 power} implies that there are exactly $2^{(k^2-k)/2}$ dissociated subsets of $\mathbb{Z}_{2^k}$ of maximum size.

\begin{prob}

Prove (or disprove) Conjecture \ref{conj on inverse tau hat pm 2 power}.

\end{prob}

Observe that, if Conjecture \ref{conj on inverse tau hat pm 2 power} is true, then every dissociated subset of $\mathbb{Z}_{2^k}$ of maximum size contains the element $2^{k-1}$.  As a modest step toward Conjecture \ref{conj on inverse tau hat pm 2 power}, we state:

\begin{conj} \label{conj on inverse tau hat pm 2 power middle}
Every dissociated subset of $\mathbb{Z}_{2^k}$ of size $k$ contains the element $2^{k-1}$.

\end{conj}

\begin{prob}

Prove (or disprove) Conjecture \ref{conj on inverse tau hat pm 2 power middle}.

\end{prob}

More generally:

\begin{prob}

For each positive integer $n$, classify all dissociated subsets of $\mathbb{Z}_n$ of maximum size $\lfloor \log_2 n \rfloor$. 
\end{prob}

There is an interest in finding the largest dissociated subsets in any subset of $G$, not just $G$ itself.  For a subset $A$ of $G$, we let
$$\dim A = \max \{ |B| \mid B \subseteq A, B \; \mbox{is dissociated} \}$$ denote the {\em dissociativity dimension} of $A$ in $G$.   (We mention in passing that there are other notions for the dimension of a set: cf.~\cite{CanHel:2015a}\index{Candela, P.}\index{Helfgott, H. A.} and \cite{SchShk:2016a}.)\index{Schoen, T.}\index{Shkredov, I. D.}  Of course,
$$\dim G = \tau \hat{_{\pm}}(G,\mathbb{N}).$$

As a (less tight) generalization of Proposition \ref{bounds for dissoc in G}, Lev and Yuster proved the following:

\begin{prop} [Lev and Yuster; cf.~\cite{LevYus:2011a}] \label{bounds for dissoc in A}\index{Lev, V. F.}\index{Yuster, R.} 
For any subset $A$ of $G$, we have
$$r_A \leq \dim A \leq \lfloor r_A \cdot \log_2 \kappa \rfloor ,$$ where $\kappa$ is the exponent of $G$ and $r_A$ is the rank of the subgroup $\langle A \rangle$ generated by $A$.  

\end{prop}

For any group $G$ of order $n$ and for any positive integer $m \leq n$, we introduce the following quantity:
$$\dim (G, m) = \min \{ \dim A \mid A \subseteq G, |A|=m \}.$$
Our question is as follows:

\begin{prob}
For every group $G$ of order $n$ and for any positive integer $m \leq n$, find $\dim (G, m)$.  

\end{prob}

Clearly, $\dim (G, 1)=0$ in every group $G$, since the only dissociated subset of $\{0\}$ is the empty-set. It is also easy to see that  $\dim (G, 2)=1$ (consider a set $\{0,g\}$ for some non-zero element $g$).  For $\dim (G, 3)$, if $G$ has exponent 3 or more, we can take the set $\{0,g,-g\}$ for any $g \in G$ of order 3 or more, so $\dim (G, 3)=1$; if $G$ is the elementary abelian 2-group, then $\dim (G, 3)=2$.  We have the following conjecture:

\begin{conj} \label{conj on dim G m}  
For all positive integers $n$ and $m \leq n$, we have $$\dim (\mathbb{Z}_n,m)=\lfloor \log_2 m \rfloor.$$

\end{conj}

For example, we have $\dim (\mathbb{Z}_{10},7)=2$, since each $4$-subset (and thus each 7-subset) of $\mathbb{Z}_n$ (when $n \geq 4$) has dimension at least 2 (only subsets of a set of the form $\{0,g,-g\}$ have dimension less than 2), and the set $$\{0,1,2,3,7,8,9\}$$ has dimension 2.

For $m=n$, Conjecture \ref{conj on dim G m} becomes Proposition \ref{tau hat pm for cyclic groups}.  Furthermore, when $n=2^k$ for some $k \in \mathbb{N}$, then for $m=n-1$, Conjecture \ref{conj on dim G m} says that $$\dim (\mathbb{Z}_{2^k},2^k-1)=k-1;$$ in other words, one can find a subset of size $2^k-1$ in $\mathbb{Z}_{2^k}$, which does not have a dissociated subset of size $k$.  Observe that, if Conjecture \ref{conj on inverse tau hat pm 2 power middle} is true, then the set $\mathbb{Z}_{2^k} \setminus \{2^{k-1}\}$ is one such set.

\begin{prob}
Prove (or disprove) Conjecture \ref{conj on dim G m}.

\end{prob}

\chapter{Sum-free sets} \label{ChapterSumfree}

Recall that for a given finite abelian group $G$, $m$-subset $A=\{a_1,\dots, a_m\}$ of $G$, $\Lambda \subseteq \mathbb{Z}$, and $H \subseteq \mathbb{N}_0$, we defined the sumset of $A$ corresponding to $\Lambda$ and $H$ as
$$H_{\Lambda}A = \{\lambda_1a_1+\cdots +\lambda_m a_m \mbox{    } |  \mbox{    }  (\lambda_1,\dots ,\lambda_m) \in \Lambda^m(H)  \}$$
where the index set $\Lambda^m(H)$ is defined as
$$\Lambda^m(H)=\{(\lambda_1,\dots ,\lambda_m) \in \Lambda^m \; |  \; |\lambda_1|+\cdots +|\lambda_m| \in H \}.$$ 

In this chapter we consider, for $G$, $H$, and $\Lambda$, {\em $H$-sum-free subsets} of $G$ over $\Lambda$; that is, subsets $A$ of $G$ for which 
$$(h_1)_{\Lambda}A \cap (h_2)_{\Lambda}A=\emptyset$$ for any two distinct elements  $h_1$ and $h_2$ of $H$.  

Let us make some preliminary observations.  First of all, if $H$ contains fewer than two elements, then every subset $A$ of $G$ is $H$-sum-free over $\Lambda$.  Second, if $H=\{0,h\}$ for some positive integer $h$, then for a subset to be $H$-sum-free is the same as it being zero-$h$-sum-free, a property we studied in Chapter \ref{ChapterZerosumfree}.  More generally, if $0 \in H$, then a set $A$ being $H$-sum-free implies that $A$ is zero-$H'$-sum-free for $H'=H \setminus \{0\}$.

Furthermore, the sum-free property is clearly a weakening of the Sidon property studied in Chapter \ref{ChapterSidon}: if linear combinations corresponding to different elements of the entire index set $\Lambda^m(H)$ are distinct, then they are certainly distinct when corresponding to distinct $h_1$ and $h_2$ of $H$.

While the sum-free property is thus closely related to properties discussed elsewhere in the book, it offers unique opportunities for the study of interesting and well known questions in additive combinatorics.

In this chapter we attempt to find $\mu_{\Lambda}(G,H)$, the maximum possible size of an $H$-sum-free set over $\Lambda$ in a given finite abelian group $G$.   If no $H$-sum-free set exists, we put $\mu_{\Lambda}(G,H) = 0.$   With this notation, our observations above can be stated as follows:
\begin{prop} \label{sumsets trivial prop}

Let $G$ and $\Lambda \subseteq \mathbb{Z}$ be arbitrary.

If $|H|\leq 1$, then $\mu_{\Lambda}(G,H) = n.$

If $0 \in H$, then $$\mu_{\Lambda}(G,H) \leq \tau_{\Lambda}(G,H \setminus \{0\});$$ in particular, if $H=\{0,h\}$, then 
$$\mu_{\Lambda}(G,H) = \tau_{\Lambda}(G,h).$$

\end{prop}

\begin{prop}
For all  $G$, $\Lambda \subseteq \mathbb{Z}$, and $H \subseteq \mathbb{N}_0$ we have
$$\mu_{\Lambda}(G,H) \geq \sigma_{\Lambda}(G,H).$$ 
\end{prop}

In the following sections we attempt to find $\mu_{\Lambda}(G,H)$ for special coefficient sets $\Lambda$.

\section{Unrestricted sumsets} \label{6maxU}

Our goal in this section is to investigate the maximum possible size of an $H$-sum-free set over the set of nonnegative integers, that is, the quantity $$\mu (G,H) =\mathrm{max} \{ |A|  \mid  A \subseteq G; h_1, h_2 \in H; h_1 \neq h_2 \Rightarrow h_1A \cap h_2A = \emptyset \}.$$  

Clearly, we always have $\mu (G,H)=0$ whenever $H$ contains two elements $h_1$ and $h_2$ whose difference is a multiple of the exponent $\kappa$ of $G$, since for any element $a \in G$, we then have $h_1a=h_2a$.  However, when distinct elements of $H$ leave different remainders when divided by $\kappa$, then $\mu  (G,H) \geq 1$: for any $a \in G$ with order $\kappa$, at least the one-element set $A=\{a\}$ will be $H$-sum-free.

It is often useful to consider $G$ of the form $G_1 \times G_2$.  (We may do so even when $G$ is cyclic if its order has at least two different prime divisors.)  It is not hard to see that, if $A_1 \subseteq G_1$ is $H$-sum-free in $G_1$, then $$A=\{(a,g)  \mid a \in A_1, g \in G_2\}$$ is $H$-sum-free  in $G$.  Indeed, if we were to have
$$(a_1,g_1)+ \cdots + (a_{h_1},g_{h_1}) = (a_1',g_1')+ \cdots + (a_{h_2}',g_{h_2}')$$
for some $h_1, h_2 \in H$, $h_1 \neq h_2$, and $(a_i,g_i),(a_i',g_i') \in A$, then the equation for the sum of the first coordinates contradicts the fact that $A_1 $ is $H$-sum-free in $G_1$.  Thus, we have the following.    

\begin{prop} \label{mufordirectsum}
For all finite abelian groups $G_1$ and $G_2$ and for all $H \subseteq \mathbb{N}_0$ we have
$$\mu(G_1 \times G_2,H) \geq \mu(G_1,H) \cdot |G_2|.$$

\end{prop}

Below we consider two special cases: when $H$ consists of two distinct positive integers (by Proposition \ref{sumsets trivial prop}, the case when one of the integers equals 0 is identical to Section \ref{5maxUfixed}), and when $H$ consists of all positive integers up to some value $s$.  The cases when $H = \mathbb{N}_0$ or $H = \mathbb{N}$, as we just mentioned, yield no $H$-sum-free sets.

\subsection{Fixed number of terms} \label{6maxUfixed}

Suppose that $k$ and $l$ are distinct positive integers; without loss of generality, we assume that $k >l$.  Sets satisfying the condition $(k A) \cap (l A) = \emptyset $ are called {\it $(k,l)$-sum-free sets}.  Here we intend  to determine the maximum value of $m$ for which $G$ contains a $(k,l)$-sum-free subset of size $m$---this value is denoted here by ${\mu}(G,\{k,l\})$.  

There are several similarities between $(k,l)$-sum-free sets and zero-$h$-sum-free sets; in fact, in some respects, one might consider a zero-$h$-sum-free set to be $(h,0)$-sum-free.  Therefore, our discussion in this section will be similar to that in Section \ref{5maxUfixed}.  This similarity has its limits, however---see, for example, our comments below before Problem \ref{prob find mu > v} and after Conjecture \ref{conj mu general G}.

Since $k$ and $l$ are positive integers with $k>l$, the first case to discuss is when $k=2$ and $l=1$; a $(2,1)$-sum-free set in $G$ is simply referred to as a \emph{sum-free set}.

For the cyclic group $\mathbb{Z}_n$, we can find explicit sum-free sets as follows.  For every $n$, the integers that are between $n/3$ and $2n/3$ form a sum-free set; more precisely, the set 
$$\left\{\left\lfloor \tfrac{n}{3} \right\rfloor, \left\lfloor \tfrac{n}{3} \right\rfloor +1, \left\lfloor \tfrac{n}{3} \right\rfloor+2, \dots, 2\left\lfloor \tfrac{n}{3} \right\rfloor-1  \right\}$$ is sum-free in $\mathbb{Z}_n$.   (The integers between $n/6$ and $n/3$, together with those between $2n/3$ and $5n/6$, with the endpoints carefully chosen, provide another sum-free set in $\mathbb{Z}_n$.)  

Like in the case of zero-3-sum-free sets, we can do better when $n$ has a prime divisor $p$ which is congruent to 2 mod 3.  We see that the set  
$$\left\{(p+1)/3+ pi_1+i_2 \mbox{    } | \mbox{    } i_1=0,1,\dots,n/p-1, \mbox{  } i_2=0,1,\dots,(p-2)/3 \right\}$$ is sum-free.

These examples show that we have 
$${\mu}(\mathbb{Z}_n, \{2,1\}) \geq \left\{
\begin{array}{ll}
\left(1+\frac{1}{p}\right) \frac{n}{3} & \mbox{if $n$ has prime divisors congruent to 2 mod 3,} \\ & \mbox{and $p$ is the smallest such divisor,}\\ \\
\left\lfloor \frac{n}{3} \right\rfloor & \mbox{otherwise;}\\
\end{array}\right.$$
and Diananda and Yap\index{Diananda, P. H.}\index{Yap, H. P.} proved in 1969 (see \cite{DiaYap:1969a}) that equality holds.  \label{mu21} Thus, using our function introduced on page \pageref{0.3.5.3 page}, we have:

\begin{thm} [Diananda and Yap; cf.~\cite{DiaYap:1969a}] \label{Diananda and Yap}\index{Diananda, P. H.}\index{Yap, H. P.} 

For all positive integers $n$, we have
$${\mu}(\mathbb{Z}_n, \{2,1\})=v_1(n,3).$$

\end{thm}  

More generally, we have the following result.

\begin{thm}[Bajnok; cf.~\cite{Baj:2009a}] \label{mu}\index{Bajnok, B.}
Suppose that $k$ and $l$ are positive integers and $k>l$.  Then we have
$$v_{k-l}(n,k+l) \leq {\mu}(\mathbb{Z}_n, \{k,l\}) \leq v_1(n,k+l).$$ 
\end{thm}  

Of course, when the lower and upper bounds coincide, we get equality:

\begin{cor} \label{k-l and n rel prime}
If $k-l$ and $n$ are relatively prime, then $${\mu}(\mathbb{Z}_n, \{k,l\}) = v_1(n,k+l).$$  In particular, for all positive integers $n$ and $l$, we have $${\mu}(\mathbb{Z}_n, \{l+1,l\}) = v_1(n,2l+1).$$ 

\end{cor} 
  We note that Corollary \ref{k-l and n rel prime} was established in \cite{HamPla:2003a} by Hamidoune and Plagne.\index{Hamidoune, Y. O.}\index{Plagne, A.}

As in Theorem \ref{zforp},  we can determine the size of the largest $(k,l)$-sum-free set in cyclic groups of prime order:

\begin{thm} \label{muforp}
The size of the largest $(k,l)$-sum-free set in  cyclic group of prime order $p$ is 
$${\mu}(\mathbb{Z}_p, \{k,l\}) =v_{k-l}(p,k+l)=\left\{
\begin{array}{ll}
0 & \mbox{if $p|(k-l)$,} \\ \\
\left \lfloor \frac{p-2}{k+l} \right \rfloor +1 & \mbox{otherwise.}\\
\end{array}\right.$$ 
\end{thm}

The value of ${\mu}(G, \{k,l\})$ for groups of composite order is not known in general.  By the following result, it suffices to consider arithmetic progressions in order to determine ${\mu}(\mathbb{Z}_n, \{k,l\})$:

\begin{thm}[Bajnok; cf.~\cite{Baj:2009a}] \label{enough for ap}\index{Bajnok, B.}
For a given divisor $d$ of $n$, let $\alpha (\mathbb{Z}_d,\{k,l\})$ be the maximum size of a $(k,l)$-sum-free arithmetic progression in $\mathbb{Z}_d$.  Then 
$${\mu}(\mathbb{Z}_n, \{k,l\}) = \max \left\{\alpha (\mathbb{Z}_d,\{k,l\}) \cdot \frac{n}{d} \mid d \in D(n) \right\}.$$

\end{thm} 

Using Theorem \ref{enough for ap}, we can compute $ {\mu}(\mathbb{Z}_n, \{3,1\})$ and $ {\mu}(\mathbb{Z}_n, \{4,1\})$:

\begin{thm} [Bajnok; cf.~\cite{Baj:2009a}] \label{Bajnok mu(3,1)}\index{Bajnok, B.}
For all positive integers $n$, we have $$ {\mu}(\mathbb{Z}_n, \{3,1\}) =v_{2}(n,4).$$

\end{thm}

\begin{thm} [Butterworth; cf.~\cite{But:2009a}]\index{Butterworth, J.}
For all positive integers $n$, we have $$ {\mu}(\mathbb{Z}_n, \{4,1\}) 
=v_{3}(n,5).$$

\end{thm}

Recall from page \pageref{v2n4} that $$v_{2}(n,4)=\left\{
\begin{array}{ll}
\left(1+\frac{1}{p}\right) \frac{n}{4} & \mbox{if $n$ has prime divisors congruent to 3 mod 4,} \\ & \mbox{and $p$ is the smallest such divisor,}\\ \\
\left\lfloor \frac{n}{4} \right\rfloor & \mbox{otherwise;}\\
\end{array}\right.$$
the value of $v_3(n,5)$ is more complicated, see page \pageref{v3n5}.

\begin{prob}
Use Theorem \ref{enough for ap} to compute $ {\mu}(\mathbb{Z}_n, \{k,l\})$ for other choices of $k$ and $l$.
\end{prob}

It may seem from our results thus far that $ {\mu}(\mathbb{Z}_n, \{k,l\})$ is given by $v_{k-l}(n,k+l).$  While this seems to be often the case, there are instances when $$ {\mu}(\mathbb{Z}_n, \{k,l\}) > v_{k-l}(n,k+l);$$ for example, $ {\mu}(\mathbb{Z}_9, \{5,2\})=2$ (e.g. the set $\{1,2\}$ is $(5,2)$-sum-free in $\mathbb{Z}_9$) while $v_3(9,7)=1$, and $ {\mu}(\mathbb{Z}_{16}, \{5,1\})=3$ (e.g. the set $\{1,2,3\}$ is $(5,1)$-sum-free in $\mathbb{Z}_{16}$) while $v_4(16,6)=2$.  (These examples point to a difference between sum-free sets and zero-sum-free sets; cf.~Conjecture \ref{zconj}.)  It is a very interesting question to find others:

\begin{prob}  \label{prob find mu > v}
Find other cases when $ {\mu}(\mathbb{Z}_n, \{k,l\}) > v_{k-l}(n,k+l).$
\end{prob}  

The general question regarding cyclic group is, of course:

\begin{prob}
Evaluate $ {\mu}(\mathbb{Z}_n, \{k,l\})$ for all positive integers $n$, $k$, and $l$.
\end{prob}

Turning to the case of noncyclic groups, we first state the following consequence of Proposition \ref{mufordirectsum}:

\begin{cor} \label{cor sumfree cyclic to non}
For every group $G$ of order $n$ and exponent $\kappa$ we have
$$ {\mu}(G, \{k,l\})  \geq {\mu}(\mathbb{Z}_{\kappa}, \{k,l\}) \cdot \frac{n}{\kappa}.$$
\end{cor}

We (somewhat hesitantly) believe that equality holds:

\begin{conj} \label{conj mu general G}
For every group $G$ of order $n$ and exponent $\kappa$ we have
$$ {\mu}(G, \{k,l\}) ={\mu}(\mathbb{Z}_{\kappa}, \{k,l\}) \cdot \frac{n}{\kappa}.$$
\end{conj}

We should point out that Conjecture \ref{conj mu general G}, if true, points to a difference in behavior between sum-free sets and zero-sum-free sets, since we definitely may have  
$$ \tau(G, h) \neq \tau(\mathbb{Z}_{\kappa}, h) \cdot \frac{n}{\kappa};$$ see, for example, Proposition \ref{Finn's constr} and Theorem \ref{tau elementary}.  If Conjecture \ref{conj mu general G} is true, it is likely to be extremely difficult to prove; we offer the challenging problem:

\begin{prob}
Prove (or disprove) Conjecture \ref{conj mu general G}.

\end{prob}

We have the following partial result:

\begin{thm} [Bajnok; cf.~\cite{Baj:2009a}]\index{Bajnok, B.}
Conjecture \ref{conj mu general G} holds whenever the exponent $\kappa$ of $G$ possesses at least one divisor $d$ that is not congruent to any integer between 1 and $\gcd (d, k-l)$ (inclusive) mod $k+l$.
\end{thm}

Thus, for example, if the exponent---or, equivalently, the order---of $G$ is divisible by 3 or has at least one (prime) divisor congruent to 2 mod 3, then 
$ {\mu}(G, \{2,1\})$ is determined.  The case when all (prime) divisors of $G$ are congruent to 1 mod 3 was unsolved for four decades, until  in 2005  in the breakthrough paper \cite{GreRuz:2005a} Green and Ruzsa proved that Conjecture \ref{conj mu general G} holds for sum-free sets:

\begin{thm} [Green and Ruzsa; cf.~\cite{GreRuz:2005a}] \label{(2,1)all}\index{Ruzsa, I.}\index{Green, B.}
Let $\kappa$ be the exponent of $G$.  Then
$$\mu(G, \{2,1\})={\mu}(\mathbb{Z}_{\kappa}, \{2,1\}) \cdot \frac{n}{\kappa}= v_{1}(\kappa,3) \cdot \frac{n}{\kappa}.$$
\end{thm}
We should mention that the proof in \cite{GreRuz:2005a} relies, in part, on a computer program.

For other $k$ and $l$,  the value of $\mu(G, \{k,l\})$ is not known in general, though we have the following bounds:

\begin{thm} [Bajnok; cf.~\cite{Baj:2009a}]  \label{Bajnok mu bounds}\index{Bajnok, B.}
Suppose that $G$ is an abelian group of order $n$ and exponent $\kappa$.  Then, for all positive integers $k$ and $l$ with $k>l$ we have
$$v_{k-l}(\kappa,k+l) \cdot \frac{n}{\kappa} \leq {\mu}(G, \{k,l\}) \leq v_1(n,k+l).$$
\end{thm}

\begin{prob}

Determine $\mu(G, \{k,l\})$ for arbitrary finite abelian groups.  In particular, find $\mu(G, \{3,1\})$ for all noncyclic groups $G$.
\end{prob}

We have the value of $\mu(G, \{k,l\})$ when $G$ is (isomorphic to) an elementary abelian 2-group.  Note that, when $k$ and $l$ have the same parity, then $ka=la$ for all $a \in \mathbb{Z}_2^r$, so $\mu(\mathbb{Z}_2^r, \{k,l\})=0$ in this case.  When $k$ and $l$ have opposite parity, then the set $\{1 \} \times \mathbb{Z}_2^{r-1}$ is $(k,l)$-sum-free in $\mathbb{Z}_2^r$, so $\mu(\mathbb{Z}_2^r, \{k,l\}) \geq 2^{r-1}$.  On the other hand, for every subset $A$ of $\mathbb{Z}_2^r$, both $kA$ and $lA$ have size at least $|A|$, so they cannot be disjoint when $A$ is $(k,l)$-sum-free and $|A| > 2^{r-1}$.  Therefore:

\begin{prop} \label{mu for 2-group}
For all positive integers $r, k$, and $l$ with $k>l$, we have
$$\mu(\mathbb{Z}_2^r, \{k,l\}) =
\left\{
\begin{array}{cl}
0 & \mbox{if $k \equiv l$ mod 2}; \\ \\
2^{r-1} & \mbox{if $k \not \equiv l$ mod 2}.
\end{array}\right.$$
\end{prop}

We can observe that, by Propositions \ref{mu for 2-group} and \ref{tau elementary}, we have $$\mu(\mathbb{Z}_2^r, \{k,l\}) =\tau(\mathbb{Z}_2^r, k+l)$$ for all parameters.

We now turn to the inverse problem of classifying all $(k,l)$-sum-free subsets $A$ of $G$ of maximum size $|A|={\mu}(G, \{k,l\})$.

Recall that arithmetic progressions play a fundamental role in providing $(k,l)$-sum-free sets (see Theorem \ref{enough for ap}) and, in particular, in providing sum-free sets.  With this in mind, we start by determining all sum-free arithmetic progressions in cyclic groups that have maximum size.

\begin{prop} \label{max sum-free AP} \label{ap max classify}
Suppose that $A$ is a sum-free arithmetic progression in $\mathbb{Z}_n$ of size $$|A|=\mu (\mathbb{Z}_n, \{2,1\})=v_1(n,3).$$  Then one of the following possibilities must hold:

\begin{enumerate}
  \item $n$ is even and $$A=\{1,3,\dots,n-1\}.$$
  \item $n$ is divisible by 3, has no prime divisors congruent to 2 mod 3, and there exists an integer $b$ relatively prime to $n$ for which
  
\begin{enumerate}
  \item $$b \cdot A=\{1,4,7,\dots,n-2\}$$ or
\item $$b \cdot A=\{n/3, n/3+1, \dots, 2n/3-1\}.$$
\end{enumerate}

\item $n$ is equal to a prime $p$ that is congruent to 2 mod 3, and there exists an integer $b$ relatively prime to $p$ for which
  
$$b \cdot A=\{(p+1)/3, (p+1)/3+1, \dots, 2(p+1)/3-1\}.$$

\item $n$ is congruent to 1 mod 3, has no prime divisors congruent to 2 mod 3, and there exists an integer $b$ relatively prime to $n$ for which
  
\begin{enumerate}
  \item $$b \cdot A=\{(n-1)/3, (n-1)/3+1, \dots, 2(n-1)/3-1\}$$ or
\item $$b \cdot A=\{(n-1)/3+1, (n-1)/3+2, \dots, 2(n-1)/3\}.$$
\end{enumerate}

\end{enumerate}    
\end{prop}

(We should add that we may assume that the integer $b$ of Case 2 (a) equals 1 or 2, or, equivalently, 1 or $-1$.)  The proof of Proposition \ref{max sum-free AP} starts on page \pageref{proof of max sum-free AP}.

To illuminate Proposition \ref{ap max classify}, we discuss an example in detail.  Recall that for each $n$, we mentioned that the elements between $n/6$ and $n/3$, together with the elements between $2n/3$ and $5n/6$, with the endpoints carefully chosen, form a sum-free set in $\mathbb{Z}_n$.  Is this an arithmetic progression, and, if it has size $v_1(n,3)$, is it included in Proposition \ref{ap max classify}?  The answers to both questions are affirmative.  For example, if $n=27$, then $v_1(n,3)=n/3=9$, and the sum-free set in question is 
$$A=\{5,6,7,8,9\} \cup \{19,20,21,22\}.$$  It is easy to verify that $A$ is sum-free in $\mathbb{Z}_{27}$, and we see that it is an arithmetic progression, since
$$A=\{5, 5+14, 5+2 \cdot 14, \dots, 5+8 \cdot 14\}.$$  We can verify that it is included in case 2 (b) of Proposition \ref{ap max classify}, by noting that $25$ and $27$ are relatively prime, and
$$25 \cdot A =\{9,10,11,12,13,14,15,16,17\}.$$

The classification of sum-free sets of maximum  size has been known for several decades for any $G$ whose order has a prime divisor congruent to 2 mod 3 or whose order is divisible by 3.  It turns out that all such sets are arithmetic progressions of cosets of a subgroup of $G$.  Namely, we have the following two results:

\begin{thm} [Diananda and Yap; cf.~\cite{DiaYap:1969a};\index{Diananda, P. H.}\index{Yap, H. P.} see also Theorem 7.8 in \cite{WalStrWal:1972a}] \label{inverse sum-free 1}\index{Wallis, W. D.}\index{Street, A. P.}\index{Wallis, J. C.} 
Suppose that the order $n$ of $G$ has prime divisors congruent to 2 mod 3, and $p$ is the smallest of them.  If $A$ is a sum-free set of size 
$$|A|=\mu (G, \{2,1\})=\left( 1 + \frac{1}{p} \right) \cdot \frac{n}{3}$$ in $G$, then there is a subgroup $H$ in $G$, so that $G/H$ is cyclic of order $p$, and $A$ is the union of $(p+1)/3$ cosets of $H$ that form an arithmetic progression.

\end{thm}  

It is worth noting that, as a special case of Theorem \ref{inverse sum-free 1}, we get cases 1 and 3 of Proposition \ref{ap max classify} when $G$ is cyclic of even order or cyclic of prime order $p$ (with $p \equiv 2$ mod 3), respectively.

\begin{thm} [Street; cf.~\cite{Str:1972a, Str:1972b}; see also Theorem 7.9 in \cite{WalStrWal:1972a}] \label{inverse sum-free 2}\index{Street, A. P.}\index{Wallis, W. D.}\index{Wallis, J. C.} 
Suppose that the order $n$ of $G$ is divisible by 3 and has no prime divisors congruent to 2 mod 3.  If $A$ is a sum-free set of size 
$$|A|=\mu (G, \{2,1\})=n/3 $$ in $G$, then there is a divisor $k$ of $n/3$ and a subgroup $H$ in $G$, so that $G/H$ is cyclic of order $3k$, and $A$ is the union of $k$ cosets of $H$ that form an arithmetic progression.

\end{thm}
This time, observe that when $G$ is cyclic, then for $k=1$ and $k=n/3$ we get cases 2 (a) and 2 (b) of Proposition \ref{ap max classify}, respectively. 

This leaves us with the task of characterizing sum-free sets of maximum size in $G$ when the order of $G$ has only prime divisors congruent to 1 mod 3.  As it turns out, in this case not all such sets are arithmetic progressions of cosets of a subgroup, though the one other type of set is not far from it.  We first present the result for cyclic groups:

\begin{thm} [Yap; cf.~\cite{Yap:1971a}] \label{inverse sum-free 3 cyclic}\index{Yap, H. P.}
Suppose that all prime divisors of $n$ are congruent to 1 mod 3.  If $A$ is a sum-free set of size 
$$|A|=\mu (\mathbb{Z}_n, \{2,1\})=(n-1)/3$$ in $\mathbb{Z}_n$, then there exists an integer $b$ relatively prime to $n$ for which 
  $$b \cdot A=\{(n-1)/3, (n-1)/3+1, \dots, 2(n-1)/3-1\},$$
$$b \cdot A=\{(n-1)/3+1, (n-1)/3+2, \dots, 2(n-1)/3\},$$ or
$$b \cdot A=\{(n-1)/3, (n-1)/3+2, (n-1)/3+3, \dots, 2(n-1)/3-1, 2(n-1)/3+1\}.$$

\end{thm}  
Observe that the third set listed is two elements short of being an arithmetic progression; note also that this possibility can only occur if $n > 7$.

The classification for noncyclic groups is  presented in the article \cite{BalPraRam:2016a} by Balasubramanian, Prakash, and Ramana.\index{Balasubramanian, R.}\index{Prakash, G.}\index{Ramana, D. S.}  Rather than stating the exact result, which is complicated, we only discuss the situation in groups of rank 2.  
\label{max sum-free classification  type 3 noncyclic}

How can one find sum-free sets of maximum size in $\mathbb{Z}_{n_1} \times \mathbb{Z}_{n_2}$?  The most obvious idea is to follow Proposition \ref{mufordirectsum}; for example, if $A$ is a maximum size sum-free set in $\mathbb{Z}_{n_2}$, then, by Theorem \ref{(2,1)all}, $$\mathbb{Z}_{n_1} \times A$$ is   
a maximum size sum-free set in $\mathbb{Z}_{n_1} \times \mathbb{Z}_{n_2}$.  The following construction is less obvious: we can verify that, with an arbitrary subgroup $H$ of $\mathbb{Z}_{n_1}$, the set $A=A_1 \cup A_2 \cup A_3$, where
\begin{eqnarray*}
A_1 & = & H \times \{(n_2-1)/3\} \\
A_2 & = & (\mathbb{Z}_{n_1} \setminus H) \times \{2 (n_2-1)/3\} \\
A_3 & = & \mathbb{Z}_{n_1} \times \{(n_2-1)/3+1, (n_2-1)/3 +2, \dots, 2(n_2-1)/3-1\}
\end{eqnarray*}
is a maximum size sum-free set in $\mathbb{Z}_{n_1} \times \mathbb{Z}_{n_2}$.  To see that $A$ is sum-free, 
observe that, considering the components in $\mathbb{Z}_{n_2}$ alone, for elements $g_1, g_2, g_3 \in A$ to satisfy $g_1+g_2=g_3$, we must have $g_1, g_2 \in A_1$ and $g_3 \in A_2$, but then the first component of $g_1+g_2$ is in $H$ while that is not the case for $g_3$.   

The following construction is similar: let $H$ again be an arbitrary subgroup of $\mathbb{Z}_{n_1}$, and let $A=A_1 \cup A_2 \cup A_3$, where
\begin{eqnarray*}
A_1 & = & H \times \{(n_2-1)/3, 2(n_2-1)/3+1\} \\
A_2 & = & (\mathbb{Z}_{n_1} \setminus H) \times \{(n_2-1)/3+1, 2 (n_2-1)/3\} \\
A_3 & = & \mathbb{Z}_{n_1} \times \{(n_2-1)/3+2, (n_2-1)/3 +3, \dots, 2(n_2-1)/3-1\}.
\end{eqnarray*}
We can again verify that $A$ is a maximum size sum-free set in $\mathbb{Z}_{n_1} \times \mathbb{Z}_{n_2}$.  

The result of Balasubramanian, Prakash, and Ramana\index{Balasubramanian, R.}\index{Prakash, G.}\index{Ramana, D. S.} in \cite{BalPraRam:2016a} says that the three types of sets just describe provide ``essentially'' the only types of sum-free sets in $\mathbb{Z}_{n_1} \times \mathbb{Z}_{n_2}$ of maximum size---see \cite{BalPraRam:2016a} for the precise statement for groups of arbitrary rank.  (We should add that the result for elementary abelian groups was proved by Street in 1971; see \cite{RheStr:1971a}\index{Rhemtulla, A. H.}\index{Street, A. P.} or Theorem 7.21 in \cite{WalStrWal:1972a}.)\index{Wallis, W. D.}\index{Street, A. P.}\index{Wallis, J. C.} 

Let us turn now to the classification of maximum-size $(k,l)$-sum-free sets for $(k,l) \neq (2,1)$.  While we may suspect that the sets come in even more variety than they did for sum-free sets, we find that this is not the case when the group is of prime order: all maximum-size sets are arithmetic progressions:

\begin{thm} [Plagne; cf.~\cite{Pla:2002a}] \label{plagne AP}\index{Plagne, A.} 
Let $p$ be a positive prime, and $k$ and $l$ be positive integers with $k>l$ and $k \geq 3$; assume also that $p$ does not divide $k-l$.  If $A$ is a $(k,l)$-sum-free set of size 
$$|A|=\mu (\mathbb{Z}_p, \{k,l\})=\left \lfloor \frac{p-2}{k+l} \right \rfloor +1$$ in $\mathbb{Z}_p$, then $A$ is an arithmetic progression.
\end{thm}  
(We note that Theorem \ref{plagne AP} was conjectured and proved in part by Bier and Chin in \cite{BieChi:2001a}.)\index{Bier, T.}\index{Chin, A. Y. M.} 

We can go a step beyond the statement of Theorem \ref{plagne AP} and determine the maximum-size $(k,l)$-sum-free sets in $\mathbb{Z}_p$ more explicitly:

\begin{thm} \label{explicit plagne AP}
Let $p$ be a positive prime, and $k$ and $l$ be positive integers with $k>l$ and $k \geq 3$; assume also that $p$ does not divide $k-l$.  We let $M$ and $r$ denote, respectively, the quotient and the remainder of $p-2$ when divided by $k+l$.  For $i=1,2,\dots, \lfloor r/2 \rfloor +1$, let $a_i$ be the (unique) solution of $$(k-l)a_i=l M +i$$ in $\mathbb{Z}_p$.  If A is a $(k,l)$-sum-free set of size 
$$|A|=\mu (\mathbb{Z}_p, \{k,l\})=M+1$$ in $\mathbb{Z}_p$, then there is an $i \in \{1,2,\dots, \lfloor r/2 \rfloor +1\}$ and an integer $b$ relatively prime to $p$ for which $$b \cdot A =\{a_i, a_i+1, \dots, a_i+M \}.$$

\end{thm}

As an example, let us consider the case of maximum-size $(4,1)$-sum-free sets in $\mathbb{Z}_{19}$.  Note that $M=\left \lfloor (19-2)/(4+1) \right \rfloor = 3$, so $\mu (\mathbb{Z}_{19}, \{4,1\})=4$, and that $r=2$ and thus $\lfloor r/2 \rfloor +1=2$.  Furthermore, solving the equations $$(4-1) \cdot a_i=1 \cdot 3+i$$ for $i=1,2$ in $\mathbb{Z}_{19}$, we find that $a_1=14$ and $a_2=8$.  Therefore, Theorem \ref{explicit plagne AP} claims that any $(4,1)$-sum-free set of maximum size 4 in $\mathbb{Z}_{19}$ is a dilate of either $\{14,15,16,17\}$ or $\{8,9,10,11\}$.  The easy proof---relying, of course, on Theorem \ref{plagne AP}---is on page \pageref{proof of explicit plagne AP}.

We do not have a classification of maximum-size $(k,l)$-sum-free sets with $k \geq 3$ in groups of composite order:

\begin{prob}
Let $n$,  $k$, and $l$ be positive integers with $k>l$ and $k \geq 3$.  Classify all $(k,l)$-sum-free sets of size $\mu (\mathbb{Z}_n, \{k,l\})$ in the cyclic group $\mathbb{Z}_n$.

\end{prob}

A particularly special $(k,l)$-sum-free set is one where no elements are ``wasted'': we say that a $(k,l)$-sum-free set $A$ in $G$ is {\em complete} if $kA$ and $lA$ partition $G$; in other words, not only do we have $kA \cap lA = \emptyset$ (i.e. $A$ is $(k,l)$-sum-free), but also $kA \cup lA=G$.  

Before going any further, a word of caution: While the definition of $A \subseteq G$ being $(k,l)$-sum-free can be given either as $kA \cap lA = \emptyset$ or as $kA-lA \subseteq G \setminus \{0\}$, the notions that for a $(k,l)$-sum-free set $kA \cup lA=G$ holds or that $kA-lA = G \setminus \{0\}$ holds are not equivalent (though both notions express some sort of ``un-wastefulness'')!  We have counter-examples in both directions: 
\begin{itemize}
  \item If $A=\{1,3,5,7,9\} \subseteq \mathbb{Z}_{10}$, then $2A=\{0,2,4,6,8\}$, so $A$ is sum-free set in $\mathbb{Z}_{10}$; $A \cup 2A=G$, thus $A$ is complete, yet $2A-A=A \neq \mathbb{Z}_{10} \setminus \{0\}$.
\item If $A=\{1,3,9\} \subseteq \mathbb{Z}_{13}$, then $2A=\{2,4,5,6,10,12\}$, so $A$ is sum-free in $\mathbb{Z}_{13}$; $2A-A=\mathbb{Z}_{13} \setminus \{0\}$, yet $A \cup 2A \neq G$ thus $A$ is not complete.
\end{itemize}  
While $(k,l)$-sum-free sets $A$ in $G$ for which $kA-lA = G \setminus \{0\}$ may also be of interest, we here define them as complete when $kA \cup lA=G$.

We may use our classification of maximum-size sum-free sets (Theorems \ref{inverse sum-free 1}, \ref{inverse sum-free 2}, and \ref{inverse sum-free 3 cyclic}) to determine which of these sets are complete sum-free sets.  Consider first the case when the order $n$ of $G$ has prime divisors congruent to 2 mod 3; let $p$ be the smallest such divisor.  From Theorem \ref{inverse sum-free 1} we see that a maximum-size sum-free set in $G$ has size $(p+1)/3 \cdot n/3$ and is of the form
$$A=(a+H) \cup (a+d+H) \cup \cdots \cup (a + (p-2)/3d+H)$$ for a subgroup $H \leq G$ of order $n/p$ and elements $a, d \in G$.   Then 
$$2A=(2a+H) \cup (2a+d+H) \cup \cdots \cup (2a + 2(p-2)/3d+H),$$ so $2A$ has size $(2p-1)/3 \cdot n/p$.  Therefore, $|A|+|2A|=n$, which (since $A$ and $2A$ are disjoint) can only be if $A \cup 2A=G$, hence $A$ is complete.

Next, assume that $n$ is divisible by 3, but has no divisors congruent to 2 mod 3.   By Theorem \ref{inverse sum-free 2}, a maximum-size sum-free set in $G$ has size $n/3$ and is of the form
$$A=(a+H) \cup (a+d+H) \cup \cdots \cup (a+(k-1)d+H)$$ for a divisor $k$ of $n/3$, a subgroup $H \leq G$ of order $n/(3k)$, and elements $a, d \in G$.   Then 
$$2A=(2a+H) \cup (2a+d+H) \cup \cdots \cup (2a + (2k-2)d+H),$$ so $2A$ has size $(2k-1) \cdot n/(3k)$.  Therefore, $|A|+|2A|=(3k-1)/(3k) \cdot n < n$, so $A$ is not complete.

Suppose now that $n$ has only prime divisors congruent to 1 mod 3.  By Theorem \ref{inverse sum-free 3 cyclic}, if $A$ is a sum-free set of maximum size 
$(n-1)/3$ in $\mathbb{Z}_n$, then there exists an integer $b$ relatively prime to $n$ for which $b \cdot A$ is equal to 
  $$A_1=\{(n-1)/3, (n-1)/3+1, \dots, 2(n-1)/3-1\},$$
$$A_2=\{(n-1)/3+1, (n-1)/3+2, \dots, 2(n-1)/3\},$$ or
$$A_3=\{(n-1)/3, (n-1)/3+2, (n-1)/3+3, \dots, 2(n-1)/3-1, 2(n-1)/3+1\}.$$
Since we have $|b \cdot A|=|A|$ and $|2 (b \cdot A)|=|2A|$, it suffices to examine sets $A_1$, $A_2$, and $A_3$.  We get  
$$2A_1=\{2(n-1)/3, 2(n-1)/3+1, \dots, 4(n-1)/3-2\},$$
$$2A_2=\{2(n-1)/3+2, 2(n-1)/3+3, \dots, 4(n-1)/3\},$$ and
$$2A_3=\{2(n-1)/3\} \cup \{ 2(n-1)/3+2, 2(n-1)/3+3, \dots, 4(n-1)/3 \} \cup \{4(n-1)/3+2\}.$$
(All elements listed are considered mod $n$, of course.)  Therefore, $2A_1$ and $2A_2$ both have size $2(n-1)/3-1$, so $A_1$ and $A_2$ are not complete.  
However, $|2A_3|=2(n-1)/3+1$, so $|A_3+|2A_3|=n$ and $A_3$ is complete.

We can summarize our findings as follows:

\begin{thm}
Let $A$ be a maximum-size sum-free set in $G$; as usual, let $|G|=n$.

\begin{enumerate}
  \item If $n$ has a prime divisor congruent to 2 mod 3, then $A$ is a complete sum-free set in $G$.
  \item If $n$ is divisible by 3 but has no prime divisors congruent to 2 mod 3, then $A$ is not a complete sum-free set in $G$.
\item If $n$ has only prime divisors congruent to 1 mod 3 and $G$ is cyclic, then $A$ is a complete sum-free set in $G$ if, and only if, there is an integer $b$ relatively prime to $n$ for which 
$$b \cdot A=\{(n-1)/3, (n-1)/3+2, (n-1)/3+3, \dots, 2(n-1)/3-1, 2(n-1)/3+1\}.$$
\end{enumerate}

\end{thm}

This leaves us with the following problem:

\begin{prob} \label{complete sum-free noncyclic}
For each noncyclic group $G$ whose order has only prime divisors congruent to 1 mod 3, determine which sum-free sets $A$ in $G$ of maximum size $|A|=(n-1)/3$ are complete.

\end{prob}
Problem \ref{complete sum-free noncyclic} is probably not difficult given the classification of the sum-sets in question that we presented above (see page \pageref{max sum-free classification  type 3 noncyclic}).

Turning to the case when $k \geq 3$, we use Theorem \ref{explicit plagne AP} to analyze maximum $(k,l)$-sum-free sets in cyclic groups of prime order $p$.  According to Theorem \ref{explicit plagne AP}, it suffices to analyze
$$A=\{a_i, a_i+1, \dots, a_i+M\}$$ where $M=\lfloor (p-2)/(k+l) \rfloor$, $r=p-2-(k+l)M$, $i=1,2,\dots, \lfloor r/2 \rfloor +1$,  and $a_i$ is the unique solution to the equation $(k-l)a_i=lM+i$ in $\mathbb{Z}_p$.

We find that 
\begin{eqnarray*}
kA & = & \{ka_i, ka_i+1, \dots, ka_i+kM\} \\
lA& = & \{la_i,la_i+1,\dots,la_i+lM\}.
\end{eqnarray*}  Now by definition, $A$ is complete when $$|kA|+|lA|=(kM+1)+(lM+1)=p,$$ and this holds if, and only if, $M=(p-2)/(k+l)$, that is, $p-2$ is divisible by $k+l$.  In this case, $r=0$, and thus $i$ can only equal 1. 
 We get:

\begin{thm}
Let $p$ be a positive prime, and $k$ and $l$ be positive integers with $k>l$ and $k \geq 3$.  Let $A$ be a maximum-size $(k,l)$-sum-free set in $\mathbb{Z}_p$.  Then $A$ is complete if, and only if, $p-2$ is divisible by $k+l$, and $A$ is a dilate of the set $$A=\{a, a+1, \dots, a+M\},$$ where $M=(p-2)/(k+l)$ and $a$ is the unique solution to the equation $(k-l)a=lM+1$ in $\mathbb{Z}_p$.

\end{thm}

As an example, we can compute explicitly that, when $p-2$ is divisible by 5, then any complete $(4,1)$-sum-free set in $\mathbb{Z}_p$ is a dilate of the set
$$\{(2p+1)/5, (2p+1)/5+1, \dots, (3p-1)/5\}.$$
The similar question is not yet solved in cyclic groups of composite order:
 
\begin{prob}
Let $n$,  $k$, and $l$ be positive integers with $k>l$ and $k \geq 3$.  Classify all complete $(k,l)$-sum-free sets of maximum size $\mu (\mathbb{Z}_n, \{k,l\})$ in the cyclic group $\mathbb{Z}_n$.

\end{prob}

Of course, a complete $(k,l)$-sum-free set need not have maximum size $\mu (\mathbb{Z}_n, \{k,l\})$.  For example, when $p$ is any (not necessarily the smallest) prime divisor of $n$ with $p \equiv 2$ mod 3, then the set
$$A=\{(p+1)/3 + i +pj \mid i = 0, 1, \dots, (p-2)/3; j=0,1,\dots,n/p-1\}$$ is a complete sum-free set in $\mathbb{Z}_n$, since   
$$2A=\{2(p+1)/3 + i +pj \mid i = 0, 1, \dots, (2p-4)/3; j=0,1,\dots,n/p-1\}$$ and thus $A \cap 2A=\emptyset$ but $A \cup 2A=\mathbb{Z}_n$. 

The work \cite{HavLev:2017a} of Haviv and Levy\index{Haviv, I.}\index{Levy, D.} generated many other examples for complete sum-free sets that also have the additional property that they are symmetric (that is, subsets $A$ of $G$ for which $A=-A$).  Their results are somewhat complicated, so we just mention the  following special case:

\begin{thm} [Haviv and Levy; cf.~\cite{HavLev:2017a}]

Let $A$ be a complete sum-free subset of $\mathbb{Z}_p$ for some prime $p$, and suppose that $A$ has size $\mu(\mathbb{Z}_p, \{2,1\})-2$.  If $p$ is sufficiently large, then there is an element $b \in \mathbb{Z}_p \setminus \{0\}$ for which:

\begin{itemize}
  \item If $p \equiv 2$ mod 3, then
  $$b \cdot A=\left \{\frac{p-5}{3}+i \mid i \in T \right \} \cup \left[\frac{p+13}{3}, \frac{2p-13}{3} \right] \cup \left \{\frac{2p+5}{3}-i \mid i \in T \right \}$$
  where $T=\{0,2,4\}$ or $T=\{0,3,4\}$;
  \item If $p \equiv 1$ mod 3, then
  $$b \cdot A=\left \{\frac{p-7}{3}+i \mid i \in T \right \} \cup \left[\frac{p+17}{3}, \frac{2p-17}{3} \right] \cup \left \{\frac{2p+7}{3}-i \mid i \in T \right \}$$
  where $T=\{0,2,4,6\}$, $T=\{0,3,5,6\},  $T=\{0,4,5,6\}, or $T=\{1,2,6,7\}$.
\end{itemize}

\end{thm}

The classification of all complete $(k,l)$-sum-free sets remains open:

\begin{prob} \label{prob complete not max}
Find all complete sum-free---or, more generally, complete $(k,l)$-sum-free---sets in finite abelian groups that do not have maximum size $\mu (\mathbb{Z}_n, \{k,l\})$.

\end{prob}

Actually, Problem \ref{prob complete not max} stays intriguing even when dropping the requirement that the set be complete: we may just ask for $(k,l)$-sum-free sets that do not have maximum size $\mu (\mathbb{Z}_n, \{k,l\})$ but are {\em maximal}---that is, they cannot be enlarged without losing the $(k,l)$-sum-free property.  One then may ask for the possible cardinalities of maximal sum-free sets:

\begin{prob}
For each group $G$ and for all positive integers $k, l$ with $k>l$, find all possible values of $m$ for which a maximal $(k,l)$-sum-free set of size $m$ exist.

\end{prob}

In particular, for positive integers $i$, we let $M_i(G,\{k,l\})$ denote the $i$-th largest cardinality of a maximal $(k,l)$-sum-free set in $G$ (for values of $i$ for which this exists), and $M (G, \{k,l\})$ denote the size of the smallest maximal $(k,l)$-sum-free set in $G$.  Note that
$$M_1(G,\{k,l\})= \mu (G, \{k,l\}).$$

We then may ask:

\begin{prob}
Find the $i$-th largest size $M_i (G, \{k,l\})$ of maximal $(k,l)$-sum-frees set in $G$.  In particular, find $M_2 (G, \{k,l\})$.
\end{prob}

\begin{prob}
Find the size $M (G, \{k,l\})$ of the smallest maximal $(k,l)$-sum-free set in $G$.
\end{prob}

As an example, consider $G=\mathbb{Z}_{11}$, where $\mu (\mathbb{Z}_{11}, \{2,1\})=4$.  The set $A=\{1,3,5\}$ is sum-free in $G$, but $A \cup \{a\}$ is not sum-free for any $a \in G \setminus A$.  We can also verify that each set of size two can be enlarged so that it stays sum-free,  hence $$M_2(\mathbb{Z}_{11}, \{2,1\})=M(\mathbb{Z}_{11}, \{2,1\}) = 3.$$ 

Clark and Pedersen have the following result in elementary abelian 2-groups:

\begin{thm} [Clark and Pedersen; cf.~\cite{ClaPed:1992a}] \label{Clark and Pedersen}\index{Clark, W. E.}\index{Pedersen, J.} 
For each positive integer $r \geq 4$, we have $$M_2(\mathbb{Z}_{2}^r, \{2,1\})=5 \cdot 2^{r-4}.$$

\end{thm}

Clark and Pedersen\index{Clark, W. E.}\index{Pedersen, J.} also computed all values of $M_i(\mathbb{Z}_{2}^r, \{2,1\})$ (that exist) and $M(\mathbb{Z}_{2}^r, \{2,1\})$ for $r \leq 6$ and possibly all for $7 \leq r \leq 10$.  From their data, we may conjecture the following:

\begin{conj} \label{conj based on Clark and Pedersen}
For each positive integer $r \geq 4$ and $2 \leq i \leq r-2$, we have $$M_i(\mathbb{Z}_{2}^r, \{2,1\})=2^{r-2}+2^{r-2-i},$$ and 
$$M(\mathbb{Z}_{2}^r, \{2,1\})=2^{\lceil r/2 \rceil}+2^{\lfloor r/2 \rfloor}-3.$$

\end{conj}

We note that our formula does not hold for $i=1$, but by Theorem \ref{Clark and Pedersen}, it holds for all $r \geq 4$ and $i=2$.  Furthermore, from Corollary 5 in \cite{ClaPed:1992a} (with $s=2$ and $t=r-2-i$) we also know that maximal sum-free sets of size $2^{r-2}+2^{r-2-i}$ exist in $\mathbb{Z}_{2}^r$ for each $r \geq 4$ and $2 \leq i \leq r-2$, and (with $s=\lfloor r/2 \rfloor$ and $t=0$) that maximal sum-free sets of size $2^{\lceil r/2 \rceil}+2^{\lfloor r/2 \rfloor}-3$ exist in $\mathbb{Z}_{2}^r$ for each $r \geq 4$.  However,  it has not been established yet that these values equal $M_i(\mathbb{Z}_{2}^r, \{2,1\})$ and $M(\mathbb{Z}_{2}^r, \{2,1\})$, respectively.

\begin{prob}
Prove Conjecture \ref{conj based on Clark and Pedersen}.

\end{prob}

We now turn to a related problem that was first investigated by Erd\H{o}s in \cite{Erd:1965a}.\index{Erd\H{o}s, P.}
Namely, we wish to examine the size of the largest $(k,l)$-sum-free subset of a given $A \subseteq G$: we denote this quantity by $\mu (G, \{k,l\}, A)$.  We can then ask for the ``worst-case scenario'': For each positive integer $m \leq n$, find the $m$-subsets $A_0$ of $G$ so that $$\mu (G, \{k,l\}, A_0) \leq \mu (G, \{k,l\}, A)$$ for all $m$-subsets $A$ of $G$; we then let $\mu (G, \{k,l\}, m)=\mu (G, \{k,l\}, A_0)$.

\begin{prob}  \label{Find mu G k,l,m}
For each group $G$ and all positive integers $k$, $l$, and $m$ (with  $k > l$ and $m \leq n$), find $\mu (G, \{k,l\}, m)$.

\end{prob}

Of course, we have $$\mu (G, \{k,l\}, n)=\mu (G, \{k,l\}),$$ but values of $\mu (G, \{k,l\}, m)$ for $m < n$ are not known (cf.~the preprint \cite{TaoVu:2016a} by Tao and Vu\index{Vu, V. H.}\index{Tao, T.} for the history of this problem and for more information).  We have the following general lower bound for the case of sum-free sets:

\begin{thm} [Alon and Kleitman; cf.~\cite{AloKle:1990a}]\index{Alon, N.} \index{Kleitman, D. J.}
For every group $G$ and positive integer $m$ (with $m \leq n$) we have $$\mu (G, \{2,1\}, m) \geq \lfloor 2m/7 \rfloor.$$
\end{thm}
We note that the bound here is tight in the sense that, by Theorem \ref{(2,1)all}, it is achieved by the group $\mathbb{Z}_7^r$ when $m=n=7^r$.

As a step toward solving Problem \ref{Find mu G k,l,m}, we may ask for the smallest value of $m$ for which $\mu (G, \{k,l\}, m)=\mu (G, \{k,l\})$:

\begin{prob}  \label{prob T}
For each group $G$ and all positive integers $k$ and $l$with  $k > l$, find the smallest value $T(G, \{k,l\})$ of $m$ for which $\mu (G, \{k,l\}, m)=\mu (G, \{k,l\})$.

\end{prob} 

When $p$ is an odd prime and $p \equiv 2$ mod 3, then by Theorem \ref{inverse sum-free 1} and Proposition \ref{ap max classify} (part 3), we know that $\mathbb{Z}_p$ has exactly $(p-1)/2$ sum-free subsets of size $\mu (\mathbb{Z}_p, \{2,1\}) =(p+1)/3$.  Furthermore, by the result of Chin\index{Chin, A. Y. M.} in \cite{Chi:2001a}, these sets form a {\em block design}; that is, every nonzero element of $\mathbb{Z}_p$ is contained in the same number (in our case, $(p+1)/6$) of these $(p-1)/2$ sets.   Since for $p > 5$, the value of $2 \cdot (p+1)/6$ is less than $(p-1)/2$, deleting two nonzero elements of $\mathbb{Z}_p$ results in a set of size $p-2$ that still contains one of the $(p-1)/2$ sum-free sets of maximum size.  Therefore, we have $$\mu (\mathbb{Z}_p, \{2,1\}, p-2)=\mu (\mathbb{Z}_p, \{2,1\}),$$ which yields:

\begin{thm} \label{block des T}
For an odd prime $p>5$ with $p \equiv 2$ mod 3, we have $$T(\mathbb{Z}_p, \{2,1\}) \leq p-2.$$

\end{thm} 

As an example, let us consider the group $\mathbb{Z}_{11}$; by Theorem \ref{block des T},  $T(\mathbb{Z}_{11}, \{2,1\})$ is at most 9.  We know that $\mu (\mathbb{Z}_p, \{2,1\})=4$ and that there are five sum-free subsets of size four:
$$\{4,5,6,7\}, \; \{8,10,1,3\}, \; \{1,4,7,10\}, \; \{5,9,2,6\}, \; \{9,3,8,2\}.$$
We see that none of these five sets are contained in $$A=\mathbb{Z}_{11} \setminus \{2,4,8\},$$ and thus we get $T(\mathbb{Z}_{11}, \{2,1\}) = 9$.  We do not know the value of $T(\mathbb{Z}_p, \{2,1\})$ in general, so, as a special case of Problem \ref{prob T},  we ask:

\begin{prob}
Find the value of $T(\mathbb{Z}_p, \{2,1\})$ for every odd prime $p$ with $p \equiv 2$ mod 3.

\end{prob}

Our next pursuit is the following related problem.
 For a given $m$-subset $A$ of a given group $G$, we may ask for the number of ordered pairs $(a_1,a_2)$ with $a_1, a_2 \in A$ for which $a_1+a_2 \in A$; we denote this quantity by $P(G,A)$.  (Note that, when $a_1+a_2 \in A$, then $a_2 +a_1 \in A$, so when $a_1$ and $a_2$ are distinct, both $(a_1,a_2)$ and $(a_2,a_1)$ contribute toward the count.)  We then may ask for the minimum value of $P(G,A)$ among all $m$-subsets $A$ of $G$; we  denote this quantity by $P(G,m)$.   Observe that, by definition, we have $P(G, m)=0$ for each $m \leq \mu (G, \{2,1\})$, but $P(G,m) \geq 1 $ for each $m \geq \mu (G, \{2,1\})+1$.  

\begin{prob}  \label{P(G,m)}

For each group $G$ and all positive integers $m$ with 
$$\mu (G, \{2,1\})+1 \leq m \leq n,$$ find $P(G, m)$.  

\end{prob}

Consider, as an example, the cyclic group $\mathbb{Z}_p$ of prime order $p \geq 3$, where $$\mu (\mathbb{Z}_p, \{2,1\})=\lfloor (p+1)/3 \rfloor.$$  As we have already observed, the ``middle'' third of the elements, namely, the set
$$\{\lfloor (p+1)/3 \rfloor ,  \lfloor (p+1)/3 \rfloor+1, \dots, 2 \lfloor (p+1)/3 \rfloor-1 \}$$ forms a sum-free set of $\mathbb{Z}_p$ of maximum size.
 Suppose that $m$ is a positive integer so that $$m > \lfloor (p+1)/3 \rfloor.$$ We enlarge our set to
$$A=\left \{ \lfloor (p+1-m)/2 \rfloor , \lfloor (p+1-m)/2 \rfloor  +1, \dots, \lfloor (p+1-m)/2 \rfloor+m-1 \right \}.$$  Then $A$ has size $m$, and it is still in the ``middle'' (when $m$ is even, it is symmetrical in that its first and last elements add to $p$, and if $m$ is odd, then it is as close to being symmetrical as possible).
A short calculation finds that $$P(\mathbb{Z}_p,A)=\left \lfloor (3m-p)^2/4 \right \rfloor.$$       
A recent result by Samotij and Sudakov says that one cannot do better:

\begin{thm} [Samotij and Sudakov; cf.~\cite{SamSud:2016a}, \cite{SamSud:2016b}]\index{Samotij, W.}\index{Sudakov, B.} 
Suppose that $p$ is an odd prime, and $A$ is a subset of $\mathbb{Z}_p$ of size $m > \lfloor (p+1)/3 \rfloor.$  Then $$P(\mathbb{Z}_p,m)=\left \lfloor (3m-p)^2/4 \right \rfloor.$$ Furthermore, $P(\mathbb{Z}_p,A)=P(\mathbb{Z}_p,m)$ holds if, and only if, there is an integer $b$, relatively prime to $p$, for which
$$b \cdot A =\left \{ \lfloor (p+1-m)/2 \rfloor , \lfloor (p+1-m)/2 \rfloor  +1, \dots, \lfloor (p+1-m)/2 \rfloor+m-1 \right \}.$$ 
\end{thm}
(We note that the statement of this result is stated erroneously in \cite{SamSud:2016a}; see \cite{SamSud:2016b} for the corrected version.)

In \cite{SamSud:2016a} the authors evaluate $P(G,m)$ for the elementary abelian group $G$, and state some partial results for other groups, but warn that finding $P(G,m)$ for all $G$ and $m$ ``would be rather difficult.''

As a generalization, for a positive integer $k$ we define
$$P(G,k,A)=|\{(a_1,\dots,a_k) \in A^k \mid a_1+\cdots a_k \in A\}|$$ and
$$P(G,k,m)= \min \{ P(G,k,A) \mid A \in G, |A|=m\}.$$

Note that we have $P(\mathbb{Z}_p,2,A)=P(G,A)$ and $P(\mathbb{Z}_p,2,m)=P(G,m)$.

\begin{prob} \label{P(G,k,m)}
Evaluate $P(G,k,m)$ for all groups $G$ and positive integers $k$ and $m$. 

\end{prob}
 Of course, partial progress towards Problem \ref{P(G,k,m)}, such as the following:

\begin{prob}
Evaluate  $P(\mathbb{Z}_p,k,m)$ for prime values of $p$ and $k>l$ with $k \geq 3$. 

\end{prob}

\subsection{Limited number of terms} \label{6maxUlimited}

Given a group $G$ and a nonnegative integer $s$, here we are interested in finding the maximum possible size of an $[0,s]$-sum-free set in $G$, that is, the quantity $$\mu (G,[0,s]) =\mathrm{max} \{ |A|  \mid  A \subseteq G; 0 \leq l < k \leq s \Rightarrow kA \cap lA = \emptyset \}.$$  

Note that in a group $G$ of exponent $\kappa$, we have $0g=\kappa g$ for any $g \in G$, so no subset of $G$ is $[0,s]$-sum-free if $s \geq \kappa$.  Furthermore, the cases of $s=0$ and $s=1$ are easy: any subset $A$ of $G$ is $[0,0]$-sum-free, and $A \subseteq G$ is $[0,1]$-sum-free if, and only if, $0 \not \in A$.  We thus have:

\begin{prop} \label{mu [0,s] triv}
If $G$ is of exponent $\kappa$ and $s \geq \kappa$, then  $\mu (G,[0,s])=0.$  Furthermore, $\mu (G, [0,0])=n$ and $\mu (G, [0,1])=n-1$.

\end{prop}

According to Proposition \ref{mu [0,s] triv}, we can assume that $2 \leq s \leq \kappa-1$.

We have two explicit constructions for $[0,s]$-sum-free sets in $\mathbb{Z}_n$.  The first comes from Zajaczkowski\index{Zajaczkowski, C.}  (see \cite{Zaj:2013a}): Let $$m=\left \lfloor \frac{n-s-1}{s^2} \right \rfloor +1,$$ and consider
$$A=\{1 + is \mid i=0,1,\dots,m-1\}$$ in $\mathbb{Z}_n$.  Note that the elements of $hA$ are congruent to $h$ mod $s$; therefore, $kA$ and $lA$ are pairwise distinct for all $1 \leq l < k \leq s$.  Furthermore,  $0 \not \in sA$, since $$s \cdot (1+(m-1)s) = s+ s^2 \cdot \left \lfloor \frac{n-s-1}{s^2} \right \rfloor \leq s+ (n-s-1) <n.$$  Therefore, $A$ is indeed  $[0,s]$-sum-free set in $\mathbb{Z}_n$.

For our second construction, let $m$ be as above, and set $$a=(s-1)(m-1)+1.$$  Then for 
$$A=\{a , a+1, \dots, a+m-1\}$$ we see that $$hA=\{ha, ha+1, \dots, ha + h(m-1)\}.$$  
Note that for each $h=0,1,2,\dots,s-1$ we have 
$$ha + h(m-1) < (h+1)a,$$ since this inequality is equivalent to 
$h(m-1) <a$, which holds as
$$h(m-1) \leq (s-1)(m-1) = a-1.$$ 
Furthermore, $$sa+s(m-1) = s((s-1)(m-1)+1) +s(m-1) = s^2 (m-1)+s = s^2 \cdot \left \lfloor \frac{n-s-1}{s^2} \right \rfloor +s <n.$$
Therefore, $A$ is $[0,s]$-sum-free set in $\mathbb{Z}_n$.

From both constructions we get:

\begin{prop} [Zajaczkowski; cf.~\cite{Zaj:2013a}] \label{Claire prop}\index{Zajaczkowski, C.} 
For all positive integers $n$ and $s$, we have $$\mu (\mathbb{Z}_n,[0,s]) \geq \left \lfloor \frac{n-s-1}{s^2} \right \rfloor +1.$$ 

\end{prop}

Recall that when the exponent (in this case, the order) of the group $G$ is $s$ or less, then no set is $[0,s]$-sum-free in $G$.  Observe that $$\left \lfloor \frac{n-s-1}{s^2} \right \rfloor +1$$ equals 0 when $n \leq s$ but is positive when $n \geq s+1$.  

In many cases, we can do better than Proposition \ref{Claire prop}.  Let $d$ be any positive divisor of $n$, and let $H$ be the subgroup of order $n/d$ in $\mathbb{Z}_n$.  Note that if $A$ is a $[0,s]$-sum-free set in $\mathbb{Z}_d$, then $A+H$ is a $[0,s]$-sum-free set in $\mathbb{Z}_n$.   
Therefore:

\begin{prop} \label{mu [0,s] factor}
For all positive integers $n$ and $s$, we have 
$$\mu (\mathbb{Z}_n,[0,s]) \geq \mu (\mathbb{Z}_d,[0,s]) \cdot \frac{n}{d}.$$
\end{prop}

Combining Propositions \ref{Claire prop} and \ref{mu [0,s] factor} yields:

\begin{cor} \label{mu [0,s] cor}  
For all positive integers $n$ and $s$, we have
$$\mu (\mathbb{Z}_n,[0,s]) \geq \max \left \{ \left( \left \lfloor \frac{d-s-1}{s^2} \right \rfloor +1 \right) \cdot \frac{n}{d} \mid d \in D(n) \right \}.$$

\end{cor}
(Zajaczkowski\index{Zajaczkowski, C.} proved a special case of this for $s=2$ in \cite{Zaj:2013a}.)  

We can analyze the bound of Corollary \ref{mu [0,s] cor} for small values of $s$, as follows.  First, note that for $s=1$, Corollary \ref{mu [0,s] cor} becomes  
$$\mu (\mathbb{Z}_n,[0,1]) \geq \max \left \{ \left( d-1 \right) \cdot \frac{n}{d} \mid d \in D(n) \right \}=n-1,$$ and from Proposition \ref{mu [0,s] triv} we know that equality holds.

For $s=2$, we get 
$$\mu (\mathbb{Z}_n,[0,2]) \geq \max \left \{  \left \lfloor \frac{d+1}{4} \right \rfloor  \cdot \frac{n}{d} \mid d \in D(n) \right \}.$$
We separate two cases.  Assume first that $n$ has no divisors that are congruent to 3 mod 4.  In this case, for all $d \in D(n)$ we have 
$$\left \lfloor \frac{d+1}{4} \right \rfloor  \cdot \frac{n}{d} \leq \frac{d}{4} \cdot \frac{n}{d}= \frac{n}{4};$$ since  
the left-hand side is an integer, it can be at most $\lfloor n/4 \rfloor$.  On the other hand, using $d=n$ we see that
$$\max \left \{  \left \lfloor \frac{d+1}{4} \right \rfloor  \cdot \frac{n}{d} \mid d \in D(n) \right \} \geq \left \lfloor \frac{n+1}{4} \right \rfloor = 
\left \lfloor \frac{n}{4} \right \rfloor,$$ 
so Corollary \ref{mu [0,s] cor} gives $$\mu (\mathbb{Z}_n,[0,2]) \geq \left \lfloor \frac{n}{4} \right \rfloor.$$  

Suppose now that $n$ has divisors that are congruent to 3 mod 4.  Then for each such divisor $d$, we get 
$$\left \lfloor \frac{d+1}{4} \right \rfloor  \cdot \frac{n}{d} = \frac{d+1}{4} \cdot \frac{n}{d}= \left(1+ \frac{1}{d}   \right) \frac{n}{4};$$
this quantity is greatest when $d$ is the smallest such divisor, namely, it is the smallest prime divisor of $n$ that is congruent to 3 mod 4.  
In summary, for $s=2$ 
Corollary \ref{mu [0,s] cor} yields (see page \pageref{v2n4}):

\begin{cor} \label{cor mu [0,2]}
For all positive integers $n$, we have
$$\mu (\mathbb{Z}_n,[0,2]) \geq v_2(n,4)=\left\{
\begin{array}{ll}
\left(1+\frac{1}{p}\right) \frac{n}{4} & \mbox{if $n$ has prime divisors congruent to 3 mod 4,} \\ & \mbox{and $p$ is the smallest such divisor;}\\ \\
\left\lfloor \frac{n}{4} \right\rfloor & \mbox{otherwise.}\\
\end{array}\right.$$

\end{cor}

We believe that equality holds in Corollary \ref{cor mu [0,2]}:

\begin{conj} \label{conj mu [0,2]}
For all positive integers $n$, we have
$$\mu (\mathbb{Z}_n,[0,2]) = v_2(n,4).$$

\end{conj}
Zajaczkowski\index{Zajaczkowski, C.} verified Conjecture \ref{conj mu [0,2]} for all $n \leq 20$ (cf.~\cite{Zaj:2013a}). 

\begin{prob}
Prove Conjecture \ref{conj mu [0,2]}.
\end{prob}   

Moving on to $s=3$, we see that Corollary \ref{mu [0,s] cor} becomes:
\begin{cor} \label{cor mu [0,3]}
For all positive integers $n$, we have
$$\mu (\mathbb{Z}_n,[0,3]) \geq \max \left \{ \left \lfloor \frac{d+5}{9} \right \rfloor  \cdot \frac{n}{d} \mid d \in D(n) \right \}.$$

\end{cor}   

(Note that the bound here is not always equal to $v_3(n,9)$.)  In particular, when $n$ only has divisors that are congruent to 0, 1, 2, or 3 mod 9, then we only get $$\mu (\mathbb{Z}_n,[0,3]) \geq \lfloor n/9 \rfloor.$$  It is not clear if we can do better:

\begin{prob}
Improve, if possible, on the  lower bound of Corollary \ref{cor mu [0,3]}.
\end{prob} 

As an example, we prove that for $n=11$, we cannot improve on the bound of Corollary \ref{cor mu [0,3]}, which yields only $$\mu (\mathbb{Z}_{11},[0,3]) \geq 1.$$  We argue, as follows.  Suppose, indirectly, that $\mathbb{Z}_{11}$ contains a 2-subset $A=\{a,b\}$ that is $[0,3]$-sum-free.  This means that the sets
\begin{eqnarray*}
0A & = & \{0\} \\
1A & = & \{a,b\} \\
2A & = & \{2a, a+b, 2b\} \\
3A & = & \{3a,2a+b,a+2b,3b\}
\end{eqnarray*}
are pairwise disjoint.  Note also that the elements listed in each set are all distinct; for example, $2a+b=3b$ would yield $2a=2b$, which in $\mathbb{Z}_{11}$ can only happen if $a=b$.  Therefore, the ten elements listed are ten distinct elements of $\mathbb{Z}_{11}$; in other words, there is a unique element $c \in \mathbb{Z}_{11}$ so that
$$[0,3]A=\mathbb{Z}_{11} \setminus \{c\}.$$
Summing the elements on the two sides, we get
$$10a+10b=(0+1+\cdots+10)-c,$$ or
$$10a+10b=-c,$$ from which $c=a+b$.  But that is a contradiction, since $a+b \in A$ and thus $c \in [0,3]A$.  Therefore,
$$\mu (\mathbb{Z}_{11},[0,3]) = 1.$$

We know even less about the cases of $s \geq 4$:

\begin{prob}
Improve, if possible, on the  lower bound of Corollary \ref{mu [0,s] cor} for $s \geq 4$.

\end{prob}

We state the following challenging problems:

\begin{prob}
Find the exact value of $\mu (\mathbb{Z}_{n},[0,s])$ for all $n$ and a given $s \geq 4$.
\end{prob}

\begin{prob}
Find the exact value of $\mu (G,[0,2])$ for noncyclic groups $G$.
\end{prob}

\begin{prob}
Find the exact value of $\mu (G,[0,s])$ for noncyclic groups $G$ and $s \geq 3$.
\end{prob}

\subsection{Arbitrary number of terms} \label{6maxUarbitrary}

Here we ought to consider $$\mu (G,H) =\mathrm{max} \{ |A|  \mid  A \subseteq G; h_1, h_2 \in H; h_1 \neq h_2 \Rightarrow h_1A \cap h_2A = \emptyset \}$$   for the case when $H$ is the set of all nonnegative or all positive integers.  However, as we have already mentioned, we have $\mu (G,H)=0$ whenever $H$ contains two elements whose difference is a multiple of the exponent of the group.

\section{Unrestricted signed sumsets} \label{6maxUS}

\subsection{Fixed number of terms} \label{6maxUSfixed}

\subsection{Limited number of terms} \label{6maxUSlimited}

\subsection{Arbitrary number of terms} \label{6maxUSarbitrary}

Here we ought to consider $$\mu_{\pm} (G,H) =\mathrm{max} \{ |A|  \mid  A \subseteq G; h_1, h_2 \in H; h_1 \neq h_2 \Rightarrow (h_1)_{\pm}A \cap (h_2)_{\pm}A = \emptyset \}$$   for the case when $H$ is the set of all nonnegative or all positive integers.  However, as we have already mentioned, we have $\mu_{\pm}(G,H)=0$ whenever $H$ contains two elements whose difference is a multiple of the exponent of the group.

\section{Restricted sumsets} \label{6maxR}

Our goal in this section is to investigate the maximum possible size of a weak $H$-sum-free set, that is, the quantity $$\mu \hat{\;} (G,H) =\mathrm{max} \{ |A|  \mid  A \subseteq G; h_1, h_2 \in H; h_1 \neq h_2 \Rightarrow h_1 \hat{\;} A \cap h_2 \hat{\;} A = \emptyset \}.$$  

Clearly, for all $G$ of order at least 2, we have  $\mu \hat{\;} (G,H) \geq 1$: for any $a \in G \setminus \{0\}$, at least the one-element set $A=\{a\}$ will be weakly $H$-sum-free.

It is important to note that Proposition \ref{mufordirectsum} does not carry through for $\mu \hat{\;} (G,H)$: for example, $\{1,2,4\}$ is weakly $\{1,2\}$-sum-free in $\mathbb{Z}_{10}$ (the sum of any two distinct elements of $A$ is in $G \setminus A$), but $\{1,2,4\} \times \mathbb{Z}_{10}$ is not weakly $\{2,1\}$-sum-free in $\mathbb{Z}_{10}^2$ since, for example, $$(1,3)+(1,4)=(2,7).$$ 

Below we consider three special cases: when $H$ consists of two distinct positive integers (by Proposition \ref{sumsets trivial prop}, the case when one of the integers equals 0 is identical to Section \ref{5maxRfixed}), when $H$ consists of all nonnegative integers up to some value $s$, and when $H=\mathbb{N}_0$.

\subsection{Fixed number of terms} \label{6maxRfixed}

The analogue of a $(k,l)$-sum-free set for restricted addition is called a weak $(k,l)$-sum-free set; namely, subsets $A$ of $G$ satisfying the condition $(k\hat{\;} A) \cap (l\hat{\;} A) = \emptyset $ are called {\it weak $(k,l)$-sum-free sets}.  Here we investigate, for a given group $G$ and positive integers $k$ and $l$ with $k>l$, the quantity $$\mu\hat{\;} (G,\{k,l\}) = \mathrm{max} \{ |A|  \mid A \subseteq G, (k\hat{\;} A) \cap (l\hat{\;} A) = \emptyset\},$$ that is, the maximum size of a weak $(k,l)$-sum-free set in $G$.

Since for $k \geq n+1$ we trivially have $\mu\hat{\;} (G,\{k,l\}) =n$, we assume that $l<k \leq n$. 

Note that, if $A$ is $(k,l)$-sum-free, then it is also weakly $(k,l)$-sum-free, so by Theorem \ref{Bajnok mu bounds}, we have the following lower bound:

\begin{prop} \label{weak sumfree vs sumfree}
Suppose that $G$ is an abelian group of order $n$ and exponent $\kappa$.  Then, for all positive integers $k$ and $l$ with $k>l$ we have
 $$\mu\hat{\;}(G,\{k,l\}) \geq \mu(G, \{k,l\}) \geq v_{k-l}(\kappa,k+l) \cdot \frac{n}{\kappa}.$$

\end{prop}

There are a variety of cases when $\mu\hat{\;}(G,\{k,l\})$ is strictly more than $\mu (G,\{k,l\})$; we discuss some of these below.

First, we investigate intervals in the cyclic group, that is, arithmetic progressions in $\mathbb{Z}_n$ whose common difference is 1.  For a fixed element $a \in \mathbb{Z}_n$ and positive integer $m \leq n$, the set $$A=\{a,a+1,\dots,a+m-1\}$$ is called an {\em interval} of length $m-1$ (we say that $A$ has size $m$ and length $m-1$).  For example, $\{3,4,5,6\}$ and $\{5,6,0,1\}$ are both intervals of length 3 (size 4) in $\mathbb{Z}_7$.

The following result exhibits a formula for the maximum size $\gamma \hat{\;}(\mathbb{Z}_n,\{k,l\})$ of a weak $(k,l)$-sum-free interval in $\mathbb{Z}_n$.  We introduce some notations: 
$$\delta=\gcd(n,k-l),$$
$$J = k^2+l^2-(k+l),$$ 
$$M= \lfloor (n+J-2)/(k+l) \rfloor ,$$
$$K= kM-J/2+1, $$
$$L= lM-J/2+1.$$
Observe that $$K+L=(k+l)\lfloor (n+J-2)/(k+l) \rfloor-J+2 \leq n.$$

\begin{prop} \label{max size interval weak sumfree}
Let $n$, $k$, and $l$ be positive integers with $l < k \leq n$, and let $\gamma \hat{\;}(\mathbb{Z}_n,\{k,l\})$ be maximum size of a weak $(k,l)$-sum-free interval in $\mathbb{Z}_n$.  With the notations just introduced,
$$\gamma \hat{\;}(\mathbb{Z}_n,\{k,l\}) = \left\{
\begin{array}{cl}
M+1  & \mbox{if $L/\delta \leq \lfloor (n-K)/\delta \rfloor$}; \\ \\
M & \mbox{otherwise.}
\end{array}
\right.$$

\end{prop}
The proof of Proposition \ref{max size interval weak sumfree} is on page \pageref{proof of max size interval weak sumfree}.

We have the following more explicit consequence of Proposition \ref{max size interval weak sumfree}  (see the corresponding Corollary \ref{cor Zlower 1} regarding weak zero-$h$-sum-free sets):

\begin{cor} \label{kllower}
For positive integers $n$, $k$, and $l$ with $l < k \leq n$ we have
$$\mu\hat{\;}(\mathbb{Z}_n,\{k,l\}) \geq  \left\lfloor \frac{n+k^2+l^2-\mathrm{gcd}(n,k-l)-1}{k+l} \right\rfloor .$$
\end{cor}
We provide two proofs: one using Proposition \ref{max size interval weak sumfree} and another that is more direct; see page \pageref{proofofkllower}.

We can use Corollary \ref{kllower} (as well as Theorem \ref{Dias Da Silva and Hamidoune}) to find the value of $\mu\hat{\;}(\mathbb{Z}_p,\{k,l\})$ for all $k$, $l$, and prime values of $p$.  By Corollary \ref{kllower}, we get
$$\mu\hat{\;}(\mathbb{Z}_p,\{k,l\}) \geq \left \lfloor \frac{p+k^2+l^2-2}{k+l} \right \rfloor.$$  To see that equality holds, we employ Theorem \ref{Dias Da Silva and Hamidoune} to conclude that for a weak $(k,l)$-sum-free set $A \subseteq \mathbb{Z}_p$ we have
$$p \geq |k \hat{\;} A| + |l \hat{\;} A| \geq (k|A|-k^2+1)+(l|A|-l^2+1),$$ from which
$$|A| \leq \frac{p+k^2+l^2-2}{k+l}$$ follows.  Therefore:

\begin{thm}
For all primes $p$ and positive integers $k$ and $l$ with $l < k \leq p$ we have
$$\mu\hat{\;}(\mathbb{Z}_p,\{k,l\}) = \left \lfloor \frac{p+k^2+l^2-2}{k+l} \right \rfloor.$$
\end{thm}

The value of $\mu\hat{\;}(\mathbb{Z}_n,\{k,l\})$ for composite $n$ is harder to find, and is known for all $n$ only for $(k,l)=(2,1)$ (the case of {\em weak sum-free sets})---we present this result next. 

By combining Proposition \ref{weak sumfree vs sumfree} and Corollary \ref{kllower}, we get: 
$$\mu\hat{\;}(\mathbb{Z}_n,\{k,l\}) \geq \max \left\{\mu(\mathbb{Z}_n, \{k,l\}), \left\lfloor \frac{n+k^2+l^2-\mathrm{gcd}(n,k-l)-1}{k+l} \right\rfloor \right\}.$$

For $k=2$ and $l=1$, this translates to
$$\mu\hat{\;}(\mathbb{Z}_n,\{2,1\}) \geq \max \left\{\mu(\mathbb{Z}_n, \{2,1\}), \left\lfloor \frac{n}{3} \right\rfloor +1 \right\}.$$ 
The value of $\mu(\mathbb{Z}_n, \{2,1\})$ is given by Corollary \ref{k-l and n rel prime} as $v_1(n,3)$; on page \pageref{mu21} we see that 
$$v_1(n,3)=\left\{
\begin{array}{ll}
\left(1+\frac{1}{p}\right) \frac{n}{3} & \mbox{if $n$ has prime divisors congruent to 2 mod 3,} \\ & \mbox{and $p$ is the smallest such divisor,}\\ \\
\left\lfloor \frac{n}{3} \right\rfloor & \mbox{otherwise.}\\
\end{array}\right.$$
We should also note that, for any prime divisor $p$ of $n$ which is congruent to 2 mod 3, we have
$$\left(1+\frac{1}{p}\right) \frac{n}{3} \geq \left\lfloor \frac{n}{3} \right\rfloor+1.$$
Indeed, if $n \geq 3p$, then 
$$\left(1+\frac{1}{p}\right) \frac{n}{3} \geq \left(1+\frac{1}{n/3}\right) \frac{n}{3}=\frac{n}{3}+1 \geq \left\lfloor \frac{n}{3} \right\rfloor+1;$$
if $n=2p$, then 
$$\left(1+\frac{1}{p}\right) \frac{n}{3} = \left(1+\frac{1}{n/2}\right) \frac{n}{3} =\frac{n-1}{3} +1=\left\lfloor \frac{n}{3} \right\rfloor+1;$$ and
if $n=p$, then 
$$\left(1+\frac{1}{p}\right) \frac{n}{3} = \left(1+\frac{1}{n}\right) \frac{n}{3} =\frac{n-2}{3} +1=\left\lfloor \frac{n}{3} \right\rfloor+1.$$
This proves that
$$\mu\hat{\;}(\mathbb{Z}_n,\{2,1\}) \geq \left\{
\begin{array}{ll}
\left(1+\frac{1}{p}\right) \frac{n}{3} & \mbox{if $n$ has prime divisors congruent to 2 mod 3,} \\ & \mbox{and $p$ is the smallest such divisor,}\\ \\
\left\lfloor \frac{n}{3} \right\rfloor +1 & \mbox{otherwise.}\\
\end{array}\right.$$

It turns out that equality holds:

\begin{thm} [Zannier; cf.~\cite{Zan:2015a}] \label{M21}\index{Zannier, U.} 
For all positive integers  we have $$\mu\hat{\;}(\mathbb{Z}_n,\{2,1\}) = \left\{
\begin{array}{ll}
\left(1+\frac{1}{p}\right) \frac{n}{3} & \mbox{if $n$ has prime divisors congruent to 2 mod 3,} \\ & \mbox{and $p$ is the smallest such divisor,}\\ \\
\left\lfloor \frac{n}{3} \right\rfloor +1 & \mbox{otherwise.}\\
\end{array}\right.$$
\end{thm}
We present Zannier's short and elegant proof on page \pageref{proof of M21}.\index{Zannier, U.}

Let us turn to the case of $k=3$ and $l=1$.  Our considerations above now yield
$$\mu\hat{\;}(\mathbb{Z}_n,\{3,1\}) \geq \max \left\{\mu(\mathbb{Z}_n, \{3,1\}), \left\lfloor \frac{n+9-\mathrm{gcd}(n,2)}{4} \right\rfloor \right\}.$$
Recall that by Theorem \ref{Bajnok mu(3,1)} and page \pageref{v2n4},\index{Bajnok, B.}
$$ \mu (\mathbb{Z}_n,\{3,1\}) =v_{2}(n,4)= \left\{
\begin{array}{ll}
\left(1+\frac{1}{p}\right) \frac{n}{4} & \mbox{if $n$ has prime divisors congruent to 3 mod 4,} \\ & \mbox{and $p$ is the smallest such divisor,}\\ \\
\left\lfloor \frac{n}{4} \right\rfloor & \mbox{otherwise.}\\
\end{array}\right.$$

Note that, unless $n$ is divisible by 4, we have $$\left\lfloor \frac{n+9-\mathrm{gcd}(n,2)}{4} \right\rfloor=\left\lfloor \frac{n}{4} \right\rfloor+2,$$ and thus for these values of $n$ we have 
$$\mu\hat{\;}(\mathbb{Z}_n,\{3,1\}) \geq \lfloor n/4 \rfloor +2.$$  
We are able to make the same claim even when $n$ is divisible by 4 as long as it is not divisible by 8, since we may apply the stronger Proposition \ref{max size interval weak sumfree} to conclude exactly this when 
$$\frac{n/4-1}{2} \leq \left \lfloor \frac{n/4-1}{2} \right \rfloor$$ holds, or, equivalently, when
$n$ leaves a remainder of 4 mod 8.

This leaves us with the cases when $n$ is divisible by 8.  Ziegler\index{Ziegler, D.} in \cite{Zie:2009a} noticed that in some subcases one still has a weak $(3,1)$-sum-free set in $\mathbb{Z}_n$ of size $n/4 +2$ (even when $n$ has no prime divisors congruent to 3 mod 4).  For example, there are several such sets of size 4 in $\mathbb{Z}_8$, including $\{0,1,2,4\}$ and $\{0,3,4,6\}$.  Furthermore, he observed that when $n$ is divisible by 16, then
$$\left\{ \frac{n}{16}, \frac{n}{16}+1, \dots,\frac{3n}{16} \right\} \cup \left\{ \frac{9n}{16}, \frac{9n}{16}+1, \dots, \frac{11n}{16} \right\}$$ 
is a weak $(3,1)$-sum-free set in $\mathbb{Z}_n$ of size $n/4 +2$.  (See a generalization of this in Proposition \ref{Hallfors prop 1} below.)    
It is tempting to conjecture that $$\mu\hat{\;}(\mathbb{Z}_n,\{3,1\}) \geq \left\lfloor n/4 \right\rfloor +2$$ always holds; however, as the computational data of Hallfors\index{Hallfors, S.}  (see below) indicates, this is false for $n=40$.  (Note that $n=24$ has a prime divisor that is congruent to 3 mod 4.)

We pose the following intriguing problem:

\begin{prob}

Suppose that $n$ is congruent to 8 mod 16 and has no prime divisors congruent to 3 mod 4.  When is $\mu\hat{\;}(\mathbb{Z}_n,\{3,1\})$ equal to $ n/4 +2$ and when is it $ n/4 +1$?

\end{prob}
As we have discussed above, $\mu\hat{\;}(\mathbb{Z}_8,\{3,1\})=8/4+2$, but $\mu\hat{\;}(\mathbb{Z}_{40},\{3,1\})=40/4+1$.  The next value of $n$ in question is $n=104$.

Some more general questions:

\begin{prob}

Find $\mu\hat{\;}(\mathbb{Z}_n,\{3,1\})$ for all values of $n$.

\end{prob}

\begin{prob}

Find $\mu\hat{\;}(\mathbb{Z}_n,\{k,1\})$ for other values of $k$.

\end{prob}

\begin{prob}

Find $\mu\hat{\;}(\mathbb{Z}_n,\{k,l\})$ for other values of $k$ and $l$.

\end{prob}

\begin{prob}

Find $\mu\hat{\;}(G,\{k,l\})$ in noncyclic groups $G$, in particular, for $G=\mathbb{Z}_k^r$.

\end{prob}

We are not aware of any general exact results  besides the ones already discussed.  However, there are a variety of constructions yielding lower bounds for $\mu\hat{\;}(\mathbb{Z}_n,\{k,l\})$, as we present next.

One such result was discovered (although for $l=1$ only) by Hallfors in \cite{Hal:2012a}:\index{Hallfors, S.} 

\begin{prop} [Hallfors; cf.~\cite{Hal:2012a}] \label{Hallfors prop 1}\index{Hallfors, S.} 
Suppose that $n$, $k$, and $l$ are positive integers so that $l<k$ and $n$ is divisible by $(k^2-l^2)(k-1)$.  Let $H$ be the subgroup of order $k-1$ in $\mathbb{Z}_n$, and set $$A=\{a,a+1,\dots,a+c\}+H,$$ where
$$a=\frac{ln}{(k^2-l^2)(k-1)} \; \; \; \mbox{and}\;  \; \; c=\frac{n}{(k+l)(k-1)}.$$
Then $A$ is a weak sum-free set in $\mathbb{Z}_n$, and thus
$$\mu\hat{\;}(\mathbb{Z}_n,\{k,l\}) \geq \frac{n}{k+l}+k-1.$$
\end{prop}

For example, if $n$ is divisible by 3, then
$$\left\{ \frac{n}{3}, \frac{n}{3}+1, \dots,\frac{2n}{3} \right\} $$
is a weak $(2,1)$-sum-free set in $\mathbb{Z}_n$ of size $n/3+1$, and if $n$ is
divisible by 16, then we get the weak $(3,1)$-sum-free set
$$\left\{ \frac{n}{16}, \frac{n}{16}+1, \dots,\frac{3n}{16} \right\} \cup \left\{ \frac{9n}{16}, \frac{9n}{16}+1, \dots, \frac{11n}{16} \right\}$$
 of size $n/4+2$ mentioned above.  The proof of Proposition \ref{Hallfors prop 1} is on page \pageref{proof of Hallfors prop 1}.  

Another construction of Hallfors\index{Hallfors, S.} from \cite{Hal:2012a} (though only presented there for the special case of $d=4$) is as follows:

\begin{prop} [Hallfors; cf.~\cite{Hal:2012a}] \label{Hallfors prop 2}\index{Hallfors, S.} 
Suppose that $n$, $k$, and $l$ are positive integers so that $l<k$.  Let $d \in D(n)$ be even, and suppose that $n \geq d(d/2-1)$ and that the remainders of $k$ and $l$ when divided by $d$ differ by $d/2$.  Let $H$ be the subgroup of order $n/d$ in $\mathbb{Z}_n$, and set $$A=(1+H) \cup H_d,$$ where
$$H_d=\{0,d,2d, \dots, (d/2-2)d\} \subseteq H$$ (with $H_2=\emptyset$).
Then $A$ is a weak sum-free set in $\mathbb{Z}_n$, and thus
$$\mu\hat{\;}(\mathbb{Z}_n,\{k,l\}) \geq n/d+d/2-1.$$
\end{prop}

The proof of Proposition \ref{Hallfors prop 2} is on page \pageref{proof of Hallfors prop 2}.  

Observe that Proposition \ref{Hallfors prop 2} has the obvious special case for $d=2$: 
If $n$ is even, $k$ and $l$ have opposite parity, and $l < k \leq n$, then the set
$$A=\{1,3,5,\dots, n-1\}$$ 
is a weak $(k,l)$-sum-free set in $\mathbb{Z}_n$, and thus $$\mu\hat{\;}(\mathbb{Z}_n,\{k,l\}) \geq n/2.$$
This, however, already follows from Proposition \ref{weak sumfree vs sumfree} and Corollary \ref{v function bounds}, since in this case we have:
$$\mu\hat{\;}(\mathbb{Z}_n,\{k,l\}) \geq \mu (\mathbb{Z}_n,\{k,l\}) \geq v_{k-l}(n,k+l) = n/2.$$ 

We continue with a set of other (rather easy) constructions.  

\begin{prop}  \label{prop weak k1 sumfree easy}
Suppose that $n=k(k-1)$.  Let $H$ be the subgroup of order $k-1$ in $\mathbb{Z}_n$, and set
$$A=(H \setminus \{0\}) \cup (1+H).$$  Then $A$ is a weak $(k,1)$-sum-free set in $\mathbb{Z}_n$, so
  $$\mu\hat{\;}(\mathbb{Z}_n,\{k,1\}) \geq 2k-3.$$
\end{prop}
 
This claim can be verified quickly, as $$k \hat{\;} A \subseteq \{2,3,\dots,k-1\}+H,$$
so $k \hat{\;} A$ is disjoint from $A$.  We should point out that Proposition \ref{prop weak k1 sumfree easy} sometimes supersedes our previous results: for example, it yields $\mu\hat{\;}(\mathbb{Z}_{20},\{5,1\}) \geq 7$ (which is the exact value, as it turns out), while $\mu (\mathbb{Z}_{20},\{5,1\}) =4$ and $\gamma \hat{\;}(\mathbb{Z}_{20},\{5,1\}) =6$.

In a similar vein:

\begin{prop} \label{prop weak k1 sumfree easy 2}
Suppose that $n=k(k^2-1)$.  Let $H$ be the subgroup of order $k$ in $\mathbb{Z}_n$, define
$$h_0=\left\{
\begin{array}{cl}
0 & \mbox{if $k$ is odd,} \\
n/2 & \mbox{if $k$ is even,}
\end{array}
\right.$$
and set
$$A=(\{1,2,\dots,k-1\}+H) \cup (k+(H \setminus \{h_0\})).$$  Then $A$ is a weak $(k,1)$-sum-free set in $\mathbb{Z}_n$, so
  $$\mu\hat{\;}(\mathbb{Z}_n,\{k,1\}) \geq k^2-1=\frac{n}{k+1}+k-1.$$
\end{prop}

To verify our claim, observe that 
\begin{eqnarray*}
k \hat{\;} A & = & \{k+k(k-1)/2 \cdot (k^2-1)\} \cup (\{k+1,k+2,\dots,(k-1)k+k-1\}+H) \\
& = & \{k +h_0\} \cup (\{k+1,k+2,\dots,k^2-1\}+H) \\
& = & \mathbb{Z}_n \setminus A,
\end{eqnarray*}
establishing our claim.

Relying on a computer program, Hallfors\index{Hallfors, S.} in \cite{Hal:2012a} exhibited the values of $\mu\hat{\;}(\mathbb{Z}_n,\{k,l\})$ for $n \leq 40$ and $1 \leq l < k \leq 4$---we provide these data (except for $(k,l)=(2,1)$ for which we have Theorem \ref{M21}) in the table below.   

\newpage 

$$\begin{array}{||c||c|c|c|c|c||} \hline \hline
n & \mu\hat{\;}(\mathbb{Z}_n,\{3,1\}) & \mu\hat{\;}(\mathbb{Z}_n,\{4,1\}) & \mu\hat{\;}(\mathbb{Z}_n,\{3,2\}) & \mu\hat{\;}(\mathbb{Z}_n,\{4,2\}) & \mu\hat{\;}(\mathbb{Z}_n,\{4,3\})   \\ \hline \hline
5 & 3 & 4 & 3 & 3 & 4 \\ \hline
6 & 3 & 4 & 3 & 4 & 4 \\ \hline
7 & 3 & 4 & 3 & 4 & 4 \\ \hline
8 & 4 & 4 & 4 & 4 & 4 \\ \hline
9 & 4 & 4 & 4 & 4 & 4 \\ \hline
10 & 4 & 5 & 5 & 4 & 5 \\ \hline
11 & 4 & 5 & 4 & 4 & 4 \\ \hline
12 & 5 & 6 & 6 & 5 & 6 \\ \hline
13 & 5 & 5 & 4 & 5 & 5 \\ \hline
14 & 5 & 7 & 7 & 5 & 7 \\ \hline
15 & 5 & 6 & 5 & 5 & 5 \\ \hline
16 & 6 & 8 & 8 & 5 & 8 \\ \hline
17 & 6 & 6 & 5 & 5 & 5 \\ \hline
18 & 6 & 9 & 9 & 6 & 9 \\ \hline
19 & 6 & 6 & 6 & 6 &  6\\ \hline
20 & 7 & 10 & 10 & 6 & 10 \\ \hline
21 & 7 & 6 & 7 & 7 &  7\\ \hline
22 & 7 & 11 & 11 & 6 & 11 \\ \hline
23 & 7 & 7 & 6 & 6 &  6\\ \hline
24 & 8 & 12 & 12 & 8 & 12 \\ \hline
25 & 8 & 8 & 7 & 7 & 6 \\ \hline
26 & 8 & 13 & 13 & 7 & 13 \\ \hline
27 & 9 & 8 & 9 & 9 & 9 \\ \hline
28 & 9 & 14 &14  & 8 & 14 \\ \hline
29 & 9 & 8 & 8 & 7 & 7 \\ \hline
30 & 10 & 15 & 15 & 10 & 15 \\ \hline
31 & 9 & 9 & 8 & 8 & 7 \\ \hline
32 & 10 & 16 & 16 & 9 & 16 \\ \hline
33 & 11 & 9 & 11 & 11 & 11 \\ \hline
34 &  10& 17 & 17 & 8 & 17 \\ \hline
35 &10  & 10 & 10 & 8 & 8 \\ \hline
36 & 12 & 18 & 18 & 12 & 18 \\ \hline
37 & 11 & 10 & 9 & 9 & 8 \\ \hline
38 & 11 & 19 & 19 & 9 & 19 \\ \hline
39 & 13 & 10 & 13 & 13 & 13 \\ \hline
40 & 11 & 20 & 20 & 11 & 20 \\ \hline
\hline \end{array}$$

Additionally, Hallfors\index{Hallfors, S.} in \cite{Hal:2011a} also computed the values of $\mu\hat{\;}(\mathbb{Z}_n,\{k,l\})$ for $n \leq 30$, $l=1$, and $5 \leq k \leq 9$, which are as follows:   

\newpage 

$$\begin{array}{||c||c|c|c|c|c||} \hline \hline
n & \mu\hat{\;}(\mathbb{Z}_n,\{5,1\}) & \mu\hat{\;}(\mathbb{Z}_n,\{6,1\}) & \mu\hat{\;}(\mathbb{Z}_n,\{7,1\}) & \mu\hat{\;}(\mathbb{Z}_n,\{8,1\}) & \mu\hat{\;}(\mathbb{Z}_n,\{9,1\})   \\ \hline \hline
5 & 4 & 5 & 5 & 5 & 5 \\ \hline
6 & 4 & 5 & 6 & 6 & 6 \\ \hline
7 & 5 & 6 & 6 & 7 & 7 \\ \hline
8 & 5 & 6 & 7 & 7 & 8 \\ \hline
9 & 5 & 6 & 7 & 8 & 8 \\ \hline
10 & 5 & 6 & 7 & 8 & 8 \\ \hline
11 & 5 & 6 & 7 & 8 & 9 \\ \hline
12 & 6 & 6 & 7 & 8 & 9 \\ \hline
13 & 6 & 6 & 7 & 8 & 9 \\ \hline
14 & 6 & 7 & 7 & 8 & 9 \\ \hline
15 & 6 & 7 & 7 & 8 & 9 \\ \hline
16 & 6 & 8 & 8 & 8 & 9 \\ \hline
17 & 6 & 7 & 8 & 8 & 9 \\ \hline
18 & 6 & 9 & 8 & 9 & 9 \\ \hline
19 & 7 & 7 & 8 & 9 & 9 \\ \hline
20 & 7 & 10 & 8 & 10 & 10 \\ \hline
21 & 7 & 8 & 8 & 9 & 10 \\ \hline
22 & 7 & 11 & 8 &11  & 10 \\ \hline
23 & 7 & 8 & 8 & 9 & 10 \\ \hline
24 & 8 & 12 & 9 & 12 & 10 \\ \hline
25 & 8 & 8 & 9 & 10 & 10 \\ \hline
26 & 8 & 13 & 9 &13  & 10 \\ \hline
27 & 9 & 9 & 9 & 10 & 10 \\ \hline
28 & 8 & 14 & 9 & 14 & 10 \\ \hline
29 & 8 & 9 & 9 & 10 & 10 \\ \hline
30 & 10 & 15 & 9 & 15 & 11 \\ \hline
\hline \end{array}$$

By carefully analyzing these data, one may be able to generate constructions that are not included in this section:

\begin{prob}
Develop new general constructions for weak $(k,l)$-sum-free sets in $\mathbb{Z}_n$, and thereby exhibit new lower bounds for $\mu \hat{\;}(\mathbb{Z}_n,\{k,l\})$.  

\end{prob}

It may be interesting to investigate the maximum size of weak $(k,l)$-sum-free sets in groups from the opposite perspective: Rather than attempting to find $\mu \hat{\;}(G,\{k,l\})$ for a given group $G$, we look for all groups $G$ that contain a weak $(k,l)$-sum-free set of a given size:

\begin{prob} \label{allnform weak sumfree}
For each $k, l, m \in \mathbb{N}$ with $l<k$, find all groups $G$ for which $\mu \hat{\;}(G,\{k,l\}) \geq m$.
\end{prob}

We can formulate two sub-problems of Problem \ref{allnform weak sumfree}:

\begin{prob}

For each  $k, l, m \in \mathbb{N}$ with $l<k$, find the least integer $f(m,\{k,l\})$ for which $\mu \hat{\;}(\mathbb{Z}_{n},\{k,l\}) \geq m$ holds for $n =f(m,\{k,l\})$.
\end{prob}

\begin{prob}

For each  $k, l, m \in \mathbb{N}$ with $l<k$, find the  least integer  $g(m,\{k,l\})$  for which $\mu \hat{\;}(\mathbb{Z}_{n},\{k,l\}) \geq m$ holds for all $n \geq g(m,\{k,l\})$.
\end{prob}

For example, the data above indicate that $f(10,\{3,1\})=30$ and $g(10,\{3,1\})=32$.  

We have the following results:
\begin{prop}  \label{all n for weak sumfree trivial}
For all positive integers $k$, $l$, and $m$, with $l<k$, we have
$$g(m,\{k,l\}) \geq f(m,\{k,l\}) \geq m,$$  with equality if, and only if, $m \leq k-1$.
\end{prop}

The fact that both $f(m,\{k,l\})$ and $g(m,\{k,l\})$ must be at least $m$ is obvious, and it is also obvious that equality holds when $m \leq k-1$.  The fact that for $m \geq k$ we have $$g(m,\{k,l\}) \geq f(m,\{k,l\}) \geq m+1$$ follows directly from the Lemma on page \pageref{lemma any j}.

We also have the values of $f(m,\{k,1\})$ and $g(m,\{k,1\})$ for $m=k$ and $m=k+1$:

\begin{prop}  \label{all n for weak sumfree m=k, k+1}
For every integer $k \geq 2$ we have  
$$f(k,\{k,1\})=g(k,\{k,1\})=\left\{
\begin{array}{cl}
k+2 & \mbox{if $k \equiv 1$ mod 4,} \\ \\
k+1 & \mbox{otherwise;}
\end{array}\right.$$  and
$$f(k+1,\{k,1\})=g(k+1,\{k,1\})=2k+2.$$
 \end{prop}

The proof starts on page \pageref{proof of all n for weak sumfree m=k, k+1}.

\begin{prob}

Find the values of (or at least good bounds for) $f(m,\{k,l\})$ and $g(m,\{k,l\})$; in particular, evaluate $f(k+2,\{k,1\})$ and $g(k+2,\{k,1\})$.

\end{prob}

As we did with $(k,l)$-sum-free sets in Section \ref{6maxUfixed}, we also investigate weak $(k,l)$-sum-free sets where no elements are ``wasted'': we say that a weak $(k,l)$-sum-free set $A$ in $G$ is {\em complete} if $k \hat{\;} A$ and $l \hat{\;} A$ partition $G$; in other words, not only do we have $k \hat{\;} A \cap l \hat{\;} A = \emptyset$ (i.e. $A$ is weakly $(k,l)$-sum-free), but also $k \hat{\;} A \cup l \hat{\;} A=G$.

Looking through our results above, we find the following complete weak $(k,l)$-sum-free sets:

\begin{itemize}

\item The set $A$ of Proposition \ref{Hallfors prop 1} is complete when $l \leq k-2$ (see page \pageref{proof of Hallfors prop 1});

\item The set $A$ of Proposition \ref{Hallfors prop 2} is complete when $$d/2-1 < l < k < n/d$$ (see page \pageref{proof of Hallfors prop 2});

\item The set $A$ in Proposition \ref{prop weak k1 sumfree easy 2} is complete (see page \pageref{prop weak k1 sumfree easy 2});

\item The set $A$ in the proof of the second statement of Proposition \ref{all n for weak sumfree m=k, k+1} is complete when $n=2k+2$ (see page \pageref{proof of all n for weak sumfree m=k, k+1}).

\end{itemize}

\begin{prob}
Find other complete weak $(k,l)$-sum-free sets in abelian groups.

\end{prob}

\subsection{Limited number of terms} \label{6maxRlimited}

\subsection{Arbitrary number of terms} \label{6maxRarbitrary}

\section{Restricted signed sumsets} \label{6maxRS}

\subsection{Fixed number of terms} \label{6maxRSfixed}

\subsection{Limited number of terms} \label{6maxRSlimited}

\subsection{Arbitrary number of terms} \label{6maxRSarbitrary}

\part{Pudding} \label{ChapterPudding}

``The proof is in the pudding!''  In this chapter we provide proofs to those of our results in the book that have not been published (and that were not presented where the results are stated).  A bit of warning: as is often the case, puddings can be messy \ldots read at your own risk!

\newpage

\addcontentsline{toc}{section}{Proof of Proposition \ref{functionsacdirect}}

\section*{Proof of Proposition \ref{functionsacdirect}} \label{proofoffunctionsacdirect}

1.  Let us define $$a'(j,k):=\sum_{i=0}^k {k \choose i} {j \choose i} 2^i.$$  Clearly, $a'(j,0)=a'(0,k)=1$; below we prove that $a'(j,k)$ also satisfies the recursion.

We have 
\begin{eqnarray*}
a'(j-1,k-1) & = & \sum_{i=0}^{k-1} {k-1 \choose i} {j-1 \choose i} 2^i \\
& = & \sum_{i=0}^{k-2} {k-1 \choose i} {j-1 \choose i} 2^i +{j-1 \choose k-1} 2^{k-1}, 
\end{eqnarray*}
and
\begin{eqnarray*}
a'(j-1,k) & = & \sum_{i=0}^{k} {k \choose i} {j-1 \choose i} 2^i \\
& = & \sum_{i=0}^{k-1} {k \choose i} {j-1 \choose i} 2^i +{j-1 \choose k} 2^k \\ 
& = &  \sum_{i=0}^{k-1} {k-1 \choose i-1} {j-1 \choose i} 2^i+\sum_{i=0}^{k-2} {k-1 \choose i} {j-1 \choose i} 2^i +{j-1 \choose k-1} 2^{k-1} +{j-1 \choose k} 2^k. 
\end{eqnarray*}

Next, we add $a'(j-1,k)$ and $a'(j-1,k-1)$.  Note that 
$${j-1 \choose k-1} 2^{k-1} +{j-1 \choose k-1} 2^{k-1}+{j-1 \choose k} 2^k={j \choose k} 2^{k},$$
and $$\sum_{i=0}^{k-2} {k-1 \choose i} {j-1 \choose i} 2^i + \sum_{i=0}^{k-2} {k-1 \choose i} {j-1 \choose i} 2^i=\sum_{i=0}^{k-2} {k-1 \choose i} {j-1 \choose i} 2^{i+1},$$ and by replacing $i$ by $i-1$, this sum becomes 
$$\sum_{i=0}^{k-1} {k-1 \choose i-1} {j-1 \choose i-1} 2^i.$$

Therefore, 
\begin{eqnarray*}
a'(j-1,k)+a'(j-1,k-1) & = & \sum_{i=0}^{k-1} {k-1 \choose i-1} {j-1 \choose i} 2^i+\sum_{i=0}^{k-1} {k-1 \choose i-1} {j-1 \choose i-1} 2^i+{j \choose k} 2^{k} \\
& = & \sum_{i=0}^{k-1} {k-1 \choose i-1} {j \choose i} 2^i+{j \choose k} 2^{k} \\
& = & \sum_{i=0}^{k} {k-1 \choose i-1} {j \choose i} 2^i \\ 
& = & \sum_{i=0}^{k} {k \choose i} {j \choose i} 2^i - \sum_{i=0}^{k} {k-1 \choose i} {j \choose i} 2^i \\
& = & a'(j,k)-a'(j,k-1).
\end{eqnarray*}  

2.  The cases of $j=0$ or $k=0$ can be verified easily (see also the comments after the statement of Proposition \ref{functionsacdirect}).  For positive $j$ and $k$ we use Proposition \ref{functionsac} and part 1 above:
\begin{eqnarray*}
c(j,k) & = & a(j,k-1)+a(j-1,k-1) \\
& = & \sum_{i=0}^k {j \choose i} {k-1 \choose i} 2^i + \sum_{i=0}^k {j-1 \choose i} {k-1 \choose i} 2^i \\
& = & \sum_{i=0}^k \left[{j \choose i}+ {j-1 \choose i} \right]{k-1 \choose i-1}  2^i, \\
\end{eqnarray*}
and
\begin{eqnarray*}
c(j,k) & = & a(j,k)-a(j-1,k) \\
& = & \sum_{i=0}^k {k \choose i} {j \choose i} 2^i - \sum_{i=0}^k {k \choose i} {j-1 \choose i} 2^i \\
& = & \sum_{i=0}^k {k \choose i} \left[ {j \choose i} -  {j-1 \choose i} \right] 2^i \\
& = & \sum_{i=0}^k {k \choose i} {j-1 \choose i-1}  2^i. \\
\end{eqnarray*}
$\Box$

\addcontentsline{toc}{section}{Proof of Proposition \ref{rankvsgen}}

\section*{Proof of Proposition \ref{rankvsgen}} \label{proofofrankvsgen}

Suppose that $G$ has invariant factorization $$\mathbb{Z}_{n_1} \times \cdots \times \mathbb{Z}_{n_r}$$ (so $2 \leq n_1$ and $n_i$ is a divisor of $n_{i+1}$ for $i=1,\dots,r-1$).  Let $p$ be any prime divisor of $n_1$, and set $$H=\langle \{p \cdot g \mid g \in G\} \rangle.$$  (Actually, $p$ would not need to be prime, any divisor $d>1$ of $n_1$ would do, including $d=n_1$.)  It is easy to see that
$$H \cong \mathbb{Z}_{n_1/p} \times \cdots \times \mathbb{Z}_{n_r/p}$$ and $$G/H \cong \mathbb{Z}_p^r.$$  Since each nonzero element of $\mathbb{Z}_p^r$ has order $p$, it cannot be generated by fewer than $r$ elements.    

Now suppose that $A=\{a_i \mid i=1,\dots,m\} \subseteq G$ has size $m$ and that $\langle A \rangle=G$.  One can readily verify that $$\overline{A}=\{a_i+H \mid i=1,\dots,m\}$$ then generates $G/H$.  Thus we have
$$m=|A| \geq |\overline{A}| \geq r,$$ as claimed.   
$\Box$

\addcontentsline{toc}{section}{Proof of Proposition \ref{sumsetidentities}}

\section*{Proof of Proposition \ref{sumsetidentities}} \label{proofofsumsetidentities}

As usual, we let $A=\{a_1,\dots,a_m\}$.    

Suppose first that $x \in h(A \cup (-A))$; we then have
\begin{eqnarray*}
x & = & \lambda_1a_1+\cdots+\lambda_ma_m+\lambda_1'(-a_1)+\cdots+\lambda_m'(-a_m) \\
& = & (\lambda_1-\lambda_1')a_1+\cdots+(\lambda_m-\lambda_m')a_m
\end{eqnarray*} for some nonnegative integers $\lambda_1, \dots, \lambda_m$ and $\lambda_1', \dots, \lambda_m'$ for which $$\lambda_1+ \cdots +\lambda_m+\lambda_1'+ \cdots +\lambda_m'=h.$$

Divide the set $\{1,2,\dots,m\}$ into two subsets, $I_1$ and $I_2$, such that $i \in I_1$ when $\lambda_i \geq \lambda_i'$ and $i \in I_2$ when $\lambda_i < \lambda_i'$.  We can then write
\begin{eqnarray*}
\sum_{i=1}^m |\lambda_i - \lambda_i'|& =&  \sum_{i \in I_1} (\lambda_i - \lambda_i') + \sum_{i \in I_2} (\lambda_i' - \lambda_i) \\ \\
& = & \sum_{i=1}^m \lambda_i + \sum_{i=1}^m \lambda_i'- 2 \left( \sum_{i \in I_1} \lambda_i' + \sum_{i \in I_2}  \lambda_i \right) \\ \\
& = & h- 2 h_0
\end{eqnarray*} 
where 
\begin{eqnarray*}
2h_0& =&  2 \left( \sum_{i \in I_1} \lambda_i' + \sum_{i \in I_2}  \lambda_i \right) \\ \\
& = & \sum_{i \in I_2}  \lambda_i + \sum_{i \in I_1} \lambda_i' + \sum_{i \in I_2}  \lambda_i + \sum_{i \in I_1} \lambda_i' \\ \\
& \leq & \sum_{i \in I_2}  \lambda_i + \sum_{i \in I_1} \lambda_i' + \sum_{i \in I_2}  \lambda_i' + \sum_{i \in I_1} \lambda_i \\ \\
& = &h.
\end{eqnarray*}  
Thus, 
$x \in (h - 2h_0)_{\pm}A .$

Conversely, let $h_0$ be a nonnegative integer (not greater than $\lfloor h/2 \rfloor$), and suppose that $x \in (h - 2h_0)_{\pm}A;$ we then have 
integers $\lambda_1, \dots, \lambda_m$ with $$x=\lambda_1a_1+\cdots+\lambda_ma_m$$ and $$|\lambda_1|+\cdots+|\lambda_m| = h-2h_0.$$ 

Similarly to above, divide the set $\{1,2,\dots,m\}$ into two subsets, $I_1$ and $I_2$, such that $i \in I_1$ when $\lambda_i  \geq 0$ and $i \in I_2$ when $\lambda_i < 0$.  We then have
\begin{eqnarray*}
x& =&  \sum_{i \in I_1} \lambda_i a_i+ \sum_{i \in I_2} \lambda_i a_i \\ \\
& = & \sum_{i \in I_1} \lambda_i a_i+ \sum_{i \in I_2} (-\lambda_i )(-a_i) \\ \\
& = & \sum_{i \in I_1} \lambda_i a_i+ \sum_{i \in I_2} (-\lambda_i )(-a_i) + h_0 a_1 + h_0 (-a_1).
\end{eqnarray*}
This expression of $x$ involves only nonnegative coefficients and thus no cancellation of terms occurs (even when terms can be combined, e.g. when $a_1=-a_1$), so for the nonnegative coefficients we have
$$\sum_{i \in I_1} \lambda_i + \sum_{i \in I_2} (-\lambda_i )+ 2h_0  = |\lambda_1|+\cdots+|\lambda_m| + 2h_0  = h.$$
Thus $x \in h(A \cup (-A))$. 
$\Box$

It may be worthwhile to mention that we can easily prove that 

$$h \hat{\;}(A \cup (-A)) \subseteq h \hat{_{\pm}} A \cup (h-2) \hat{_{\pm}} A \cup (h-4) \hat{_{\pm}} A \cup \cdots.$$

The proof is essentially the same as the corresponding proof above: note that with the additional assumption that each coefficient $\lambda_i$ and $\lambda_i'$ is 0 or 1 implies that the coefficients $\lambda_i-\lambda_i'$ in the expression of $x$ above are all in $\{-1,0,1\}$, and thus $x \in (h-2h_0)\hat{_{\pm}} A$.

\addcontentsline{toc}{section}{Proof of Proposition \ref{hA and h pm A with sym}} 

\section*{Proof of Proposition \ref{hA and h pm A with sym}} \label{proof of hA and h pm A with sym}

The claim is obvious for $h=0$ and $h=1$, since $$0_{\pm}A=\{0\}=0A,$$ and $$1_{\pm}A=A \cup -A=A=1A.$$  For $h \geq 2$, we show that $(h-2)_{\pm}A \subseteq h_{\pm}A$, from which our proposition follows by Proposition \ref{sumsetidentities}.

To show that $(h-2)_{\pm}A \subseteq h_{\pm}A$,  let $A=\{a_1,\dots,a_m\}$ and $$g=\lambda_1a_1+ \cdots + \lambda_ma_m$$ with coefficients $\lambda_1, \dots, \lambda_m \in \mathbb{Z}$ for which $$|\lambda_1|+ \cdots  + |\lambda_m|=h-2,$$ so that $g \in (h-2)_{\pm}A$.  Without loss of generality, assume that $\lambda_1 \geq 0$.  Since $A=-A$, there is an index $i \in \{1,\dots,m\}$ so that $-a_1=a_i$.  We consider three cases.

If $i=1$, then $2a_1=0$, and thus $$g=(\lambda_1+2)a_1+ \lambda_2a_2+ \cdots + \lambda_ma_m,$$ so $g \in h_{\pm}A$.  

If $i \neq 1$, w.l.o.g. $i=2$, and $\lambda_2 \geq 0$, then  $a_1+a_2=0$, and thus 
$$g=(\lambda_1+1)a_1+ (\lambda_2+1)a_2+ \lambda_3a_3+ \cdots + \lambda_ma_m,$$ so $g \in h_{\pm}A$.

If $i \neq 1$, w.l.o.g. $i=2$, and $\lambda_2 < 0$, then  
$$g=(\lambda_1+1)a_1+ (-\lambda_2+1)a_2+ \lambda_3a_3+ \cdots + \lambda_ma_m,$$ so $g \in h_{\pm}A$ again. $\Box$

\addcontentsline{toc}{section}{Proof of Proposition \ref{v for prime power}}

\section*{Proof of Proposition \ref{v for prime power}} \label{proof of v for prime power}

The claim is true for $h=1$, so we assume that $h \geq 2$.

Since $g$ is not divisible by $p$, we have $\gcd(p^k,g)=1$, so $g$ is relatively prime with all divisors of $p^r$.

Suppose first that $p \equiv 1$ mod $h$.  We then have
\begin{eqnarray*}
v_g(p^r,h) & = & \max \left\{ \left( \left \lfloor \frac{p^k-2}{h} \right \rfloor +1  \right) \cdot \frac{p^r}{p^k}  \mid 1 \leq k \leq r \right\} \\
& = & \max \left\{ \left(  \frac{p^k-1-h}{h}  +1  \right) \cdot p^{r-k}  \mid 1 \leq k \leq r \right\} \\
& = & \max \left\{   \frac{p^k-1}{h}   \cdot p^{r-k}  \mid 1 \leq k \leq r \right\} \\
& = & \max \left\{   \frac{p^r-p^{r-k}}{h}    \mid 1 \leq k \leq r \right\} \\
& = & \frac{p^r-1}{h}.
\end{eqnarray*}

Next, let's assume that $p \equiv 0$ mod $h$; this can only happen if $p=h$.  We then have
\begin{eqnarray*}
v_g(p^r,p) & = & \max \left\{ \left( \left \lfloor \frac{p^k-2}{p} \right \rfloor +1  \right) \cdot \frac{p^r}{p^k}  \mid 1 \leq k \leq r \right\} \\
& = & \max \left\{ \left(  \frac{p^k-p}{p}  +1  \right) \cdot p^{r-k}  \mid 1 \leq k \leq r \right\} \\
& = & p^{r-1} \\
& = & \left(  \left \lfloor \frac{p-2}{p} \right \rfloor  +1  \right) \cdot p^{r-1}.
\end{eqnarray*}

For our main case, assume that $p \equiv i$ mod $h$ for some $2 \leq i \leq h-1$.  Since $p \in D(p^r)$, we must have
$$v_g(p^r,p) \geq \left(  \left \lfloor \frac{p-2}{h} \right \rfloor  +1  \right) \cdot p^{r-1};$$
so, to prove our claim, we must show that
$$\left(  \left \lfloor \frac{p-2}{h} \right \rfloor  +1  \right) \cdot p^{r-1} \geq \left(  \left \lfloor \frac{p^k-2}{h} \right \rfloor  +1  \right) \cdot p^{r-k} $$ holds for every $2 \leq k \leq r$.  (We will, in fact, show that strict inequality holds.)  We consider two subcases: when $p<h$ and when $p \geq h$.

If $p < h$, then for every $2 \leq k \leq r$ we have
\begin{eqnarray*}
\left(  \left \lfloor \frac{p-2}{h} \right \rfloor  +1  \right) \cdot p^{r-1} & = & p^{r-1} \\
& = & \frac{p^r +p^{r-2} \cdot p}{p+1} \\
& > & \frac{p^r +p^{r-k} \cdot (p-1)}{p+1} \\ 
& = & \left( \frac{p^k -2}{p+1} + 1 \right) \cdot p^{r-k} \\ 
& \geq & \left(  \left \lfloor \frac{p^k-2}{h} \right \rfloor  +1  \right) \cdot p^{r-k},
\end{eqnarray*}
as claimed.

Similarly, if $p \geq h$, then for every $2 \leq k \leq r$ we have
\begin{eqnarray*}
\left(  \left \lfloor \frac{p-2}{h} \right \rfloor  +1  \right) \cdot p^{r-1} & \geq &  \left(   \frac{p-(h-1)}{h}   +1  \right) \cdot p^{r-1}  \\
& = & \frac{p+1}{h} \cdot p^{r-1} \\
& = & \frac{p^k+p^{k-1}}{h} \cdot p^{r-k} \\
& > & \frac{p^k+h-2}{h} \cdot p^{r-k} \\
& \geq & \left(  \left \lfloor \frac{p^k-2}{h} \right \rfloor  +1  \right) \cdot p^{r-k},
\end{eqnarray*}
as claimed.  $\Box$

\addcontentsline{toc}{section}{Proof of Proposition \ref{nusmall}}

\section*{Proof of Proposition \ref{nusmall}} \label{proofofnusmall}

The first claim is quite straight-forward:
\begin{eqnarray*}
v_1(n,1) & = & \max \left\{ \left( \left \lfloor \frac{d-1-1}{1} \right \rfloor +1  \right) \cdot \frac{n}{d}  \mid d  \in D(n) \right\} \\
& = & \max \left\{ n- \frac{n}{d}  \mid d \in D(n) \right\} \\
& = & n- \frac{n}{n} \\
& = & n-1.
\end{eqnarray*}

For $v_1(n,2)$ we have
\begin{eqnarray*}
v_1(n,2) & = & \max \left\{ \left( \left \lfloor \frac{d-1-1}{2} \right \rfloor +1  \right) \cdot \frac{n}{d}  \mid d \in D(n) \right\} \\
& = & \max \left\{ \left \lfloor \frac{d}{2} \right \rfloor \cdot \frac{n}{d}  \mid d \in D(n) \right\}. \\
\end{eqnarray*}

Now if $n$ is odd, then all its divisors are odd, and we have
\begin{eqnarray*}
v_1(n,2) & = & \max \left\{ \frac{d-1}{2} \cdot \frac{n}{d} \ d \in D(n) \right\} \\
& = & \max \left\{ \frac{n}{2} -  \frac{n}{2d} \; d \in D(n) \right\} \\
& = & \frac{n}{2}- \frac{n}{2n} \\
& = &  \frac{n-1}{2},
\end{eqnarray*} 
as claimed.

If $n$ is even, then it has both even and odd divisors.  For each even divisor $d$ of $n$, we get
$$\left \lfloor \frac{d}{2} \right \rfloor \cdot \frac{n}{d}=\frac{n}{2};$$ while for odd divisors $d$, we have $$\left \lfloor \frac{d}{2} \right \rfloor \cdot \frac{n}{d} < \frac{n}{2}.$$  Therefore, if $n$ is even, we have
$$v_1(n,2)=\frac{n}{2}.$$

The proof of the third claim is similar.
$\Box$

\addcontentsline{toc}{section}{Proof of Theorem \ref{nu}}

\section*{Proof of Theorem \ref{nu}} \label{proofofnu}

Suppose that $d$ is a positive divisor of $n$, and let $i$ be the remainder of $d$ when divided by $h$.  
We define the function
$$f(d)= \left(  \left \lfloor \frac{d-1-\mathrm{gcd}(d,g)}{h} \right \rfloor +1 \right) \cdot \frac{n}{d}.$$

We first prove the following.

{\bf Claim 1:}  We have
$$ f(d) =
\left\{
\begin{array}{ll}
\frac{n}{h} \cdot \left(1+ \frac{h-i}{d} \right) & \mbox{if } \; \mathrm{gcd}(d,g) < i; \\ \\
\frac{n}{h}\cdot \left( 1-\frac{h}{d} \right) & \mbox{if } \; h|d \; \mbox{and} \; g=h; \\ \\
\frac{n}{h} \cdot \left( 1-\frac{i}{d} \right) & \mbox{otherwise. }
\end{array}\right.$$

{\em Proof of Claim 1.}  We start with
$$ \left \lfloor \frac{d-1-\mathrm{gcd}(d,g)}{h} \right \rfloor=\frac{d-i}{h}+\left \lfloor \frac{i-1-\mathrm{gcd}(d,g)}{h} \right \rfloor.$$
We investigate the maximum and minimum values of the quantity $\left \lfloor \frac{i-1-\mathrm{gcd}(d,g)}{h} \right \rfloor$.

For the maximum, we have 
$$\left \lfloor \frac{i-1-\mathrm{gcd}(d,g)}{h} \right \rfloor \leq \left \lfloor \frac{(h-1)-1-1}{h} \right \rfloor \leq 0,$$ with equality if, and only if, $i-1-\mathrm{gcd}(d,g) \geq 0$; that is, $\mathrm{gcd}(d,g) < i$.

For the minimum, we get
$$\left \lfloor \frac{i-1-\mathrm{gcd}(d,g)}{h} \right \rfloor \geq \left \lfloor \frac{0-1-g}{h} \right \rfloor \geq 
\left \lfloor \frac{0-1-h}{h} \right \rfloor =-2,$$
with equality if, and only if, $i=0$, $\mathrm{gcd}(d,g)=g$, and $g=h$; that is, $h|d$ and $g=h$.  

The proof of Claim 1 now follows easily.  $\Box$

{\bf Claim 2:}  Using the notations as above, assume that $\mathrm{gcd}(d,g) \geq i$.  Then
$$ f(d) \leq
\left\{
\begin{array}{ll}
\frac{n}{h}  & \mbox{if } \;  g \neq h; \\ \\
\frac{n-1}{h}& \mbox{if } \;  g = h.
\end{array}\right.$$

{\em Proof of Claim 2.}  By Claim 1, we have 
$$ f(d) =
\left\{
\begin{array}{ll}
\frac{n}{h}\cdot \left( 1-\frac{h}{d} \right) & \mbox{if } \; h|d \; \mbox{and} \; g=h; \\ \\
\frac{n}{h} \cdot \left( 1-\frac{i}{d} \right) & \mbox{otherwise. }
\end{array}\right.$$
Therefore, $$f(d) \leq \frac{n}{h}.$$    
Furthermore, unless $i=0$ and $g \neq h$, we have 
$$f(d) \leq \frac{n}{h}\cdot \left( 1-\frac{1}{d} \right) \leq \frac{n}{h}\cdot \left( 1-\frac{1}{n} \right) = \frac{n-1}{h}.$$
$\Box$

{\bf Claim 3:}  For all $g$, $h$, and $n$ we have 
$$ v_g(n,h) \geq
\left\{
\begin{array}{ll}
\left \lfloor \frac{n}{h} \right \rfloor & \mbox{if } \; g \neq h; \\ \\
\left \lfloor \frac{n-1}{h} \right \rfloor & \mbox{if } \; g=h.
\end{array}\right.$$

{\em Proof of Claim 3.}  We first note that 
\begin{eqnarray*}
v_g(n,h) &=& \max \left\{ \left( \left \lfloor \frac{d-1-\mathrm{gcd} (d, g)}{h} \right \rfloor +1  \right) \cdot \frac{n}{d}  \mid d \in D(n) \right\} \\ \\
& \geq & \left( \left \lfloor \frac{n-1-\mathrm{gcd} (n, g)}{h} \right \rfloor +1  \right) \cdot \frac{n}{n} \\ \\
& = & \left \lfloor \frac{n-1-\mathrm{gcd} (n, g)}{h} \right \rfloor +1 \\ \\
& \geq & \left \lfloor \frac{n-1-g}{h} \right \rfloor +1 \\ \\
& = & \left \lfloor \frac{n+(h-g-1)}{h} \right \rfloor.
\end{eqnarray*}
The claim now follows, since $h-g-1 \geq 0$, unless $g=h$ in which case $h-g-1=-1$.   $\Box$

We are now ready for the proof of Theorem \ref{nu}.

{\em Proof of Theorem \ref{nu}.}  If $I = \emptyset$, then by Claims 2 and 3 we have 
$$ v_g(n,h) =
\left\{
\begin{array}{ll}
\left \lfloor \frac{n}{h} \right \rfloor & \mbox{if } \; g \neq h; \\ \\
\left \lfloor \frac{n-1}{h} \right \rfloor & \mbox{if } \; g=h
\end{array}\right.$$

Suppose now that $I \neq \emptyset$.  Here we call a positive divisor $d$ of $n$ good, if its remainder mod $h$ is larger than $\gcd(d,g)$.  Since $I \neq \emptyset$, $n$ has at least one good divisor; we let $D_I$ be the collection of good divisors of $n$.  By Claim 1, we have 
$$v_g(n,h)  = \max \{ f(d) \; | \; d \in D_I \}.$$  Suppose now that $d$ and $d'$ are two elements of $D_I$ that both leave a remainder of $i$ mod $h$.  If $d'>d$, then by Claim 1, we have
$$f(d)=  \frac{n}{h} \cdot \left(1+ \frac{h-i}{d} \right) > \frac{n}{h} \cdot \left(1+ \frac{h-i}{d'} \right) =f(d'),$$ 
from which our result follows.  $\Box$

\addcontentsline{toc}{section}{Proof of Proposition \ref{nu pm small}}

\section*{Proof of Proposition \ref{nu pm small}} \label{proofofnu pm small}

For a divisor $d$ of $n$, let
$$g_d(n,h)=\left( 2 \cdot \left \lfloor \frac{d-2}{2h} \right \rfloor +1 \right) \cdot \frac{n}{d}.$$ 

Regarding $h=1$, we see that 
$$g_d(n,1)=\left( 2 \cdot \left \lfloor \frac{d-2}{2} \right \rfloor +1 \right) \cdot \frac{n}{d} = \left( d-\delta \right) \cdot \frac{n}{d}= \left( 1 - \frac{\delta}{d} \right) \cdot n,$$
where $\delta=1$ when $d$ is even and $\delta=2$ when $d$ is odd.  This quantity is maximized with $d=n$, from which the result for $v_{\pm}(n,1)$ follows.

Now for $h=2$ and even $n$, we have $2 \in D(n)$, and thus
$$\frac{n}{2}=g_2(n,2) \leq v_{\pm}(n,2) \leq v_1(n,2)=\frac{n}{2},$$ so $v_{\pm}(n,2)=n/2$.
When $n$ is odd, then each $d \in D(n)$ is odd too, so
$$v_{\pm}(n,2) = \max \left \{ g_d(n,2)\right \} \leq  \max \left \{ \left( 2 \cdot  \frac{d-3}{4}  +1 \right) \cdot \frac{n}{d} \right \} = 
\max \left \{ \left( 1-\frac{1}{d} \right) \cdot \frac{n}{2} \right \}.$$
Therefore, when $n \equiv 3$ mod 4, then
$$\frac{n-1}{2}=g_n(n,2) \leq v_{\pm}(n,2) \leq \left( 1-\frac{1}{n} \right) \cdot \frac{n}{2}=\frac{n-1}{2},$$ so equality holds throughout.  When $n \equiv 1$ mod 4, then we get
$$\frac{n-3}{2}=g_n(n,2) \leq v_{\pm}(n,2),$$ and there is no $d \in D(n)$ for which $g_d(n,2)$ is larger than $(n-3)/2$, since if $d \neq n$, then $d \leq n/3$, and thus
$$\left( 1-\frac{1}{d} \right) \cdot \frac{n}{2} \leq \left( 1-\frac{1}{n/3} \right) \cdot \frac{n}{2} = \frac{n-3}{2}.$$

For $h=3$ and even $n$, we get $v_{\pm}(n,3)=n/2$ as we did for $h=2$.  
When $n$ is odd but divisible by 3, then $3 \in D(n)$, so
$$\frac{n}{3}=g_3(n,3) \leq v_{\pm}(n,3) = \max \left \{ g_d(n,3)\right \} \leq  \max \left \{ \left( 2 \cdot  \frac{d-3}{6}  +1 \right) \cdot \frac{n}{d} \right \} = 
 \frac{n}{3},$$ and thus $v_{\pm}(n,3)=n/3$.  When $n$ is odd and is not divisible by 3, then each $d \in D(n)$ leaves a remainder of 1 or 5 mod 6, so 
$$v_{\pm}(n,3) = \max \left \{ g_d(n,3)\right \} \leq  \max \left \{ \left( 2 \cdot  \frac{d-5}{6}  +1 \right) \cdot \frac{n}{d} \right \} = 
\max \left \{ \left( 1-\frac{2}{d} \right) \cdot \frac{n}{3} \right \}.$$
Therefore, when $n \equiv 5$ mod 6, then
$$\frac{n-2}{3}=g_n(n,3) \leq v_{\pm}(n,3) \leq \left( 1-\frac{2}{n} \right) \cdot \frac{n}{3}=\frac{n-2}{3},$$ so equality holds throughout.  When $n \equiv 1$ mod 6, then we get
$$\frac{n-4}{2}=g_n(n,3) \leq v_{\pm}(n,3),$$ and there is no $d \in D(n)$ for which $g_d(n,3)$ is larger than $(n-4)/3$, since if $d \neq n$, then $d \leq n/5$, and thus
$$\left( 1-\frac{2}{d} \right) \cdot \frac{n}{3} \leq \left( 1-\frac{2}{n/5} \right) \cdot \frac{n}{3} = \frac{n-10}{3}.$$
This completes our proof.  $\Box$

\addcontentsline{toc}{section}{Proof of Proposition \ref{nu pm h=4}}

\section*{Proof of Proposition \ref{nu pm h=4}} \label{proofofnu pm h=4}

For a divisor $d$ of $n$, let
$$g_d(n,4)=\left( 2 \cdot \left \lfloor \frac{d-2}{8} \right \rfloor +1 \right) \cdot \frac{n}{d}.$$

As before, for even $n$, we have $2 \in D(n)$, and thus
$$\frac{n}{2}=g_2(n,4) \leq v_{\pm}(n,4) \leq v_1(n,4)=\frac{n}{2},$$ so $v_{\pm}(n,4)=n/2$.

Suppose that $n$ is odd but has one or more divisors $d$ congruent to 3 mod 8, then for all such $d$ we have
$$ g_{d}(n,4) = \left( 2 \cdot  \frac{d-3}{8}  +1 \right) \cdot \frac{n}{d} = \left ( 1+\frac{1}{d} \right) \frac{n}{4},$$ which is maximal when $d$ is minimal.  For any other $d \in D(n)$ we have $$ g_{d}(n,4) \leq  \left( 2 \cdot  \frac{d-5}{8}  +1 \right) \cdot \frac{n}{d} < \frac{n}{4},$$ so our claim holds in this case.

Suppose now that $n$ is odd and has no divisors congruent to 3 mod 8; in particular, $3 \not \in D(n)$ and thus $n/3 \not \in D(n)$, so the second largest divisor of $n$ is at most $n/5$.
We have
$$g_d(n,4) \leq \left( 2 \cdot \frac{d-5}{8}  +1 \right) \cdot \frac{n}{d} = \left( 1- \frac{1}{d}   \right) \cdot \frac{n}{4}.$$
We consider the cases of $n \equiv 5$ mod 8, $n \equiv 7$ mod 8, and $n \equiv 1$ mod 8 separately.   Let $d_0 \in D(n)$ be such that $v_{\pm}(n,4)=g_{d_0}(n,4)$.

If $n \equiv 5$ mod 8, then 
$$\frac{n-1}{4} = 2 \cdot \left \lfloor \frac{n-2}{8} \right \rfloor +1 = g_n(n,4) \leq g_{d_0}(n,4) \leq \left( 1- \frac{1}{d_0}   \right) \cdot \frac{n}{4} \leq \left( 1- \frac{1}{n}   \right) \cdot \frac{n}{4} =\frac{n-1}{4},$$ 
so equality holds throughout.

If $n \equiv 7$ mod 8, then $d_0=n$ and thus
$$v_{\pm}(n,4)=g_n(n,4)= 2 \cdot \left \lfloor \frac{n-2}{8} \right \rfloor +1 =\frac{n-3}{4},$$ since if $d_0 \neq n$, then $d_0 \leq n/5$, so
$$g_{d_0}(n,4) \leq \left( 1- \frac{1}{d_0}   \right) \cdot \frac{n}{4} \leq \left( 1- \frac{1}{n/5}   \right) \cdot \frac{n}{4} =\frac{n-5}{4}.$$ 

Suppose, finally, that $n \equiv 1$ mod 8.  
If $d_0=n$, we get
$$v_{\pm}(n,4)=g_n(n,4)= 2 \cdot \left \lfloor \frac{n-2}{8} \right \rfloor +1 =\frac{n-5}{4}.$$  If $d_0 \neq n$, then $d_0 \leq n/5$.  If $d_0=n/5$, we get
$$v_{\pm}(n,4)=g_{n/5}(n,4) \leq \left( 1 - \frac{1}{n/5} \right) \cdot \frac{n}{4} =\frac{n-5}{4}.$$
If $d_0 \leq n/7$, then  
$$v_{\pm}(n,4) = g_{d_0}(n,4) \leq \left( 1- \frac{1}{d_0}   \right) \cdot \frac{n}{4} \leq \left( 1- \frac{1}{n/7}   \right) \cdot \frac{n}{4} =\frac{n-7}{4}.$$  
Therefore, 
$$v_{\pm}(n,4)=\frac{n-5}{4}=2\cdot \left \lfloor \frac{n-2}{8} \right \rfloor +1.$$
That completes our proof.  $\Box$

\addcontentsline{toc}{section}{Proof of Theorem \ref{h-critical of n}}

\section*{Proof of Theorem \ref{h-critical of n}} \label{proof of h-critical of n}

We need to prove that, for $m=v_1(n,h)$, we have $u(n,m,h)<n$ but $u(n,m+1,h) \geq n$.  Let $d_0 \in D(n)$ be such that $$v_1(n,h)=\max \left\{ \left( \left \lfloor \frac{d-2}{h} \right \rfloor +1  \right) \cdot \frac{n}{d}   \mid d \in D(n) \right\}=\left( \left \lfloor \frac{d_0-2}{h} \right \rfloor +1  \right) \cdot \frac{n}{d_0}.$$    

To establish the first inequality, simply note that $u(n,m,h) \leq f_{n/d_0}(m,h)$ where
$$f_{n/d_0}(m,h)=\left(h \cdot \left ( \left \lfloor \frac{d_0-2}{h} \right \rfloor +1 \right )    -h+1  \right) \cdot \frac{n}{d_0}=
\left(h \cdot \left \lfloor \frac{d_0-2}{h} \right \rfloor +1  \right) \cdot \frac{n}{d_0} \leq (d_0-1) \cdot \frac{n}{d_0} <n. $$

For the second inequality, we must prove that, for any $d \in D(n)$, we have $f_d(m+1,h) \geq n$; that is, 
$$h \cdot \left \lceil \frac{ 
\left( \left \lfloor \frac{d_0-2}{h} \right \rfloor +1 \right )  \cdot \frac{n}{d_0} +1}{d} 
\right \rceil  -h+1 \geq \frac{n}{d}.$$
 But $n/d \in D(n)$, so by the choice of $d_0$, we have
$$\left( \left \lfloor \frac{d_0-2}{h} \right \rfloor +1 \right )  \cdot \frac{n}{d_0} \geq \left( \left \lfloor \frac{n/d-2}{h} \right \rfloor +1 \right )  \cdot \frac{n}{n/d},$$ and thus
\begin{eqnarray*}
h \cdot \left \lceil \frac{ 
\left( \left \lfloor \frac{d_0-2}{h} \right \rfloor +1 \right )  \cdot \frac{n}{d_0} +1}{d} 
\right \rceil  -h+1
& \geq & h \cdot \left \lceil 
\left( \left \lfloor \frac{n/d-2}{h} \right \rfloor +1 \right )  +\frac{1}{d} 
\right \rceil  -h+1 \\
& = & h \cdot  
\left( \left \lfloor \frac{n/d-2}{h} \right \rfloor +2 \right )  
  -h+1 \\
& \geq & h \cdot  
\left(  \frac{n/d-2-(h-1)}{h}  +2 \right )  
  -h+1 \\
& = & \frac{n}{d}.
\end{eqnarray*}
Our proof is complete.  $\Box$

\addcontentsline{toc}{section}{Proof of Proposition \ref{uhatforh=2}}

\section*{Proof of Proposition \ref{uhatforh=2}} \label{proofofuhatforh=2}

From the formulae on page \pageref{fhatsmall}  we already have the following values:
$$f\hat{_1}(n,m,2)= \min\{n,2m-3\};$$ and
$$f\hat{_2}(n,m,2)= \left\{\begin{array}{ll}
\min\{n,2m-4\}  & \mbox{if $m$ is even,}\\ \\
\min\{n,2m-3\} & \mbox{if $m$ is odd.} 
\end{array}\right.$$
For $d \geq 3$, we have
$$f\hat{_d}(n,m,2)=
\left\{\begin{array}{ll}
\min\{n, f_d, 2m-3\}  & \mbox{if $d \not| m-1$,}\\ \\
\min\{n, 2m-2\} & \mbox{otherwise.}  
\end{array}\right.$$

To prove Proposition \ref{uhatforh=2}, we suppose first that $m$ is odd.  Then 
$$f\hat{_1}(n,m,2)= f\hat{_2}(n,m,2)=\min\{n,2m-3\}.$$
We consider three subcases: when $u(n,m,2)=n$, when $2m-3 \leq u(n,m,2) <n$, and when $u(n,m,2)<2m-3$.  (Note that we always have $u(n,m,2) \leq n$.)

If $u(n,m,2)=n$, then for every $d \in D(n)$ we have $f_d \geq n$, and thus for $d \geq 3$ we have
$$f\hat{_d}(n,m,2)=
\left\{\begin{array}{ll}
\min\{n, 2m-3\}  & \mbox{if $d \not| m-1$,}\\ \\
\min\{n, 2m-2\} & \mbox{otherwise.}  
\end{array}\right.$$  Therefore, $$u\hat{\;}(n,m,2)=\min\{ f\hat{_d}(n,m,h)  \mid d \in D(n) \}=\min\{n, 2m-3\}=\min\{u(n,m,2), 2m-3\},$$ as claimed.

Similarly, if $2m-3 \leq u(n,m,2) <n$, then for every $d \in D(n)$ we have $f_d \geq 2m-3$, and thus for $d \geq 3$ we have
$$f\hat{_d}(n,m,2)=
\left\{\begin{array}{ll}
\min\{n, 2m-3\}  & \mbox{if $d \not| m-1$,}\\ \\
\min\{n, 2m-2\} & \mbox{otherwise.}  
\end{array}\right.$$  Therefore, $$u\hat{\;}(n,m,2)=\min\{n, 2m-3\}=2m-3=\min\{u(n,m,2), 2m-3\},$$ as claimed.

Finally, if $u(n,m,2)<2m-3$, then let us choose a $d_0 \in D(n)$ for which $u(n,m,2)=f_{d_0}$.  Note that $f_{d_0} <2m-3$ but, as a quick computation shows, $f_1=2m-1$ and $f_2=2m$, so $d_0 \geq 3$.  Furthermore, whenever $d|m-1$, we get $f_d=2m+d-2>2m-3$, so $d_0 \not|m-1$, and thus $f\hat{_{d_0}}=f_{d_0}$.  

We will verify that $u\hat{\;}(n,m,2)=f_{d_0}$.  Clearly, $u\hat{\;}(n,m,2) \leq f\hat{_{d_0}}=f_{d_0}$, so we only need to show that there is no $d \in D(n)$ for which $f\hat{_d} < f_{d_0}$.  Since $f_{d_0} <2m-3$, this could only happen if $f\hat{_d} =f_d$, but that is not possible as $f_d \geq f_{d_0}$.
Thus, again we get, $$u\hat{\;}(n,m,2)=f_{d_0}=u(n,m,2)=\min\{u(n,m,2), 2m-3\},$$ as claimed.

The case when $m$ is even can be analyzed similarly. $\Box$

\addcontentsline{toc}{section}{Proof of Proposition \ref{uhatforh=3}} 

\section*{Proof of Proposition \ref{uhatforh=3}} \label{proofofuhatforh=3}

We can carry out a similar analysis for the case of $h=3$, where we may assume that $4 \leq m \leq n$.  If $d \geq 4$ and $k \geq 3$, then 
$$f \hat{_d}(n,m,3)=\min\{n,f_d,3m-8\}.$$  For the remaining cases we have the following.
$$\begin{array}{|c|c|c|c||c|} \hline
d & k & r & \delta_d & f \hat{_d}(n,m,3) \\ \hline \hline
1 & 1 & 1 & 0 & \min\{n,3m-8\} \\ \hline
2 & 1 & 1 & 0 & \min\{n,3m-8\} \\ \hline
2 & 2 & 1 & 0 & \min\{n,3m-8\} \\ \hline
3 & 1 & 3 & 0 & \min\{n,3m-8\} \\ \hline
3 & 2 & 3 & 0 & \min\{n,3m-8\} \\ \hline
3 & 3 & 3 & 2 & \min\{n,3m-10\} \\ \hline
\geq 4 & 1 & 3 & d-5 & \min\{n,3m-3-d\} \\ \hline
\geq 4 & 2 & 3 & -2 & \min\{n,3m-6\} \\ \hline
\end{array}$$
The case of $d \geq 4$ and $k=1$ is interesting: it implies that $m-1$ is divisible by $d$; thus, to minimize $3m-3-d$, we select $d=\gcd(n,m-1)$.  Note also that when $3|m$ then $\gcd(n,m-1) \neq 6$.  Therefore, 
$$u\hat{\;} (n,m,3)=\left\{
\begin{array}{ll}
\min\{u(n,m,3), 3m-3-\gcd(n,m-1)\} & \mbox{if} \; \gcd(n,m-1) \geq 8; \\ \\
\min\{u(n,m,3), 3m-10\} & \mbox{if} \; \gcd(n,m-1) = 7, \; \mbox{or} \\
& \gcd(n,m-1) \leq 5, \; 3|n, \; \mbox{and} \; 3|m; \\ \\
\min\{u(n,m,3), 3m-9\} & \mbox{if} \; \gcd(n,m-1) = 6; \\ \\
\min\{u(n,m,3), 3m-8\} & \mbox{otherwise.} 
\end{array}
\right.$$

\addcontentsline{toc}{section}{Proof of Proposition \ref{uhat lessthan u}}

\section*{Proof of Proposition \ref{uhat lessthan u}} \label{proofofuhat lessthan u}

Choose $d_0 \in D(n)$ so that $f_{d_0}=u(n,m,h)$.  We will show that then $f\hat{_{d_0}} \leq f_{d_0}$; this will then imply our claim since
$$u\hat{\;}(n,m,h) =\min\{ f\hat{_d}  \mid d \in D(n) \} \leq f\hat{_{d_0}} \leq f_{d_0} = u(n,m,h).$$

Suppose, indirectly, that $f\hat{_{d_0}} > f_{d_0}$;  we then have $h> \min\{d_0-1,k\}$ and 
$$f\hat{_{d_0}} =\min\{n, hm-h^2+1-\delta_{d_0}\},$$ otherwise we would have
$$f\hat{_{d_0}} =\min\{n, f_{d_0} , hm-h^2+1\} \leq f_{d_0},$$ a contradiction.  

Therefore, $$f\hat{_{d_0}} \leq hm-h^2+1-\delta_{d_0}.$$

Assume first that $h < d_0$.  Then $k<h=r <d_0$, so $\delta_{d_0} = (d_0-r)(r-k)-(d_0-1)$, and thus
\begin{eqnarray*}
h m-hk + d_0 & = &  \left( h \cdot  \frac{m-k+d_0}{d_0} -h+1 \right) \cdot d_0 \\
& = & \left( h \left \lceil \frac{m}{d_0} \right \rceil -h+1 \right) \cdot d_0 \\
 & = & f_{d_0} \\
& < &  f\hat{_{d_0}} \\
& \leq & hm-h^2+1-\delta_{d_0} \\
& = & hm-h^2- (d_0-r)(r-k)+d_0 \\
& = & hm-h^2- (d_0-h)(h-k)+d_0 \\
& = & hm-hk+d_0-d_0(h-k) \\
& < & hm-hk+d_0,
\end{eqnarray*}
which is a contradiction.

Assume now that $h \geq d_0$.  We then have $h \geq k$ since $k \leq d_0$.  Furthermore, we see from the definition of $\delta$ that $\delta \geq 1-(d-1)$ in all cases; therefore, $-\delta_{d_0} \leq d_0-2$.  We now see that 
\begin{eqnarray*}
h m-hk + d_0 & = &  \left( h \cdot  \frac{m-k+d_0}{d_0} -h+1 \right) \cdot d_0 \\
& = & \left( h \left \lceil \frac{m}{d_0} \right \rceil -h+1 \right) \cdot d_0 \\
 & = & f_{d_0} \\
& < &  f\hat{_{d_0}} \\
& \leq & hm-h^2+1-\delta_{d_0} \\
& \leq& hm-h^2+d_0 -1\\
& \leq & hm-hk+d_0 -1\\
\end{eqnarray*}
which is a contradiction. $\Box$

\addcontentsline{toc}{section}{Proof of Proposition \ref{uhat(n,m,h)extreme}}

\section*{Proof of Proposition \ref{uhat(n,m,h)extreme}} \label{proofofuhat(n,m,h)extreme}

Let $d \in D(n)$.  If $h \leq \min\{d-1,k\}$, then 
$$f\hat{_d}=\min\{n,f_d,hm-h^2+1\}.$$ Clearly, $n \geq m$, and we can easily see that $hm-h^2+1 \geq m$ holds as well, with equality if, and only if, $h=1$ or $h=m-1$.  Furthermore, from the argument for Proposition \ref{u(n,m,h)extreme} we see that we also have $f_d \geq m$ with equality if, and only if, $h=1$ or $m=d$.  Therefore, when $h \leq \min\{d-1,k\}$,  we have $f\hat{_d} \geq m$ with equality if, and only if, $h=1$, $h=m-1$, or $m=d$. 

Suppose now that  $h > \min\{d-1,k\}$, in which case
$$f\hat{_d}=\min\{n,hm-h^2+1-\delta_d\}.$$

If $\delta_d=0$, then $f\hat{_d}\geq m$ with equality if, and only if, $h=1$ or $h=m-1$.  This leaves three cases, according to the definition of $\delta_d$.

Assume first that $r=k=d$.  Note that, if $d=1$, then $\delta_d=0$, so we assume that $d \geq 2$; this then also implies that $h \geq 2$.  Furthermore, $r=d$ implies that $h \geq d$ and $k=d$ implies that $m-h \geq d$.  Therefore, we have
\begin{eqnarray*}
hm-h^2+1-\delta_d-m &=& (h-1)(m-h-1)-(d-1) \\
&\geq &(h-1)(m-h-1)-(h-1) \\
&=&(h-1)(m-h-2) \\
&\geq & (h-1)(d-2) \\
& \geq &0,
\end{eqnarray*} 
with equality if, and only if, $h=d$, $m-h=d$, and $d=2$; that is, $h=2$, $m=4$, and $d=2$. 

Next, assume that $r<k$; note that this then implies that $r \leq d-1$.  We now have
$$hm-h^2+1-\delta_d-m =(h-1)(m-h-1)-(k-r)r+(d-1)=(h-1)(m-h)-(k-r)r-(h-d).$$
Here $$m-h=(cd+k)-(qd+r)=(c-q)d+(k-r) \geq k-r,$$ so 
\begin{eqnarray*}
hm-h^2+1-\delta_d-m & \geq & (h-1)(k-r)-(k-r)r-(h-d) \\
&=& (h-r-1)(k-r)-(h-d) \\
&\geq & (h-r-1)-(h-d) \\
&=& d-1-r \\
&\geq &0.
\end{eqnarray*}
We see that equality holds if, and only if, $m-h=k-r$, $k-r=1$, and $r=d-1$; that is, $h=qd+d-1$, and $m=h+1=qd+d$.

In our last case, we have $k<r<d$.  In this case, we must have $d \geq 3$ and  
$$1 \leq m-h = (cd+k)-(qd+r)=d(c-q)+(k-r)<d(c-q),$$ so $c-q \geq 1$ and thus $$m-h=d(c-q-1)+(d+k-r) \geq d+k-r.$$
Thus,
\begin{eqnarray*}
hm-h^2+1-\delta_d-m &=& (h-1)(m-h-1)-(d-r)(r-k)+(d-1)\\
& \geq & (h-1)(d+k-r)-(d-r)(r-k)+(d-1) \\
&=& (d-r)(h-r+k-1)+k(h-1)+(d-1).
\end{eqnarray*}
Here each term is nonnegative, but $d -1 >0$.  

Therefore, we have proved that $f\hat{_d} \geq m$; since $d \in D(n)$ was arbitrary, this implies that $u\hat{\;} (n,m,h) \geq m$.

If $h=1$ or $h=m-1$, then with $d=n$ we get $h \leq \min\{d-1,k\}$, and thus $u\hat{\;} (n,m,h) \leq f\hat{_d} =m$.

If $m$ is a divisor of $n$, then we may choose $d=m$, with which again $h \leq \min\{d-1,k\}$, and thus $u\hat{\;} (n,m,h) \leq f\hat{_d} =m$.

If $h=2$, $m=4$, and $n$ is even, we may pick $d=2$, with which $h > \min\{d-1,k\}$, and so $u\hat{\;} (n,m,h) \leq f\hat{_d} =m$ holds again.

If none of these occur, then $f\hat{_d}> m$ for all $d \in D(n)$ and thus $u\hat{\;} (n,m,h) > m$.  $\Box$

\addcontentsline{toc}{section}{Proof of Proposition \ref{nu m=h cyclic}}

\section*{Proof of Proposition \ref{nu m=h cyclic}}  \label{proof of nu m=h cyclic}

We first show that we can write every odd integer between $-2^m$ and $2^m$ in the form
$\Sigma_{i=0}^{m-1} \pm 2^i$.  There are exactly $2^m$ such odd integers, so we just need to verify that no two of the $2^m$ signed sums yield the same integer. Suppose we have
$$\Sigma_{i=0}^{m-1} \lambda_i \cdot 2^i = \Sigma_{i=0}^{m-1} \lambda_i' \cdot 2^i$$ with coefficients $\lambda_i$ and $\lambda_i'$ all from the set $\{-1,1\}$, and assume that $j$ is the largest index for which $\lambda_i \neq \lambda_i'$; without loss of generality, let $\lambda_j=1$ and $\lambda_j'=-1$.  But then 
$$\Sigma_{i=0}^{m-1} \lambda_i \cdot 2^i \geq 1+\Sigma_{i=j+1}^{m-1} \lambda_i \cdot 2^i,$$ and
$$\Sigma_{i=0}^{m-1} \lambda_i' \cdot 2^i \leq -1+\Sigma_{i=j+1}^{m-1} \lambda_i \cdot 2^i,$$ which is a contradiction.

Now if $n$ is odd, then our $2^m$ odd integers yield exactly $\min \{n,2^m\}$ elements in $\mathbb{Z}_n$; since by Proposition \ref{nuupperRSfixed} we cannot hope for more, we have $$\nu\hat{_\pm}(\mathbb{Z}_n,m,h) = \mathrm{min} \left\{ n, 2^m \right\}.$$  If $n$ is even, our integers yield only $\min \{n/2,2^m\}$ elements in $\mathbb{Z}_n$, which is less than the upper bound in Proposition \ref{nuupperRSfixed} when $n<2^{m+1}$.  But note that, whenever $n$ is even, no $m$-subset $A=\{a_1,\dots,a_m\}$ of $\mathbb{Z}_n$ has a restricted $m$-fold signed sumset of size more than $n/2$: indeed, the signed sums in $$m \hat{_\pm}A=\{\pm a_1\pm \cdots \pm a_m\}$$
all have the same parity.  $\Box$

\addcontentsline{toc}{section}{Proof of Theorem \ref{no perfect bases s=2,3}}

\section*{Proof of Theorem \ref{no perfect bases s=2,3}} \label{proof of no perfect bases s=2,3}

We consider the cases of $s=2$ and $s=3$ separately.

Suppose that $A$ is a perfect 2-basis of size $m$ in $G$.  Then we must have $$n={m+2 \choose 2}.$$

Clearly, $0 \in A-A$; we show that $A-A$ has size $m(m-1)+1$: in other words, if $a_1$ and $a_2$ are distinct elements of $A$ and $a_3$ and $a_4$ are distinct elements of $A$ so that $$a_1-a_2=a_3-a_4,$$ then $a_1=a_3$ and $a_2=a_4$.  Indeed: our hypothesis implies that $$a_1+a_4=a_2+a_3,$$ so if $A$ is a perfect 2-basis, then $a_1 \in \{a_2,a_3\}$ and $a_4 \in \{a_2,a_3\}$, from which our claim follows since $a_1 \neq a_2$ and $a_3 \neq a_4$.

Furthermore, $A-A$ and $A$ are disjoint, since if we were to have $a_1, a_2, a_3 \in A$ for which $$a_1-a_2=a_3,$$ then $$a_1=a_2+a_3,$$ contradicting the assumption that $A$ is a perfect 2-basis.

Therefore,
$$|(A-A) \cup A|=m(m-1)+1+m \leq n = {m+2 \choose 2},$$ from which $m \leq 3$ follows.  We will rule out $m=2$ and $m=3$, as follows.

If $m=2$, then $n={m+2 \choose 2}=6$, so $G=\mathbb{Z}_6$.  Suppose that $A=\{a,b\} \subseteq \mathbb{Z}_6$.  Then 
$$[0,2]A=\{0,a,b,2a,2b,a+b\}=\mathbb{Z}_6,$$
so $a$, $b$, and $a+b$ must all be odd, which is impossible.

Similarly, if $m=3$, then $n={m+2 \choose 2}=10$, so $G=\mathbb{Z}_{10}$.  Suppose that $A=\{a,b\} \subseteq \mathbb{Z}_{10}$.  Then 
$$[0,2]A=\{0,a,b,c,2a,2b,2c,a+b,a+c,b+c\}=\mathbb{Z}_{10},$$
so exactly one of $a$, $b$, $c$, $a+b$, $a+c$, or $b+c$ is even, which is impossible.

Now we turn to $s=3$.  Suppose that $A$ is a perfect 3-basis of size $m$ in $G$.  Then we must have $$n={m+3 \choose 3}.$$

(Our argument here is similar to the proof of Theorem 3 in \cite{BraRuiTru:2012a}.)  Note that the set $2A-A$ is the disjoint union of the sets
$$B=\{a_1+a_2-a_3 \mid a_1, a_2, a_3 \in A; a_1 \neq a_2 \neq a_3 \neq a_1 \},$$
$$C=\{2a_1-a_2 \mid a_1, a_2 \in A, a_1 \neq a_2 \},$$ and $A$ since, for example, $$a_1+a_2-a_3=2a_1'-a_2'$$ would imply $$a_1+a_2+a_2'=2a_1'+a_3,$$ from which $a_3 \in \{a_1,a_2,a_2'\}$ so $a_3=a_2'$, and thus $a_1=a_2$, which is a contradiction.  Furthermore, $B$ and $C$ each have `maximum' cardinality; that is, if $$a_1+a_2-a_3=a_1'+a_2'-a_3',$$ then $\{a_1,a_2\}=\{a_1',a_2'\}$ and $a_3=a_3'$; similarly for $C$.  Therefore, 
$$|2A-A|=|B|+|C|+|A|=m(m-1)(m-2)/2+m(m-1)+m.$$  

We can similarly see that $2A-A$ must also be disjoint from $A-A$, and, as above, 
 we have $$|A-A|=m(m-1)+1.$$  

Therefore, we have $$m(m-1)(m-2)/2+m(m-1)+m+ m(m-1)+1 \leq n = {m+3 \choose 3},$$ from which $m \leq 3$.  

We can easily rule out $m=2$ by hand: If $A=\{a,b\}$ were to be a perfect 3-basis in $G$, then $n=10$ and thus $G=\mathbb{Z}_{10}$; with $$[0,3]A=\{0,a,b,2a,2b,a+b,3a,3b,2a+b,a+2b\}$$ we see that (i) if both $a$ and $b$ are odd, then 6 elements of $[0,3]A$ are odd; (ii) if exactly one of $a$ or $b$ is odd, then 4 elements of $[0,3]A$ are odd; and (iii) if both $a$ and $b$ are even, then no element of $[0,3]A$ is odd.  The case of $m=3$ leads to $G=\mathbb{Z}_{20}$ or $G=\mathbb{Z}_2 \times \mathbb{Z}_{10}$; we ruled out these by the computer program \cite{Ili:2017a}.  $\Box$

\addcontentsline{toc}{section}{Proof of Proposition \ref{nform=2}}

\section*{Proof of Proposition \ref{nform=2}} \label{proofofpropnform=2}

We will prove that for each $2 \leq n \leq \lfloor \tfrac{f^2+6f+5}{4} \rfloor$, the set $A=\{1, a \}$ with $a=\lfloor \tfrac{f+3}{2} \rfloor$ has folding number at most $f$ in $\mathbb{Z}_n$.

What we need to verify is that for every $g \in \mathbb{Z}_n$, we have $g \in \cup_{h=0}^f hA$.  We may assume that $g$ is a nonnegative integer which is less than $n$.

By the Division Theorem, there exist (unique) integers $q$ and $r$ with $0 \leq r \leq a-1$ for which $g=q \cdot a+r \cdot 1$.  We will prove that $q+r \leq f$, from which our claim follows.

Case 1: $f$ is even.  In this case, we have $a=\frac{f+2}{2}$ and $n \leq \frac{f^2+6f+4}{4}$.  We separate several subcases.

Subcase 1.1: $q \leq a-1$.  Then $q+r \leq 2a-2 \leq f$, as claimed.

Subcase 1.2: $q=a$ and $r \leq a-2$.  Again, $q+r \leq 2a-2 \leq f$, as claimed.

Subcase 1.3: $q=a$ and $r = a-1$.  Then we have $$g=a^2+a-1=\frac{f^2+6f+4}{4}\geq n,$$ which cannot happen, since we assumed that $g$ is less than $n$.

Subcase 1.4: $q \geq a+1$.  Then we have $g \geq a^2+a$, which is a contradiction as in Subcase 1.3.

Case 2: $f$ is odd.  In this case, we have $a=\frac{f+3}{2}$ and $n \leq \frac{f^2+6f+5}{4}$.  We again separate several subcases.

Subcase 2.1: $q \leq a-2$.  Then $q+r \leq 2a-3 \leq f$, as claimed.

Subcase 2.2: $q=a-1$ and $r \leq a-2$.  Again, $q+r \leq 2a-3 \leq f$, as claimed.

Subcase 2.3: $q=a-1$ and $r = a-1$.  Then we have $$g=a(a-1)+a-1=\frac{f^2+6f+5}{4} \geq n,$$ which cannot happen, since we assumed that $g$ is less than $n$.

Subcase 2.4: $q \geq a$.  Then we have $g \geq a^2$, which is a contradiction as in Subcase 2.3.
$\Box$

\addcontentsline{toc}{section}{Proof of Proposition \ref{perfex}}

\section*{Proof of Proposition \ref{perfex}} \label{proofofperfex} 

As we noted below the statement of Proposition \ref{perfex}, here we only need to prove that the set $\{1,2s+1\}$ is a perfect $s$-spanning set in $\mathbb{Z}_n$ for $n=2s^2+2s+1$.

Let $n=2s^2+2s+1$.  We will prove that for every integer $g$ between $-s$ and $n-s-1$, inclusive, one can find integers $\lambda_1$ and $\lambda_2$ so that $|\lambda_1|+|\lambda_2| \leq s$ and $$g \equiv \lambda_1 \cdot 1+\lambda_2 \cdot (2s+1) \; \mbox{mod} \; n.$$

Let $k$ be the largest odd integer for which $g>k(s+1)$; $k$ is then uniquely defined by the inequalities
$$k (s+1) +1 \leq g \leq (k+2)(s+1).$$
Note that $g \geq -s$ implies that $k \geq -1$, and $g \leq 2s^2+s$ implies that $k \leq 2s-1$, so we have $$-1 \leq k \leq 2s-1.$$
We consider two cases.

Case 1: Suppose that $$k(s+1)+1 \leq g \leq k(s+1)+(2s-k)=(k+2)s.$$
(Note that $2s-k \geq 1$.)  Let $$\lambda_1=g-(k+1)s-\frac{k+1}{2},$$ and $$\lambda_2=\frac{k+1}{2}.$$  We can readily verify that $$g = \lambda_1 \cdot 1+\lambda_2 \cdot (2s+1);$$ we need to prove that $|\lambda_1|+|\lambda_2| \leq s$ holds as well.

Note that $k \geq -1$, so $$|\lambda_2|=\frac{k+1}{2}.$$  Since $\lambda_1$ is a linear function of $g$, $|\lambda_1|$ achieves its maximum value over the interval for $g$ at one (or both) of the endpoints; therefore, we need to evaluate
$$\left| k(s+1)+1-(k+1)s-\frac{k+1}{2} \right|$$ and 
$$\left| (k+2)s-(k+1)s-\frac{k+1}{2} \right|;$$ since $s \geq \frac{k+1}{2}$, both expressions equal $s-\frac{k+1}{2}$.  
Therefore, $|\lambda_1|+|\lambda_2| \leq s$, as claimed.

Case 2:  Suppose now that $$(k+2)s+1 \leq g \leq (k+2)(s+1)=(k+2)s+(k+2).$$
(Note that $k+2 \geq 1$.)  Let $$\lambda_1=g-(k+2)s-\frac{k+3}{2},$$ and $$\lambda_2=\frac{k+1}{2}-s.$$  We can verify that $$\lambda_1 \cdot 1+\lambda_2 \cdot (2s+1)=g-(2s^2+2s+1),$$ and thus $\lambda_1 \cdot 1+\lambda_2 \cdot (2s+1) \equiv g$ mod $n$.  Next, we prove that $|\lambda_1|+|\lambda_2| \leq s$ holds as well.

First, since $s \geq \frac{k+1}{2}$, we have $$|\lambda_2| = s-\frac{k+1}{2}.$$  Since $\lambda_1$ is a linear function of $g$, $|\lambda_1|$ achieves its maximum value over the interval for $g$ at one (or both) of the endpoints; therefore, we need to evaluate
$$\left| (k+2)s+1-(k+2)s-\frac{k+3}{2} \right|$$ and 
$$\left| (k+2)(s+1)-(k+2)s-\frac{k+3}{2} \right|;$$ since $k \geq -1$, both expressions equal $\frac{k+1}{2}$.  
Therefore, we again have $|\lambda_1|+|\lambda_2| \leq s$, as claimed.  $\Box$

\addcontentsline{toc}{section}{Proof of Proposition \ref{p}}

\section*{Proof of Proposition \ref{p}} \label{proofofp}

We verify that for every $n \leq 2s^2+2s+1$, we have
$$[0,s]_{\pm}\{s,s+1\}=\{\lambda_1s+\lambda_2(s+1) \in \mathbb{Z}_n \mbox{    } | \mbox{    }  \lambda_1,\lambda_2 \in \mathbb{Z}, |\lambda_1|+|\lambda_2| \leq s \}=\mathbb{Z}_n.$$ 

Considering the elements of $[0,s]_{\pm}\{s,s+1\}$ in $\mathbb{Z}$ (rather than $\mathbb{Z}_n$), we see that the elements of 
$$\Sigma=\{\lambda_1s+\lambda_2(s+1) \in \mathbb{Z} \mbox{    } | \mbox{    }  \lambda_1,\lambda_2 \in \mathbb{Z}, |\lambda_1|+|\lambda_2| \leq s \}$$ 
 lie in the interval $[-(s^2+s),(s^2+s)]$.  Since the index set $$I_{\pm}(2,[0,s])=\{(\lambda_1,\lambda_2) \mbox{    }  | \mbox{    }  \lambda_1,\lambda_2 \in \mathbb{Z}, |\lambda_1|+|\lambda_2| \leq s \}$$ contains exactly $2s^2+2s+1$ elements, it suffices to prove that no integer in $[-(s^2+s),(s^2+s)]$ can be written as an element of $\Sigma$ in two different ways.  

For that, suppose that $$\lambda_1s+\lambda_2(s+1)=\lambda_1's+\lambda_2'(s+1)$$ for some $(\lambda_1,\lambda_2) \in I_{\pm}(2,[0,s])$ and $(\lambda_1',\lambda_2') \in I_{\pm}(2,[0,s])$; w.l.o.g., we can assume that $\lambda_2 \geq \lambda_2'$.  Our equation implies that $\lambda_2 - \lambda_2'$ is divisible by $s$ and is, therefore, equal to 0, $s$, or $2s$.  

If $\lambda_2-\lambda_2'=2s,$ then $\lambda_2=s$ and $\lambda_2'=-s$, which can only happen if $\lambda_1=\lambda_1'=0$; this case leads to a contradiction with our equation.  Assume, next, that $\lambda_2-\lambda_2'=s.$  Our equation then yields $\lambda_1-\lambda_1'=-s-1.$  In this case, we have $\lambda_1 \leq 0$, $\lambda_1' \geq 0$, $\lambda_2 \geq 0$, and $\lambda_2' \leq 0$, and we see that
$$|\lambda_1|+|\lambda_2|=-\lambda_1+\lambda_2=(s+1-\lambda_1')+(s+\lambda_2')=2s+1-(|\lambda_1'|+|\lambda_2'|) \geq s+1,$$ but that is a contradiction.

This leaves us with the case that $\lambda_2=\lambda_2'$, which implies $\lambda_1=\lambda_1'$, as claimed, and therefore the set $\{s,s+1\}$ is $s$-spanning in $\mathbb{Z}_n$.
$\Box$

\addcontentsline{toc}{section}{Proof of Proposition \ref{Z2 X Z2k s spanning}}

\section*{Proof of Proposition \ref{Z2 X Z2k s spanning}} \label{proof of Z2 X Z2k s spanning}

Suppose first that $s \geq 2$, $k \leq 3s-4$, and $k \neq 2s-1$.  We will show that with
$$A=\{(0,1),(1,s-1)\},$$ we have $[0,s]_{\pm}A=\mathbb{Z}_2 \times \mathbb{Z}_{2k}$.  

Let $x$ and $y$ be integers.  Observe that
\begin{enumerate}[(i)]
  \item if $0 \leq x \leq s$, then $$(0,x)=x \cdot (0,1) \in [0,s]_{\pm}A;$$ 
\item if $s \leq x \leq 3s-4$, then $$(0,x)=2 \cdot (1,s-1)+(x-2s+2) \cdot (0,1) \in [0,s]_{\pm}A;$$

  \item if $0 \leq y \leq 2s-2$, then $$(1,y)=1 \cdot (1,s-1)+(y-s+1) \cdot (0,1) \in [0,s]_{\pm}A;$$
\item if $2s \leq y \leq 4s-6$, then $$(1,y)=3 \cdot (1,s-1)+(y-3s+3) \cdot (0,1) \in [0,s]_{\pm}A.$$

\end{enumerate}

Therefore, by (i) and (ii), $(0,x) \in [0,s]_{\pm}A$ for all $0 \leq x \leq 3s-4$, thus, since $2k \leq 6s-8$, by symmetry, $[0,s]_{\pm}A$ contains all group elements whose first component is 0.

If $k \leq 2s-2$, then by (iii) and by symmetry, $[0,s]_{\pm}A$ contains all group elements whose first component is 1, and we are done.

Suppose now that $2s \leq k \leq 3s-4$.  By (iii), $(1,y) \in [0,s]_{\pm}A$ for all $0 \leq y \leq 2s-2$, so, by symmetry, $(1,y) \in [0,s]_{\pm}A$ for all $2k-2s+2 \leq y \leq 2k$.  Note that
$$2k-2s+2 \leq 4s-6,$$ so combining this with (iv), we get that $(1,y) \in [0,s]_{\pm}A$ for all $2s \leq y \leq 2k$.  We are under the assumption that $2s \leq k$, however, so by symmetry,  again $[0,s]_{\pm}A$ contains all group elements whose first component is 1, and $[0,s]_{\pm}A=\mathbb{Z}_2 \times \mathbb{Z}_{2k}$

This completes all cases of Proposition \ref{Z2 X Z2k s spanning}, except for $k=2s-1$ with $s \geq 5$.  We will show that, in this case, for
$$A=\{(0,1),(1,s+1)\},$$ we have $[0,s]_{\pm}A=\mathbb{Z}_2 \times \mathbb{Z}_{2k}$.  

Let $x$ and $y$ be integers.  Observe that
\begin{enumerate}[(i)]
  \item if $0 \leq x \leq s$, then $$(0,x)=x \cdot (0,1) \in [0,s]_{\pm}A;$$ 
\item if $s+1 \leq x \leq s+3$, then $$(0,x)=(-2) \cdot (1,s+1)+(x-2s+4) \cdot (0,1) \in [0,s]_{\pm}A;$$

\item if $s+4 \leq x \leq 3s$, then $$(0,x)=2 \cdot (1,s+1)+(x-2s-2) \cdot (0,1) \in [0,s]_{\pm}A;$$
\item if $0 \leq y \leq 1$, then $$(1,y)=(-3) \cdot (1,s+1)+(y-s+5) \cdot (0,1) \in [0,s]_{\pm}A;$$
  \item if $2 \leq y \leq 2s$, then $$(1,y)=1 \cdot (1,s+1)+(y-s-1) \cdot (0,1) \in [0,s]_{\pm}A.$$

\end{enumerate}
As explanation for (ii), we should add that, since $k=2s-1$, we have 
$$-2(s+1)+(x-2s+4)=-2k+x \equiv x \; \mbox{mod} \; 2k,$$ and, since $s \geq 5$, 
$$|-2|+|x-2s+4| \leq s$$ always holds when $s+1 \leq x \leq s+3$.  The explanation for (iv) is similar.  

The result now follows by symmetry, as above.
$\Box$

\addcontentsline{toc}{section}{Proof of Proposition \ref{rhoZ_nUfixed}}

\section*{Proof of Proposition \ref{rhoZ_nUfixed}} \label{proof of rhoZ_nUfixed}

We will prove a slightly more general result.  

For positive integers $n,m,h, a$, and $b$, consider the geometric progression $$A=\{a,a  b, a  b^2, \dots, a  b^{m-1}\}$$ in the cyclic group $\mathbb{Z}_n$.  (By this we mean that each integer in the set is considered mod $n$; we will talk about integers and their values mod $n$ interchangeably.)  Suppose that $b \geq h \geq 2$ and $$n \geq  hab^{m-1}.$$  We first note that with these conditions, $A$ has size $m$: indeed, we have $$1 \leq a < ab < \cdots < ab^{m-1} < n.$$  Furthermore, we show that two $h$-fold sums of $A$ can only equal if they correspond to the same element of the index set $$\mathbb{N}_0^m(h)=\{(\lambda_1,\dots ,\lambda_m) \in \mathbb{N}_0^m \mbox{    } | \mbox{    } \lambda_1+\cdots +\lambda_m =h \}.$$  First, observe that, for nonnegative integer coefficients $\lambda_1, \dots, \lambda_m$ with $$\lambda_1+\cdots +\lambda_m=h,$$ we have
$$\lambda_1a+\lambda_2ab+\cdots+\lambda_mab^{m-1} \geq \lambda_1a+\lambda_2a+\cdots+\lambda_ma=ha >0$$ and
$$\lambda_1a+\lambda_2ab+\cdots+\lambda_mab^{m-1} \leq \lambda_1ab^{m-1}+\lambda_2ab^{m-1}+\cdots+\lambda_mab^{m-1}=hab^{m-1} \leq n,$$
so linear combinations
$$\lambda_1a+\lambda_2ab+\cdots+\lambda_mab^{m-1} $$ and $$\lambda_1'a+\lambda_2'ab+\cdots+\lambda_m'ab^{m-1} $$ in $hA$ can only equal in $\mathbb{Z}_n$ if they also equal in $\mathbb{Z}$.  But 
$$\lambda_1a+\lambda_2ab+\cdots+\lambda_mab^{m-1}=\lambda_1'a+\lambda_2'ab+\cdots+\lambda_m'ab^{m-1}$$ in $\mathbb{Z}$ implies that $\lambda_1-\lambda'_1$ is divisible by $b$.  Assume, without loss of generality, that $\lambda_1 \geq \lambda_1'$; then $$0 \leq \lambda_1-\lambda'_1 \leq h.$$  In fact, we cannot have $\lambda_1-\lambda'_1 = h$, since that could only happen if $\lambda_1=h$ and $\lambda_1'=0$, but that would imply that 
$$ha=\lambda_2'ab+\cdots+\lambda_m'ab^{m-1},$$ which is impossible since
$$\lambda_2'ab+\cdots+\lambda_m'ab^{m-1} \geq \lambda_2'ab+\cdots+\lambda_m'ab=hab>ha.$$
Therefore, $$0 \leq \lambda_1-\lambda'_1<h \leq b,$$ so the only way for $\lambda_1-\lambda'_1$ to be divisible by $b$ is if $\lambda_1=\lambda_1'$.

This, in turn, implies that $\lambda_i=\lambda_i'$ for every $i=2,3,\dots,m$ as well; therefore, each element of the index set  yields a different element.
Thus we have shown that, under the conditions of $b \geq h \geq 2$ and
$$n \geq hab^{m-1},$$  the geometric progression $$A=\{a,a b, a  b^2, \dots, a  b^{m-1}\}$$ has size $m$ and is a $B_h$ set in the cyclic group $\mathbb{Z}_n$.  (We should note that there is nothing special about $A$ being a geometric progression: all that's needed is that each element in the progression is much larger than the previous one.)

By choosing $a=1$ and $b=h$, we arrive at Proposition \ref{rhoZ_nUfixed}.  $\Box$

\addcontentsline{toc}{section}{Proof of Proposition \ref{rhoZ_nUSfixed}}

\section*{Proof of Proposition \ref{rhoZ_nUSfixed}} \label{proof of rhoZ_nUSfixed}

Again we will prove a slightly more general result.

We consider, for positive integers $n,m,h, a$, and $b$, the geometric progression $$A=\{a,a  b, a  b^2, \dots, a  b^{m-1}\}$$ in the cyclic group $\mathbb{Z}_n$.  Here we assume that $b \geq 2h \geq 2$ and $$n >  2hab^{m-1}.$$  Again we see that, with these conditions, $A$ has size $m$.  Similarly as before, we see that for integer coefficients $\lambda_1, \dots, \lambda_m$ with $$|\lambda_1|+\cdots +|\lambda_m|=h,$$ we have
$$\lambda_1a+\lambda_2ab+\cdots+\lambda_mab^{m-1} \geq -|\lambda_1|ab^{m-1}-|\lambda_2|ab^{m-1}-\cdots-|\lambda_m|ab^{m-1}=-hab^{m-1},$$ and
$$\lambda_1a+\lambda_2ab+\cdots+\lambda_mab^{m-1} \leq |\lambda_1|ab^{m-1}+|\lambda_2|ab^{m-1}+\cdots+|\lambda_m|ab^{m-1}=hab^{m-1},$$
so linear combinations
$$\lambda_1a+\lambda_2ab+\cdots+\lambda_mab^{m-1} $$ and $$\lambda_1'a+\lambda_2'ab+\cdots+\lambda_m'ab^{m-1} $$ in $h_{\pm}A$ can only equal in $\mathbb{Z}_n$ if they also equal in $\mathbb{Z}$.  

But 
$$\lambda_1a+\lambda_2ab+\cdots+\lambda_mab^{m-1}=\lambda_1'a+\lambda_2'ab+\cdots+\lambda_m'ab^{m-1}$$ in $\mathbb{Z}$ implies that $\lambda_1-\lambda'_1$ is divisible by $b$.  Assume, without loss of generality, that $\lambda_1 \geq \lambda_1'$; then $$0 \leq \lambda_1-\lambda'_1 \leq 2h.$$  In fact, we cannot have $\lambda_1-\lambda'_1 = 2h$, since that could only happen if $\lambda_1=h$ and $\lambda_1'=-h$, but that would imply that 
$$ha=-ha,$$ which is impossible.  Therefore, $$0 \leq \lambda_1-\lambda'_1<2h \leq b,$$ so the only way for $\lambda_1-\lambda'_1$ to be divisible by $b$ is if $\lambda_1=\lambda_1'$.

This, in turn, implies that $\lambda_i=\lambda_i'$ for every $i=2,3,\dots,m$ as well; therefore, each element of the index set  yields a different element in the $h$-fold signed sumset, and $$|h_{\pm} A|=|\mathbb{Z}^m(h)|=c(h,m).$$   In summary, we have shown that, under the conditions $b \geq 2h \geq 2$ and $$n  > 2hab^{m-1},$$ the geometric progression $$A=\{a,a b, a  b^2, \dots, a  b^{m-1}\}$$ in the cyclic group $\mathbb{Z}_n$ has size $m$, and its $h$-fold signed sumset $h_{\pm} A$ has size $$|\mathbb{Z}^m(h)|=c(h,m)$$ and thus it is a $B_h$ set over $\mathbb{Z}$ in $\mathbb{Z}_n$.  (We should note that there is nothing special about $A$ being a geometric progression: all that's needed is that each element in the progression is much larger than the previous one.)  

By choosing $a=1$ and $b=2h$, we arrive at Proposition \ref{rhoZ_nUSfixed}.  $\Box$

\addcontentsline{toc}{section}{Proof of Proposition \ref{nupmmaxm=2}}

\section*{Proof of Proposition \ref{nupmmaxm=2}} \label{proofofnupmmaxm=2}

We first provide a proof  for the case when $h$ is even: we will show that, in this case, $A=\{1,h\}$ has an $h$-fold signed sumset of size exactly $4h$.  

We start by observing that the relevant index set is $\mathbb{Z}^2(h)$, which consists of the $4h$ points 
$$\{(0, \pm h), (\pm 1, \pm (h-1)), (\pm 2, \pm (h-2)), \dots, (\pm (h-1), \pm 1), (\pm h,0)\};$$ we may rewrite this layer of the integer lattice as 
$$\mathbb{Z}^2(h)=I \cup J \cup (-I) \cup (-J)$$ where
$$I=\{(i , h-i )  \mid i=0,1,2,\dots,h\}$$ and
$$J=\{(-j ,h-j)  \mid j=1,2,\dots,h-1\}.$$
(Note that $|I|+|J|+|-I|+|-J|=(h+1)+(h-1)+(h+1)+(h-1)=4h.$)

Therefore, we get
$h_{\pm}A= S \cup T \cup (-S) \cup (-T)$, where
$$S=\{i \cdot 1 + (h-i) \cdot h  \mid i=0,1,2,\dots,h\}$$ and
$$T=\{-j \cdot 1 +(h-j) \cdot h  \mid j=1,2,\dots,h-1\}.$$

We need to prove that the four sets $S$, $T$, $-S$, and $-T$ are pairwise disjoint in $\mathbb{Z}_n$.  
First, we prove that $S$ and $T$ are disjoint.  Suppose, indirectly, that this is not so and we have indices $i \in \{0,1,2,\dots,h\}$ and $j \in \{1,2,\dots,h-1\}$ for which
$$i \cdot 1 + (h-i) \cdot h=-j \cdot 1 +(h-j) \cdot h.$$  Rearranging yields that
$$j \cdot (h+1) - i \cdot (h-1)=0$$ in $\mathbb{Z}_n$.  
Now the integer $j \cdot (h+1) - i \cdot (h-1)$ is at least $$1 \cdot(h+1)-h \cdot (h-1)=-h^2+2h+1 >-n$$ and at most $$(h-1) \cdot (h+1)-0 \cdot (h-1)=h^2-1<n,$$ so if $$j \cdot (h+1) - i \cdot (h-1)=0$$ in $\mathbb{Z}_n,$ then $$j \cdot (h+1) - i \cdot (h-1)=0$$ in $\mathbb{Z}$.
Rearranging again, we get $$i +j= (i-j) \cdot h,$$ which implies that $i+j$ is divisible by $h$; since $1 \leq i+j \leq 2h-1$, this can only occur if $i+j=h$.  Our equation then yields $i-j=1$, adding these two equations results in $2i=h+1$, which is a contradiction as the left-hand side is even and the right-hand side is odd.  Therefore, we proved that $S \cap T = \emptyset$.

Observe that the elements of $S$ are between $h$ and $h^2$ (inclusive), and the elements of $T$ are between $1$ and $h^2-h-1$ (inclusive); in particular, all of them are between $1$ and $h^2$.  Therefore, if $S$ and $T$ are disjoint in $\mathbb{Z}_n$ and $n >2h^2$, then the four sets $S$, $T$, $-S$, and $-T$ must be pairwise disjoint.  This completes our proof for the case when $h$ is even.  

Next, we turn to the case when $h$ is odd, and show that $A=\{2,h\}$ has an $h$-fold signed sumset of size exactly $4h$.  
Again, we write $h_{\pm}A$ as $S \cup T \cup (-S) \cup (-T)$, but this time with
$$S=\{i \cdot 2 + (h-i) \cdot h  \mid i=0,1,2,\dots,h\}$$ and
$$T=\{-j \cdot 2 +(h-j) \cdot h  \mid j=1,2,\dots,h-1\}.$$
As above, we can show that the four sets $S$, $T$, $-S$, and $-T$ must be pairwise disjoint.
$\Box$

\addcontentsline{toc}{section}{Proof of Proposition \ref{|hA|=m}}

\section*{Proof of Proposition \ref{|hA|=m}} \label{proofof|hA|=m}

Since the claim is obviously true for $h=1$, we will assume that $h \geq 2$.

Suppose first that $A=a+H$ for some $a \in G$ and $H \leq G$; we then have $|A|=|H|=m$.  Note that
$hA=ha+H$, since $$a_1+\cdots+a_h=(a+h_1)+\cdots+(a+h_h) \in ha+H$$ and
$$ha+h_0=(a+0)+\cdots +(a+0)+(a+h_0) \in hA.$$
Therefore, $$|hA|=|ha+H|=|H|=|A|.$$  

Conversely, suppose that $|hA|=|A|$.  Let $H$ be the stabilizer of $(h-1)A$; that is,
$$H=\{g \in G  \mid g+(h-1)A=(h-1)A\}.$$  Then $H \leq G$.  Choose any $a \in A$; we will show that $A=a+H$.

Consider the set $A'=A-a$.  Then $|A'|=m$ and $0 \in A'$, and therefore
$$(h-1)A =\{0\}+(h-1)A \subseteq A'+(h-1)A.$$
But then
$$|hA|=|hA-a|=|A'+(h-1)A|\geq |(h-1)A| \geq |(h-2) \cdot a+A|=|A|;$$
since we assumed $|hA|=|A|$, equality must hold throughout, and thus $$A'+(h-1)A = (h-1)A,$$ so $A' \subseteq H$ by definition.  This means that $A \subseteq a+H$.  This implies that
$$|a+H| \geq |A|=|hA| =|(h-1)A+A| \geq |(h-1)A|=|H+(h-1)A| \geq |H|=|a+H|.$$  Therefore, equality holds throughout and $|a+H|=|A|$ and thus $a+H=A$.
$\Box$

\addcontentsline{toc}{section}{Proof of Theorem \ref{|hA|=hm-h+1<p}}

\section*{Proof of Theorem \ref{|hA|=hm-h+1<p}} \label{proofof|hA|=hm-h+1<p}

Suppose, first, that $A$ is an arithmetic progression (AP) of length $m$; that is, $$A=\{a+i g  \mid i=0,1,\dots,m-1\}$$  for some $a, g \in \mathbb{Z}_p$.  Then  
$$hA=\{ha+i g  \mid i=0,1,\dots,hm-h\},$$ and thus
$$|hA|=\min\{p,hm-h+1\}=hm-h+1,$$ by assumption.

For the other direction, we will use the Cauchy--Davenport Theorem\index{Cauchy, A--L.}\index{Davenport, H.}  (see \cite{Cau:1813a} and \cite{Dav:1935a}) and Vosper's Theorem\index{Vosper, A. G.} (see \cite{Vos:1956a} or Theorem 2.7 in \cite{Nat:1996a}).  They both refer to the sum $$A+B=\{a+b  \mid a \in A, b \in B \}$$ of subsets $A$ and $B$ in $\mathbb{Z}_p$, and can be stated, as follows.

\noindent {\bf Theorem (Cauchy--Davenport Theorem)} \label{CDthm}\index{Cauchy, A--L.}\index{Davenport, H.}  
{\em Suppose that $p$ is a prime, and let $A, B \subseteq \mathbb{Z}_p$ for which $|A+B| < p$.  Then $|A+B| \geq |A|+|B|-1.$}

\noindent {\bf Theorem (Vosper's Theorem)}\index{Vosper, A. G.} 
{\em Suppose that $p$ is a prime, and let $A, B \subseteq \mathbb{Z}_p$ for which $|A+B| < p$.  Then $|A+B|=|A|+|B|-1$ if, and only if, at least one of the following three conditions holds:
\begin{enumerate}[i]
\item $|A|= 1$ or $|B|=1$;
\item there is an element $c \in \mathbb{Z}_p$ for which $(A+B) \cup \{c\}=\mathbb{Z}_p$ and $\mathbb{Z}_p \setminus B=\{c\}-A$; or
\item $A$ and $B$ are APs with a common difference.
\end{enumerate}}

Suppose now that $$|hA|=hm-h+1<p.$$  We first observe that for any positive integer $k$ and any $a \in A$, we have
$$(k-1)A +\{a\} \subseteq kA,$$ and thus $|(k-1)A| \leq |kA| $.  Therefore, $|kA|<p$ for $k=1,2,\dots,h$, and we 
can repeatedly apply the Cauchy--Davenport Theorem\index{Cauchy, A--L.}\index{Davenport, H.} to get
\begin{eqnarray*}
hm-h+1 &=& |hA|\\
& = & |A+(h-1)A| \\
& \geq & |A|+|(h-1)A|-1 \\
&\geq&  2|A|+|(h-2)A|-2 \\
& \vdots & \\
&\geq& (h-1)|A|+|A|-(h-1) \\
& = & hm-h+1.
\end{eqnarray*} 
Therefore, we must have equality throughout; in particular, $|kA|=km-k+1$ for $k=1,2,\dots,h$.  

In fact, to complete our proof, we will only need the fact that $|2A|=2m-1<p$.  If $m=1$, then we are done since if $A=\{a\}$ for some $a \in \mathbb{Z}_p$, then $A$ is an AP (of length 1).  Assume then that $m \geq 2$.  Using the ``only if '' part of Vosper's Theorem\index{Vosper, A. G.} for $A=B$ (which we are allowed as the hypotheses hold), we see that conditions i and ii cannot occur as condition ii would imply that $p=2m$ and thus $m=1$.  Thus condition iii holds, proving our claim. $\Box$

\addcontentsline{toc}{section}{Proof of Theorem \ref{thm Gryn |hA|=p}}

\section*{Proof of Theorem \ref{thm Gryn |hA|=p}}  \label{proof of thm Gryn |hA|=p}

Recall that the stabilizer of a nonempty subset $S$ in $G$ is the set $$H=\{g \in G \mid g+S=S\}.$$  It is easy to verify that the stabilizer of $S$ is a subgroup of $G$, and we have $$H+S=S;$$ that is, $S$ is the union of some cosets of $H$.

We shall use the following famous result of Kneser (cf.~\cite{Kne:1953a} or Theorem 4.5 in \cite{Nat:1996a}):

\noindent {\bf Theorem (Kneser's Theorem)}\index{Kneser, M.}
{\em Suppose that $A$ is an $m$-subset of a finite abelian group $G$, $h$ is a positive integer, and $H$ is the stabilizer subgroup of $hA$ in $G$.  We then have}
$$|hA| \geq h|A+H|-(h-1)|H|.$$

Suppose now that $A$ is an $m$-subset of $G$ and that $hA$ has size $p$, where $p$ is the smallest prime divisor of the order $n$ of $G$.  We also assume that $m \leq p < hm-h+1$.  Note that if the stabilizer $H$ of $hA$ were to have order 1, then by Kneser's Theorem, we would get $$p=|hA| \geq h|A+H|-(h-1)|H| \geq hm-(h-1),$$ contradicting our assumption that $p < hm-h+1$.  Therefore, $|H|>1$, and by the definition of $p$, we have $|H| \geq p$.  Since $H$ is the stabilizer of $hA$, $hA$ is the union of some cosets of $H$; with $|H| \geq p$ and $|hA|=p$ this is only possible if $|H|=p$ and $hA$ equals a coset of $H$.

Note that if $A$ would not lie in a single coset of $H$, then it would have elements $a_1$ and $a_2$ for which the cosets $a_1+H$ and $a_2+H$ are different, but then the cosets $(h-1)a_1+a_2+H$ and $ha_1+H$ would be different too, so $hA$ could not lie in a single coset of $H$, which is a contradiction.  $\Box$

\addcontentsline{toc}{section}{Proof of Theorem \ref{thm |hA|=hm-h+1}}

\section*{Proof of Theorem \ref{thm |hA|=hm-h+1}} \label{proof of thm |hA|=hm-h+1}

Suppose, first, that $A$ is an arithmetic progression (AP) of length $m$; that is, $$A=\{a+i g  \mid i=0,1,\dots,m-1\}$$  for some $a, g \in G$.  Then  
$$hA=\{ha+i g  \mid i=0,1,\dots,hm-h\},$$ and thus
$$|hA|=\min\{p,hm-h+1\}=hm-h+1,$$ by assumption.

For the other direction, we will use the Corollary on page 74 of Kemperman's work \cite{Kem:1960a}:

\noindent {\bf Theorem (Kemperman; cf.~\cite{Kem:1960a})}  \label{Kemp thm}\index{Kemperman, J. H. B.} 
{\em Let $p$ be the minimum prime divisor of $n$.   Suppose that $m \geq 2$ and $$|2A| \leq \min\{p-2, 2m-1\}$$ for some $A \subseteq G$.  Then $A$ is an AP.}

Note that our assumption of $hm-h+1 <p$ implies that $2m-1<p$ and, since $p$ is odd, we then have $2m+1 \leq p$ or $p-2 \geq 2m-1$.  Hence, the conditions of Kemperman's Theorem\index{Kemperman, J. H. B.} above are met.  
The proof of Theorem \ref{thm |hA|=hm-h+1} then follows as the proof of Theorem \ref{|hA|=hm-h+1<p} above.

\addcontentsline{toc}{section}{Proof of Theorem \ref{Matzke limited prime}}

\section*{Proof of Theorem \ref{Matzke limited prime}} \label{proof of Matzke limited prime}

We can evaluate $u_{\pm}(p,m,[0,s])$ for odd prime values of $p$: we have
$$D_1(p)=\left\{
\begin{array}{cll}
\{1,p\} & \mbox{if} & m \; \mbox{is odd}, \\ \\
\{p\} & \mbox{if} & m \; \mbox{is even};
\end{array} \right.$$
and
$$D_2(p)=\left\{
\begin{array}{cll}
\emptyset & \mbox{if} & m \; \mbox{is odd}, \\ \\
\{1\} & \mbox{if} & m \; \mbox{is even}.
\end{array} \right.$$
Therefore, when $m$ is odd, we get
$$u_{\pm}(p,m,[0,s])=\min\{f_1(m,s),f_p(m,s)\}=\min\{p,sm-s+1\}=\min\{p,2s  \lfloor m/2 \rfloor +1\},$$ and 
when $m$ is even, we get
$$u_{\pm}(p,m,[0,s])=\min\{f_1(m+1,s),f_p(m,s)\}=\min\{p,sm+1\}=\min\{p,2s  \lfloor m/2 \rfloor +1\}.$$  This means that, by Theorem \ref{matzke upper for [0,s]}, we have
$$\rho_{\pm} (\mathbb{Z}_p,m,[0,s]) \leq \min\{p,2s  \lfloor m/2 \rfloor +1\}.$$

Now we prove that for every $m$-subset $A$ of $\mathbb{Z}_p$, we have 
$$|[0,s]_{\pm} A| \geq \min \{p, 2s \lfloor m/2 \rfloor +1\}.$$

When $m$ is odd, this follows immediately from the Cauchy--Davenport Theorem, since
$$|[0,s]_{\pm} A| \geq |sA| \geq \min \{p, sm-s +1\}=\min \{p, 2s \lfloor m/2 \rfloor +1\}.$$

Observe that when $m$ is even, then $A$ is a proper subset of $A \cup (-A) \cup \{0\}$; therefore, by Propositions 3.3 and 3.4, we get
$$|[0,s]_{\pm} A| = |s(A \cup (-A) \cup \{0\}) |  \geq \min \rho (\mathbb{Z}_p,m+1,s).$$
Our claim again follows from the Cauchy--Davenport Theorem, since
$$\rho (\mathbb{Z}_p,m+1,s)  =\min \{p, s (m+1)-s +1\}=\min \{p, 2s \lfloor m/2 \rfloor +1\}.$$
Our proof is now complete.  
$\Box$

\addcontentsline{toc}{section}{Proof of Proposition \ref{rhohat m=4}}

\section*{Proof of Proposition \ref{rhohat m=4}} \label{proof of rhohat m=4}

Let us write an arbitrary 4-subset $A$ of $G$ in the form
$$A=\{a,a+d_1,a+d_2,a+d_3\}.$$ Here $d_1$, $d_2$, and $d_3$ are distinct nonzero elements of $G$.

We first prove the following.

{\bf Proposition} \label{h=2 m=4 size}
{\em Suppose that $d_1$, $d_2$, and $d_3$ are distinct nonzero elements of $G$, and let $A=\{a,a+d_1,a+d_2,a+d_3\}.$
\begin{enumerate}

\item If none of $d_1$, $d_2$, or $d_3$ equals the sum of the other two, then $|2 \hat{\;} A|=6$.

\item Suppose (w.l.o.g.) that $d_3=d_1+d_2$.
\begin{enumerate}
\item If both $d_1$ and $d_2$ have order 2, then $|2 \hat{\;} A|=3$;
\item if exactly one of $d_1$ or $d_2$ have order 2, then $|2 \hat{\;} A|=4$;
\item if neither of $d_1$ or $d_2$ have order 2, then $|2 \hat{\;} A|=5$.
\end{enumerate}

\end{enumerate}}

{\em Proof of Proposition:}  Since $2a+d_1$, $2a+d_2$, and $2a+d_3$ are distinct, 
$$2 \hat{\;} A=\{2a+d_1, 2a+d_2, 2a+d_3, 2a+d_1+d_2, 2a+d_1+d_3, 2a+d_2+d_3\}$$ has size 3, 4, 5, or 6, exactly when three, two, one, or none of the equations
$$d_3=d_1+d_2, \; \; \; d_2=d_1+d_3, \; \; \; \mbox{or} \; d_1=d_2+d_3$$
hold, respectively.  This proves 1.  

Assume now that $d_3=d_1+d_2$, and that both $d_1$ and $d_2$ have order 2.  Then $$d_1+d_3=d_1+(d_1+d_2)=d_2;$$ similarly, $d_2+d_3=d_1$.  This proves 2 (a).

If $d_3=d_1+d_2$, and (w.l.o.g.) $d_1$ has order 2 but $d_2$ does not, then we still have $d_1+d_3=d_2$, but $$d_2+d_3=d_2+(d_1+d_2) \neq d_1.$$  

Finally, if neither of $d_1$ or $d_2$ have order 2, then $d_2 \neq d_1+d_3$ and $d_1 \neq d_2+d_3$. $\Box$

{\em Proof of Proposition \ref{rhohat m=4}.}  By our proposition above, it suffices to prove that 
$$\rho\hat{\;} (G,4,2) \leq
\left\{\begin{array}{ll}
3  & \mbox{if $\; |\mathrm{Ord}(G,2)| \geq 2$,}\\ 
4  & \mbox{if $\; |\mathrm{Ord}(G,2)| =1$,}\\ 
5  & \mbox{if $\; |\mathrm{Ord}(G,2)| =0$.}  
\end{array}\right.$$

Suppose, first, that $|\mathrm{Ord}(G,2)| \geq 2$, and let $d_1$ and $d_2$ be two distinct elements of $\mathrm{Ord}(G,2)$.  Then $d_3=d_1+d_2$ (which is also of order 2) is nonzero, and distinct from $d_1$ or $d_2$.  Thus, by our proposition above, the set $A=\{0,d_1,d_2,d_3\}$ has $|2 \hat{\;} A|=3$.

Suppose now that $d_1$ is the unique element of $G$ of order 2; let $d_2$ be any other nonzero element of $G$ (exists since $n \geq m=4$).  Again, $d_3=d_1+d_2$ is nonzero, and is distinct from $d_1$ or $d_2$, and thus, by our proposition above, the set $A=\{0,d_1,d_2,d_3\}$ has $|2 \hat{\;} A|=4$.

Finally, assume that $\mathrm{Ord}(G,2) = \emptyset$, and let $d_1$ and $d_2$ be distinct nonzero elements of $G$; assume further that $d_2 \neq -d_1$ (this is possible since $n \geq m=4$).  As before, $d_3=d_1+d_2$ is nonzero, and is distinct from $d_1$ or $d_2$, and thus, by our proposition above, the set $A=\{0,d_1,d_2,d_3\}$ has $|2 \hat{\;} A|=5$. $\Box$

\addcontentsline{toc}{section}{Proof of Theorem \ref{uhattheorem}}

\section*{Proof of Theorem \ref{uhattheorem}} \label{proofofuhattheorem}

We first state and prove the following lemma.

{\bf Lemma} \label{lemma any j}
{\em Suppose that $d$ and $t$ are positive integers with $t \leq d-1$, and let $j \in \mathbb{Z}_d$.  Then there is a $t$-subset $J=\{j_1,\dots,j_t\}$ of $\mathbb{Z}_d$ for which $$j_1+j_2+\cdots+j_t=j.$$}

Note that the restriction of $t \leq d-1$ is necessary: for $t=d$ the only $j \in \mathbb{Z}_d$ for which such a set exists is, of course, $$j=0+1+\cdots+(d-1)=\frac{d(d-1)}{2}=\left\{
\begin{array}{ll}
\frac{d}{2} & \mbox{if $d$ is even}; \\ \\
0 & \mbox{if $d$ is odd}.
\end{array}\right.$$

{\em Proof of lemma:}  Write $j$ as $$j=\frac{t^2-t}{2}+j_0 \; \mbox{mod}\; d$$ where $j_0=0,1,\dots,d-1$.  Note that $1 \leq t \leq d-1$.  We will separate two cases: when $j_0$ is less than $t$ and when it is not.

If $0 \leq j_0 \leq t-1$, let $$J=\{0,1,\dots,t-1,t\} \setminus \{t-j_0\}.$$  Then $|J|=t$, and the elements of $J$ add up to $$\frac{t^2+t}{2}-(t-j_0)=j.$$

If $t \leq j_0 \leq d-1$, take $$J=\{1,\dots,t-1\} \cup \{j_0\}.$$  Again, $|J|=t$, and the elements of $J$ add up to $j.$
$\Box$

Before we turn to the proof of Proposition \ref{uhattheorem}, let us recall our notations and perform some computations.  We write

$$m=dc+k \; \mbox{with} \; c=\left \lceil \frac{m}{d}  \right \rceil -1,$$
and
$$h=dq+r \; \mbox{with} \; q=\left \lceil \frac{h}{d}  \right \rceil -1.$$
Note that $ c \geq 0$, $q \geq 0$, $1 \leq k \leq d$, and $1 \leq r \leq d$.

Recall that we have set $$A=A_d(n,m)=\bigcup_{i=0}^{c-1} (i+H) \cup \left\{c + j \cdot \frac{n}{d}  \mid j=0,1,2,\dots,k-1 \right\}.$$ Here $$\bigcup_{i=0}^{c-1} (i+H)=\emptyset$$ when $c=0$ (that is, when $m \leq d$), but $$\left\{c + j \cdot \frac{n}{d}  \mid j=0,1,2,\dots,k-1 \right\} \neq \emptyset.$$

Note that every element of $h \hat{\;} A$ is of the form 
$$(i_1+i_2+\cdots+i_h)+(j_1+j_2+\cdots+j_h) \cdot \frac{n}{d}$$
with $i_1,\dots,i_h \in \{0,1,\dots,c\}$ and $j_1,\dots,j_h \in \{0,1,\dots,d-1\}$, with the added conditions that when any of the $i$-indices equals $c$, the corresponding $j$-index is at most $k-1$, and that when two $i$-indices are equal, the corresponding $j$-indices are distinct.

Clearly, the least value  of $i_1+\dots+i_h$ is 
$$i_{\mathrm{min}}=d(0+1+\cdots+(q-1))+rq=q \cdot \frac{h+r-d}{2}.$$ 

To compute the largest value $i_{\mathrm{max}}$ of $i_1+\dots+i_h$, we consider four cases depending on whether $r >k$ or not and whether $q=0$ or not.  

First, when $q=0$ and $r>k$, then $r=h$ and $1 \leq h-k =r-k<r \leq d$, so it is easy to see that
$$i_{\mathrm{max}}=kc+(h-k)(c-1)=h(c-1)+k.$$

In the case when $r >k$ and $q \geq 1$, we write $h$ as $h=k+dq+(r-k)$; thus
\begin{eqnarray*}
i_{\mathrm{max}} & = & k c+[(c-1)+(c-2)+\cdots+(c-q)]d+(r-k)(c-q-1) \\
& = & k c+q c d-\frac{q(q+1)}{2}d+r c-r q-r-k c+k q+k \\
& = & qcd+rc-q \cdot \frac{dq+d+2r}{2}-r+kq+k \\
& = & hc-q \cdot \frac{h+r-d}{2}-dq-r+kq+k \\
& = & h(c-1)-q \cdot \frac{h+r-d}{2}+kq+k.
\end{eqnarray*} 

Next, when  $q=0$ and $r \leq k$, then $r=h$ and $h =r \leq k$, so we have
$$i_{\mathrm{max}}=hc=h(c-1)+r.$$

Finally, in the case when $q \geq 1$ and $r \leq k$, then $0 < d-k+r  \leq d$; we write $h$ as $h=k+d(q-1)+(d-k+r)$ and thus (using our result from the first case above)
\begin{eqnarray*}
i_{\mathrm{max}} & = & k c+[(c-1)+(c-2)+\cdots+(c-q+1)]d+(d+r-k)(c-q) \\
& = & k c+[(c-1)+(c-2)+\cdots+(c-q+1)+(c-q)]d+(r-k)(c-q-1)+(r-k) \\
& = & h(c-1)-q \cdot \frac{h+r-d}{2}+kq+r.
\end{eqnarray*} 
All four cases can be summarized by the formula
$$i_{\mathrm{max}}=h(c-1)-q \cdot \frac{h+r-d}{2}+kq+ \min\{r,k\}.$$

Clearly, when $m \leq d$, then $i_{\mathrm{min}}=i_{\mathrm{max}}=0$.  But when $m>d$, we can verify that $i_{\mathrm{max}} > i_{\mathrm{min}}$, as follows.  We have
\begin{eqnarray*}
i_{\mathrm{max}}-i_{\mathrm{min}} & = &  h(c-1) - q \cdot (h+r-d) +kq + \min\{r,k\} \\ 
& = &  h \cdot (c-q-1)+q(d-r+k)+\min\{r,k\};
\end{eqnarray*}
this quantity is positive when $c \geq q+1$.  Note that we must have $c- q =\left \lceil \frac{m}{d}  \right \rceil - \left \lceil \frac{h}{d}  \right \rceil \geq 0$, thus the only remaining case is when $c=q$, in which case $k-r=m-h>0$ and $q>0$ (since $q=c$ and $m>d$), so we now have
$$i_{\mathrm{max}}-i_{\mathrm{min}}=-h+q(d-r+k)+\min\{r,k\}=-h+qd+q(k-r)+r=q(k-r) > 0.$$

Obviously, $i=i_1+i_2+\cdots+i_h$ can assume the value of any integer between these two bounds, and thus $h \hat{\;} A$ lies in exactly $$\min \left\{\frac{n}{d},i_{\mathrm{max}}-i_{\mathrm{min}}+1 \right\}$$ cosets of $H$. 

{\em Proof of Proposition \ref{uhattheorem}.}  We can easily check that the result holds for $h=1$, so below we assume that $2 \leq h \leq m-1$.  We will separate the rest of the proof into several cases.  In the first two cases, we have $h \leq k \leq d$, and thus $q=0$ and $h=r$, so we have $i_{\mathrm{min}}=0$ and $i_{\mathrm{max}}=hc$.

{\bf Claim 1}: If $h \leq k$ and $h<d$, then $|h \hat{\;} A| = \min\{n, hcd+d, hm-h^2+1\}$.

We here recognize the quantity $f_d=hcd+d$.

{\em Proof of Claim 1}:  Note that the assumptions, using our lemma above, imply that
$$h \hat{\;} A =\bigcup_{i=0}^{hc-1} (i+H) \cup \left\{hc+j \cdot \frac{n}{d}  \mid j=\frac{h(h-1)}{2}, \dots, h(k-1)-\frac{h(h-1)}{2} \right\}.$$

We first consider the case when $c=0$ or, equivalently, when $m \leq d$, i.e. when $m=k$.  In this case 
$$h \hat{\;} A = \left\{j \cdot \frac{n}{d}  \mid j=\frac{h(h-1)}{2}, \dots, h(m-1)-\frac{h(h-1)}{2} \right\},$$ and therefore
$$|h \hat{\;} A| = \min\{d, hm-h^2+1\}=\min\{n, d, hm-h^2+1\}=\min\{n, hcd+d, hm-h^2+1\},$$ as claimed.

Assume now that $c \geq 1$ (iff $m >d$, iff $m>k$).  If we also have $$hc-1 \geq \frac{n}{d}-1,$$  then $$|h \hat{\;} A|=n.$$  But note that, when $k \geq h$ and $$hc-1 \geq \frac{n}{d}-1,$$ then we also have
$$hm-h^2+1=hcd+hk-h^2+1 \geq hcd+1 >n$$ and
$$hcd+d \geq n+d >n,$$ and thus
$$|h \hat{\;} A|=n=\min\{n, hcd+d, hm-h^2+1\},$$
as claimed.

If, on the other hand, $c \geq 1$ and $$hc-1 \leq \frac{n}{d}-2,$$  then
$$|h \hat{\;} A| = hcd+ \min\{d, hk-h^2+1\}=\min\{ hcd+d, hcd+hk-h^2+1\}=\min\{ hcd+d, hm-h^2+1\}.$$ 
Note that, under the assumption that  $$hc-1 \leq \frac{n}{d}-2,$$ we have
$$\min\{ hcd+d, hm-h^2+1\} \leq hcd+d \leq n,$$ and therefore again we get
$$|h \hat{\;} A|=n=\min\{n, hcd+d, hm-h^2+1\},$$
as claimed.  This completes the proof of Claim 1.

{\bf Claim 2}: If $h=k=d$, then $|h \hat{\;} A| = \min\{n, hm-h^2-h+2\}$.

{\em Proof of Claim 2}:  First, we note that $h<m=dc+k=dc+h$, so $c \geq 1$.  In the case of $h=k=d$, our lemma above cannot be used for the coset $i+H$ when $i=0$; we in fact now have  
$$h \hat{\;} A =\left\{\frac{d(d-1)}{2} \cdot \frac{n}{d} \right\} \cup \bigcup_{i=1}^{dc-1} (i+H) \cup \left\{dc+ \frac{d(d-1)}{2} \cdot \frac{n}{d}  \right\}.$$

If we have $$dc-1 \geq \frac{n}{d},$$ then clearly $|h \hat{\;} A|=n$, but we also have
$$hm-h^2-h+2=h(dc+h)-h^2-h+2=h(dc-1)+2 \geq n+2>n,$$
so 
$$|h \hat{\;} A| = \min\{n, hm-h^2-h+2\},$$ as claimed.

Assume now that 
 $$dc-1 \leq \frac{n}{d}-2.$$
In this case
$$|h \hat{\;} A| =1+(dc-1)d+1 = h(cd+h)-h^2-h+2=hm-h^2-h+2;$$ furthermore, 
$$1+(dc-1)d+1 \leq n-2d+2  \leq n,$$  thus
$$|h \hat{\;} A| = \min\{n, hm-h^2-h+2\},$$ as claimed.

This leaves us with the case of $$dc-1 = \frac{n}{d}-1,$$  when we have
\begin{eqnarray*}
h \hat{\;} A &=& \left\{\frac{d(d-1)}{2} \cdot \frac{n}{d} \right\} \cup \bigcup_{i=1}^{dc-1} (i+H) \cup \left\{dc+ \frac{d(d-1)}{2} \cdot \frac{n}{d}  \right\} \\ \\
&=& \left\{\frac{d(d-1)}{2} \cdot \frac{n}{d} \right\} \cup \bigcup_{i=1}^{\frac{n}{d}-1} (i+H) \cup \left\{\frac{n}{d}+ \frac{d(d-1)}{2} \cdot \frac{n}{d}  \right\} \\ \\
&=& \left\{j \cdot \frac{n}{d} \mid j= \frac{d(d-1)}{2}, \frac{d(d-1)}{2}+1 \right\} \cup \bigcup_{i=1}^{\frac{n}{d}-1} (i+H),
\end{eqnarray*}
and $d=h \geq 2$, so
$$|h \hat{\;} A|=2+\left(\frac{n}{d}-1 \right) \cdot d=2+(dc-1)d= cdh-h+2=(cd+h)h-h^2-h+2=hm-h^2-h+2.$$
Furthermore, $2+\left(\frac{n}{d}-1 \right) \cdot d=n-(d-2) \leq n,$ so again we have
$$|h \hat{\;} A| = \min\{n, hm-h^2-h+2\},$$ completing the proof of Claim 2.

We can also observe that, since we always have $k \leq d$, Claims 1 and 2 cover all possibilities under the assumption $h \leq k$.

{\bf Claim 3}: If $h>k$, $r \neq d$, and $r \neq k$, then $$|h \hat{\;} A| = \min\{n, hcd+d-(h-r)(h+r-k)-d \cdot \max\{0,r-k\}\}.$$

{\em Proof of Claim 3}:  First, observe that by our lemma above, the three conditions imply that 
$$h \hat{\;} A=\bigcup_{i=i_{\mathrm{min}}}^{i_{\mathrm{max}}} (i+H).$$
Therefore, we just need to prove that
$$(i_{\mathrm{max}}-i_{\mathrm{min}}+1)\cdot d= hcd+d-(h-r)(h+r-k)-d \cdot \max\{0,r-k\}.$$
Indeed, from our computations above we get
\begin{eqnarray*}
i_{\mathrm{max}}-i_{\mathrm{min}}+1 & = & h(c-1)-q \cdot(h+r-d)+kq+ \min\{r,k\} +1\\
& = & hc - h +qd -q \cdot(h+r-k)+ \min\{r,k\} +1\\
& = & hc - r -q \cdot(h+r-k)+ \min\{r,k\} +1\\
& = & hc+1 -\frac{h-r}{d} \cdot(h+r-k)+ \min\{0,k-r\},
\end{eqnarray*}
from which our result follows.

{\bf Claim 4}: If $h>k$, $r=d$, and $r \neq k$, then $$|h \hat{\;} A| = \min\{n, hcd+d-(h-r)(h+r-k)-d \cdot \max\{0,r-k\}-(d-1)\}.$$

{\em Proof of Claim 4}:  The only difference between the conditions here and the conditions for Claim 3 is that this time we have 
$$h \hat{\;} A=\{ x_{\mathrm{min}} \} \cup \bigcup_{i=i_{\mathrm{min}}+1}^{i_{\mathrm{max}}} (i+H)$$
where $x_{\mathrm{min}}$ equals the sum of the $h$ elements of the set
$$\bigcup_{i=0}^{q} (i+H).$$  
Therefore, 
$$|h \hat{\;} A|=
\left\{\begin{array}{cll}
n  & \mbox{if} & i_{\mathrm{max}}-i_{\mathrm{min}} \geq \frac{n}{d}, \\  \\
(i_{\mathrm{max}}-i_{\mathrm{min}})\cdot d +1  & \mbox{if} & i_{\mathrm{max}}-i_{\mathrm{min}} \leq \frac{n}{d} -1 
\end{array}\right.$$
or, equivalently,
$$|h \hat{\;} A|=\min\{n, (i_{\mathrm{max}}-i_{\mathrm{min}})\cdot d +1 \}.$$
Our claim then follows, since, as can be seen from the proof of Claim 3,
$$(i_{\mathrm{max}}-i_{\mathrm{min}})\cdot d +1=(i_{\mathrm{max}}-i_{\mathrm{min}}+1)\cdot d -(d-1)=hcd+d-(h-r)(h+r-k)-d \cdot \max\{0,r-k\}-(d-1).$$

{\bf Claim 5}: If $h>k$, $r \neq d$, and $r = k$, then $$|h \hat{\;} A| = \min\{n, hcd+d-(h-r)(h+r-k)-d \cdot \max\{0,r-k\}-(d-1)\}.$$

{\em Proof of Claim 5}:  The only difference between the conditions here and the conditions for Claim 3 is that this time we have 
$$h \hat{\;} A= \bigcup_{i=i_{\mathrm{min}}}^{i_{\mathrm{max}}-1} (i+H) \cup \{ x_{\mathrm{max}} \}$$
where $x_{\mathrm{max}}$ equals the sum of the $h$ elements of the set
$$\bigcup_{i=c-q}^{c-1} (i+H) \cup \left\{c + j \cdot \frac{n}{d}  \mid j=0,1,2,\dots,k-1 \right\}.$$  Our claim then follows as in Claim 5.

{\bf Claim 6}: If $h>k$, $r = d$, and $r = k$, then 
$$|h \hat{\;} A| = 
\left\{
\begin{array}{cl}
\min\{n,n-(d-2)\}  & \mbox{if} \; m = \frac{n}{h}+h, \\  \\
\min\{n, hcd+d-(h-r)(h+r-k)-d \cdot \max\{0,r-k\}-2(d-1)\}  & \mbox{otherwise.} 
\end{array}\right.
$$

{\em Proof of Claim 6}:  This time we get
$$h \hat{\;} A=\{ x_{\mathrm{min}} \} \cup \bigcup_{i=i_{\mathrm{min}}+1}^{i_{\mathrm{max}}-1} (i+H) \cup \{ x_{\mathrm{max}} \}$$
where $x_{\mathrm{min}}$ and $x_{\mathrm{max}}$ were defined in Claims 4 and 5, respectively.

With $r=k=d$, we can simplify the expression we got for $i_{\mathrm{max}}-i_{\mathrm{min}}+1$ in the proof of Claim 3:
\begin{eqnarray*}
i_{\mathrm{max}}-i_{\mathrm{min}}+1 & = & hc+1 -\frac{h-r}{d} \cdot(h+r-k)+ \min\{0,k-r\} \\
& = & hc+1-\frac{(h-d)h}{d} \\
& = & \frac{h(cd+d)-h^2 +d}{d} \\
& = & \frac{hm-h^2 +d}{d}.
\end{eqnarray*} 
Thus we see that, if $m \neq \frac{n}{h}+h,$ then $i_{\mathrm{max}}-i_{\mathrm{min}}-1 \neq \frac{n}{d}-1$.  Moreover, if $m > \frac{n}{h}+h,$ then $i_{\mathrm{max}}-i_{\mathrm{min}}-1 \geq \frac{n}{d}$, and thus
\begin{eqnarray*}
|h \hat{\;} A| & = & n \\
& =&  \min\{n, (i_{\mathrm{max}}-i_{\mathrm{min}}-1) \cdot d +2 \}  \\
& =&  \min\{n, (i_{\mathrm{max}}-i_{\mathrm{min}}+1) \cdot d -2(d-1) \}  \\
& = & \min\{n, hcd+d-(h-r)(h+r-k)-d \cdot \max\{0,r-k\}-2(d-1)\} .
\end{eqnarray*}
If, on the other hand, $m < \frac{n}{h}+h,$ then $i_{\mathrm{max}}-i_{\mathrm{min}}-1 \leq \frac{n}{d}-2$, and thus
\begin{eqnarray*}
|h \hat{\;} A| & = & (i_{\mathrm{max}}-i_{\mathrm{min}}-1) \cdot d +2 \\
& =&  \min\{n, (i_{\mathrm{max}}-i_{\mathrm{min}}-1) \cdot d +2 \}  \\
& =&  \min\{n, (i_{\mathrm{max}}-i_{\mathrm{min}}+1) \cdot d -2(d-1) \}  \\
& = & \min\{n, hcd+d-(h-r)(h+r-k)-d \cdot \max\{0,r-k\}-2(d-1)\} .
\end{eqnarray*}
This leaves us with the case of $m=\frac{n}{h}+h$.  In this case
\begin{eqnarray*}
h \hat{\;} A & =& \{ x_{\mathrm{min}} \} \cup \bigcup_{i=i_{\mathrm{min}}+1}^{i_{\mathrm{max}}-1} (i+H) \cup \{ x_{\mathrm{max}} \} \\
& =& \bigcup_{i=i_{\mathrm{min}}+1}^{i_{\mathrm{min}}+\frac{n}{d}-1} (i+H) \cup \{x_{\mathrm{min}},  x_{\mathrm{max}} \}.
\end{eqnarray*}

A simple computation shows that, denoting the sum of the elements in a subset $S$ of $\mathbb{Z}_n$ by $\sum S$, we have
$$x_{\mathrm{min}} = \sum \bigcup_{i=0}^{q} (i+H) = \frac{dq(q+1)}{2} + \frac{d(d-1)(q+1)}{2} \cdot \frac{n}{d}$$
and
$$x_{\mathrm{max}} =\sum \bigcup_{i=c-q}^{c} (i+H)  = cd(q+1) - \frac{dq(q+1)}{2} + \frac{d(d-1)(q+1)}{2} \cdot \frac{n}{d}.$$
But
$$cd=m-d=\frac{n}{h}+h-d=\frac{n}{d(q+1)}+dq,$$
thus
$$x_{\mathrm{max}} = \frac{dq(q+1)}{2} + \left( \frac{d(d-1)(q+1)}{2} +1 \right) \cdot \frac{n}{d},$$
showing that $x_{\mathrm{min}} = x_{\mathrm{max}}$ if, and only if, $d=1$.
Therefore, when $d \geq 2$, we get
$$|h \hat{\;} A|=\left(\frac{n}{d}-1 \right)d+2=n-d+2,$$ and when $d=1$ we get $$|h \hat{\;} A|=\left(\frac{n}{d}-1 \right)d+1=n.$$
This completes the proof of Claim 6.  

Now it is an easy exercise to verify that in all cases, $|h \hat{\;} A|$ is as claimed in the statement of Theorem \ref{uhattheorem}.

\addcontentsline{toc}{section}{Proof of Proposition \ref{rho hat for chi hat -1}}

\section*{Proof of Proposition \ref{rho hat for chi hat -1}} \label{proof of rho hat for chi hat -1}

Our task is to prove that, when $n$ is not a power of 2, then
$$\rho \hat{\;} (G,m_0,2) \leq n-2,$$ where $$m_0=\frac{n+|\mathrm{Ord}(G,2)|+1}{2}.$$

We use induction on the rank of $G$.  When $G$ is cyclic, the claim follows easily from Corollary \ref{cor h=2}; for the sake of completeness---and since the rest of our proof is constructive---we exhibit a subset $A$ of $\mathbb{Z}_n$ of size $m_0$ for which $2 \hat{\;} A$ is of size $n-2$.  

When $n$ is odd, we have $m_0=(n+1)/2$, and the set $$A=\{0,1,2,\dots, (n-1)/2 \}$$ has restricted sumset
$$2 \hat{\;} A = \{1,2,\dots, n-2\},$$  as claimed.

Suppose now that $n$ is even, and set $n=2^k \cdot d$ where $k \geq 1$ and $d \geq 3$ is odd.  In this case, $$m_0=n/2+1=2^{k-1} \cdot d+1.$$
Consider the set
$$A=B_d(n,m_0;(d+1)/2, (d+1)/2, 1, (d-1)/2),$$ defined on page \pageref{def of B_d sets}: in our case, this set becomes
$$A=B' \cup \bigcup_{i=1}^{2^{k-1}-1} \{i+ j \cdot 2^k \mid j=0,1,\dots,d-1\} \cup B'',$$
with
$$B'=\{j \cdot 2^k \mid j=0,1,\dots,(d-1)/2\},$$ and
$$B''=\{2^{k-1}+j \cdot 2^k \mid j=(d-1)/2, (d+1)/2, \dots, d-1\}.$$
Then $$|A|=(d+1)/2+(2^{k-1}-1) \cdot d+(d+1)/2 =2^{k-1} \cdot d+1=m_0,$$ and 
$$2 \hat{\;} A = \mathbb{Z}_n \setminus \{0, (d-1) \cdot 2^k\},$$
since $$2 \hat{\;} B'=\{j \cdot 2^k \mid j=1,2, \dots,d-2\}$$ and
$$2 \hat{\;} B''=\{2^k+j \cdot 2^k \mid j=d, d+1, \dots, 2d-3\}=2 \hat{\;} B'.$$
This completes the case when $G$ is cyclic.

Assume now that our claim holds for groups of rank $r-1$, and let $G$ be a group of rank $r \geq 2$.   Suppose further that $G$ is of type $(n_1, \dots, n_r)$; we write 
$$G=\mathbb{Z}_{n_1} \times G_2$$ where $G_2$ has order $n/n_1$ and rank $r-1$.  Note that $n_1$ is a divisor of $n/n_1$, so if $n$ is not a power of 2, then $n/n_1$ is not a power of 2 either.  Therefore, by our inductive assumption, $G_2$ contains a subset $A_2$ of size $$|A_2|=\frac{n/n_1+|\mathrm{Ord}(G_2,2)|+1}{2}$$ for which $$2 \hat{\;} A_2 = G_2 \setminus X$$ for some $X \subseteq G_2$ of size $|X| \geq 2$.

We separate two cases depending on the parity of $n_1$.

Observe that when $n_1$ is odd, then $$|\mathrm{Ord}(G,2)|=|\mathrm{Ord}(G_2,2)|,$$ so $$m_0=\frac{n+|\mathrm{Ord}(G_2,2)|+1}{2}.$$   Set
$$A=\left( \{0\} \times A_2 \right) \cup \left( \{1,2,\dots,(n_1-1)/2\} \times G_2 \right).$$  The size of $A$ is then
$$|A|=|A_2|+ (n_1-1)/2  \cdot |G_2| = \frac{n/n_1+|\mathrm{Ord}(G_2,2)|+1}{2} + \frac{(n_1-1) \cdot n/n_1}{2} = m_0.$$
We can also see that $$2 \hat{\;} A \subseteq (\mathbb{Z}_{n_1} \times G_2) \setminus (\{0\} \times X),$$ and hence $$|2 \hat{\;} A| \leq n-|X| \leq n-2,$$ as claimed.  

Suppose now that $n_1$ is even.  In this case, we can easily see that
$$|\mathrm{Ord}(G,2)|=2 \cdot |\mathrm{Ord}(G_2,2)|+1,$$ so $$m_0=\frac{n+2 \cdot|\mathrm{Ord}(G_2,2)|+2}{2}=\frac{n}{2} + |\mathrm{Ord}(G_2,2)|+1.$$
This time, we set
$$A=\left( \{0, n_1/2\} \times A_2 \right) \cup \left( \{1,2,\dots,n_1/2-1\} \times G_2 \right).$$
We again have
$$|A|=2 \cdot |A_2| + (n_1/2-1) \cdot |G_2| = (n/n_1+|\mathrm{Ord}(G_2,2)|+1) + (n/2-n/n_1) = m_0$$
and $$2 \hat{\;} A \subseteq (\mathbb{Z}_{n_1} \times G_2) \setminus (\{0\} \times X),$$ and thus $$|2 \hat{\;} A| \leq n-|X| \leq n-2,$$ which completes our proof.  $\Box$

\addcontentsline{toc}{section}{Proof of Theorem \ref{rho hat m=n}}

\section*{Proof of Theorem \ref{rho hat m=n}} \label{proof of rho hat m=n}

Our claim is already established for cyclic groups by, for example, Proposition \ref{rhohat(G,m,h)>=mcyclic}.   

Turning to noncyclic groups, let us first consider the elementary abelian 2-group $\mathbb{Z}_2^r$, which we write as $G_1 \times \mathbb{Z}_2$ with $G_1=\mathbb{Z}_2^{r-1}$.  The result for $h=2$ has been delivered in Proposition \ref{rho hat Z_2^r}, so assume that $h \geq 3$.  Let $g_1 \in G_1$ and $g_2 \in \mathbb{Z}_2$ be arbitrary; we need to prove that there are $h$ pairwise distinct elements in $G_1 \times \mathbb{Z}_2$ that add to $(g_1,g_2)$.

Noting that we have $h \leq 2^{r-1}$, we choose $h-1$ arbitrary elements $a_1, \dots, a_{h-1}$ in $G_1$.  Set $$a=g_1-(a_1+\cdots+a_{h-1}).$$ (To make the proof more transparent, we used subtraction here, though of course in $G_1$ it is equivalent to addition.)  If $a$ is distinct from $a_i$ for all $1 \leq i \leq h-1$, then $$(a,g_2), (a_1,0), \dots, (a_{h-1},0)$$ are $h$ distinct elements of $G_1 \times \mathbb{Z}_2$ that add to $(g_1,g_2)$.

In the case when $a$ equals one of $a_1, \dots, a_{h-1}$, say $a=a_1$, then, noting also that $h-1 \geq 2$, we see that 
$$(a_1,0), (a_1,1), (a_2,g_2-1), (a_3,0), \dots, (a_{h-1},0)$$ are $h$ distinct elements of $G_1 \times \mathbb{Z}_2$ that add to $(g_1,g_2)$. 

Next, we consider the group $\mathbb{Z}_2 \times \mathbb{Z}_{\kappa}$ with $\kappa \geq 4$ even.  Let $g_1 \in \mathbb{Z}_2$ and $g_2 \in \mathbb{Z}_{\kappa}$; we need to find $h$ distinct elements of $\mathbb{Z}_2 \times \mathbb{Z}_{\kappa}$ that add to $(g_1,g_2)$.  If $h < \kappa$, then this follows easily from the result for cyclic groups.  Given that $h \leq n/2$, the remaining case is that $h=\kappa$, which needs more attention.  We separate two cases depending on $\kappa$ mod 4.

Suppose first that $\kappa$ is divisible by 4, and let $c=\kappa/4$. Since $$c \leq \kappa/2-1,$$  we have pairwise distinct integers $a_1, \dots, a_c$ that are all between 1 and $\kappa/2-1$, inclusive.  
We consider two subcases, as follows: if $g_2=0$, let $$A=\{(0,a_i),(0,-a_i),(1,a_i),(1,-a_i) \mid i=1,\dots, c-1\} \cup \{(1,a_c),(g_1,-a_c),(0,0),(1,0)\}.$$
If $g_2 \neq 0$, then one of the integers $a_1, \dots, a_c$ may equal $g_2$ or $-g_2$; we will assume that none of $a_1, \dots, a_{c-1}$ equals $\pm g_2$, and we let
$$A=\{(0,a_i),(0,-a_i),(1,a_i),(1,-a_i) \mid i=1,\dots, c-1\} \cup \{(0,a_c),(0,-a_c),(1,g_2),(g_1-1,0)\}.$$ 
In both subcases, $A$ consists of $h=\kappa$ pairwise distinct elements that add to $(g_1,g_2)$.

Now suppose that $\kappa \equiv 2$ mod 4, and let $c=(\kappa-2)/4$. This time, $$c \leq \kappa/2-2,$$  so we have pairwise distinct integers $a_1, \dots, a_{c+1}$ that are all between 1 and $\kappa/2-1$, inclusive.  
Again we consider two subcases: if $g_2=0$, let $$A=\{(0,a_i),(0,-a_i),(1,a_i),(1,-a_i) \mid i=1,\dots, c\} \cup \{(0,a_{c+1}),(g_1,-a_{c+1})\}.$$
If $g_2 \neq 0$, then one of the $c+1$ integers may equal $g_2$ or $-g_2$; we will assume that none of $a_1, \dots, a_c$ equals $\pm g_2$, and we let
$$A=\{(0,a_i),(0,-a_i),(1,a_i),(1,-a_i) \mid i=1,\dots, c\} \cup \{(0,0),(g_1,g_2)\}.$$ 
In both subcases, $A$ consists of $h=\kappa$ pairwise distinct elements that add to $(g_1,g_2)$.

This completes the cases when $G$ is cyclic, the elementary abelian 2-group, or is of the form $\mathbb{Z}_2 \times \mathbb{Z}_{\kappa}$, so we assume that $G \cong G_1 \times \mathbb{Z}_{\kappa}$, where $G_1$ is of order at least three, $\kappa$ is the exponent of $G$, and $\kappa \geq 3$.  We will also write
$$h=c \kappa+b$$ where $$0 \leq b \leq \kappa-1.$$  Furthermore, since $|G_1| \geq 3$ and $h \leq |G_1| \cdot \kappa/2$, we also have $$c \leq |G_1|-2.$$
Let $g_1 \in G_1$ and $g_2 \in \mathbb{Z}_{\kappa}$; we need to find $h$ distinct elements of $G_1 \times \mathbb{Z}_{\kappa}$ that add to $(g_1,g_2)$.  We consider two cases depending on whether $b$ is positive or not.

If $b \geq 1$, we let $A_1$ be any $c$ distinct elements of $G_1 \setminus \{0,g_1\}$ (possible since $c \leq |G_1|-2$), and let $A_2$ be any $b$ distinct elements of $\mathbb{Z}_{\kappa}$ that add to $$g_2-\kappa (\kappa-1) c/2$$ (possible since $1 \leq b \leq \kappa-1$).  Let $a$ be any element of $A_2$.  Then
$$(A_1 \times \mathbb{Z}_{\kappa}) \cup (\{0\} \times (A_2 \setminus \{a\})) \cup \{(g_1,a)\}$$ is a set of $h$ distinct elements of $G$, and its elements add to $(g_1,g_2)$.  (Note that, since $\kappa$ is the exponent of $G$, $\kappa a_1=0$ for any $a_1 \in A_1$.)

Suppose now that $b=0$; we need to separate several subcases.  If $g_1 \neq 0$, we let $A_1 \subset G_1$ to be any $c-1$ distinct elements of $G_1 \setminus \{0,g_1\}$, and let $A_2$ to be the $\kappa-1$ distinct elements of $\mathbb{Z}_{\kappa}$ that add to $$g_2-\kappa (\kappa-1) (c-1)/2.$$  Then
$$(A_1 \times \mathbb{Z}_{\kappa}) \cup (\{0\} \times A_2) \cup \{(g_1,0)\}$$ is a set of $h$ distinct elements of $G$ and its elements add to $(g_1,g_2)$. 

Next, suppose that $b=0$, $g_1=0$, and $\kappa$ is odd.  Let $a$ be any nonzero element of $G_1$, choose $A_1$ to be any $c-1$ distinct elements of $G_1 \setminus \{0,a,-a\}$ (note that this is possible even when $a \neq -a$, since $c-1 \leq |G_1|-3$); furthermore, let $$A_2=\{1, 2,3,\dots, (\kappa -1)/2\}.$$  Then
$$(A_1 \times \mathbb{Z}_{\kappa}) \cup (\{a\} \times A_2) \cup (\{-a\} \times (-A_2)) \cup \{(0,g_2) \}$$ is a set of $h$ distinct elements of $G$ and its elements add to $(0,g_2)$.  (Note that, since $\kappa-1$ is even, $\kappa (\kappa-1) (c-1)/2)$ equals zero in $\mathbb{Z}_{\kappa}$.) 

The case when $b=0$, $g_1=0$, $\kappa$ is even, and $g_2 \neq t$ where $$t=\kappa (\kappa-1) (c-1)/2$$ is very similar: again we let $a$ be any nonzero element of $G_1$, choose $A_1$ to be any $c-1$ distinct elements of $G_1 \setminus \{0,a,-a\}$, but now we set $$A_2=\{1,2,\dots, \kappa/2-1\}.$$  Then
$$(A_1 \times \mathbb{Z}_{\kappa}) \cup (\{0\} \times A_2) \cup (\{0\} \times (-A_2))\cup \{(a,0), (-a,g_2-t)\}$$ is a set of $h$ distinct elements of $G$ and its elements add to $(0,g_2)$.  This very construction works even if $g_2 =t$ as long as $a$ has order at least three in $G_1$.

This leaves us with the case when $b=0$, $g_1=0$, $\kappa$ is even, $g_2=t$, and $G_1$ is an elementary abelian 2-group.  Our construction is again similar, but we start with two distinct nonzero elements $a_1$ and $a_2$ of $G_1$.  We then let $A_1$ be any $c-1$ distinct elements of $G_1 \setminus \{0,a_1,a_2\}$, and set $$A_2=\{1,2,\dots, \kappa/2-2\}.$$  Then
$$(A_1 \times \mathbb{Z}_{\kappa}) \cup (\{0\} \times A_2) \cup (\{0\} \times (-A_2))\cup \{(a_1,0), (a_2,0), (a_1,\kappa/2), (a_2, \kappa/2\}$$ is a set of $h$ distinct elements of $G$ and its elements add to $(0,g_2)$.  $\Box$

\addcontentsline{toc}{section}{Proof of Proposition \ref{m=3 zero-free min size}}

\section*{Proof of Proposition \ref{m=3 zero-free min size}} \label{proof of m=3 zero-free min size}

Let $A=\{a,b,c\}$ be a weakly zero-sum-free 3-subset of $G$; we then have
$$\Sigma^*  A = \{a,b,c,a+b,a+c,b+c,a+b+c\}.$$

We can observe right away that the elements $a$, $b$, $c$, and $a+b+c$ must be pairwise distinct (otherwise one of $a+b$, $a+c$, or $b+c$ would be $0$).  Furthermore, none of $a+b$, $a+c$, or $b+c$ can equal $a+b+c$.  Therefore, the size of $\Sigma^*  A$ is 4, 5, 6, or 7, depending on how many of the equations $a+b=c$, $a+c=b$, and $b+c=a$ hold.  Note that if all three of them hold, then $a+b+c=0$.  Therefore, the size of $\Sigma^*  A$ is 5, 6, or 7.

Suppose first that $n$ is odd.  In this case, $|\Sigma^*  A| \geq 6$, since otherwise we would have, wlog, $a+b=c$ and $a+c=b$, which would imply that $2a=0$, and that is impossible if $n$ is odd.  On the other hand, if $c=a+b$, then clearly $|\Sigma^*  A| = 6$.

Now if $n \in \{1,3,5\}$, then $G$ has no weakly zero-sum-free 3-subsets.  In the case when $n \geq 7$ and odd, we show that we can always find elements $a \in G$ and $b \in G$ so that $A=\{a, b, a+b\}$ is weakly zero-sum-free.  If $G$ is isomorphic to $G_1 \times G_2$, then for any $g_1 \in G_1 \setminus \{0\}$ and  $g_2 \in G_2 \setminus \{0\}$, we can take $a=(g_1,0)$ and $b=(0,g_2)$, since then   
$$0 \not \in \Sigma^*  A = \{(g_1,0), (0,g_2), (g_1,g_2), (2g_1,g_2), (g_1,2g_2),(2g_1,2g_2)\}.$$
This leaves us with cyclic groups of (prime) order at least 7, in which case we take $a=1$ and $b=2$, for which 
$$0 \not \in \Sigma^*  A = \{1,2,3,4,5,6\}.$$

Suppose next that $n$ is even and that the exponent $\kappa$ of $G$ is at least 5.  In this case, let $a$ be an order 2 element of $G$ and $b$ to be an order $\kappa$ element of $G$; we also set $c=a+b$.  In this case,   
$$\Sigma^*  A = \{a,b,a+b,2b, a+2b\}.$$  None of these elements are $0$; for example, $a+2b=0$ would imply that $4b=-2a=0$, contradicting the fact that $b$ has order at least 5.

Suppose now that $G$ has exponent 4 and order at least 8.  In this case, we can write $G$ as $G_1 \times \mathbb{Z}_4$; let $g_1$ be an element of $G_1$ of order 2.  Choosing $a=(g_1,0)$, $b=(0,1)$ and $c=(g_1,1)$ gives us
$$\Sigma^*  A = \{(g_1,0),(0,1),(g_1,1),(0,2),(g_1,2)\},$$ so $A$ is weakly zero-sum-free.

We are left with the group $\mathbb{Z}_2^r$, for which we prove that $\Sigma^*  A$ has size 7 for every weakly zero-sum-free subset $A=\{a,b,c\}$ of size 3.  Indeed, none of the equations $a+b=c$, $a+c=b$, and $b+c=a$ hold: for example, $a+b=c$ implies that $a+b+c=2c=0$.  Finally, observe that when $r \geq 3$, then letting $a$, $b$, and $c$ denote the elements of $\mathbb{Z}_2^r$ that consist of $r-1$ zero components, with the 1 component being in three different places, the set $A=\{a,,c\}$ is weakly zero-sum-free.  This completes our proof.  $\Box$

\addcontentsline{toc}{section}{Proof of Theorem \ref{From Klopsch Lev}}

\section*{Proof of Theorem \ref{From Klopsch Lev}} \label{proof of From Klopsch Lev}

We will use the notations of \cite{KloLev:2009a}.  First, note that
$$\widehat{\chi} (G,[0,s]) = s^+_{s+1}(G)+1.$$  Therefore, we need to prove that $$s^+_{\rho}(G) \leq \widehat{v}(n, \rho -1)= \max \left \{\left( \left \lfloor \frac{d-2}{\rho-1} \right \rfloor +1 \right) \cdot \frac{n}{d} \; \mid \; d \in D(n), d \geq \rho+1\right\}$$ for every integer $\rho \geq 2$. 

By Lemma 2.3 in \cite{KloLev:2009a}, we have
$$s^+_{\rho}(G) = \max \{ |H| \cdot t^+_{\rho} (G/H) \; \mid \; H \leq G, H \neq G\}.$$
Let $H  \leq G, H \neq G$ be such that $$s^+_{\rho}(G) = |H| \cdot t^+_{\rho} (G/H),$$ and suppose that $|G/H|=d$
and that $G/H$ has invariant factorization $$G/H \cong \mathbb{Z}_{d_1} \times \cdots \times \mathbb{Z}_{d_r}.$$ 

By Theorem 2.1 in \cite{KloLev:2009a}, we have
$$\mathrm{diam}^+(G/H) = \Sigma_{i=1}^r (d_i-1).$$  Since the number of elements in $\mathbb{Z}_{d_1} \times \cdots \times \mathbb{Z}_{d_r}$ with exactly one nonzero coordinate equals  $\Sigma_{i=1}^r (d_i-1)$ and this number is clearly at most the number of nonzero elements in $G/H$, we have 
$$\mathrm{diam}^+(G/H) \leq |G/H|-1=d-1.$$
In particular, when $d \leq \rho$, we have $\mathrm{diam}^+(G/H) < \rho$.  
According to page 27 in \cite{KloLev:2009a}, we then have $t^+_{\rho} (G/H)=0$, so $s^+_{\rho}(G)=0$, from which our claim trivially follows. 

Assume now that $d \geq \rho +1$.  By Proposition 2.8 in \cite{KloLev:2009a}, we have
$$t^+_{\rho} (G/H) \leq \left \lfloor \frac{d-2}{\rho-1} \right \rfloor +1,$$ and so
\begin{eqnarray*}
s^+_{\rho}(G) &\leq&  \left( \left \lfloor \frac{d-2}{\rho-1} \right \rfloor +1 \right) \cdot \frac{n}{d} \\ \\
& \leq& \widehat{v}(n, \rho -1),
\end{eqnarray*} as claimed. $\Box$

\addcontentsline{toc}{section}{Proof of Proposition \ref{fixed h crit nonzero}}

\section*{Proof of Proposition \ref{fixed h crit nonzero}} \label{proof of fixed h crit nonzero}

Clearly, $$\chi \hat{\;} (G^*,h) \leq {\chi} \hat{\;} (G,h),$$ since if there is a subset $A$ of $G \setminus \{0\}$ with $h \hat{\;} A \neq G$, then $|A|+1$ provides a lower bound for ${\chi} \hat{\;} (G,h)$.

For the other direction, let $B$ be a subset of $G$ of size ${\chi} \hat{\;} (G,h)-1$ for which $h \hat{\;} B \neq G$.  Since ${\chi} \hat{\;} (G,h) \leq n$, we have $|B| \leq n-1$, and so $|-B| \leq n-1$ as well.  Let $g \in G \setminus (-B)$.  Then $A=g+B$ has size ${\chi} \hat{\;} (G,h)-1$, and $A \subseteq G \setminus \{0\}$, since $0 \in A=g+B$ would imply that $0=g+b$ for some $b \in B$, and thus $g=-b \in -B$, a contradiction.  But $h \hat{\;} A$ and $h \hat{\;} B$ have the same size, so we conclude that $h \hat{\;} A \neq G$, from which 
$$\chi \hat{\;} (G^*,h) \geq {\chi} \hat{\;} (G,h)$$ follows. $\Box$

\addcontentsline{toc}{section}{Proof of Lemma \ref{lemma with prime sqrt}}

\section*{Proof of Lemma \ref{lemma with prime sqrt}} \label{proof of lemma with prime sqrt}

Suppose first that $k$ is odd, in which case $h=(k+1)/2$, and our desired inequality becomes
$$\lfloor (2n-4)/(k+1) \rfloor +(k+1)/2 \leq k$$ or
$$\lfloor (2n-4)/(k+1) \rfloor  \leq (k-1)/2.$$
Since $(k-1)/2$ is an integer, it suffices to prove that
$$(2n-4)/(k+1)  < (k-1)/2+1$$ or, equivalently, that
$$4n-8 < (k+1)^2.$$
But $$(k+1)^2 = (\lfloor 2 \sqrt{n-2} \rfloor +1)^2> (2 \sqrt{n-2} )^2=4n-8,$$ as claimed.

Assume now that $k$ is even, in which case $h=k/2$, and our desired inequality becomes
$$\lfloor (2n-4)/k \rfloor +k/2 \leq k$$ or
$$\lfloor (2n-4)/k \rfloor  \leq k/2.$$
Since $k/2$ is an integer, it suffices to prove that
$$(2n-4)/k  < k/2+1$$ or, equivalently, that
$$4n-7 < (k+1)^2.$$
As above, we see that $$(k+1)^2 > 4n-8,$$  so we just need to rule out the possibility that $$(k+1)^2=4n-7.$$  This is indeed impossible, as the left-hand side is the square of an odd integer and thus congruent to 1 mod 8, while the right-hand side, since $n$ is odd, is congruent to 5 mod 8.  $\Box$

\addcontentsline{toc}{section}{Proof of Theorem \ref{four inverse probs}}

\section*{Proof of Theorem \ref{four inverse probs}}  \label{proof of four inverse probs}

Clearly, any P1-set is a P3-set, and any P2-set is a P4-set.

 Suppose now that $A$ is a P3-set; we will prove that $0 \in A$ and $A \setminus \{0\}$ is a P2-set, as follows.  If we were to have $0 \not \in A$, then, since by Theorem \ref{four quantities} $A$ has size 
$$|A|={\chi} \hat{\;} (G,\mathbb{N})-1={\chi} \hat{\;} (G^*,\mathbb{N}),$$  we have $\Sigma^* A=G$, but that contradicts our assumption that $A$ is a P3-set.  Therefore, $0 \in A$, and, using Theorem \ref{four quantities} again,
$$|A \setminus \{0\}|=|A|-1={\chi} \hat{\;} (G,\mathbb{N})-2={\chi} \hat{\;} (G^*,\mathbb{N}_0)-1.$$ Furthermore, since $0 \in A$, we have $\Sigma^* A =\Sigma A$, but then $\Sigma (A \setminus \{0\}) \neq G$, since $\Sigma^* A \neq G$.  This proves that $A \setminus \{0\}$ is a P2-set. 

Finally, Theorem \ref{four quantities} immediately implies that if $0 \in A$ and $A \setminus \{0\}$ is a P2-set, then $A$ is a P1-set.  This completes our proof. 
$\Box$

\addcontentsline{toc}{section}{Proof of Theorem \ref{size n/2 not spanning}}

\section*{Proof of Theorem \ref{size n/2 not spanning}} \label{proof of size n/2 not spanning}   

Our proof will follow that of Theorem 3.1 in \cite{GaoHamLlaSer:2003a} for the case of $|S|=n/2-1$.  We will need two lemmas, also found (though one stated with a crucial condition missing) in \cite{GaoHamLlaSer:2003a}:

{\bf Lemma 1 (Gao, Hamidoune, Llad\'o, and Serra; cf.~\cite{GaoHamLlaSer:2003a} Lemma 2.7)} \label{lemma 1 for gao et al}\index{Hamidoune, Y. O.}\index{Gao, W.}\index{Llad\'o, A. S.}\index{Serra, O.}  
{\em Let $S \subseteq G \setminus \{0\}$, $\langle S \rangle = G$, and $|S| \geq 14$.  Suppose further that there is no proper subgroup $H$ of $G$ for which $$|S \cap H| \geq |S|-1.$$  Then
$$|\Sigma S | \geq \min \{n-3, 3|S|-3\}.$$}

{\bf Lemma 2 (Gao, Hamidoune, Llad\'o, and Serra;\index{Hamidoune, Y. O.}\index{Gao, W.}\index{Llad\'o, A. S.}\index{Serra, O.}  cf.~\cite{GaoHamLlaSer:2003a} Lemma 2.9)} \label{lemma 2 for gao et al}
{\em Let $S \subseteq G \setminus \{0\}$.  Suppose further that $H$ is a subgroup of $G$ of prime index $p$, for which $$|S \cap H| \geq p-1.$$  Then
$$\Sigma (S \setminus H) +H=G.$$}

{\em Proof of Theorem \ref{size n/2 not spanning}}:  If $A$ is a subgroup of order $n/2$, then clearly $\Sigma A \neq G$, so we only need to prove the converse.  We can check the claim for all groups of order $16$ using the computer program \cite{Ili:2017a}, so we will assume that $n \geq 18$.

Let us assume, indirectly, that there is a subset $A$ of $G$ of size $|A|=n/2$ that is not a subgroup of $G$ but for which $\Sigma A \neq G$.  Since $A$ is not a subgroup of $G$, $ A  \subset \langle A \rangle$, which implies that $\langle A \rangle=G$.

Observe that by Theorem \ref{Diderrich and Mann}, if $|A|=n/2$ and $\Sigma A \neq G$, then $0 \in A$.  Let $S=A \setminus \{0\}$.  Then  
$$\langle S \rangle = \langle S \cup \{0\} \rangle = \langle A \rangle =G.$$  Since $\Sigma A \neq G$ implies that $\Sigma S \neq G$, we may choose an element $g \in G \setminus \Sigma S$.  

{\bf Claim}: For every proper subgroup $H$ of $G$ we have
$$|S \cap H| \leq n/4.$$

{\em Proof of Claim}:  Suppose that our claim is false; we then have a proper subgroup $H$ of $G$ for which  
$$|S \cap H| \geq \lceil (n+2)/4 \rceil.$$  The index of $H$ in $G$ then is either 2 or 3.  Since $\langle S \rangle =G$, $S$ cannot be a subset of $H$; furthermore, if the index of $H$ in $G$ is 3, then $S \setminus H$ cannot consist of a single element.  Therefore, $H$ is a subgroup of $G$ of order $p \in \{2,3\}$ so that
$$|S \setminus H| \geq p-1.$$  Thus, by Lemma 2 above, $$\Sigma (S \setminus H) +H=G.$$ But then $$\Sigma (S \cap H) \neq H,$$ otherwise we would have
$$G=\Sigma (S \setminus H) + \Sigma (S \cap H) = \Sigma S ,$$ contradicting our assumption.  We fix an element $h \in H \setminus \Sigma (S \cap H).$ 

Since for $n \geq 16$, $|S \cap H| \geq \lceil (n+2)/4 \rceil \geq 5$, we see that we can find an element $s \in S \cap H$ so that $S_1=(S \cap H) \setminus \{s\}$ either has size at least 5, or has size 4 but is not of the form $\{\pm a, \pm 2a\}$ for any $a \in S$.  

Now $$|S_1 \cup \{0\}|=|S \cap H|  > n/4  \geq |H|/2,$$ so $$\langle S_1 \rangle = \langle S_1 \cup \{0\} \rangle = H,$$ and, 
therefore, by Theorem \ref{Ham min sumset size},   
$$|\Sigma S_1| \geq \min \{|H|-1, 2|S_1|\}.$$
Here $$|H|-1 \leq n/2-1 = 2 \cdot ((n+2)/4-1) \leq 2 \cdot (|S \cap H| -1) = 2 \cdot |S_1|,$$ so
$$|\Sigma S_1| \geq |H|-1.$$
This is impossible, however, since we know of at least two distinct elements of $H$ that are not in $\Sigma S_1$: $h$ and $h-s$.  Indeed, $h \in H \setminus \Sigma (S \cap H),$ so $h \in H \setminus \Sigma S_1;$  for the same reason, if $h-s \in \Sigma S_1$, then $h=s+(h-s) \in \Sigma (S \cap H)$, a contradiction.  This proves our Claim.

Since $|S|=n/2-1>2$, we can select distinct elements $s_1, s_2 \in S$ so that $s_1+s_2 \neq 0$.  Set $S_2=S \setminus \{s_1,s_2\}$.  We find that $\langle S_2 \rangle=G$, otherwise $\langle S_2 \rangle$ would be a proper subgroup of $G$ for which 
$$|S \cap \langle S_2 \rangle | \geq |S_2| =n/2-3 >n/4,$$ contradicting our Claim above.  By the same Claim, when $n \geq 18$, then for each proper subgroup $H$ of $G$, we have
$$|S_2 \cap H| \leq |S \cap H| \leq \lfloor n/4 \rfloor \leq n/2-5  = |S_2|-2.$$ 
We can then apply Lemma 1 above, and get
$$|\Sigma S_2| \geq \min \{n-3, 3|S_2|-3\}=\min\{n-3, 3(n/2-3)-3\}=n-3.$$  This is impossible; recall that $g \not \in \Sigma S$, so none of the four pairwise distinct elements $g, g-s_1, g-s_2,$ or $g-s_1-s_2$ are in $\Sigma S_2$.
$\Box$

\addcontentsline{toc}{section}{Proof of Theorem \ref{size n/3 not spanning}}

\section*{Proof of Theorem \ref{size n/3 not spanning}} \label{proof of size n/3 not spanning}   

Our proof will follow that of Theorem 3.2 in \cite{GaoHamLlaSer:2003a} for the case of $|S|=n/3+1$.  Besides the two lemmas in the proof of Theorem \ref{size n/2 not spanning} above, we will also need the following from \cite{GaoHamLlaSer:2003a}:

{\bf Lemma (Gao, Hamidoune, Llad\'o, and Serra; cf.~\cite{GaoHamLlaSer:2003a} Lemma 2.8)}\index{Hamidoune, Y. O.}\index{Gao, W.}\index{Llad\'o, A. S.}\index{Serra, O.} 
{\em Let $X$ be a generating subset of $G$ that is asymmetric (i.e.~$X \cap -X=\emptyset$), and suppose that $| \Sigma X| \leq n/2$.  Then there is a proper subset $V$ of $X$ so that
$$|\Sigma X | \geq 4|V|+\frac{(|X|+|V|+5)(|X|-|V|-1)-2}{4}.$$}

{\em Proof of Theorem \ref{size n/3 not spanning}}:  If $A$ is of the form specified, then clearly $\Sigma A \neq G$, so we only need to prove the converse.  Since 27 and 45 are the only odd values of $n$ up to 62 for which $n/3$ is a composite integer, we may assume that $n \geq 63$.  

Let us assume that $A$ is a subset of $G$ so that $|A|=n/3+1$ and $\Sigma A \neq G$.  Then $\langle A \rangle=G$; furthermore, by Theorem \ref{Gao and Hamidoune}, we see that $0 \in A$.  Let $S=A \setminus \{0\}$; then $S$ has size $n/3$.  

{\bf Claim 1}: There is an asymmetric subset $X$ of $S$ of size $(n-3)/6$ so that $|\Sigma X| \leq (n-1)/2$.    

{\em Proof of Claim 1}: Let $s$ be an arbitrary element of $S$, and consider $S \setminus \{s\}$.  Observe that $S \setminus \{s\}$ can be partitioned as $X_1 \cup X_2$ so that $X_1$ and $X_2$ are both asymmetric and have size $(n-3)/6$.  We just have to prove that at least one of $\Sigma X_1$ or $\Sigma X_2$ has size at most $(n-1)/2$.  If this were not the case, then for any $g \in G \setminus \Sigma S$, we would have $$|\Sigma X_1|+|g-\Sigma  X_2|=|\Sigma  X_1|+| \Sigma  X_2| \geq n+1,$$ so $\Sigma X_1$ and $g-\Sigma X_2$ cannot be disjoint and thus $$g \in \Sigma (X_1 \cup X_2) \subseteq \Sigma S,$$ a contradiction.  This proves Claim 1.

Let $X$ be a set specified by Claim 1, and set $H=\langle X \rangle$.  

{\bf Claim 2}: $|H| = n/3$.  
 
{\em Proof of Claim 2}: Since $X$ is asymmetric, we must have $$|H| \geq 2|X|+1 = n/3;$$ since $n$ is odd, we only need to prove that $H \neq G$.

If $H=G$, then by the Lemma above, there is a proper subset $V$ of $X$ so that
$$|\Sigma X | \geq 4|V|+\frac{(|X|+|V|+5)(|X|-|V|-1)-2}{4}$$ and thus, by Claim 1,
$$\frac{n-1}{2} \geq 4|V|+\frac{((n-3)/6+|V|+5)((n-3)/6-|V|-1)-2}{4}.$$  
However, we find that, when $n \geq 63$, there is no value of $0 \leq |V| \leq (n-3)/6-1$ for which this inequality holds. (The minimum value of the right-hand side must occur either at $|V|=0$ or at $|V|=(n-3)/6-1$, but both possibilities yield values more than $(n-1)/2$.) 
This proves our claim.

{\em Proof of Theorem \ref{size n/3 not spanning}}:  Let $p$ be the smallest prime divisor of $n/3$.  Since 
$$|X|=\frac{n-3}{6} \geq \frac{n}{9}+2 = \frac{3n+9p^2-9p+(p-3)(n-9p)}{9p} \geq \frac{n/3}{p}+p-1,$$ by Theorems \ref{Mann and Wou}--\ref{Freeze et al}, we get $\Sigma X=H$.  Therefore,
$$\Sigma (S \setminus H) +H =  \Sigma (S \setminus \Sigma X) +\Sigma X \subseteq \Sigma (S \setminus X) +\Sigma X=\Sigma S \neq G,$$ so by Lemma 2 from the proof of Theorem \ref{size n/2 not spanning}, $|S \setminus H| \leq 1$.  But $|S|=|H|=n/3$ and $0 \in H \setminus S$, so $|S \setminus H|=1$, and our claim follows.  $\Box$

{\em Remark}:  We did not use in this proof that $n/3$ is composite, so as long as $n \geq 63$, our claim holds.

\addcontentsline{toc}{section}{Proof of Theorem \ref{thmz}}

\section*{Proof of Theorem \ref{thmz}} \label{proof of thm thmz}

If $n|h$, then $\tau (\mathbb{Z}_n,h)=0$ and, since for each $d \in D(n)$ we have $\gcd (d,h)=d$, 
$$v_h (n,h)= \max \left\{ \left ( \left \lfloor \frac{d-1- \gcd (d,h) }{h} \right \rfloor  +1 \right) \cdot \frac{n}{d} \mid d \in D(n) \right \}=0$$
as well.

Suppose now that $h$ is not divisible by $n$, and let $d \in D(n)$.  If $d|h$, then 
$$\left ( \left \lfloor \frac{d-1- \gcd (d,h) }{h} \right \rfloor  +1 \right) \cdot \frac{n}{d}=0;$$ otherwise (and there is at least one such $d$) $$d-1-\gcd(d,h) \geq 0.$$  Let $$c=\left \lfloor \frac{d-1- \gcd (d,h) }{h} \right \rfloor.$$
Note that
$$1 \leq \gcd(d,h) \leq \gcd (d,h) +hc \leq d-1.$$

Choose integers $a$ and $b$ for which
$$\gcd(d,h)=ha+db.$$

Let $H$ be the subgroup of $\mathbb{Z}_n$ of index $d$, and set
$$A=\bigcup_{i=a}^{a+c} (i+H);$$ then
$A$ has size $(c+1)n/d.$  We also see that $A$ is zero-$h$-sum-free in $\mathbb{Z}_n$, since
$$hA= \bigcup_{i=ha}^{ha+hc} (i+H) = \bigcup_{i=\gcd(d,h)-db}^{\gcd(d,h)-db+hc} (i+H)=\bigcup_{i=\gcd(d,h)}^{\gcd(d,h)+hc} (i+H) \subseteq \bigcup_{i=1}^{d-1} (i+H).$$

Thus we have constructed a zero-$h$-sum-free set of size $(c+1)n/d$ in $\mathbb{Z}_n$ for each $d \in D(n)$, which proves our lower bound.  $\Box$

\addcontentsline{toc}{section}{Proof of Proposition \ref{zer2hisBh}}

\section*{Proof of Proposition \ref{zer2hisBh}} \label{proofofzer2hisBh}

Suppose that $A=\{a_1,\dots,a_m\}$ is a $B_h$ set over $\mathbb{Z}$ in $G$.  Furthermore, suppose that $\lambda_1, \dots, \lambda_m$ are integers for which $$0=\lambda_1a_1+\cdots+\lambda_ma_m$$ and $$|\lambda_1|+\cdots+|\lambda_m| =2h.$$ 

Let $k$ be the smallest index for which 
$$h < |\lambda_1|+\cdots+|\lambda_k|;$$ w.l.o.g., we may assume that $\lambda_k >0$.  Let us write $\lambda=|\lambda_1|+\cdots+|\lambda_{k-1}|$ (if $k=1$, then simply let $\lambda=0$).

We may then write our original equation as
$$\lambda_1a_1+ \cdots +\lambda_{k-1}a_{k-1} + (h- \lambda)a_k=(h-\lambda-\lambda_k)a_k-\lambda_{k+1}a_{k+1}-\cdots-\lambda_ma_m$$ (with the understanding that if $k=1$ or $k=m$ then the terms before or after the one with $a_k$ in it vanish).

We first note that the left-hand side above is a signed sum of exactly $h$ terms, since 
$$|\lambda_1|+\cdots+|\lambda_{k-1}|+(h- \lambda)=h.$$
Before calculating the number of terms on the right-hand side, note that
$$h-\lambda-\lambda_k=h-\lambda-|\lambda_k|=h-(|\lambda_1|+\cdots+|\lambda_k|)<0,$$ so
the right-hand side consists of 
\begin{eqnarray*}
|h-\lambda-\lambda_k|+|\lambda_{k+1}|+\cdots+|\lambda_m| & = & |\lambda_1|+\cdots+|\lambda_k|-h+|\lambda_{k+1}|+\cdots+|\lambda_m| \\
& = & |\lambda_1|+\cdots+|\lambda_m|-h \\
& = & 2h - h = h
\end{eqnarray*} terms.  Thus, the equation 
$$\lambda_1a_1+ \cdots +\lambda_{k-1}a_{k-1} + (h- \lambda)a_k=(h-\lambda-\lambda_k)a_k-\lambda_{k+1}a_{k+1}-\cdots-\lambda_ma_m$$
has exactly $h$ terms on each side; since $A$ is a $B_h$ set over $\mathbb{Z}$, we get that $\lambda_i=0$ for each $i=1,\dots,m$, and therefore $A$ is a zero-$2h$-sum-free set over $\mathbb{Z}$ in $G$.
$\Box$

\addcontentsline{toc}{section}{Proof of Proposition \ref{Bhiszero4}}

\section*{Proof of Proposition \ref{Bhiszero4}} \label{proofofBhiszero4}

Let $A=\{a_1,\dots,a_m\}$ be a zero-$h$-sum-free set over $\mathbb{Z}$ in $G$ for some positive integer $h$ which is divisible by 4, and suppose that 
$$\lambda_1a_1+\cdots+\lambda_ma_m=\lambda_1'a_1+\cdots+\lambda_m'a_m,$$ where the coefficients are integers and $$|\lambda_1|+\cdots+|\lambda_m|=|\lambda_1'|+\cdots+|\lambda_m'|=2.$$  We will prove that we must have $\lambda_i=\lambda_i'$ for every $1 \leq i \leq m$.

It doesn't take long to verify, that the only possible values of
$$|\lambda_1-\lambda_1'|+\cdots+|\lambda_m-\lambda_m'|$$
are 0, 2, or 4.

If it is 0, then we are done.  If it is 2 or 4, then consider the sum  
$$\tfrac{h}{2}(\lambda_1-\lambda_1')a_1+\cdots+\tfrac{h}{2}(\lambda_m-\lambda_m')a_m$$ or
$$\tfrac{h}{4}(\lambda_1-\lambda_1')a_1+\cdots+\tfrac{h}{4}(\lambda_m-\lambda_m')a_m,$$ respectively.  Since the sum of the absolute values of the coefficients equals $h$ in both cases, the expressions cannot be zero, and therefore we cannot have
$$\lambda_1a_1+\cdots+\lambda_ma_m=\lambda_1'a_1+\cdots+\lambda_m'a_m,$$ which is a contradiction.  $\Box$

\addcontentsline{toc}{section}{Proof of Proposition \ref{greedy zero-even}}

\section*{Proof of Proposition \ref{greedy zero-even}} \label{proof of greedy zero-even}

First, we note that, for any $k \in \mathbb{N}$, the equation $k \cdot x=0$ has exactly $\gcd(k,n)$ solutions in $\mathbb{Z}_n$: indeed, the equation is equivalent to $kx$ being divisible by $n$ and thus to $(k / \gcd (k,n)) \cdot x$ being divisible by $n / \gcd (k,n)$,  which happens if, and only if, $x$ itself is divisible by $n / \gcd (k,n)$.  Consequently, the equation $k \cdot x=0$ has at most $k$ solutions in $\mathbb{Z}_n$.

We use the sequences $a(j,k)$ and $c(j,k)$ defined in Section \ref{0.2.3}.

Suppose that $$n > a(0,m-1)+a(1,m-1)+ \cdots + a(h-1,m-1).$$
We construct the set $A=\{a_1, \dots, a_m\}$ recursively, as follows.  Since $a(j,k)\geq 1$ for all $j,k$ implies $n>h$, we can find an element $a_1 \in \mathbb{Z}_n$ for which $h \cdot a_1 \neq 0$.  

Suppose now that we have already found the $(m-1)$-subset $A_{k-1}=\{a_1, \dots, a_{k-1}\}$ in $\mathbb{Z}_n$ so that $A_{k-1}$ is zero-$h$-sum-free over $\mathbb{Z}$ in $\mathbb{Z}_n$.  To find an element $a_k \in \mathbb{Z}_n$ so that $A_k=A_{k-1} \cup \{a_k\}$ is zero-$h$-sum-free over $\mathbb{Z}$ in $\mathbb{Z}_n$, we must have that none of 
$$h \cdot a_k, \; (h-1) \cdot a_k + g_1, \; (h-2) \cdot a_k +  g_2, \; \dots, \; 1 \cdot a_k +  g_{h-1}$$ equals zero for any $g_1 \in 1_{\pm}A_{k-1}, g_2 \in 2_{\pm}A_{k-1}, \dots, g_{h-1} \in (h-1)_{\pm}A_{k-1}$.   
This rules out at most
$$h+(h-1) \cdot |1_{\pm}A_{k-1}| + (h-2) \cdot |2_{\pm}A_{k-1}| + \cdots +  1\cdot |(h-1)_{\pm}A_{k-1}|$$ elements.  This quantity is at most
$$h \cdot c(0,k-1) + (h-1) \cdot c(1,k-1)+ (h-2) \cdot c(2,k-1) + \cdots + 1 \cdot c(h-1,k-1)$$ which, using Proposition \ref{functionsac}, can be rewritten as  
$$a(0,k-1)+a(1,k-1)+ \cdots + a(h-1,k-1).$$
Therefore, when $k \leq m$, we can find the desired element $a_k$.  

In particular, for $h=4$, it is sufficient that 
\begin{eqnarray*}
n & > &  a(0,m-1)+a(1,m-1)+a(2,m-1)+a(3,m-1) \\
& = & 4+12(m-1)+16 {m-1 \choose 2} + 8 {m-1 \choose 3} \\
& = & \frac{4}{3} m^3 + \frac{8}{3}m.
\end{eqnarray*}
$\Box$

\addcontentsline{toc}{section}{Proof of Proposition \ref{zetapmhodd}}

\section*{Proof of Proposition \ref{zetapmhodd}} \label{proofofzetapmhodd}

Let $A$ be the set of integers (viewed as elements of $G = \mathbb{Z}_n$) that are strictly between $\frac{h-1}{2}\frac{n}{h}$ and $\frac{h+1}{2}\frac{n}{h}$.  In particular, when $n$ is divisible by $h$, we take
$$A=\left\{ \frac{h-1}{2} \frac{n}{h} +1, \frac{h-1}{2} \frac{n}{h} +2, \dots, \frac{h+1}{2} \frac{n}{h} -1 \right\}.$$ When $n$ is not divisible by $h$, then neither are $\frac{h-1}{2}n$ and $\frac{h+1}{2}n$, since $h$ is relatively prime to $\frac{h-1}{2}$ and $\frac{h+1}{2}$; so, if $n$ is not divisible by $h$, we have  
$$A=\left\{ \left \lfloor \frac{h-1}{2} \frac{n}{h}\right \rfloor + 1, \left \lfloor \frac{h-1}{2} \frac{n}{h}\right \rfloor+2, \dots, \left \lfloor \frac{h+1}{2} \frac{n}{h} \right \rfloor \right\}.$$

We prove that $0 \not \in h_{\pm}A$ by showing that, for every $k=0,1,\dots,h$, the set $kA-(h-k)A$ does not contain 0.  Indeed, the smallest integer in $kA-(h-k)A$ is greater than $$k \cdot \frac{h-1}{2}\frac{n}{h} - (h-k) \cdot \frac{h+1}{2}\frac{n}{h}=\left(k-\frac{h+1}{2} \right)n$$ and less than $$k \cdot \frac{h+1}{2}\frac{n}{h} - (h-k) \cdot \frac{h-1}{2}\frac{n}{h}=\left(k-\frac{h-1}{2} \right)n,$$ or two consecutive multiples of $n$, and therefore $0 \not \in kA-(h-k)A$.

It remains to be verified that the size of $A$ equals $2  \left \lfloor \frac{n+h-2}{2h} \right \rfloor$.

This is easy to see if $n$ is divisible by $h$, since then the size of $A$ is clearly  
$$\frac{h+1}{2}\frac{n}{h}-\frac{h-1}{2}\frac{n}{h}-1=2 \left(  \frac{n+h}{2h}-1 \right)= 2 \left( \left \lfloor \frac{n+h}{2h} \right \rfloor -1 \right)=2  \left \lfloor \frac{n+h-2}{2h} \right \rfloor.$$

When $n$ is not divisible by $h$, then 
\begin{eqnarray*}
|A| & = & \left \lfloor \frac{h+1}{2} \frac{n}{h} \right \rfloor - \left \lfloor \frac{h-1}{2} \frac{n}{h}\right \rfloor \\
& = & \left( \frac{n+1}{2} + \left \lfloor \frac{n-h}{2h} \right \rfloor \right) - \left( \frac{n-1}{2} - \left \lfloor \frac{n-h}{2h} \right \rfloor -1 \right) \\
& = & 2 \left( \left \lfloor \frac{n-h}{2h} \right \rfloor+1 \right) \\ 
& = & 2 \left \lfloor \frac{n+h}{2h} \right \rfloor ;
\end{eqnarray*}
since neither $n+h$ nor $n+h-1$ is divisible by $2h$ (the first would imply that $h|n$ and the second that $n+h$ is odd), this equals $2  \left \lfloor \frac{n+h-2}{2h} \right \rfloor$, as claimed.  
$\Box$

\addcontentsline{toc}{section}{Proof of Proposition \ref{q}} \label{proofofq}

\section*{Proof of Proposition \ref{q}} \label{proofofq}

For statement 2 (a), we prove that for $t=2s$, the set $\{s,s+1\}$ is $2s$-independent in $\mathbb{Z}_n$ for all $n \geq 2s^2+2s+1$.  Suppose that for the integer $$a=\lambda_1 \cdot s + \lambda_2 \cdot (s+1)$$ we have $a=0$ in $\mathbb{Z}_n$ for some integer coefficients with $|\lambda_1|+|\lambda_2| \leq 2s$; w.l.o.g., assume that $\lambda_2 \geq 0$.  Then  $$-2s^2 \leq a \leq 2s^2+2s,$$ so $a=0$ in $\mathbb{Z}_n$ means that $a=0$ in $\mathbb{Z}$.  Therefore, $\lambda_2$ is divisible by $s$, hence $\lambda_2=0$,  $\lambda_2=s$, or $\lambda_2=2s$.  In the first and last cases, we have $\lambda_1=0$, and we are done.  The second case cannot happen, since $$\lambda_1 \cdot s + s \cdot (s+1)=0$$ implies that $\lambda_1=-(s+1)$, contradicting $|\lambda_1|+|\lambda_2| \leq 2s$.

Similarly, for statement 2 (b) (ii), we show (somewhat clarifying Miller's proof) that the set $\{s,s+1\}$ is $(2s-1)$-independent in $\mathbb{Z}_n$ for $n \geq 2s^2+s.$  Suppose that  for the integer $$a=\lambda_1 \cdot s + \lambda_2 \cdot (s+1)$$ we have $a=0$ in $\mathbb{Z}_n$ for some integer coefficients with $|\lambda_1|+|\lambda_2| \leq 2s-1$; w.l.o.g., assume that $\lambda_2 \geq 0$.  Then  $$-(2s^2-s) \leq a \leq 2s^2+s-1,$$ so $a=0$ in $\mathbb{Z}_n$ means that $a=0$ in $\mathbb{Z}$.  Therefore, $\lambda_2$ is divisible by $s$, hence $\lambda_2=0$ or  $\lambda_2=s$.  In the first case, we have $\lambda_1=0$, and we are done.  The second case cannot happen, since $$\lambda_1 \cdot s + s \cdot (s+1)=0$$ implies that $\lambda_1=-(s+1)$, contradicting $|\lambda_1|+|\lambda_2| \leq 2s-1$.

Next, we consider 2 (b) (iii), and prove that $\{1,2s-1\}$ is $(2s-1)$-independent in $\mathbb{Z}_n$ for $n=2s^2$.   Assume that  $$a=\lambda_1 \cdot 1 + \lambda_2 \cdot (2s-1)=0$$ in $\mathbb{Z}_n$ for some integer coefficients with $|\lambda_1|+|\lambda_2| \leq 2s-1$; w.l.o.g., assume that $\lambda_2 \geq 0$.   Then for the integer $a$ we have $$-(2s-1) \leq a \leq 4s^2-4s+1,$$ so $a=0$ in $\mathbb{Z}_n$ means that either $a=0$ in $\mathbb{Z}$ or $a=2s^2$ in $\mathbb{Z}$.  If $a=0$, then $\lambda_1$ is divisible by $2s-1$, hence $\lambda_1=0$ which implies that $\lambda_2=0$ and we are done, or $|\lambda_1|=2s-1$, which again implies that $\lambda_2=0$ which can only happen if $\lambda_1=0$ as well, a contradiction.

If $a=2s^2$ in $\mathbb{Z}$, then $\lambda_1-\lambda_2$ must be divisible by $2s$, which can only happen if it is 0, since 
$$|\lambda_1-\lambda_2| \leq |\lambda_1|+|\lambda_2| \leq 2s-1.$$  
But solving the system $$\lambda_1 \cdot 1 + \lambda_2 \cdot (2s-1)=2s^2$$ and $$\lambda_1-\lambda_2=0$$ yields $\lambda_1=\lambda_2=s$, contradicting $|\lambda_1|+|\lambda_2| \leq 2s-1$.

For statement 2 (b) (iv), we assume that $s$ and $n$ are both even, $n \geq 2s^2$, and prove that the set $\{s-1,s+1\}$ is $(2s-1)$-independent in $\mathbb{Z}_n$.  Suppose that $$a=\lambda_1 \cdot (s-1) + \lambda_2 \cdot (s+1)=0$$ in $\mathbb{Z}_n$ for some integer coefficients with $|\lambda_1|+|\lambda_2| \leq 2s-1$; w.l.o.g., assume that $\lambda_2 \geq 0$.

Note that every integer has the same parity as its absolute value. Therefore,
$$\lambda_1 \cdot (s-1) + \lambda_2 \cdot (s+1) \equiv |\lambda_1| \cdot (s-1) +|\lambda_2| \cdot (s+1) \equiv (|\lambda_1| +|\lambda_2|) \cdot (s+1) \; \mbox{mod} \; 2.$$   So if $|\lambda_1|+|\lambda_2|$ is odd and $s$ is even, then $a$ is odd and, since $n$ is even, $a \neq 0$ in $\mathbb{Z}_n$.  In particular, $0 \not \in (2s-1)_{\pm} \{s-1,s+1\}$.  It remains to be shown that $0 \not \in [1,2s-2]_{\pm} \{s-1,s+1\}$. 

For $|\lambda_1|+|\lambda_2| \leq 2s-2$,  we have $$-(2s^2-4s+2) \leq a \leq 2s^2-2,$$ so $a=0$ in $\mathbb{Z}_n$ means that $a=0$ in $\mathbb{Z}$.  Therefore, $\lambda_1$ is divisible by $s+1$ and $\lambda_2$ is divisible by $s-1$, hence either one of them equals 0, in which case they both do, or 
$$|\lambda_1|+|\lambda_2| \geq (s+1)+(s-1)=2s,$$ a contradiction.

For statement 2 (b) (v), we assume that $s$ is odd, $n$ is congruent to 2 mod 4, $n \geq 2s^2$, and prove that the set $\{s-2,s+2\}$ is $(2s-1)$-independent in $\mathbb{Z}_n$.  Suppose that $$a=\lambda_1 \cdot (s-2) + \lambda_2 \cdot (s+2)=0$$ in $\mathbb{Z}_n$ for some integer coefficients with $|\lambda_1|+|\lambda_2| \leq 2s-1$; w.l.o.g., assume that $\lambda_2 \geq 0$.

Like in the previous case, we can show that if $|\lambda_1|+|\lambda_2|$ is odd and $s$ is odd, then $a$ is odd and, since $n$ is even, $a \neq 0$ in $\mathbb{Z}_n$.  Therefore, $0 \not \in (2s-1)_{\pm} \{s-2,s+2\}$ and $0 \not \in (2s-3)_{\pm} \{s-2,s+2\}$.  We will now show that $0 \not \in (2s-2)_{\pm} \{s-2,s+2\}$.

Assume first that $\lambda_1 \geq 0$.  In this case,
$$0 < (2s-2)(s-2) \leq a \leq (2s-2)(s+2) <2n,$$ so $a =0$ in $\mathbb{Z}_n$ implies that $a=n$.  But then
$$n=\lambda_1 \cdot (s-2) + \lambda_2 \cdot (s+2) = (\lambda_1 + \lambda_2)  \cdot (s-2) + 4 \lambda_2 =(2s-2) \cdot (s-2) + 4 \lambda_2  \equiv 0 \; \mbox{mod} \; 4,$$ a contradiction, since $n$ is congruent to 2 mod 4. 

If $\lambda_1 < 0$, then 
$$-n< -(2s-2)(s-2) \leq a \leq (-1) (s-2) +(2s-3)(s+2) < n.$$  Therefore, $a=0$ in $\mathbb{Z}$, and since $s-2$ and $s+2$ are relatively prime, this implies that $\lambda_1$ is divisible by $s+2$ and $\lambda_2$ is divisible by $s-2$, hence either one of them equals 0, in which case they both do, or 
$$|\lambda_1|+|\lambda_2| \geq (s+2)+(s-2)=2s,$$ a contradiction.

The proof that $0 \not \in [1,2s-4]_{\pm} \{s-2,s+2\}$ is similar: 
for $|\lambda_1|+|\lambda_2| \leq 2s-4$,  we have $$-n< -(2s-4)(s-2) \leq a \leq (2s-4)(s+2) < n$$ so $a=0$ in $\mathbb{Z}$, and we can complete the proof as above.
$\Box$

\addcontentsline{toc}{section}{Proof of Proposition \ref{Zlower}}

\section*{Proof of Proposition \ref{Zlower}} \label{proofofZlower}

Let $n$ and $h$ be fixed positive integers; note that $$\overline{h}=\frac{h^2-h-2}{2}$$ is an integer.  
Let $d=\mathrm{gcd}(n,h)$, and set $r$ equal to the remainder of $\overline{h}$ when divided by $d$.  
Furthermore, let $$c=\left \lfloor \frac{n+h^2-r-2}{h} \right \rfloor.$$
Note that $h \leq n-1$, so $$c \leq (n+h^2-2)/h < (n+h^2-1)/h \leq (n+(h-1)n)/h = n;$$ we need to prove that $\tau\hat{\;}(\mathbb{Z}_n,h) \geq c$.

Since $d$ is a divisor of $\overline{h}-r$, we can find integers $a$ and $b$ for which $$a \cdot h - b \cdot n =-(\overline{h}-r);$$ furthermore, we may assume w.l.o.g. that $b \geq 0$.

Now define the set $A$ as $$A=\{a,a+1,\dots,a+c-1\}.$$  (Throughout this proof, we will talk about integers and elements of the group $\mathbb{Z}_n$ interchangeably: when an integer is considered as an element of $\mathbb{Z}_n$, we regard it mod $n$; on the other hand, an element $k$ of $\mathbb{Z}_n$ is treated as the integer $k$ in $\mathbb{Z}$.)  

Clearly, $|A|=c$.  We will show that $A$ is a weak zero-$h$-sum-free set in $\mathbb{Z}_n$ by showing that every element of $h \hat{\;} A$ is strictly between $bn$ and $(b+1)n$.

The smallest element of $h \hat{\;} A$ is $$ha+\frac{h(h-1)}{2}=h a+\overline{h}+1=b n-(\overline{h}-r)+\overline{h}+1=bn+1+r \geq bn+1.$$

For the largest element of $h \hat{\;} A$ we have
\begin{eqnarray*}
h a+h c-\frac{h(h+1)}{2} & = &  (h a -bn) +bn +hc - \frac{h(h+1)}{2} \\
& = & -(\overline{h}-r) +bn +hc - \frac{h(h+1)}{2} \\
& = & bn + hc -h^2+r+1 \\
& \leq & bn +(n+h^2-r-2) -h^2+r+1 \\
& = & (b+1)n-1,
\end{eqnarray*}
as claimed. $\Box$

\addcontentsline{toc}{section}{Proof of Proposition \ref{Zlower2}}

\section*{Proof of Proposition \ref{Zlower2}} \label{proofofZlower2}

Let $n$ and $h$ be fixed positive integers, and assume that $h^2|n$.  We let $\epsilon$ equal 0 if $h$ is even and 1 if $h$ is odd.  We define the set
$$A=\bigcup_{i=0}^{h-1} \left\{ (ih+\epsilon)\cdot \frac{n}{h^2} + j  \mid j=0,1,\dots, \frac{n}{h^2} \right\}.$$
We will prove that $A$ is zero-$h$-sum-free in $\mathbb{Z}_n$.

First, observe that every element of $h \hat{\;} A$ can be written as 
\begin{eqnarray*}
b & = & ((i_1+\cdots+i_h)h+\epsilon h) \cdot \frac{n}{h^2} +(j_1+\cdots+j_h) \\
& = & ((i_1+\cdots+i_h)+\epsilon ) \cdot \frac{n}{h} +(j_1+\cdots+j_h)
\end{eqnarray*}
 where $$\{i_1,\dots,i_h\} \subseteq \{0,1,\dots,h-1\}$$ and $$\{j_1,\dots,j_h\} \subseteq \left\{0,1,\dots,\frac{n}{h^2} \right\}.$$

If $b$ were to be divisible by $n$ (that is, $b=0$ in $\mathbb{Z}_n$), then it would need to be divisible by $\frac{n}{h}$.  Since the first term above is always divisible by $\frac{n}{h}$, the second term must be as well.  However, since we have $$0 \leq j_1+\cdots+j_h \leq h \cdot \frac{n}{h^2}=\frac{n}{h},$$ this can only happen if either $j_1=\cdots=j_h=0$ or $j_1=\cdots=j_h=\frac{n}{h^2}$.  In either case, all $j$ values are equal.

Therefore, since $b \in h \hat{\;} A$, all $i$ values must be distinct, and thus $$\{i_1,\dots,i_h\} = \{0,1,\dots,h-1\}.$$  Thus we must have
\begin{eqnarray*}
b & = & ((i_1+\cdots+i_h)h+\epsilon h) \cdot \frac{n}{h^2} +(j_1+\cdots+j_h) \\
& = & ((0+1+\cdots+(h-1))+\epsilon ) \cdot \frac{n}{h} +(j_1+\cdots+j_h) \\
& = & \left(\frac{h(h-1)}{2}+\epsilon \right) \cdot \frac{n}{h} +(j_1+\cdots+j_h).
\end{eqnarray*}
We will show that $b$ is not divisible by $n$.  We have four cases.

Case 1: $h$ is even and $j_1=\cdots=j_h=0$.  In this case we get $$b=\frac{h(h-1)}{2} \cdot \frac{n}{h} = \frac{h-1}{2} \cdot n,$$ and this is not divisible by $n$ since $\frac{h-1}{2}$ is not an integer.

Case 2: $h$ is odd and $j_1=\cdots=j_h=0$.  In this case we get $$b=\left(\frac{h(h-1)}{2}+1\right) \cdot \frac{n}{h} = \left(\frac{h-1}{2} +\frac{1}{h} \right) \cdot n,$$ and this is not divisible by $n$ since $\frac{h-1}{2}$ is an integer but  $\frac{1}{h}$ is not (if $h >1$).

Case 3: $h$ is even and $j_1=\cdots=j_h=\frac{n}{h^2}$.  In this case we get $$b=\frac{h(h-1)}{2} \cdot \frac{n}{h} + \frac{n}{h} = \left(\frac{h}{2} -\frac{h-2}{2h} \right) \cdot n,$$ and this is not divisible by $n$ since $\frac{h}{2}$ is an integer, but $\frac{h-2}{2h}$ is not an integer as it is strictly between 0 and $\frac{1}{2}$ (if $h>2$).

Case 4: $h$ is odd and $j_1=\cdots=j_h=\frac{n}{h^2}$.  In this case we get $$b=\left(\frac{h(h-1)}{2}+1\right) \cdot \frac{n}{h} + \frac{n}{h} = \left(\frac{h-1}{2} +\frac{2}{h} \right) \cdot n,$$ and this is not divisible by $n$ since $\frac{h-1}{2}$ is an integer but  $\frac{2}{h}$ is not (if $h >2$).
$\Box$

\addcontentsline{toc}{section}{Proof of Theorem \ref{Zfor n even h odd}}

\section*{Proof of Theorem \ref{Zfor n even h odd}} \label{proof of Zfor n even h odd}

By Proposition \ref{tau hat bounds} and Theorem \ref{Roth and Lempel even} we have 
$$\tau\hat{\;}(G,h) \leq \chi \hat{\;}(G,h) -1 = \left \{
\begin{array}{cl}
n-1 & \mbox{if} \; h=1; \\ \\
n/2 & \mbox{if} \; 3 \leq h \leq n/2-2; \\ \\
n/2+1 & \mbox{if} \; h=n/2-1; \\ \\
h+1 & \mbox{if} \; n/2 \leq h \leq  n-2; \\ \\
n-1 & \mbox{if} \; h=n-1.
\end{array}
\right.$$ 
Therefore, it suffices to prove that a weakly zero-$h$-sum-free set of the desired size exists.

For $h=1$ we take, of course, all nonzero elements of $\mathbb{Z}_n$, and for $3 \leq h \leq n/2-2$ (indeed, for any odd $h$) we take all odd elements of $\mathbb{Z}_n$.  

When $h=n/2-1$, then, since $h$ is odd, $\gcd(n,h)=1$, so by Corollary \ref{cor Zlower 1}, we get 
$$\tau\hat{\;}(\mathbb{Z}_n,h) \geq  \left\lfloor \frac{n+h^2- 2}{h} \right\rfloor =\left\lfloor \frac{2h+2+h^2- 2}{h} \right\rfloor =h+2=\frac{n}{2}+1.$$

The result for $n/2 \leq h \leq  n-2$ follows from Proposition \ref{Zlower11}.

Finally, when $h=n-1$, we again have $\gcd(n,h)=1$, with which Corollary \ref{cor Zlower 1} yields
$$\tau\hat{\;}(\mathbb{Z}_n,h) \geq  \left\lfloor \frac{n+h^2- 2}{h} \right\rfloor =\left\lfloor \frac{h+1+h^2- 2}{h} \right\rfloor =h=n-1,$$ which completes our proof.  $\Box$ 

We should also note that some of our cases carry through for even $h$ as well.

\addcontentsline{toc}{section}{Proof of Proposition \ref{Yin Z2}}

\section*{Proof of Proposition \ref{Yin Z2}} \label{proof of Yin Z2}

Since no two distinct elements of $\mathbb{Z}_2^r$ add to zero, $\{1\} \times \mathbb{Z}_2^{r-1}$ is a weak zero-$[1,3]$-free set in $\mathbb{Z}_2^r$, and thus
$$\tau \hat{\;} (\mathbb{Z}_2^r,[1,3]) \geq 2^{r-1}.$$

To prove the reverse inequality, suppose that$A$ is a weak zero-$[1,3]$-free set in $\mathbb{Z}_2^r$.  Let $a \in A$ be arbitrary, and let $$B=a+(A \setminus \{a\}).$$    Then $|B|=|A|-1$; furthermore, $0 \not \in A$, and $0 \not \in B$.

Furthermore, $A$ and $B$ are disjoint, since if we had an $a' \in A \setminus \{a\}$ for which $a+a'=a''$ for some $a'' \in A$, then $a+a'+a''$ would be $0$, which can only happen if $a'' \in \{a,a'\}$; however, this cannot be, as $0 \not \in A$.

Therefore, $A$ and $B$ are disjoint subsets of $\mathbb{Z}_2^r \setminus \{0\}$, from which $$|A|+|B|=|A|+|A|-1 \leq 2^r-1,$$ and thus $|A| \leq 2^{r-1}$, as claimed. $\Box$

\addcontentsline{toc}{section}{Proof of Proposition \ref{collins}}

\section*{Proof of Proposition \ref{collins}} \label{proof of collins}

Note that
$$\left \lfloor \frac{2n+h^2-3}{2h} \right \rfloor \leq \left \lfloor \frac{2n+h^2}{2h} \right \rfloor \leq \left \lfloor \frac{2n+hn}{2h} \right \rfloor \leq \left \lfloor \frac{2hn}{2h} \right \rfloor = n.$$

Let $$2n+h^2-3=4hq+r$$ with $$0 \leq r < 4h.$$  Then 
$$\left \lfloor \frac{2n+h^2-3}{2h} \right \rfloor =
\left\{
\begin{array}{cl}
2q & \mbox{if $0 \leq r < 2h$}, \\ \\
2q +1 & \mbox{if $2h \leq r < 4h$}.
\end{array}
\right.$$

We treat the cases of $0 \leq r < 2h$ and $2h \leq r < 4h$ separately.

Suppose first that $0 \leq r < 2h$; set $$A=\{ (n+1)/2-q, (n+1)/2-q+1, \dots, (n+1)/2+q-1\}.$$
As we showed above, $2q \leq n$, so we have $|A|=2q$.

For each integer $0 \leq k \leq h$, let us denote by $h(k) \hat{_{\pm}} A$ the collection of all signed sums of $h$ distinct terms of $A$ in which exactly $k$ terms are added and $h-k$ are subtracted.  (For example, $h(h) \hat{\;} A = h \hat{\;} A$, and $h(0) \hat{_{\pm}} A = -h \hat{\;} A$.)
Considering the elements of $A$ as integers (rather than elements of $\mathbb{Z}_n$), we will compute $a_{\min}(k)$ and $a_{\max}(k)$, the smallest and largest elements of $h(k) \hat{_{\pm}} A$, respectively.  

For $a_{\min}(k)$ we get
\begin{eqnarray*}
a_{\min}(k) & = & \left( k \cdot ((n+1)/2-q) + k(k-1)/2 \right) \\
& & - \left( (h-k) \cdot ((n+1)/2+q) - (h-k)(h-k+1)/2 \right) \\
& = & k(n+1)-h(n+1)/2-hq+\tfrac{k(k-1)+(h-k)(h-k+1)}{2} \\
& \geq & k(n+1)-h(n+1)/2-\tfrac{2n+h^2-3}{4}+\tfrac{k(k-1)+(h-k)(h-k+1)}{2}\\
& = & (k-(h+1)/2)n + ((h-2k)^2+3)/4 \\
& > & (k-(h+1)/2)n.
\end{eqnarray*}

Similarly, 
\begin{eqnarray*}
a_{\max}(k) & = & \left( k \cdot ((n+1)/2+q) - k(k+1)/2 \right) \\
& & - \left( (h-k) \cdot ((n+1)/2-q) +(h-k)(h-k-1)/2 \right) \\
& = & k(n+1)-h(n+1)/2+hq-\tfrac{k(k+1)+(h-k)(h-k-1)}{2} \\
& \leq & k(n+1)-h(n+1)/2+\tfrac{2n+h^2-3}{4}-\tfrac{k(k+1)+(h-k)(h-k-1)}{2}\\
& = & (k-(h+1)/2+1)n - ((h-2k)^2+3)/4 \\
& < & (k-(h+1)/2+1)n.
\end{eqnarray*}

Therefore, all elements of $h(k) \hat{_{\pm}} A$ lie strictly between two consecutive multiples of $n$, thus $0 \not \in h(k) \hat{_{\pm}} A$.  Since this holds for every $0 \leq k \leq h$ and $$h \hat{_{\pm}} A = \cup_{k=0}^h h(k) \hat{_{\pm}} A,$$ we have $0 \not \in h \hat{_{\pm}} A$.

Suppose now that $2h \leq r < 4h$; set $$B=\{ (n+1)/2-q, (n+1)/2-q+1, \dots, (n+1)/2+q\}.$$
As we showed above, $2q+1 \leq n$, so we have $|B|=2q+1$.

Since $B=A \cup \{(n+1)/2+q\}$, we have
\begin{eqnarray*}
b_{\min}(k) & = & a_{\max}(k)-k \\
& = & \left( k \cdot ((n+1)/2-q) + k(k-1)/2 \right) \\
& & - \left( (h-k) \cdot ((n+1)/2+q) - (h-k)(h-k+1)/2 \right) -k \\
& = & k(n+1)-h(n+1)/2-hq+\tfrac{k(k-1)+(h-k)(h-k+1)}{2} -k \\
& \geq & k(n+1)-h(n+1)/2-\tfrac{2n+h^2-3-2h}{4}+\tfrac{k(k-1)+(h-k)(h-k+1)}{2} -k \\
& = & (k-(h+1)/2)n + ((h-2k+1)^2+2)/4 \\
& > & (k-(h+1)/2)n.
\end{eqnarray*}

Similarly, 
\begin{eqnarray*}
a_{\max}(k) & = & a_{\max}(k)+k \\ 
& = & \left( k \cdot ((n+1)/2+q) - k(k+1)/2 \right) \\
& & - \left( (h-k) \cdot ((n+1)/2-q) +(h-k)(h-k-1)/2 \right) +k \\
& = & k(n+1)-h(n+1)/2+hq-\tfrac{k(k+1)+(h-k)(h-k-1)}{2} +k \\
& \leq & k(n+1)-h(n+1)/2+\tfrac{2n+h^2-3-2h}{4}-\tfrac{k(k+1)+(h-k)(h-k-1)}{2}+k \\
& = & (k-(h+1)/2+1)n - ((h-2k+1)^2+2)/4 \\
& < & (k-(h+1)/2+1)n.
\end{eqnarray*}

Again, all elements of $h(k) \hat{_{\pm}} B$ lie strictly between two consecutive multiples of $n$, thus $0 \not \in h(k) \hat{_{\pm}} B$, and, therefore, $0 \not \in h \hat{_{\pm}} B$.  This completes our proof.  $\Box$

\addcontentsline{toc}{section}{Proof of Proposition \ref{prob tau hat pm G, 1,t n small}}

\section*{Proof of Proposition \ref{prob tau hat pm G, 1,t n small}} \label{proof of prob tau hat pm G, 1,t n small}

First, we prove that if $A=\{a_1,\dots,a_t\}$ is a $t$-subset of $\mathbb{Z}_n$, then $0 \in [1,t] \hat{_{\pm}}  A$.  Observe that since $|\Sigma A| \leq n < 2^t$, we must have $$\epsilon_1, \dots, \epsilon_t, \delta_1, \dots, \delta_t \in \{0,1\}$$ for which 
$$\epsilon_1a_1+ \cdots+ \epsilon_ta_t= \delta_1a_1+ \cdots+ \delta_ta_t,$$ but for at least one index $1 \leq i \leq t$, $\epsilon_i \neq \delta_i$.  But then $$1 \leq |\epsilon_1-\delta_1|+ \cdots+ |\epsilon_t-\delta_t| \leq t,$$ and $$0=(\epsilon_1-\delta_1)a_1+ \cdots+ (\epsilon_t-\delta_t)a_t \in [1,t] \hat{_{\pm}}  A,$$ as claimed.  

Second, we prove that for the set $$A=\{1,2,\dots,2^{t-2}\} \subseteq \mathbb{Z}_n,$$ we have $0 \not\in [1,t] \hat{_{\pm}}  A$.  Suppose that we have
$$\epsilon_0 , \epsilon_1 , \dots , \epsilon_{t-2} \in \{-1,0,1\}$$ for which
$$\epsilon_0 \cdot 1 + \epsilon_1 \cdot 2+ \cdots + \epsilon_{t-2} \cdot 2^{t-2}=0$$ in $\mathbb{Z}_n$.  Since
the left-hand side has absolute value at most $2^{t-1}-1<n$, the same equation must hold in $\mathbb{Z}$ as well.  Now if there is an index $0 \leq i \leq t-2$ for which $\epsilon_i \neq 0$, then let
$$k=\min\{i \mid \epsilon_i \neq 0\}.$$  Dividing our equation (in $\mathbb{Z}$) by $2^k$ we get
$$\epsilon_k + \epsilon_{k+1} \cdot 2+ \cdots + \epsilon_{t-2} \cdot 2^{t-2-k}=0,$$ which is a contradiction, since the left-hand side is odd.  Therefore, $\epsilon_i=0$  for all $0 \leq i \leq t-2$, proving our claim.
$\Box$

\addcontentsline{toc}{section}{Proof of Proposition \ref{max sum-free AP}}

\section*{Proof of Proposition \ref{max sum-free AP}} \label{proof of max sum-free AP}

Let $$A=\{a,a+d,\dots,a+(m-1)d\}$$ be a sum-free arithmetic progression in $\mathbb{Z}_n$ of size $$m=v_1(n,3)=\left\{
\begin{array}{ll}
\left(1+\frac{1}{p}\right) \frac{n}{3} & \mbox{if $n$ has prime divisors congruent to 2 mod 3,} \\ & \mbox{and $p$ is the smallest such divisor,}\\ \\
\left\lfloor \frac{n}{3} \right\rfloor & \mbox{otherwise.}\\
\end{array}\right.$$
Note that $m \geq (n-1)/3$ (since in the case when $n \equiv 2$ mod 3, $n$ must have a prime divisor $p \equiv 2$ mod 3).

Let $\delta=\gcd(d,n)$, $d=d_1 \delta$, $n=n_1 \delta$.  Since $d_1$ and $n_1$ are relatively prime, there are integers $b$ and $c$ for which
$$b \cdot d_1=1+c \cdot n_1.$$  Then $$b \cdot d=\delta + c \cdot n,$$ so
$$b \cdot A=\{b  a, b   a+\delta,\dots,b   a+(m-1)\delta\}.$$  Note that the fact that $A$ is sum-free implies that $b \cdot A$ is as well.

Observe that $b \cdot A$ has size $m$, since otherwise we would have an $i \in \{1,2,\dots,m-1\}$ for which $i \delta=0$ in $\mathbb{Z}_n$, but then $i d=0$ in $\mathbb{Z}_n$ for the same $i$, contradicting $|A|=m$.  Furthermore, $b \cdot A$ has size at most $n_1=n/\delta$, since it is contained in a coset of the subgroup of $G$ that has order $n_1$.  Thus $m \leq n_1$; but since $m \geq (n-1)/3$, we get $n/\delta \geq (n-1)/3$.  Therefore, either $\delta \leq 3$, or $\delta=4$ and $n=4$, but the latter case is impossible since for even $n$ we have $m=n/2$, contradicting $m \leq n_1$.

We consider the cases $\delta=3$, $\delta=2$, and $\delta=1$ separately.

When $\delta=3$, we have $n_1=n/3$, so $n/3 \geq m \geq (n-1)/3$ implies that $m=n/3$ and thus $b \cdot A$ is a full coset of the subgroup of $G$ that has order $n/3$.  This coset cannot be the subgroup itself, since then $0 \in A$, contradicting sum-freeness.  Therefore, 
$$b \cdot A=\{1,4,7, \dots, n-2\}$$ or $$b \cdot A=\{2,5,8, \dots, n-1\},$$ but note that in the second case we have 
$$(2b) \cdot A=\{1,4,7, \dots, n-2\}.$$
Since $m=v_1(n,3)=n/3$, we see that $n$ cannot have any prime divisors congruent to 2 mod 3; in particular, $n$ is odd.  
We also see that $b$ and $n$ are relatively prime, since $1 \in b \cdot A$ or $-1 \in b \cdot A$ . Since $n$ is odd, $2b$ and $n$ must be relatively prime too.  This yields case 2 (a).

If $\delta=2$, then $n$ is even and thus $m=v_1(n,3)=n/2$.  Therefore, $$b \cdot A = ba+\{0,2,\dots,n-2\}=\{1,3,\dots,n-1\}.$$   
Since then $b$ is relatively prime to $n$, it has an inverse, with which we get $$A = a+\{0,2,\dots,n-2\}=\{1,3,\dots,n-1\}.$$
This yields case 1.

Finally, assume that $\delta=1$.  In this case, 
$$b \cdot A=\{b  a, b   a+1,\dots,b   a+(m-1)\},$$ and so
$$2(b \cdot A)=b \cdot A+b \cdot A=\{2b  a, 2b   a+1,\dots,2b   a+2(m-1)\},$$ and thus
$$2(b \cdot A)-b \cdot A=\{ba-m+1, ba-m, ba-m+1,\dots, ba+2m-2\}.$$
Therefore, we see that $b \cdot A$ is sum-free, that is, $0 \not \in 2(b \cdot A)-b \cdot A$, exactly when there is an integer $k$ for which
$$kn+1 \leq ba-m+1$$ and $$ba+2m-2 \leq (k+1)n-1,$$  or, equivalently, 
$$m \leq ba-kn \leq n-2m+1.$$
Therefore, $m \leq (n+1)/3$.  We separate three subcases.

If $n \equiv 2$ mod 3, then $$m=v_1(n,3)=\left(1+\frac{1}{p}\right) \frac{n}{3},$$ where $p$ is the smallest prime divisor of $n$ with $p \equiv 2$ mod 3.
So $m \leq (n+1)/3$ is possible only when $$m=(n+1)/3=\left(1+\frac{1}{p}\right) \frac{n}{3},$$ which gives $n=p$.  So $m=(p+1)/3$, $$(p+1)/3 \leq ba-kp \leq n-2(p+1)/3+1,$$ so
$$ba-kp=(p+1)/3,$$ and thus $ba=(p+1)/3$ in $\mathbb{Z}_p$.  This yields case 3.

If $3|n$, then $m \leq (n+1)/3$ implies that $m \leq n/3$, but we also know that $m \geq (n-1)/3$, so $m=n/3$.  Therefore, $n$ cannot have a prime divisor congruent to 2 mod 3.  We have $$n/3 \leq ba-kn \leq n-2n/3+1,$$ so $ba-kn \in \{n/3, n/3+1\}$ and thus $ba=n/3$ or $ba=n/3+1$ in $\mathbb{Z}_n$, with which
$$b \cdot A=\{n/3, n/3+1, \dots, 2n/3-1\}$$ or
$$b \cdot A=\{n/3+1, n/3+2, \dots, 2n/3\}.$$  Note that 
$$\{n/3, n/3+1, \dots, 2n/3-1\}=-\{n/3+1, n/3+2, \dots, 2n/3\},$$ so the second possibility is superfluous.
This yields case 2 (b).

Finally, suppose that $n \equiv 1$ mod 3.  Then $m \leq (n+1)/3$ implies that $m \leq (n-1)/3$, but we also know that $m \geq (n-1)/3$, so $m=(n-1)/3$.
Therefore, $n$ cannot have a prime divisor congruent to 2 mod 3.  We have $$(n-1)/3 \leq ba-kn \leq n-2(n-1)/3+1,$$ so $$ba-kn \in \{(n-1)/3, (n-1)/3+1, (n-1)/3+2\}$$ and thus $ba=(n-1)/3$, $ba=(n-1)/3+1$, or $ba=(n-1)/3+2$ in $\mathbb{Z}_n$, with which
$$b \cdot A=\{(n-1)/3, (n-1)/3+1, \dots, 2(n-1)/3-1\},$$ 
$$b \cdot A=\{(n-1)/3+1, (n-1)/3+2, \dots, 2(n-1)/3\},$$ or
$$b \cdot A=\{(n-1)/3+2, (n-1)/3+3, \dots, 2(n-1)/3+1\}.$$  Note that 
$$\{(n-1)/3, (n-1)/3+1, \dots, 2(n-1)/3-1\}=-\{(n-1)/3+2, (n-1)/3+3, \dots, 2(n-1)/3+1\},$$ so the third possibility is superfluous.
This yields case 4.  $\Box$

\addcontentsline{toc}{section}{Proof of Theorem \ref{explicit plagne AP}}

\section*{Proof of Theorem \ref{explicit plagne AP}} \label{proof of explicit plagne AP}

Let $A$ be a $(k,l)$-sum-free set of size $M+1$ in $\mathbb{Z}_p$.  According to Theorem \ref{plagne AP}, $A$ is an arithmetic progression (of common difference $d$); since $p$ is prime, (multiplying by the inverse of $d$ in $\mathbb{Z}_p$) we may assume that $$A=\{a, a+1, \dots, a+M\}$$ for some $a \in \mathbb{Z}_p$.  Then
$$0 \not \in kA-lA=\{(k-l)a-lM, (k-l)a-lM+1, \dots, (k-l)a+kM\},$$ so there is an integer $c$ for which
$$cp+1 \leq (k-l)a-lM$$ and $$(k-l)a+kM \leq (c+1)p-1.$$ Combining these two inequalities we get $$lM+1 \leq (k-l)a-cp \leq p-kM-1$$ or, equivalently,
$$(k-l)a=lM+i$$ in $\mathbb{Z}_p$ for some $$1 \leq i \leq (p-kM-1)-(lM+1)+1=r+1.$$
Thus $A$ is (a dilate of) one of the $r+1$ sets
$$\{a_i,a_i+1,\dots,a_i+M\}.$$  Due to symmetry, we don't need $\lfloor (r+1)/2 \rfloor$ of these choices; more precisely, for any $i$ with $1 \leq i \leq \lfloor (r+1)/2 \rfloor$, $A_i=-A_{r+2-i}$, since $$a_i=-(a_{r+2-i}+M).$$  To see this, note that this equation holds if, and only if, 
$$(k-l)a_i=-(k-l)(a_{r+2-i}+M).$$ This equation holds, since
$$-(k-l)(a_{r+2-i}+M)=-(lM+r+2-i)-(k-l)M=lM+i=(k-l)a_i$$ in $\mathbb{Z}_p$.  Therefore, it suffices to assume that $i \leq r+1-\lfloor (r+1)/2 \rfloor=\lfloor r/2 \rfloor +1$.  $\Box$

\addcontentsline{toc}{section}{Proof of Proposition \ref{max size interval weak sumfree}}

\section*{Proof of Proposition \ref{max size interval weak sumfree}} \label{proof of max size interval weak sumfree}

Let $$A=[a,a+m-1]=\{a,a+1,\dots,a+m-1\}$$ be an interval of length $m-1$ (size $m$) in $\mathbb{Z}_n$. 
Note that if $A$ is an interval, then so is $h \hat{\;} A $ for all integers $h$ with $1 \leq h \leq m$; in particular,
$$h \hat{\;} A=[ha+h(h-1)/2, h(a+m-1)-h(h-1)/2].$$
Therefore, for positive integers $k$ and $l$ with $l < k \leq M+1$, $k \hat{\;} A - l \hat{\;} A$ is also an interval: its ``smallest'' element is 
$$ka+k(k-1)/2-l(a+m-1)+l(l-1)/2,$$ and its ``largest'' element equals 
$$k(a+m-1)-k(k-1)/2-la-l(l-1)/2;$$  using the notation $$J=k(k-1)+l(l-1)=k^2+l^2-k-l,$$
we have 
$$k \hat{\;} A - l \hat{\;} A=[ka-l(a+m-1)+J/2, k(a+m-1)-la-J/2].$$
Now $A$ is a weak $(k,l)$-sum-free set in $\mathbb{Z}_n$ if, and only if, $0 \not \in k \hat{\;} A - l \hat{\;} A$, and this occurs if, and only if, there is an integer $b$ for which
$$ka-l(a+m-1)+J/2 \geq bn+1$$ and
$$k(a+m-1)-la-J/2 \leq (b+1)n-1,$$ or, equivalently, 
$$l(m-1)-J/2+1 \leq (k-l)a -bn \leq n-k(m-1)+J/2-1.$$  
We can divide by $\delta=\gcd(n,k-l)$ and write
$$\frac{l(m-1)-J/2+1}{ \delta} \leq \frac{k-l}{\delta} a -\frac{n}{\delta} b \leq \frac{ n-k(m-1)+J/2-1 }{ \delta}.$$
Note that $(k-l)/\delta$ and $n/\delta$ are relatively prime integers, so any integer can be written as their linear combination (with integer coefficients). 
We thus see that $\mathbb{Z}_n$ contains a weak $(k,l)$-sum-free interval of size $m$ if, and only if, there is an integer $C$ for which
$$\frac{l(m-1)-J/2+1}{ \delta} \leq C \leq \frac{ n-k(m-1)+J/2-1 }{ \delta}.$$

Therefore, a necessary condition for the existence of a weak $(k,l)$-sum-free interval of size $m$ is that 
$$\frac{l(m-1)-J/2+1}{ \delta} \leq \frac{ n-k(m-1)+J/2-1 }{ \delta},$$ from which
$$m \leq \left \lfloor \frac{n+J-2}{k+l} \right \rfloor +1=M+1$$ follows.  On the other hand, a sufficient condition for the existence of a weak $(k,l)$-sum-free interval of size $m$ is that 
$$\frac{l(m-1)-J/2+1}{ \delta} +1 \leq \frac{ n-k(m-1)+J/2-1 }{ \delta},$$ for which it suffices that $m \leq M$ since then
$$m \leq  \frac{n+J-2}{k+l}  \leq \frac{n+J-2-\delta}{k+l} +1,$$ from which our condition follows.
Therefore, $\gamma \hat{\;} (\mathbb{Z}_n,\{k,l\})$ equals either $M$ or $M+1$. 

To see when $\gamma \hat{\;} (\mathbb{Z}_n,\{k,l\})=M+1$, note that the existence 
of an integer $C$ for which
$$\frac{L}{ \delta} = \frac{lM-J/2+1}{ \delta} \leq C \leq \frac{ n-kM+J/2-1 }{ \delta} = \frac{n-K}{ \delta}$$ is equivalent to saying that
$$\frac{L}{ \delta} \leq \left \lfloor  \frac{n-K}{ \delta} \right \rfloor.$$
This completes our proof.  $\Box$

\addcontentsline{toc}{section}{Proof of Corollary \ref{kllower}}

\section*{Proof of Corollary \ref{kllower}} \label{proofofkllower}

We provide two proofs: one using Proposition \ref{max size interval weak sumfree} and another that is more direct.

{\em Proof I}: Here we use Proposition \ref{max size interval weak sumfree}; it suffices to prove that $$\gamma \hat{\;} (\mathbb{Z}_n,\{k,l\}) \geq 
\left \lfloor \frac{n+k^2+l^2-\gcd(n,k-l)-1}{k+l} \right \rfloor.$$

Recall the notations 
$$\delta=\gcd(n,k-l),$$
$$J = k^2+l^2-(k+l),$$ 
$$M= \lfloor (n+J-2)/(k+l) \rfloor ,$$
$$K= kM-J/2+1, $$
$$L= lM-J/2+1.$$
Since $$M+1=\left \lfloor \frac{n+k^2+l^2-2}{k+l} \right \rfloor \geq \left \lfloor \frac{n+k^2+l^2-\gcd(n,k-l)-1}{k+l} \right \rfloor,$$ by Proposition \ref{max size interval weak sumfree}
we may assume that $$L/\delta > \lfloor (n-K)/\delta \rfloor,$$ and we need to show that
$$M \geq \left \lfloor \frac{n+k^2+l^2-\gcd(n,k-l)-1}{k+l} \right \rfloor$$ in this case.  

Let $R$ be the remainder of $n+k^2+l^2-2$ when divided by $k+l$; we can then write
$$n+k^2+l^2-2 = (k+l)(M+1) +R,$$
and so
$$n-K=(k+l)(M+1)+R-(k^2+l^2-2)-(kM-J/2+1)=L+R.$$
Therefore, our hypothesis that 
$$L/\delta > \lfloor (n-K)/\delta \rfloor$$
can only hold if $R < \delta -1$.
But then
$$\left \lfloor \frac{n+k^2+l^2-\gcd(n,k-l)-1}{k+l} \right \rfloor=M+1+\left \lfloor \frac{R-\delta+1}{k+l} \right \rfloor \leq M+1 + \left \lfloor \frac{-1}{k+l} \right \rfloor=M,$$ as claimed.  $\Box$

{\em Proof II}: 
We set $\delta=\mathrm{gcd}(n,k-l)$ and $$c=\left \lfloor \frac{n+k^2+l^2-\delta-1}{k+l} \right \rfloor.$$  We need to prove that $\mu \hat{\;} (\mathbb{Z}_n, \{k,l\}) \geq c$.  

Note that $c<n$ since
$$c \leq  \frac{n+k^2+l^2-2}{k+l} = n - \frac{n(k+l-1)-(k^2+l^2-2)}{k+l},$$ and, since $0 < l < k \leq n$, we have
\begin{eqnarray*}
n(k+l-1)-(k^2+l^2-2) & \geq & k(k+l-1)-(k^2+l^2-2) \\
& = & k(l-1)-(l^2-2) \\
& > & k(l-1)-(l^2-1) \\
& = & (k-l-1)(l-1) \\
& \geq & 0.
\end{eqnarray*}

Now $$\overline{h}=\frac{k^2+l^2-k+l-2}{2}$$ is an integer; let $r$ be the remainder of $\overline{h}-lc$ when divided by $\delta$; we then have $0 \leq r \leq \delta-1$.
 
Since $\delta$ is a divisor of $\overline{h}-lc-r$, we can find integers $a$ and $b$ for which $$a \cdot (k-l) - b \cdot n =-(\overline{h}-lc-r);$$ furthermore, we may assume w.l.o.g. that $b \geq 0$.    

Now define the set $A$ as $$A=\{a,a+1,\dots,a+c-1\}.$$  (Throughout this proof, we will talk about integers and elements of the group $\mathbb{Z}_n$ interchangeably: when an integer is considered as an element of $\mathbb{Z}_n$, we regard it mod $n$; on the other hand, an element $t$ of $\mathbb{Z}_n$ is treated as the integer $t$ in $\mathbb{Z}$.)  

Clearly, $|A|=c$; we will show that $A$ is a weak $(k,l)$-sum-free set in $\mathbb{Z}_n$ by showing that every element of $k \hat{\;} A - l\hat{\;} A$ is strictly between $bn$ and $(b+1)n$.

The smallest element of $k \hat{\;} A - l\hat{\;} A$ is 
\begin{eqnarray*}
ka+\frac{k(k-1)}{2}-la-lc+\frac{l(l+1)}{2} & = & (k-l) a-lc+\overline{h}+1 \\
& = & b n-\overline{h}+r+\overline{h}+1 \\
& = & bn+1+r \\
& \geq & bn+1.
\end{eqnarray*}

For the largest element of $k \hat{\;} A - l\hat{\;} A$ we have 
\begin{eqnarray*}
k a+k c-\frac{k(k+1)}{2} -la- \frac{l(l-1)}{2}& = &  (k-l) a +kc - \frac{k^2+l^2+k-l}{2} \\
& = & b n+lc -\overline{h}+r +kc - \frac{k^2+l^2+k-l}{2}\\
& = & b n+(k+l)c - k^2-l^2+1+ r\\
& \leq & b n+(n+k^2+l^2-\delta-1) - k^2-l^2+1+ r \\
&= & (b+1) n  - \delta +r \\
& \leq & (b+1) n  - 1.
\end{eqnarray*}  
Therefore, no element of $k \hat{\;} A - l\hat{\;} A$ is 0 in $\mathbb{Z}_n$.  $\Box$

\addcontentsline{toc}{section}{Proof of Theorem \ref{M21}}

\section*{Proof of Theorem \ref{M21}} \label{proof of M21}

Let $A$ be a weak sum-free set in $\mathbb{Z}_n$.  We claim that
$$|A| \leq \left\{
\begin{array}{ll}
\left(1+\frac{1}{p}\right) \frac{n}{3} & \mbox{if $n$ has prime divisors congruent to 2 mod 3,} \\ & \mbox{and $p$ is the smallest such divisor,}\\ \\
\left\lfloor \frac{n}{3} \right\rfloor +1 & \mbox{otherwise.}\\
\end{array}\right.$$

If $A=\{0\}$, the claim is obviously true.  If $A$ contains a nonzero element $a$, then it cannot contain 0, since otherwise $0+a=a$ would contradict $A$ being weakly sum-free.  So we may assume that $0 \not \in A$.

By Theorem \ref{Diananda and Yap}, the claim is clear when $A$ is sum-free, so we may assume that there is an element $a_0 \in A$ for which $2a_0 \in A$.  We have $a_0\neq 0$ and thus $a_0 \neq 2a_0$.

Let $$A_1=a_0+(A \setminus \{a_0\})=\{a_0 + a \mid a \in A \setminus \{a_0\} \}$$ and $$A_2=2a_0+(A \setminus \{a_0,2a_0\})=\{2a_0 + a \mid a \in A \setminus \{a_0,2a_0\} \}.$$
Note that $A_1 \subseteq 2 \hat{\;} A$ and $A_2 \subseteq 2 \hat{\;} A$, so $A$ is disjoint from both $A_1$ and $A_2$.  Furthermore, $A_1$ and $A_2$ are disjoint too, since otherwise we would have elements $a_1 \in A \setminus \{a_0\}$ and $a_2 \in A \setminus \{a_0,2a_0\}$ for which
$$a_0+a_1=2a_0+a_2,$$ but then $$a_1=a_0+a_2,$$ contradicting that $A$ and $A_1$ are disjoint.  Since $A$, $A_1$, and $A_2$ are pairwise disjoint, we have
$$|A|+|A_1|+|A_2| = 3|A|-3 \leq n,$$ from which $|A| \leq \lfloor n/3 \rfloor +1$, proving our claim.  $\Box$

\addcontentsline{toc}{section}{Proof of Proposition \ref{Hallfors prop 1}}

\section*{Proof of Proposition \ref{Hallfors prop 1}} \label{proof of Hallfors prop 1}

Since $$A=\{a,a+1,\dots,a+c\}+H$$ consists of $c+1$ consecutive cosets of $H$, and each coset is of size $k-1$, we have
$$k \hat {\;} A =\{ka+1, ka+2, \dots, ka+kc-1\}+H,$$ and
$$l \hat {\;} A \subseteq \{la, la+1, \dots, la+lc\}+H,$$with equality holding if $l \leq k-2$; however, for $l=k-1$, only a single element from each of the cosets $la+H$ and $la+lc+H$ is in $l \hat {\;} A$.

Using our specified values of $a$ and $c$, we find that
$$la+lc=\frac{l^2n}{(k-1)(k^2-l^2)} + \frac{ln}{(k-1)(k+l)}=\frac{kln}{(k-1)(k^2-l^2)}=ka,$$ and 
\begin{eqnarray*}
ka+kc-1 & = & \frac{kln}{(k-1)(k^2-l^2)} + \frac{kn}{(k-1)(k+l)}-1 \\
& = & \frac{k^2n}{(k-1)(k^2-l^2)}-1\\
& = & \frac{n}{k-1}+\frac{l^2n}{(k-1)(k^2-l^2)}-1\\
& = & \frac{n}{k-1}+la-1.
\end{eqnarray*}
Therefore, $k \hat {\;} A \cap l \hat {\;} A =\emptyset$, as claimed.  We also see that, in the case of $l \leq k-2$, we have $k \hat {\;} A \cup  l \hat {\;} A =\mathbb{Z}_n$, so $A$ is a complete weak $(k,l)$-sum-free set in $\mathbb{Z}_n$.  $\Box$

\addcontentsline{toc}{section}{Proof of Proposition \ref{Hallfors prop 2}}

\section*{Proof of Proposition \ref{Hallfors prop 2}} \label{proof of Hallfors prop 2}

Using our notations, for any positive integer $h$, we have:
$$h \hat {\;} A \subseteq \{h_1, h_1+1, \dots, h_2\}+H,$$ where $$h_1=\max\{0, h- d/2+1\},$$ and $$h_2=\min\{h, n/d\}.$$ 
Therefore,
$$k \hat {\;} A \subseteq \{k-d/2+1, k-d/2+2, \dots, k\}+H,$$ and 
$$l \hat {\;} A \subseteq \{l-d/2+1, l-d/2+2, \dots, l\}+H.$$
So both $k \hat {\;} A$ and $l \hat {\;} A$ are within $d/2$ consecutive cosets of $H$; furthermore, these $d$ cosets are all distinct, since $$(k-d/2)+H=l+H.$$  This proves that $k \hat {\;} A$ and $l \hat {\;} A$ are disjoint.

Furthermore, we see that $$k \hat {\;} A \cup l \hat {\;} A = \mathbb{Z}_n$$ if, and only if, the inequalities $k < n/d$, $k > d/2-1$, $l < n/d$, and $l > d/2-1$ all hold; so $A$ is a complete weak $(k,l)$-sum-free set in $\mathbb{Z}_n$ if, and only if,   
$$d/2-1 < l < k < n/d.$$
$\Box$

\addcontentsline{toc}{section}{Proof of Proposition \ref{all n for weak sumfree m=k, k+1}}

\section*{Proof of Proposition \ref{all n for weak sumfree m=k, k+1}} \label{proof of all n for weak sumfree m=k, k+1}

By Proposition \ref{all n for weak sumfree trivial}, we may assume that $n \geq k+1$.  For $i=0,1,2,\dots,k$, consider the subsets 
$$A_i=\{0,1,2,\dots,k\} \setminus \{i\}$$ of $\mathbb{Z}_n$.  

Consider first the case when $n=k+1$ and $k \equiv 1$ mod 4.  Note that every $k$-subset of $\mathbb{Z}_n$ is of the form $A_i$ for some $i=0,1,2,\dots,k$; to prove our claim, we will show that $$s_i =k \hat {\;} A=0+1+2+\cdots+k-i \in A_i$$ for each $i$.  Indeed, if for some $i$ we were to have   
$$s_i=k(k+1)/2-i=i$$ in $\mathbb{Z}_n$, then $$k(k+1)/2=2i,$$ which is impossible, since $n$ is even, $2i$ is even, but $k(k+1)/2$ is odd.

We will now show that in all other cases, $\mathbb{Z}_n$ contains a $k$-subset that is weakly $(k,1)$-sum-free.  Assuming that there is no $i$ for which $A_i$ is weakly $(k,1)$-sum-free yields that, for each $i$, 
$s_i \in A_i$, so $s_i \in \{0,1,\dots, k\}$ and thus $$\{s_k,s_{k-1},\dots, s_1, s_0\} \subseteq \{0,1,2,\dots,k\}.$$  Observe that $$(s_k,s_{k-1},\dots, s_1, s_0)=(s_k, s_k+1, \dots, s_k+k),$$ so we must have $s_i=k-i$ for each $i$.  In particular, $$s_0=k(k+1)/2=k.$$

We consider three cases.  Assume first that $n \geq k+3$.  Then for
$$B=\{0,1,2,\dots,k, k+1,k+2\} \setminus \{1,k, k+1\}=A_0 \cup \{0,k+2\} \setminus \{1,k\}$$ we have
$$k \hat {\;} B=s_0+(k+2)-1-k=k+1.$$ Since $k+1 \not \in B$, $B$ is weakly $(k,1)$-sum-free in $\mathbb{Z}_n$. 

Assume now that $n=k+2$.  In this case, we have
$$6=(k^2-k)-(k+2)(k-3)=k^2-k=2 \cdot k(k+1)/2 -2k=2 s_0-2s_0=0.$$  This can only happen in $\mathbb{Z}_{k+2}$ if $k+2=6$ and thus $k=4$.  However, in this case, $$A_2=\{0,1,3,4\}$$ is weakly $(4,1)$-sum-free in $\mathbb{Z}_6$.

Next, assume that $n=k+1$ and that $k$ is even.  As we have seen above, we have $$s_{k/2}=k-k/2=k/2,$$ so $A_{k/2}$ is weakly $(k,1)$-sum-free in $\mathbb{Z}_n$. 

Finally, assume that $n=k+1$ and $k \equiv 3$ mod 4.  In this case, we have $$s_{(k+1)/4} = k(k+1)/2-(k+1)/4=(k+1) \cdot (k-1)/2 +(k+1)/4=(k+1)/4,$$ which means that $A_{(k+1)/4}$ is weakly $(k,1)$-sum-free in $\mathbb{Z}_n$.

To prove our second claim, first note that if $n \leq 2k+1$ and $|A|=k+1$, then $A$ and $k \hat{\;} A$ cannot be disjoint, since $|k \hat{\;} A|=k+1$.  Suppose then that $n \geq 2k+2$.  We will consider two cases depending on the parity of $k$.

When $k$ is even, we set
$$A=\{1\} \cup \{\pm 2, \pm 4, \dots, \pm k\}.$$
Then $|A|=k+1$, and 
$$k \hat{\;} A = 1 - A = \{0\} \cup \{-1\} \cup \{\pm 3, \pm 5, \dots, \pm (k-1)\} \cup \{k+1\}.$$
Since $k+1 < n-k$, $A$ is weakly $(k,1)$-sum-free in $\mathbb{Z}_n$.

When $k$ is odd, we let
$$A=\{0,-1,-2,4\} \cup \{\pm 5, \pm 7, \dots, \pm k\}.$$  (We simply let $A=\{0,-1,-2,4\}$ if $k=3$.)
Again, $|A|=k+1$, and 
$$k \hat{\;} A = 1 - A =\{1,2, \pm 3, -4\} \cup \{\pm 6, \pm 8, \dots, \pm (k-1)\} \cup \{k+1\},$$
and $A$ is weakly $(k,1)$-sum-free in $\mathbb{Z}_n$.

We can also observe that if $n=2k+2$, then $A$ is a complete weak $(k,1)$-sum-free set.  $\Box$

\backmatter

\newpage
\pagebreak

\addcontentsline{toc}{chapter}{Author Index}

\printindex

\end{document}